\documentclass[11pt,twoside,reqno]{report}
\usepackage[ngerman, american]{babel}
\usepackage{amsmath}
\usepackage{amssymb,amsfonts,amsthm,amscd}
\usepackage{exscale}
\usepackage[DIV13, headinclude]{typearea}
\usepackage{graphicx}
\usepackage{mathrsfs}
\usepackage{local}

\usepackage{url}

\begin{document}

\begin{titlepage}
\title{On Rho invariants of fiber bundles}
\author{Michael Bohn\\[12pt] \small\it Department of Mathematics, Bonn
University\\ \small\it Endenicher Allee 60, D-53115 Bonn, Germany\\ \small
\textrm{e-mail: mbohn@math.uni-bonn.de}\\[3ex]
\small Dissertation \\
\small published under\\ \small \url{http://hss.ulb.uni-bonn.de/diss_online/}
}
\date{Bonn, 2009}
\maketitle
\end{titlepage}

\pagenumbering{roman}\setcounter{page}{1}
\selectlanguage{american}

\cleardoublepage

\chapter*{Summary}
\addcontentsline{toc}{chapter}{Summary}

The content of this thesis is a detailed investigation of Rho
invariants of the total spaces of fiber bundles. The main idea is
to use adiabatic limits of Eta invariants to obtain a formula for
Rho invariants that separates the contribution coming from the
fiber and the one coming from the base. An adiabatic metric on a
fiber bundle rescales the metric of the base manifold in such a
way that the geometry of the fiber bundle approaches a product
situation. Concerning the Eta invariant, this process has
received a far-reaching treatment in the literature. For this
reason, one concern of this thesis is to formulate the technical
aspects of local index theory for families of Dirac operator in
terms of the odd signature operator, and place known results in a
context which permits the treatment of Rho invariants.

The resulting formula expresses the Rho invariant as a sum of
three terms, each of which is of a very different nature. First
of all, a higher dimensional analog of the Rho invariants of the
fiber has to be integrated over the base. This term is of a local
nature on the base, but contains global spectral information
about the fiber. The next term is essentially a Rho invariant of
the base, where the underlying flat connection is defined on the
bundle of cohomology groups of the fiber. Lastly, there is a
purely topological term, which can be computed from the spectral
sequence of the fiber bundle. Together, this formula casts the
Rho invariant of the total space into a form which incorporates
the structure of the fiber bundle in a satisfactory way.

The main concern of this thesis is, however, to use this
theoretical formula to compute Rho invariants for explicit
classes of fibered 3-manifolds. More precisely, we consider
principal $S^1$-bundles over closed, oriented surfaces as well as
mapping tori with fiber a closed, oriented surface. For the first
class of examples, one can compute $\U(1)$-Rho invariants without
using this general formula. In particular, this yields the
opportunity to compare the different approaches and test the
systematical advantage of the general formula.

For 3-dimensional mapping tori, the presented theory can also be
used for explicit computations. We first consider the case that
the monodromy map is of finite order. In this case, a general
formula for Rho invariants can be derived. To investigate a
further interesting class of mapping tori, we consider $\U(1)$-Rho
invariants in the case that the fiber is a 2-dimensional torus.
Here, hyperbolic monodromy maps deserve particular attention. When
discussing them, the logarithm of a generalized Dedekind Eta
function naturally appears. A satisfactory formula for $\U(1)$-Rho
invariants of hyperbolic mapping tori can then be deduced from a
transformation formula for these Eta functions.

\pagestyle{myheadings}
\markboth{\textsc{Contents}}{\textsc{Contents}} \tableofcontents

\renewcommand{\thechapter}{\arabic{chapter}}
\renewcommand{\thesection}{\arabic{chapter}.\arabic{section}}
\renewcommand{\thesubsection}{\arabic{chapter}.\arabic{section}.\arabic{subsection}}

\numberwithin{equation}{chapter}

\markboth{\textsc{Introduction}}{\textsc{Introduction}}
\renewcommand{\thetheorem}{\arabic{theorem}}
\renewcommand{\theequation}{\arabic{equation}}

\chapter*{Introduction}

\addcontentsline{toc}{chapter}{Introduction}

\noindent\textbf{The Rho Invariant for Closed Manifolds.} In their
famous series of articles \cite{APS1,APS2, APS3}, Atiyah, Patodi
and Singer established an index theorem for manifolds with
boundary. Part of their motivation was to find a generalization
of Hirzebruch's Signature Theorem to manifolds with boundary and
give a differential geometric explanation for the signature
defect.

We recall this briefly. Let $W$ be a closed, oriented 4-manifold,
and let $\Sign(W)$ be its signature. Then the Hirzebruch's
Signature Theorem states that
\begin{equation}\label{SignIntro}
\Sign(W) = \frac 13 \int_W p_1(TW).
\end{equation}
Here, $p_1(TW)$ is the first Pontrjagin form, and since $W$ is
closed, it is a characteristic class independent of the connection
used to compute it. Let us now assume that $W$ has a boundary
$M$. If $\ga:\pi_1(W)\to\U(k)$ is a unitary representation of the
fundamental group, one defines a twisted signature $\Sign_\ga(W)$
using cohomology groups with local coefficients. Then an
application of the signature formula \eqref{SignIntro} and its
twisted version to the closed double $W\cup_M -W$ shows that the
difference
\[
\Sign_\ga(W)- k\cdot\Sign(W)
\]
depends only on the topology of the boundary $\pd W = M$ as well
as the restriction of $\ga$ to $\pi_1(M)$. This signature defect
is in general non-trivial. Moreover, the signature itself fails
to be multiplicative under finite coverings of manifolds with
boundary. Both observations show that the signature of a manifold
with boundary is in general not expressible as in
\eqref{SignIntro}.

The Atiyah-Patodi-Singer Index Theorem for manifolds with
boundary identifies the correction term in great generality. For
a formally self-adjoint elliptic differential operator $D$ of
first order, acting on sections of a vector bundle over a closed
manifold $M$, one defines the Eta function
\begin{equation}
\eta(D,s):= \sum_{0\neq \gl\in\spec(D)}
\frac{\sgn(\gl)}{|\gl|^s},\quad \Re(s)\text{ large}.
\end{equation}
The function $\eta(D,s)$ admits a meromorphic continuation to the
whole $s$-plane, and it is a remarkable fact that $s=0$ is not a
pole. The \emph{Eta invariant} $\eta(D)$ is defined as this
finite value. Then a special case of the Atiyah-Patodi-Singer
Index Theorem for a compact, oriented 4-manifold $W$ with
boundary $M$ is
\begin{equation}\label{APSIntro}
\Sign_\ga(W) = \frac k 3 \int_W p_1(TW,\nabla^g) -
\eta(B_A^{\ev}).
\end{equation}
Here, $p_1(TW,\nabla^g)$ is the first Pontrjagin form, computed
with respect to a metric $g$ in product form near the boundary,
$A$ is a flat $\U(k)$-connection over $M$ whose holonomy
coincides with $\ga|_{\pi_1(M)}$, and---most
importantly---$B_A^{\ev}$ is the \emph{odd signature operator} on
$M$. It is defined on differential forms of even degree with
values in the flat bundle $E$ as
\[
B_A^{\ev}\go = (-1)^{p}(*d_A-d_A*)\go, \quad \go\in \gO^{2p}(M,E).
\]
A simple consequence of \eqref{APSIntro} is that
\begin{equation}\label{SignDefIntro}
\Sign_\ga(W)- k\cdot\Sign(W) =\eta(B_A^{\ev})-k\cdot\eta(B^{\ev}).
\end{equation}
The right hand side of \eqref{SignDefIntro} is called the
\emph{Rho invariant} $\rho_A(M)$. It is defined for every odd
dimensional manifold $M$ and flat unitary connection $A$ over $M$,
without any reference to a bounding manifold. Moreover,
$\rho_A(M)$ turns out to be independent of the choice of the
metric used in defining the involved odd signature operators.
Therefore, it is an intrinsically defined smooth invariant of
$M$, which extends the signature defect.

Since the Eta invariants appearing in the definition of
$\rho_A(M)$ are non-local spectral invariants, it is difficult to
compute Rho invariants directly without using property
\eqref{SignDefIntro}. However, one cannot always find a bounding
manifold in such a way that the flat connection extends.
Therefore, intrinsic methods to compute Rho invariants are in
demand. The concern of this thesis is to investigate an intrinsic
approach to this problem in the case that the manifold $M$ is the
total space of an oriented fiber bundle of closed
manifolds.\\

\noindent\textbf{Adiabatic Limits of Eta Invariants.} In a
remarkable paper of Witten \cite{Wit85}, Eta invariants appeared
in the interpretation of anomalies in physics. Associated to a
family of Dirac operators is a determinant line bundle over the
parameter space, first described by in \cite{Qui85a} by Quillen.
It comes equipped with a natural connection, defined in terms of
the Ray-Singer analytic torsion \cite{RS71}. The topology of this
line bundle encodes the obstruction to defining in a consistent
way a regularized determinant associated to the family of Dirac
operators. In the physicists terminology, there is no ``local
anomaly'' if Quillen's connection on the determinant bundle is
flat. However, the bundle might still not be trivial, and this
``global anomaly'' is encoded in the holonomy of Quillen's
connection. Witten suggests an interpretation of this holonomy
using the Eta invariant. For example, a family of spin Dirac
operators is naturally associated to a fiber bundle over the
parameter space $B$ whose fiber $F$ is a closed spin manifold.
Here, we also restrict to the case that $F$ is even dimensional.
Pulling back this structure using a closed loop $c:S^1\to B$ leads
to a fiber bundle $F\hookrightarrow M\to S^1$. Now Witten
considers an \emph{adiabatic metric}, that is, a family of
submersion metrics of the form
\begin{equation}\label{AdiabMetIntro}
g_\eps = \frac {g_{S^1}}{\eps^2}\oplus g_v,
\end{equation}
where $g_v$ is a metric on the vertical tangent bundle of $M$.
Then, if $D_\eps$ denotes the Dirac operator associated to
\eqref{AdiabMetIntro} on the total space $M$, Witten suggests that
the holonomy of Quillen's connection around the loop $c:S^1\to B$
is given by
\[
\lim_{\eps\to 0}\exp\big(2\pi i \eta(D_\eps)\big).
\]
The mathematical treatment of this holonomy theorem is due to
Bismut-Freed \cite{BF2} and Cheeger \cite{Che87}.

Motivated by this geometric interpretation, Bismut and Cheeger
\cite{BC89} gave a formula for the \emph{adiabatic limit} of the
Eta invariant, and generalized it to fiber bundles
$F\hookrightarrow M\to B$ with higher dimensional base spaces.
Using ideas of Bismut's local index theory for families
\cite{B86}, they construct a differential form $\widehat \eta$ on
$B$, whose value at each point $x\in B$ depends only on global
information of the fiber over $x$. Under the assumption that the
fiberwise Dirac operator is invertible, they prove that
\begin{equation}\label{BCIntro}
\lim_{\eps\to 0} \eta(D_\eps) = \int_B \widehat
A(TB,\nabla^{g_B}) \widehat \eta,
\end{equation}
where $A(TB,\nabla^{g_B})$ is the Hirzebruch $\widehat A$-form of
$B$. Very roughly, a consequence of \eqref{BCIntro} is that the
adiabatic limit of the Eta invariant is a simpler object than the
Eta invariant itself, since it is local on the base. The main
matter of this thesis is to analyze in which way this effect can
be used to simplify explicit computations of Rho invariants of
fiber bundles. \\

\noindent\textbf{Dai's Adiabatic Limit Formula.} The Rho
invariants we are considering are associated to the odd signature
operator. Here, the kernel of the vertical operator forms a
vector bundle $\sH_v^\bullet(M)\to B$ whose fiber over each point
is isomorphic to the cohomology of $F$. Therefore, the
invertibility hypothesis leading to \eqref{BCIntro} is too
restrictive. Fortunately, the result of Bismut and Cheeger has
been generalized by Dai \cite{Dai91} to a setting which applies
in particular to the case we are interested in. The bundle
$\sH_v^\bullet(M)\to B$ of vertical cohomology groups is endowed
with a natural flat connection $\nabla^{\sH_v}$. Using this, one
associates a twisted odd signature operator $D_B\otimes
\nabla^{\sH_v}$ over the base. Then Dai proves the following, very
remarkable adiabatic limit formula.

\begin{theorem}[Dai]\label{DaiIntro}
Let $F\hookrightarrow M\to B$ be an oriented fiber bundle of
closed manifolds with odd dimensional total space, endowed with a
submersion metric. Let $B_\eps^{\ev}$ be the family of odd
signature operators on $M$ associated to an adiabatic metric. Then
\[
\lim_{\eps\to 0}\eta(B_{\eps}^{\ev}) = 2^{[\frac{b+1}2]}\int_B
\widehat L(TB,\nabla^{g_B})\wedge \widehat \eta + \lfrac 12
\eta\big(D_B\otimes \nabla^{\sH_v}\big) + \gs.
\]
Here, $b$ is the dimension of $B$, and the differential forms
$\widehat L(TB,\nabla^B)$ and $\widehat \eta$ are the Hirzebruch
$\widehat L$-form and the Eta form of Bismut-Cheeger,
respectively. Moreover, $\gs$ is a topological invariant computed
from the Leray-Serre spectral sequence.
\end{theorem}
We will give more details on the terms appearing here in the main
body of this thesis. However, we already want to stress that
$\eta\big(D_B\otimes \nabla^{\sH_v}\big)$ and $\gs$ are of a very
different nature than the integral of the Eta form. Where the
latter is local on the base and contains spectral information
about the fiber, the twisted Eta term is a spectral invariant of
the base which contains cohomological information of the fiber.
Moreover, since $\gs$ arises from the Leray-Serre spectral
sequence it is purely cohomological. In this respect, Theorem
\ref{DaiIntro} is a very satisfactory decomposition of the
adiabatic limit into contributions coming from the base and the
fiber, respectively.\\

\noindent\textbf{Concern of this Thesis.} Since the treatment in
\cite{Dai91} is more general than what we have stated in Theorem
\ref{DaiIntro}, Dai's adiabatic limit formula continues to hold
for the odd signature operator twisted by a flat connection $A$
over $M$. We will see that there are natural analogs $\widehat
\eta_A$, $D_B\otimes \nabla^{\sH_{A,v}}$ and $\gs_A$ of the
quantities appearing in Theorem \ref{DaiIntro}. As the Rho
invariant $\rho_A(M)$ is independent of the metric, it is
immediate, that with respect to every adiabatic metric,
\[
\rho_A(M) = \lim_{\eps\to 0}\eta(B_{A,\eps}^{\ev}) - k \cdot
\lim_{\eps\to 0}\eta(B_{\eps}^{\ev}).
\]
Then, Theorem \ref{DaiIntro} yields
\begin{theorem}\label{RhoGenIntro}
Let $A$ be a flat $\U(k)$-connection over $M$. Then with respect
to every submersion metric
\[
\begin{split}
\rho_A(M) = 2^{[\frac{b+1}2]}\int_B \widehat L(TB,&\nabla^B)\wedge
\big(\widehat \eta_A -k\cdot \widehat \eta\big)\\ & + \lfrac 12
\eta\big(D_B\otimes \nabla^{\sH_{A,v}}\big) - \lfrac k2
\eta\big(D_B\otimes \nabla^{\sH_v}\big) + \gs_A - k\cdot \gs.
\end{split}
\]
\end{theorem}
Now, the main matter of this thesis is to investigate how this
rather straightforward consequence of Theorem \ref{DaiIntro} can
be used for explicit computations of Rho invariants. Due to the
technical nature of local families index theory, our first
concern is to assemble the building blocks we need, and specialize
many known results to the case of the odd signature operator. The
motivation here is certainly not to exhibit new results, but to
present the theory in such a way that it becomes accessible for a
treatment of Rho invariants of fiber bundles. Therefore, our
perspective will usually be a geometric one rather than discussing
analytical difficulties. For a discussion of these aspects, we
usually refer to the wide variety of literature. Nevertheless, we
include proofs of some folklore results, for instance a fibered
version of the Hodge decomposition theorem, and a result about
how to achieve that the mean curvature of a fiber bundle vanishes.

Apart from the theoretical discussion, our true focus is on
explicit examples. We will examine two important classes of
fibered 3-manifolds in detail.\\

\noindent\textbf{Circle Bundles over Surfaces.} The simplest
class of fiber bundles for which a discussion of Rho invariants
is meaningful, is given by principal $S^1$-bundles over Riemann
surfaces. Nevertheless, this family of manifolds already deserves
some attention as it is a model for two important classes of
manifolds, namely 3-dimensional Seifert fibrations and higher
dimensional principal bundles.

In this spirit, Nicolaescu \cite{Nic98, Nic99} has analyzed the
Seiberg-Witten equations of Seifert manifolds, and parts of our
discussion are influenced by his work. Given a closed, oriented
surface $\gS$, and an oriented principal bundle
$S^1\hookrightarrow M\to \gS$, we will see that we can represent
every flat $\U(1)$-connection $A$ over $M$ by pulling back a line
bundle of degree $k$ over $\gS$. In terms of this data, the Rho
invariant associated to $A$ is given as follows, see Theorem
\ref{RhoCircBund}.

\begin{theorem}\label{CircIntro}
Assume that $l\neq 0$, and that $q_0\in [0,1)$ is such that
$k/l\equiv q_0\mod \Z$. Then
\[
\rho_A(M) = 2l (q_0^2-q_0) +\sgn(l).
\]
If $l=0$, so that the fiber bundle is trivial, all Rho invariants
vanish.
\end{theorem}

We shall include two proofs of this result. The first one in
Chapter \ref{ChapBasic} uses only basic considerations about the
geometry of fiber bundles, and the second one in Chapter
\ref{ChapAbst} shows how Theorem \ref{RhoGenIntro} can be used
for this class of examples.\\

\noindent\textbf{3-dimensional Mapping Tori.} The second family
of manifolds we will consider are fiber bundles
$\gS\hookrightarrow M \to S^1$, where $\gS$ is again a closed,
oriented surface. A manifold of this type is determined by an
element $f$ of the mapping class group of $\gS$. Due to the rich
algebraic structure encoded in the latter, we have not attempted
to treat the class of 3-dimensional mapping tori in full
generality.

What we shall do instead, is to assume first that the monodromy
$f$ of the mapping torus $M$ is of finite order. Under this
assumption, the formula of Theorem \ref{RhoGenIntro} for the Rho
invariant of a flat $\U(k)$-connection $A$ over $M$ reduces to
\[
\rho_A(M) = \lfrac 12 \eta\big(D_B\otimes \nabla^{\sH_{A,v}}\big)
- \lfrac k2 \eta\big(D_B\otimes \nabla^{\sH_v}\big).
\]
In Theorem \ref{RhoFiniteMapTor} we shall give an expression of
the right hand side in terms of Hodge-de-Rham cohomology of
$\gS$, thus obtaining a cohomological formula for the Rho
invariant. Since the precise statement would need a longer
preamble, we refer Chapter \ref{3dimMapTor} for details.

After this we will consider $\U(1)$-Rho invariants of a mapping
torus $T^2_M$, where the fiber is the 2-dimensional torus $T^2$,
and $M\in \SL_2(\Z)$---the mapping class group of $T^2$. One
naturally has to distinguish between three cases, depending of
whether $M$ is elliptic, parabolic or hyperbolic. The first two
cases are rather special, and we shall not discuss them here. The
explicit formul{\ae} for the corresponding Rho invariants can be
found in Theorem \ref{RhoFiniteTorusBundle} and Theorem
\ref{RhoParab}, respectively.

The case of a hyperbolic monodromy matrix requires more
background material. Here we will use ideas of Atiyah
\cite{Ati87}, who gives a far-reaching treatment of the untwisted
Eta invariant for mapping tori with fiber $T^2$. In particular,
he uses the relation to Hirzebruch's signature defect to show
that for a hyperbolic element $M=\left(\begin{smallmatrix} a &b\\
c&d \end{smallmatrix}\right)\in \SL_2(\Z)$, the Eta invariant of
the odd signature operator with respect to a natural metric on
$T^2_M$, is given by
\begin{equation}\label{AtiIntro}
\eta(B^{\ev}) = \frac {a+d}{3c} - 4\sgn(c) s(a,c) -
\sgn\big(c(a+d)\big).
\end{equation}
Here, $s(a,c)$ is the \emph{Dedekind sum}
\begin{equation}\label{DedSumIntro}
s(a,c) = \sum_{k=1}^{|c|-1} P_1\big(\lfrac
{ak}{c}\big)P_1\big(\lfrac k{c}\big),
\end{equation}
where for $x\in \R$,
\[
P_1(x)= \begin{cases}\hphantom{B_n\big(x} 0, &\text{if $x\in
\Z$},\\ x- [x]-\lfrac 12, &\text{if $x\notin \Z$.}
\end{cases}
\]
Atiyah also relates the Eta invariant to the classical
\emph{Dedekind Eta function}, which is defined for a point $\gs$
in the upper half plane as
\[
\boldsymbol{\eta}(\gs):= q_\gs^{\frac 1{24}}\prod_{n=1}^\infty
(1-q_\gs^n),\quad  q_\gs:= e^{2\pi i \gs}.
\]
One can define a natural logarithm of $\boldsymbol{\eta}$, and the
transformation property of $\log \boldsymbol{\eta}$ under the
action of elements of $\SL_2(\Z)$ has a long history, going back
to Dedekind \cite{Ded}. For an element
$M=\left(\begin{smallmatrix} a &b\\ c&d
\end{smallmatrix}\right)\in \SL_2(\Z)$ with $c\neq 0$, this
transformation formula states that
\begin{equation}\label{EtaTransIntro}
\log \boldsymbol{\eta}(M\gs) - \log \boldsymbol{\eta}(\gs) =
\frac 12 \log\Big(\frac{c\gs +d}{\sgn(c)i}\Big) + \pi
i\Big(\frac{a+d}{12c} - \sgn(c)s(a,c)\Big).
\end{equation}
Here, the logarithm on the right hand side is the standard branch
on $\C\setminus\R^-$, and $s(a,c)$ is the Dedekind sum
\eqref{DedSumIntro}. Atiyah's explanation of the relation between
\eqref{EtaTransIntro} and the formula \eqref{AtiIntro} makes
essential use of the idea of taking the adiabatic limit.

Motivated by this, we will study in detail the expression
$\int_{S^1} \widehat \eta_A$ appearing in Theorem
\ref{RhoGenIntro} for the case of a flat $\U(1)$-connection over a
hyperbolic mapping torus over $S^1$. Using ideas related to
Kronecker's second limit formula, we will cast it into a form,
where the logarithm of a generalized Dedekind Eta function
naturally appears. We will then employ a transformation formula
due to Dieter \cite{Die59}, to obtain the analog of
\eqref{EtaTransIntro} for the logarithm of this generalized
Dedekind Eta function. From this we shall deduce the main result
of Chapter \ref{3dimMapTor}, see Theorem \ref{RhoHyp}.

\begin{theorem}\label{RhoHypIntro}
Let $M=\left(\begin{smallmatrix} a& b\\ c&
d\end{smallmatrix}\right)\in \SL_2(\Z)$ be hyperbolic, and let
$(\nu_1,\nu_2)\in \R^2\setminus \Z^2$ with $\nu_1\in [0,1)$
satisfy
\[
\begin{pmatrix} m_1 \\ m_2 \end{pmatrix} = (\Id -M^t)
\begin{pmatrix}
\nu_1\\ \nu_2 \end{pmatrix}\in \Z^2.
\]
This defines a flat connection $A$ over $T^2_M$, and
\[
\begin{split}
\rho_A(T^2_M) & = \lfrac {2(a+d)- 4}c ( \nu_1^2 -\nu_1) -
4\sum_{k=1}^{|c|-r} P_1\big(\lfrac{dk}{c}\big)  +
\sgn\big(c(a+d)\big) -
\sgn(c)\gd(\nu_1)\big(1-\gd(\lfrac{m_1}c)\big)
\\  & \qquad\quad  -2 P_1\big(\lfrac
{dm_1}{c}\big) -2\gd(\nu_1)\Big( P_1\big(\lfrac{m_1}{c}\big) -
P_1\big(\lfrac{dm_1}{c}\big)\Big),
\end{split}
\]
where $r\in \{0,\ldots |c|-1\}$ is such that $m_1\equiv r\mod c$,
and $\gd$ is the characteristic function of $\R\setminus \Z$.
\end{theorem}

Although this formula might appear to be somewhat involved, it is
satisfactory in two ways. First of all, the involved terms are
easy to compute for explicit choices of
$M=\left(\begin{smallmatrix} a& b\\ c& d\end{smallmatrix}\right)$
and $(\nu_1,\nu_2)$. Secondly, we shall see that it contains
previous computations of Chern-Simons invariants by Freed and Vafa
\cite{FreVaf} as a special case. In this respect, the author hopes
that a possible generalization to $\SU(2)$-connections will
reprove results of Kirk and Klassen \cite{KK90} and might shed a
new light on Jeffrey's conjecture \cite{Jef92} concerning the
spectral flow associated to twisted odd signature operators on a
mapping torus of the form considered here.\\

\noindent\textbf{Outline of this Thesis.} We end the introduction
with a very brief outline of this thesis. We will keep this
rather short since the beginning of each chapter contains a more
detailed outline of its contents.

\begin{itemize}
\item Chapter \ref{Signature} is a survey of results from index
theory that we will need. In particular, we shall introduce the
signature of a manifold, discuss its relation to index theory,
and recall the Atiyah-Singer Index Theorem in its cohomological
version for geometric Dirac operators. Then we introduce the Eta
and Rho invariant, and discuss how they appear in the index
theorem for manifolds with boundary. We also place some emphasis
on variation formul{\ae} and sketch how they are related to local
index theory.
\item In Chapter \ref{ChapBasic} we start with the discussion of fiber
bundles. We will introduce the geometric setup, paying close
attention to the structure of the odd signature operator. Then,
we shall encounter the basic idea of adiabatic limits and use
this to give an elementary proof of Theorem \ref{CircIntro} above.
\item Chapter \ref{ChapAbst} contains the main theoretical part of this
thesis. Here, the main objective is to introduce all quantities
appearing in Theorem \ref{DaiIntro}. After discussing the bundle
of vertical cohomology groups in some detail, we will give a
heuristic derivation of Theorem \ref{DaiIntro}. For this, we also
include a short survey of local families index theory. All this
discussion will lead to Theorem \ref{RhoGenIntro}, which we will
then use to reprove Theorem \ref{CircIntro} in a more abstract
way.
\item The content of Chapter \ref{3dimMapTor} is the discussion
of 3-dimensional mapping tori along the lines we have already
outlined above.
\item For the reader's convenience, and to keep our discussion more
self-contained, we have also included a couple of appendices,
which contain material that we freely use, but that would lead to
far afield if discussed in the main body of this thesis.
\begin{itemize}
\item Appendix \ref{CharClass} contains a discussion of Chern-Weil
theory in a way which is particularly well-suited for the
applications we need. Moreover, we have included some aspects
concerning Chern-Simons invariants, as they appear throughout our
discussion.
\item Since the Rho invariants we are interested in depend only on
the gauge equivalence class of the involved flat connection, we
include some remarks concerning the moduli space of flat
connections in Appendix \ref{ModApp}. We start giving  some
details on the relation to representations of the fundamental
group. Moreover, the moduli space of flat connections over a
mapping torus is discussed, since we need this in Chapter
\ref{3dimMapTor}. We end this appendix with a brief survey of the
moduli space of holomorphic line bundles over a Riemann surface,
which is an ingredient for discussing flat $\U(1)$-connections
over principal $S^1$-bundles over surfaces.
\item Appendix \ref{Comp} contains some computations. On the one
hand, we need explicit values of basic Eta and Zeta functions to
which we reduce most computations in the main body of the thesis.
On the other hand, we shall discuss the Dedekind sum in
\eqref{DedSumIntro} and its generalization, establishing a
relation among them which we need to prove Theorem
\ref{RhoHypIntro}.
\item Finally, Appendix \ref{AppEtaVar} includes some more
analytical details concerning the heat operator. Specifically, we
will give some remarks concerning families of heat operators, and
derive the variation formula for the Eta invariant.
\end{itemize}
\end{itemize}

\clearpage

\chapter*{Acknowledgements}

\pagestyle{empty}

\addcontentsline{toc}{chapter}{Acknowledgements}

Foremost, I want to express my deep gratitude to Prof. Dr.
Matthias Lesch, under whose supervision I chose this topic and
worked on this thesis. Without his interest and insight into
various branches of mathematics, I would not have learnt to
appreciate the beauty of index theory. I am particularly thankful
for his confidence in me with which he helped to overcome any dry
spell. I am also deeply indebted to Prof. Dr. Paul Kirk for
hosting an inspiring one-year stay at the Indiana University in
Bloomington. Many ideas in this thesis are influenced by numerous
discussions in which he explained to me the elegant way he thinks
about mathematics. Furthermore, I would like to thank Prof. Dr.
Matthias Kreck, Prof. Dr. Jens Frehse and Prof. Dr. Albrecht
Klemm for joining my thesis committee. Special thanks are also due
to Christian Frey and Benjamin Himpel, who were my office mates
in K\"{o}ln and Bonn. Both were always available for giving their
mathematical advice as well as encouraging me with their
friendship.

Moreover, I would like to thank the SFB/TR12, the Department of
Mathematics at Indiana University, and the Bonn International
Graduate School in Mathematics, for providing funding during the
period when this thesis was written. Especially the financial
support for attending various conferences and workshops is
gratefully acknowledged.

I shall not end without expressing my deep gratitude to my family
and friends for their encouragement and love on which I could
always rely. I cannot be thankful enough for all the support my
parents, Lothar and Monika Bohn, offered me to pursue my studies.

\cleardoublepage

\pagestyle{myheadings}

\renewcommand{\chaptermark}[1]%
   {\markboth{\textsc{\thechapter.\ #1}}{}}
\renewcommand{\sectionmark}[1]%
   {\markright{\textsc{\thesection.\ #1}}}

\renewcommand{\thetheorem}{\arabic{chapter}.\arabic{section}.\arabic{theorem}}
\renewcommand{\theequation}{\arabic{chapter}.\arabic{equation}}

\cleardoublepage \pagenumbering{arabic} \setcounter{page}{1}

\chapter{The Signature Operator and the Rho Invariant}\label{Signature}

In this chapter we will survey some results from index theory,
which we need in our discussion later on. The objective is not to
give a systematic treatment but merely to introduce the objects we
are interested in and to fix notation. For this reason we shall
include proofs only if they enrich the discussion and do not lead
too far afield.

We start with a brief discussion of the signature of a closed
manifold, placing emphasis on the generalized version where the
intersection form is associated to a flat unitary vector bundle.
Introducing the signature operator relates the signature to the
index of an elliptic operator, and this leads us to a discussion
of the main facts concerning the heat equation on closed
manifolds. Although our focus is on the signature operator and its
odd dimensional analog, we present the Atiyah-Singer Index Theorem
in its version for geometric Dirac operators. This is because many
ideas in the later chapters are influenced by local index theory
which is more transparent when formulated in terms of Clifford
modules and Dirac operators.

Then we will introduce the object which constitutes the main topic
of this thesis---namely the Eta invariant of an elliptic operator
on a closed manifold. Variation formul{\ae} for Eta invariants
will play a prominent role in the discussion in the next chapters.
Therefore, we discuss this topic in some detail. Most notably, the
behaviour under the variation of a flat twisting connection leads
naturally to the first appearance of a Rho invariant as the
difference of certain Eta invariants.

After this we describe---briefly leaving the realm of closed
manifolds---how the Eta invariant arises as a correction term in
the index theorem for a manifold with boundary. From then on the
focus will be on the case that the elliptic operator in question
is the odd signature operator. From the signature theorem for
manifolds with boundary we derive some general and well-known
properties of the Eta invariant. In particular, the relation of
Rho invariants to Chern-Simons invariants will be exhibited.

We close the general discussion of this chapter with a short
outline of how local index theory methods for odd dimensional
manifolds can be used to obtain important properties of Rho
invariants without referring to the Atiyah-Patodi-Singer Index
Theorem. Our interest in this is not only of a purely
academic nature, as similar ideas are the ones underlying local
families index theory, which we will need in the context of fiber
bundles in Chapter \ref{ChapAbst}.

\clearpage

\section{The Signature of a Manifold}

\subsection{Intersection Forms and Local Coefficients}

Let $M$ be a closed, oriented and connected manifold of dimension
$m$. On the real cohomology groups there exists the
\emph{intersection pairing}
\[
H^p(M,\R)\times H^{m-p}(M,\R)\to \R,\quad (a,b)\mapsto
\Scalar{a\cup b}{[M]},
\]
where $\scalar{.}{.}$ is the Kronecker pairing, $\cup$ is the cup
product, and $[M]$ denotes the fundamental class of $M$ determined
by the orientation. Expressing $H^p(M,\R)$ in terms of de Rham
cohomology groups, the intersection pairing is induced by
\[
\gO^p(M)\times \gO^{m-p}(M)\mapsto \R,\quad (\ga,\gb)\mapsto
\int_M\ga\wedge\gb.
\]
As a consequence of the Poincar\'e duality theorem, the above
pairing is non-degenerate. In particular, if $\dim M=m$ is even,
there is a non-degenerate bilinear form
\begin{equation}\label{IntersectionForm}
Q:H^{m/2}(M,\R)\times H^{m/2}(M,\R)\to \R,
\end{equation}
which is called the \emph{intersection form} of $M$. If $(m\equiv
0 \mod 4)$ the intersection form is symmetric, and if $(m\equiv
2\mod 4)$ it is skew. Recall that the signature of a symmetric
form is the number of positive minus the number of negative
eigenvalues.

\begin{convention}\label{SignConv}
We also use the convention that the signature of a skew form is
the number of positive imaginary eigenvalues minus the number of
negative imaginary ones.
\end{convention}

\begin{dfn}
Let $M$ be a closed, oriented and connected manifold of even
dimension $m$. Then the \emph{signature} of $M$ is defined as
\[
\Sign(M):=\Sign(Q).
\]
\end{dfn}

\begin{remark*}
In topology, the intersection form is usually considered as a
form over $\Z$. If one uses cohomology groups with integer
coefficients, one has to divide out the torsion subgroup of
$H^{m/2}(M,\Z)$ to get a non-degenerate form. Moreover, note that
if we want to work with complex coefficients, we have to extend
$Q$ anti-linearly in, say, the first variable to get a (skew)
Hermitian form. Then the signature is also well-defined and agrees
with the signature of the underlying real form. Note, however,
that if a skew form comes from a form over $\R$ it has zero
signature since its eigenvalues come in conjugate pairs.
\end{remark*}

\noindent\textbf{Cohomology with Local Coefficients.} We will
also need a twisted version of the intersection form. For this we
briefly recall the construction of cohomology groups with local
coefficients. We refer to \cite[Ch. 5]{DavKir} for more details
and proofs.

Let $M$ be a connected manifold, not necessarily closed, and let
$\widetilde M$ be the universal cover of $M$. Let $\pi=\pi_1(M)$
be the fundamental group of $M$, and let $\C[\pi]$ denote the
group algebra of $\pi$. The fundamental group acts from the right
on $\widetilde M$, so that the cellular chain groups
$C_p(\widetilde M)$ are $\C[\pi]$ right modules in a natural way.
This is because a cell decomposition of $M$ induces via lifting
of cells a cell decomposition of $\widetilde M$ which is
compatible with the action of $\pi$ on $\widetilde M$.

Now let $\ga:\pi \to \U(k)$ be a unitary representation, and let
\[
C^p(M,E_\ga):=\Hom_{\C[\pi]}\big(C_p(\widetilde M), \C^k).
\]
Here, the action of $\C[\pi]$ on $\C^k$ is by matrix
multiplication $x\mapsto \ga(x)^{-1}$. The differential on
cochains on $\widetilde M$ turns $C^\bullet(M,E_\ga)$ into a
complex. As for untwisted cellular cohomology, the homology of
this complex does not depend on the particular cell decomposition
and the chosen lifts.

\begin{dfn}
Let $M$ be a compact, connected manifold, and
$\ga:\pi_1(M)\to\U(k)$ a representation. Then the homology of
$C^\bullet(M,E_\ga)$ is denoted by $H^\bullet(M,E_\ga)$ and is
called the cohomology of $M$ with \emph{local coefficients} given
by $\ga$.
\end{dfn}

If we now assume that $M$ is also closed and oriented, then there
exists a non-degenerate pairing
\begin{equation}\label{TwistIntPair}
H^p(M,E_\ga)\times H^{m-p}(M,E_\ga) \to \C
\end{equation}
induced by the cup product on the cohomology of $\widetilde M$ and
the scalar product on $\C^k$. If $M$ is of even dimension $m$ this
yields a bilinear form
\[
Q_\ga: H^{m/2}(M,E_\ga)\times H^{m/2}(M,E_\ga) \to \C,
\]
which we call the \emph{twisted intersection form.}

\begin{dfn}
Let $M$ be a closed, oriented and connected manifold of even
dimension $m$, and let $\ga:\pi_1(M)\to\U(k)$ be a unitary
representation of the fundamental group of $M$. Then the
\emph{twisted signature} of $M$ is defined as
\[
\Sign_\ga(M):=\Sign(Q_\ga),
\]
where we use again the convention that the signature of a skew
form is the number of positive imaginary eigenvalues minus the
number of negative imaginary ones.
\end{dfn}

\begin{remark*}
As we will see soon, the twisted signatures we have just defined
give no new topological information for a closed manifold $M$.
However, their version for manifolds with boundary are a
non-trivial generalization of topological importance.
Nevertheless, we have included the definition here to keep the
discussion parallel.
\end{remark*}

\noindent\textbf{Local Coefficients and Flat Bundles.} Twisted
cohomology groups can also be defined in terms of differential
forms. Let $E\to M$ be a Hermitian vector bundle of rank $k$,
endowed with a unitary connection $A$. This gives rise to a
twisted version of the exterior differential
\[
d_A:\gO^p(M,E)\to \gO^{p+1}(M,E), \quad d_A(\go\otimes e)=d\go
\otimes e + \go\wedge Ae.
\]
The square of $d_A$ is given by exterior multiplication with the
curvature $F_A\in\gO^2\big(M,\End(E)\big)$. Therefore, we get a
complex $\big(\gO^\bullet(M,E), d_A\big)$ precisely if $A$ is
flat.

\begin{dfn}
Let $E\to M$ be a Hermitian vector bundle of rank $k$, endowed
with a unitary flat connection $A$. Then we denote the homology of
$\big(\gO^\bullet(M,E), d_A\big)$ by $H^\bullet(M,E_A)$ and call
it the cohomology of $M$ \emph{with values in the flat bundle
$(E,A)$}.
\end{dfn}

We briefly sketch the relation between cohomology with local
coefficients and cohomology with values in a flat bundle. As a
general reference, we refer to \cite[Sec. 5.5]{Ram}. In addition,
we have included a detailed discussion concerning the equivalence
of flat connections and representations of the fundamental group
in Appendix \ref{FlatConn}. The language there is in terms of
principal bundles, but the translation to Hermitian vector
bundles is done without effort.

Let $M$ be a connected manifold, not necessarily closed, and let
$E$ be a flat unitary bundle over $M$ with connection $A$. As
explained in Appendix \ref{FlatConn} we can lift every closed
loop $c$ in $M$ horizontally to $E$ with respect to $A$. The
starting point and the end point of the lifted loop lie in the
same fiber of $E$ and since $A$ is unitary they differ by the
action of an element in $\U(k)$. This construction gives rise to
the holonomy representation of the based loop group of $M$, see
\eqref{HolRepDef} on p. \pageref{HolRepDef}. Since $A$ is flat,
the holonomy representation depends only on homotopy classes.
Thus, we obtain a representation
\[
\hol_A:\pi_1(M)\to \U(k),
\]
which is precisely the object we need to define cohomology with
local coefficients.

Conversely, let us start with a representation $\ga:\pi_1(M)\to
\U(k)$ of the fundamental group. Via the action of $\pi_1(M)$ as
the group of deck transformations we may interpret the universal
cover $\widetilde M$ as a $\pi_1(M)$-principal bundle over $M$.
The representation $\ga$ defines an associated vector bundle
\[
E_\ga= \widetilde M\times_\ga \C^k \to M.
\]
Since $\ga$ is unitary, one can define a natural Hermitian metric
on $E_\ga$. Moreover, the trivial connection on $\widetilde
M\times\C^k$ descends to a unitary, flat connection $A_\ga$ on
$E_\ga$.

When taking suitable equivalence classes, the above constructions
are inverses of each other, see Appendix \ref{FlatConn}. Then
there is the following twisted version of the de Rham Theorem,
see for example the discussion in \cite[p. 154]{Ram}.

\begin{prop}\label{TwistedDeRhamThm}
Let $M$ be a connected manifold, and let $E$ be a Hermitian vector
bundle with flat connection $A$. If $\ga:\pi_1(M)\to \U(k)$ is the
holonomy representation of $A$, then there is a natural
isomorphism
\[
H^\bullet(M,E_A)\xrightarrow{\cong} H^\bullet(M,E_\ga).
\]
Moreover, if $M$ is closed, the twisted intersection pairing of
\eqref{TwistIntPair} corresponds under this isomorphism to the
bilinear form on $H^\bullet(M,E_A)$ induced by
\[
\gO^p(M,E)\times \gO^{m-p}(M,E)\mapsto \C,\quad (\go,\eta)\mapsto
\int_M \langle \go\wedge\eta \rangle,
\]
The notation $\langle\go\wedge\eta\rangle$ is shorthand for
taking the exterior product in the differential form part and
pairing with the Hermitian metric in the bundle part.
\end{prop}
Having the above canonical isomorphism in mind, we will henceforth
not carefully distinguish between $H^\bullet(M,E_A)$ and
$H^\bullet(M,E_\ga)$. Similarly, when concerned with the
intersection form, we will also write $Q_A$ and $\Sign_A(M)$ if
the focus is on a flat connection rather than a representation of
the fundamental group.

\subsection{Twisted Signature Operators}

Historically, one of the starting points of index theory is the
observation that the signature can be described as the index of
an elliptic operator. Before we can define the signature
operator, we need to fix some notation and conventions. Since
they are of a purely linear algebraic nature, we formulate them
for an oriented vector space $V$ of dimension $m$ over $\R$ which
plays the role of the cotangent space $T^*_xM$.\\

\noindent\textbf{Algebraic Preliminaries.} Let $V$ be an
Euclidean vector space with scalar product $g$, and let
$\gL^\bullet V_\C$ denote the complexified exterior algebra. We
endow it with the Hermitian metric $g_h$ given by extending $g$
antilinearly in the first variable. Then
\begin{equation*}
g_h(\ga,\gb)\vol(g) = \overline \ga\wedge *\gb,\quad \ga,\gb\in
\gL^\bullet V_\C,
\end{equation*}
where $\vol(g)$ is the volume element given by $g$ and the
orientation of $V$ and $*$ is the complex linear extension of the
Hodge $*$ operator. $V$ acts on $\gL^\bullet V$ via exterior
multiplication
\begin{equation*}\label{ExtDef}
\emu(v)\ga = v\wedge\ga,\quad \ga\in \gL^\bullet V.
\end{equation*}
Using the metric $g$, one defines an interior multiplication by
requiring that
\begin{equation*}\label{IntDef}
\imu(v)w = g(v,w)\quad\text{and}\quad \imu(v)(\ga\wedge\gb) =
\imu(v)(\ga)\wedge\gb + (-1)^{|\ga|}\ga\wedge \imu(v)\gb,
\end{equation*}
for every $v,w\in V$ and $\ga,\gb\in \gL^\bullet V$ with $\ga$
homogeneous of degree $|\ga|$. For $v\in V$ we extend $\imu(v)$
and $\emu(v)$ complex linearly to $\gL^\bullet V_\C$, and define
the \emph{Clifford multiplication}
\begin{equation}\label{CliffDef}
c: V\to \End_\C\big(\gL^\bullet V_\C\big),\quad
c(v):=\emu(v)-\imu(v).
\end{equation}
Then, for all $v\in V$,
\[
c(v)^2 = -g(v,v),\quad c(v)^* = -c(v).
\]
This means that $c$ extends to a complex representation of the
Clifford algebra\footnote{Recall (e.g. from \cite[Sec. 1.3]{BGV})
that the \emph{Clifford algebra} $\cl(V,g)$ of an $\R$ vector
space $V$ with a metric $g$ is the $\R$ algebra generated by $V$
and the relations \[ c(v)c(w)+c(w)c(v) = -2g(v,w),\quad v,w\in
V.\]}
\[
c:\cl(V,g) \to \End_\C\big(\gL^\bullet V_\C\big).
\]
Equivalently, this is a representation of the complexified
Clifford algebra $\cl_\C(V)$. We define the \emph{symbol map}
\begin{equation}\label{SymbolMapDef}
\boldsymbol{\gs}: \cl_\C(V) \to \gL^\bullet V_\C ,\quad a\mapsto
c(a)1.
\end{equation}
The symbol map $\boldsymbol{\gs}$ is an isomorphism of vector
spaces, but certainly not of algebras. The inverse
$\boldsymbol{\gs}^{-1}$ is called the \emph{quantization map}.
Using this we define the \emph{chirality operator}
\begin{equation}\label{TauDef}
\tau : = i^{[\frac{m+1}{2}]}\;c\circ
 \boldsymbol{\gs}^{-1}\big(\vol(g)\big) : \gL^\bullet
V_\C \to \gL^\bullet V_\C.
\end{equation}
Here, $[\frac{m+1}{2}]$ denotes the integral part of
$\frac{m+1}{2}$. The following result is straightforward, see
\cite[Prop. 3.58]{BGV}.

\begin{lemma}\label{TauProp}
The chirality operator $\tau$ satisfies
\[
\tau^2 = \Id\quad\text{and}\quad \tau^* = \tau = \tau^{-1}.
\]
On $\gL^p V_\C$ it is explicitly given by
\begin{equation}\label{TauExpl}
\tau = (-1)^{\frac{p(p-1)}2 +mp}\,i^k *_p,
\end{equation}
where $k:=[\frac{m+1}{2}]$ and $*_p$ is the complex linear Hodge
$*$ operator on $\gL^p V_\C$. Moreover,
\[
\tau\circ c(v) = (-1)^{m+1} c(v)\circ \tau\quad\text{and}\quad
\tau\circ \imu(v) = (-1)^m \emu(v)\circ \tau.
\]
\end{lemma}

\begin{convention}\label{ComplexConv}
From now on we will always drop the subscripts $\C$. Thus, for a
real vector space $V$, we will use $\gL^\bullet V$ to denote the
complexified exterior algebra, and $\cl(V)$ will denote the
complexified Clifford algebra.
\end{convention}

\noindent\textbf{The Signature an the Index.} Now let $M$ be an
oriented manifold of dimension $m$, endowed with a Riemannian
metric $g$. We fix a Hermitian vector bundle $E\to M$ of rank
$k$, endowed with a unitary connection $A$. Let
$\gO^\bullet(M,E)$ denote differential forms with values in $E$.
The Riemannian metric $g$ and the bundle metric on $E$ define an
$L^2$ scalar product on $\gO^\bullet(M,E)$. With respect to this,
the formal adjoint of the twisted exterior differential
\[
d_A: \gO^\bullet(M,E)\to \gO^{\bullet+1}(M,E)
\]
is given in terms of the chirality operator $\tau=\tau_M$ from
\eqref{TauDef} by
\begin{equation}\label{d^tDef}
d_A^t=(-1)^{m+1}\tau\circ d_A\circ \tau,
\end{equation}
see \cite[Prop. 3.58]{BGV}. Here, $\tau_M$ acts only on the
differential form part. Note that we are using that $A$ is
unitary. Now \eqref{d^tDef} implies that the twisted de Rham
operator $d_A+d^t_A$ satisfies
\begin{equation}\label{deRhamTauComm}
\tau(d_A+d^t_A) =(-1)^{m+1}(d_A+d^t_A)\tau.
\end{equation}
Let us assume from now on that $m$ is even. Since $\tau$ is an
involution, we may decompose $\gO^\bullet(M,E)$ into the $\pm 1$
eigenspaces of $\tau$,
\[
\gO^\bullet(M,E)=\gO^+(M,E)\oplus \gO^-(M,E).
\]
It follows from \eqref{deRhamTauComm} that we can decompose
\[
d_A+d_A^t=\begin{pmatrix} 0 & D_A^-\\ D_A^+
&0\end{pmatrix},\quad\text{where }
D_A^+:\gO^+(M,E)\to\gO^-(M,E),\quad D_A^-= (D_A^+)^t.
\]

\begin{dfn}\label{SignOpDef}
Let $M$ be an even dimensional, oriented Riemannian manifold, and
let $E\to M$ be a Hermitian vector bundle with a unitary
connection $A$. Then
\[
D^+_A:\gO^+(M,E)\to\gO^-(M,E)
\]
is called the \emph{twisted signature operator} of $M$ twisted by
$A$.
\end{dfn}

Since the signature operator is an elliptic operator of first
order and $M$ is a closed manifold, we have a well-defined index
problem. The name \emph{signature operator} is only justified if
$A$ is a flat connection, since then the index of $D_A^+$ is
indeed the twisted signature $\Sign_A(M)$.

\begin{prop}\label{Sign=Ind}
Let $M$ be a closed, oriented and connected Riemannian manifold
of even dimension $m$. Let $E$ be a Hermitian vector bundle,
endowed with a unitary flat connection $A$. Then
\[
\Sign_A(M) = \ind(D_A^+).
\]
\end{prop}

\begin{proof}
Although this result can be found in many textbooks, we include
its proof as related arguments will appear again. First of all,
since $(D^+_A)^t=D^-_A$,
\begin{equation}\label{Sign=Ind:1}
\ind(D^+_A)=\dim\big( \ker D^+_A\big) - \dim \big(\ker D^-_A\big)
\end{equation}
To identify $\Ker D^\pm_A$ in cohomological terms consider the
twisted Laplacian
\begin{equation}\label{TwistLaplace}
\gD_A=(d_A+d^t_A)^2: \gO^\bullet(M,E)\to \gO^\bullet(M,E).
\end{equation}
The twisted version of the Hodge isomorphism identifies
$H^p(M,E_A)$ with the space of harmonic forms
\[
\sH^p(M,E_A):=\big(\ker \gD_A\big)\cap \gO^p(M,E).
\]
It follows from \eqref{deRhamTauComm} that the chirality operator
$\tau$ commutes with $\gD_A$, hence it induces an involution on
$\sH^\bullet(M,E_A)$. On $\sH^{m/2}(M,E_A)$ the chirality
operator is $\tau=i^k*$, where $k:= m^2/4$. Using this the
intersection form can be expressed in terms of harmonic forms as
\[
Q_A\big(\ga,\gb\big)=\int_M\langle \ga\wedge\gb\rangle =
\LScalar{\ga}{i^k\tau\gb},\quad \ga,\gb\in \sH^{m/2}(M,E_A).
\]
Therefore,
\[
\Sign(Q_A)= \Sign(i^k\tau|_{\sH^{m/2}}),
\]
where we use the same convention as before regarding the signature
of a skew endomorphism. We deduce that
\[
\Sign(Q_A)= \dim \sH^{m/2}(M,E_A)^+ - \dim \sH^{m/2}(M,E_A)^-.
\]
If $p\neq {m/2}$, there are isomorphisms
\begin{equation*}
\gF^\pm: \sH^p\xrightarrow{\cong} \big(\sH^p\oplus
\sH^{m-p}\big)^\pm,\quad \gF^\pm(\ga):=\lfrac12(\ga\pm\tau \ga).
\end{equation*}
From this and from the fact that $\Ker(d_A+d^t_A)=\Ker\gD_A$ we
find
\[
\dim\big( \ker D^\pm_A\big) = \dim \big(\sH^{m/2}(M,E_A)^\pm\big)
+ \sum_{p<m/2} \dim\big(\sH^p(M,E_A)\big)
\]
Thus, all terms in \eqref{Sign=Ind:1} are cancelled except the
one in the middle degree so that
\[
\ind(D^+_A)= \dim \sH^{m/2}(M,E_A)^+ - \dim \sH^{m/2}(M,E_A)^- =
\Sign(Q_A)\qedhere
\]
\end{proof}

\begin{remark*}
We have already pointed out that for closed manifolds the twisted
signatures do not carry interesting topological information. If
the flat twisting bundle is trivial, this can be deduced from the
above result: Let $A$ be a flat connection on the trivial vector
bundle $E=M\times \C^k$, and let $D_A^+$ be the associated
signature operator. Furthermore, let $D_{\oplus k}^+$ denote the
signature operator associated to the trivial connection on $E$.
Clearly, $D_A^+ - D_{\oplus k}^+$ is an operator of order 0. On a
closed manifold adding a 0-order perturbation to an elliptic
operator of first order does not change the index. This is because
it is a compact perturbation of a Fredholm operator in the
appropriate Hilbert space setting. Therefore,
\[
\Sign_A(M)=\ind(D_A^+)=\ind(D_{\oplus k}^+) = k\cdot\ind(D^+) =
k\cdot \Sign(M).
\]
This means that the only new information encoded in the twisted
signature is the rank of the twisting bundle. We will see in
Corollary \ref{SignMultClosed} below that this is also true for
flat twisting bundles which are topologically non-trivial.
\end{remark*}

\section{Dirac Operators and the Atiyah-Singer Index
Theorem}\label{AtiyahSinger}

The famous index theorem equates the index of an elliptic operator
over a closed manifold $M$ with the integral over certain
characteristic classes over $M$. In this section we briefly recall
the definitions occurring in the index theorem for Dirac type
operators.

\subsection{The Index and the Heat Equation}

We first recall some facts about the spectral theory of formally
self-adjoint elliptic operators on closed manifolds, see e.g.
\cite[Sec.'s. 1.3 \& 1.6]{Gil}.

\begin{dfn}\label{SelfAdjEllDef}
Let $M$ be a Riemannian manifold, and let $E\to M$ be a Hermitian
vector bundle. We denote by
\begin{equation*}
\sP_{s,e}^d=\sP_{s,e}^d(M,E)
\end{equation*}
the space of  formally self-adjoint elliptic differential
operators of order $d$.
\end{dfn}


\begin{theorem}\label{EllOpSpec}
Let $D\in \sP_{s,e}^d(M,E)$, and assume that $M$ closed.
\begin{enumerate}
\item The operator $D$ extends to an unbounded self-adjoint
operator in $L^2(M,E)$ with domain the Sobolev space
$L^2_d(M,E)$. It defines Fredholm operators
\[
D:L^2_{s+d}(M,E)\to L^2_s(M,E),\quad s\in \R,
\]
with Fredholm index independent of $s$.
\item There exists a constant $C$ such that for
all $\gf\in C^\infty(M,E)$
\begin{equation}\label{Garding}
\|\gf \|_{L_d^2} \le C\big(\|\gf\|_{L^2}+ \|D \gf\|_{L^2}\big).
\end{equation}
\item The spectrum $\spec(D)$ is a discrete subset of $\R$
consisting of eigenvalues with finite multiplicities. There is an
orthonormal basis of $L^2(M,E)$ consisting of smooth eigenvectors.
\item If we define
\[
N(\gl) := \#\bigsetdef{\gl_n\in\spec(D)}{|\gl_n|\le \gl},\quad
\gl\ge 0,
\]
then for some constant $C>0$
\begin{equation}\label{EigenvalueAsymp}
N(\gl) \sim C \gl^{m/d},\quad \text{as } \gl\to\infty.
\end{equation}
\item If $D$ has a positive definite leading symbol, then
the spectrum $\spec(D)$ is bounded from below.
\end{enumerate}
\end{theorem}

\noindent\textbf{The Heat Kernel.} We now specialize to a second
order operator, $H\in \sP_{s,e}^2(M,E)$, and assume that $H$ has
positive definite leading symbol. Clearly, the model we are
having in mind is a generalized Laplacian, i.e., an operator
$H\in \sP_{s,e}^2(M,E)$ such that its principal symbol $\gs(H)$
satisfies
\[
\gs(H)(x,\xi) = |\xi|_g^2\id_{E_x},\quad \xi\in T^*_xM.
\]
Recall that the heat equation with initial condition $\gf\in
L^2(M,E)$ in terms of $H$ is the partial differential equation
\begin{equation}\label{HeatEqn}
\big(\lfrac d{dt} + H\big) \gf(t) = 0,\quad t\ge0,\quad
\gf(0)=\gf.
\end{equation}
Formally, if $\{\gl_n\}_{n\ge -n_0}$ denotes the set of
eigenvalues of $H$ with eigenvectors $\gf_n$, the solution to
\eqref{HeatEqn} is
\[
\begin{split}
\gf(t)  = e^{-t H} \gf  &= \sum_{n\ge -n_0} e^{-t\gl_n}
\gf_n\LScalar{\gf_n}{\gf } \\
&= \int_M \sum_{n\ge -n_0} e^{-t\gl_n}
\gf_n(x)\scalar{\gf_n(y)}{\gf(y)}\vol_M(y).
\end{split}
\]
Thus, $e^{-t H}$ is an integral operator with kernel
\begin{equation}\label{HeatKernel}
k_t(x,y)=e^{-tH}(x,y)= \sum_{n\ge -n_0} e^{-t\gl_n}
\gf_n(x)\otimes \gf_n(y)^* \in C^\infty\big(M\times M,E\boxtimes
E^*\big).
\end{equation}
Here, for vector bundles $E\to M$ and $F\to N$, we employ the
standard notation
\begin{equation}\label{FiberProd}
E\boxtimes F:= \pi_M^*E\otimes \pi_N^*F\to M\times N,
\end{equation}
where $\pi_M$ and $\pi_N$ are the natural projections. The formal
expression \eqref{HeatKernel} can be made precise using the
following basic estimate.

\begin{lemma}\label{BasicKernelEst}
Let $0<\gl_{0}\le \gl_1\le \ldots$ denote the positive
eigenvalues of $H$, and let $\gl_{-1}\le\ldots\le\gl_{-n_0}$
denote the finite number of eigenvalues of $H$ which are less or
equal than 0. Let $\{\gf_n\}_{n\ge -n_0}$ be a basis of smooth
eigenvectors and for $N\in\N$ consider
\[
k_t^N(x,y):=\sum_{n=-n_0}^N e^{-t\gl_n}\gf_n(x)\otimes \gf_n^*(y)
\in C^\infty\big(M\times M,E\boxtimes E^*\big).
\]
Then for every $k\in \N$ and $t_0>0$ there exists a constant $C$
such that for every $N\in\N$ and $t\ge t_0$ we can estimate
\begin{equation*}
\Big\|k_t^N(x,y)- \sum_{n< 0}e^{-t\gl_n}\gf_n(x)\otimes
\gf_n^*(y) \Big\|_{C^k} \le C e^{-t\gl_{0}/2}.
\end{equation*}
\end{lemma}

\begin{proof}
Sobolev embedding \cite[Lem. 1.3.5]{Gil} and the elliptic
estimate \eqref{Garding} imply that for $l>k+  m/2$ there exist
constants $C_1, C_2$ such that for all $n\ge 0$
\[
\|\gf_n\|_{C^k}\le C_1 \|\gf_n\|_{L_l^2}\le C_2\big(
\|\gf_n\|_{L^2} + \|H^{l/2}\gf_n\|_{L^2}\big) = C_2\big(1+
\gl_n^{l/2}\big).
\]
Thus, for some other constants $C_1$ and $C_2$,
\[
\big\|e^{-t\gl_n}\gf_n(x)\otimes \gf_n^*(y)\big\|_{C^k} \le C_1
e^{-t\gl_n} \big(1+\gl_n^l \big) \le C_2 e^{-t\gl_n/2} (1+t^{-l}),
\]
where we have used that for $x=\gl_nt>0$
\[
x^l\le C_l\, e^{x/2}\quad\text{and thus,}\quad x^l e^{-x} \le
C_l\, e^{-x/2}.
\]
Now for each $n\ge 0$ we have $\gl_n\ge\gl_{0}$ and thus for $t\ge
t_0$
\[
e^{-t\gl_n/2} (1+t^{-l}) \le e^{-t_0\gl_n/2} (1+t_0^{-l})
e^{-(t-t_0)\gl_{0}/2}= C e^{-t_0\gl_n/2} e^{-t\gl_{0}/2},
\]
where the constant $C$ depends only on $t_0$ and $l$. Putting the
pieces together we find that there exists a constant $C$ such
that for every $N\ge 0$ and $t\ge t_0$
\[
\Big\|\sum_{n=0}^N e^{-t\gl_n}\gf_n(x)\otimes
\gf_n^*(y)\Big\|_{C^k} \le C e^{-t\gl_{0}/2} \sum_{n=0}^N
e^{-t_0\gl_n/2}
\]
Now the eigenvalue asymptotics \eqref{EigenvalueAsymp} shows that
$\sum_{n=0}^N e^{-t_0\gl_n/2}$ is absolutely convergent for
$N\to\infty$. Absorbing this to the constant, we get the desired
result.
\end{proof}

The above result shows that the kernels $k_t^N(x,y)$ converge for
$N\to\infty$ to the expression \eqref{HeatKernel}, uniformly with
respect to all $C^k$. Hence, we can define $e^{-tH}$ for $t>0$ by
\[
(e^{-tH}\gf)(x) = \int_M k_t(x,y)\gf(y)\vol_M(y),\quad \gf\in
L^2(M,E).
\]
In Appendix \ref{MoreHeat} we give an expression for $e^{-tH}$
using the spectral theorem. It is then easy to check that the
collection $e^{-tH}$ forms a strongly continuous semi-group in
each Sobolev space $L_s^2(M,E)$, i.e.,
\[
e^{-(s+t)H}= e^{-sH}e^{-tH},\quad\text{and}\quad \lim_{t\to
0}\|e^{-tH}\gf - \gf\|_{L_s^2}=0\quad\text{for each $\gf\in
L_s^2(M,E)$},
\]
see also Proposition \ref{HeatSolve}. Moreover, $e^{-tH}$ is
smooth in $t$ and does indeed solve the heat equation
\eqref{HeatEqn}. In addition, each $e^{-tH}$ is trace class, and
\[
\Tr e^{-tH} =\int_M\tr_E\big[k_t(x,x)\big]\vol_M(x),
\]
where $\tr_E:C^\infty\big(M,\End(E)\big)\to C^\infty(M)$ denotes
the fiberwise trace.

More generally, given an auxiliary differential operator
$D:C^\infty(M,E)\to C^\infty(M,E)$ of order $d\ge 0$, the uniform
bound of Lemma \ref{BasicKernelEst} ensures that we can apply $D$
under the integral to get
\[
De^{-tH}\gf(x) = D\Big(\int_M k_t(x,y)\gf(y)\vol_M(y)\Big) =
\int_M D_xk_t(x,y)\gf(y)\vol_M(y),
\]
where $D_x$ means applying $D$ with respect to the $x$ variable.
Thus, the operator $De^{-tH}$ has a smooth kernel
\[
D_x k_t(x,y)\in C^\infty\big(M\times M,E\boxtimes E^*\big),
\]
so that $De^{-tH}$ is trace class. The estimate in Lemma
\ref{BasicKernelEst} implies the following basic estimate on
$\Tr(De^{-tH})$.

\begin{prop}\label{BasicTraceProp}
For every $t_0>0$ there exists a constant $C$ such that for all
$t\ge t_0$
\begin{equation*}
\big|\Tr\big(De^{-tH}P_{(0,\infty)}\big)\big| \le C e^{-t\gl_0/2},
\end{equation*}
where $\gl_0$ is the smallest positive eigenvalue of $H$ and
$P_{(0,\infty)}$ is the spectral projection of $H$ associated to
the interval $(0,\infty)$.
\end{prop}

\begin{proof}
The kernel of $e^{-tH}P_{(0,\infty)}$ is given by
\[
\tilde k_t(x,y) := k_t(x,y)- \sum_{n< 0}e^{-t\gl_n}\gf_n(x)\otimes
\gf_n^*(y).
\]
Then we find that for $t\ge t_0$
\[
\begin{split}
\big|\Tr\big(De^{-tH}P_{(0,\infty)}\big)\big| &=
\Big|\int_M\tr_E\big[D_x\tilde k_t(x,y)\big]_{y=x}\vol_M(x)\Big|\\
&\le \rk(E)\vol(M)\, \big\|D_x\tilde k_t(x,y)\big\|_{C^0}\\ &\le
C_1\, \big\|\tilde k_t(x,y)\big\|_{C^d} \le C_2\, e^{-t\gl_0/2},
\end{split}
\]
where in the last line we have used that $D$ is a differential
operator of order $d$ and then Lemma \ref{BasicKernelEst}.
\end{proof}

\noindent\textbf{The McKean-Singer Formula.} We now turn our
attention to first order differential operators. To be able to
restrict to the formally self-adjoint case, we use the following
construction: Let $D^+:C^\infty(M,E^+)\to C^\infty(M,E^-)$ be an
elliptic differential operator of first order, acting between
Hermitian vector bundles $E^+$ and $E^-$. We define $E:=E^+\oplus
E^-$, and consider
\begin{equation}\label{GradedEllOp}
D:=\begin{pmatrix} 0 & D^-\\ D^+ &0\end{pmatrix}:C^\infty(M,E)\to
C^\infty(M,E) ,\quad \text{where}\quad D^-:= (D^+)^t.
\end{equation}
Certainly, an operator of this form is formally self-adjoint and
elliptic.

\begin{dfn}
Let $E\to M$ be a Hermitian vector bundle endowed with a splitting
$E=E^+\oplus E^-$.
\begin{enumerate}
\item Let $\gs:E\to E$ be the involution on $E$ given by
$\gs|_{E^\pm}=\pm\id$. Then $\gs$ is called the \emph{grading
operator} of $E$.
\item  An operator $D\in\sP_{s,e}^1(M,E)$ is called
\emph{$\Z_2$-graded} if
\[
\{D,\gs\} = D\gs +\gs D = 0.
\]
Note that---unless stated otherwise---we are using commu\-tators
and anti-commu\-tators in an \emph{ungraded} sense. Clearly, $D$
is $\Z_2$-graded if and only if it is of the form
\eqref{GradedEllOp}.
\item For $T\in C^\infty\big(M,\End(E)\big)$, we define
\emph{fiberwise supertrace} of $T$ as
\[
\str_E (T) := \tr_E(\gs T)\in C^\infty(M).
\]
\item If $D\in \sP_{e,s}^1(M,E)$ is $\Z_2$-graded, and $M$ is
closed, then the \emph{heat supertrace} associated to $D$ is
defined as
\[
\Str (e^{-tD^2}) := \Tr\big(\gs e^{-tD^2}\big).
\]
\end{enumerate}
\end{dfn}

Note that in (iv) the operator $H=D^2$ is positive and splits as
$H=H_+\oplus H_-$, where $H_\pm$ are formally self-adjoint,
positive operator as well. Thus, we are in the situation
considered before so that the respective heat traces exist, and
\[
\Str (e^{-tD^2})= \Tr(e^{-tH_+}) - \Tr(e^{-tH_-}) =
\Tr(e^{-t(D^+)^tD^+}) - \Tr(e^{-tD^+(D^+)^t}).
\]
Now, as an elliptic operator on a closed manifold $D^+$ has a
well-defined Fredholm index
\[
\ind(D^+) = \dim\big(\ker D^+\big) - \dim\big(\ker
(D^+)^t\big)=\Str(P_0),
\]
where $P_0:L^2(M,E)\to \ker(D)$ is the orthogonal projection on
the kernel of $D$. The famous McKean-Singer formula \cite{McKSin}
relates this index and the heat supertrace.

\begin{theorem}[McKean-Singer]\label{McKeanSinger}
Let $M$ be a closed manifold, $E\to M$ a Hermitian vector bundle,
and let $D\in \sP_{e,s}^1(M,E)$ be $\Z_2$-graded. Then for all
$t>0$
\[
\ind(D^+) = \Str (e^{-tD^2}).
\]
\end{theorem}
\begin{proof}
Since $D^2$ has no negative eigenvalues, we have
$P_0+P_{(0,\infty)}=\Id$. Then the estimate in Proposition
\ref{BasicTraceProp} implies that there exist constants $c$ and
$C$ such that for large $t$
\[
\big|\Str (e^{-tD^2})-\Str(P_0)\big| = \big| \Str
(e^{-tD^2}P_{(0,\infty)})\big| \le Ce^{-ct}.
\]
Thus,
\[
\lim_{t\to\infty} \Str (e^{-tD^2}) = \Str(P_0) = \ind(D^+).
\]
It remains to check that $\Str (e^{-tD^2})$ is independent of
$t>0$. For this, we note that the heat equation yields
\[
\lfrac d{dt} e^{-tD^2} = - D^2 e^{-tD^2}.
\]
The basic trace estimate then implies that $\Str(e^{-tD^2})$ is
differentiable in $t$ with
\[
\lfrac d{dt} \Str(e^{-tD^2}) = - \Str(D^2 e^{-tD^2}).
\]
However, since $D$ anti-commutes with $\gs$, we infer from the
trace property that
\[
\Str(D^2 e^{-tD^2}) =  \Tr(\gs D^2 e^{-tD^2}) = -\Tr (D\gs D
e^{-tD^2}) = - \Tr(\gs e^{-t D^2}D^2)=-\Str (D^2e^{-tD^2}),
\]
so that indeed $\Str (D^2e^{-tD^2})=0$.
\end{proof}

\noindent\textbf{Heat Kernel Asymptotics.} So far, the treatment
of the heat kernel has been of a functional analytic nature. The
only input are the basic properties of elliptic differential
operators on closed manifolds as in Theorem \ref{SelfAdjEllDef}.
However, one important missing piece is the analysis of the heat
trace as $t\to 0$, which requires further work. We summarize the
following from \cite[Sec.'s 1.8 \& 1.9]{Gil}.

\begin{theorem}\label{HeatTrace}
Let $M$ be a closed manifold, and let $H\in\sP_{s,e}^2(M,E)$ have
positive definite leading symbol. Let $D:C^\infty(M,E)\to
C^\infty(M,E)$ be an auxiliary differential operator of order
$d\ge 0$, and let $k_t(x,y)$ denote the kernel of $De^{-tH}$.
\begin{enumerate}
\item There exists an asymptotic expansion
\[
k_t(x,x) \sim \sum_{n=0}^\infty
t^{\frac{n-m-d}{2}}e_n(x),\quad\text{as }t\to 0,
\]
with $e_n\in C^\infty\big(M,\End(E)\big)$ such that $e_n(x)$ is
locally computable from the total symbols of $H$ and $D$ near
$x$. If $n+d$ is odd, then $e_n=0$.
\item The trace of $De^{-tH}$ admits an asymptotic expansion
\[
\Tr(De^{-tH})=\int_M\tr_E\big[k_t(x,x)\big]\vol_M(x) \sim
\sum_{n=0}^\infty t^{\frac{n-m-d}{2}}a_n(D,H) ,\quad\text{as }t\to
0.
\]
The asymptotic expansion can be differentiated in $t$, and the
$a_n$ are given by
\[
a_n(D,H) = \int_M \tr_E\big[ e_n\big]\vol_M.
\]
\end{enumerate}
\end{theorem}

\noindent\textbf{The Index Density.} As a consequence of the
McKean-Singer formula and the asymptotic expansion of the heat
trace, we get the following result of \cite{McKSin} and
\cite{ABP}.

\begin{theorem}\label{LocIndPrep}
Let $e_n(x)$ be the coefficient appearing in the asymptotic
expansion of $e^{-tD^2}(x,x)$ as in Theorem \ref{HeatTrace}. Then
\[
\int_M \str_E[e_n]\vol_M=  \begin{cases} \ind(D^+), &\text{\rm if
}
n=\dim M, \\
\hphantom{\ind}0,&\text{\rm if }n< \dim M.
\end{cases}
\]
\end{theorem}

\begin{remark*}\quad\nopagebreak
\begin{enumerate}
\item For aparent reasons, the differential form $\str_E[e_m]\vol_M$
with $m=\dim M$ is called the {\em index density} of $D^+$.
\item As mentioned in Theorem \ref{HeatTrace}, the coefficients
$e_n$ vanish if $n$ is odd. This implies that the index of $D^+$
vanishes if $M$ is odd dimensional.
\end{enumerate}
\end{remark*}

We also want to note an important consequence of the local nature
of the $e_n$. Assume that two operators $D\in \sP_{e,s}^1(M,E)$
and $D'\in \sP_{e,s}^1(M,E')$ are \emph{locally equivalent}. This
means that for every $x\in M$, there exists a neighbourhood $U$ of
$x$ and a local isometry $\gF: E|_U \to E'|_U$ such that
\begin{equation}\label{LocalEquiv}
\gF\circ D\circ \gF^{-1} = D'\quad\text{over $U$}.
\end{equation}
Since the sections $e_n$ and $e_n'$ in the respective asymptotic
expansions of the heat kernels are computable from the total
symbols over $U$, relation \eqref{LocalEquiv} implies that
$\str_E[e_n]|_U = \str_{E'}[e_n']|_U$. This has the following
consequence.

\begin{cor}\label{IndCoverTwist}
Let $M$ be a closed manifold, $E\to M$ a Hermitian vector bundle,
and let $D\in \sP_{e,s}^1(M,E)$ be a $\Z_2$-graded operator.
\begin{enumerate}
\item If $\pi:\widehat M\to M$ is a $k$-fold regular cover,
and
\[
\widehat D: C^\infty(\widehat M,\pi^*E)\to C^\infty(\widehat
M,\pi^*E)
\]
is the natural lift of $D$ to $\widehat M$, then
\[
\ind(\widehat D^+) = k\cdot \ind(D^+).
\]
\item If $A$ is a flat connection on a Hermitian vector bundle $F$
of rank $k$ over $M$, and if $D_A$ is the operator $D$ twisted by
$F$ and $A$, then
\[
\ind\big(D_A^+\big) = k\cdot \ind(D^+).
\]
\end{enumerate}
\end{cor}

\begin{proof}
For (i) one notes that the $\widehat e_n$ are the lifts to
$\widehat M$ of the $e_n$. Thus,
\[
\int_{\widehat M} \str_{\pi^*E}[\widehat e_n]\vol_{\widehat M}=
\int_{\widehat M} \pi^*\big(\str_E[e_n]\vol_M\big) = k\cdot
\int_M \str_E[e_n]\vol_M,
\]
where we have used that $\vol(\widehat M)=k\cdot \vol(M)$. For
(ii) note that choosing local trivializations for $F$ which are
parallel with respect to $A$, one finds that the flat bundle $F$
is locally isomorphic to $M\times \C^k$ endowed with the trivial
connection. This yields that $D_A$ is locally equivalent to
$D^{\oplus k}$.
\end{proof}

As we have seen in Proposition \ref{Sign=Ind}, the signature of a
closed manifold equals the index of an elliptic differential
operator of first order. Hence, Corollary \ref{IndCoverTwist}
proves our earlier assertion that the twisted signatures do not
contain new information other than the rank of the twisting
bundle.

\begin{cor}\label{SignMultClosed}
Let $M$ be a closed, even dimensional manifold.
\begin{enumerate}
\item If $\pi:\widehat M\to M$ is a $k$-fold regular cover,
then
\[
\Sign(\widehat M) = k\cdot\Sign(M)
\]
\item If $A$ is a flat connection on a Hermitian vector bundle $E$
of rank $k$ over $M$, then
\[
\Sign_A(M) = k\cdot\Sign(M).
\]
\end{enumerate}
\end{cor}

\subsection{Geometric Dirac Operators and the Local Index
Theorem}\label{GeomDirac}

Theorem \ref{LocIndPrep} is the starting point for local index
theory. Since the coefficients $e_n$ are---in principal---locally
computable, a strategy to prove the Atiyah-Singer Index Theorem
is to identify the index density $\str_E[e_m]\vol_M$ with a
Chern-Weil representative of an appropriate characteristic class.
Note, however, that Chern-Weil classes are expressions in the
curvature, whereas the $e_n$ a priori contain higher order
derivatives of the connection. Nevertheless, this strategy works
for an important class of differential operators, which we
describe briefly.\\

\noindent\textbf{Geometric Dirac Operators.} Let $(M,g)$ be an
oriented, Riemannian manifold, and let $E\to M$ be a Hermitian
vector bundle. $E$ is called a \emph{Clifford module} if it is
endowed with a bundle map $c: T^*M\to \End(E)$ such that for every
$\xi\in T^*M$
\begin{equation}\label{Clifford}
c(\xi)^2=-|\xi|^2_g\id_E \quad\text{and}\quad c(\xi)^*=-c(\xi).
\end{equation}
A Clifford module is called $\Z_2$-graded if $E$ is $\Z_2$-graded
and if $c(\xi)$ is an odd element of $\End(E)$ for all $\xi\in
T^*M$. Provided that $E$ is endowed with suitable connection, we
can construct a natural first order differential operator:

\begin{dfn}\label{GeomDiracOp}
Let $E$ be a Clifford module over a Riemannian manifold $(M,g)$.
\begin{enumerate}
\item A connection $\nabla^E$ on $E$ which is compatible with the
metric is called a {\em Clifford connection} if for all $e\in
C^\infty(M,E)$ and $\xi\in\gO^1(M)$
\[
\big[\nabla^E,c(\xi)\big]e= \nabla^E\big(c(\xi)e\big) -
c(\xi)\nabla^Ee = c(\nabla^g\xi)e,
\]
where $\nabla^g$ is the Levi-Civita connection acting on forms.
\item If $\nabla^E$ is a Clifford connection, we define the
associated \emph{geometric Dirac operator} by
\[
D:= c\circ \nabla^E: C^\infty(M,E)\to C^\infty(M,E).
\]
Here, we are viewing the Clifford structure as a bundle map
$c:T^*M\otimes E\to E$.
\item A geometric Dirac operator is called {\em $\Z_2$-graded} if $E$
is a $\Z_2$-graded Clifford module, and the Clifford connection
respects the splitting $E=E^+\oplus E^-$.
\end{enumerate}
\end{dfn}

\begin{remark}\label{GeomDiracOpRem}
One often defines a \emph{Dirac operator} to be a formally
self-adjoint elliptic operator whose square is a generalized
Laplacian. It is straightforward to check that geometric Dirac
operators have this property. However, not every Dirac operator is
a geometric one. In Section \ref{LocalFamIndex} we sketch how to
associate Dirac operators to generalized connections, in
particular Clifford {\em superconnections}, thereby obtaining a
larger class of Dirac operators, see also \cite[Prop. 3.42]{BGV}.
\end{remark}

\noindent\textbf{Canonical Structures on Clifford Modules.} As
pointed out, the index density is a purely local object. The
reason for studying geometric Dirac operators is that Clifford
modules have a canonical local structure which we now describe
briefly. For proofs we refer to \cite[Sec.'s 3.2 \& 3.3]{BGV}.

We assume from now on that $m=\dim M$ is even. Some aspects of
the odd dimensional case are contained in Section
\ref{Hirzebruch} and Section \ref{VarLocalInd} in a special case.
Let us further assume for the moment that $M$ is \emph{spin}.
Without going into the details of the definition and the
topological restrictions that this imposes on $M$, we note that it
implies that there exists a unique irreducible Clifford module
$S$ of rank $2^{m/2}$ over $M$, called the {\em spinor module}.
It follows from the representation theory of Clifford algebras
that
\begin{equation}\label{EndS}
\End(S) = \cl(T^*M).
\end{equation}
The Clifford module $S$ is naturally $\Z_2$-graded, and
\eqref{EndS} is an isomorphism of $\Z_2$-graded algebras. Here,
the grading on $\cl(T^*M)$ is the one induced via the symbol map
from the even/odd grading on differential forms. Now, every
$\Z_2$-graded Clifford module $E$ can be decomposed as
$E=S\otimes W$, where $W$ carries a trivial Clifford structure.
Moreover,
\begin{equation}\label{TwistBundle}
\End(W)=\End_{\cl}(E):=\Bigsetdef{T\in \End(E)}{[T,c(\ga)]_s =
0\quad\text{for all }\ga\in T^*M},
\end{equation}
and
\begin{equation}\label{TwistBundleDecomp}
\End(E)= \cl(T^*M)\otimes_s \End_{\cl}(E).
\end{equation}
Here, $[.,.]_s$ and $\otimes_s$ are the commutator and the tensor
product in the $\Z_2$-graded sense. If $\gs$ denotes the grading
operator of $E$, there exists a unique decomposition
\begin{equation}\label{TwistGradOp}
\gs=\tau\otimes \gs_W\in \cl(T^*M)\otimes_s \End_{\cl}(E),
\end{equation}
where $\tau:=i^{m/2}c(\vol_M)$ is the chirality operator, and
$\gs_W$ is a grading operator on $W$. Now, if $T\in \End(W)$,
then one verifies that
\begin{equation}\label{TwistTrace}
\str_W(T)=\tr_W[\gs_W T] = \frac{1}{\rk(S)}
\tr_E\big[(\tau^2\otimes\gs_W) T\big] =
\frac{1}{2^{m/2}}\str_E[\tau T].
\end{equation}
Moreover, the spinor module $S$ comes equipped with a canonical
Clifford connection $\nabla^S$ which is induced by the Levi-Civita
connection $\nabla^g$ via \eqref{EndS}. From this one gets a 1-1
correspondence between Clifford connections $\nabla^E$ on $E$ and
connections of the form $\nabla^S\otimes 1+ 1\otimes \nabla^W$,
where $\nabla^W$ is a Hermitian connection on $W$. The curvature
$F_{\nabla^W}$ satisfies
\begin{equation}\label{TwistCurv}
F_{\nabla^W}= F_{\nabla^E} - R^S\in
\gO^2\big(M,\End_{\cl}(E)\big),
\end{equation}
where for any orthonormal frame $\{e_i\}$ of $TM$
\begin{equation}\label{CliffordCurv}
R^S := \lfrac 18 g\big(R^g(e_i,e_j)e_k,e_l\big) e^i\wedge
e^j\otimes c(e^k)c(e^l).
\end{equation}
Here, $R^g\in \gO^2\big(M,\End(TM)\big)$ is the curvature tensor
of the Levi-Civita connection.

We now note that the right hand sides of \eqref{TwistBundle},
\eqref{TwistTrace} and \eqref{TwistCurv} can be defined globally
on $M$ without referring to the spinor module $S$. Thus, we can
introduce corresponding objects also in the case that $M$ is not
spin. In particular, the definition of $\End_{\cl}(E)$ in
\eqref{TwistBundle} is meaningful for every Clifford module $E$.

\begin{dfn}
Let $M$ be an $m$-dimensional manifold, where $m$ is even, and let
$E$ be a $\Z_2$-graded Clifford module over $M$, endowed with a
Clifford connection $\nabla^E$.
\begin{enumerate}
\item Let $R^S$ be defined as in \eqref{CliffordCurv}. Then we can
decompose
\[
F_{\nabla^E} = R^S + F^{E/S},\quad  \text{where}\quad F^{E/S}\in
\gO^2\big(M,\End_{\cl}(E)\big).
\]
The term $F^{E/S}$ is called the \emph{twisting curvature} of $E$.
\item If $T\in C^\infty\big(M,\End_{\cl}(E)\big)$, then its
\emph{relative supertrace} is defined as
\[
\str_{E/S}(T):= \frac 1{2^{m/2}}\str_E(\tau T).
\]
\item The \emph{relative Chern character form} of $E$ is given by
\begin{equation}\label{RelChern}
\ch_{E/S}(E,\nabla^E):= \str_{E/S}\big[\exp\big(\lfrac{i}{2\pi}
F^{E/S}\big)\big]\in \gO^{ev}(M).
\end{equation}
\end{enumerate}
\end{dfn}

\begin{remark*}
It follows from \eqref{TwistTrace} and \eqref{TwistCurv} that if
$M$ is spin so that $E=S\otimes W$, then
\[
\ch_{E/S}(E,\nabla^E) = \str_W\big[\exp\big(\lfrac{i}{2\pi}
F^W\big)\big].
\]
In particular, if $W$ is ungraded this coincides with the Chern
character form of $W$ as in Definition \ref{ChernCharDef}. In
general, if $W=W^+\oplus W^-$, the relative Chern character is the
difference of the Chern characters of $W^+$ and $W^-$, see also
\eqref{ChernSuperConnRem}.
\end{remark*}

\noindent\textbf{The Local Index Theorem.} We can now state a
version of the local index theorem as in \cite[Thm. 4.2]{BGV}.

\begin{theorem}[Patodi, Gilkey]\label{LocalIndex}
Let $M$ be a closed, oriented Riemannian manifold of even
dimension $m$, and let $E\to M$ be a $\Z_2$-graded Clifford module
with Clifford connection $\nabla^E$ and Dirac operator $D$. Let
$e_n(x)$ be the coefficient appearing in the asymptotic expansion
of $e^{-tD^2}(x,x)$ as in Theorem \ref{HeatTrace}. Then
\begin{equation}\label{LocalIndexEq}
\str_E[e_n]\vol_M(x) = \begin{cases} \big(\widehat
A(TM,\nabla^g)\wedge \ch_{E/S}(E,\nabla^E)\big)_{[m]}, &\text{\rm
if }
n=\dim M, \\
\hphantom{ \big(\widehat A(TM,\nabla^g)}\,0,&\text{\rm if }n< \dim
M,
\end{cases}
\end{equation}
where the Hirzebruch $\widehat A$-form is as in Definition
\ref{PontLAHatDef}, and $(\ldots)_{[m]}$ means taking the $m$-form
part of a differential form.
\end{theorem}

\begin{remark}\label{GetlzerLocInd}
There is a stronger version of the local index theorem due to E.
Getzler, see \cite{Get83} and \cite[Thm. 4.1]{BGV}, which we also
want to recall. Let $\boldsymbol{\gs}:\cl(T^*M)\to \gL^\bullet
T^*M$ be the symbol map \eqref{SymbolMapDef}, and use this to
endow $\cl(T^*M)$ with a $\Z$-grading. With respect to this let
$\cl_n(T^*M)$ be the subbundle of $\cl(T^*M)$ of elements of
degree $\le n$. Then it can be shown that
\[
e_n \in C^\infty\big(M,\cl_n(T^*M)\otimes_s\End_{\cl}(E)\big).
\]
The stronger version of the local index theorem is the formula
\[
(4\pi)^{m/2}\sum_{n\le m}\boldsymbol{\gs}(e_n) =
{\det}^{1/2}\left(\frac{R^g/2}{\sinh(R^g/2)}\right)\wedge
\exp(-F^{E/S}),
\]
where $R^g$ is the Riemann curvature tensor. For the definition
of the right hand side, see Appendix \ref{CharClass}, in
particular \eqref{SqrtDet} and Definition \ref{PontLAHatDef}. Now
the supertrace of elements in $\cl(T^*M)$ vanishes away from
degree $m$. Hence, Theorem \ref{LocalIndex} follows by computing
the supertrace of $\tau$ and taking into account the powers of
$\lfrac{i}{2\pi}$ appearing in our definition of $\widehat A$ and
$\ch_{E/S}$. In Section \ref{VarLocalInd}, we will sketch a proof
of a variation formula for the Eta invariant based on this more
general local index theorem.
\end{remark}

A direct consequence of the local index theorem \ref{LocalIndex}
and Theorem \ref{LocIndPrep} is the famous Atiyah-Singer Index
Theorem for geometric Dirac operators in its cohomological
version.

\begin{theorem}[Atiyah-Singer]\label{IndThm}
Let $M$ be a closed, oriented Riemannian manifold of even
dimension $m$, and let $E\to M$ be a $\Z_2$-graded Clifford module
with Clifford connection $\nabla^E$ and Dirac operator $D$. Then
\[
\ind(D^+) = \int_M \widehat A(TM,\nabla^g)\wedge
\ch_{E/S}(E,\nabla^E).
\]
In terms of characteristic classes,
\[
\ind(D^+) = \big\langle\widehat A(TM)\cup
\ch_{E/S}(E),[M]\big\rangle.
\]
\end{theorem}

\subsection{Hirzebruch's Signature Theorem}\label{Hirzebruch}

A special case of the Atiyah-Singer Index Theorem is one of its
predecessors, the Hirzebruch Signature Theorem. It arises if the
geometric Dirac operator in question is a twisted signature
operator. Therefore, we now collect some details about the
structure of the exterior algebra as a Clifford module.\\

\noindent\textbf{Clifford Structures on the Exterior Algebra.} Let
$M$ be an $m$-dimensional closed, oriented Riemannian manifold.
For the moment we do not assume that $m$ is even. We consider the
Clifford structure
\[
c:T^*M\to\End\big(\gL^\bullet T^*M\big),\quad c(\xi)=
\emu(\xi)-\imu(\xi),
\]
see \eqref{CliffDef}. The Levi-Civita connection $\nabla^g$
acting on forms is a Clifford connection, and we get a geometric
Dirac operator
\[
d+d^t=c\circ \nabla^g:\gO^\bullet(M)\to\gO^\bullet(M).
\]
There are two natural gradings on $\gL^\bullet T^*M$, one given
by the even/odd grading and one given by the chirality operator
$\tau$. We know from Lemma \ref{TauProp} that if $m$ is even,
Clifford multiplication is odd with respect to both gradings.
However, if $m$ is odd, Clifford multiplication commutes with
$\tau$ so that in this case we do not get a $\Z_2$-graded
Clifford module. Moreover, if $M$ is spin and even dimensional,
then
\[
\gL^* T^*M \cong \cl(T^*M) = \End(S) = S\otimes S^*,
\]
which means that the twisting bundle is isomorphic to the dual
bundle $S^*$ of $S$. This motivates the following

\begin{dfn}\label{RightCliffDef}
Let $M$ be a manifold of dimension $m$, not necessarily even. We
define a \emph{transposed Clifford multiplication}
\begin{equation*}
\widehat c: T^*M\to \End\big(\gL^\bullet T^*M\big), \quad
\widehat c(\xi):= \emu(\xi)+\imu(\xi).
\end{equation*}
\end{dfn}

The transposed Clifford multiplication has the following
properties.

\begin{lemma}\label{RightCliffLem}
Let $\{e^i\}$ be a local orthonormal frame for $TM$.
\begin{enumerate}
\item With the obvious abbreviations, we have
\[
\widehat c^i \widehat c^j + \widehat c^j\widehat c^i =
2\gd^{ij},\quad\text{and}\quad c^i \widehat c^j + \widehat c^j
c^i = 0.
\]
\item If we define $\widehat \tau:= i^{[\frac{m+1}2]}\widehat
c^1\ldots \widehat c^m$, then
\[
\widehat \tau^2 = (-1)^m,\quad\text{and}\quad \widehat \tau =
\tau\circ (-1)^\nu,
\]
where $\nu:\gL^\bullet T^*M\to \N$ is the number operator given by
$\nu(\go)=k$ if $\go\in \gL^kT^*M$.
\item Let $\widehat \cl(T^*M)$ denote the subbundle of $\End\big(\gL^\bullet
T^*M\big)$ generated by transposed Clifford multiplication. Then
\[
\widehat \cl(T^*M) = \End_{\cl}\big(\gL^\bullet T^*M\big),
\]
where $\End_{\cl}$ is defined as in \eqref{TwistBundle} with
respect to the even/odd grading.
\end{enumerate}
\end{lemma}

We shall not include the proof which is straightforward but a bit
tedious. Part (i) and (iii) of Lemma \ref{RightCliffLem} can be
found in \cite[p. 144]{BGV}. Part (ii) can be easily proved by
induction on $m$. However, we want to point out that part (iii)
implies that we have an isomorphism of $\Z_2$-graded algebras,
\begin{equation}\label{FormsSuperDecomp}
\End\big(\gL^\bullet T^*M\big) = \cl(T^*M) \otimes_s \widehat
\cl(T^*M),
\end{equation}
which is \eqref{TwistBundleDecomp} translated to the case at hand.
Moreover, part (ii) of Lemma \ref{RightCliffLem} shows that the
decomposition \eqref{TwistGradOp} in the case at hand is
\[
(-1)^\nu = \tau\otimes \widehat \tau.
\]

\begin{remark*}
In the case that $m$ is even, one might expect that the
decomposition \eqref{TwistBundleDecomp} for $\End\big(\gL^\bullet
T^*M\big)$ with respect to the $\tau$-grading is given by
\eqref{FormsSuperDecomp} together with the grading operator
$\tau\otimes 1$. However, with respect to this, the endomorphism
$\widehat c(\ga)$ for $\ga \in T^*M$ is even, and this is
incompatible with part (i) of Lemma \ref{RightCliffLem}. To stay
in the $\Z_2$-graded formalism, one would have to consider yet
another kind of Clifford multiplication, namely
\[
\widetilde c:= \widehat c\circ (-1)^{\nu}:T^*M \to
\End\big(\gL^\bullet T^*M\big).
\]
This generates a subalgebra $\widetilde \cl(T^*M)$ of
$\End\big(\gL^\bullet T^*M\big)$ of purely even degree with
respect to $\tau$ so that
\[
\End\big(\gL^\bullet T^*M\big) = \cl(T^*M) \otimes_s \widetilde
\cl(T^*M).
\]
Fortunately, in the discussion to follow, we are interested only
in elements of $\widetilde \cl(T^*M)$, respectively $\widehat
\cl(T^*M)$, which are of even with respect to the even/odd
grading. For elements of this form, $\widetilde c$ and $\widehat
c$ coincide up to sign. More precisely, if $\{e_i\}$ is a local
frame for $TM$, then for all $k\le m/2$
\[
\widetilde c^{i_1}\ldots \widetilde c^{i_{2k}} = (-1)^k \widehat
c^{i_1}\ldots \widehat c^{i_{2k}}.
\]
Hence, even if it is incorrect from a formal point of view, we
use the transposed Clifford multiplication $\widehat c$ also in
the case that $\gL^\bullet T^*M$ is graded by $\tau$.
\end{remark*}

\noindent\textbf{Traces of the Exterior Algebra.} Let
$\pi_0:\cl(T^*M)\to \C$ be the projection onto the subalgebra
$\C\subset \cl(T^*M)$. One easily verifies that
$\big[\cl(T^*M),\cl(T^*M)\big]\cap \C = \{0\}$, so that we can
define a trace on $\cl(T^*M)$ by
\begin{equation*}\label{CliffordTraceDef}
\tr_{\cl}:=  2^{m/2} \pi_0 : \cl(T^*M)\to \C,
\end{equation*}
see \cite[Thm. 1.8]{Get83}. In the same way, we get a trace
$\widehat {\tr_{\cl}}$ on $\widehat\cl(T^*M)$. Note that in the
case that $m$ is even, the natural supertrace of \cite[Prop.
3.21]{BGV} is given by
\[
\str_{\cl}:= \tr_{\cl}\circ \tau: \cl(T^*M)\to \C.
\]

\begin{prop}\label{TraceSuperDecomp}
With respect to the decomposition \eqref{FormsSuperDecomp} we have
\[
\tr_{\gL^\bullet} = \tr_{\cl}\otimes \widehat {\tr_{\cl}},
\]
where $\tr_{\gL^\bullet}$ is the natural trace on
$\End(\gL^\bullet T^*M)$.
\end{prop}

\begin{remark*}
We include a proof, since the treatment in \cite{BGV} considers
only the case that $m$ is even. There are some non-trivial sign
difficulties involved, since \eqref{FormsSuperDecomp} involves
the graded tensor product, whereas $\tr_{\cl}$ and $\widehat
{\tr_{\cl}}$ are traces rather than supertraces. Yet, for elements
of pure degree,
\begin{equation}\label{SuperTraceComp}
\begin{split}
\tr_{\cl}\otimes \widehat {\tr_{\cl}}\big((a\otimes \hat
a)(b\otimes \hat b)\big) &= (-1)^{|\hat a| |b|} \tr_{\cl}(ab)
\widehat {\tr_{\cl}}(\hat a \hat b)\\
&= (-1)^{|\hat a| |b|} \tr_{\cl}(ba)\widehat {\tr_{\cl}} (\hat b
\hat a)\\ &= (-1)^{|\hat a| |b| + |\hat b||a|} \tr_{\cl}\otimes
\widehat {\tr_{\cl}}\big((b\otimes \hat b)(a\otimes \hat a)\big).
\end{split}
\end{equation}
Moreover, for $\tr_{\cl}(ab)$ and $\widehat {\tr_{\cl}}(\hat a
\hat b)$ to be non-zero it is necessary that $|a|=|b|$ and $|\hat
a|=|\hat b|$. In this case the sign in \eqref{SuperTraceComp} is
always $+1$. Hence, $\tr_{\cl}\otimes \widehat {\tr_{\cl}}$ is
indeed a trace.
\end{remark*}

\begin{proof}[Proof of Proposition \ref{TraceSuperDecomp}]
Since the assertion is local, it suffices to consider an
$m$-dim\-en\-sio\-nal Euclidean vector space $V$. Let us first
consider the case $V=\R$, and let $e$ be a unit vector. Then
\[
\gL^\bullet \R = \C\oplus \C e \cong \C^2.
\]
Consider the following elements of $\End(\C^2)$
\[
c:= \begin{pmatrix}0 &-1\\ 1 &0\end{pmatrix},\quad \widehat
c:=\begin{pmatrix}0 &1\\ 1 &0\end{pmatrix},\quad\text{and}\quad
n:= \begin{pmatrix}1 &0\\ 0 &-1\end{pmatrix}.
\]
Then $\cl(\R)\subset \End(\C^2)$ is the algebra generated by $\Id$
and $c$, and $\widehat \cl(\R)$ is generated by $\Id$ and
$\widehat c$. Moreover, $c\widehat c = -\widehat c c = -n$, which
implies that all monomials in $c$ and $\widehat c$ have vanishing
trace except $c^0\widehat c^0 = \Id$. In this case,
\[
\tr(c^0\widehat c^0)=\tr(\Id)=2 = \tr_{\cl}(c^0)\widehat{
\tr_{\cl}}(\widehat c^0),
\]
which yields the claimed formula in the case that $V=\R$.

We now assume that the claim holds for and $m$-dimensional vector
space $V$, and want to prove it for $V\oplus \R$. Let $\{e_i\}$
be an orthonormal basis for $V$, and let $e$ be a unit vector in
$\R$. For an ordered multi-index $A=(i_1<\ldots<i_k)$ we consider
the following sets of generators
\[
e_A:= e_{i_1}\wedge\ldots\wedge e_{i_k}\in \gL^\bullet V,\quad
c_A:= c(e_{i_1})\ldots c(e_{i_k})\in \cl(V),
\]
and for $\ga \in \{0,1\}$,
\[
e_{A,\ga}:= e_A\wedge e^{\ga}\in \gL^\bullet (V\oplus\R),\quad
c_{A,\ga} := c_A c^\ga \in \cl(V\oplus \R),
\]
where $c:=c(e)$. In the same way we define $\widehat c_A\in
\widehat \cl(V)$ and $\widehat c_{A, \ga}\in \widehat \cl(V\oplus
\R)$. Then a short computation shows that
\[
\begin{split}
\tr_{\gL^\bullet(V\oplus \R)} (c_{A,\ga}\widehat c_{B,\gb}) &=
\sum_{|C|\le m}\sum_{\gamma\in\{0,1\}}
\Scalar{c_{A,\ga}\widehat c_{B,\gb}e_{C,\gamma}}{e_{C,\gamma}}\\
&= \sum_{|C|\le m}\sum_{\gamma\in\{0,1\}} (-1)^{|B||\ga|}
\Scalar{(c_A\widehat c_B)(c^\ga\widehat c^\gb
)e_{C,\gamma}}{e_{C,\gamma}}\\
&= \sum_{|C|\le m}(-1)^{|B||\ga|+|C|(|\ga|+|\gb|)}
\Scalar{(c_A\widehat c_B)e_C}{e_C}\tr_{\gL^\bullet
\R}(c^\ga\widehat c^\gb).
\end{split}
\]
From the case $m=1$ we know that $\tr_{\gL^\bullet
\R}(c^\ga\widehat c^\gb)=0$ if $(\ga,\gb)\neq (0,0)$. In this
case the factor in the last line above is $+1$ and so
\begin{equation}\label{TraceSuperDecomp:1}
\tr_{\gL^\bullet(V\oplus \R)} (c_{A,\ga}\widehat c_{B,\gb}) =
\tr_{\gL^\bullet V}(c_A\widehat c_B) \tr_{\gL^\bullet \R}
(c^\ga\widehat c^\gb).
\end{equation}
Now by induction we have
\begin{equation}\label{TraceSuperDecomp:2}
\tr_{\gL^\bullet V}= \tr_{\cl(V)}\otimes
\widehat{\tr_{\cl(V)}},\quad \tr_{\gL^\bullet \R}=
\tr_{\cl(\R)}\otimes \widehat{\tr_{\cl(\R)}}.
\end{equation}
Moreover, as in \eqref{SuperTraceComp} one checks that with
respect to
\[
\cl(V\oplus \R)\otimes_s \widehat \cl(V\oplus \R) \cong
\big(\cl(V)\otimes_s\cl(V)\big)\otimes_s\big(\widehat
\cl(\R)\otimes_s\widehat\cl(\R)\big)
\]
one has
\[
\tr_{\cl(V\oplus \R)}\otimes \widehat{\tr_{\cl(V\oplus \R)}} =
\big(\tr_{\cl(V)}\otimes \widehat{\tr_{\cl(V)}}\big)\otimes
\big(\tr_{\cl(\R)}\otimes \widehat{\tr_{\cl(\R)}}\big),
\]
which together with \eqref{TraceSuperDecomp:1} and
\eqref{TraceSuperDecomp:2} proves the assertion for $V\oplus \R$.
\end{proof}

\noindent\textbf{Local Index Density and the Signature Theorem.}
As in the above proof let $V:=T^*_xM$ for some $x\in M$, and let
$R$ be an element in the Lie algebra $\cso(V)\subset \End(V)$. Let
$V_\C$ be the complexification of $V$. Then $iR\in \End(V_\C)$ is
a self-adjoint endomorphism, and we can define $\cosh(iR)\in
\End(V_\C)$ via the spectral theorem. Since the eigenvalues of
$iR$ are real, and $\cosh$ is a positive function on $\R$, we can
define
\[
{\det}^{1/2}\big(\cosh(iR)\big):=\sqrt{\det\big(\cosh(iR)\big)}.
\]
It follows from the spectral theorem that
\[
{\det}^{1/2}\big(\cosh(iR)\big) = \exp\big(\lfrac 12\tr\big[\log
\cosh(iR) \big]\big),
\]
which agrees with the definition in \eqref{SqrtDet}. Note,
however, that the context here is slightly different since we are
considering elements in $\cso(V)$ whereas in \eqref{SqrtDet} we
are considering elements of the algebra $(\gL^{\ev}\C^m)\otimes
\End(V_\C)$.

We then have the following version of \cite[Lem. 4.5]{BGV}

\begin{lemma}\label{FormsTwistChern}
Let $V$ be an $m$-dimensional oriented Euclidean vector space, and
let $R\in \cso(V)$. Define
\[
\widehat R^S := -\lfrac 14 \scalar{Re_i}{e_j} \widehat
c^i\widehat c^j\in \widehat \cl(V^*),
\]
where $\{e_i\}$ is any orthonormal basis for $V$, and $\widehat
c^i=\widehat c(e^i)$. Then
\[
\widehat{\tr_{\cl}}\big[\exp(i\widehat R^S)\big] = 2^{m/2}
{\det}^{1/2}\big(\cosh(iR/2)\big).
\]
\end{lemma}

\begin{proof}
Let $k\in \N$ be such that $m=2k$ or $m=2k+1$. Since
$R\in\cso(V)$ we can find an orthonormal basis $\{e_i\}$ such that
\begin{equation}\label{FormsTwistChern:1}
R(e_{2j-1}) = \gt_j e_{2j},\quad R(e_{2j}) = -\gt_j e_{2j-1},\quad
j=1,\ldots,k,\quad\text{and}\quad R(e_{2k+1})=0,
\end{equation}
where the last condition has to be considered as empty if $m$ is
even. Then
\[
\widehat R^S = -\lfrac 12\sum_{j\le k} \gt_j \widehat
c^{\,2j-1}\widehat c^{\,2j}.
\]
Since $c^{\,2i-1}\widehat c^{\,2i}$ and $c^{\,2j-1}\widehat
c^{\,2j}$ commute for $i\neq j$, one finds
\[
\exp(i\widehat R^S) = \prod_{j\le k} \exp\big((-i\gt_j/2) \widehat
c^{\,2j-1}\widehat c^{\,2j}\big) =\prod_{j\le
k}\Big(\cosh(\gt_j/2)-i\sinh(\gt_j/2)\widehat c^{\,2j-1}\widehat
c^{\,2j}\Big),
\]
where the last equality follows from the relation $(\widehat
c^{\,2j-1}\widehat c^{\,2j})^2 = -1$. By definition of the trace
on $\widehat\cl(V^*)$, we find
\[
\widehat{\tr_{\cl}}\big[\exp(i\widehat R^S)\big] =
2^{m/2}\prod_{j\le k}\cosh(\gt_j/2).
\]
On the other hand, for all $z\in \C$
\[
\cosh\left[\begin{pmatrix}0 &-z\\ z & 0 \end{pmatrix}\right]=
\cos(z)
\begin{pmatrix}1 &0\\ 0 & 1 \end{pmatrix},\quad\cosh(0)=1,
\]
from which it follows that
\[
\sqrt{\det\big(\cosh(iR/2)\big)} = \prod_{j\le k}
\sqrt{\cosh(\gt_j/2)^2} = \prod_{j\le k}\cosh(\gt_j/2).\qedhere
\]
\end{proof}

The next result can be found in \cite[p. 145]{BGV}. Although the
treatment there is only for $m$ even, one verifies without effort
that it holds for odd $m$ as well.

\begin{lemma}\label{FormsTwistCurv}
Let $(M,g)$ be an oriented Riemannian manifold with Riemann
curvature tensor $R^g\in\gO^2\big(M,\End(TM)\big)$. In a local
orthonormal frame $\{e_i\}$ for $TM$ write
\[
R^g = \lfrac12 R_{kij}^l e^i\wedge e^j \otimes (e_l\otimes
e^k),\quad\text{and}\quad R_{lkij}=g_{ln} R_{kij}^n =
g\big(R(e_i,e_j)e_k,e_l\big).
\]
Define $R^S\in \gO^2\big(M,\cl(T^*M)\big)$ as in
\eqref{CliffordCurv}, and
\[
\widehat R^S:= - \lfrac 18 R_{lkij}e^i\wedge e^j\otimes \widehat
c^k \widehat c^l \in \gO^2\big(M,\widehat \cl(T^*M)\big).
\]
Then, the curvature $R^{\gL^\bullet T^*M}$ of the induced
connection on $\gL^\bullet T^*M$ decomposes as
\[
R^{\gL^\bullet T^*M} = R^S +\widehat R^S.
\]
\end{lemma}

Lemma \ref{FormsTwistChern} extends to the case $R\in \gL^2
V\otimes \cso(V)$, where \eqref{SqrtDet} is used to define the
right hand side, see \cite[pp. 144--146]{BGV}. Then Lemma
\ref{FormsTwistChern} and Lemma \ref{FormsTwistCurv} imply the
following

\begin{prop}\label{SignIndDens}
Let $(M,g)$ be an oriented Riemannian manifold of dimension $m$,
and let $R^g$ be the Riemann curvature tensor. Then
\begin{equation*}
\widehat A(TM,\nabla^g)\wedge
\widehat{\tr_{\cl}}\big[\exp(\lfrac{i}{2\pi}\widehat R^S)\big] =
2^{m/2}\widehat L(TM,\nabla^g).
\end{equation*}
with the Hirzebruch $\widehat L$-form, see Definition
\ref{PontLAHatDef}.
\end{prop}

From Proposition \ref{SignIndDens} and Theorem \ref{IndThm}, we
obtain the index theorem for twisted signature operators.

\begin{theorem}[Atiah-Singer, Hirzebruch]\label{TwistSignIndThm}
Let $M$ be a closed, oriented Riemannian manifold of even
dimension $m$. Let $E\to M$ be a Hermitian vector bundle with
connection $A$.
\begin{enumerate}
\item The index of the twisted signature operator $D_A^+$ is
given by
\[
\ind(D_A^+)= 2^{m/2} \int_M \widehat L(TM,\nabla^g)\wedge
\ch(E,A).
\]
\item If $A$ is a flat connection, and $E$ has rank $k$, then
\begin{equation}\label{SignThm}
\Sign_A(M)= k \int_M L(TM,\nabla^g) = k\cdot \big\langle
L(TM),[M]\big\rangle,
\end{equation}
where $L(TM,\nabla^g)$ is the Hirzebruch $L$-form as in
\eqref{L-L-Hat}.
\end{enumerate}
\end{theorem}


\section{Manifolds with Boundary and the Eta Invariant}

\subsection{The Eta Function}

Let $(M,g)$ be a closed, oriented Riemannian manifold of
dimension $m$. Let $E\to M$ be a Hermitian vector bundle, and let
\[
D:C^\infty(M,E)\to C^\infty(M,E)
\]
be a formally self-adjoint elliptic differential operator of
first order, i.e., $D\in \sP_{s,e}^1(M,E)$. As $M$ is closed, the
growth of the eigenvalues of $D$ are controlled by
\eqref{EigenvalueAsymp}. This allows us to make the following
definition.

\begin{dfn}\label{EtaFctnDef}
The \emph{Eta function} of $D$ is defined for $s\in\C$ with $\Re
(s)>m$ as
\begin{equation*}
\eta(D,s):= \sum_{0\neq \gl\in\spec(D)} \frac{\sgn(\gl)}{|\gl|^s}.
\end{equation*}
\end{dfn}

Via a Mellin transform, the Eta function is related to the heat
operator $e^{-tD^2}$ in the following way
\begin{equation}\label{EtaMellin}
\eta(D,s) = \frac{1}{\gG\big(\lfrac{s+1}{2}\big)} \int_0^\infty
\Tr\big(De^{-tD^2}\big) t^{\frac{s-1}2}dt,\quad \Re(s)>m,
\end{equation}
where $\gG(s)$ is the Gamma function\footnote{See Appendix
\ref{CompGamma} for some facts about the Gamma function which we
will use freely.},
\[
\gG(s) := \int_0^\infty e^{-t}t^{s-1}dt,\quad \Re(s)>0.
\]
Note that for \eqref{EtaMellin} to exist we are using that
$De^{-tD^2}$ is trace class, and that there exist constants $c$
and $C$ such that for large $t$
\[
|\Tr(De^{-tD^2})|\le Ce^{-ct}.
\]
This follows from Proposition \ref{BasicTraceProp}, because using
the notation introduced there, we have
$\Tr(De^{-tD^2})=\Tr\big(De^{-tD^2}P_{(0,\infty)}\big)$. Thus,
\eqref{EtaMellin} yields a holomorphic function in the half plane
$\Re(s)>m$.

As we have noted in Theorem \ref{HeatTrace}, there is an
asymptotic expansion
\begin{equation}\label{EtaAsymp}
\Tr\big(De^{-tD^2}\big) \sim \sum_{n=0}^\infty
t^{\frac{n-m-1}{2}}a_n(D) ,\quad\text{as }t\to 0,
\end{equation}
where $a_n(D)$ is an integral over a quantity locally computable
from the total symbol of $D$. Substituting the asymptotic
expansion into \eqref{EtaMellin} and dividing the integration
into $\int_0^1+\int_1^\infty$ one easily verifies that for each
$N\ge 0$
\begin{equation}\label{EtaCont}
\eta(D,s) = \frac{1}{\gG\big(\lfrac{s+1}{2}\big)}
\Big(\sum_{n=0}^N \frac{2a_n(D)}{n-m+s} + h_N(s) \Big),
\end{equation}
where $h_N(s)$ is holomorphic in the half plane $\Re(s)>
m-(N+1)$. Since $\gG\big(\lfrac{s+1}{2}\big)^{-1}$ is an entire
function, one can use this to deduce

\begin{prop}\label{EtaMerom}
The Eta function $\eta(D,s)$ extends uniquely to a meromorphic
function on the whole plane with possible simple poles for $s\in
\setdef{m-n}{n\in \N}$.
\end{prop}

\noindent\textbf{Regularity at $\boldsymbol{s=0}$.} We note that
$\gG\big(\lfrac{s+1}{2}\big)^{-1}$ has no zeros to cancel the
possible poles. This is an important difference between the Eta
function and the Zeta function of, say, a generalized Laplacian,
see e.g. \cite[Prop. 9.35]{BGV}. Therefore, the following result
is very remarkable. For references we refer to Remark
\ref{EtaRegRem} below.

\begin{theorem}[Atiyah-Patodi-Singer, Gilkey]\label{EtaReg}
Let $D$ be a formally self-adjoint elliptic differential operator
of first order on a closed, Riemannian manifold $M$. Then the Eta
function $\eta(D,s)$ has no pole at $s=0$. If $D$ is a geometric
Dirac operator, $\eta(D,s)$ is holomorphic for $\Re(s)>-1/2$.
\end{theorem}

\begin{dfn}\label{EtaDef}
Using the result of Theorem \ref{EtaReg}, we can define the
\emph{Eta invariant} of $D$ as
\begin{equation*}
\eta(D):=\eta(D,0).
\end{equation*}
Moreover, we define the \emph{$\xi$-invariant} and the
\emph{reduced $\xi$-invariant} by
\begin{equation*}
\xi(D):= \frac{\eta(D)+\dim(\ker D)}{2},\quad \text{and}\quad
[\xi(D)]:= \xi(D) \mod\Z.
\end{equation*}
\end{dfn}

\begin{remark}\label{EtaRegRem}
To further stress the non-triviality of Theorem \ref{EtaReg}, we
want to give some historical remarks.
\begin{enumerate}
\item Atiyah, Patodi and Singer first deduced the regularity of
the Eta function at 0 from their proof of the index theorem for
elliptic differential operators of first order on manifolds with
boundary, see \cite[Thm. 3.10]{APS1}. The improved regularity for
geometric Dirac operators is \cite[Thm. 4.2]{APS1}.
\item Later, the same authors generalized the result to
pseudo-differential operators of arbitrary order on closed, odd
dimensional manifolds using $K$-theoretic arguments and the
regularity results for geometric Dirac operators, see \cite[Thm.
4.5]{APS3}.
\item In \cite{Gil81}, Gilkey was able to generalize Theorem
\ref{EtaReg} to the case of even dimensional manifolds, again for
the much larger class of formally self-adjoint elliptic
pseudo-differential operators, see also \cite[Sec. 3.8]{Gil}.
\item Gilkey \cite{Gil79} also initiated the study of the
regularity of the \emph{local Eta function} of an elliptic
operator $D$,
\begin{equation}\label{LocalEta}
\eta(D,s,x):= \frac{1}{\gG\big(\lfrac{s+1}{2}\big)} \int_0^\infty
\tr\big( k_t(x,x)\big) t^{\frac{s-1}2}dt,
\end{equation}
where $k_t(x,y)$ is the kernel of $De^{-tD^2}$. Studying various
examples, Gilkey found that $\eta(D,s,x)$ is in general not
regular at $s=0$.
\item Later Bismut and Freed were able to refine the result of
Atiyah-Patodi-Singer for geometric Dirac operators. They showed
using local index theory techniques, that the local Eta function
$\eta(D,s,x)$ of a geometric Dirac operator is holomorphic for
$\Re(s)>-2$, see \cite[Thm. 2.6]{BF2}. Their result implies that
for a geometric Dirac operator one can define the Eta invariant
directly by
\begin{equation}\label{EtaDirect}
\eta(D) = \frac1{\sqrt\pi}\, \int_0^\infty t^{-1/2}
\Tr\big(De^{-tD^2}\big) dt.
\end{equation}
\end{enumerate}
\end{remark}

\begin{lemma}\label{EtaProp}
Let $M$ and $N$ be closed, oriented Riemannian manifolds, with
Hermitian vector bundles $E$ over $M$ and $F$ over $N$. Let $D\in
\sP_{s,e}^1(M,E)$ and $B\in \sP_{s,e}^1(N,F)$.
\begin{enumerate}
\item Assume that there exists and isometry $\gf:M\to N$, and a
unitary bundle map $\gF:E\to F$ covering $\gf$ such that
\[
\gF\circ D=   B\circ \gF.
\]
Then $\eta(D)=\eta(B)$. In particular, if $M=N$, $E=F$ and
$\{D,\gF\}=0$, then $\eta(D)=0$.
\item If $M=N$, then
\[
\eta(D\oplus B) = \eta(D) + \eta(B).
\]
\item Assume that $D$ is $\Z_2$-graded with grading operator $\gs$ on
$E$. Consider the operator
\[
D\otimes 1 + \gs\otimes B \quad\text{on}\quad
C^\infty\big(M\times N,E\boxtimes F\big),
\]
with the fiber product $E\boxtimes F$ as in \eqref{FiberProd}.
Then
\[
\eta\big(D\otimes 1 + \gs\otimes B\big) = \ind(D^+)\cdot \eta(B).
\]
\end{enumerate}
\end{lemma}

\begin{proof}[Sketch of proof]
Part (i) and (ii) of the above result are immediate for
$\eta(D,s)$ for $\Re(s)$ large, since the whole spectrum has the
respective properties. By meromorphic continuation, they continue
to hold for $s=0$. We sketch a proof of (iii).

First note that $D\otimes 1$ and $\gs \otimes B$ anti-commute as
operators on $C^\infty\big(M\times N,E\boxtimes F\big)$.
Therefore,
\[
\big(D\otimes 1 + \gs\otimes B\big) e^{-t(D\otimes 1 + \gs\otimes
B)^2} = (D\otimes 1)e^{-t(D^2\otimes 1+1\otimes B^2)}+ (\gs\otimes
B) e^{-t(D^2\otimes 1+1\otimes B^2)}.
\]
To compute traces one may choose an orthonormal basis of
$L^2\big(M\times N,E\boxtimes F\big)$ of the form $\{\gf_i\otimes
\psi_j\}$ with $\gf_i\in C^\infty(M,E)$ and $\psi_j\in
C^\infty(N,F)$. Then one easily finds that
\begin{multline*}
\Tr\Big((D\otimes 1)e^{-t(D^2\otimes 1+1\otimes B^2)}\Big) +
\Tr\Big((\gs\otimes
B) e^{-t(D^2\otimes 1+1\otimes B^2)}\Big)\\
=\Tr_{L^2(M,E)}\big(De^{-tD^2}\big)
\Tr_{L^2(N,F)}\big(e^{-tB^2}\big) \\
+ \Str_{L^2(M,E)}\big(e^{-tD^2})\Tr_{L^2(N,F)}\big(B
e^{-tB^2}\big).
\end{multline*}
Since $D$ is $\Z_2$-graded, $\gs De^{-tD^2}= -De^{-tD^2}\gs$ and
thus,
\[
\Tr\big(De^{-tD^2}\big) = \Tr\big(\gs^2 De^{-tD^2}\big)= -
\Tr\big(\gs D e^{-tD^2}\gs\big) = -\Tr\big(De^{-tD^2}\big),
\]
where we use the trace property in the last equality. Therefore,
$\Tr\big(De^{-tD^2}\big) = 0$. Moreover, Theorem
\ref{McKeanSinger} asserts that for all $t>0$,
\[
\Str\big(e^{-tD^2}) = \ind(D^+).
\]
We conclude that for $\Re(s)$ large,
\[
\begin{split}
\eta\big(D\otimes 1 + \gs\otimes B,s\big) &=
\frac{1}{\gG\big(\lfrac{s+1}{2}\big)} \int_0^\infty
\Str\big(e^{-tD^2})\Tr\big(B e^{-tB^2}\big) t^{\frac{s-1}2}dt\\
&= \ind(D^+)\cdot \eta(B,s).
\end{split}
\]
By meromorphic continuation, part (iii) follows.
\end{proof}

\noindent\textbf{The Rho Function.} We now want to use the local
nature of the coefficients $a_n(D)$ appearing in \eqref{EtaCont}
to introduce the Rho function---respectively, the Rho invariant.

\begin{dfn}\label{RhoFctnDef}
Let $D\in \sP^1_{s,e}(M,E)$, and let $A$ be a flat connection on a
Hermitian vector bundle of rank $k$. Denote by $D_A$ the operator
$D$ twisted by $A$, and use the notation $D^{\oplus k}$ for the
operator $D$ twisted by the trivial flat bundle $\C^k$. We then
define the \emph{Rho function} of $D_A$ as
\[
\rho(D_A,s) = \eta(D_A,s) - \eta(D^{\oplus k},s),\quad s\in\C.
\]
Moreover, we define the \emph{Rho invariant} of $D_A$ as
\[
\rho(D_A):=\rho(D_A,0).
\]
\end{dfn}

From Theorem \ref{EtaReg} we know that the meromorphic functions
$\eta(D_A,s)$ and $\eta(D^{\oplus k},s)$ have no pole in 0. Thus,
the Rho invariant is well-defined. However, unlike in the case for
the individual Eta invariants, this already follows from the
local nature of heat trace asymptotics.

\begin{prop}\label{RhoIntProp}
Let $D_A$ and $D^{\oplus k}$, where $A$ is a flat
$\U(k)$-connection. Then the Rho function is holomorphic on the
whole plane, and for all $s\in \C$, we have
\begin{equation}\label{RhoMellin}
\rho(D_A,s) =  \frac{1}{\gG\big(\lfrac{s+1}{2}\big)} \int_0^\infty
\big[\Tr\big(D_Ae^{-tD_A^2}\big)-\Tr\big(D^{\oplus
k}e^{-t(D^{\oplus k})^2}\big)\big] t^{\frac{s-1}2}dt,
\end{equation}
in particular,
\[
\rho(D_A) = \frac1{\sqrt\pi}\, \int_0^\infty t^{-1/2}
\big[\Tr\big(D_Ae^{-tD_A^2}\big)-\Tr\big(D^{\oplus
k}e^{-t(D^{\oplus k})^2}\big)\big] dt.
\]
\end{prop}

\begin{proof}
For $\Re(s)>m$, we already know from \eqref{EtaMellin} that
$\rho(D_A,s)$ is holomorphic and that \eqref{RhoMellin} is the
correct formula. Moreover, as for the Eta function, we can split
up the integral and use that $\int_1^\infty$ extends to a
holomorphic function on $\C$. Concerning $\int_0^1$, we now make
the following observation: As we have already seen in the proof
of Corollary \ref{IndCoverTwist}, the operators $D_A$ and
$D^{\oplus k}$ are locally equivalent. Then Theorem
\ref{HeatTrace} implies that all coefficients of the asymptotic
expansions of $\Tr\big(D_Ae^{-tD_A^2}\big)$ and $\Tr\big(D^{\oplus
k}e^{-t(D^{\oplus k})^2}\big)$ as $t\to 0$ agree, since they are
local in the total symbols of the involved operators. Thus, for
all $N$ there exists a constant $C$ such that as $t\to 0$
\[
\big|\Tr\big(D_Ae^{-tD_A^2}\big)-\Tr\big(D^{\oplus
k}e^{-t(D^{\oplus k})^2}\big)\big| \le Ct^N.
\]
This shows that the integral
\[
\int_0^1 \big[\Tr\big(D_Ae^{-tD_A^2}\big)-\Tr\big(D^{\oplus
k}e^{-t(D^{\oplus k})^2}\big)\big] t^{\frac{s-1}2}dt
\]
exists and defines a holomorphic function for all $s\in \C$ with
$\Re(s)>-(N+1)$. Continuing in this way, one finds that
$\rho(D_A,s)$ is holomorphic for $s\in \C$.
\end{proof}

\subsection{Variation of the Eta Invariant}

While the Eta invariant is a spectral invariant and thus encodes
global information about the manifold and the operator, its
deformation theory turns out to be expressible in terms of heat
trace asymptotics, thus being a local quantity. We have included
further details in Appendix \ref{AppEtaVar}, and summarize only
briefly what we need here.\\

\noindent\textbf{Families of Operators.} Let $M$ be a closed
manifold of dimension $m$, and let $E\to M$ be a Hermitian vector
bundle. Let $U\in \R^p$ be open, and let $(D_u)_{u\in U}$ be a
$p$-parameter family of differential operators on $C^\infty(M,E)$
of order $d$. Choosing local frames for $E$ we can write locally
\begin{equation}\label{LocalChart}
D_u = \sum_{|\ga|\le d} a_\ga(x,u) \pd^{\ga}_x,
\end{equation}
where $\ga=(\ga_1,\ldots,\ga_m)$ is a multi-index,
$x=(x_1,\ldots,x_m)$ is a local coordinate chart on $M$, and the
$a_\ga(x,u)$ are matrix valued functions.

\begin{dfn}\label{SmoothFamily}
The $p$-parameter family $(D_u)_{u\in U}$ is called \emph{smooth}
if for every local frame and coordinate chart, the functions
$a_\ga(x,u)$ as in \eqref{LocalChart} are smooth, jointly in $x$
and $u$.
\end{dfn}

We also need a functional analytic consequence of the geometric
notion of smoothness we have just defined. Let $s, s'\in\R$ and
denote by $\sB(L_s^2,L_{s'}^2)$ the space of bounded linear
operators $L_s^2\to L_{s'}^2$ endowed with the operator norm
$\|.\|_{s,s'}$. The proof of the following result is
straightforward.

\begin{lemma}\label{SmoothFamilyMap}
Let $(D_u)_{u\in U}$ be a $p$-parameter family of differential
operator of order $d\ge 0$ which is smooth in the sense of
Definition \ref{SmoothFamily}. Then for all $s\in \R$, we have a
smooth map
\[
U\to\sB(L_{s+d}^2,L_s^2),\quad u\mapsto D_u.
\]
\end{lemma}

For more details on the following results we refer to Appendix
\ref{AppEtaVar}, in particular Proposition \ref{EtaDerApp} and
Proposition \ref{RedEtaDerApp}.

\begin{prop}
Let $(D_u)_{u\in \R}$ be a smooth one-parameter family of
formally self-adjoint elliptic operators of first order on
$C^\infty(M,E)$, and let $a_m(\lfrac{dD_u}{du},D_u^2)$ denote the
constant term in the asymptotic expansion of
\begin{equation}\label{EtaVarAsymp}
\sqrt t\,\Tr\big(\lfrac{dD_u}{du}e^{-tD_u^2}\big),\quad\text{as
}t\to 0.
\end{equation}
Then the following holds.
\begin{enumerate}
\item Assume that $\dim(\ker D_u)$ is constant. Then the meromorphic
extension of $\eta(D_u,s)$ is continuously differentiable in $u$,
and
\begin{equation}\label{EtaDer}
\lfrac{d}{du}\eta(D_u) = - \lfrac2{\sqrt\pi}\,
a_m(\lfrac{dD_u}{du},D_u^2).
\end{equation}
\item Without the assumption on $\ker(D_u)$, the reduced
$\xi$-invariant $[\xi(D_u)]\in \R/\Z$ is continuously
differentiable in $u$, and
\begin{equation}\label{RedEtaDer}
\lfrac d{du}[\xi(D_u)] = - \lfrac1{\sqrt\pi}\,
a_m(\lfrac{dD_u}{du},D_u^2).
\end{equation}
\end{enumerate}
\end{prop}

An immediate consequence of the local nature of
$a_m(\lfrac{dD_u}{du},D_u^2)$ is the following result, see
\cite[Thm. 3.3]{APS2}.

\begin{cor}\label{RhoPrep}
Let $(D_u)_{u\in \R}$ be a smooth one-parameter family of
operators in $\sP^1_{s,e}(M,E)$, and let $A$ be a flat connection
on a Hermitian vector bundle of rank $k$. Denote by
$(D_{A,u})_{u\in\R}$ the one-parameter family of operators
obtained by twisting with $A$.
\begin{enumerate}
\item If the kernels of $D_u$ and $D_{A,u}$ are of constant
dimensions, then the Rho invariant $\rho(D_{A,u})$ is independent
of $u$.
\item In the general case, only the reduced Rho invariant,
\begin{equation*}
\big[\rho(D_{A,u})\big]:=\big[\xi(D_{A,u})\big] - k\cdot
\big[\xi(D_u)\big],
\end{equation*}
is independent of $u$.
\end{enumerate}
\end{cor}

\begin{proof}
As in the proof of Proposition \ref{RhoIntProp}, the operators
$D_{A,u}$ and $D_u^{\oplus k}$ are locally equivalent smooth
families in the sense of \eqref{LocalEquiv}. Since the twisting
connection is independent of $u$, the families
$\big(\frac{dD_u^{\oplus k }}{du}\big)_{u\in\R}$ and
$\big(\frac{d D_{A,u} }{du}\big)_{u\in\R}$ are locally equivalent
as well. Now, the $a_n$ in the asymptotic expansion of
\eqref{EtaVarAsymp} can be computed locally from the total
symbols of the involved operators, and so
\[
a_m\big(\lfrac{dD_u^{\oplus k }}{du},(D_u^{\oplus k})^2\big) =
a_m \big(\lfrac{dD_{A,u}}{du},(D_{A,u})^2\big).
\]
Since $\eta(D^{\oplus k}_u) = k\cdot\eta(D_u)$, the result
follows from \eqref{EtaDer} respectively \eqref{RedEtaDer}.
\end{proof}

\begin{remark*}
One can also use Proposition \ref{RhoIntProp} and proceed as in
Appendix \ref{AppEtaVar} to differentiate under the integral to
proof part (i) in Corollary \ref{RhoPrep}. This has the advantage
that it is a bit less involved than the proof of the variation
formula of the individual Eta invariant, since one does not have
to deal with meromorphic continuations. However, the main steps
remain the same.
\end{remark*}

\noindent\textbf{Spectral Flow.} The smoothness of $[\xi(D_u)]$
shows that the discontinuities of $\xi(D_u)$ as $u$ varies are
only integer jumps. Heuristically, this is due to the fact that
the Eta invariant is a regularized signature so that whenever an
eigenvalue of $D_u$ crosses $0$, it changes by an integer multiple
of 2. This can be made more precise using the notion of
\emph{spectral flow}, which we now introduce briefly.

For a smooth one-parameter family of formally self-adjoint
elliptic operators $(D_u)_{u\in[a,b]}$, it can be shown that the
associated family of compact resolvents varies smoothly with $u$
in the operator norm on $L^2(M,E)$, see the proof of Theorem
\ref{DuhamelThm} for some related ideas. This implies that the
eigenvalues of $D_u$ can be arranged in such a way that they vary
continuously with $u$. In particular, we can find a partition
$a=u_0<u_1<\ldots<u_n=b$ such that for each $i\in\{1,\ldots,n\}$
there is $c_i>0$ with
\begin{equation}\label{SFDefHelp}
c_i\notin \bigcup_{u\in[u_{i-1},u_i]}\spec(D_u).
\end{equation}
For $u\in [u_{i-1},u_i]$ denote by $P_{[0,c_i]}(u)$ the
finite-rank spectral projection associated to eigenvalues in the
interval $[0,c_i]$. We then define the spectral flow in the
spirit of \cite{Phi}.

\begin{dfn}\label{SFDef}
Let $(D_u)_{u\in[a,b]}$ be a smooth one-parameter family in
$\sP_{s,e}^d(M,E)$, and let $a=u_0<u_1<\ldots<u_n=b$ be a
partition such that there exist $c_i$ as in \eqref{SFDefHelp}.
Then the \emph{spectral flow} between $D_a$ and $D_b$ is defined
as
\[
\SF(D_u)_{u\in[a,b]} := \sum_{i=1}^n \rk\big(
P_{[0,c_i]}(u_i)\big) -  \rk\big(P_{[0,c_i]}(u_{i-1})\big).
\]
\end{dfn}

%
%

\begin{figure}[htbp]
\centering
\includegraphics[width=0.6\linewidth]{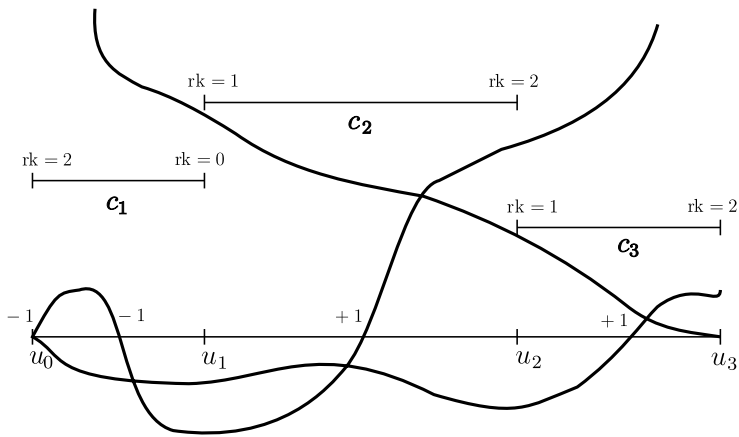}
\caption{Spectral Flow}\label{Fig:SF}
\end{figure}

Without going into further details, we note that the spectral
flow is well-defined and independent of the choices made.
Moreover, as indicated in Figure \ref{Fig:SF}, there is a built-in
convention of how to count zero eigenvalues at the end points.
For more details and generalizations we refer to \cite{BLP, L05,
Phi}. In the situation at hand, there is also the approach as in
\cite{RobSal} using Kato's selection theorem.\\

\noindent\textbf{The Variation Formula.} With the notion of
spectral flow at hand, we can now state a result on the variation
of the Eta invariant, due to Atiyah, Patodi and Singer in
\cite{APS3}. We sketch a proof in Corollary \ref{EtaDiffSFApp},
see also \cite[Lem. 3.4]{KL04}.

\begin{prop}\label{EtaDiffSF}
Let $(D_u)_{u\in [a,b]}$ be a smooth one-parameter family of
operators in $\sP_{e,s}^1(M,E)$, where $M$ is closed. Then
\begin{equation*}
\xi(D_b) - \xi (D_a) = \SF (D_u)_{u\in[a,b]} + \int_a^b \lfrac
d{du}[\xi(D_u)]du.
\end{equation*}
\end{prop}

\subsection{The Atiyah-Patodi-Singer Index Theorem}

To relate the Eta invariant to the index theorem, we have to
leave the realm of closed manifolds briefly.\\

\begin{figure}[htbp]
\centering
\includegraphics[width=0.5\linewidth]{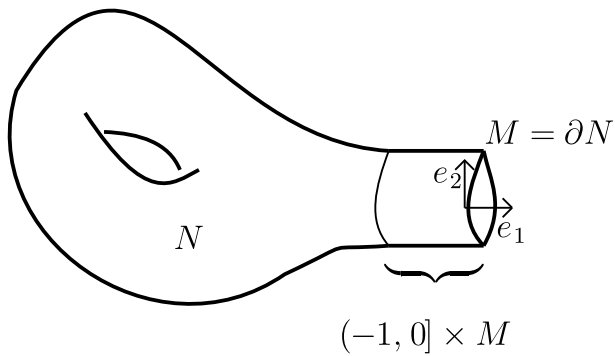}
\caption{Collar and orientation convention}\label{Fig:Collar}
\end{figure}

\noindent\textbf{Product Structures at the Boundary.} Let $N$ be
compact manifold with boundary $\pd N=M$. We equip $N$ with a
metric $g_N$ of \emph{product type} near the boundary, i.e., we
assume that a collar of the boundary is isometric to
$(-1,0]\times M$ endowed with the metric $g=dt^2+ g_M$, where
$g_M$ is a metric on $M$. If $N$ is oriented, we get an induced
orientation on $(-1,0]\times M$ which we use to orient $M$. This
is the \emph{outward normal first} convention, see Figure
\ref{Fig:Collar}. Atiyah, Patodi and Singer in \cite{APS1} use a
different convention, which explains some sign differences.

\begin{dfn}\label{ProductOpDef}
A $\Z_2$-graded, formally self-adjoint elliptic differential
operator of first order $D:C^\infty(N,E)\to C^\infty(N,E)$ is
called \emph{in product form} if the following holds.
\begin{enumerate}
\item On a collar of $M$
\[
E^+|_{(-1,0]\times M} = \pi^* E_M,
\]
where $\pi:(-1,0]\times M\to M$ is the projection, and $E_M$ is a
Hermitian vector bundle over $M$.
\item If $\gamma$ is the bundle isomorphism $E^+|_M\to E^-|_M$,
given by applying the symbol of $D^+$ to the outward normal unit
vector, then
\begin{equation*}
D^+ = \gamma\big(\lfrac d{dt} - D_M\big),\quad\text{over}\quad
(-1,0]\times M,
\end{equation*}
where $D_M$ is a formally self-adjoint elliptic differential
operator of first order on $E_M$, called the \emph{tangential
operator}.
\end{enumerate}
\end{dfn}

An operator of the above form can be extended canonically to an
operator on the \emph{closed double}
\[
X = N\cup_M -N.
\]
Here, $-N$ denotes $N$ with the reversed orientation. One uses
$\gamma$ to glue $E|_{\pd N}$ to $E|_{\pd (-N)}$ to get a bundle
$E_X\to X$, and the operator $D$ extends naturally over $X$. For
more details and a detailed proof of the following we refer to
\cite[Thm. 9.1]{BW}.

\begin{prop}
If $D$ is in product form, then the natural extension
\[
D_X:C^\infty(X,E_X)\to C^\infty(X,E_X)
\]
is invertible and extends $D$ in the sense that
$D_X|_{C^\infty(N,E)} = D$.
\end{prop}

\noindent\textbf{Boundary Conditions and the Index Theorem.} To
formulate an index theorem for $D$ one needs to introduce
suitable boundary conditions in order to render $D$ Fredholm. As
observed by Atiyah and Bott \cite{AB64}, most geometric Dirac
operators do not admit local boundary conditions that set up a
well-posed boundary value problem. Atiyah, Patodi and Singer
\cite{APS1} solved this by introducing the following global
boundary projection.

\begin{dfn}\label{APSProj}
Let $D\in \sP_{s,e}^1(N,E)$ be $\Z_2$-graded and in product form,
with tangential operator $D_M$. The \emph{Atiyah-Patodi-Singer
projection}
\begin{equation*}
P_{\ge}(D_M): C^\infty(M,E_M)\to C^\infty(M,E_M)
\end{equation*}
is defined as the spectral projection onto the subspace spanned by
the eigenvectors of $D_M$ corresponding to eigenvalues $\ge 0$.
\end{dfn}

Then the index theorem for manifolds with boundary in \cite[Thm.
3.10]{APS1} reads

\begin{theorem}[Atiyah-Patodi-Singer]\label{APSIndThm}
Let $N$ be compact, oriented Riemannian manifold of dimension $n$
with boundary $\pd N=M$, and let
\[
D:C^\infty(N,E)\to C^\infty(N,E)
\]
be a $\Z_2$-graded, formally self-adjoint elliptic differential
operator of first order. Assume that the metric and $D$ are in
product form in a collar of $M$. Let
\[
C^\infty\big(N,E^+;P_\ge\big):= \bigsetdef{\gf\in
C^\infty(N,E^+)}{P_{\ge}(D_M)(\gf|_M)=0},
\]
where $P_{\ge}(D_M)$ is the projection of Definition
\ref{APSProj}. Then
\[
D^+: C^\infty\big(N,E^+;P_\ge\big) \to C^\infty\big(N,E^-\big)
\]
has a natural Fredholm extension with
\begin{equation}\label{GenAPS}
\ind\big(D^+;P_\ge\big) = \int_N \str_E[\widetilde e_n]\vol_N
-\xi(D_M).
\end{equation}
Here, $\xi(D_M)$ is the $\xi$-invariant as in Definition
\ref{EtaDef}, and $\widetilde e_n(x)$ is the coefficient of the
constant term in the asymptotic expansion as $t\to 0$ of the heat
kernel $e^{-t D^2_X}(x,x)$ associated to the extension $D_X$ to
the closed double $X$.
\end{theorem}

\begin{remark*}
We are following \cite{APS1} and use the operator $D_X$ on the
closed double to define the index density. Clearly, this is an ad
hoc method which allows to use the asymptotic expansion of the
heat trace in its version for closed manifolds as in Theorem
\ref{HeatKernel}. However, the trace expansion can be formulated
in a much more abstract functional analytic setting to
incorporate the case of manifolds with boundary, see \cite{BL98}
for an expository account.
\end{remark*}

If $\dim N$ is even and $D$ is a $\Z_2$-graded geometric Dirac
operator, the index density can be made explicit using the local
index theorem for closed manifolds. One obtains the index theorem
for geometric Dirac operators on manifolds with boundary, see
\cite[Thm. 4.2]{APS1}.

\begin{theorem}[Atiyah-Patodi-Singer]\label{APSIndThmDirac}
Let $N$ be compact, oriented Riemannian manifold of even dimension
$n$ with boundary $\pd N=M$, and let $D:C^\infty(N,E)\to
C^\infty(N,E)$ be a $\Z_2$-graded, geometric Dirac operator which
is in product form near $M$. Then the index of
\[
D^+: C^\infty\big(N,E^+;P_\ge\big) \to C^\infty\big(N,E^-\big)
\]
is given by
\begin{equation}\label{GeomAPS}
\ind\big(D^+;P_\ge\big) = \int_N \widehat A(TN,\nabla^g)\wedge
\ch_{E/S}(E,\nabla^E) -\xi(D_M),
\end{equation}
where $\widehat A(TN,\nabla^g)$ is the $\widehat A$-form of
Definition \ref{PontLAHatDef} with respect to the Levi-Civita
connection $\nabla^g$ and $\ch_{E/S}$ is the realtive Chern
character \eqref{RelChern}.
\end{theorem}

\section{The Atiyah-Patodi-Singer Rho Invariant}

In Proposition \ref{Sign=Ind} we have seen that the signature of
closed manifolds equals the index of the signature operator. One
motivation that led Atiyah, Patodi and Singer to the discovery of
Theorem \ref{APSIndThm} was the search for a generalization of
this to the case for manifolds with boundary. We thus briefly
recall the definition of the signature of a manifold with
boundary.

\subsection{The Signature of Manifolds with Boundary}\label{SignBound}

Let $N$ be a compact, oriented and connected manifold of
dimension $n$ with boundary $\pd N$. Then $\pd N$ is closed and
naturally oriented, but we allow that it consists of several
connected components. Let $\ga:\pi_1(N)\to \U(k)$ be a
representation of the fundamental group. Via the map induced by
the inclusion $\pd N\hookrightarrow N$, the representation $\ga$
restricts to $\pi_1(\pd N)$. There is a relative version
$H^\bullet(N,\pd N ,E_\ga)$ of cohomology with local coefficients,
whose construction we will not describe in detail. We only note
that the machinery of algebraic topology---like cup and cap
products, Poincar\'e duality, and exact sequence of
pairs---extends to this context\footnote{However, one has to be
careful with the Mayer-Vietoris sequence and excision, where one
has to use van Kampen's Theorem to relate the representations of
the involved fundamental groups.}. In particular, there is a
relative intersection pairing
\[
H^p(N,E_\ga)\times H^{n-p}(N,\pd N,E_\ga)\to \C,\quad
(a,b)\mapsto \Scalar{a\cup b}{[N,\pd N]},
\]
where the fundamental class $[N,\pd N]$ is the generator of
$H_n(N,\pd N)$ determined by the orientation. Via the de Rham
isomorphism, the relative cohomology groups are isomorphic to de
Rham cohomology with compact support in the interior of $N$,
\[
H^p(N,\pd N,E_\ga)\cong H^p_c(N,E_\ga).
\]
Here, we also use $E_\ga$ to denote the flat vector bundle
determined by $\ga$. With respect to this identification, the
intersection pairing is induced by
\[
\gO^p(N,E_\ga)\times \gO^{n-p}_c(N,E_\ga)\mapsto \C,\quad
(\go,\eta)\mapsto \int_N \langle\go\wedge\eta\rangle.
\]
Now assume that $n$ is even. There is a natural map
$H^{n/2}(N,\pd N,E_\ga)\to H^{n/2}(N,E_\ga)$ which we can combine
with the intersection pairing to define a \emph{twisted
intersection form}
\[
Q_\ga: H^{n/2}(N,\pd N,E_\ga)\times H^{n/2}(N,\pd N,E_\ga) \to \C.
\]
As for closed manifolds, the intersection form is skew for
$(n\equiv 2 \mod 4)$ and symmetric for $(n\equiv 0\mod 4)$, but it
is in general degenerate.

\begin{dfn}
Let $N$ be a compact, connected manifold with boundary of even
dimension $n$, and let $\ga:\pi_1(N)\to\U(k)$ be a representation
of the fundamental group . Then the \emph{twisted signature} of
$M$ is defined as
\[
\Sign_\ga(N):=\Sign(Q_\ga),
\]
where we use the earlier convention regarding the signature of a
skew endomorphism.
\end{dfn}

\begin{remark*}
To turn $Q_\ga$ into a non-degenerate form, one needs to restrict
to
\[
\widehat H^{n/2}(N,E_\ga):=\im\big( H^{n/2}(N,\pd N,E_\ga)\to
H^{n/2}(N,E_\ga)\big).
\]
Then Poincar\'e duality ensures that we get a non-degenerate form
\[
\widehat Q_\ga: \widehat H^{n/2}(N,E_\ga)\times \widehat
H^{n/2}(N,E_\ga) \to \C.
\]
Since this process just eliminates the radical of $Q_\ga$, it is
immediate that
\[
\Sign Q_\ga=\Sign \widehat Q_\ga.
\]
\end{remark*}

\subsection{Twisted Odd Signature Operators}

To state the relation of the signature on a manifold with
boundary to the index of the signature operator, we first need to
understand the structure of the signature operator near the
boundary. Therefore, we now consider the model case of a
cylinder, and derive the formula for the odd signature operator.

Let $M$ be a closed, oriented manifold of odd dimension $m$ and
consider the cylinder $N:= \R\times M$. We use the natural
splitting $T^*N = \R\oplus T^*M$ to orient $N$ and make the
identification
\[
\gF:C^\infty\big(\R,\gO^\bullet(M)\big)\oplus
\C^\infty\big(\R,\gO^\bullet(M)\big) \xrightarrow{\cong}
\gO^\bullet(N),\quad \gF(\go_0,\go_1):= dt\wedge \go_0 + \go_1,
\]
where $t$ denotes the $\R$ coordinate. Clearly,
\begin{equation}\label{dCyl}
d_N\gF\big(\go_0,\go_1) = dt\wedge\big(\pd_t\go_1 - d_M\go_0\big)
+ d_M\go_1 = \gF\big(\pd_t\go_1 - d_M\go_0,d_M\go_1\big).
\end{equation}
We now endow $N$ with a metric of product form $g = dt^2 + g_M$,
and denote by $\tau_N$ and $\tau_M$ the chirality operators on
$\gO^\bullet(N)$ and $\gO^\bullet(M)$, respectively. Then one
checks that
\[
\tau_N \gF\big(\go_0,\go_1) =
\gF\big(\tau_M\go_1,\tau_M\go_0\big).
\]
From this we obtain isomorphisms
\[
\gF_\pm: C^\infty\big(\R,\gO^\bullet(M)\big)\xrightarrow{\cong}
\gO^\pm(N),\quad \go\mapsto \gF(\go,\pm\tau_M\go).
\]
Let $D^+_N:\gO^+(N)\to \gO^-(N)$ be the signature operator. Then a
short computation using \eqref{dCyl} and the formula
\eqref{d^tDef} for the adjoint differential yields
\[
\gF_-^{-1}\circ D^+_N\circ \gF_+ = \tau_M\big(\pd_t -\tau_Md_M
-d_M\tau_M\big).
\]
The same continues to hold if $E\to M$ is a Hermitian vector
bundle, and when we twist with a unitary connection $A$ on
$\pi^*E$ in \emph{temporal gauge}, i.e., a connection of the form
$A:=\pi^*a$. Here, $\pi:N\to M$ is the natural projection and $a$
is a unitary connection on $E$. We summarize what we have
observed so far.

\begin{prop}\label{SignOpCylinder}
Let $D_A^+$ be the signature operator on the cylinder $N=\R\times
M$ twisted by a unitary connection $A=\pi^* a$ in temporal gauge.
Then $D^+_A$ is isometric to
\begin{equation*}
\tau_M \big(\pd_t -
B_a\big):C^\infty\big(\R,\gO^\bullet(M,E)\big)\to
C^\infty\big(\R,\gO^\bullet(M,E)\big),
\end{equation*}
where
\begin{equation}\label{OddSignTot}
B_a := \tau_M\big(d_a+d^t_a) = \tau_M d_a + d_a\tau_M.
\end{equation}
\end{prop}

\begin{dfn}\label{OddSignDef}
Let $a$ be a unitary connection over an oriented Riemannian
manifold $M$ of odd dimension. The operator
\[
B_a^{\ev}:=B_a|_{\gO^{\ev}(M,E)}:\gO^{\ev}(M,E)\to \gO^{\ev}(M,E)
\]
is called the \emph{odd signature operator} on $M$ twisted by $a$.
\end{dfn}

\begin{remark}\label{OddSignRem}\quad\nopagebreak
\begin{enumerate}
\item Note that the operator $B_a$ does indeed preserve the
even/odd grading,
\[
B_a=B_a^{\ev}\oplus B_a^{\odd} : \gO^{\ev}(M,E)\oplus
\gO^{\odd}(M,E) \to \gO^{\ev}(M,E)\oplus \gO^{\odd}(M,E).
\]
Moreover, $B_a^{\ev}$ and $B_a^{\odd}$ are conjugate via
$\tau_M$. The kernels of $B_a$ and $B_a^{\ev}$ are of a
topological nature, since
\[
\ker(B_a) = \ker(d_a+d_a^t) = \sH^\bullet(M,E_a),\quad
\ker(B_a^{\ev}) = \sH^{\ev}(M,E_a).
\]
\item It is straightforward to check that the odd signature
operator is a geometric Dirac operator in the sense of Definition
\ref{GeomDiracOp}. Here, one defines the Clifford structure by
\begin{equation}\label{OddCliffordMult}
c^{\ev}:T^*M\otimes \gL^{\ev} T^*M\to \gL^{\ev} T^*M,\quad
c^{\ev}(\xi)\go := \tau_M\big(\xi\wedge \go-\imu(\xi)\go\big).
\end{equation}
Paying close attention to the various identifications made, one
also verifies that the map
\[
\tau_M: \gL^\bullet T^*M\to \gL^\bullet T^*M
\]
corresponds to Clifford multiplication
\[
c_N(dt): \gL^+T^*N|_M\to  \gL^-T^*N|_M.
\]
Hence, Proposition \ref{SignOpCylinder} shows that $D_A^+$ is of
product form in the sense of Definition \ref{ProductOpDef}.
\item To give an explicit formula, let $m=2k-1$. Then one can
check, using the formula \eqref{TauExpl} for $\tau_M$, that for
all $\go\in\gO^p$
\[
B_a\go = -i^{k+p(p+1)}\big((-1)^{p}*d_a-d_a*\big)\go.
\]
In particular for $p=2q$,
\begin{equation*}
B^{\ev}_a\go= i^k(-1)^{q+1}(*d_a-d_a*)\go.
\end{equation*}
\end{enumerate}
\end{remark}

\subsection{The Signature Theorem for Manifolds with Boundary}

Having identified the tangential operator, the APS projection
sets up a well-defined index problem for the signature operator
on manifolds with boundary. Moreover, the Atiyah-Patodi-Singer
Index Theorem for geometric Dirac operators \ref{APSIndThmDirac}
and the index theorem for the twisted signature operator in
Theorem \ref{TwistSignIndThm} calculates its index. If the
twisting bundle is flat, this index is related to the signature
for manifolds with boundary as introduced in Section
\ref{SignBound}. However, this relation is more difficult to
establish than in Proposition \ref{Sign=Ind} for the case that
the manifold is closed. We only state the result and refer to
\cite[Sec. 4]{APS1} for the proof. A concise discussion also be
found in \cite[Sec. 23]{BW}.

\begin{theorem}[Atiyah-Patodi-Singer]\label{SignIndThmBound}
Let $N$ be a compact, oriented Riemannian manifold of even
dimension $n$ with boundary $\pd N=M$, and let $E\to N$ be a
Hermitian vector bundle of rank $k$ with a unitary connection $A$.
Assume that the metric is in product form, and that $A$ is in
temporal gauge $A=\pi^*a$ on a collar of $M$. Let $D_A^+$ be the
twisted signature operator on $N$, and let
\[
P_{\ge}(a):\gO^\bullet(M,E|_M)\to  \gO^\bullet(M,E|_M)
\]
be the APS projection in Definition \ref{APSProj} of the
tangential operator $B_a$. Then
\[
\ind\big(D^+_A;P_\ge(a)\big) = 2^{n/2} \int_N
\widehat L(TN,\nabla^g)\wedge \ch(E,A) - \xi(B_a).
\]
Moreover, if $A$ is flat, then
\[
\ind\big(D^+_A;P_\ge(a)\big) = \Sign_A(N) - \lfrac 12 \dim(\ker
B_a),
\]
and therefore,
\[
\Sign_A(N) = k\cdot  \int_N L(TN,\nabla^g) - \eta(B_a^{\ev}).
\]
\end{theorem}

\begin{remark*}
In the last equation, the occurrence of $B_a^{\ev}$ stems from the
relation
\[
\xi(B_a) = \eta(B_a^{\ev}) +\lfrac 12 \dim(\ker B_a),
\]
which follows from Remark \ref{OddSignRem} (i).
\end{remark*}

\noindent\textbf{Rho Invariants.} Motivated by Theorem
\ref{SignIndThmBound}, Atiyah, Patodi and Singer \cite{APS2}
introduced the Rho invariant, which we have already briefly
considered in Proposition \ref{RhoIntProp} and Corollary
\ref{RhoPrep}. We now treat the Rho invariant associated to the
odd signature operator in more detail.

\begin{dfn}\label{RhoDef}
Let $M$ be a closed, oriented Riemannian manifold of odd dimension
$m$, and let $A$ be a flat unitary connection on a Hermitian
bundle $E$ of rank $k$. Then the \emph{Rho invariant of $A$} is
defined as
\begin{equation*}
\rho_A(M):= \rho(B^{\ev}_A) =\eta(B^{\ev}_A)-k\cdot \eta(B^{\ev}).
\end{equation*}
\end{dfn}

We have the following immediate consequences of what we have
discussed so far.

\begin{prop}\label{RhoProp} \quad\nopagebreak
\begin{enumerate}
\item If $A'$ is a flat connection, unitarily equivalent to $A$,
then $\rho_{A'}(M) = \rho_A(M)$. In particular, the Rho invariant
depends only on the holonomy representation
\[
\hol_A : \pi_1(M)\to \U(k).
\]
For this reason, we also use the notation $\rho_\ga(M)$ if the
focus is on representations of the fundamental group.
\item The Rho invariant is independent of the metric used to
define $\eta(B^{\ev}_A)$ and $\eta(B^{\ev})$. Therefore, it is a
smooth invariant of $M$ and $A$.
\item If $N$ is a compact, oriented manifold with boundary $\pd
N=M$, and the representation $\ga:\pi_1(M)\to \U(k)$ extends to a
unitary representation $\gb:\pi_1(N)\to\U(k)$, then
\[
\rho_\ga(M)= \Sign_{\gb}(N) - k\cdot\Sign(N).
\]
\end{enumerate}
\end{prop}

\begin{proof}
Part (i) follows from the fact that the Eta invariant does not
change, if we transform with a unitary bundle isomorphism. Part
(ii) is a consequence of Corollary \ref{RhoPrep} (i), since the
dimensions of the kernels of $B^{\ev}_\ga$ and $B^{\ev}$ are
independent of the metric. Part (iii) follows immediately from
the signature formula of Theorem \ref{SignIndThmBound}.
\end{proof}

\begin{remark}\label{RhoRem}\quad\nopagebreak
\begin{enumerate}
\item Part (i) of Proposition \ref{RhoProp} shows that the Rho
invariant can be interpreted as a map, defined on the moduli space
of flat connections
\[
\rho: \cM\big(M,\U(k)\big)\to \R,\quad [A]\mapsto \rho_A(M).
\]
In the explicit examples in Section \ref{S1Bundles} and Chapter
\ref{3dimMapTor} we will use this and consider particularly
well-suited representatives for gauge equivalence classes of flat
connections.
\item Proposition \ref{RhoProp} (iii) gives a negative answer to the
question of whether Corollary \ref{SignMultClosed} continues to
hold for the signature of manifolds with boundary. As mentioned
in the introduction, an explanation of this \emph{signature
defect} was one of the motivations leading to the discovery of
Theorem \ref{SignIndThmBound}, and the Rho invariant is indeed an
intrinsic characterization of this.
\item Rho invariants do in general give non-trivial invariants.
As an easy example we consider $M=S^1$. We view $S^1$ as a subset
of $\C$, endowed with the metric of length $2\pi$. A flat
$\U(1)$-connection on the trivial line bundle is determined by
its holonomy $e^{2\pi ia}\in \U(1)$ with $a\in \R$. The
corresponding odd signature operator is easily seen to be
\[
B_a^{\ev}= -i(\sL_e - ia): C^\infty(S^1)\to C^\infty(S^1),
\]
where with $\gf\in C^\infty(S^1)$ and $z\in S^1$,
\[
\sL_e\gf(z) = \lfrac d{dt}\big|_{t=0} \gf(z e^{it}).
\]
Therefore, $B_a^{\ev}\gf = \gl \gf$ if and only if
\[
\gl+ a\in \Z\quad\text{and}\quad \gf(z)= z^{\gl+a}\gf(1).
\]
This implies
\[
\eta(B_a^{\ev},s) = \sum_{\begin{smallmatrix} n \in\Z \\
n\neq a
\end{smallmatrix}}\frac{\sgn(n-a)}{|n -a|^s}, \quad\Re(s)>1.
\]
In Proposition \ref{ZEtaCalc} (i), we have included a computation
of the value of the meromorphic continuation of this expression
at 0. The result is, see also Definition \ref{PeriodicBernoulli},
\begin{equation*}
\eta(B_a^{\ev}) = 2 P_1(a) = \begin{cases}\quad \,\,0, &\text{for
$a\equiv 0\mod\Z$},\\ 2a_0-1,&\text{for $a_0\in (0,1)$ and
$a\equiv a_0\mod\Z$.}
\end{cases}
\end{equation*}

\item We want to point out that the Rho invariant is a true
extension of the signature defect. More explicitly, there are Rho
invariants which cannot be calculated using the formula of
Proposition \ref{RhoProp} (iii). As an example, consider a
compact oriented surface $\gS$ with one boundary component $\pd
\gS = S^1$. Then the fundamental class of $S^1$ is a commutator
in $\pi_1(\gS)$, see for example the discussion in Section
\ref{S1BundlesModSpace}. This implies that a non-trivial
$\U(1)$-representations of $\pi_1(S^1)$ cannot extend to a
representation of $\pi_1(\gS)$.
\end{enumerate}
\end{remark}

\section{Rho Invariants and Local Index Theory}\label{RhoLocalInd}

\subsection{Relation to Chern-Simons Invariants}

As seen in Proposition \ref{RhoProp}, the Rho invariant can be
computed in a purely topological way if the representation $\ga$
extends over a bounding manifold. However, we have already
pointed out that this situation is often too restrictive.
Therefore, intrinsic methods to compute Rho invariants are of
great interest. One observation for applying topological tools is
the relation of Rho invariants to Chern-Simons invariants. We
refer to Appendix \ref{CharClass} for definitions and basic
properties of Chern-Simons invariants, and make the relation more
precise now. The following result goes back to \cite[Sec.
4]{APS2}, see also \cite[Sec. 7]{KKR}.

\begin{prop}\label{EtaConnVar}
Assume that $M$ is a closed, oriented Riemannian manifold of odd
dimension $m$. Let $A_t$ be a smooth path of connections on a
fixed Hermitian vector bundle $E\to M$. Then the reduced
$\xi$-invariant satisfies
\[
\int_0^1 \lfrac d{dt}[\xi(B_{A_t})]dt =
2^{\frac{m+1}2}\int_M\widehat L(TM,\nabla^g)\wedge \cs(A_0,A_1),
\]
where $\cs(A_0,A_1)$ is the transgression form of the Chern
character, see Definition \ref{TransgressDef}.
\end{prop}

We postpone the proof and mention some consequences.

\begin{cor}\label{RhoConnVar}
Let $A_t$ be a smooth path of connections on $E$.
\begin{enumerate}
\item The following variation formula holds
\[
\xi(B_{A_1})- \xi(B_{A_0})= \SF(B_{A_t})_{t\in[0,1]} +
2^{\frac{m+1}2}\int_M\widehat L(TM,\nabla^g)\wedge \cs(A_0,A_1).
\]
\item If $A_0$ and $A_1$ are flat, and $M$ is 3-dimensional, then
\[
\rho_{A_1}(M)  = \rho_{A_0}(M)+ 4 \CS(A_0,A_1) \mod \Z.
\]
Here, $\CS(A_0,A_1)$ is the Chern-Simons invariant associated to
the Chern character, see Definition \ref{CSDef} and
\eqref{ChernSimonsExpl}.
\item Assume that $A_t$ is a path of flat connections
and that either $A_0$ and $A_1$ reduce to $\SU(k)$-connections or
$(\dim M \equiv 3\mod 4)$. Then
\[
\eta(B^{\ev}_{A_1})- \eta(B^{\ev}_{A_0}) =
2\SF(B^{\ev}_{A_t})_{t\in [0,1]} - \dim\ker B^{\ev}_{A_1} +
\dim\ker B^{\ev}_{A_0}.
\]
In particular,
\[
\rho_{A_1}(M)= \rho_{A_0}(M)\mod \Z.
\]
\end{enumerate}
\end{cor}

\begin{remark}\label{RhoConnVarRem}\quad\nopagebreak
\begin{enumerate}
\item Corollary \ref{RhoConnVar} (ii) shows that on a 3-manifold
the reduction mod $\Z$ of the Rho invariant is up to a constant
the Chern-Simons invariant of the corresponding flat connection.
Therefore, the reduced Rho invariant is the integral over local
invariants of the connections. Now, the unreduced version is
essentially this ``local'' contribution plus a spectral flow term
which encodes ``global'' topological information.

\item Under the hypothesis of part (iii) the reduction mod $\Z$
of $\rho_{A_t}(M)$ is constant. In other words, the Chern-Simons
invariant is constant on connected components of the moduli space
of flat connections. Therefore, unreduced Rho invariants
associated to one-parameter families of flat connections have only
integer jumps, which occur precisely at the points where the rank
of the twisted cohomology groups changes. Independently,
Farber-Levine \cite{FL} and Kirk-Klassen \cite{KK94, KK97} have
developed powerful methods to compute this spectral flow term in
purely cohomological terms.
\end{enumerate}
\end{remark}

We will now give a proof of Proposition \ref{EtaConnVar} based on
the signature theorem for manifolds with boundary in Theorem
\ref{SignIndThmBound}. We shall also sketch a different proof
using Getzler's approach to local index theory in Section
\ref{VarLocalInd} below.

\begin{figure}[t]
\centering
\includegraphics[width=0.45\linewidth]{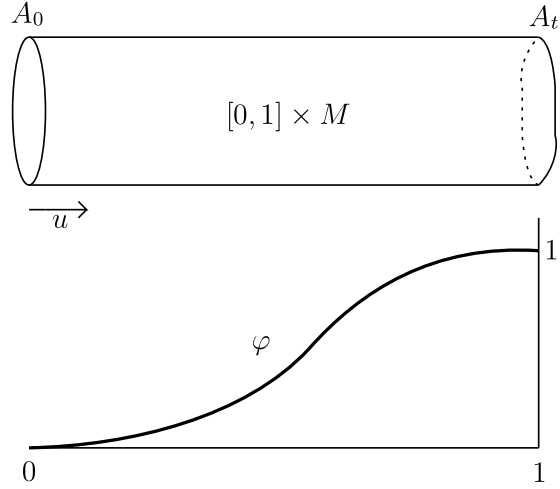}
\caption{Cylinder $N$ and the cutoff function
$\gf$}\label{Fig:CSRel}
\end{figure}

\begin{proof}[Proof of Proposition \ref{EtaConnVar}]
Consider the cylinder $N:=[0,1]\times M$, endowed with the product
metric $g_N= du^2+ \pi^*g$. Here, we are using $u$ to denote the
coordinate on $[0,1]$, and $\pi:N\to M$ is the natural projection.
We endow $\pi^*E$ with the connection
\[
\widetilde A_t:= du\wedge \lfrac d{du} + \pi^* A_{t\gf(u)},
\]
where $\gf:[0,1]\to[0,1]$ is a smooth function such that for some
$\eps>0$,
\begin{equation}\label{ProdStrucCutoff}
\gf(u) = \begin{cases} 0, &\text{if $u<\eps$},\\ 1, &\text{if
$u>1-\eps$},
\end{cases}
\end{equation}
see Figure \ref{Fig:CSRel}. For fixed $t$, the connection
$\widetilde A_t$ is in temporal gauge on a collar of $\pd N$.
Therefore, we can apply Theorem \ref{SignIndThmBound} for each
$t$ to the signature operator $D_{\widetilde A_t}^+$ and conclude
that
\[
\xi(B_{A_t})- \xi(B_{A_0}) = 2^{\frac{m+1}2} \int_N \widehat
L(TN,\nabla^{g_N})\wedge \ch(\pi^*E,\widetilde A_t)\mod\Z.
\]
Since $g_N$ is the product metric on $N$ it is straightforward to
check that
\[
\widehat L(TN,\nabla^{g_N})=\pi^*\widehat L(TM,\nabla^g).
\]
Then, if $\int_{N/M}$ denotes integration along the fiber, see
Proposition \ref{FiberInt} below, we observe that
\[
\begin{split}
\int_N \widehat L(TN,\nabla^{g_N})\wedge \ch(\pi^*E,\widetilde
A_t) &= \int_M \widehat L(TM,\nabla^g)\wedge \int_{N/M}
\ch(\pi^*E,\widetilde A_t)\\
&= \int_M \widehat L(TM,\nabla^g)\wedge \cs(A_{t\gf(u)}).
\end{split}
\]
Here, we have used Lemma \ref{TransgressCyl} in the second line to
replace $\int_{N/M} \ch(\pi^*E,\widetilde A_t)$ with the
transgression form of the Chern character computed with respect
to the path $u\mapsto A_{t\gf(u)}$. If we choose a different path,
the result will differ by an exact form on $M$, see Proposition
\ref{CSPathIndep} for a proof. This allows us to remove the
function $\gf$ and use the path $u\mapsto A_{tu}$. Then
\[
\cs(A_{t\gf(u)}) = \cs(A_0,A_t)\mod d\gO^{\ev}(M).
\]
Since $M$ is closed, this shows that
\[
\xi(B_{A_t})- \xi(B_{A_0}) = 2^{\frac{m+1}2} \int_M \widehat
L(TM,\nabla^g)\wedge\cs(A_0,A_t) \mod\Z,
\]
which implies Proposition \ref{EtaConnVar}.
\end{proof}

\begin{proof}[Proof of Corollary \ref{RhoConnVar}]
Part (i) is an immediate consequence of Proposition
\ref{EtaConnVar} and Proposition \ref{EtaDiffSF}. Next, we note
that
\begin{equation}\label{EtaXiRel}
\xi(B_{A_1})- \xi(B_{A_0}) = \eta(B_{A_1}^{\ev})-
\eta(B_{A_0}^{\ev}) +\dim\ker(B_{A_1}^{\ev}) -
\dim\ker(B_{A_0}^{\ev}).
\end{equation}
Using part (i) and reducing mod $\Z$ one finds that
\[
\eta(B_{A_1}^{\ev})- \eta(B_{A_0}^{\ev}) =
2^{\frac{m+1}2}\int_M\widehat L(TM,\nabla^g)\wedge
\cs(A_0,A_1)\mod\Z.
\]
Now, if we assume that $m=3$, the $\widehat L$-form equals 1, and
we obtain part (ii). Concerning part (iii), we now consider the
connection $\widetilde A:=dt\wedge\lfrac d{dt}+A_t$ on the
cylinder $[0,1]\times M$. Since we are assuming that $A_t$ is a
path of \emph{flat} connections, we have
\[
F_{\widetilde A}= dt\wedge \lfrac d{dt} A_t,\quad\text{and}\quad
\exp\big(dt\wedge \lfrac d{dt} A_t \big) = 1+ dt\wedge \lfrac
d{dt} A_t.
\]
This implies
\[
\begin{split}
\int_M\widehat L(TM,\nabla^g)\wedge \cs(A_0,A_1) &=
\int_{[0,1]\times M}\pi^*\widehat L(TM,\nabla^g)\wedge
\tr_E\big[\exp(\lfrac i{2\pi}F_{\widetilde A})\big] \\
&=  \int_M\widehat L(TM,\nabla^g)\wedge \lfrac i{2\pi}
\tr_E[A_1-A_0].
\end{split}
\]
If $A_0$ and $A_1$ are $\SU(k)$-connections, then
$\tr_E[A_1-A_0]=0$, so that the integrand vanishes. On the other
hand, if $(m\equiv 3\mod 4)$, the integrand has no degree $m$
part so that the integral vanishes again. In both cases, part
(iii) follows from part (i) and \eqref{EtaXiRel}.
\end{proof}

\noindent\textbf{Variation of the Metric.} As already pointed out
in Remark \ref{RhoConnVarRem}, an important tool for studying the
Rho invariant is studying its variation under deformations of the
flat connection. However, the moduli space of flat connections is
often discrete or at least disconnected. Therefore, it is not
always possible to find a path of flat connections which joins
given endpoints. However, in many cases one can deform the
geometry of the underlying manifold in such a way that twisted
Eta invariants become computable. The main concern of this thesis
are Rho invariants of fiber bundles, and as already explained in
the introduction there are powerful methods to deform the
geometry to a much simpler situation. One underlying result is the
following, which is an analog of Proposition \ref{EtaConnVar},
see \cite[Sec. 2]{APS2}.

\begin{prop}\label{EtaMetricVar}
Assume that $g_t$ is a smooth path of Riemannian metrics on a
closed, odd dimensional manifold $M$, and denote by $\nabla^{g_t}$
the associated one-parameter family of Levi-Civita connections.
Moreover, let $A$ be a flat connection on a Hermitian vector
bundle $E\to M$ of rank $k$. Then
\[
\eta(B^{\ev}_A,g_1)- \eta(B^{\ev}_A,g_0) = k\cdot\int_M
TL(\nabla^{g_0},\nabla^{g_1}),
\]
where $TL(\nabla^{g_0},\nabla^{g_1})$ is the transgression form
of the $L$-class of $TM$, see Remark \ref{ChernWeilRem}.
\end{prop}

\begin{proof}
As before, let $N:=[0,1]\times M$ be the cylinder, endowed with
the bundle $\pi^*E$ and the connection $\widetilde A:= dt\wedge
\lfrac{d}{dt}+ \pi^*A$. In contrast to the proof of Proposition
\ref{EtaConnVar}, the connection $\widetilde A$ is flat, since
$A$ is independent of $t$. Let $\ga:\pi_1(N)\to \U(k)$ denote the
holonomy representation of $\widetilde A$. Recall that
$\Sign_\ga(N)$ is defined using the homomorphism
\[
H^\bullet(N,M,E_\ga)\to H^\bullet(N,E_\ga).
\]
However, in the case at hand, this map is trivial as $M$ is a
deformation retract of $N$, and $\ga$ is compatible with the
natural retraction. Therefore, $\Sign_\ga(N)=0$.

Now let $\gf$ be a cutoff function as in \eqref{ProdStrucCutoff},
and endow $N$ with the metric
\[
g_N:= dt^2+ \pi^*g_{\gf(t)}.
\]
Let $\nabla^{g_N}$ be the Levi-Civita connection associated to
$g_N$. Since $g_N$ is in product form near the boundary and
$\Sign_A(N)=0$, Theorem \ref{SignIndThmBound} yields
\[
\eta(B^{\ev}_A,g_1)- \eta(B^{\ev}_A,g_0) =  k\cdot  \int_N
L(TN,\nabla^{g_N}).
\]
Moreover, $\nabla^N$ is in temporal gauge on a collar of the
boundary, so that we can deduce from Proposition \ref{CSProp} and
Remark \ref{CSMultSeq} that
\[
\int_N L(TN,\nabla^{g_N}) = \int_M
TL(\nabla^{g_0},\nabla^{g_1}).\qedhere
\]
\end{proof}

\begin{remark}
In the next chapter we will also encounter a one-parameter family
of connections $(\nabla^t)_{t\in [0,1]}$ on $TM$ which is not
associated to a family of Riemannian metrics. Nevertheless, we can
study the associated family of generalized odd signature operators
\[
D_t := c^{\ev}\circ\nabla^t:\gO^{\ev}(M)\to \gO^{\ev}(M),
\]
where $c^{\ev}$ is Clifford multiplication as defined in
\eqref{OddCliffordMult}. Certainly, the operators $D_t$ will in
general not be formally self-adjoint, even if all $\nabla^t$ are
compatible with the metric. Without going into detail, we note
that a certain restriction on the torsion tensor of $\nabla^t$
guarantees that we get formally self-adjoint operators. Then we
can use the variation formula of Proposition \ref{EtaDiffSF},
\[
\xi(D_1) - \xi(D_0) = \SF(D_t)_{t\in[0,1]} +\int_0^1 \lfrac
d{dt}\big[\xi(D_t\big]dt.
\]
However, since $(D_t)_{t\in[0,1]}$ is not a family of geometric
Dirac operators, the local variation can not be identified using
Theorem \ref{SignIndThmBound}. Yet, if we are interested in Rho
invariants, it follows from Corollary \ref{RhoPrep} that for any
flat $\U(k)$-connection $A$, the local variations of $(D_t^{\oplus
k})_{t\in [0,1]}$ and $(D_{A,t})_{t\in [0,1]}$ agree. In
particular,
\begin{equation}\label{GenDiracEtaVar}
\rho(D_{A,1})= \rho(D_{A,0}) + \SF(D_{A,t})_{t\in[0,1]} - k\cdot
\SF(D_t)_{t\in[0,1]}.
\end{equation}
\end{remark}

\subsection{The Variation Formula and Local Index
Theory}\label{VarLocalInd}

In the remainder of this chapter we sketch a proof of Proposition
\ref{EtaConnVar} based on local index theory techniques. Mainly,
this is because the underlying ideas will be helpful in the
discussion of Rho invariants of fiber bundles in Chapter
\ref{ChapAbst}. Aside from that, a proof of the variation formula
in Proposition \ref{EtaConnVar}, which does not rely on the index
theorem for manifolds with boundary, underlines the intrinsic
nature of Rho invariants.\\

\noindent\textbf{The setup.} Assume that $M$ is a closed,
oriented Riemannian manifold of odd dimension $m$, and let $A_u$
be a smooth path of connections on a fixed Hermitian vector
bundle $E\to M$. We want to have an explicit formula for the
variation $\lfrac d{du}[\xi(B_{A_u})]$ of the reduced
$\xi$-invariant, where
\[
B_{A_u}=\tau \big(d_{A_u}+d_{A_u}^t\big):\gO^\bullet(M,E)\to
\gO^\bullet(M,E).
\]
Proposition \ref{RedEtaDerApp} shows that
\[
\lfrac d{du}[\xi(B_{A_u})] = - \lfrac1{\sqrt\pi}\, a_m(B_{A_u}),
\]
where $a_m(B_{A_u})$ is the constant term in the asymptotic
expansion of
\[
\sqrt
t\,\Tr\big(\lfrac{dB_{A_u}}{du}e^{-tB_{A_u}^2}\big),\quad\text{as
}t\to 0.
\]
For brevity we will also use the common notation
\begin{equation}\label{LIMDef}
a_m(B_{A_u}) = \LIM_{t\to 0} \sqrt
t\,\Tr\big(\lfrac{dB_{A_u}}{du}e^{-tB_{A_u}^2}\big).
\end{equation}
In the case at hand, it is immediate that
\[
\lfrac{dB_{A_u}}{du} = \tau c\big(\lfrac{d}{du}A_u\big).
\]
Now, to get a formula for $a_m(B_{A_u})$, we can fix $u$. To keep
the notation short, we thus extract the following setup with
which we will work

\begin{dfn}\label{OddHeatKernel}
Let $A$ be a connection on $E$, and let $\dot A\in
\gO^1\big(M,\End(E)\big)$. Define
\[
D_A:= d_A+ d_A^t:\gO^\bullet(M,E)\to \gO^\bullet(M,E),
\]
and denote by
\[
k_t(x,x)\in C^\infty\big(M,\End(\gL^\bullet T^*M\otimes E)\big)
\]
the restriction to the diagonal of the kernel
\[
k_t(x,y):= \Big(\sqrt t c(\dot A)e^{-t D^2_A}\Big)(x,y).
\]
\end{dfn}

With this notation our goal is now to compute
\[
\LIM_{t\to 0}\, \tr_{\gL^\bullet T^*M\otimes E}\big[\tau
k_t(x,x)\big].
\]
According to \eqref{FormsSuperDecomp} we can decompose
\[
\End(\gL^\bullet T^*M\otimes E) = \cl(T^*M)\otimes_s
\End_{\cl}\big(\gL^\bullet T^*M\otimes E\big),
\]
where
\[
\End_{\cl}\big(\gL^\bullet T^*M\otimes E\big)=
\widehat\cl(T^*M)\otimes \End(E).
\]
The local index theory proof of Proposition \ref{EtaConnVar} will
follow from the following odd dimensional version of \cite[Thm.
4.1]{BGV}, see also Remark \ref{GetlzerLocInd}.

\begin{theorem}\label{LocalIndOdd}
Let $M$ be a Riemannian manifold of odd dimension $m$, and let
$k_t(x,x)$ be as in Definition \ref{OddHeatKernel}. There is an
asymptotic expansion
\[
k_t(x,x) \sim (4\pi t)^{-\frac m2} \sum_{n=0}^\infty
t^{\frac{n+1}2} k_n(x),\quad\text{as }t\to 0,
\]
such that
\[
k_n(x)\in C^\infty\Big(M,\cl_{n+1}(T^*M)\otimes
\End_{\cl}\big(\gL^\bullet T^*M\otimes E\big)\Big).
\]
With respect to the symbol map $\boldsymbol{\gs}:\cl(T^*M)\to
\gL^\bullet T^*M$, one has
\[
\sum_{n=0}^{m-1} \boldsymbol{\gs}(k_n) =
{\det}^{1/2}\left(\frac{R^g/2}{\sinh(R^g/2)}\right)\wedge
\big(\dot A\wedge \exp(-F)\big),
\]
where $R^g$ is the Riemann curvature tensor, and $F$ is the
twisting curvature, defined in a local orthonormal frame $\{e_i\}$
for $TM$ by
\[
F:= F_{\nabla^{A,g}} - \lfrac 18 g\big(R^g(e_i,e_j)e_k,e_l\big)
e^i\wedge e^j\otimes c^k c^l \in
\gO^2\big(M,\End_{\cl}\big(\gL^\bullet T^*M\otimes E\big)\big).
\]
Here, $\nabla^{A,g}$ is the connection on $\gL^\bullet T^*M\otimes
E$ induced by the Levi-Civita connection and the connection $A$
on $E$.
\end{theorem}

Before we sketch how Theorem \ref{LocalIndOdd} can be proved along
parallel lines as in \cite[Ch. 4]{BGV}, we deduce the following
consequence which gives an alternative proof of Proposition
\ref{EtaConnVar}.

\begin{prop}\label{EtaConnVarLoc}
Assume that $M$ is a closed, oriented Riemannian manifold of odd
dimension $m$. Let $A_u$ be a smooth path of connections on a
fixed Hermitian vector bundle $E\to M$. Then the reduced
$\xi$-invariant associated to the family $B_{A_u} =
\tau\big(d_{A_u}+d_{A_u}^t\big)$ satisfies
\[
\lfrac d{du}[\xi(B_{A_u})]  = 2^{\frac {m+1}2}\int_M \widehat
L(TM,\nabla^g)\wedge \lfrac {i}{2\pi} \tr_E\big[\lfrac d{du}A_u
\wedge \exp(\lfrac {i}{2\pi}F_{A_u})\big],
\]
where $F_{A_u}$ is the curvature or $A_u$.
\end{prop}

\begin{proof}
For fixed $u$ we use the notation of Definition
\ref{OddHeatKernel}, and let $A:=A_u$ and $\dot A:= \lfrac
d{du}A_u$. Then Proposition \ref{RedEtaDerApp} shows that with the
kernel $k_t(x,x)$
\begin{equation}\label{EtaConnVarLoc:0}
\lfrac d{du}[\xi(B_{A_u})] = -\frac 1{\sqrt \pi} \LIM_{t\to 0}
\int_M \tr_{\gL^\bullet T^*M\otimes E}\big[\tau k_t(x,x)\big]
\vol_M.
\end{equation}
It follows from Theorem \ref{LocalIndOdd} that
\[
\LIM_{t\to 0} \int_M \tr_{\gL^\bullet T^*M\otimes E}\big[\tau
k_t(x,x)\big] \vol_M = (4\pi)^{-m/2} \int_M \tr_{\gL^\bullet
T^*M\otimes E}\big[\tau k_{m-1}\big] \vol_M.
\]
As in Proposition \ref{TraceSuperDecomp} we can decompose
\[
\tr_{\gL^\bullet T^*M\otimes E} = \tr_{\cl}\otimes \widehat
{\tr_{\cl}}\otimes \tr_E.
\]
Now, the definition of $\tr_{\cl}$ is in such a way that for
$a\in \cl(T^*M)$
\begin{equation*}\label{EtaConnVarLoc:1}
\tr_{\cl}(\tau a) = \begin{cases}\hphantom{2,} 0, &\text{if $a\in
\cl_{m-1}(T^*M)$},\\ 2^{m/2},&\text{if $a=\tau$.}
\end{cases}
\end{equation*}
This implies that for all $\gk \in
C^\infty\big(M,\End(T^*M\otimes E)\big)$
\[
\int_M \tr_{\gL^\bullet T^*M\otimes E}\big[\tau \gk\big] \vol_M =
2^{m/2}\int_M i^{-\frac{m+1}2}(\widehat {\tr_{\cl}}\otimes
\tr_E)\big[\boldsymbol{\gs}(\gk)\big],
\]
where the factor $i^{-\frac{m+1}2}$ arises from the fact that
$\boldsymbol{\gs}(\tau) = i^{\frac{m+1}2}\vol_M$. Then, we can use
Theorem \ref{LocalIndOdd} again to infer that
\[
\begin{split}
(4\pi)^{-m/2} &\int_M \tr_{\gL^\bullet T^*M\otimes E}\big[\tau
k_{m-1}\big] \vol_M = \\ & \sqrt{2\pi}(2\pi i)^{-\frac{m+1}2}
\int_M {\det}^{1/2}\left(\frac{R^g/2}{\sinh(R^g/2)}\right)\wedge
(\widehat {\tr_{\cl}}\otimes \tr_E) \big[\dot A\wedge
\exp(-F)\big].
\end{split}
\]
Now, we decompose $F=\widehat R^S+ F_A$, where $\widehat R^S$ is
the twisting curvature of $\gL^\bullet T^*M$ as in Lemma
\ref{FormsTwistCurv}, and $F_A\in \gO^2\big(M,\End(E)\big)$ is
the curvature of $A$. Then one computes
\[
\begin{split}
(2\pi i)^{-\frac{m+1}2}& \int_M
{\det}^{1/2}\left(\frac{R^g/2}{\sinh(R^g/2)}\right)\wedge
(\widehat {\tr_{\cl}}\otimes \tr_E) \big[\dot A\wedge
\exp(-F)\big]\\ &= - \int_M \widehat A(TM,\nabla^g)\wedge
(\widehat {\tr_{\cl}}\otimes \tr_E) \big[\lfrac{i}{2\pi}\dot
A\wedge \exp(\lfrac{i}{2\pi} F)\big]\\ &= - \int_M \widehat
A(TM,\nabla^g)\wedge\widehat{\tr_{\cl}}\big[\exp\big(\lfrac
i{2\pi}\widehat R^S\big)\big]\wedge\tr_E \big[\lfrac{i}{2\pi}\dot
A\wedge \exp\big(\lfrac{i}{2\pi} F_A\big)\big]\\
&= -2^{m/2} \int_M \widehat L(TM,\nabla^g)\wedge\tr_E
\big[\lfrac{i}{2\pi}\dot A\wedge \exp\big(\lfrac{i}{2\pi}
F_A\big)\big],
\end{split}
\]
where we have used Proposition \ref{SignIndDens} in the last
line. Using \eqref{EtaConnVarLoc:0}, the proof of Proposition
\ref{EtaConnVarLoc} is finished.
\end{proof}

\noindent\textbf{Getzler's Rescaling.} We now want to motivate why
Theorem \ref{LocalIndOdd} can be proved in the same way as the
local index theorem \cite[Thm. 4.1]{BGV}. Since it is also basic
for the considerations in the next chapters, we first want to
extract one of the main ideas of Getzler's approach in
\cite{Get86}. This is to consider an appropriate rescaling of the
Riemannian metric.

Let $M$ be a closed manifold, and let $g$ be a Riemannian metric.
For $t>0$ consider the rescaled metric $g_t :=t^{-1} g$. We define
a rescaled Clifford multiplication by
\begin{equation}\label{RescaledClifford}
c_t: \big(T^*M,g_t\big)\to \big(\gL^\bullet T^*M, g\big),\quad
c_t(\xi) := \sqrt t\big(\emu(\xi)-\imu(\xi)\big),
\end{equation}
where $\imu(\xi)$ denotes inner multiplication by $\xi$ with
respect to the fixed metric $g$ on $\gL^\bullet T^*M$. One easily
checks that
\begin{equation}\label{RescaledCliffProp}
c_t(\xi)^2 = -t\cdot |\xi|_{g}^2 =
-|\xi|_{g_t}^2,\quad\text{and}\quad c_t(\xi)^* = -c_t(\xi)\quad
\text{w.r.t. $g$}.
\end{equation}
This means that $c_t$ defines a Clifford structure for
$(T^*M,g_t)$ on the bundle $\big(\gL^\bullet T^*M, g\big)$. Since
the Levi-Civita connection $\nabla^g$ is invariant under
rescaling with a constant parameter, it defines a Clifford
connection for every $t$,
\[
\big[\nabla^g,c_t(\xi)\big] = c_t\big(\nabla^g\xi\big),\quad
\xi\in \gO^1(M).
\]
We thus get a family of de Rham operators $D_t$, defined in a
local orthonormal frame $\{e_j\}$ for $(TM,g)$ by
\begin{equation}\label{RescaledDeRhamOp}
D_t:=c_t(e^j)\nabla^g_{e_j} = \sqrt t(d+d^t):
\gO^\bullet(M)\to\gO^\bullet(M).
\end{equation}
Here, $d^t$ is the adjoint differential with respect to the fixed
metric $g$.

\begin{remark*}
It might be confusing that the de Rham differential $d$ is also
rescaled, although it is defined without using a metric. This is
due to the fact that we are fixing $g$ as a reference metric on
$\gL^\bullet T^*M$ while varying the metric on $T^*M$.
\end{remark*}

We wish to make this more precise, and take $g_t$ as a metric on
$\gL^\bullet T^*M$ rather than the fixed metric $g$. Consider
\begin{equation}\label{RescaledClifford:alt}
\widetilde c_t: \big(T^*M,g_t\big)\to \big(\gL^\bullet T^*M,
g_t\big),\quad \widetilde c_t(\xi) := \emu(\xi)-\imu_t(\xi),
\end{equation}
where now, $\imu_t(\xi)$ is inner multiplication with respect to
$g_t$. Note that $\imu_t(\xi)$ is related to inner multiplication
with respect to $g$ via $\imu_t(\xi) = t\cdot \imu(\xi)$. From
this, one deduces that
\begin{equation}\label{RescaledCliffProp:alt}
\widetilde c_t(\xi)^2 = -t\cdot |\xi|_{g}^2 =
-|\xi|_{g_t}^2,\quad\text{and}\quad c_t(\xi)^* = -c_t(\xi)\quad
\text{w.r.t. $g_t$}.
\end{equation}
Hence, $\widetilde c_t$ defines a Clifford structure for
$(T^*M,g_t)$ on the bundle $\big(\gL^\bullet T^*M, g_t\big)$. The
relation between $c_t$ and $\widetilde c_t$ is as follows.

\begin{lemma}\label{GetzlerRescalingPrep}
There is an isometry of vector bundles given by
\[
\gd_t: \big(\gL^\bullet T^*M,g_t\big) \to \big(\gL^\bullet
T^*M,g\big),\quad \gd_t(\ga):= (\sqrt t)^{|\ga|} \ga,
\]
where $\ga$ is homogenous of degree $|\ga|$. Moreover, the
Clifford structures $c_t$ and $\widetilde c_t$ are related by
\[
\widetilde c_t = \gd_t^{-1} \circ c_t \circ \gd_t.
\]
\end{lemma}

The proof is straightforward and is left to the reader. We also
note that if $d^t_{g_t}$ denotes the adjoint differential with
respect to $g_t$, then
\[
\sqrt t(d+d^t) = \gd_t\circ\big(d+d^t_{g_t}\big)\circ\gd_t^{-1}.
\]
This explains in more detail why the de Rham differential in
\eqref{RescaledDeRhamOp} is rescaled.

Now, one of the ideas underlying Getzler's approach in
\cite{Get86} can be expressed in the simple identity
\begin{equation}\label{GetzlerApproach}
e^{-tD^2}=e^{-sD_t^2}|_{s=1},
\end{equation}
where $D$ is the de Rham operator with respect to $g$, and $D_t$
is the rescaled de Rham operator \eqref{RescaledDeRhamOp}. The
deep insight behind \eqref{GetzlerApproach} is that the
asymptotic expansion as $t\to 0$ of the kernel of $e^{-tD^2}$ can
be related to the Euclidean heat kernel since---very
roughly---the metric $g_t$ converges locally to a Euclidean metric
as $t\to 0$.\\

\noindent\textbf{Remarks on the Proof of Theorem
\ref{LocalIndOdd}.} We shall make the above idea only a bit more
precise, and refer to \cite[Sec. 4.3]{BGV} for more details. Let
$x\in M$, and let $V:=T_xM$. We choose a ball $U\subset V$ of
radius less than the injectivity radius of $M$, so that
$\exp_x:U\to M$ parametrizes a normal neighbourhood of $x$ in
$M$. Using parallel transport along geodesic rays with respect to
the flat connection $A$, we can identify the bundle $E|_U$ with
the trivial bundle $U\times E_x$. Hence, for $y\in U$ we can
consider
\begin{equation}\label{LocalHeatKernel}
h(t,y):= e^{-tD_A^2}\big(\exp_x(y),x\big) \in \End(\gL^\bullet
V^*\otimes E_x),
\end{equation}
and
\begin{equation}\label{LocalOddHeatKernel}
k(t,y):= \Big(\sqrt t c(\dot A)e^{-t
D^2_A}\Big)\big(\exp_x(y),x\big) \in \End(\gL^\bullet V^*\otimes
E_x).
\end{equation}
Using the symbol map $\boldsymbol{\gs}:\cl(V^*)\to \gL^\bullet
V^*$, we get sections
\[
\boldsymbol{\gs}(h), \boldsymbol{\gs}(k) \in
C^\infty\big(\R^+\times U, \gL^\bullet V^*\otimes
\End_{\cl}(\gL^\bullet V^*\otimes E_x)\big).
\]
Now Getzler's rescaling method can be described as follows: We
fix $U$ as the coordinate space, but replace the metric on $M$ by
$g_t:=t^{-1}g$. This implies that the new system of normal
coordinates is given by
\[
\exp_x\circ \sqrt t: U\to M,\quad y\mapsto \exp(\sqrt t y).
\]
On $\End_{\cl}(\gL^\bullet V^*\otimes E_x)$, we fix the reference
metric given by $g$ and the metric on $E$, but we use the
rescaled metric on $\cl(V^*)$. According to Lemma
\ref{GetzlerRescalingPrep} this means that the symbol map has to
be replaced by
\[
\boldsymbol{\gs}_t:= \gd_t^{-1}\circ \boldsymbol{\gs}: \cl(V^*)\to
\gL^\bullet V^*.
\]
One checks that if $\xi\in V^*\subset \cl(V^*)$ and $a\in
\cl(V^*)$, then
\begin{equation}\label{GetzlerRescaleEff}
\boldsymbol{\gs}_t(\xi\cdot a) = \lfrac 1{\sqrt t} \xi\wedge
\boldsymbol{\gs}_t(a) - \sqrt t\imu(\xi) \boldsymbol{\gs}_t(a),
\end{equation}
where $\imu(\xi)$ is interior multiplication with respect to the
fixed metric $g$.

\begin{dfn}\label{RescaledHeatKernel}
Let $s\in \R^+$ be an auxiliary parameter as in
\eqref{GetzlerApproach}. Then the \emph{rescaled heat kernel} is
defined as
\[
h(s,t,y):=t^{m/2}\boldsymbol{\gs}_t\big(h(st,\sqrt t y)\big)\in
\gL^\bullet V^*\otimes \End_{\cl}(\gL^\bullet V^*\otimes E_x).
\]
\end{dfn}

As remarked in \cite[p. 155]{BGV} the extra factor $t^{m/2}$
enters because $h(t,y)$ is a density in the $y$ variable. Now
Getzler's local index theorem can be reformulated as follows, see
\cite[Thm. 4.21]{BGV}.

\begin{theorem}\label{LocalInd:Alt}
The limit $\lim_{t\to 0}h(s,t,y)$ exists. For $y=0$ and $s=1$,
\[
\lim_{t\to 0}h(1,t,0) = (4\pi)^{-m/2}
{\det}^{1/2}\left(\frac{R^g_x/2}{\sinh(R^g_x/2)}\right)\wedge
\exp(-F_x),
\]
where $R^g_x$ is the Riemann curvature tensor at $x$ and $F$ is
the twisting curvature, defined in Theorem \ref{LocalIndOdd}.
\end{theorem}

\begin{remark*}
In \cite{BGV}, Theorem \ref{LocalInd:Alt} is proved under the
assumption that $m$ is even. However, a close examination of all
intermediate steps shows that this restriction is not necessary.
In particular, \cite[Prop. 4.12]{BGV} continues to hold for odd
$m$ as well. The argument is similar to \eqref{FormsTwistChern:1}
in the proof of Lemma \ref{FormsTwistChern}.
\end{remark*}

Using \eqref{LocalOddHeatKernel} we now define as in Definition
\ref{RescaledHeatKernel}
\[
k(s,t,y):=t^{m/2} \boldsymbol{\gs}_t\big(k(st,\sqrt t y)\big)\in
\gL^\bullet V^*\otimes \End_{\cl}(\gL^\bullet V^*\otimes E_x).
\]
For $y\in U$ let $\dot A_y\in \gL^1V^*\otimes \End(E_x)$ denote
the pullback of $\dot A$ to $U$ at the point $y$. Then
\eqref{GetzlerRescaleEff} shows that
\[
k(s,t,y) = \dot A_{\sqrt t y}\wedge h(s,t,y)  - t\imu\big(\dot
A_{\sqrt t y}\big)h(s,t,y),
\]
and Theorem \ref{LocalInd:Alt} implies

\begin{cor}\label{LocalIndOdd:Alt}
The limit $\lim_{t\to 0}k(s,t,y)$ exists. For $y=0$ and $s=1$,
\[
\lim_{t\to 0}k(1,t,0) = (4\pi)^{-m/2}
{\det}^{1/2}\left(\frac{R^g_x/2}{\sinh(R^g_x/2)}\right)\wedge
\big(\dot A_x\wedge \exp(-F_x)\big).
\]
\end{cor}

On the other hand---as in the even dimensional case---the rescaled
kernel satisfies
\[
\lim_{t\to 0}k(1,t,0) = (4\pi)^{-\frac m2} \sum_{n=0}^{m-1}
\boldsymbol{\gs}\big(k_n(x)\big),
\]
where the $k_n$ are the coefficients appearing in the asymptotic
expansion in Theorem \ref{LocalIndOdd}. This together with
Corollary \ref{LocalIndOdd:Alt} finishes the outline of the proof
of Theorem \ref{LocalIndOdd}.

\cleardoublepage
\chapter{Rho Invariants of Fiber Bundles, Basic
Considerations}\label{ChapBasic}

In this chapter we start to work on the main topic of this thesis.
Our concern is to investigate how the structure of an oriented
fiber bundle of closed manifolds can be used to analyze Rho
invariants of its total space.

For this reason, we will first give a detailed summary of some
geometric preliminaries as they appear in the theory of Riemannian
submersions and in Bismut's local index theory for families.
Since we are dealing with the odd signature operator, our emphasis
is to understand the structure of the space of differential forms.
In particular, we obtain descriptions of the exterior
differential, the adjoint differential and the Levi-Civita
connection, which account for the special situation arising in
the context of the total space of a fiber bundle.

After having established how the odd signature operator can be
expressed in terms of a submersion metric, we describe the idea
of an \emph{adiabatic metric} on a fiber bundle. This is the main
tool on which our discussion of Rho invariants of fiber bundles
relies. Very roughly, the idea is to rescale the metric on the
base manifold in order to deform the geometry of the fiber bundle
to an ``almost product'' situation. Using the variation formula,
we shall see that the Eta invariant of the odd signature operator
has a well-defined limit under this process, which is called the
\emph{adiabatic limit}. As the Rho invariant is independent of
the underlying metric, we can then replace the Eta invariants in
its definition by adiabatic limits.

With this idea in mind, we will analyze the first class of
examples, namely principal circle bundles over closed surfaces.
Thanks to the low dimensions of fiber and base as well as the
enhanced symmetry provided by the principal bundle structure, one
can compute Rho invariants without having to use the more advanced
theory of Chapter \ref{ChapAbst}. Nevertheless, adiabatic metrics
will already play a prominent role in the discussion.

\clearpage

\section{Fibered Calculus}\label{FiberedCalc}

\subsection{Connections on Fiber Bundles.}\label{FiberedConn} Let
$F\hookrightarrow M\xrightarrow{\pi} B$ be an oriented fiber
bundle, where all manifolds are assumed to be closed, connected
and oriented. Let $T^vM:=\ker \pi_*$ be the vertical subbundle of
$TM$. Then $T^vM$ is \emph{involutive}, i.e., if $[.,.]$ denotes
the Lie bracket on $C^\infty(M,TM)$, we have
\[
[U,V] \in T^vM,\quad U,V\in C^\infty(M,T^vM).
\]
We now assume that $M$ is endowed with a \emph{connection}, i.e.,
a vertical projection
\[
P^v:TM\to T^vM.
\]
This induces a splitting
\[
TM = T^hM\oplus T^vM, \quad T^hM:=\ker P^v,
\]
and $\pi_*:T^hM\to \pi^*TB$ is an isomorphism. In the following we
usually identify $T^hM$ and $\pi^*TB$ via this isomorphism. Given
a vector field $X\in C^\infty(B,TB)$ we can use the connection to
lift $X$ horizontally to a vector field $X^h\in C^\infty(M,T^hM)$.
We will frequently use the following easy result, see \cite[Lem.
10.7]{BGV}.

\begin{lemma}\label{HorizontalVerticalLie}
Let $V\in C^\infty(M,T^vM)$ be a vertical vector field on $M$, and
let $X\in C^\infty(B,TB)$ be a vector field on $B$. Then
\[
[X^h,V]\in C^\infty(M,T^vM).
\]
\end{lemma}

The horizontal distribution $T^hM$ is in general not involutive.
The following quantity measures the failure of being so.

\begin{dfn}\label{CurvFiber}
The \emph{curvature form} of the connection $P^v$ is given by
\begin{equation*}
\gO\in C^\infty(M, \gL^2 T^hM^*\otimes T^vM),\quad \gO(X,Y):= -
P^v\big([X,Y]\big),
\end{equation*}
where $X,Y\in C^\infty(M,T^hM)$.
\end{dfn}

\noindent\textbf{Riemannian Connections on Fiber Bundles.} We now
assume that the fiber bundle is equipped with a \emph{submersion
metric}. This means that with respect to the splitting
$TM=\pi^*TB\oplus T^vM$ induced by the connection, we take a
metric of the form
\begin{equation}\label{g:split}
g= \pi^*g_B\oplus g_v,
\end{equation}
where $g_v$ is a Riemannian metric on $T^vM$, and $g_B$ is a
metric on $B$. We will frequently write $g= g_B\oplus g_v$, the
pullback and the identification $TM=\pi^*TB\oplus T^vM$ being
understood.

\begin{remark*}
Note that our point of view is somewhat reversed to the situation
encountered in differential geometry. Recall, see e.g. \cite[pp.
212--214]{ON}, that a \emph{Riemannian submersion} is defined  as
a submersion $\pi:(M,g)\to (B,g_B)$ such that the push-forward
$\pi_*:(\ker \pi_*)^\perp\to TB$ is an isometry. Then one deduces
that $M$ is a fiber bundle over $B$, endowed with a natural
connection given by the orthogonal projection onto $\ker \pi_*$.
Moreover, via the isomorphism $(\ker \pi_*)^\perp\cong \pi^*TB$
the metric $g$ is of the form $g_B\oplus g_v$. In contrast to
this point of view, we start with a fiber bundle, choose a
connection and then endow the vertical and the horizontal bundles
with metrics. This is because we often want to fix only a
vertical metric.
\end{remark*}

Given a vertical metric $g_v$, the following result provides a
natural connection $\nabla^v$ on $T^vM$, see \cite[Prop.
10.2]{BGV}.

\begin{prop}\label{CanVerticalConn}
Let $g_v$ be a metric on $T^vM$, and let $P^v$ be a vertical
projection. Then there is a natural connection $\nabla^v$ on
$T^vM$ defined by
\[
\nabla^v := P^v\circ \nabla^g\circ P^v,
\]
where $\nabla^g$ is the Levi-Civita connection associated to a
metric of the form \eqref{g:split}. The connection $\nabla^v$ is
independent of the choice of the metric $g_B$ on $B$. It is
compatible with $g_v$ and torsion free when restricted to $T^vM$,
i.e.,
\[
\nabla^v_UV - \nabla^v_VU = [U,V],\quad U,V\in C^\infty(M,T^vM).
\]
\end{prop}

\begin{remark}\label{NablaNatRem}
The connection $\nabla^v$ can be thought of as a family of
Levi-Civita connections parametrized by $B$, when we consider
every fiber $F\subset M$ endowed with the metric induced by
$g_v|_F$. This might give a good intuition why $\nabla^v$ is
canonically associated to $g_v$ and independent of $g_B$.
\end{remark}

\begin{dfn}\label{MeanCurvFiber}
Let $g_v$ be a vertical Riemannian metric. Then we define the
\emph{Weingarten map}
\[
W\in C^\infty\big(M,T^hM^*\otimes \End(T^vM)\big),\quad W_X(V) :=
\nabla_X^vV - P^v\big([X,V]\big),
\]
where $V\in C^\infty(M,T^vM)$ and $X\in C^\infty(M,T^hM)$. The
\emph{mean curvature} of $\nabla^v$ with respect to the vertical
projection $P^v$ is defined as
\begin{equation*}
k_v=k_v(P^v,g_v)\in C^\infty(M,T^hM^*),\quad k_v(X):=\tr_v(W_{X}),
\end{equation*}
where $X\in C^\infty(M,T^hM)$ and $\tr_v: \End(T^vM)\to \C$ is the
fiberwise vertical trace.
\end{dfn}

\begin{remark*}\quad\nopagebreak
\begin{enumerate}
\item If $\{e_i\}$ is a local orthonormal frame for $T^vM$, we
have
\begin{equation}\label{MeanCurvRem}
\tr_v(W_X) = \sum_i g_v\big(\nabla_X^ve_i - [X,e_i],e_i \big) =
-\sum_i g_v\big([X,e_i],e_i \big).
\end{equation}
This is because $\nabla^v$ is a metric connection so that
\[
g_v(\nabla^v_Xe_i,e_i) = Xg_v(e_i,e_i) - g_v(e_i,\nabla^v_Xe_i) =
-g_v(\nabla^v_Xe_i,e_i).
\]
\item In the literature, one finds different conventions of
how to define the Weingarten map. First of all in Riemannian
geometry, one usually defines it as the negative of what we have
defined. Moreover, one often normalizes the mean curvature by a
factor of $(\dim F)^{-1}$. We chose to follow the conventions of
\cite[Sec. 10.1]{BGV}.
\end{enumerate}
\end{remark*}

Associated to the metric $g_B$ on $B$, there is the Levi-Civita
connection $\nabla^B$. Apart from the Levi-Civita connection
$\nabla^g$ on $TM$ we can thus form the direct sum connection
\begin{equation}\label{SplitConn}
\nabla^\oplus := \pi^*\nabla^B\oplus \nabla^v
\end{equation}
with respect to the splitting $TM=\pi^*TB\oplus T^vM$. Note that
if $X^h$ is a horizontal lift and if $V$ is vertical, then
$\nabla^\oplus_V X^h =0$. Clearly, the connection $\nabla^\oplus$
preserves the metric $g=g_B\oplus g_v$. However, it is not
necessarily torsion free. Hence, it does usually not coincide
with the Levi-Civita connection of $M$. Following \cite[Sec.
I]{B86} we introduce the following natural tensors which measure
the difference of $\nabla^\oplus$ and $\nabla^g$.

\begin{dfn}\label{SThetaDef}
Let $\nabla^\oplus$ be defined as in \eqref{SplitConn} and let
$\nabla^g$ be the Levi-Civita connection associated to the metric
\eqref{g:split}. Define
\begin{equation*}
S:= \nabla^g - \nabla^\oplus\in C^\infty\big(M,\End(TM)\big),
\end{equation*}
and let $\gt$ be its metric contraction
\begin{equation*}
\gt(X)(Y,Z):= g\big(S(X)Y,Z\big),\quad X,Y,Z\in C^\infty(M,TM).
\end{equation*}
\end{dfn}
Since both connections $\nabla^\oplus$ and $\nabla^g$ preserve
the metric, the tensor $\gt$ is antisymmetric in $Y$ and $Z$.
Therefore, it is a section of $T^*M\otimes \gL^2T^*M$. Also note
that the above tensors are related to the \emph{O'Neill tensors}
of the Riemannian submersion, see \cite[pp. 212--214]{ON}. The
following result describes all non-trivial components of $\gt$,
see \cite[Sec. 4 (a)]{BC89}.

\begin{prop}\label{PropFiber}
The tensor $\gt$ is independent of the chosen metric $g_B$ on
$B$. Moreover, if $X$, $Y$ are horizontal vector fields, and
$U,V$ are vertical vector fields, then
\begin{equation*}
\gt(X)(V,Y) = \gt(V)(X,Y) = \lfrac12 g\big(\gO(X,Y),V\big),\quad
\gt(U)(V,X) =  g\big(\nabla^g_U V, X\big),
\end{equation*}
where $\gO$ is the curvature of the fiber bundle, see Definition
\ref{CurvFiber}.
\end{prop}

We also note the following formul{\ae} for $\gt$ and the mean
curvature $k_v$.

\begin{lemma}\label{MeanCurvFiber:alt}
Let $g_v$ be a vertical Riemannian metric. Then for $U$,$V$
vertical and $X$ horizontal
\[
\gt(U)(V,X) = -\lfrac 12 \sL_X (g_v)(U,V),
\]
where $\sL_X$ denotes the Lie derivative in the $X$ direction. If
$k_v\in C^\infty(M,T^hM^*)$ is the mean curvature of $g_v$, then
\[
k_v(X) = \lfrac 12 \tr_v \big[\sL_X (g_v)\big] = -
\sum_i\gt(e_i)(e_i,X),
\]
where $\{e_i\}$ is an arbitrary local orthonormal frame for
$T^vM$.
\end{lemma}

\begin{proof}
Let $U,V$ be vertical and $X$ horizontal. Then by definition of
the Lie derivative
\begin{equation}\label{MeanCurvFiber:alt:1}
\begin{split}
\sL_X (g_v)(U,V) &= X g_v(U,V) - g_v\big([X,U],V\big) -
g_v\big(U,[X,V]\big)\\
&= g_v\big(\nabla^v_XU - [X,U],V\big) +
g_v\big(U,\nabla^v_XV-[X,V]\big),
\end{split}
\end{equation}
where we have used that $\nabla^v$ is a metric connection. This
shows that
\[
\tr_v\big[\sL_X (g_v)\big] = 2 \sum_i g_v\big(\nabla^v_Xe_i -
[X,e_i],e_i\big) = 2 \tr_v (W_X) = 2\,k_v(X).
\]
Now choose a metric $g_B$ on $B$ and endow $M$ with the submersion
metric $g=g_B\oplus g_v$. By definition of $\nabla^v$, we can
replace $g_v$ and $\nabla^v$ in \eqref{MeanCurvFiber:alt:1} with
$g$ and $\nabla^g$, respectively. Since $\nabla^g$ is torsion
free, we find
\[
\begin{split}
\sL_X (g_v)(U,V) &= g\big(\nabla^g_XU - [X,U],V\big) +
g\big(U,\nabla^g_XV-[X,V]\big)\\
&= g\big(\nabla^g_UX,V\big) + g\big(U,\nabla^g_VX\big)\\
&= g\big(X,- \nabla^g_UV- \nabla^g_VU\big) = -2
g\big(X,\nabla^g_UV\big) + g\big(X,[U,V]\big).
\end{split}
\]
Now, $[U,V]$ is vertical, so that we find using Proposition
\ref{PropFiber} that
\[
\sL_X (g_v)(U,V) = -2 g\big(X,\nabla^g_UV\big) = -
2\gt(U)(V,X).\qedhere
\]
\end{proof}

\begin{convention}\label{FramesConv}
At this point it is convenient to introduce a convention regarding
local computations. We will always denote by $\{e_i\}$ a local,
oriented frame for $T^vM$ and by $\{f_a\}$ a local, oriented frame
for $TB$. The horizontal lifts to a frame for $T^hM$ will be
denoted with the same letters. Upper indices denote the dual
frames and the summation convention will be understood. Indices
$a,b,c,\ldots$ always refer to horizontal directions and
$i,j,k,\ldots$ to vertical ones. If we have chosen metrics $g_v$
and $g_B$, we will always choose local \emph{orthonormal} frames.
Whenever we do not want to distinguish horizontal and vertical
directions, we use the notation $\{E_I\}$ for the frame
$\{f_1,f_2,\ldots,e_1,e_2,\ldots\}$ with uppercase indices
$I,J,K,\ldots$. Moreover, if $\nabla$ is a connection on a vector
bundle over $M$, we will use the abbreviations $\nabla_a$,
$\nabla_i$, $\nabla_I$ for $\nabla_{f_a}$, $\nabla_{e_i}$,
$\nabla_{E_I}$.
\end{convention}

As an example regarding this convention, we write the tensor $\gt$
of Definition \ref{SThetaDef} as
\[
\gt = \lfrac 12 \gt_{IJK}E^I\otimes E^J\wedge E^K.
\]
Then Proposition \ref{PropFiber} shows that if we distinguish
vertical and horizontal direction, the functions $\gt_{IJK}$
satisfy the relations
\begin{equation}\label{ThetaCoord}
\begin{split}
\gt_{a i b} &= - \gt_{a b i} = \gt_{i ab } = \lfrac 12
g\big(\gO_{ab},e_i\big), \\
\gt_{ija} &= - \gt_{ia j} = \gt_{jia} =
g\big(\nabla_i^ge_j,f_a\big),\\
\gt_{a jk}& = \gt_{ijk}=\gt_{abc} = 0.
\end{split}
\end{equation}

\noindent\textbf{Another Natural Connection.} The connection
$\nabla^v$ associated to a vertical Riemannian metric is not the
only natural connection on $T^vM$. For $X\in C^\infty(M,T^hM)$
consider the vertical projection of the Lie derivative $\sL_X$,
i.e.,
\[
\sL_X^v: T^vM\to T^vM,\quad \sL_X^v(V):= P^v[X,V].
\]
Note that Lemma \ref{HorizontalVerticalLie} yields that if $X$ is
a horizontal lift, then $[X,V]$ is automatically vertical. For
general $X\in C^\infty(M,T^hM)$, $V\in C^\infty(M,T^vM)$ and
$\gf\in C^\infty(M)$ we have
\[
\sL_X^v(\gf V)= P^v\big(\gf[X,V]+ (X\gf)V\big) = \gf\sL_X^v(V) +
(X\gf)V
\]
and, since $X$ is horizontal,
\[
\sL^v_{\gf X}(V) = P^v\big(\gf[X,V]- (V\gf)X \big) =
\gf\sL_X^v(V).
\]
Therefore, we can define a connection on $T^vM$ as follows:

\begin{dfn}\label{CanVerticalConn2}
Let $\nabla^v$ be the natural connection associated to a vertical
metric $g_v$. We define the connection $\widetilde\nabla^v$ on
$T^vM$ by
\begin{equation*}
\widetilde\nabla^v_U := \nabla^v_U,\quad U\in
C^\infty(M,T^vM),\quad\text{and}\quad
\widetilde\nabla^v_X:=\sL_X^v,\quad X\in C^\infty(M,T^hM).
\end{equation*}
If $g_B$ is a metric on $B$ with Levi-Civita connection
$\nabla^B$, we also define a connection on $TM=\pi^*TB\oplus TM$
via
\[
\widetilde\nabla^\oplus :=\pi^*\nabla^B\oplus \widetilde\nabla^v.
\]
\end{dfn}

Clearly, the definition of $\widetilde\nabla^v_X$ for horizontal
$X$ is independent of any choice of metric. However, for vertical
$U$ we cannot use the Lie derivative to define
$\widetilde\nabla^v_U$ which is why the connection $\nabla^v_U$
enters the definition.

\begin{lemma}\quad\label{CanVerticalConnTor}
\begin{enumerate}
\item For $X\in C^\infty(M,T^hM)$
\[
\widetilde\nabla^v_X (g_v) = \sL_X(g_v),
\]
i.e., $\widetilde\nabla^v_X$ does in general not preserve the
metric $g_v$. Similarly $\widetilde\nabla^\oplus$, does in general
not preserve the metric $g=g_B\oplus g_v$.
\item The torsion of $\widetilde\nabla^\oplus$ coincides with the
curvature of the fiber bundle,
\[
T({\widetilde\nabla^\oplus}) = \gO \in
C^\infty\big(M,\gL^2T^hM^*\otimes T^vM\big).
\]
\end{enumerate}
\end{lemma}
\begin{proof}
Part (i) is clear by definition. For part (ii) we use local
orthonormal frames $\{f_a\}$ and $\{e_i\}$ for $TB$ and $T^vM$
according to the convention on p. \pageref{FramesConv}. Then one
computes
\[
\begin{split}
\widetilde\nabla^\oplus_ie_j- \widetilde\nabla^\oplus_je_i &=
\nabla^v_ie_j - \nabla^v_je_i = P^v\big(\nabla^g_ie_j -
\nabla^g_je_i\big) = P^v[e_i,e_j] =[e_i,e_j], \\
\widetilde\nabla^\oplus_ae_i- \widetilde\nabla^\oplus_i f_a
&= \widetilde\nabla^v_ae_i = P^v [f_a,e_i] = [f_a,e_i],\quad\text{and} \\
\widetilde\nabla^\oplus_af_b- \widetilde\nabla^\oplus_bf_a &=
\pi^*\big(\nabla^B_af_b- \nabla^B_bf_a \big) = \pi^*[f_a,f_b] =
[f_a,f_b]^h = [f_a,f_b] + \gO(f_a,f_b).
\end{split}
\]
Note that in the first and third row we have used that $\nabla^g$
respectively $\nabla^B$ are torsion free, and in the second row
that $[f_a,e_i]$ is vertical as $f_a$ is a horizontal lift, see
Lemma \ref{HorizontalVerticalLie}.
\end{proof}

\subsection{Differential Forms on a Fiber Bundle.}

Let $F\hookrightarrow M\xrightarrow{\pi} B$ be an oriented fiber
bundle as above. The natural exact sequence
\[
T^vM\hookrightarrow TM\to \pi^*TB
\]
translates to exterior bundles as
\[
\pi^*\gL^pT^*B\hookrightarrow \gL^pT^*M\to \gL^p(T^vM)^*.
\]
In terms of differential forms, the pullback of forms gives a
natural inclusion
\begin{equation*}
\pi^*: \gO^p(B)\to \gO^p(M).
\end{equation*}

\begin{dfn}\label{HorBasVertDef}
Let $\go\in \gO^\bullet(M)$, and let $\imu(.)$ denote inner
multiplication with a vector field. Then $\go$ is called
\begin{enumerate}
\item \emph{horizontal} if $\imu(V)\go=0$ for all $V\in
C^\infty(M,T^vM)$,
\item \emph{basic} if $\go=\pi^*\ga$ for some
$\ga\in\gO^\bullet(B)$, and
\item \emph{vertical} if $\imu(X)\go=0$ for all $X\in
C^\infty(M,T^hM)$.
\end{enumerate}
We denote by $\gO^\bullet_v(M)$ the algebra of vertical
differential forms, and by $\gO^{\bullet}_h(M)$ the algebra of
horizontal differential forms.
\end{dfn}

\newpage

\begin{remark*}\quad
\begin{enumerate}
\item Note that the definition of vertical forms requires the choice of
a vertical projection $P^v:TM\to T^vM$, whereas horizontal and
basic forms are defined independently of such a choice. Moreover,
$P^v$ gives an identification
\[
\gF:C^\infty\big(M,\gL^q (T^vM)^*\big)\to \gO_v^q(M),\quad
\gF(\go)(X_1,\ldots,X_q) := \go(P^vX_1,\ldots,P^vX_q).
\]
We will usually suppress this isomorphism from the notation and
identify a section of $(T^vM)^*$ with a vertical differential
form.
\item A horizontal form $\go$ can always be written (non-uniquely) as a
sum
\[
\go = \sum_i \gf_i(\pi^*\ga_i),\quad \gf_i\in C^\infty(M),\quad
\ga_i\in \gO^\bullet(B).
\]
Thus, on a somewhat formal level,
\[
\gO^\bullet_h(M) \cong \pi^*\gO^\bullet(B)\otimes C^\infty(M),
\]
where the tensor product is over $\pi^*C^\infty(B)$.
\item More generally, we can decompose every differential form
$\go$ on $M$ as
\[
\go = \sum_i \ga_i\wedge\gb_i,\quad \ga_i\in\gO^\bullet_h(M),\quad
\gb_i\in\gO^\bullet_v(M),
\]
and so
\begin{equation}\label{FormsSplit}
\gO^k(M) \cong \bigoplus_{p+q=k}\pi^*\gO^p(B)\otimes_s
\gO^q_v(M)=: \bigoplus_{p+q=k}\gO^{p,q}(M).
\end{equation}
Here, $\otimes_s$ is a graded tensor product over
$\pi^*C^\infty(B)$. We will often write $\otimes$ instead of
$\otimes_s$, keeping in mind that there is a grading involved.
\item So far we have only implicitly remarked on our orientation
convention. We use the \emph{basis first} orientation which can
be described as follows: In \eqref{FormsSplit} we consider the
case $k=\dim M$, $p=\dim B$ and $q=\dim F$. Let
$\vol_B(g_B)\in\gO^p(B)$ and $\vol_F(g_v)\in \gO^q_v(M)$ be
oriented volume forms associated to metrics $g_B$ and $g_v$. Then
we orient $M$ according to the prescription
\begin{equation}\label{OrientConvent}
\vol_M(g) = \pi^*\big(\vol_B(g_B)\big)\wedge \vol_F(g_v).
\end{equation}
\end{enumerate}
\end{remark*}

\noindent\textbf{Integration Along the Fiber.} If the fiber bundle
is endowed with a connection and a vertical metric, there is a
natural right-inverse for the map $\pi^*: \gO^p(B)\to \gO^p(M)$.
We recall the following well-known facts, see e.g. \cite[Sec.
1.6]{BT}.

\begin{prop}\label{FiberInt}
There is a natural homomorphism of $C^\infty(B)$ modules,
\[
\int_{M/B}: \gO^{\bullet,\dim F}(M) \to \gO^\bullet(B),
\]
called \emph{``integration along the fiber''}, which is uniquely
defined by the property that
\[
\int_B \ga\wedge\Big( \int_{M/B} \go\Big) = \int_M \pi^*\ga\wedge
\go,\quad \ga\in \gO^{\dim B-k}(B),\quad \go\in \gO^{k,\dim F}(M).
\]
Moreover, for $\ga$ and $\go$ as above,
\[
d_B \int_{M/B}\go =  \int_{M/B} d_M\go\quad\text{and}\quad
\ga\wedge \int_{M/B}\go= \int_{M/B}\pi^*\ga\wedge\go.
\]
\end{prop}

\begin{dfn}\label{BasicProj}
Let $\vol_F(g_v)\in \gO^{\dim F}_v(M)$ be the vertical volume form
associated to a vertical metric, and let
\[
v_F(g_v):=\int_{M/B}\vol_F(g_v)\in C^\infty(B)
\]
be the function which associates to a point $y\in B$ the volume of
the fiber over $y$. Then we define the \emph{basic projection} on
horizontal forms as
\[
\Pi_B: \gO^\bullet_h(M)\to \gO^\bullet(B),\quad
\Pi_B(\go):=\frac1{v_F(g_v)} \int_{M/B}\go\wedge \vol_F(g_v).
\]
\end{dfn}

Here, the normalization factor enters since we want
$\Pi_B(\pi^*\ga) = \ga$ for every $\ga\in\gO^\bullet(B)$. We also
want to point out that if we allow conformal changes of the
vertical metric, we can easily achieve that $v_F(g_v)=1$. This is
the content of the following simple result.

\begin{lemma}\label{VolumeNorm}
Let $n:=\dim F$, and let $g_v$ be a vertical metric. Define
\[
u:= \lfrac 1 n \log\big(\pi^* v_F(g_v)\big)\in C^\infty(M).
\]
Then the metric $\widetilde g_v:=e^{-2u}g_v$ has unit volume along
the fibers, i.e., $v_F(\widetilde g_v)=1$.
\end{lemma}

\begin{proof}
Let $\{e_i\}$ be a local, oriented orthonormal frame for
$(T^vM,g_v)$. Let $u$ be defined as above, and let $\widetilde
e_i:= \exp(u) e_i$. Then $\{\widetilde e_i\}$ is a local,
oriented orthonormal frame with respect to the metric $\widetilde
g_v$, and
\[
\vol_F(\widetilde g_v) = \widetilde e^1\wedge\ldots\wedge
\widetilde e^n = (e^{-n u }) e^1\wedge\ldots\wedge e^n =
\big(\pi^*v_F(g_v)\big)^{-1}\vol_F(g_v).
\]
This yields that $v_F(\widetilde g_v)=1$.
\end{proof}

One might expect, that there is a canonical description of the
kernel of $\Pi_B$. This is indeed true. However, the
corresponding result is not completely straightforward. We give a
proof in Chapter \ref{ChapAbst}, see Proposition
\ref{BasicProjKer}. For the time being we need a better
understanding of the calculus for differential forms on fiber
bundles.

\subsection{The Exterior Differential of a Fiber Bundle.}\label{FiberedExtDiff}

Let $F\hookrightarrow M\xrightarrow{\pi} B$ be an oriented fiber
bundle of closed manifolds as before. Let $E\to M$ be a Hermitian
vector bundle which admits a flat connection $A$. We denote by
$\gO^q_v(M,E)$ the space of vertical $E$-valued $q$-forms, i.e.,
the space of sections of $\gL^q T^vM^*\otimes E$. The canonical
connection $\nabla^v$ in Proposition \ref{CanVerticalConn}
together with $A$ induces a natural connection
\[
\nabla^{A,v}: \gO^\bullet_v(M,E) \to \gO^1(M)\otimes
\gO^\bullet_v(M,E).
\]
Moreover, we get a vertical differential
\[
d_{A,v} : \gO^q_v(M,E)\to \gO^{q+1}_v(M,E), \quad d_{A,v} =
e^i\wedge \nabla^{A,v}_i,
\]
where $\{e_i\}$ is any local orthonormal frame for $T^vM$. As in
\eqref{FormsSplit}, we can split the space of $E$-valued
$k$-forms on $M$ as
\[
\gO^k(M,E) =  \bigoplus_{p+q=k}\gO^{p,q}(M,E).
\]
The vertical differential then extends to $\gO^\bullet(M,E)$ by
requiring that
\[
d_{A,v}(\ga\otimes \go) = (-1)^p \ga\otimes d_{A,v}\go,\quad
\ga\otimes \go\in \gO^{p,q}(M,E).
\]
On the other hand, the connection $A$ defines a total exterior
differential $d_A$ on $\gO^\bullet(M,E)$. It inherits a bigrading
\[
d_A = \sum_{i+j=1} d_{ij},\quad\text{where}\quad
d_{ij}:\gO^{p,q}(M,E)\to\gO^{p+i,q+j}(M,E)\quad \text{for all
$p,q$}.
\]

We now want to describe this in terms of the data introduced in
Section \ref{FiberedConn}. Let $\widetilde\nabla^v$ be the
connection on $T^vM$ as in Definition \ref{CanVerticalConn2}. It
induces a connection on vertical differential forms, which we
denote by the same letter. Similarly, we obtain a connection
$\widetilde \nabla^{\oplus}$ on $\gL^\bullet T^*M$, and using
$A$, we define $\widetilde\nabla^{A,\oplus}$ on $\gL^\bullet
T^*M\otimes E$.

\begin{remark*}
There is a subtlety concerning the action on vertical
differential forms of the connection $\widetilde\nabla^v_{X^h}$,
if $X^h$ is the horizontal lift of a vector field $X$ on $B$.
Since on $T^vM$, the action of $\widetilde\nabla^v_{X^h}$ is
given by the Lie derivative, one might expect that the same is
true for its action on forms. However, it is in general not true
that $\sL_{X^h}\go$ is automatically vertical for a vertical
differential form $\go$, compare with Lemma
\ref{HorizontalVerticalLie}. For example, if $\go\in\gO_v^1(M)$,
then the Cartan formula yields
\[
\sL_{X^h}(\go)(Y^h) = (\imu(X^h)\circ d\go)(Y^h) = -
\go\big([X^h,Y^h]\big) = \go\big(\gO(X,Y)\big).
\]
This is in general non-zero, so that $\sL_{X^h}\go$ is in general
not a vertical form. Thus, $\widetilde\nabla^v_{X^h}$ agrees in
general only with the vertical projection $\sL_{X^h}^v$ of the
Lie derivative.
\end{remark*}

As always let $\{f_a\}$ be a local frame for $TB$, and write $\gO
= \lfrac 12 f^a\wedge f^b\otimes \gO_{ab}$ for the curvature of
the fiber bundle.

\begin{prop}\label{dSplit}
The total exterior differential $d_A$ on $\gO^\bullet(M,E)$ splits
as
\begin{equation*}
d_A = d_{A,v} + d_{A,h} + \imu(\gO),
\end{equation*}
where for $\go\in \gO^\bullet(M,E)$
\[
d_{A,h}\go = f^a\wedge \widetilde
\nabla^{A,\oplus}_a\go,\quad\text{and}\quad \imu(\gO)\go = \lfrac
12 f^a\wedge f^b\wedge \imu(\gO_{ab})\go.
\]
Here, $\imu(\gO_{ab})$ denotes interior multiplication with
$\gO_{ab}\in C^\infty(M,T^vM)$.
\end{prop}

For convenience we sketch a proof, although the result is well
known, see \cite[Prop. 10.1]{BGV} or \cite[Prop. 3.4]{BL95}.
Before we do so, let us point out that there is no $d_{-1,2}$
contribution to $d_A$, which is due to the fact that the vertical
distribution is integrable, see \cite[p. 58]{Mo}.

\begin{proof}[Proof of Proposition \ref{dSplit}]
We assume for simplicity that $A$ is the trivial connection on
the trivial line bundle. Recall, e.g. from \cite[Prop. 1.22]{BGV},
that if $\nabla$ is a torsion free connection on $TM$, then $d$
can be expressed in terms of $\nabla$ as the composition
\begin{equation}\label{d=ExtTor}
\gO^\bullet(M) \xrightarrow{\nabla} \gO^1(M)\otimes \gO^\bullet(M)
\xrightarrow{\emu\circ} \gO^{\bullet+1}(M).
\end{equation}
Here, the second arrow means contraction with exterior
multiplication. Let $g$ be a metric of the form $g_B\oplus g_v$.
It follows from Lemma \ref{CanVerticalConnTor}, that we can
define a torsion free\footnote{Clearly, we could also use the
Levi-Civita connection associated to the metric $g$. However, its
explicit formula is more complicated (see also Remark
\ref{CanVertConnDiff} below).} connection $\nabla$ on $TM$ by
\[
\nabla_X:= \widetilde\nabla^\oplus_X - \lfrac12 \gO(X,.),\quad
X\in C^\infty(M,TM).
\]
Let $\{f_a\}$ and $\{e_i\}$ be local orthonormal frames for $TB$
and $T^vM$ respectively. Since $d$ satisfies the Leibniz rule it
suffices to compute $de^i$ and $df^a$. From \eqref{d=ExtTor} and
the definition of $\widetilde\nabla^\oplus$ we get
\[
\begin{split}
de^i & = e^j\wedge \nabla_je^i + f^a\wedge\nabla_ae^i\\
&= e^j\wedge\nabla^v_je^i + f^a\wedge\big(- e^i(\nabla_ae_k) e^k
-e^i(\nabla_af_b)f^b \big) \\
&= d_ve^i + f^a\wedge\big(- e^i([f_a,e_k])e^k -
e^i(-\lfrac12 \gO_{ab})f^b \big) \\
&= d_ve^i + f^a\wedge \widetilde \nabla^\oplus_ae^i + \lfrac12
 e^i(\gO_{ab}) f^a\wedge f^b.
\end{split}
\]
This is the required formula for $de^i$. On the other hand, since
$\nabla_j$ acts trivially on basic forms, and $\gO(f_b,f_c)$ has
no horizontal component, one easily finds that
\[
df^a = f^b\wedge\nabla_bf^a = f^b\wedge \widetilde\nabla^\oplus_b
f^a.\qedhere
\]
\end{proof}

Since $A$ is flat we have $d_A^2=0$. This implies the following.

\begin{cor}\label{dSquare}
Let $\{.,.\}$ denote the anti-commutator of two operators in the
ungraded sense. Then
\[
\begin{split}
d_{A,v}^2 = \imu(\gO)^2=0,&\quad d_{A,h}^2 +
\big\{d_{A,v},\imu(\gO) \big\}=0,\\ \big\{d_{A,v}, d_{A,h} \big\}
&=\big\{d_{A,h}, \imu(\gO) \big\}=0.
\end{split}
\]
\end{cor}

\noindent\textbf{More on the Mean Curvature.} From Proposition
\ref{dSplit}, we can also deduce the following formula which
relates the mean curvature with the differential of the vertical
volume form, see \cite[Lem. 10.4]{BGV}. In the theory of
foliations, this is known as \emph{Rummler's formula}, see e.g.
\cite[p. 38]{Ton}.

\begin{prop}\label{Rummler}
Let $\vol_F(g_v)$ be the volume form associated to a vertical
metric. Let $k_v\in \gO^{1,0}(M)$ be the mean curvature form. Then
\begin{equation*}
d_M \vol_F(g_v) = k_v\wedge \vol_F(g_v) + \imu(\gO)\vol_F(g_v).
\end{equation*}
\end{prop}

\begin{proof}
Since $\vol_F(g_v)$ has maximal vertical degree, Proposition
\ref{dSplit} implies that
\[
d_M \vol_F(g_v) = d_h \vol_F(g_v) + \imu(\gO) \vol_F(g_v).
\]
If $\{f_a\}$ and $\{e_j\}$ are local frames for $TB$ and $T^vM$,
we compute
\[
\begin{split}
d_h \vol_F(g_v) &= f^a \wedge \widetilde \nabla_a^v\big(
\vol_F(g_v) \big) = f^a \wedge \widetilde \nabla_a^v(e^j)\wedge
\imu(e_j) \big( \vol_F(g_v) \big) \\
&= - f^a \wedge e^j\big([f_a,e_k]\big)e^k\wedge \imu(e_j) \big(
\vol_F(g_v) \big)\\
&= - f^a \wedge e^j\big([f_a,e_k]\big)\gd^k_j\wedge
\vol_F(g_v) \\
&=- \sum_j g_v\big([f_a,e_j],e_j\big)f^a\wedge \vol_F(g_v).
\end{split}
\]
Now, \eqref{MeanCurvRem} identifies the last line with $k_v\wedge
\vol_F(g_v)$.
\end{proof}

\begin{cor}\label{BasicMeanCurv}
Let $\vol_F(g_v)$ be the volume form associated to a vertical
metric $g_v$, and let $v_F$ be the volume of the fiber as in
Definition \ref{BasicProj}. Then the basic projection of the mean
curvature form is given by
\[
\Pi_B(k_v) = d_B\log(v_F) \in \gO^1(B).
\]
\end{cor}

\begin{proof}
We differentiate $v_F$ and use Proposition \ref{Rummler} to find
that
\[
d_B v_F =d_B\int_{M/B} \vol_F(g_v)= \int_{M/B} d_M \vol_F(g_v) =
\int_{M/B} k_v\wedge \vol_F(g_v) = v_F\, \Pi_B(k_v).\qedhere
\]
\end{proof}

Corollary \ref{BasicMeanCurv} shows that the basic projection of
the mean curvature form gives a trivial element in the cohomology
of the base. Moreover, it vanishes if the metric $g_v$ has
constant volume along the fiber. As we have seen in Lemma
\ref{VolumeNorm}, this can be achieved by a conformal deformation
of the vertical metric. In Section \ref{VertCohom} we will
transfer a result of \cite{Dom98} from the theory of foliations
to the situation at hand and show that one can always deform the
vertical projection and the vertical metric of the fiber bundle
in such a way that not only the basic projection but the mean
curvature form itself vanishes.

\subsection{The Levi-Civita Connection on Forms and the Adjoint
Differential}

To study the de Rham operator on a fiber bundle $F\hookrightarrow
M\xrightarrow{\pi} B$, we also need to understand the adjoint
differential $d_A^t$, where $A$ is a flat connection on a
Hermitian vector bundle $E\to M$. For this we will use the local
formula
\[
d_A^t = -\imu(E^I)\nabla^{A,g}_{I},
\]
where $\{E_I\}$ is a local orthonormal frame for $TM$, and
$\nabla^{A,g}$ is the Levi-Civita connection on forms twisted by
$A$. We want to use this formula to split $d_A^t$ in terms of its
bidegrees with respect to the decomposition \eqref{FormsSplit}.
For this we need to relate the Levi-Civita connection
$\nabla^{A,g}$ on forms with the connection $\nabla^{A,\oplus}$.
Recall that in Definition \ref{RightCliffDef} we have introduced
a transposed Clifford as
\[
\widehat c: T^*M \to \End(\gL^\bullet T^*M),\quad \widehat c(\xi)
= \emu(\xi) + \imu(\xi).
\]

\begin{lemma}\label{LCDiff}
Let $E\to M$ be a Hermitian bundle which admits a flat connection
$A$. Then, on $\gL^\bullet T^*M\otimes E$, the difference of
$\nabla^{A,g}$ and $\nabla^{A,\oplus}$ is given by
\[
\nabla^{A,g} = \nabla^{A,\oplus} +\lfrac 12 \big(c(\gt)-\widehat
c(\gt)\big),
\]
where $\gt$ is the tensor defined in Definition \ref{SThetaDef}.
\end{lemma}
\begin{proof}
Let $\{E_I\}$ be a local orthonormal frame of $TM$ with dual
coframe $\{E^I\}$. Then by definition of $S$ and $\gt$,
\[
\begin{split}
(\nabla^{A,g}_I - \nabla^{A,\oplus}_I)  &= (\nabla^{A,g}_I -
\nabla^{A,\oplus}_I)(E^J)\wedge \imu(E_J) = -
E^J\big(S(E_I)E_K\big)\emu(E^K)\imu(E_J)\\ &= -
\gt_{IKJ}\emu(E^K)\imu(E^J).
\end{split}
\]
On the other hand,
\[
\begin{split}
\lfrac 12 \big(c(\gt_I)-\widehat c(\gt_I)\big) &= \lfrac 14
\gt_{IJK}\big(c(E^J)c(E^K) - \widehat c(E^J)\widehat c(E^K)\big)\\
&= -\lfrac 12
\gt_{IJK}\big(\emu(E^J)\imu(E^K) - \imu(E^J)\emu(E^K)\big)\\
&= - \lfrac 12\gt_{IKJ} \emu(E^K)\imu(E^J) +\lfrac12 \gt_{IJK}
\emu(E^K)\imu(E^J) = - \gt_{IKJ} \emu(E^K)\imu(E^J),
\end{split}
\]
where in the last line we have first renamed $J$ and $K$ in the
first summand and then used antisymmetry of $\gt_{IJK}$ in $J$
and $K$ for the second summand, see \eqref{ThetaCoord}.
\end{proof}

\begin{remark}\label{CanVertConnDiff}
One could use the above result to give a different proof of
Proposition \ref{dSplit} by writing out locally
\[
d_A = E^I\wedge\nabla^{A,g}_I =  E^I\wedge\nabla^{A,\oplus}_I -
\gt_{IKJ} E^I\wedge E^K\wedge\imu(E^J)
\]
Splitting this into horizontal and vertical contributions, and
using \eqref{ThetaCoord}, one verifies that
\begin{equation*}\label{dSplit:2}
d_A = d_{A,v} + f^a\wedge \big(\nabla^{A,\oplus}_a -
\gt_{kja}\,e^k\wedge \imu(e^j)\big) + \imu(\gO).
\end{equation*}
Comparing this with Proposition \ref{dSplit} we see that in
particular,
\begin{equation}\label{NablaNat}
\widetilde  \nabla^v_a = \nabla^v_a - \gt_{kja}\,e^k\wedge
\imu(e^j).
\end{equation}
A more invariant description of the term occurring here is as
follows: Define a tensor field $B\in C^\infty\big(M,T^hM^*\otimes
\End(\gL^\bullet T^vM^*)\big)$ by requiring that
\[
g_v\big(\ga,B(X)\gb\big) = \sL_X^v(g_v)(\ga,\gb),\quad X\in
C^\infty(M,T^hM),\quad\ga,\gb\in \gO^\bullet_v(M).
\]
Then Lemma \ref{MeanCurvFiber:alt} (or a direct computation)
easily implies that
\begin{equation}\label{NablaNatInv}
B(X) = X^a\gt_{kja}\,e^k\wedge \imu(e^j).
\end{equation}
\end{remark}

We now have the following analog of Proposition \ref{dSplit} for
the adjoint differential.

\begin{prop}\label{d^tSplit}
Let $E\to M$ be a Hermitian vector bundle which admits a flat
connection $A$. Then the twisted adjoint differential splits as
\begin{equation*}
d_A^t = d_{A,v}^t + d_{A,h}^t + \imu(\gO)^t,
\end{equation*}
where the terms are given in local orthonormal frames $\{e_i\}$
and $\{f_a\}$ by
\[
\begin{split}
d_{A,v}^t &= -\imu(e^i)\circ \nabla^{A,\oplus}_i:
\gO^{p,q}\to \gO^{p,q-1},\\
d_{A,h}^t &= -\imu(f^a)\circ\big(\nabla^{A,\oplus}_a
+ B(f_a) +   k_v(f_a)\big): \gO^{p,q}\to \gO^{p-1,q},\\
\imu(\gO)^t &= - \lfrac 12 \imu(f^a)\imu(f^b)\emu(\gO_{ab}) :
\gO^{p,q}\to \gO^{p-2,q+1}.
\end{split}
\]
Here, $B(f_a)$ is defined as in \eqref{NablaNatInv}, $k_v$ is the
mean curvature form, and $\emu(\gO_{ab})$ denotes exterior
multiplication with the dual of $\gO_{ab}\in C^\infty(M,T^vM)$.
\end{prop}

\begin{remark}\label{d^tSplit:alt}
The definition of the connection $\widetilde \nabla^v$ and the
relation \eqref{NablaNat} between $\nabla^v$ and $\widetilde
\nabla^v$ shows that we can write alternatively
\[
\begin{split}
d_{A,v}^t &= -\imu(e^i)\circ \widetilde\nabla^{A,\oplus}_i:
\gO^{p,q}\to \gO^{p,q-1},\\
d_{A,h}^t &= -\imu(f^a)\circ\big(\widetilde \nabla^{A,\oplus}_a +
2 B(f_a) + k_v(f_a)\big): \gO^{p,q}\to \gO^{p-1,q}.
\end{split}
\]
\end{remark}

\begin{proof}[Proof of Proposition \ref{d^tSplit}]
For convenience we drop again the reference to the flat
connection. According to Lemma \ref{LCDiff} and the local formula
for $d^t$ we have
\[
d^t = -\imu(E^I)\nabla^\oplus_I +
\imu(E^I)\gt_{IJK}\emu(E^J)\imu(E^K).
\]
Splitting this into horizontal and vertical parts and checking
bidegrees one finds
\[
d_v^t = -\imu(e^j)\nabla^\oplus_j +
\imu(e^j)\gt_{jab}\emu(f^a)\imu(f^b) +
\imu(f^a)\gt_{abj}\emu(f^b)\imu(e^j) = -\imu(e^j)\nabla^\oplus_j.
\]
Here, we have used the relation $\gt_{abj} = -\gt_{jab}$, see
\eqref{ThetaCoord}. Similarly,
\[
\begin{split}
d_h^t &= -\imu(f^a)\nabla^\oplus_a +
\imu(e^j)\gt_{jka}\emu(e^k)\imu(f^a)=
-\imu(f^a)\big(\nabla^\oplus_a - \gt_{jka}
\imu(e^j)\emu(e^k)\big)\\ &= -\imu(f^a)\Big(\nabla^\oplus_a +
\gt_{kja} \big(\emu(e^k)\imu(e^j) - \gd^{jk}\big)\Big) =
-\imu(f^a)\big(\nabla^{A,\oplus}_a + B(f_a) + k_v(f_a)\big),
\end{split}
\]
where we have used symmetry of $\gt_{jka}$ in $j$ and $k$, the
definition of $B(f_a)$, and Lemma \ref{MeanCurvFiber:alt}. Lastly,
one has
\[
\imu(\gO)^t = \imu(f^a)\gt_{ajb}\emu(e^j)\imu(f^b) = - \lfrac 12
\imu(f^a)\imu(f^b)\emu(\gO_{ab}).\qedhere
\]
\end{proof}

\noindent\textbf{Partial de Rham Operators.} Having established
the description of the adjoint differential in analogy to
Proposition \ref{dSplit}, we can now split the twisted de Rham
operator on $M$.

\begin{dfn}\label{PartDeRhamOpDef}
We use the abbreviations
\begin{equation*}
D_{A,v} := d_{A,v}+ d_{A,v}^t,\quad D_{A,h} := d_{A,h}+
d_{A,h}^t,\quad\text{and}\quad T:=\imu(\gO)+\imu(\gO)^t.
\end{equation*}
The operators $D_{A,v}$ and $D_{A,h}$ are called the
\emph{vertical} respectively \emph{horizontal} twisted de Rham
operator on $M$.
\end{dfn}

\begin{remark}\label{PartDeRhamOpRem}\quad\nopagebreak
\begin{enumerate}
\item Tautologically, the de Rham operator on $M$ splits as
\begin{equation}\label{DeRhamSplit}
D_A = D_{A,v} + D_{A,h} + T:\gO^\bullet(M,E)\to \gO^\bullet(M,E).
\end{equation}
\item The vertical de Rham operator is a first order differential
operator \emph{acting fiberwise}, i.e., $[D_{A,v},\pi^*\gf]=0$
for all $\gf\in C^\infty(B)$. This means roughly, that it can be
thought of as a smooth family of first order elliptic differential
operators on $\gL^\bullet T^*F\otimes E|_F$ parametrized by $B$,
see Definition \ref{SmoothFamily}. We refer to Section
\ref{VertCohom} for some more details. Clearly, an analogous
statement cannot be formulated for the horizontal de Rham
operator $D_{A,h}$, unless the horizontal distribution is
integrable.

\item We want to add some remarks about the effect the splitting
\eqref{DeRhamSplit} has on the spectrum of $D_A$. Let us assume
that $A$ is the trivial connection. The appearance of the mean
curvature form in the formula for $d_h^t$ shows that in general,
$D_h$ will not restrict to an operator on basic forms, see
Definition \ref{HorBasVertDef}. Instead, if $D_B$ denotes the de
Rham operator on $B$,
\begin{equation}\label{deRhamIntertwine}
D_h (\pi^*\ga) = \pi^*(D_B\ga) - \imu(k_v) \pi^*\ga,\quad \ga\in
\gO^\bullet(B),
\end{equation}
where $\imu(k_v)$ denotes interior multiplication with the mean
curvature form. However---as already pointed out---we will see in
Section \ref{VertCohom} that upon changing the vertical metric
and the horizontal distribution, we can achieve that $k_v$
vanishes. In this case \eqref{deRhamIntertwine} shows that
eigenforms of $D_B$ lift to eigenforms of $D_h$. Since $D_v$
vanishes on basic forms, this produces eigenforms of $D_v+D_h$. In
particular, $\spec(D_B)\subset \spec(D_v+D_h)$. However, the full
de Rham operator on $M$ is given by \eqref{DeRhamSplit}, and $T$
will in general not act trivially on $\pi^*\gO^\bullet(B)$. This
should give a hint at why---even in the case that $k_v$
vanishes---the relation between the spectrum of $D_M$ and the
spectra of $D_h$ and $D_v$ is non-trivial. We refer to
\cite[Ch.'s 3\&4]{GLP} for a detailed study of related questions.
\end{enumerate}
\end{remark}

\noindent\textbf{Some Commutator Relations.} The explicit
description of the adjoint differential has the following
consequence, see also \cite[Prop. 3.1]{AK01} for a generalization.

\begin{prop}\label{dvdhCommutator}
Let $\{f_a\}$ be local orthonormal frame for $TB$, and define a
bundle endomorphism
\begin{equation}\label{KDef}
K:= -\imu(f^a)\circ\big(2B(f_a) +   k_v(f_a)\big): \gL^\bullet
T^*M\to \gL^\bullet T^*M.
\end{equation}
Then
\[
d_{A,v}d_{A,h}^t + d_{A,h}^td_{A,v} = d_{A,v}K + Kd_{A,v}.
\]
\end{prop}

\begin{proof}
If $\{e_i\}$ is a local orthonormal frame for $T^vM$, one checks
that
\[
\widetilde\nabla^{A,\oplus}_a\circ d_{A,v} = d_{A,v}\circ
\widetilde\nabla^{A,\oplus}_a.
\]
Since $\widetilde \nabla^{A,\oplus}_i f^a = 0$, this implies
\[
d_{A,v}\circ\big(\imu(f^a)\circ \widetilde\nabla^{A,\oplus}_a\big)
+ \big( \imu(f^a)\circ\widetilde\nabla^{A,\oplus}_a\big)\circ
d_{A,v}=0.
\]
Now Remark \ref{d^tSplit:alt} yields that
\[
d_{A,h}^t = -\imu(f^a)\circ\big(\widetilde\nabla^{A,\oplus}_a +
2B(f_a) +   k_v(f_a)\big).
\]
Then, with $K$ as in \eqref{KDef}, one easily verifies that indeed
\[
d_{A,v}d_{A,h}^t + d_{A,h}^td_{A,v} = d_{A,v}K + Kd_{A,v}.\qedhere
\]
\end{proof}

\begin{cor}\label{FirstOrderFiber}
The anti-commutator $\{D_{A,v},D_{A,h}\}$ is a first order
differential operator acting fiberwise.
\end{cor}

\begin{proof}
According to Corollary \ref{dSquare} and the corresponding
statement for the formal adjoints, we have
\[
\{d_{A,h},d_{A,v}\} =0\quad\text{and}\quad
\{d_{A,h}^t,d_{A,v}^t\}=0.
\]
This implies
\[
\{D_{A,v},D_{A,h}\} = \{d_{A,v},d_{A,h}^t\} +
\{d_{A,v}^t,d_{A,h}\} = \{d_{A,v},K\} + \{K^t,d_{A,v}^t\},
\]
which is $C^\infty(B)$ linear and thus a first order differential
operator acting fiberwise.
\end{proof}

\section{Rho Invariants and Adiabatic Metrics}\label{RhoAdiabaticSec}

\subsection{The Odd Signature Operator}

As before let $F\hookrightarrow M\xrightarrow{\pi} B$ be an
oriented fiber bundle of closed manifolds, and let $E\to M$ be a
Hermitian vector bundle endowed with a unitary flat connection
$A$. We endow $T^vM$ with a metric $g_v$ and $B$ with a
Riemannian metric $g_B$, and consider the associated submersion
metric $g:=g_B\oplus g_v$.

If $\dim M$ is odd the odd signature operator on $M$ twisted by
$A$ is given in terms of the partial de Rham operators introduced
in Definition \ref{PartDeRhamOpDef} as
\begin{equation}\label{SignatureSplit}
B_A^{\ev} = \tau_M D_{A,v}+ \tau_M D_{A,h}+ \tau_M T
:\gO^{\ev}(M,E)\to \gO^{\ev}(M,E),
\end{equation}
where $\tau_M$ is the chirality operator on the total space $M$
of the fiber bundle. It is useful to identify $\gO^{\ev}(M,E)$ in
terms of the splitting $TM = \pi^*TB \oplus T^vM$. From
\eqref{FormsSplit} we see that
\[
\gO^{\ev}(M,E) = \sum_{p+q\equiv 0\,(2)} \pi^*\gO^p(B)\otimes
\gO^q_v(M,E).
\]
Using this identification, we define
\begin{equation*}
\gF: \gO^{\ev}(M,E)\to \pi^*\gO^\bullet(B)\otimes
\gO^\bullet_v(M,E),\quad \gF(\ga\otimes \go) = \ga^e\otimes \go^e
+\tau_M(\ga^o\otimes \go^o),
\end{equation*}
where $\ga^{e/o}$ and $\go^{e/o}$ refer to the even/odd degree
parts. Since it is straightforward, we skip the proof of the
following result.

\begin{lemma}\label{OddSignSplit}
Assume that $M$ is odd dimensional.
\begin{enumerate}
\item
If the fiber $F$ is even dimensional, then $\gF$ gives rise to an
isometry
\[
\gF:\gO^{\ev}(M,E)\xrightarrow{\cong}\pi^*\gO^{\ev}(B)\otimes
\gO^\bullet_v(M,E),
\]
and the odd signature operator is equivalent to
\[
\gF\circ B_A^{\ev}\circ \gF^{-1} = D_{A,v}+ \tau_M D_{A,h}+ T
\]
\item If $F$ is odd dimensional, then
\[
\gF:\gO^{\ev}(M,E)\xrightarrow{\cong}\pi^*\gO^\bullet(B)\otimes
\gO^{\ev}_v(M,E),
\]
and
\[
\gF\circ B_A^{\ev}\circ \gF^{-1} = \tau_M D_{A,v}+ D_{A,h}+
\tau_M T.
\]
\end{enumerate}
\end{lemma}

\noindent\textbf{The Vertical Chirality Operator.} For a more
explicit formula for the odd signature one needs to understand
how the chirality operator splits with respect to
\eqref{FormsSplit}. In the general setting at hand we will not
give a detailed account but add some remarks which will be used in
the examples below.

\begin{dfn}\label{VertChirDef}
Let $M$ be endowed with a vertical metric $g_v$, and let $n:=\dim
F$. Then the \emph{vertical chirality operator}
$\tau_v:\gO^q_v(M)\to \gO^{n-q}_v(M)$ is defined with respect to
an oriented, orthonormal frame $\{e_i\}$ for $T^vM$ by
\[
\tau_v = i^{[\frac{n+1}{2}]}c_v(e^1)\cdot\ldots\cdot c_v(e^{n}).
\]
Here, $c_v: T^vM^*\to \End\big(\gL^\bullet T^vM^*\big)$ is the
vertical Clifford multiplication,
\[
c_v(\xi) = \emu(\xi)-\imu(\xi),\quad \xi\in \gO^1_v(M).
\]
\end{dfn}

The vertical Clifford multiplication extends naturally to vertical
differential forms, and up to the normalization factor, $\tau_v$
is Clifford multiplication with the vertical volume form. In
particular, it is independent of the chosen frame. We also recall
the convention \eqref{OrientConvent} that if $g=g_B\oplus g_v$, we
orient $M$ using
\[
\vol_M(g) = \pi^*\big(\vol_B(g_B)\big)\wedge \vol_F(g_v).
\]

\begin{lemma}\label{TauSplit}
Let $\tau_B$ be the chirality operator on $\gO^\bullet(B)$, and
let $(\pi^*\ga)\wedge \go\in \gO^{p,q}(M)$.
\begin{enumerate}
\item Assume that $F$ is even dimensional. Then
\[
\tau_M (\pi^*\ga\wedge \go) = \pi^*(\tau_B\ga)\wedge \tau_v \go.
\]
\item If $F$ is odd dimensional, then
\[
\tau_M (\pi^*\ga\wedge \go) = (-1)^{p}\cdot
\begin{cases}
\pi^*(\tau_B\ga)\wedge \tau_v \go ,&\text{if $B$ is even dimensional}, \\
-i\cdot \pi^*(\tau_B \ga)\wedge \tau_v \go,&\text{if $B$ is odd
dimensional}.
\end{cases}
\]
\end{enumerate}
\end{lemma}

The proof is a bit tedious but straightforward and shall be
skipped.

\begin{remark*}
In Proposition \ref{d^tSplit} we have described the adjoint
differential $d_A^t$ using the local formula $d_A^t =
-\imu(E^I)\circ \nabla^{A,g}_I$. However, as in \eqref{d^tDef},
we also have the description
\begin{equation*}
d_A^t = (-1)^{m+1}\tau_M\circ d_A\circ \tau_M,
\end{equation*}
where $m=\dim M$. Using this together with Lemma \ref{TauSplit}
and Proposition \ref{dSplit}, one could give a different proof of
Proposition \ref{d^tSplit}. Clearly, the main point is then to
compute $\tau_M d_{A,h} \tau_M$, which amounts to proof that
\begin{equation}\label{LieTauComm}
\tau_v \big[\widetilde \nabla_X^v,\tau_v\big] = 2  B(X) +
k_v(X),\quad X\in C^\infty(M,T^hM).
\end{equation}
Conversely, \eqref{LieTauComm} can be verified using Proposition
\ref{d^tSplit} and \eqref{NablaNat}.
\end{remark*}


\subsection{Adiabatic Metrics on Fiber Bundles}\label{AdiabaticMetrics}

In a similar way as in Section \ref{VarLocalInd}, we now want to
rescale the metric on the fiber bundle $F\hookrightarrow M
\xrightarrow{\pi}B$. Yet, the important difference is that we
only rescale the metric on the base manifold. In order avoid
square roots of $\eps$, we are using $\eps^2$ rather than $\eps$
to rescale the metric.

\begin{dfn}
Let $g_B$ be a metric on $B$ and $g_v$ be a vertical metric. For
$\eps>0$ we define the \emph{adiabatic metric}
\begin{equation}\label{AdiabaticMetric}
g_\eps:= \lfrac1{\eps^2}g_B\oplus g_v.
\end{equation}
\end{dfn}

Associated to each $g_\eps$, we have a Levi-Civita connection
$\nabla^{g_\eps}$. Note that unlike in the case of a single
manifold, the family $\nabla^{g_\eps}$ is not independent of
$\eps$ since we only scale the base metric. However, the direct
sum connection $\nabla^\oplus$ is independent of $\eps$ since
both, $\nabla^B$ and $\nabla^v$ are, see Proposition
\ref{CanVerticalConn}. Similarly, the tensor $\gt$ as in
Definition \ref{SThetaDef} is independent of $\eps$.\\

\noindent\textbf{Adiabatic Families of Odd Signature Operators.}
Now let $E\to M$ be a flat unitary bundle with connection $A$, and
let $\nabla^{A,g_\eps}$ and $\nabla^{A,\oplus}$ be the induced
connections on $\gL^\bullet T^*M\otimes E$. We can use Lemma
\ref{LCDiff} to write
\begin{equation}\label{AdiabLCDiff}
\nabla^{A,g_\eps} = \nabla^{A,\oplus} +\lfrac12
\big(c_\eps(\gt)-\widehat c_\eps(\gt)\big).
\end{equation}
Here, Clifford multiplication is defined with respect to the
fixed metric $g=g_B\oplus g_v$ on $\gL^\bullet T^*M$, i.e.,
\[
c_\eps(f^a) = \eps c(f^a),\quad c_\eps(e^i) = c(e^i),\quad
\widehat c_{\eps}(f^a)=\eps \widehat c(f^a),\quad \widehat
c_{\eps}(e^i)=\widehat c(e^i),
\]
compare with \eqref{RescaledClifford}.

\begin{lemma}\label{AdiabaticClifford}
For each $\eps>0$, the connection $\nabla^{A,g_\eps}$ on
$\gL^\bullet T^*M\otimes E$ is compatible with the fixed metric
$g=g_B\oplus g_v$. Moreover, it is a Clifford connection with
respect to $c_\eps$, i.e.,
\[
\big[\nabla^{A,g_\eps},c_\eps(\xi)\big] =
c_\eps\big(\nabla^{g_\eps}\xi\big),\quad \xi\in \gO^1(M).
\]
\end{lemma}
\begin{proof}[Sketch of proof]
Lemma \ref{AdiabaticClifford} is basically \cite[Prop.
10.10]{BGV}. The main observation there is that
$\nabla^{A,\oplus}$ is compatible with $g$ and satisfies
\begin{equation*}\label{SplitConnClifford}
\big[\nabla^{A,\oplus},c_\eps(\xi)\big] =
c_\eps\big(\nabla^\oplus\xi\big),\quad \xi\in \gO^1(M).
\end{equation*}
On the other hand, according to \eqref{AdiabLCDiff}, we need to
consider
\[
c_\eps\big(\gt(X)\big)-\widehat c_\eps\big(\gt(X)\big)\in
C^\infty\big(M,\End(\gL^\bullet T^*M\otimes E)\big),\quad X\in
C^\infty(M,TM).
\]
Since $c_\eps$ and $\widehat c_\eps$ are defined with respect to
the fixed metric $g$, and since $\gt(X)$ is a 2-form, one finds
that $c_\eps\big(\gt(X)\big)$ and $\widehat
c_\eps\big(\gt(X)\big)$ are self-adjoint with respect to $g$.
This implies that $\nabla^{A,g_\eps}$ is compatible with the
metric. Lemma \ref{RightCliffLem} (i) shows that for $\xi\in
\gO^1_v(M)$
\[
\big[c_\eps(\gt)-\widehat c_\eps(\gt),c_\eps(\xi)\big] =
\big[c_\eps(\gt),c_\eps(\xi)\big]
\]
As in \cite[Prop. 10.10]{BGV} one then finds that
\[
\big[\nabla^{A,\oplus} + \lfrac 12 c_\eps(\gt),c_\eps(\xi)\big] =
c_\eps\big(\nabla^{g_\eps}\xi\big),
\]
which proves that $\nabla^{A,g_\eps}$ is indeed a Clifford
connection.
\end{proof}

\newpage

\begin{remark}\label{DeRhamRescaleRem}\quad
\begin{enumerate}
\item Again, it might be confusing that all connections
$\nabla^{A,g_\eps}$ are compatible with the fixed metric $g$ on
$\gL^\bullet T^*M$. This is due to the fact that we have defined
$\nabla^{A,g_\eps}$ in such a way, that it already incorporates
the isometry of Lemma \ref{GetzlerRescalingPrep}, which in the
case at hand takes the form
\[
\gd_\eps: \big(\gO^\bullet(M,E),g_\eps\big) \to
\big(\gO^\bullet(M,E),g\big),\quad \gd_\eps(\pi^*\ga\wedge \go):=
\eps^{|\ga|}\pi^*\ga\wedge \go.
\]
\item We also want to point out that the chirality operator $\tau_M$ on
$\big(\gO^\bullet(M,E),g\big)$ does not change with $\eps$.
Indeed, it is immediate that $\vol_M(g_\eps)= \eps^{-\dim
B}\vol_M(g)$ from which it follows that
\[
c_{\eps}\big(\vol_M(g_\eps) \big) = c\big(\vol_M(g)\big)
\]
\end{enumerate}
\end{remark}

\begin{dfn}\label{AdiabaticFamily}
Let $g_\eps$ be an adiabatic metric on $M$, and assume that
$m=\dim M$ is odd. We define the \emph{adiabatic family} of odd
signature operators as
\begin{equation}\label{AdiabaticSignatureSplit}
B_{A,\eps}^{\ev} := \tau_M D_{A,v}+ \eps\cdot \tau_M D_{A,h}+
\eps^2\cdot \tau_M T:\gO^{\ev}(M,E)\to \gO^{\ev}(M,E).
\end{equation}
\end{dfn}
The definition is in such a way that $B_{A,\eps}^{\ev}$ is given
by Clifford contraction of $\nabla^{A,g_\eps}$ with respect to the
Clifford multiplication $\tau_M\circ c_\eps$. Thus, all
$B_{A,\eps}^{\ev}$ are geometric Dirac operators on
$\gO^{\ev}(M,E)$ which are formally self-adjoint with respect to
the $L^2$-structure induced by the fixed reference metric $g$.
The $\eps$ factors occur since each
horizontal Clifford variable is scaled with $\eps$.\\

\noindent\textbf{Adiabatic Limit of the Eta Invariant.} The
family of operators $B_{A,\eps}^{\ev}$ converges pointwise to
$\tau_M D_{A,v}$, which is not an elliptic operator. Therefore,
the following result is remarkable, see also \cite[Prop.
4.3]{BC89}.

\begin{prop}\label{AdiabaticLimitExist}
Let $g_\eps$ be an adiabatic metric on the total space of a fiber
bundle $F\hookrightarrow M\xrightarrow{\pi} B$, and assume that
$m=\dim M$ is odd. Let $A$ be a flat $\U(k)$-connection, and let
$\eta(B_{A,\eps}^{\ev})$ be the family of Eta invariants
associated to the adiabatic family of odd signature operators
$B_{A,\eps}^{\ev}$. Then the \emph{``adiabatic limit of the Eta
invariant''} exists in $\R$. More precisely,
\[
\lim_{\eps\to 0} \eta(B_{A,\eps}^{\ev}) = \eta(B_A^{\ev}) +
k\cdot \int_M TL(\nabla^g,\nabla^\oplus),
\]
where $TL(\nabla^g,\nabla^\oplus)$ is the transgression form of
the $L$-class with respect to the connection $\nabla^g$ and
$\nabla^\oplus$ on $TM$.
\end{prop}
\begin{proof}
Fix $\eps\in (0,1)$. Then Proposition \ref{EtaMetricVar} shows
that
\begin{equation}\label{AdiabaticLimitExist1}
\eta(B_{A,\eps}^{\ev}) = \eta(B_A^{\ev}) + k\cdot \int_M
TL(\nabla^g,\nabla^{g_\eps}).
\end{equation}
Moreover, we deduce from Proposition \ref{CSProp} that
\[
\int_M TL(\nabla^g,\nabla^{g_\eps}) = \int_M
TL(\nabla^g,\nabla^\oplus) + \int_M
TL(\nabla^\oplus,\nabla^{g_\eps}).
\]
As in Definition \ref{SThetaDef} consider
\[
S= \nabla^g - \nabla^\oplus\quad\text{and}\quad S^\eps =
\nabla^{g_\eps}-\nabla^\oplus.
\]
It is straightforward to check that
\[
P^v S^\eps = P^vS\quad\text{and}\quad P^hS^\eps = \eps^2 P^hS,
\]
where $P^{v/h}:TM\to T^{v/h}M$ is the vertical, respectively
horizontal, projection of the fiber bundle. Hence,
\[
\lim_{\eps\to 0} \nabla^{g_\eps} = \nabla^\oplus + P^vS,
\]
and the limit is uniform in $\eps$. Therefore,
\[
\lim_{\eps\to 0} \int_M TL(\nabla^\oplus,\nabla^{g_\eps}) = \int_M
TL(\nabla^\oplus,\nabla^\oplus + P^vS).
\]
Now, a consequence of Proposition \ref{PropFiber} is that for
fixed $X\in C^\infty(M,TM)$, the only non-trivial component of
$P^vS(X)$ is
\[
P^vS(X): T^h M \to T^v M.
\]
In particular, $P^vS(X)$ and all its powers are trace-free which
implies that
\[
TL(\nabla^\oplus,\nabla^\oplus + P^vS)=0.
\]
Hence, we can take the limit in \eqref{AdiabaticLimitExist1} and
get
\[
\lim_{\eps\to 0} \eta(B_{A,\eps}^{\ev}) = \eta(B_A^{\ev}) +
k\cdot \int_M TL(\nabla^g,\nabla^\oplus).\qedhere
\]
\end{proof}

\begin{remark*}
So far Proposition \ref{AdiabaticLimitExist} is not of particular
value for explicit computations of $\eta(B_A^{\ev})$. First of
all, we do not yet know anything about the adiabatic limit
$\lim_{\eps\to 0} \eta(B_{A,\eps}^{\ev})$. However, in Chapter
\ref{ChapAbst} we will describe how powerful methods of local
families index theory give an alternative expression of the
adiabatic limit in more topological terms. Another aspect worth
mentioning is that the Chern-Simons term $\int_M
TL(\nabla^g,\nabla^\oplus)$ can be very difficult to compute, see
e.g. \cite{Nic99} for very explicit computations in the case of
circle bundles over surfaces.
\end{remark*}

Concerning Rho invariants we already know at this point that the
transgression term does not play a role. This is because the
transgression term is the same for $B_A^{\ev}$ and the untwisted
odd signature operator $B^{\ev}$, since $A$ is flat. Moreover,
according to Proposition \ref{RhoProp} the Rho invariant is
independent of the metric, so that we have the freedom of
choosing particular well-suited vertical and horizontal metrics.
Summarizing these observations, we obtain the following result,
which is the underlying idea for our discussion of Rho invariants
of fiber bundles.

\begin{cor}\label{RhoAdiabatic}
With respect to all adiabatic metrics $g_\eps$ on $M$ we have
\[
\rho_A(M) = \lim_{\eps\to 0} \eta(B_{A,\eps}^{\ev}) - k\cdot
\lim_{\eps\to 0} \eta(B_{\eps}^{\ev})
\]
\end{cor}

\section[The U$(1)$-Rho Invariant for $S^1$-Bundles over Surfaces]{
The $\boldsymbol{\U(1)}$-Rho Invariant for Principal
$\boldsymbol{S^1}$-Bundles over Surfaces}\label{S1Bundles}

In this section we will see how the idea of Corollary
\ref{RhoAdiabatic} is already helpful for explicit computations,
even without employing more abstract theory we will encounter in
Chapter \ref{ChapAbst}. We give an elementary computation of Rho
invariants for the simple but already non-trivial example of a
principal $S^1$-bundle over a closed surface. We content
ourselves with the $\U(1)$-Rho invariant since all phenomena
related to adiabatic limits appear. Some parts of our discussion
are borrowed from \cite{Nic99}. The setup there is the Spin$^c$
Dirac operator which is, however, closely related to a twisted
odd signature operator.

Before we can start with the discussion of the odd signature
operator on a principal $S^1$-bundle over a Riemannian surface,
we need an explicit description of flat $\U(1)$-connections.

\subsection[The U$(1)$-Moduli Space]{The
$\boldsymbol{\U(1)}$-Moduli Space}\label{S1BundlesModSpace}

Let $\gS$ be a closed, oriented surface of genus $g$, and let
$S^1\hookrightarrow M\xrightarrow{\pi}\gS$ be an oriented
principal circle bundle. Since $H^2(\gS,\Z)=\Z$, such a bundle is
classified up to isomorphism by its degree $l\in \Z$. Given that,
there is a very explicit construction, which we describe now. \\

\noindent\textbf{Topological Description.} Let $\D\subset \gS$ be
an embedded disc, and let $\gS_0:=\gS\setminus \D$. Clearly,
$H^2(\D,\Z)=\{0\}$, and the long exact cohomology sequence of the
pair $(\gS_0,\pd\gS_0)$ implies that $H^2(\gS_0,\Z)=\{0\}$ as
well. Since principal $S^1$-bundles are classified by their first
Chern class, this shows that the restriction of
$S^1\hookrightarrow M\xrightarrow{\pi}\gS$ to $\D$ and $\gS_0$ is
trivializable. Fixing an identification $\pd\D=-\pd\gS_0=S^1$ as
oriented manifolds, we conclude that---up to isomorphism---the
bundle $\pi:M\to \gS$ is given by a glueing function of the form
\begin{equation}\label{CircleGlue}
\gf:\pd(\D\times S^1)\to \pd(\gS_0\times S^1),\quad \gf(z,\gl) =
(z,z^{-l}\gl),
\end{equation}
where $z\in S^1=\pd\D=-\pd\gS_0$. We want to use this description
to determine the fundamental group of $M$. For elements $a,b\in
\pi_1(\gS)$ we write $[a,b]=b^{-1}a^{-1}ba$, which according to
our convention means to first follow the path $a$, then $b$ and
then the same again with the orientations reversed.

\begin{lemma}\label{CircleBundleFund}
Let $S^1\hookrightarrow M\xrightarrow{\pi}\gS$ be an oriented
principal circle bundle of degree $l\in \Z$. Then the fundamental
group of $M$ has the presentation
\[
\pi_1(M)=\Big\langle a_1,b_1,\ldots,a_g,b_g,\gamma\Big |
\prod_{i=1}^g [a_i,b_i]=\gamma^l,\;\text{\rm $\gamma$
central}\Big\rangle,
\]
where $a_1,b_1,\ldots,a_g,b_g$ are lifts to $M$ of the standard
generators of $\pi_1\gS$ and $\gamma$ is the homotopy class of the
$S^1$-fiber.
\end{lemma}

\begin{proof}
Let $c$ be the homotopy class of $\pd\gS_0$. Then the canonical
generators of $\pi_1(\gS_0)$ are the ones indicated in Figure
\ref{Fig:Surf}. It is well known, see e.g. \cite[Sec.
III.3.5]{FarKra}, that
\begin{figure}[htbp]
\centering
\includegraphics[width=0.6\linewidth]{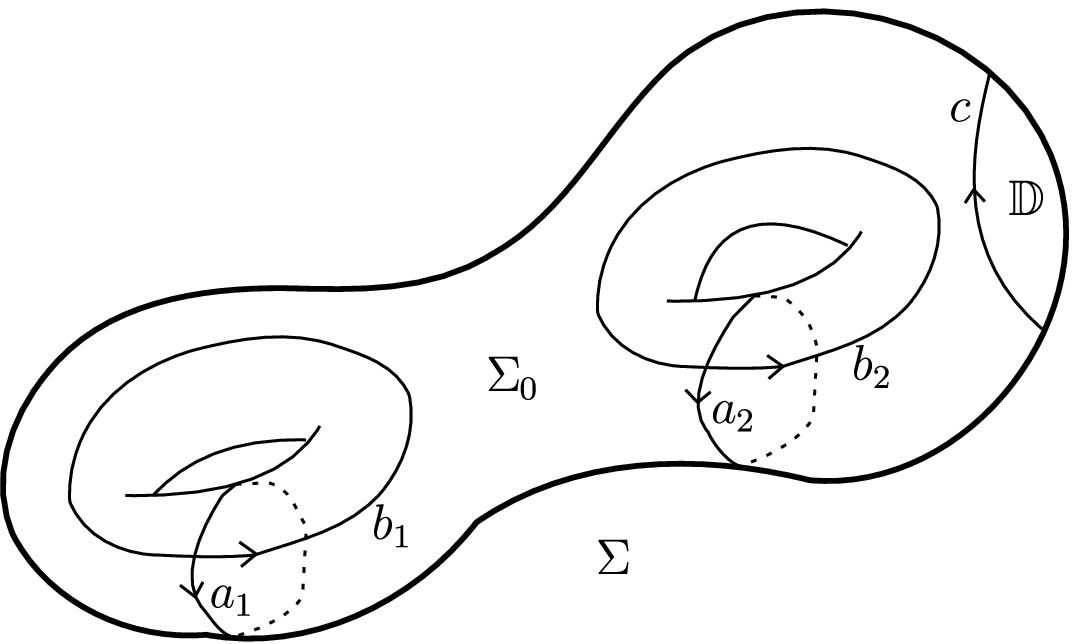}
\caption{Generators of $\pi_1(\gS_0)$}\label{Fig:Surf}
\end{figure}
\[
\pi_1(\gS_0) = \Big\langle a_1,b_1,\ldots,a_g,b_g,c\,\Big |
\prod_{i=1}^g [a_i,b_i]=c^{-1} \Big\rangle.
\]
Write
\[
\pi_1\big(\pd(\D\times S^1)\big) = \big\langle \tilde
c,\tilde\gamma\,\big| [\tilde\gamma,\tilde c] =1\big\rangle,\quad
\pi_1\big(\pd(\gS_0\times S^1)\big) =\big\langle c,\gamma\,\big|
[\gamma,c]=1 \big\rangle.
\]
Note that $\tilde c$ is annihilated under the inclusion
$\pd\D\hookrightarrow \D$. Moreover, the map \eqref{CircleGlue}
induces a map on fundamental groups
\[
\gf_*: \pi_1\big(\pd(\D\times S^1)\big)\to
\pi_1\big(\pd(\gS_0\times S^1)\big),\quad \gf_*(\tilde\gamma)=
\gamma,\quad \gf_*(\tilde c)= \gamma^{-l}c^{-1}.
\]
Van Kampen's Theorem now shows that
\[
\pi_1(M)=\Big\langle a_1,b_1,\ldots,a_g,b_g,c,\gamma\,\Big |
\prod_{i=1}^g [a_i,b_i]=c^{-1},\; \forall_i\;
[a_i,\gamma]=[b_i,\gamma]=1, c^{-1}=\gamma^l \Big\rangle
\]
which by cancelling $c$ coincides with the claimed presentation.
\end{proof}

Since $H_1=\pi_1/[\pi_1,\pi_1]$, it follows immediately from the
above Lemma that
\begin{equation}\label{H1CircleBundle}
H_1(M,\Z)= H_1(\gS,\Z)\oplus \Z_l,
\end{equation}
where we set $\Z_l=\Z$ if $l=0$. As the first homology group
$H_1(\gS,\Z)$ is equal to $\Z^{2g}$, we deduce that
\[
\Hom\big(H_1(\gS,\Z),\U(1)\big) = \U(1)^{2g}.
\]
The long exact coefficient sequence shows that this
$2g$-dimensional torus can be identified with
$H^1(\gS,\R)/H^1(\gS,\Z)$. From \eqref{H1CircleBundle} one can
now determine the moduli space of $\U(1)$-representations, which
is the topological version of moduli space of flat Hermitian line
bundles, see Proposition \ref{FlatModuliRepVar}.

\begin{lemma}\label{ModuliAlgebraic}
Let $M\to\gS$ be an oriented principal circle bundle of degree
$l$. Then the moduli space of flat Hermitian line bundles on $M$
is given by
\[
\cM\big(M,\U(1)\big) \cong \begin{cases}\U(1)^{2g}\times \Z_l,
&\textup{if }l\neq 0,\\ \U(1)^{2g+1}, &\textup{if }l=0.\end{cases}
\]
\end{lemma}

\begin{remark}
Note that in the case $l\neq 0$ it follows from Poincar\'e
duality $H_1(M,\Z) \cong H^2(M,\Z)$ and \eqref{H1CircleBundle}
that $\Tor H^2(M,\Z)=\Z_l$. Hence there are flat line bundles
which are topologically non-trivial. This corresponds to the above
decomposition of $\cM\big(M,\U(1)\big)$ into $l$ different
components. We have included some more details on flat line
bundles which are topologically non-trivial in Appendix
\ref{FlatConn}, see in particular Lemma \ref{TorsionChernTop}.
\end{remark}

\noindent\textbf{Flat Line Bundles over $\boldsymbol{M}$.} We now
need a description of $\cM\big(M,\U(1)\big)$ in terms of flat line
bundles. We will see that every flat line bundle over the total
space $M$ arises as the pullback of a line bundle on the base
$\gS$. Since $M\to \gS$ is a principal $S^1$-bundle, there exists
an associated Hermitian line bundle $L\to \gS$ which clearly has
to play a particular role. Much of the discussion to follow is
inspired by \cite[Sec. 3.3]{Nic98}, although we include some more
details and put more emphasis on the explicit description of the
$\U(1)$-moduli space.

We will work with respect to a fixed connection on the principal
$S^1$-bundle. For this we first endow $\gS$ with a Riemannian
metric $g_\gS$ of unit volume. As noted in Appendix
\ref{HolomAspects}, this amounts to fixing a complex structure on
$\gS$. We identify the Lie algebra of $S^1$ with $i\R$. Let $e$
be the vector field on $M$ associated to the $S^1$-action,
\[
e|_p = \lfrac d{dt}\big|_{t=0} p\cdot e^{it},\quad p\in M.
\]
A connection on the principal bundle $\pi:M\to \gS$ is a 1-form
$i\go\in \gO^1(M,i\R)$ such that
\[
\go(e) = 1\quad\text{and}\quad R_{e^{it}}^*\go = \go,
\]
where $R_{e^{it}}$ denotes right-multiplication, compare with
\eqref{ConnDef}. Let $F_\go\in \gO^2(\gS,i\R)$ be the curvature of
$i\go$. Since the cohomology class of $\lfrac i{2\pi}F_\go$
represents the rational first Chern class of the bundle $\pi:M\to
\gS$ we can choose $\go$ in such a way that
\begin{equation}\label{CurvDeg}
-\lfrac1{2\pi}d\go = \lfrac i{2\pi}\pi^*F_\go =  l\cdot
\pi^*\vol_\gS.
\end{equation}

Let $L\to \gS$ be the line bundle associated to the principal
bundle structure. The connection $\go$ induces a natural
connection $A_\go$ on $L$. We write $L_\go$ for the line bundle
$L$ endowed with this particular connection $A_\go$. As explained
in Appendix \ref{HolomAspects}, this is the same as fixing a
holomorphic structure on $L$. The following simple observation
relies only on the principal bundle structure and not on the
particular structure group $\U(1)$ or the dimension of the base.

\begin{lemma}\label{TrivialPullBack}
The pullback $\pi^*L_\go\to M$ is canonically trivial and the
pullback connection $\pi^*A_\go$ satisfies
\[
\pi^*A_\go = d_M+i\go.
\]
\end{lemma}

\begin{proof}
Recall that the associated bundle $L_\go\to \gS$ is defined by
the pullback diagram
\[
\begin{CD}
M\times \C @>>> M\\
@VVV @V{\pi}VV \\
M\times \C/\sim @>>> \gS
\end{CD}
\]
where $(p,v)\sim (pz,z^{-1}v)$ for all $z\in S^1$. Since
$\pi^*L_\go$ is given by the same pullback diagram, we
tautologically get $M\times \C = \pi^*L_\go$. Under this
identification,
\[
\pi^* C^\infty(\gS,L_\go) = \bigsetdef{\gf:M\to \C}{ \gf(pz) =
z^{-1} \gf(p)}.
\]
The pullback connection $\pi^*A_\go$ acts on equivariant
functions $\gf$ as
\[
(\pi^*A_\go)_X \gf = X^h\gf = d_M\gf(X) - X^v\gf,
\]
where $X^{h/v}$ denotes the horizontal/vertical projection of $X$
with respect to the connection $i\go$. The latter is explicitly
given by
\[
X^v_p = \lfrac{d}{dt}\big|_0
p\cdot\exp\big(ti\go(X_p)\big)=\go(X)\,e|_p.
\]
If $\gf$ is an equivariant function, then $e\gf  = -i\gf$.
Therefore,
\[
X^v_p\gf =  -i\go(X_p)\gf(p).
\]
Extending by the Leibniz rule to all functions on $M$, we get
\[
\pi^*A_\go = d_M + i\go.\qedhere
\]
\end{proof}

Now let $L_A\to \gS$ be an arbitrary Hermitian line bundle of
degree $k$ with a holomorphic structure given by a unitary
connection $A$, see Appendix \ref{HolomAspects}. It follows from
Proposition \ref{PicModRel} that upon transforming $A$ with a
complex gauge transformation $f\in \cG^c$ we may---and
will---assume in the following that
\begin{equation}\label{IntegralConnection}
\lfrac i{2\pi}F_A = k\vol_\gS.
\end{equation}
In general, to achieve this, we really need to transform with a
complexified gauge transformation and not just a unitary one, see
Proposition \ref{PicModRel}.

\begin{lemma}\label{FlatPullBack}
Assume that $l\neq 0$, and let $q:=k/l$. Then the connection
\begin{equation*}
A_q:= \pi^*A - q\, i\go\quad\text{on $\pi^*L_A$}
\end{equation*}
is flat. Moreover, the holonomy of $A_q$ along the $S^1$-fiber
$\gamma$ is given by
\[
\hol_{A_q}(\gamma) = \exp(2\pi i q).
\]
\end{lemma}
\begin{proof}
By functoriality, we have $F_{\pi^*A} = \pi^*F_A$. Thus, it
follows from assumptions \eqref{CurvDeg} and
\eqref{IntegralConnection} that the curvature of $A_q$ satisfies
\[
F_{A_q} = \pi^*F_A - q\,i d \go = -2\pi i\big(k- ql
\big)\pi^*\vol_\gS = 0.
\]
To compute the holonomy, let $p\in M$ be arbitrary and let
$\gamma(t) = p\cdot\exp(it)$ with $t\in[0,2\pi]$ parametrize the
fiber containing $p$. Clearly, $\hol_{\pi^*A}(\gamma) = 0$ and
$\go_{\gamma(t)}(\overset{.}{\gamma}(t)) =1$. Therefore,
\[
\hol_{A_q}(\gamma) = \exp\Big(-\int_\gamma - q\,i \go \Big) =
\exp\Big(q\,i  \int_0^{2\pi}
\go_{\gamma(t)}(\overset{.}{\gamma}(t))dt \Big)= \exp(2\pi i
q).\qedhere
\]
\end{proof}

\noindent\textbf{The Moduli Space of Flat Line Bundles.} After
this preparation, we can now give the geometric description of
$\cM\big(M,\U(1)\big)$. Recall that $\Pic(\gS)$ denotes the Picard
group of holomorphic line bundles over $\gS$, see Definition
\ref{PicDef}.

\begin{prop}\label{CircleBundleModuli}
Let $l\neq 0$ and let $\cM\big(M,\U(1)\big)$ be the moduli space
of flat line bundles over $M$. Then $\go$ induces a natural
surjection
\[
\pi^*: \Pic(\gS)\to \cM\big(M,\U(1)\big),\quad [L_A]\mapsto
\big[\pi^*{L_A,A_q}\big],
\]
where $A_q$ is defined as in Lemma \ref{FlatPullBack}. There is a
natural $\Z$-action on $\Pic(\gS)$ given by $(L_A,k)\mapsto
L_A\otimes L_\go^{\otimes k}$, and with respect to this,
\[
\Pic(\gS)/\Z \cong \cM\big(M,\U(1)\big).
\]
\end{prop}

\begin{proof}
Let $L_A\to \gS$ be a holomorphic line bundle. Assume that $B$ is
another unitary connection on $L$, satisfying condition
\eqref{IntegralConnection} and inducing an equivalent holomorphic
structure, i.e.,
\[
B = A + u^{-1}du,\quad \text{for some $u\in \cG^c$.}
\]
Condition \eqref{IntegralConnection} means in particular that
$F_A=F_B$, which implies that in fact $u\in \cG$. From this we
obtain that
\[
B_q = A_q +\pi^*(u^{-1}du),\quad \text{for some $u\in \cG$.}
\]
i.e., the flat connections $B_q$ and $A_q$ on $\pi^*L_A$ are
equivalent. This shows that the map in Proposition
\ref{CircleBundleModuli} is well-defined.

To verify that it is surjective, let $L\to M$ be a flat Hermitian
line bundle with connection $\tilde A$. Let $\gamma$ denote the
generator of the $S^1$-fiber in $\pi_1(M)$. Then Lemma
\ref{CircleBundleFund} shows that $\gamma^l$ is a commutator. It
follows that for some $k\in \Z$,
\begin{equation}\label{CircleBundleModuli:1}
\hol_{\tilde A}(\gamma) = \exp(2\pi i k/l).
\end{equation}
Now let $L_A\to\gS$ be an arbitrary holomorphic line bundle of
degree $k$. We infer from \eqref{CircleBundleModuli:1} and Lemma
\ref{FlatPullBack} that $\pi^*L_A\otimes L^{-1}$, endowed with
the connection $A_q\otimes 1 - 1\otimes \tilde A$, is a flat line
bundle on $M$ with trivial holonomy along the fiber $\gamma$.
This easily implies that it is equivalent to the pullback
$\pi^*\C_B$ of the trivial line bundle over $\gS$ endowed with a
flat connection $B$. Thus, as line bundles with connection,
\[
L = \pi^*(L_A\otimes \C_B),
\]
which proves surjectivity of the map in Proposition
\ref{CircleBundleModuli}.

As we have seen in Lemma \ref{TrivialPullBack}, the pullback
$\pi^*L_\go$ with connection $\pi^*A_\go -i\go$ is the trivial
flat line bundle. Using this one observes that the map $\pi^*$ is
invariant under the natural $\Z$-action on $\Pic(\gS)$. Assume
now that $\pi^*L_A = \pi^*L_B$ for two holomorphic line bundles
over $\gS$ of degree $k$ and $m$ respectively. Since their
holonomies along $\gamma$ agree, it follows that $m-k=nl$ for
some $n\in\Z$. Thus,
\[
\pi^*(L_B\otimes L_A^{-1}) = \pi^*L_\go^{\otimes n}\quad
\text{and}\quad \pi^*B=\pi^*A + n\cdot i\go.
\]
We deduce that $L_B  = L_A\otimes L_\go^{\otimes n}$ as
holomorphic line bundles, which is what we needed to prove.
\end{proof}

\begin{remark}\label{Gysin}
Recall that the $S^1$-bundle $\pi:M\to \gS$ gives rise to the
Gysin sequence
\[
...\rightarrow H^0(\gS)  \xrightarrow{\cup c} H^2(\gS)
\xrightarrow{\pi^*} H^2(M) \xrightarrow{\pi_*} H^1(\gS)
\rightarrow 0,
\]
see \cite[Prop. 14.33]{BT}. Here, $c=c(M)\in H^2(\gS)$ is the
first Chern class (or Euler class) of the oriented $S^1$-bundle.
If we are assuming that $l\neq 0$, the map $H^0(\gS)
\xrightarrow{\cup c} H^2(\gS)$ gives an isomorphism in de Rham
cohomology. This implies that for cohomology with integer
coefficients, the map $\pi^*:H^2(\gS,\Z)\to H^2(M,\Z)$ appearing
in the Gysin sequence surjects onto the torsion subgroup of
$H^2(M,\Z)$. It is related to the map $\pi^*$ of Proposition
\ref{CircleBundleModuli} by the following diagram
\[
\begin{CD}
\Pic(\gS) @>>{c_1}> H^2(M,\Z)\\
@V{\pi^*}VV @V{\pi^*}VV \\
\cM\big(M,\U(1)\big)@>>{c_1}> \Tor\big(H^2(M,\Z)\big).
\end{CD}
\]
Note also that the first Chern class $c_1$ is equivariant with
respect to the $\Z$-action on $\Pic(\gS)$,
\[
c_1(L_A\otimes L_\go^{\otimes k}) = c_1(L_A) + k\cdot c(M).
\]
Using Proposition \ref{PicModRel} one can now interpret the above
diagram as the geometric version of Lemma \ref{ModuliAlgebraic}
in the case $l\neq 0$.
\end{remark}

The structure result Proposition \ref{CircleBundleModuli} excludes
the case that the circle bundle is of degree $l=0$, i.e.,
isomorphic to $\gS\times S^1$. However, a geometric description
in this case is easy to find directly. As in Remark \ref{RhoRem}
(iii), a flat line bundle $L_q$ over $S^1$ is the trivial line
bundle endowed with the connection $d-qz^{-1}dz$ for some
$q\in\R$. Here, we view $S^1$ as a subset of $\C$, and $z^{-1}dz$
expresses the Maurer-Cartan form of $S^1$. Clearly, $L_q$ and
$L_{q'}$ are unitarily equivalent if and only if $q-q'\in 2\pi
i\Z$. Without effort one verifies the following result.

\begin{lemma}\label{ModuliTrivialBundle}
If $M=\gS\times S^1$ is the trivial circle bundle over $\gS$, then
\[
\cM\big(\gS,\U(1)\big)\times \cM\big(S^1,\U(1)\big)\cong
\cM\big(M,\U(1)\big).
\]
Here, the isomorphism is given by
\[
\big([L_A],[L_q]\big)\mapsto \big[L_A\boxtimes L_q\big],
\]
where $L_A\boxtimes L_q$ is the fiber product defined in
\eqref{FiberProd}, endowed with its natural connection.
\end{lemma}

\subsection{The Odd Signature Operator.} We now want to identify
the odd signature operator on the total space of a principal
circle bundle over a closed, oriented surface. Certainly, the
underlying principal bundle structure will play an important
role, and many features generalize to arbitrary principal bundles
with compact structure group. However, we will not give many
comments about these
generalizations. \\

\noindent\textbf{Fibered Calculus on $\boldsymbol{M}$.} To start,
we need to identify some of the quantities defined in Section
\ref{FiberedCalc} in the case at hand. Let $i\go$ be a connection
on the principal $S^1$-bundle $\pi:M\to \gS$. Since the vector
field $e$ associated to the $S^1$-action gives a trivialization
of the vertical tangent bundle, we get a vertical projection
\[
P^v:TM\to T^vM,\quad X\mapsto \go(X)e.
\]
With respect to this, the curvature $\gO$ in the sense of
Definition \ref{CurvFiber} is related to the curvature of $i\go$
by
\[
\gO(X^h,Y^h)=-\go\big([X^h,Y^h]\big)e = d\go(X^h,Y^h)e =
-iF_\go(X,Y) e,
\]
where $X^h$ and $Y^h$ are horizontal lifts of vector fields $X,Y$
on $\gS$. In particular, when we fix a metric $g_\gS$ of unit
volume and require that $\go$ satisfies \eqref{CurvDeg}, we have
\begin{equation}\label{CurvFiberCirc}
\gO(X^h,Y^h) = -2\pi l \vol_\gS(X,Y) e.
\end{equation}
We now endow $T^vM$ with the vertical metric
$g_v:=\go\otimes\go$, and consider the submersion metric
$g=g_\gS\oplus g_v$.

\begin{lemma}\label{CircleBundleCalc}\quad\nopagebreak
\begin{enumerate}
\item With respect to the trivialization given by $e$, the
canonical connection $\nabla^v$ on $T^vM$ is the trivial
connection, i.e.,
\[
\nabla^v_Xe = 0,\quad X\in C^\infty(M,TM).
\]
\item If $X\in C^\infty(\gS,T\gS)$, we have
\[
\sL_{X^h}(e) = [X^h,e] = 0,\quad\text{and}\quad \sL^v_{X^h}(g_v)
=0.
\]
In particular, the connection $\widetilde\nabla^v$ from Definition
\ref{CanVerticalConn2} agrees with $\nabla^v$, and the mean
curvature $k_v$ as well as the tensor $B$ in \eqref{NablaNatInv}
vanish.
\end{enumerate}
\end{lemma}

\begin{proof}
The connection $\nabla^v$ is compatible with $g_v$. Hence,
\[
0= g_v(\nabla_X^ve,e) + g_v(e,\nabla_X^ve),
\]
which yields $\nabla_X^ve =0$. This proves (i). Since $i\go$ is a
connection, we have $R_{e^{it}}^*\go = \go$. This implies that
the metric $g_v=\go\otimes \go$ is invariant under the flow
associated to the vector field $e$. Since $\pi^*g_\gS$ is
constant along the fiber, we find that for all vector fields $X$
on $\gS$
\[
0 = \sL_e(g)(X^h,e) = g\big([e,X^h],e\big) + g\big(X^h,[e,e]\big).
\]
As $[e,e]=0$ we conclude that $g\big([e,X^h],e\big)=0$. This
implies that $[e,X^h]=0$, because  Lemma
\ref{HorizontalVerticalLie} ensures that the vector field
$[e,X^h]$ is vertical. In particular, since
\[
\sL_{X^h}^vg_v(e,e) = X^h\big(g_v(e,e)\big) -
2g_v\big([X^h,e],e\big),
\]
we deduce that the vertical Lie derivative of $g_v$ vanishes.
Using its very definition, we see that the tensor $B$ is indeed
trivial. Moreover, we know from Lemma \ref{MeanCurvFiber:alt} that
the mean curvature is given by the trace of $\sL_{X^h}^v(g_v)$.
Thus, it is also is zero. Moreover, using part (i) we find that
the derivations $\sL_{X^h}$ and $\nabla_{X^h}^v$ agree on $e$.
Since both satisfy the Leibniz rule they are necessarily equal.
Hence, by definition, the connection $\widetilde\nabla^v$ agrees
with $\nabla^v$.
\end{proof}

\noindent\textbf{Rho Invariants for Trivial Circle Bundles.}
Before we continue with the general discussion, we assume that
$l=0$ so that $M=\gS\times S^1$. We endow $M$ with the natural
connection $i\go = z^{-1}dz$ given by the Maurer-Cartan form on
$S^1$. Choose a flat line bundle $L\to M$, i.e.,
\[
L = L_A\boxtimes L_q\to \gS\times S^1,
\]
where $L_A$ and $L_q$ are flat line bundles over $\gS$
respectively $S^1$, see Lemma \ref{ModuliTrivialBundle}. We
identify
\[
\gO^{\ev}(\gS\times S^1,L) = \gO^\bullet(\gS,L_A)\otimes
C^\infty(S^1).
\]
Using Lemma \ref{OddSignSplit} and Lemma \ref{TauSplit}, we can
write the odd signature operator as
\[
B_{A,q}:=B^{\ev}_{A,q} = \tau_\gS\otimes B_q+ D_A\otimes 1,
\]
where $D_A$ is the twisted de Rham operator on $\gS$ and $B_q$ is
the odd signature operator on $S^1$,
\begin{equation*}
B_q=-i(\sL_e - iq) :C^\infty(S^1)\to C^\infty(S^1).
\end{equation*}
Hence, $B_{A,q}$ is of the form considered in Lemma \ref{EtaProp}
(iii). According to the Hirzebruch Signature Theorem, the index
of $D_A^+$ vanishes for all flat connections $A$ on $\gS$ and so
\begin{equation}\label{EtaTrivCirc}
\eta(B_{A,q}) = \ind(D_A^+)\cdot\eta(B_q)=0.
\end{equation}
Therefore, all Rho invariants for the trivial circle bundle
$\gS\times S^1$ vanish.\\

\noindent\textbf{The Structure of $\boldsymbol{B_{A,q}}$ in the
General Case.} We now assume that $l\neq 0$. Let $L_A\to\gS$ be a
line bundle of degree $k$ endowed with a Hermitian connection $A$
which satisfies the condition of \eqref{IntegralConnection},
\[
\lfrac i{2\pi} F_A = k\cdot\vol_\gS.
\]
We endow the pullback $L:=\pi^*L_A\to M$ with the flat connection
of Lemma \ref{FlatPullBack}, i.e.,
\[
A_q= \pi^*A - iq\, \go, \quad q:=k/l.
\]
Since $L$ is the pullback of $L_A$, we alter the identification
\eqref{FormsSplit} slightly to
\[
\gO^\bullet(M,L) = \pi^*\gO^\bullet(\gS,L_A)\otimes
\gO^\bullet_v(M).
\]
As in Lemma \ref{dSplit} we write the twisted de Rham operator as
\[
d_{A_q} = d_{q,v} + d_{A,h} + \imu(\gO),
\]
where
\begin{equation}\label{DiffCircleBundle}
d_{q,v} = d_v - iq\, \emu(\go),\quad d_{A,h} =
\emu(f^a)\widetilde \nabla^{A_q,\oplus}_a =  (\pi^*d_A )\otimes 1
+\emu(f^a)\otimes \widetilde \nabla_a^v.
\end{equation}
As always $\emu(.)$ denotes exterior multiplication and
$\{f_1,f_2\}$ is a local orthonormal frame for $T\gS$.

To describe the odd signature operator, we split the space of
$L$-valued differential forms of even degree as in Lemma
\ref{OddSignSplit},
\begin{equation*}
\gO^{\ev}(M,L)= \pi^*\gO^\bullet(\gS,L_A)\otimes C^\infty(M).
\end{equation*}

\begin{prop}\label{OddSignCirc}
With respect to the above identification, the odd signature
operator is given by
\[
B_{A,q}= \tau_\gS\otimes B_{q,v} + D_{A,h} + \tau_MT,
\]
where the individual terms are
\[
B_{q,v} = -i(\sL_e -iq), \quad D_{A,h} = D_A\otimes 1 +
c(f^a)\otimes \sL_{f_a}^v,
\]
and
\[
\tau_M T =\begin{cases} \quad 0 &\text{\rm on } \gO^{0,0}\oplus
\gO^{1,0},\\ -2\pi l &\text{\rm on } \gO^{2,0}.
\end{cases}
\]
Moreover, we have the (anti-)commutator relations
\begin{equation}\label{DvDh:anticommute}
\big[1\otimes B_{q,v},D_{A,h}\big] = 0\quad \text{and}\quad
\big\{\tau_\gS\otimes B_{q,v},D_{A,h} \big \}=0.
\end{equation}
\end{prop}

\begin{proof}
Let $\ga\in \pi^*\gO^p(\gS,L_A)$ and $\gf\in C^\infty(M)$. Then,
\[
(\tau_M D_{q,v})(\ga\wedge \gf) = (-1)^p \tau_M\big(
\ga\wedge(D_{q,v}\gf)\big) =
(\tau_\gS\ga)\wedge(\tau_vD_{q,v}\gf),
\]
where we have used Lemma \ref{TauSplit} in the last equality.
Now, checking the factors of $i$ in Definition \ref{VertChirDef}
one finds that $\tau_v(\go) = -i$. Thus,
\[
(\tau_vD_{q,v})\gf = (\tau_v d_{q,v})\gf = \tau_v\big(d_v -
iq\go\big) \gf = \tau_v(\go)\big(\sL_e - iq\big)\gf = -i
\big(\sL_e - iq\big)\gf.
\]
According to Lemma \ref{OddSignSplit} the horizontal part of
$B_{A,q}$ coincides with the horizontal de Rham operator
\[
D_{A,h}= d_{A,h} + d_{A,h}^t = \emu(f^a) \widetilde
\nabla^{A,\oplus}_a -\imu(f^a)\big(\widetilde\nabla^{A,\oplus}_a
+ 2B(f_a) + k_v(f_a)\big).
\]
Here, we have used Proposition \ref{d^tSplit} and
\eqref{NablaNat}. Hence, we deduce from Lemma
\ref{CircleBundleCalc} that
\[
D_{A,h} = c(f^a) \widetilde\nabla^{A,\oplus}_a = D_A\otimes 1 +
c(f^a)\otimes \widetilde \nabla_a^v.
\]
For the last term appearing in the formula for $B_{A,q}$ note that
\[
T(\ga\otimes \gf) = \emu(f^1)\emu(f^2)\imu(\gO_{12}) -
\imu(f^1)\imu(f^2)\emu(\gO_{12}),
\]
where \eqref{CurvFiberCirc} shows that $\gO_{12}= -2\pi l e$.
Therefore, $T$ is non-zero only on $\gO^{2,0}$, and
\[
\begin{split}
(\tau_M T)(\vol_\gS\wedge \gf) &= \tau_M \big(-2\pi
l\imu(f^1)\imu(f^2)\vol_\gS\wedge \gf\go\big)\\ &= -2\pi
l\big(\tau_\gS(1)\wedge \gf \tau_v(\go)\big) = -2\pi l
\vol_\gS\wedge \gf.
\end{split}
\]
In this computation we have used Lemma \ref{TauSplit} to express
$\tau_M$ in terms of $\tau_\gS$ and $\tau_v$. Also note that
$\tau_\gS(1) = i\vol_\gS$ and $\tau_v(\go) = -i$.

Lemma \ref{CircleBundleCalc} (ii) implies that the bundle
endomorphism $K$ as defined in \eqref{KDef} vanishes. Thus, we
deduce from Proposition \ref{dvdhCommutator} that
\begin{equation}\label{DvDh:anticommute:1}
d_{q,v}d_{A,h}^t + d_{A,h}^td_{q,v} = 0,\quad\text{and}\quad
D_{q,v}D_{A,h} +  D_{A,h}D_{q,v} =0.
\end{equation}
Also $D_A$ and $c(f^a)$ anti-commute with $\tau_\gS$, since $\gS$
is even dimensional. This yields the relations in
\eqref{DvDh:anticommute}.
\end{proof}

\noindent\textbf{The Spectrum of $\boldsymbol{B_{q,v}}$.} The
vertical odd signature operator $B_{q,v}$ is not elliptic, since
its principal symbol vanishes in all directions orthogonal to the
fiber. Thus, we do not know much about the spectrum of $B_{q,v}$
by employing the general theory. However, due to the
$S^1$-symmetry, we can determine its eigenvalues by hand.

\begin{remark*}
Before we state the next result, recall that $L_\go\to\gS$ denotes
the line bundle associated to $M$ endowed with the connection
$A_\go$ induced by $\go$. As we have seen in Lemma
\ref{TrivialPullBack}, the pullback $\pi^*L_\go\to M$ is
isomorphic to the trivial line bundle. Under this identification,
a function $\gf\in C^\infty(M)$ is the pullback of a section
$s_\gf \in C^\infty(\gS,L_\go)$ if and only if
\begin{equation}\label{PullBackSection}
\gf(p\cdot z) = z^{-1}\gf(p),\quad p\in M,\quad z\in S^1.
\end{equation}
Moreover, the Lie derivative is related to the connection $A_\go$
via
\begin{equation}\label{PullBackDer}
s_{(X^h\gf)} = A_\go(X)s_\gf.
\end{equation}
We refer to the proof of Lemma \ref{TrivialPullBack} for more
details.
\end{remark*}

\begin{lemma}\label{VerticalEigenspaces}
Assume that $l\neq 0$, and let $L_A$ be a holomorphic line bundle
over $\gS$ of degree $k$. Moreover, let $q:=k/l$, and let
$L=\pi^*L_A$ be the associated flat line bundle over $M$. Then
\[
\ker \big( B_{q,v} -\gl\big) \neq \{0\}\quad\text{if and only
if}\quad \gl+ q\in \Z.
\]
Moreover, if $\gl+ q\in \Z$, then
\[
\ker \big( B_{q,v} -\gl\big) \cong \pi^*C^\infty(\gS,
L_{B_\gl}),\quad \text{where}\quad L_{B_\gl}:=L_A\otimes
L_\go^{-(\gl+q)}.
\]
The operator $D_{A,h}$ restricted to $\gO^\bullet(\gS)\otimes
\ker \big( B_{q,v} -\gl\big)$ corresponds under the above
isomorphism to
\[
D_{B_\gl}: \gO^\bullet(\gS,L_{B_\gl})\to
\gO^\bullet(\gS,L_{B_\gl}),
\]
where $B_\gl = A\otimes 1 + 1\otimes -(\gl+q)A_\go$ is the natural
connection on $L_{B_\gl}$.
\end{lemma}

\begin{proof}
Let $\gf\in\ker (B_{q,v}-\gl)$. According to Proposition
\ref{OddSignCirc} this means that
\[
\sL_e\gf = i\big(q+ \gl\big)\gf .
\]
For $t\in \R$ and $p\in M$ let $\widehat\gf_t(p):= \gf(p\cdot
e^{it})$. Then, since $e$ is the vector field generated by the
$S^1$-action,
\[
\lfrac d{dt}  \widehat\gf_t = \sL_e(\widehat\gf_t) =
i(q+\gl)\cdot\widehat\gf_t,\quad\text{i.e.,}\quad \widehat\gf_t =
e^{i(q+\gl)t}\gf.
\]
This implies that $q+\gl\in\Z$ or $\gf =0$. Rewriting the result
in terms of $z=e^{it}$ we see that
\[
\gf(p\cdot z) = z^{q+\gl}\cdot \gf(p),\quad z\in S^1.
\]
As in \eqref{PullBackSection} this means that we can identify
$\gf$ with a section
\[
\gf \in \pi^*C^\infty\big(\gS,L_A\otimes
L_\go^{-q-\gl}\big)=\pi^*C^\infty(\gS,L_{B_\gl}).
\]
Tracing the proof backwards shows that conversely every such
element gives an eigenvector of $B_{q,v}$. The assertion about
$D_{A,h}$ easily follows from \eqref{PullBackDer} and Proposition
\ref{OddSignCirc}.
\end{proof}

\begin{remark}\label{VerticalEigenspacesDecomp}
Without going into details, we want to mention that we have
actually determined the full spectrum of $B_{q,v}$. We note
without proof that $B_{q,v}$ is essentially self-adjoint in
$\pi^*\gO^\bullet(\gS,L_A)\otimes C^\infty(M)$ and that
\eqref{DvDh:anticommute} implies that it commutes with the
formally self-adjoint elliptic operator $D_h+B_{q,v}$. This
suffices to guarantee that $\spec(B_{q,v})$ consists only of
eigenvalues---though, with infinite multiplicities. Then Lemma
\ref{VerticalEigenspaces} implies that
\[
\spec(B_{q,v}) = \setdef{\gl}{\gl+ q\in \Z}.
\]
Moreover, as in the case of eigenvalues with finite multiplicity,
we can decompose
\[
\pi^*\gO^\bullet(\gS,L_A)\otimes C^\infty(M) =
\bigoplus_{\gl\in\spec(B_{q,v})} \gO^\bullet(\gS)\otimes \ker
\big( B_{q,v} -\gl\big) \cong
\bigoplus_{\gl\in\spec(B_{q,v})}\gO^\bullet(\gS,L_{B_\gl}).
\]
We also want to note that this decomposition is essentially the
decomposition of the infinite dimensional $S^1$-module
$\pi^*\gO^\bullet(\gS,L_A)\otimes C^\infty(M)$ into its
irreducible components. A similar situation should occur for
general Lie groups.
\end{remark}

\subsection{The Eta Invariant of the Truncated Odd Signature
Operator.} The fact that the commutators in
\eqref{DvDh:anticommute} in Proposition \ref{OddSignCirc} are zero
allows us to give an elementary computation of Eta invariants,
see \cite[App. C]{Nic99} for a related treatment. However, the
Eta invariant of the full signature operator is not directly
tractable. Therefore, we introduce the following.

\begin{dfn}
Let $L_A\to\gS$ be a line bundle of degree $k$, and let
$L:=\pi^*L_A$ be the corresponding flat line bundle over $M$. We
call the operator
\[
B_{A,q}^\oplus:= \tau_\gS\otimes
B_{q,v}+D_{A,h}\quad\text{on}\quad
\pi^*\gO^\bullet(\gS,L_A)\otimes C^\infty(M)
\]
the \emph{truncated odd signature operator} twisted by $L$.
\end{dfn}

\begin{remark}
The connection $\nabla^\oplus$ from \eqref{SplitConn} together
with $A_q$ induces a connection $\nabla^{A_q,\oplus}$ on
$\gL^\bullet T^*M\otimes L$. Then the truncated odd signature
operator is given by Clifford contraction of
$\nabla^{A_q,\oplus}$. Therefore, it is almost as good as a
geometric Dirac operator. However, $\nabla^{A_q,\oplus}$ is not a
Clifford connection since it is compatible with $\nabla^\oplus$
and not with the Levi-Civita connection $\nabla^g$. As remarked
earlier an operator of this type is in general not formally
self-adjoint. However, in the situation at hand, $B_{A,q}^\oplus$
is clearly formally self-adjoint, since $B_{q,v}$ and $D_{A,h}$
are.
\end{remark}

Since $B_{A,q}^\oplus$ is an formally self-adjoint elliptic
differential operator on a closed manifold, its Eta function is
well-defined and for $\Re(s)$ large,
\[
\eta(B_{A,q}^\oplus,s) = \frac{1}{\gG\big(\lfrac{s+1}{2}\big)}
\int_0^\infty
\Tr\big[B_{A,q}^\oplus\exp\big(-t(B_{A,q}^\oplus)^2\big)\big]
t^{\frac{s-1}2}dt.
\]
Moreover, Theorem \ref{EtaReg} implies that the meromorphic
extension of $\eta(B_{A,q}^\oplus,s)$ has no pole in 0. Our
strategy is now to compute the Eta invariant of the truncated odd
signature operator explicitly and determine its kernel, see
Proposition \ref{EtaTruncSign} and Proposition
\ref{TruncSignOpKern}. Then in Section \ref{S1BundlesAdiab} we
will use these results to determine the Rho invariant of the full
odd signature operator $B_{A,q}$. For this we want to use the
variation formula of Proposition \ref{EtaDiffSF} for Eta
invariants so that we will need to understand the spectral flow
between $B_{A,q}^\oplus$ and $B_{A,q}$. However, the difference
$B_{A,q} - B_{A,q}^\oplus$ might be ``too large'' compared to
$B_{A,q}^\oplus$ to get a good control of the spectrum near zero.
As we will see the solution to this problem is to study an
adiabatic metric. After this short digression on our strategy let
us now investigate the truncated odd signature operator.

\begin{prop}\label{EtaTruncSign}
Let $S^1\hookrightarrow M\xrightarrow{\pi} \gS$ be a principal
circle bundle of degree $l\neq 0$. Let $L_A\to\gS$ be a line
bundle over $\gS$ of degree $k$, and let $L:=\pi^*L_A$ be the
corresponding flat line bundle over $M$. Then, with $q=k/l$,
\[
\eta(B_{A,q}^\oplus) =  2lP_2(q),
\]
where $P_2$ is the second periodic Bernoulli function, i.e., if
$q-[q]=q_0$ with $q_0\in [0,1)$ and $[q]\in \Z$, then
\[
P_2(q) = q_0^2- q_0+ \lfrac 16,
\]
see Definition \ref{PeriodicBernoulli}. In particular,
$\eta(B_{A,q}^\oplus)$ is independent of the metric $g_\gS$, and
the connection $A$ involved in its definition.
\end{prop}

\begin{proof}
Formula \eqref{DvDh:anticommute} in Proposition \ref{OddSignCirc}
shows that $\tau_\gS\otimes B_{q,v}$ anti-commutes with
$D_{A,h}$. Hence, we can split
\[
B_{A,q}^\oplus e^{-t(B_{A,q}^\oplus)^2} = D_{A,h} e^{-tD_{A,h}^2
-t(\tau_\gS\otimes B_{q,v})^2} + \big(\tau_\gS\otimes
B_{q,v}\big) e^{-tD_{A,h}^2 -t(\tau_\gS\otimes B_{q,v})^2}.
\]
Since $D_{A,h}$ anti-commutes with $\tau_\gS$ one finds as in the
proof of Lemma \ref{EtaProp} that
\[
\Tr\big[D_{A,h} e^{-tD_{A,h}^2 -t(\tau_\gS\otimes
B_{q,v})^2}\big] =0.
\]
Now let $\gl\in\spec(B_{q,v})$. Then according to Lemma
\ref{VerticalEigenspaces}, the operator
\[
\big(\tau_\gS\otimes B_{q,v}\big) e^{-tD_{A,h}^2
-t(\tau_\gS\otimes B_{q,v})^2}\quad \text{on}\quad
\gO^\bullet(\gS)\otimes \ker \big( B_{q,v} -\gl\big)
\]
is unitarily equivalent to
\begin{equation*}\label{EtaTrunc:1}
\big(\tau_\gS e^{-t D_{B_\gl}^2}\big)\gl e^{-t\gl^2}:
\gO^\bullet(\gS,L_{B_\gl})\to \gO^\bullet(\gS,L_{B_\gl}),
\end{equation*}
where $B_\gl = A -(\gl+q)A_\go$. It follows from the McKean-Singer
formula and the index theorem for the signature operator that
\begin{equation*}\label{EtaTrunc:2}
\Tr\big[\tau_\gS e^{ -t(D_{B_\gl})^2}
\big]=\ind\big(D_{B_\gl}^+\big) = 2\big(k -l(\gl+q)\big) = -2l\gl.
\end{equation*}
Hence, using the decomposition from Remark
\ref{VerticalEigenspacesDecomp} we find
\[
\begin{split}
\Tr\big[B_{A,q}^\oplus\exp\big(-t(B_{A,q}^\oplus)^2\big)\big] &=
\sum_{\gl\in\spec(B_{q,v})} \Tr\big[\tau_\gS e^{
-t(D_{B_\gl})^2} \big]\gl e^{-t\gl^2} \\
&= -2l \sum_{\gl\in\spec(B_{q,v})} \gl^2 e^{-t\gl^2},
\end{split}
\]
Hence, we find that for $\Re(s)$ large,
\[
\begin{split}
\eta(B_{A,q}^\oplus,s) &= \frac{1}{\gG\big(\lfrac{s+1}{2}\big)}
\int_0^\infty
\Tr\big[B_{A,q}^\oplus\exp\big(-t(B_{A,q}^\oplus)^2\big)\big]
t^{\frac{s-1}2}dt\\
&= \frac{1}{\gG\big(\lfrac{s+1}{2}\big)} \int_0^\infty -2l
\sum_{\gl\in\spec(B_{q,v})} \gl^2 e^{-t\gl^2}
t^{\frac{s-1}2}dt\\
&= -2l \sum_{\gl\in\spec(B_{q,v})} |\gl|^{1-s}
\frac{1}{\gG\big(\lfrac{s+1}{2}\big)} \int_0^\infty e^{-x}
x^{\frac{s-1}2} dx\\
&= -2l \sum_{\gl\in\spec(B_{q,v})} |\gl|^{1-s},
\end{split}
\]
where we have substituted $x=t\gl^2$. Also note that
interchanging summation and integration can be justified by the
large and small time estimates on $\sum_{\gl\in\spec(B_{q,v})}
\gl^2 e^{-t\gl^2}$, see Proposition \ref{BasicTraceProp} and
Theorem \ref{HeatTrace}. Now, since
\[
\spec(B_{q,v}) = \bigsetdef{\gl\in\R}{ \gl+ q\in \Z},
\]
we find that
\[
\eta(B_{A,q}^\oplus,s) = -2l \sum_{\begin{smallmatrix} n \in\Z \\
n\neq q
\end{smallmatrix}}|n -q|^{1-s}, \quad
\Re(s)>1.
\]
We have included a computation of the value at $s=0$ in
Proposition \ref{ZEtaCalc} (ii). The result is the claimed formula
\[
\eta(B_{A,q}^\oplus) = 2l P_2(q).\qedhere
\]
\end{proof}

\begin{prop}\label{TruncSignOpKern}
The kernel of the truncated odd signature operator is given by
\[
\ker(B_{A,q}^\oplus) = \ker(1\otimes B_{q,v})\cap \ker(D_{A,h})
\cong
\begin{cases}
\hphantom{H^\bullet(} \{0\}, &\text{\rm if $q\notin \Z$}, \\
H^\bullet(\gS,L_B), &\text{\rm if $q\in \Z$}.
\end{cases}
\]
Here, $L_B$ is the trivial line bundle endowed with the flat
connection $B=A-qA_\go$. Moreover, if $q\in\Z$, and if $g$ denotes
the genus of $\gS$, then
\[
H^\bullet(\gS,L_B) \cong \begin{cases}\quad\,\C\oplus
\C^{2g}\oplus \C, &\text{if $B$ is the trivial connection}, \\
\{0\}\oplus \C^{2g-2}\oplus \{0\}, &\text{otherwise.}
\end{cases}
\]
\end{prop}

\begin{proof}
Since $B_{A,q}^\oplus$ is formally self adjoint, we have
\[
\ker(B_{A,q}^\oplus) = \ker(B_{A,q}^\oplus)^2 = \ker\big(1\otimes
B_{q,v}^2+ D_{A,h}^2\big),
\]
where we have used that $\tau_\gS\otimes B_{q,v}$ anti-commutes
with $D_{A,h}$. Since both, $1\otimes B_{q,v}$ and $D_{A,h}$, are
formally self-adjoint we get
\[
\ker\big(1\otimes B_{q,v}^2+ D_{A,h}^2\big) = \ker(1\otimes
B_{q,v})\cap \ker(D_{A,h}).
\]
Now Lemma \ref{VerticalEigenspaces} shows that
\begin{equation}\label{TruncSignOpKern:1}
\ker(1\otimes B_{q,v}) \cong \begin{cases}
\hphantom{\gO^\bullet(} \{0\}, &\text{\rm if $q\notin \Z$}, \\
\gO^\bullet(\gS,L_B), &\text{\rm if $q\in \Z$},
\end{cases}
\end{equation}
where $L_B = L_A\otimes L_\go^{-q}$, endowed with the connection
$B=A\otimes 1 +1\otimes qA_\go$. If $q\notin\Z$ the proof is
finished, so that we assume from now on that $q\in\Z$. Since
$L_A$ is of degree $k$ and $L_\go$ of degree $l$, we find that
$L_A\cong L_\go^{k/l}=L_\go^q$. Thus, $L_B$ is isomorphic to the
trivial line bundle and $B$ is flat. Moreover, Lemma
\ref{VerticalEigenspaces} identifies the restriction of the
operator $D_{A,h}$ to $\ker(1\otimes B_{q,v})$ with the de Rham
operator $D_B$ on $\gO^\bullet(\gS,L_B)$. Using this and
\eqref{TruncSignOpKern:1} we deduce from the Hodge-de-Rham
isomorphism that
\[
\ker(1\otimes B_{q,v})\cap \ker(D_{A,h}) \cong \ker(D_B) \cong
H^\bullet(\gS,L_B).
\]
Now if $B$ is isomorphic to the trivial connection, we have the
well-known cohomology groups of a surface
\[
H^\bullet(\gS) \cong \C\oplus \C^{2g}\oplus \C,
\]
where $g$ is the genus of $\gS$. In the case that $B$ is
non-trivial, the index theorem for the twisted de Rham operator
shows that $\ind(D_B)$ is independent of $B$ as long as $B$ is
flat. Hence,
\[
\sum_{p} (-1)^p \dim H^p(\gS,L_B) =  \ind(D_B) = \ind(D) =
\chi(\gS) = 2-2g.
\]
Moreover, Poincar\'e duality shows that
\[
\dim H^0(\gS,L_B) = \dim H^2(\gS,L_B).
\]
However, if $B$ is a non-trivial flat connection, then
\[
H^0(\gS,L_B) = \{0\}.
\]
One way to see this is as follows\footnote{See also Corollary
\ref{HarmonicFormsTorus} below for a different proof of this
fact.}: Let $\gb: \pi_1(\gS)\to \U(1)$ be the holonomy
representation of $B$. Then according to Proposition
\ref{TwistedDeRhamThm} and \cite[Prop. 5.14]{DavKir},
\[
H^0(\gS,L_B) = \bigsetdef{z\in \C}{ \gb(c)z = z\text{ for all
$c\in \pi_1(\gS)$}}.
\]
Therefore, for $H^0(\gS,L_B)$ is trivial unless $\gb\equiv 1$.
Putting these observations together, we find that for non-trivial
$B$ we have indeed
\[
H^\bullet(\gS,L_B) \cong \{0\}\oplus \C^{2g-2}\oplus
\{0\}.\qedhere
\]
\end{proof}

\subsection{Adiabatic Metrics and the Spectral Flow}\label{S1BundlesAdiab}

After having calculated the Eta invariant for the truncated
signature operator $B_{A,q}^\oplus$, we turn our attention to the
Rho invariant of the odd signature operator $B_{A,q}$. As before,
let $S^1\hookrightarrow M\xrightarrow{\pi}\gS$ be a circle bundle
of degree $l\neq 0$ over a closed surface $\gS$. We let $g_\gS$
be a metric on $\gS$ of unit volume, and $g_v$ be the metric on
$T^vM$ such that the vector field $e$ has length 1. For $\eps>0$
we consider the adiabatic metric \eqref{AdiabaticMetric},
\[
g_\eps:= \lfrac1{\eps^2}g_\gS\oplus g_v,
\]
and let $\nabla^{g_\eps}$ be the Levi-Civita connection
associated to $g_\eps$. For each $t\in [0,1]$ and $\eps>0$ we
define a connection on $TM$ by
\[
\nabla^{\eps,t}:= (1-t)\nabla^\oplus + t\nabla^{g_\eps},
\]
where $\nabla^\oplus$ is the direct sum connection of
\eqref{SplitConn}, which is independent of the scaling parameter
$\eps$.

Now let $L_A\to\gS$ be a holomorphic line bundle of degree $k$,
and let $L:=\pi^*L_A$ be the corresponding flat line bundle over
$M$. Contracting the connection
\[
\nabla^{\eps,t}\otimes 1 + 1\otimes A_q \quad\text{on}\quad
\pi^*\gL^\bullet(T^*\gS)\otimes \pi^*L_A
\]
with the natural Clifford multiplication, we obtain a 2-parameter
family of formally self-adjoint elliptic operators
\[
B_{A,q}^{\eps,t} := \tau_\gS\otimes B_{q,v} + \eps D_{A,h} +
t\eps^2 \tau_MT\quad\text{on}\quad
\pi^*\gO^\bullet(\gS,L_A)\otimes C^\infty(M),
\]
which connects the truncated odd signature operator with the full
odd signature operator associated to $g_\eps$.

\begin{prop}\label{FullSignOpKern}
There exists $\eps_0$ such that for all $\eps<\eps_0$ the
following holds.
\begin{enumerate}
\item If $q\notin \Z$, then for all $t\in[0,1]$
\[
\ker\big(B_{A,q}^{\eps,t}) =\{0\},\quad\text{and}\quad
\SF\big(B^{\eps,t}_{A,q}\big)_{t\in[0,1]} =0.
\]
\item If $q\in\Z$ and $A_q$ is non-trivial, then for all $t\in[0,1]$
\[
\ker\big(B_{A,q}^{\eps,t}) \cong \C^{2g-2},\quad \text{and}\quad
\SF\big(B^{\eps,t}_{A,q}\big)_{t\in[0,1]}=0.
\]
\item For the trivial connection we have
\[
\ker\big(B^{\eps,t}\big)\cong \begin{cases} \C^{2g+1}, &\text{\rm
if $t\neq 0$},\\ \C^{2g+2}, &\text{\rm if $t= 0$}.
\end{cases}\quad\text{and}\quad \SF\big(B^{\eps,t}\big)_{t\in[0,1]}
=\begin{cases} \,\,0,&\text{\rm if $l<0$},\\ -1, &\text{\rm if
$l>0$.}\end{cases}
\]
\end{enumerate}
\end{prop}

\begin{proof}
To keep the notation simple we abbreviate
\[
B_{A,q}^{\eps,t} = \tau_\gS\otimes B_{q,v} + \eps D_{A,h} +
t\eps^2 \tau_MT = :B +\eps D + t\eps^2 S.
\]
According to \eqref{DvDh:anticommute} we have $\{B,D\}=0$, and so
\begin{equation}\label{FullSignOpKern:1}
\big(B_{A,q}^{\eps,t}\big)^2 = B^2 +\eps^2 D^2 +\eps^4t^2 S^2
+\eps^2t\{B,S\} + \eps^3t \{D,S\}.
\end{equation}
By definition the operators $B$, $D$ and $S$ are formally
self-adjoint. Thus,
\[
(\eps^{1/2} B + \eps^{3/2}t S)^2 \ge 0\quad\text{and}\quad
(\eps^{3/2} D + \eps^{3/2}tS)^2\ge 0.
\]
This implies
\[
\eps^2 t \{B,S\}\ge -\eps B^2 - \eps^3t^2 S^2
\quad\text{and}\quad \eps^3t \{D,S\}\ge -\eps^3 D^2 -\eps^3t^2
S^2.
\]
Using this in \eqref{FullSignOpKern:1} we can estimate that for
$\eps<1/2$
\[
\big(B_{A,q}^{\eps,t}\big)^2 \ge \lfrac{\eps^2}2\, \big(B^2+D^2 -
4\eps t^2 S^2\big).
\]
Now $B^2+D^2$ is an elliptic operator, and since $B$ and $D$ are
both formally self-adjoint, it has non-negative spectrum. Hence,
its non-zero eigenvalues are bounded from below by some $\gl>0$.
Moreover, $S$ is an operator of order 0 so that $S^2$ is a
bounded operator. Letting $\eps_0< \min\{\frac 12,
\frac{\gl}{4||S^2||}\}$ we find that for all $\eps<\eps_0$ and
$t\neq 0$
\[
\ker\big(B^2+D^2 - 4\eps t^2 S^2\big) = \ker B\cap\ker D\cap\ker
S,
\]
where we have used that $\ker(B^2+D^2)=\ker B\cap\ker D$.

We now switch back to the usual notation. Using Proposition
\ref{TruncSignOpKern} we can reformulate what we have observed so
far:
\begin{enumerate}
\item If $q\notin \Z$ and $\eps<\eps_0$, then
\[
\ker(B_{A,q}^{\eps,t}) =\{0\}.
\]
\item If $q\in\Z$ and $\eps<\eps_0$, then for all $t\neq 0$
\[
\ker(B_{A,q}^{\eps,t}) = \ker(\tau_MT)\cap
H^\bullet(\gS,L_B),\quad L_B=L_A\otimes L_\go^{-q}.
\]
\end{enumerate}
Thus, in the case $q\notin\Z$ the proof is finished and we assume
henceforth that $q\in\Z$. From Proposition \ref{OddSignCirc} we
know that
\[
\tau_M T =\begin{cases} \quad 0 &\text{\rm on } \gO^{0,0}\oplus
\gO^{1,0},\\ -2\pi l  &\text{\rm on } \gO^{2,0}.
\end{cases}
\]
On the other hand, we have seen in Proposition
\ref{TruncSignOpKern} that
\[
H^\bullet(\gS,L_B) \cong \begin{cases}\quad\,\C\oplus
\C^{2g}\oplus \C, &\text{if $B$ is the trivial connection}, \\
\{0\}\oplus \C^{2g-2}\oplus \{0\}, &\text{otherwise.}
\end{cases}
\]
Therefore, since $l\neq 0$,
\[
\ker(\tau_MT)\cap H^\bullet(\gS,L_B) \cong
\begin{cases}\quad\,\C\oplus \C^{2g}\oplus \{0\}, &\text{if $B$
is the trivial connection}, \\ \{0\}\oplus \C^{2g-2}\oplus \{0\},
&\text{otherwise.}
\end{cases}
\]
This implies part (ii) and the first assertion of (iii). To finish
the proof we still have to compute the spectral flow in the case
that $B$ is the trivial connection. Since the kernels of
$B^{\eps,t}$ are of constant dimension for $t\neq 0$ the only
possible spectral flow contribution is at $t=0$. As we have seen
$\ker( B^{\eps,0})$ contains $H^2(\gS)$ as a summand, whereas
$\ker( B^{\eps,t})$ for $t\neq0$ does not. Now $H^2(\gS)$ is
spanned by the cohomology class of $\vol_\gS$ and we can
explicitly compute that
\begin{equation}\label{CircSmallEV}
B^{\eps,t}(\vol_\gS) = \big(\tau_\gS B_v + \eps D_h + t\eps^2
\tau_MT\big)(\vol_\gS) = -2\pi l\eps^2 t \vol_\gS.
\end{equation}
According to the convention in Definition \ref{SFDef} of how to
count eigenvalues at the endpoints we find that
\[
\SF\big(B^{\eps,t}\big)_{t\in[0,1]} =\begin{cases}
\,\,0,&\text{\rm if $l<0$},\\ -1, &\text{\rm if
$l>0$.}\end{cases}\qedhere
\]
\end{proof}

\begin{remark*}
Note that \eqref{CircSmallEV} also shows that the odd signature
operator associated to the metric $g_\eps$ has a non-trivial
eigenvalue $-2\pi l\eps^2$ which is of order $\eps^2$.
Eigenvalues of this type play a special role in the general
adiabatic limit formula of \cite{Dai91}. We will discuss this in
more detail in Sections \ref{Dai:Gen} and \ref{SF:SpecSeq}, see
in particular Proposition \ref{TopCorrTermCirc}.
\end{remark*}

We now have collected all ingredients to compute the $\U(1)$-Rho
invariant for circle bundles over surfaces.

\begin{theorem}\label{RhoCircBund}
Let $S^1\hookrightarrow M\xrightarrow{\pi} \gS$ be a principal
circle bundle of degree $l\neq 0$. Let $L_A\to\gS$ be a line
bundle over $\gS$ of degree $k$, and let $L:=\pi^*L_A$ be the
corresponding flat line bundle over $M$. Write $q:=k/l$ and
assume that the flat connection $A_q= \pi^*A - iq\go$ is not the
trivial connection. Then
\[
\rho_{A_q}(M) = 2l\big(P_2(q)-\lfrac 16\big) + \sgn(l).
\]
If $M=\gS\times S^1$ is the trivial circle bundle, then all Rho
invariants vanish.
\end{theorem}

\begin{proof}
The Rho invariant associated to the odd signature operator is
independent of the metric. In particular,
\[
\rho_{A_q}(M) = \eta(B_{A,q}^{\eps,1})- \eta(B^{\eps,1}),\quad
\eps>0.
\]
Hence, the variation formula \eqref{GenDiracEtaVar} implies that
for all $\eps>0$
\begin{multline*}
\rho_{A_q}(M) = \eta(B_{A,q}^{\eps,0}) - \eta(B^{\eps,0}) +
2\SF\big(B_{A,q}^{\eps,t}\big)_{t\in[0,1]} -
2\SF\big(B^{\eps,t}\big)_{t\in[0,1]} \\ -\dim(\ker
B_{A,q}^{\eps,1}) +\dim(\ker B^{\eps,1}) +\dim(\ker
B_{A,q}^{\eps,0}) -\dim(\ker B^{\eps,0}).
\end{multline*}
As we have seen in Proposition \ref{EtaTruncSign}, the Eta
invariant associated to the truncated odd signature operator does
not change if the metric on the base is rescaled. Thus,
\[
\eta(B_{A,q}^{\eps,0}) - \eta(B^{\eps,0}) = 2l\big(P_2(q)-\lfrac
16\big).
\]
From Propositions \ref{TruncSignOpKern} and \ref{FullSignOpKern}
we see that if $A_q$ is non-trivial
\[
\dim(\ker B_{A,q}^{\eps,1}) = \dim(\ker B_{A,q}^{\eps,0}).
\]
On the other hand, in the untwisted case
\[
\dim(\ker B^{\eps,1}) - \dim(\ker B^{\eps,0}) = (2g+1) -(2g+2) =
-1.
\]
Lastly, we have seen in Proposition \ref{FullSignOpKern} that for
$\eps$ small enough
\[
\SF\big(B_{A,q}^{\eps,t}\big)_{t\in[0,1]} = 0\quad\text{and}\quad
\SF\big(B^{\eps,t}\big)_{t\in[0,1]} =\begin{cases}\,\, 0,&\text{if
$l<0$},\\ -1, &\text{if $l>0$.}
\end{cases}
\]
Putting all pieces together we find that
\[
\rho_{A_q}(M) = 2l\big(P_2(q)-\lfrac 16\big) + \sgn(l).
\]
The triviality of Rho invariants for $\gS\times S^1$ follows from
\eqref{EtaTrivCirc}.
\end{proof}

\cleardoublepage
\chapter{Rho Invariants of Fiber Bundles, Abstract
Theory}\label{ChapAbst}

This chapter forms the main theoretical part of the thesis. After
having introduced the idea of adiabatic metrics on fiber bundles
and seen their effect in the computation of the Rho invariant for
principal circle bundles, we now want to describe how powerful
tools of local index theory lead to a general formula for the
adiabatic limit of Eta invariants. Since there exists a wide
range of literature on this subject, the ideas presented here are
not new. Nevertheless, we give a detailed account, including some
proofs if feasible.

The treatment starts with the bundle of vertical cohomology
groups over the base of the fiber bundle. To relate it to the
kernel of the vertical de Rham operator, we discuss a fibered
version of the Hodge decomposition theorem. As a byproduct of
this we can prove a result about how to achieve that the mean
curvature of a fiber bundle vanishes. Continuing with the main
line of argument, we give a detailed discussion of the natural
flat connection that exists on the bundle of vertical cohomology
groups. It is precisely this topological nature of the kernel of
the vertical de Rham operator which will make the adiabatic limit
formula accessible for computational purposes. In particular, we
will need to discuss a version of the odd signature operator on
the base twisted by the bundle of vertical cohomology groups.

As the formulation of the general adiabatic limit formula relies
on Bismut's local index theory for families, we continue with a
brief survey of the main constructions and necessary results. In
particular, we include a short discussion of superconnections and
associated Dirac operators. Returning to the context of fiber
bundles, we introduce the Bismut superconnection and recall how
it appears in the local index theorem for families.

With these notions at hand, we will give a heuristic derivation
of the adiabatic limit formula for families of odd signature
operators. Then, referring to the literature for rigorous proofs,
we finally state the general adiabatic limit formula for the Eta
invariant due to Dai. One of the terms appearing there has a
topological interpretation in terms of the Leray-Serre spectral
sequence, and we discuss this briefly.

We finish this chapter using the adiabatic limit formula to
derive again the formula for the Rho invariant of a principal
$S^1$-bundle over a closed surface. Although we have obtained the
formula already in the last chapter, it is illuminating to
observe how the abstract theory leads to a shorter a more
conceptual proof.

\clearpage

\section{The Bundle of Vertical Cohomology Groups}\label{VertCohom}

\subsection{The Vertical de Rham Operator}\label{VertDeRham}

In Remark \ref{PartDeRhamOpRem} we have pointed out without
further comments that there is a relationship between differential
operators acting fiberwise and families of differential operators
in the sense of Definition \ref{SmoothFamily}. We want to make
this a bit more precise now. Let $F\hookrightarrow
M\xrightarrow{\pi} B$ be an oriented fiber bundle, where as before
all manifolds are assumed to be closed, connected and oriented,
and let $E\to M$ be a Hermitian over $M$.

\begin{dfn}\label{FiberEllOpDef}
Let $D:\gO^\bullet(M,E)\to \gO^\bullet(M,E)$ be a differential
operator. Then we call $D$ a \emph{fiberwise} differential
operator, if
\[
[D,\pi^*\gf] = 0,\quad \text{for all }\gf\in C^\infty(B).
\]
We call $D$ \emph{fiberwise elliptic} if in addition its
principal symbol
\[
\gs(D)(x,\xi): E_x\to E_x,\quad x\in M,
\]
is invertible for every non-vanishing $\xi\in T^v_xM^*$.
\end{dfn}

Certainly, if $T^vM$ is endowed with a metric, and $A$ is a flat
connection on $E$, the vertical de Rham operator $D_{A,v}$ as in
Definition \ref{PartDeRhamOpDef} is a fiberwise elliptic operator
in the sense of this definition.\\

\noindent\textbf{Local Trivializations and Families.} To relate
fiberwise differential operators with families of differential
operators as in Definition \ref{SmoothFamily}, we describe a
particular way to construct local trivializations of the fiber
bundle, see \cite[Lem. 1.3.3]{GLP}.

\begin{lemma}\label{FiberTrivialization}
Let $g=g_B\oplus g_v$ be a submersion metric as in Section
\ref{FiberedCalc}. Let $y\in B$ and $F:=\pi^{-1}(y)$. Then for
every sufficiently small geodesic neighbourhood $U$ around $y$
there exists an isomorphism of fiber bundles
\[
\gF:U\times F \to \pi^{-1}(U),
\]
such that for all $(u,x)\in U\times F$ and every vector $v\in
T_{y}U\subset T_{(y,x)}U\times F$,
\[
\gF(y,x) = x, \quad \pi\circ \gF(u,x) = u,\quad \text{and}\quad
\gF_* v = v^h,
\]
where $v^h$ refers to the horizontal lift of $v$.
\end{lemma}

\begin{proof}
Let $b=\dim B$, and let $U\subset B$ be a geodesic ball centered
in $y$. We identify $U$ with an open ball in $\R^b$ in such a way
that $y=0$. We can use the horizontal projection $P^h:TM\to T^hM$
to lift the coordinate vector fields $\pd_a$ to horizontal vector
fields $\pd_a^h$ on $\pi^{-1}(U)$. Identifying a point $u\in U$
with the vector field $u^a\pd_a$ we get a vector field $u^h=
u^a\pd_a^h$ on $\pi^{-1}(U)$, and hence a flow
\[
\gF_t(u,.):M\to M,\quad u\in U.
\]
It follows from the construction that for small $t$, the flow
$\gF_t(u,.)$ maps the fiber $F$ diffeomorphically onto the fiber
$\pi^{-1}(tu)$. Moreover, for all $x\in F$ we have
$\gF_{st}(u,x)=\gF_t(su,x)$, so that we can choose $U$ small
enough to define
\[
\gF:U\times F \to \pi^{-1}(U),\quad (u,x)\mapsto \gF_1(u,x).
\]
The claimed properties all follow immediately from this
definition.
\end{proof}


Using a fiber bundle chart $\gF:U\times F\to \pi^{-1}(U)$ of the
form just described, we can transfer all geometric structures on
$\pi^{-1}(U)$ to $U\times F$: First of all, it is straightforward
to check that
\[
\gF^*T^vM  = U\times TF.
\]
Therefore, the pullback $\gF^*(g_v)$ of the vertical metric is
the same as a family of Riemannian metrics $g_{F,u}$ on $F$. In
the same way the restriction of $\nabla^v$ to $T^vM$ induces a
family $\nabla^{F,u}$ of covariant derivatives on $TF$, and
Proposition \ref{CanVerticalConn} shows that each $\nabla^{F,u}$
is the Levi-Civita connection on $F$ with respect to the metric
$g_{F,u}$, compare with Remark \ref{NablaNatRem}. Similarly, we
pull the horizontal distribution $T^hM$ back to $U\times F$ and
use this to identify
\begin{equation}\label{VertFormsLocChart}
\gO_v^\bullet\big(\pi^{-1}(U)\big) \cong
C^\infty\big(U,\gO^\bullet(F)\big),\quad\text{and}\quad
\gO^{p,q}\big(\pi^{-1}(U)\big) \cong \gO^p\big(U,\gO^q(F)\big).
\end{equation}
Note, however, that $\gF^*T^hM$ will in general not coincide with
$TU\times F$, unless the curvature $\gO$ of the fiber bundle is
trivial.

\begin{remark*}
We want to give a note about the definition of
$\gO^p\big(U,\gO^q(F)\big)$. A naive way---which is sufficient
for our purposes---is to define elements of
$C^\infty\big(U,\gO^q(F)\big)$ to be locally of the form
\[
\sum_{|I|=q} f_I(y,x)dx^I,\quad y\in U,
\]
with local coordinates $x_i$ for $F$, multi indices $I$, and
smooth functions $f_I(y,x)$ satisfying the appropriate
transformation laws with respect to changes of the coordinate
chart. Similarly one treats elements of
$\gO^p\big(U,\gO^q(F)\big)$. From a more invariant perspective,
one could endow $\gO^q(F)$ with its natural Fr\'{e}ch\'{e}t
topology and consider smooth maps with respect to this.
\end{remark*}

Under the identification of \eqref{VertFormsLocChart}, a vertical
differential operator $D$ on $\gO^\bullet(M)$ can be written over
$U\times F$ as
\[
D = \sum_j K_j(u)\otimes D_j(u),
\]
where each $K_j(u)$ is a bundle endomorphism of $\gL^\bullet T^*U$
and each $D_j(u)$ is a smooth $b$-parameter family of
differential operators on $\gO^\bullet(F)$ in the sense of
Definition \ref{SmoothFamily}.\\

\noindent\textbf{Vertical de Rham Operators.} We use the above
digression to give a description of the vertical de Rham
operator. For this we first need to incorporate a bundle $E$ over
$M$, endowed with a flat connection $A$. Let $\pi_F:U\times F\to
F$ be the projection onto the second factor, and denote by $E|_F$
be the restriction of $E$ to the fiber $F$.

\begin{lemma}\label{FamilyFlatConn}
There exists a natural lift of the bundle isomorphism
\[
\gF: U\times F\to \pi^{-1}(U)
\]
to an isomorphism of flat Hermitian vector bundles
\[
\gF_E:\pi_F^*(E|_F)\to E|_{\pi^{-1}(U)}.
\]
\end{lemma}

\begin{proof}
Let $u\in U$, and let $F_u$ be the fiber over $u$. Then
$\gF(u,.)$ maps $F$ diffeomorphically to $F_u$. Using parallel
transport with respect to $A$ along the flow lines of $\gF_t(u,.)$
we can lift this to an isomorphism
\[
\gF_E(u,.): E|_F\to E|_{F_u}.
\]
Now, since $A$ is a flat Hermitian connection, the bundles $E|_F$
and $E|_{F_u}$ are naturally endowed with flat Hermitian
connections induced by $A$. Since we are using parallel transport
with respect to $A$, a locally constant, unitary frame for $E|_F$
will be mapped by $\gF_E$ to a locally constant, unitary frame for
$E|_{F_u}$. This implies that $\gF_E$ is, in fact, an isomorphism
of flat Hermitian bundles.
\end{proof}

In a similar way as in \eqref{VertFormsLocChart}, we can use Lemma
\ref{FamilyFlatConn} to identify
\begin{equation}\label{VertFamilyIdentify}
\gO^{p,q}\big(\pi^{-1}(U),E|_{\pi^{-1}(U)}\big) \cong
\gO^p\big(U,\gO^q(F,E|_F)\big),
\end{equation}
where $E|_F$ is endowed with a fixed flat Hermitian connection
$A_F$. We then let $D_{A_F,u}$ be the family of de Rham operators
on $\gO^\bullet(F,E|_F)$ associated to the metric $g_{F,u}$ and
the flat connection $A_F$. Then under the identification
\eqref{VertFamilyIdentify} we can write
\begin{equation}\label{VertDeRhamFamily}
D_{A,v} = (-1)^p\otimes D_{A_F,u}\quad \text{on }
\gO^p\big(U,\gO^\bullet(F,E|_F)\big).
\end{equation}

\begin{prop}\label{VertDeRhamKer}
The $C^\infty(B)$-module $\ker(D_{A,v})\cap \gO^\bullet_v(M,E)$
is isomorphic to the space of smooth sections of a vector bundle,
which we denote by $\sH_{A,v}^\bullet(M) \to B$. Moreover,
\[
\ker(D_{A,v}) \cong \gO^\bullet\big(B,\sH_{A,v}^\bullet(M)\big).
\]
\end{prop}

\begin{proof}[Sketch of proof]
It suffices to show that the assertion is true locally, i.e.,
that for sufficiently small open subsets $U\subset M$
\[
\ker(D_{A,v})\cap \gO^\bullet\big(\pi^{-1}(U),E\big)
\]
is isomorphic to the space of differential forms over $U$ with
values in a vector bundle. For this let $U\subset B$ be as in
Lemma \ref{FiberTrivialization} such that $\pi^{-1}(U)\cong
U\times F$, and write $D_{A,v}$ as in \eqref{VertDeRhamFamily}.
Since $D_{A,v}$ acts as $\pm\Id$ on $\gO^\bullet(U)$, it suffices
to consider $D_{A_F,u}$ acting on
\[
\gO^\bullet_v\big(\pi^{-1}(U),E\big) \cong
C^\infty\big(U,\gO^\bullet(F,E|_F)\big).
\]
For fixed $u$, the Hodge-de-Rham theorem for $D_{A_F,u}$ implies
that $\ker(D_{A_F,u})$ is isomorphic to $H^\bullet(F,E_A|_F)$,
where $E_A|_F$ is short for $E|_F$ endowed with the flat
connection $A_F$. Since we know from Lemma \ref{FamilyFlatConn}
that $A_F$ does not vary with $u$, we infer that
$\dim\ker(D_{A_F,u})$ is constant for $u\in U$. Hence, we are
precisely in the situation of Proposition
\ref{ProjGreenSmooth}---respectively Remark \ref{SmoothDepRem}.
Therefore, the family of projections
\[
P_u: \gO^\bullet(F,E|_F) \to \ker\big(D_{A_F,u}\big),\quad u\in U,
\]
is a smooth family of finite rank smoothing operators. Using this,
it is straightforward to check that the collection
\[
\sH_{A,v}^\bullet(U):=\bigcup_{u\in U}\ker\big(D_{A_F,u}\big) \to
U
\]
forms a smooth vector bundle over $U$, see \cite[Lem. 9.9]{BGV}
for a detailed proof. Then the assertion of Proposition
\ref{VertDeRhamKer} easily follows.
\end{proof}

\subsection{Vertical Hodge Decomposition}

We now want to use Proposition \ref{VertDeRhamKer} to prove the
following fibered version of the de Hodge decomposition theorem.

\begin{theorem}\label{FiberedHodge}
Let $E\to M$ be a Hermitian vector bundle over the total space of
an oriented fiber bundle of closed manifolds $F\hookrightarrow
M\xrightarrow{\pi} B$. Assume that $E$ admits a flat connection
$A$. With respect to every submersion metric, there is an
$L^2$-orthogonal splitting of smooth forms
\[
\begin{split}
\gO^\bullet(M,E) &= \ker(D_{A,v}) \oplus \im(D_{A,v}) \\
&= \big(\ker d_{A,v}\cap \ker d_{A,v}^t\big) \oplus \im(d_{A,v})
\oplus \im(d_{A,v}^t).
\end{split}
\]
Moreover, the splitting is independent of the chosen metric $g_B$
on $B$.
\end{theorem}

\begin{proof}
We start with a local consideration. With the same notation as in
the proof of Proposition \ref{VertDeRhamKer} we consider the
family of de Rham operators $D_{A_F,u}$ on $\gO^\bullet(F,E|_F)$.
We know from the proof of Proposition \ref{VertDeRhamKer} that
$\dim\ker(D_{A_F,u})$ is constant for $u\in U$, so that the
family of projections $P_u$ onto the kernels depends smoothly on
$u$. According to Proposition \ref{ProjGreenSmooth} (ii) the same
is true for the family of Green's operators,
\[
G_u:\gO^\bullet(F,E|_F)\to \gO^\bullet(F,E|_F).
\]
Let $\go_u\in C^\infty\big(U, \gO^\bullet(F,E|_F)\big)$. For fixed
$u$ we can decompose
\[
\go_u = P_u\go_u + D_u\circ G_u\circ(\Id-P_u)\go_u,
\]
and both summands depend smoothly on $u$ as $P_u$ and $G_u$ do
so. Writing $D_{A,v}|_{\pi^{-1}(U)}$ for the restriction of
$D_{A,v}$ to $\gO^\bullet\big(\pi^{-1}(U),E|_{\pi^{-1}(U)}\big)$,
one readily concludes that
\begin{equation}\label{LocalHodgeDec}
\gO^\bullet\big(\pi^{-1}(U),E|_{\pi^{-1}(U)}\big) =
\ker\big(D_{A,v}|_{\pi^{-1}(U)}\big) \oplus
\im\big(D_{A,v}|_{\pi^{-1}(U)}\big).
\end{equation}

Now let $\{\gf_i\}$ be a partition of unity on $B$, subordinate to
a finite covering $B=\bigcup_i U_i$, such that
\eqref{LocalHodgeDec} holds for every
$\gO^\bullet\big(\pi^{-1}(U_i),E|_{\pi^{-1}(U_i)}\big)$. For
$\go\in \gO^\bullet(M,E)$ and every $i$ we can decompose
\[
(\pi^*\gf_i) \go = \ga_i + D_{A,v}\gb_i,\quad \ga_i\in
\ker\big(D_{A,v}|_{\pi^{-1}(U_i)}\big),\quad \gb_i\in
\gO^\bullet\big(\pi^{-1}(U_i),E|_{\pi^{-1}(U_i)}\big).
\]
Then, since $D_{A,v}$ is $C^\infty(B)$ linear,
\[
\go = \sum_i(\pi^*\gf_i)\ga_i+ \sum_i(\pi^*\gf_i)D_{A,v}\gb_i =
\sum_i(\pi^*\gf_i)\ga_i+ D_{A,v}
\Big(\sum_i(\pi^*\gf_i)\gb_i\Big),
\]
which is a decomposition of $\go$ in terms of $\ker(D_{A,v})
\oplus \im(D_{A,v})$. Clearly, the decomposition is
$L^2$-orthogonal, since $D_{A,v}$ is formally self-adjoint. The
equalities
\[
\ker(D_{A,v})= \big(\ker d_{A,v}\cap \ker d_{A,v}^t\big)
\quad\text{and}\quad \im(D_{A,v})= \im(d_{A,v}) \oplus
\im(d_{A,v}^t)
\]
follow as in the unparametrized case. Finally, the assertion that
the vertical Hodge decomposition is independent of the metric
$g_B$ on $B$ is immediate from the fact that $D_{A,v}$ is
independent of $g_B$.
\end{proof}

\subsection{Vanishing Mean Curvature}

Before we continue the discussion of the bundle
$\sH_{A,v}^\bullet(M) \to B$, we want to give an interesting
application of the fibered Hodge decomposition theorem. The
corresponding result for foliations is \cite[Thm. 4.18]{Dom98}.
However, in the case of fiber bundles, the proof can be
simplified considerably, and the author of this thesis is not
aware of a corresponding treatment in the literature. This
subsection is not essential for the line of thoughts in the later
sections. Yet, it might be helpful in more complicated examples.

\begin{theorem}\label{VanishingMeanCurv}
Let $g_v$ be a metric on $T^vM$ of unit volume. Then there exists
a vertical projection $P^v:TM\to T^vM$ such that the associated
mean curvature form $k_v(g_v,P^v)$ vanishes.
\end{theorem}

Since we have seen in Lemma \ref{VolumeNorm} that we can deform a
vertical metric conformally to a metric of unit volume, Theorem
\ref{VanishingMeanCurv} implies

\begin{cor}\label{VanishingMeanCurvCor}
Every oriented fiber bundle of closed manifolds admits a
connection and a vertical metric such that the mean curvature form
vanishes.
\end{cor}

Before we give the proof of Theorem \ref{VanishingMeanCurv}, we
extract the part where we will use the vertical Hodge
decomposition of Theorem \ref{FiberedHodge}. Recall that we have
introduced the basic projection $\Pi_B$ in Definition
\ref{BasicProj}.

\begin{prop}\label{BasicProjKer}
There is an $L^2$-orthogonal splitting of smooth horizontal forms,
\[
\gO^\bullet_h(M) = \pi^*\gO^\bullet(B)\oplus
d_v^t\big(\gO^{\bullet,1}(M)\big).
\]
Moreover, the kernel of the basic projection is given by
\[
\ker \Pi_B = d_v^t\big(\gO^{\bullet,1}(M)\big).
\]
\end{prop}

\begin{proof}
First of all, as the 0th cohomology group of the fiber consists
only of constant functions, one deduces from Proposition
\ref{VertDeRhamKer} that
\[
\ker D_v\cap \gO^\bullet_h(M) = \pi^*\gO^\bullet(B).
\]
Since $\gO^\bullet_h(M)\perp \im(d_v)$, the vertical Hodge
decomposition in Theorem \ref{FiberedHodge} yields
\[
\gO^\bullet_h(M) = \pi^*\gO^\bullet(B)\oplus
d_v^t\big(\gO^{\bullet,1}(M)\big).
\]
For the second assertion we note that for all
$\ga\in\gO^\bullet(B)$ and $\go\in \gO^\bullet_h(M)$
\[
\Scalar{\ga}{\Pi_B(\go)}= v_F^{-1} \int_{M/B}\scalar
{\pi^*\ga}{\go}\vol_F,
\]
where $v_F$ is the function which associates to a point $y\in B$
the volume of the fiber over $y$, see Definition \ref{BasicProj}.
This implies that $\ker\Pi_B\perp\pi^*\gO^\bullet(B)$. On the
other hand, as in the unparametrized case, one finds that for
$\go\in \gO^{\bullet,1}(M)$,
\[
\Pi_B(d^t_v\go) = v_F^{-1} \int_{M/B} d^t_v(\go)\wedge \vol_F= 0.
\]
Therefore,
\[
d_v^t\big(\gO^{\bullet,1}(M)\big)\subset \ker \Pi_B,
\]
which finishes the proof.
\end{proof}

\begin{proof}[Proof of Theorem \ref{VanishingMeanCurv}]
Let $g_v$ be a vertical metric such that $v_F(g_v)=1$. We choose
an arbitrary vertical projection $P^v:TM\to T^vM$, and let $g$ be
a submersion metric on $M$ satisfying
\[
T^vM^\perp = \ker P^v\quad\text{and}\quad g\big|_{T^vM\times T^vM}
= g_v.
\]
According to Corollary \ref{BasicMeanCurv}, the assumption that
$v_F(g_v)=1$ implies that the basic projection of the mean
curvature $k_v$ vanishes. Form Proposition \ref{BasicProjKer} we
deduce that there exists $\eta\in \gO^{1,1}(M)$ such that
\[
d_v^t\eta = k_v \in\gO^{1,0}(M).
\]
Define $h\in C^\infty(M,T^*M\otimes T^*M)$ by
\[
h(X,Y):= \eta(P^hX,P^vY) + \eta(P^hY,P^vX),\quad X,Y\in
C^\infty(M,TM).
\]
Here, $P^h=\Id-P^v$ is the horizontal projection. Then $h$ is a
symmetric 2-tensor, and $h(X,X)=0$ for all $X\in C^\infty(M,TM)$.
Thus, we can define a new metric on $TM$ by letting
\[
\widetilde g:= g + h.
\]
Note that the restriction of $\widetilde g$ to $T^vM$ still
coincides with $g_v$. Let $\widetilde P^{v/h}$ denote the
vertical respectively horizontal projection associated to
$\widetilde g$. Then, if $\{e_i\}$ is any local orthonormal frame
for $T^vM$ with respect to $g_v$, we have
\[
\widetilde P^v(X) = P^v(X) + \sum_i \eta(X,e_i)e_i, \quad
\widetilde P^h(X) = P^h(X) - \sum_i \eta(X,e_i)e_i.
\]
Let $\widetilde k_v$ denote the mean curvature form associated to
$g_v$ and $\widetilde P^v$, and let $\{f_a\}$ be a local
orthonormal frame for $T^hM$ with respect to the original metric
$g$. Then, according to formula \eqref{MeanCurvRem} for the mean
curvature,
\[
\widetilde k_v(f_a) = -\sum_i g_v\big([\widetilde
P^hf_a,e_i],e_i\big) = k_v(f_a) + \sum_{ij}
g_v\big([\eta(f_a,e_j)e_j,e_i],e_i\big).
\]
Now, using standard arguments involving the Lie bracket and the
fact that $\nabla^v$ is metric and torsion-free as a connection on
$T^vM$, one gets
\[
\begin{split}
g_v\big([\eta(f_a,e_j)e_j,e_i],e_i\big) &= - e_i
\big[\eta(f_a,e_j)\big]g_v(e_j,e_i) + \eta(f_a,e_j)
g_v\big([e_j,e_i],e_i\big)\\
&= - e_i \big[\eta(f_a,e_i)\big] + \eta(f_a,e_j)
g_v(e_j,\nabla^v_{e_i}e_i).
\end{split}
\]
On the other hand, according to Proposition \ref{d^tSplit},
\[
d_v^t\eta(f_a) = \sum_{i}(\nabla^{\oplus}_{e_i}\eta)(f_a,e_i) =
\sum_{i}\Big(e_i \big[\eta(f_a,e_i)\big] -
\eta(\nabla_{e_i}^\oplus f_a,e_i)
-\eta(f_a,\nabla^\oplus_{e_i}e_i)\Big).
\]
Since $\nabla_{e_i}^\oplus f_a=0$ and $\nabla^\oplus_{e_i}e_i =
\nabla^v_{e_i}e_i$ we conclude that
\[
\sum_{ij} g_v\big([\eta(f_a,e_j)e_j,e_i],e_i\big)=
-d_v^t\eta(f_a).
\]
Employing the definition of $\eta$ we have thus achieved that
\[
\widetilde k_v = k_v - d_v^t\eta = 0.
\]
This shows that the mean curvature associated to $\widetilde P^v$
and $g_v$ vanishes.
\end{proof}

\begin{remark*}
The statement of Theorem \ref{VanishingMeanCurv} for foliations is
not true without changes. The underlying reason is that Corollary
\ref{BasicMeanCurv} does not generalize, i.e., the basic
projection of the mean curvature does not necessarily give a
trivial cohomology class. Note that the definition of cohomology
requires extra work for the possibly singular leaf space of a
foliation. But even if one uses the \emph{basic cohomology} as
the correct substitute, Corollary \ref{BasicMeanCurv} does not
carry over, and one finds topological obstructions to the
vanishing of the mean curvature form. In the language of
foliation theory, Corollary \ref{VanishingMeanCurvCor} asserts
that the vertical distribution of a fiber bundles is a
\emph{taut} foliation. For a detailed discussion of the aspects
mentioned here, in particular the difference between \emph{tense}
and taut foliations, we refer to \cite{Dom98} and references
given therein.
\end{remark*}

\subsection{A Flat Connection on the Bundle of Vertical Cohomology
Groups}

Let $E\to M$ be a Hermitian vector bundle over the total space of
the fiber bundle $F\hookrightarrow M\xrightarrow{\pi} B$, and
assume that $E$ admits a flat connection $A$. As we have seen in
Corollary \ref{dSquare} there is a vertical differential
\[
d_{A,v}: \gO^\bullet_v(M,E) \to \gO^{\bullet+1}_v(M,E),\quad
d_{A,v}^2=0.
\]
Hence, we can form the quotient $\ker d_{A,v}/\im d_{A,v}$. If
$M$ is endowed with a vertical metric, it is an immediate
consequence of the fibered Hodge decomposition theorem that
\begin{equation}\label{VertHodgeIsom}
\ker d_{A,v}/\im d_{A,v} \cong \ker(D_{A,v})\cap
\gO^\bullet_v(M,E),
\end{equation}
which is the space of sections of the bundle
$\sH_{A,v}^\bullet(M) \to B$ of Proposition \ref{VertDeRhamKer}.
This is not surprising, since if $F\subset M$ is a fiber of
$\pi:M\to B$ we can restrict $d_{A,v}$ as in Section
\ref{VertDeRham} to an operator
\[
d_{A,v}|_F: \gO^\bullet(F,E|_F)\to \gO^{\bullet+1}(F,E|_F).
\]
Then $\ker (d_{A,v}|_F)/\im (d_{A,v}|_F)$ is just the de Rham
cohomology of $F$ with values in the flat bundle $E|_F$, so that
$\ker d_{A,v}/\im d_{A,v}$ is roughly the union of the cohomology
groups of all fibers. These observations and Proposition
\ref{VertDeRhamKer} imply

\begin{prop}\label{VertCohomBundleDef}
The space $\ker d_{A,v}/\im d_{A,v}$ is isomorphic to the space of
sections of a finite rank vector bundle $H_{A,v}^\bullet(M)\to
B$, which we call the \emph{``bundle of vertical cohomology
groups''}. Its fiber over a point $y\in B$ is isomorphic to the
de Rham cohomology of $\pi^{-1}(y)$ with values in the flat bundle
$E|_{\pi^{-1}(y)}$.
\end{prop}

\begin{remark*}\quad\nopagebreak
\begin{enumerate}
\item Although the bundles $H_{A,v}^\bullet(M)$ and
$\sH_{A,v}^\bullet(M)$ are isomorphic, we will usually not
identify them, since the latter is only defined when we have
chosen a vertical metric. This is why we have introduced the term
``bundle of vertical cohomology groups'' only here rather than
already in Proposition \ref{VertDeRhamKer}.
\item As for the bundle $\sH^\bullet_{A,v}(M)$ we consider also
differential forms on $B$ with values in the bundle of vertical
cohomology groups. Then it is immediate that
\begin{equation}\label{VertCohomValuedForms}
\gO^p\big(B,H^q_{A,v}(M)\big) \cong
\frac{\ker\big(d_{A,v}:\gO^{p,q}(M,E)\to
\gO^{p,q+1}(M,E)\big)}{\im\big(d_{A,v}:\gO^{p,q-1}(M,E)\to
\gO^{p,q}(M,E)\big)}.
\end{equation}
\end{enumerate}
\end{remark*}

\noindent\textbf{The Natural Flat Connection.} Recall from
Corollary \ref{dSquare}, that on $\gO^\bullet(M,E)$,
\[
d_{A,h}^2 + \big\{d_{A,v},\imu(\gO) \big\}=0,\quad\text{and}\quad
\big\{d_{A,v}, d_{A,h} \big\}=\big\{d_{A,h}, \imu(\gO) \big\}=0,
\]
where $\gO$ is the curvature of the fiber bundle. Now the fact
that $d_{A,h}$ anti-commutes with $d_{A,v}$ implies that $d_{A,h}$
descends to a well-defined map
\[
\Bar d_{A,h}: \gO^{p}\big(B,H_{A,v}^\bullet(M)\big) \to
\gO^{p+1}\big(B,H_{A,v}^\bullet(M)\big).
\]
Moreover, if $\go \in \gO^{p,q}(M,E)$ satisfies $d_{A,v}\go=0$,
then it follows from the relation
\[
d_{A,h}^2\go = -  d_{A,v}\circ \imu(\gO)\go
\]
that $d_{A,h}^2\go$ is a $d_{A,v}$-exact element of
$\ker\big(d_{A,v}:\gO^{p+2,q}(M,E)\to \gO^{p+2,q+1}(M,E)\big)$.
This implies that $(\Bar d_{A,h})^2 =0$. In other words, we have
found a natural flat connection on the bundle of vertical
cohomology groups which is induced by $d_{A,h}$.

\begin{dfn}\label{FlatVertCohom}
We denote by
\[
\nabla^{H_{A,v}}: C^\infty\big(B,H_{A,v}^\bullet(M)\big) \to
\gO^1\big(B,H_{A,v}^\bullet(M)\big)
\]
the flat connection defined by $\Bar d_{A,h}$. More precisely, for
all $X\in C^\infty(B,TB)$ and $\go \in \gO^\bullet_v(M,E)$ with
$d_{A,v}\go=0$ we define
\[
\nabla^{H_{A,v}}_X [\go] := \big[\imu(X^h) d_{A,h} \go\big]\in
C^\infty\big(B,H_{A,v}^\bullet(M)\big).
\]
\end{dfn}

\noindent\textbf{Relation to the Leray-Serre Spectral Sequence.}
We want to point out that we have just constructed the term
$(E_1^{\bullet,\bullet},d_1)$ of the spectral sequence associated
to the complex $\big(\gO^\bullet(M,E),d_A\big)$. To explain
this---and also for later use---we make a short digression on the
Leray-Serre spectral sequence.

Recall, e.g. from \cite[Sec 2.2]{McC}, that a complex endowed
with a decreasing filtration gives rise to a spectral sequence.
The appropriate filtration in the case at hand is the \emph{Serre
filtration} given by
\begin{equation}\label{SerreFilt}
F^k \gO^\bullet := \sum_{p\ge k} \gO^{p,\bullet}(M,E).
\end{equation}
Then, if $b=\dim B$ and $d=d_A$, we have
\[
\{0\}=F^{b+1} \gO^\bullet \subset F^b \gO^\bullet  \subset \ldots
\subset  F^0 \gO^\bullet  = \gO^\bullet ,\quad d\big(F^k
\gO^\bullet \big)\subset F^k \gO^\bullet .
\]
Note that the latter follows from Proposition \ref{dSplit}, since
each of the terms appearing in $d_A$,
\[
d_A = d_{A,v}+d_{A,h}+ \imu(\gO),
\]
preserves $F^k \gO^\bullet$. In the same way one verifies that the
de Rham cohomology $H^\bullet(M,E_A)$ inherits a filtration
defined by
\begin{equation}\label{CohomFilt}
F^kH^n:= \im\big(H^n(F^k\gO^\bullet,d)\to H^n(M,E_A)\big).
\end{equation}
One now constructs a spectral sequence as follows, see the proof
of \cite[Thm. 2.6]{McC}. For $r\in\N$ define
\begin{equation}\label{SpecSeqDef}
\begin{split}
Z_r^{p,q} &:= F^p \gO^{p+q} \cap
d ^{-1}\big(F^{p+r}\gO^{p+q+1} \big),\\
B_r^{p,q} &:= F^p\gO^{p+q}  \cap d  \big(F^{p-r}\gO^{p+q-1}
\big),\\
E_r^{p,q} &:=
\frac{Z_r^{p,q}}{Z_{r-1}^{p+1,q-1}+B_{r-1}^{p,q}},\quad
E_0^{p,q}:= \frac{F^p\gO^{p+q}}{F^{p+1}\gO^{p+q}}\cong
\gO^{p,q}(M,E).
\end{split}
\end{equation}
Then the differential $d$ naturally defines on each bigraded
module $E_r^{\bullet,\bullet}$ a differential $d_r$ of bidegree
$(r,1-r)$ in such a way that
\[
E_{r+1}^{p,q} \cong \frac{\ker \big(d_r:E_r^{p,q}\to
E_r^{p+r,q+1-r}\big)}{\im \big(d_r:E_r^{p-r,q+r-1}\to
E_r^{p+r,q+1-r}\big)}.
\]
The general theory of spectral sequences now implies the
following, see again \cite[Thm. 2.6]{McC}.

\begin{theorem}\label{SpecSeqThm}
The spectral sequence $\big(E_r^{\bullet,\bullet},d_r\big)$
collapses for $r=b+1$ and converges to $H^\bullet(M,E_A)$. More
precisely, for all $r\ge b+1$
\[
E_r^{p,q} \cong \frac{F^pH^{p+q}}{F^{p+1}H^{p+q}},
\]
where $F^pH^\bullet$ is the filtration of $H^\bullet(M,E)$ given
by \eqref{CohomFilt}.
\end{theorem}

Now we can interpret the bundle of vertical cohomology groups in
terms of the Leray-Serre spectral sequence. Certailnly, the term
$E_0^{\bullet,\bullet}$ in \eqref{SpecSeqDef} coincides with
$\gO^{\bullet,\bullet}(M,E)$. Moreover, one easily verifies that
the natural construction of the differential in the proof of
\cite[Thm. 2.6]{McC} coincides with $d_{A,v}$. Thus,
\[
\big(E_0^{\bullet,\bullet},d_0\big) =
\big(\gO^{\bullet,\bullet}(M,E),d_{A,v}\big).
\]
In the discussion following \eqref{VertCohomValuedForms} we have
constructed a natural differential $\Bar d_{A,h}$ on the
cohomology of $\big(E_0^{\bullet,\bullet},d_0\big)$, and again,
one can easily check that it coincides with the differential on
$E_1^{\bullet,\bullet}$ abstractly constructed from
\eqref{SpecSeqDef}. Without giving more details we summarize that

\begin{lemma}\label{SpecSeqLem}
The Leray-Serre spectral sequence satisfies
\[
E_1^{p,q} \cong \big(\gO^p\big(B,H^q_{A,v}(M)\big),\Bar
d_{A,h}\big) \quad\text{and}\quad E_2^{p,q} \cong
H^p\big(\gO^\bullet\big(B,H^q_{A,v}(M)\big).
\]
\end{lemma}

\subsection{Twisting with the Bundle of Vertical Cohomology
Groups}\label{FlatHodge}

We can also use the vertical Hodge decomposition of Theorem
\ref{FiberedHodge} to give a Hodge theoretic description of the
flat connection $\nabla^{H_{A,v}}$. We fix a vertical metric $g_v$
and use this to identify $H_{A,v}^\bullet(M)$ with
$\sH_{A,v}^\bullet(M)$ using the vertical Hodge-de-Rham
isomorphism in \eqref{VertHodgeIsom}. From Section
\ref{FiberedExtDiff} we know that $\gO^\bullet_v(M,E)$ is endowed
with the natural connection $\widetilde\nabla^{A,v}$, induced by
the vertical Lie derivative and the connection $A$. Then we have
the following, see also \cite[Prop 3.14]{BL95}.

\begin{prop}\label{FlatConnHodge}
Under the vertical Hodge-de-Rham isomorphism the flat connection
$\nabla^{H_{A,v}}$ coincides with the connection defined by
\[
\nabla^{\sH_{A,v}}_X := P_{\ker(D_{A,v})}\circ
\widetilde\nabla^{A,v}_{X^h} ,\quad X\in C^\infty(B,TB).
\]
\end{prop}

\begin{proof}
For convenience we drop the reference to the flat connection $A$.
Denote by
\[
\Psi: \ker d_v\cap\gO^\bullet_v(M) \to
C^\infty\big(B,H_v^\bullet(M)\big)
\]
the quotient map. Then, according to Definition
\ref{FlatVertCohom},
\[
\nabla_X^{H_v} \big(\Psi(\go)\big) = \Psi\circ\imu(X^h)\circ
d_h(\go),\quad \go\in \ker d_v\cap\gO^\bullet_v(M),\quad X\in
C^\infty(B,TB).
\]
Using Proposition \ref{VertCohomBundleDef} and Theorem
\ref{FiberedHodge}, we can explicitly describe the isomorphism
$\sH_v^\bullet(M)\cong H_v^\bullet(M)$ in terms of sections by
the composition
\[
C^\infty\big(B,\sH_v^\bullet(M)\big) = \ker
D_v\cap\gO^\bullet_v(M)\hookrightarrow \ker
d_v\cap\gO^\bullet_v(M) \xrightarrow{\Psi}
C^\infty\big(B,H_v^\bullet(M)\big).
\]
This implies that $\nabla_X^{H_v} \Psi(\go)\in
C^\infty\big(B,H_v^\bullet(M)\big)$ corresponds to
\[
P_{\ker D_v}\circ \imu(X^h)\circ d_h\big(P_{\ker D_v}\go\big)\in
C^\infty\big(B,\sH_v^\bullet(M)\big).
\]
Finally, Proposition \ref{dSplit} shows that on $\gO^\bullet_v(M)$
\[
\imu(X^h)\circ d_h = \widetilde\nabla_{X^h}^v,
\]
from which we obtain the claimed formula.
\end{proof}

\noindent\textbf{Metrics on the Bundle of Vertical Cohomology
Groups.} The $C^\infty(B)$-module $\gO^\bullet_v(M,E)$ is endowed
with the pairing
\begin{equation}\label{InfiniteMetric}
(\go,\eta)_{M/B}:= \int_{M/B} \scalar{\go}{\eta}\vol_F(g_v)\in
C^\infty(B) ,\quad \go,\eta\in \gO^\bullet_v(M,E),
\end{equation}
where the scalar product in the integrand is induced by $g_v$
together with the Hermitian metric on $E$.

\begin{dfn}\label{CohomBundleMet}
Let $g_v$ be a vertical metric. We define
\[
\scalar{\go}{\eta}_{\sH_{A,v}}:= (\go,\eta)_{M/B},\quad
\go,\eta\in C^\infty\big(B,\sH_{A,v}^\bullet(M)\big),
\]
and, if $\tau_v$ is the vertical chirality operator,
\[
Q_{A,v}(\go,\eta):= \scalar{\go}{\tau_v\eta}_{\sH_{A,v}}.
\]
We also use $\scalar{.}{.}_{\sH_{A,v}}$ and $Q_{A,v}$ to the
corresponding objects induced by the vertical Hodge-de-Rham
isomorphism on $H_{A,v}^\bullet(M)$.
\end{dfn}


Clearly, $\scalar{.}{.}_{\sH_{A,v}}$ is a Hermitian metric on the
vector bundle $H_{A,v}^\bullet(M)\to B$, which through the
vertical Hodge-de-Rham isomorphism depends on the vertical metric
$g_v$. In contrast, $Q_{A,v}$ is independent of $g_v$ since it is
related to the vertical intersection form via
\begin{equation}\label{VerticalIntersection}
Q_{A,v}\big([\ga],[\gb]\big) = i^k \int_{M/B} \langle
\ga\wedge\gb\rangle,\quad [\ga],[\gb]\in H_{A,v}^\bullet(M),
\end{equation}
where $k$ depends only on the degrees of $\ga$ and $\gb$.
Furthermore, one easily checks that $Q_{A,v}$ is an indefinite
Hermitian form with signature
\[
\Sign(Q_{A,v})= \rk\big(\sH_{A,v}^+(M)\big) -
\rk\big(\sH_{A,v}^-(M)\big).
\]
Here, $\sH_{A,v}^\pm(M)$ denotes the $\pm 1$ eigenbundle of
$\tau_v$. This implies that $Q_{A,v}$ has signature 0 unless the
dimension of the fiber is divisible by 4 in which case
$\Sign(Q_{A,v})=\Sign(F)$.

\begin{prop}\label{FlatSignBundle}
The flat connection $\nabla^{H_{A,v}}$ is compatible with the
indefinite Hermitian metric $Q_{A,v}$. It is compatible with the
Hermitian metric $\scalar{.}{.}_{\sH_{A,v}}$ if and only if for
all $X\in C^\infty(M,T^hM)$
\begin{equation}\label{FlatSignBundle:1}
2 P_{\ker D_v}\circ B(X)\circ P_{\ker D_v}  +   k_v(X)  =0,
\end{equation}
where $B(X)$ is the tensor as in \eqref{NablaNatInv}, and $k_v$ is
the mean curvature form.
\end{prop}

\begin{proof}
For the first part we use the description
\eqref{VerticalIntersection} for $Q_{A,v}$. Let $\ga,\gb\in
\gO^\bullet_v(M,E)$ be chosen in such a way that $\langle
\ga\wedge\gb\rangle$ is of maximal vertical degree. Then
\[
d_B \int_{M/B} \langle \ga\wedge\gb\rangle = \int_{M/B} d_h
\langle \ga\wedge\gb\rangle = \int_{M/B} \langle
d_{A,h}\ga\wedge\gb\rangle + (-1)^{|\ga|} \langle \ga\wedge
d_{A,h}\gb\rangle.
\]
Hence, if $d_{A,v}\ga = d_{A,v}\gb =0$ we have
\[
d_B Q_{A,v}\big([\ga],[\gb]\big) =
Q_{A,v}\big([d_{A,h}\ga],[\gb]\big) + (-1)^{|\ga|}
Q_{A,v}\big([\ga],[d_{A,h}\gb]\big),
\]
so that, according to Definition \ref{FlatVertCohom},
\[
X Q_{A,v}\big([\ga],[\gb]\big) = Q_{A,v}\big(\nabla_X^{H_{A,v}}
[\ga],[\gb]\big) + Q_{A,v}\big([\ga],\nabla_X^{H_{A,v}}
[\gb]\big),\quad X\in C^\infty(B,TB).
\]
This shows that $\nabla^{H_{A,v}}$ is indeed compatible with
$Q_{A,v}$. Now let $g_v$ be a vertical metric, and let $\go,
\eta\in C^\infty\big(B,\sH_{A,v}^\bullet(M)\big)$. Then
\[
d_B \int_{M/B} \scalar{\go}{\eta} \vol_F(g_v) = \int_{M/B}
d_h\big(\scalar{\go}{\eta}\big)\wedge \vol_F(g_v) + \int_{M/B}
\scalar{\go}{\eta}\, d_h\vol_F(g_v).
\]
It follows from Proposition \ref{Rummler} that
\[
d_h\vol_F(g_v) = k_v\wedge \vol_F(g_v).
\]
Since the connection $\nabla^{A,v}$ is compatible with the metric
on $\gO^\bullet_v(M,E)$ we know that for all $X\in C^\infty(B,TB)$
\[
X^h \scalar{\go}{\eta} = \Scalar{\nabla^{A,v}_{X^h}\go}{\eta} +
\Scalar{\go}{\nabla^{A,v}_{X^h}\eta}.
\]
From \eqref{NablaNat} and \eqref{NablaNatInv} we then deduce
\[
X^h \scalar{\go}{\eta} =
\Scalar{\widetilde\nabla^{A,v}_{X^h}\go}{\eta} +
\Scalar{\go}{\widetilde \nabla^{A,v}_{X^h}\eta} +
\Scalar{B(X^h)\go}{\eta} + \Scalar{\go}{B(X^h)\eta}.
\]
Now, $B(X^h)$ is easily seen to be self-adjoint with respect to
the metric $\scalar{.}{.}$ on $\gO^\bullet_v(M,E)$. Putting all
pieces together, we find that for the connection
$\nabla_X^{\sH_{A,v}}$ of Proposition \ref{FlatConnHodge}
\[
\begin{split}
X \scalar{\go}{\eta}_{\sH_{A,v}} &=
\Scalar{\nabla_X^{\sH_{A,v}}\go}{\eta}_{\sH_{A,v}} +
\Scalar{\go}{\nabla_X^{\sH_{A,v}}\eta}_{\sH_{A,v}}\\& \qquad +
\Scalar{\go}{2B(X^h)\eta}_{\sH_{A,v}} +
\Scalar{\go}{k_v(X^h)\eta}_{\sH_{A,v}},
\end{split}
\]
which proves that $\nabla^{\sH_{A,v}}$ is compatible with
$\scalar{.}{.}_{\sH_{A,v}}$ if and only if
\eqref{FlatSignBundle:1} holds. Since the metrics as well as the
connections $\nabla^{H_{A,v}}$ and $\nabla^{\sH_{A,v}}$ on
$H^\bullet_{A,v}(M)$ and $\sH^\bullet_{A,v}(M)$ coincide under
the vertical Hodge-de-Rham isomorphism, the proof of Proposition
\ref{FlatSignBundle} is finished.
\end{proof}

\begin{remark*}
If we denote by $p$ and $q$ the maximal ranks of subbundles of
$H_{A,v}^\bullet(M)$ on which $Q_{A,v}$ is positive respectively
negative definite, we can rephrase the first part of Proposition
\ref{FlatSignBundle} by saying that $\nabla^{H_{A,v}}$ is a flat
$\U(p,q)$-connection.\footnote{Recall that $\U(p,q)$ denotes the
isometry group of the quadratic from
\[
\sum_{j=1}^p |z_j|^2 - \sum_{j=p+1}^{p+q} |z_j|^2.
\]} The choice of a vertical metric
reduces the structure group of $H_{A,v}^\bullet(M)$ to the
subgroup $\U(p)\times \U(q)$. However, the connection does not
necessarily reduce to a flat $\U(p)\times \U(q)$-connection, the
geometric obstruction being \eqref{FlatSignBundle:1}. As we have
seen in Theorem \ref{VanishingMeanCurv} we can always arrange
that the mean curvature form vanishes. For arbitrary fiber
bundles, the tensor $B(X)$ is, however, a non-trivial
obstruction. It would be interesting to find a topological
condition which guarantees that there exists a vertical metric
such that \eqref{FlatSignBundle:1} holds.
\end{remark*}

\begin{dfn}\label{CohomTwistSignOp}
Let $D_{A,v}$ and $D_{A,h}$ be the vertical and horizontal de
Rham operators as in Definition \ref{PartDeRhamOpDef}. If $\dim
M$ is odd, we define the odd signature operator on $B$ with
\emph{values in the bundle of vertical cohomology groups},
\[
D_B\otimes
\nabla^{\sH_{A,v}}:\gO^\bullet\big(B,\sH_{A,v}^\bullet(M)\big) \to
\gO^\bullet\big(B,\sH_{A,v}^\bullet(M)\big),
\]
by
\[
D_B\otimes \nabla^{\sH_{A,v}}:= P_{\ker D_{A,v}}\circ \tau_M
D_{A,h}\circ P_{\ker D_{A,v}}.
\]
Here, $\tau_M$ is the chirality operator associated to a fixed
submersion metric on $M$.
\end{dfn}

\begin{remark}\label{CohomTwistSignOpRem}\quad\nopagebreak
\begin{enumerate}
\item Certainly, $D_B\otimes \nabla^{\sH_{A,v}}$ is a formally
self-adjoint elliptic differential operator and thus has a
well-defined Eta invariant. This will play an important role in
Dai's general adiabatic limit formula for the Eta in Section
\ref{GenAdiabaticLimit}.
\item We note that if $\dim B$ is odd, and \eqref{FlatSignBundle:1}
is satisfied, then $D_B\otimes \nabla^{\sH_{A,v}}$ is actually
isometric to two copies of the odd signature operator on $B$
twisted by $\nabla^{\sH_{A,v}}$. This is because we have not
restricted to forms on the base of even degree, compare with
Remark \ref{OddSignRem} (i).
\end{enumerate}
\end{remark}

\subsection[Eta Invariants of U$(p,q)$-Connections]{Eta
Invariants of
$\boldsymbol{\U(p,q)}$-Connections}\label{IndefiniteMetric}

Before we continue with the general discussion, we briefly want
to digress on the Eta invariant of the operator $D_B\otimes
\nabla^{\sH_{A,v}}$ introduced above. Without any effort, we can
treat the more general case that $E\to B$ is a complex vector
bundle, endowed with an indefinite Hermitian metric $Q$ and a
connection $\nabla$, not necessarily flat, but compatible with
$Q$. We choose a splitting $E=E^+\oplus E^-$ into subbundles
where $Q$ is positive respectively negative definite, and define
$\tau_E$ to be $\pm\id_E$ on $E^\pm$. Then we can define a
Hermitian metric on $E$ via
\[
h(e,f):= Q(e,\tau f),\quad e,f\in E,
\]
compare with Definition \ref{CohomBundleMet}. Note that the
splitting $E=E^+\oplus E^-$ is orthogonal with respect to $h$. We
now define an $\End(E)$-valued 1-form on $B$
\[
\go^{\nabla,\tau_E}(X):= \tau_E \big[\nabla_X,\tau_E\big] =
\tau_E\circ \nabla_X\circ \tau_E - \nabla_X,\quad X\in
C^\infty(B,TB).
\]
Then we have the following simple result.

\begin{lemma}\label{NablaCompPart}
For all $X\in C^\infty(B,TB)$, the endomorphism
$\go^{\nabla,\tau_E}(X)$ is self-adjoint with respect to $h$. It
interchanges the subbundles $E^+$ and $E^-$. Moreover, the
connection
\[
\nabla^u:= \nabla +\lfrac 12 \go^{\nabla,\tau_E},
\]
is unitary with respect to $h$.
\end{lemma}

\begin{proof}
Let $e,f\in C^\infty(B,E)$. Since $\nabla_X$ is compatible with
$Q$, one verifies---using in particular that $\tau^2_E=\id_E$ and
that $Q\circ \tau_E = Q$,
\[
\begin{split}
h\big(\go^{\nabla,\tau_E}(X)e,f\big) &=
Q\big(\tau_E\big[\nabla_X,\tau_E\big]e,\tau_E f\big) =
Q\big(\nabla_X (\tau_E e),f\big) - Q\big(\nabla_Xe,\tau_E f\big)\\
&= -Q\big(\tau_Ee,\nabla_X f \big) + Q\big(e,\nabla_X(\tau_E
f)\big) = Q\big(e,\big[\nabla_X,\tau_E\big]f\big)\\
&= h\big(e,\go^{\nabla,\tau_E}(X)f\big).
\end{split}
\]
Hence, $\go^{\nabla,\tau_E}(X)$ is self-adjoint with respect to
$h$. Now, let $P_{E^\pm} := \lfrac 12 (\id_E\pm \tau_E)$ denote
the projection onto $E^\pm$. Then one easily obtains that
\[
P_{E^+}\circ \nabla_X \circ P_{E^-} + P_{E^-}\circ \nabla_X\circ
P_{E^+}= -\lfrac12 \go^{\nabla,\tau_E}(X).
\]
On the one hand, this implies that $\go^{\nabla,\tau_E}(X)$
interchanges the subbundles $E^+$ and $E^-$. By definition of
$\nabla^u$, we can deduce on the other hand, that $\nabla^u$
preserves $E^+$ and $E^-$ from which it easily follows that
$\nabla^u$ is unitary with respect to $h$.
\end{proof}

\begin{remark*}
In the case that $E=\sH_{A,v}^\bullet(M)$ is the bundle of
vertical cohomology groups, $Q=Q_v$ is the vertical intersection
form and $\nabla=\nabla^{H_{A,v}}$ is the natural flat
connection, the 1-form $\go^{\nabla,\tau_E}$ is precisely the
1-form appearing in \eqref{FlatSignBundle:1}, compare also with
\eqref{LieTauComm}. This gives a more abstract explanation of
Proposition \ref{FlatSignBundle}.
\end{remark*}

\noindent\textbf{The Odd Signature Operator with values in
$\boldsymbol{E}$.} To define the analog of $D_B\otimes
\nabla^{\sH_{A,v}}$ in the case at hand, we choose a metric $g_B$
on $B$, and let $\tau_B$ the associated chirality operator on
$\gO^\bullet(B,E)$. Let $b:=\dim B$, and extend $\tau_E$ to
$\gO^\bullet(B,E)$ by requiring that
\begin{equation}\label{TauEDef}
\tau_E (\ga\otimes e) = (-1)^{p (b+1)} \ga\otimes \tau_E e,\quad
\ga\in \gO^p(B),\quad e\in C^\infty(B,E).
\end{equation}
Checking signs one finds that $\tau_B\tau_E = \tau_E\tau_B$. We
then define
\[
\tau:= \tau_B\tau_E: \gO^\bullet(B,E)\to \gO^{b-\bullet}(B,E),
\]
which takes the place of the total chirality operator $\tau_M$,
compare with Lemma \ref{TauSplit}. Note that more explicitly, if
$\ga\otimes e\in \gO^p(B,E)$, then
\[
\tau(\ga\otimes e) =  \tau_B\big((-1)^{p (b+1)} \ga\otimes \tau_E
e\big)  = (-1)^{p (b+1)} (\tau_B\ga)\otimes \tau_E e.
\]
Then the analog of $D_B\otimes \nabla^{\sH_{A,v}}$ is given by
\begin{equation}\label{CohomTwistSignOp:Gen}
D_B\otimes \nabla:= \tau d_{\nabla} + d_{\nabla} \tau:
\gO^\bullet(B,E)\to \gO^\bullet(B,E),
\end{equation}
where $d_\nabla$ is the exterior differential on $B$ twisted by
the connection $\nabla$ on $E$. We also define
\[
D_B\otimes \nabla^u := \tau d_{\nabla^u} + d_{\nabla^u} \tau,
\]
and denote by $\nabla^{u,\pm}$ the restriction of $\nabla^u$ to
$E^\pm$.

\begin{lemma}\label{SignCompPart}
With respect to the splitting $E=E^+\oplus E^-$, the operator
$D_B\otimes \nabla^u$ is of the form
\begin{equation*}
D_B\otimes \nabla^u = \begin{pmatrix} D_B^{\nabla^{u,+}} &0\\ 0 &
-D_B^{\nabla^{u,-}}
\end{pmatrix},
\end{equation*}
where
\[
D_B^{\nabla^{u,\pm}}:= \tau_B (-1)^{b+1} d_{\nabla^{u,\pm}} +
d_{\nabla^{u,\pm}} \tau_B,\quad b:=\dim B.
\]
Moreover, if we define
\[
V:= D_B\otimes \nabla - D_B\otimes \nabla^u,
\]
then $V$ is a self-adjoint operator on $\gO^\bullet(B,E)$ of
order 0 which interchanges $\gO^\bullet(B,E^+)$ and
$\gO^\bullet(B,E^-)$. In particular, $D_B\otimes \nabla$ and
$D_B\otimes \nabla^u$ are formally self-adjoint.
\end{lemma}

\begin{proof}
It follows from Lemma \ref{NablaCompPart} that $\tau_E$ commutes
with $\nabla^u$. Using the sign convention in \eqref{TauEDef} it
is immediate that
\[
\tau_E d_{\nabla^u} = (-1)^{b+1} d_{\nabla^u} \tau_E.
\]
Hence,
\[
\tau d_{\nabla^u} + d_{\nabla^u} \tau = \big((-1)^{b+1}\tau_B
d_{\nabla^u} + d_{\nabla^u} \tau_B\big)\tau_E,
\]
which proves the first assertion. The other assertions are a
simple consequence of the corresponding properties of
$\go^{\nabla,\tau_E}$ in Lemma \ref{NablaCompPart}, since by
definition
\[
d_{\nabla^u} = d_\nabla +\lfrac 12
\emu\big(\go^{\nabla,\tau_E}\big),
\]
where $\emu(.)$ is exterior multiplication.
\end{proof}

\noindent\textbf{Difference of Eta Invariants.} Roughly, Lemma
\ref{SignCompPart} asserts that $D\otimes \nabla$ is the direct
sum of two geometric Dirac operator plus a lower order
perturbation which interchanges the twisting bundles. This leads
to a simple relation between the Eta invariants of $D\otimes
\nabla$ and $D\otimes \nabla^u$. The following result is a
reformulation of \cite[Thm.'s 2.7 \& 2.35]{BC89}. We formulate it
in terms of the $\xi$-invariant, see Definition \ref{EtaDef}.

\begin{theorem}\label{CohomTwistEtaDiff}
As before, let $V=D_B\otimes \nabla - D_B\otimes \nabla^u$. Then
\[
\xi(D_B\otimes \nabla) - \xi (D_B\otimes \nabla^u) = \SF
\big(D_B\otimes \nabla^u+ xV\big)_{x\in[0,1]}.
\]
\end{theorem}

\begin{remark*}
By comparison with the variation formula of Corollary
\ref{EtaDiffSFApp}, we see that Theorem \ref{CohomTwistEtaDiff}
asserts that the contribution coming from the variation of the
reduced $\xi$-invariant vanishes, i.e.,
\[
\int_0^1 \lfrac d{dx}[\xi(D_x)]dx =0.
\]
In fact, this is precisely what Bismut and Cheeger prove, see
\cite[Lem. 2.11]{BC89}. Recall from Proposition
\ref{RedEtaDerApp} that
\[
\lfrac d{dt}[\xi(D_x)]dx =  - \lfrac1{\sqrt\pi}\, a_n(V,D_x^2),
\]
where $a_n(V,D_x^2)$ is the constant term in the asymptotic
expansion of
\[
\sqrt t\,\Tr\big(Ve^{-tD_x^2}\big),\quad\text{as }t\to 0.
\]
In the case that $\dim B$ is even, Theorem \ref{HeatTrace} shows
that there are no half integer powers of $t$ in the asymptotic
expansion of $\Tr\big(Ve^{-tD_x^2}\big)$, so that Theorem
\ref{CohomTwistEtaDiff} is a consequence of the general theory
for elliptic operators. However, the odd dimensional case
requires considerably more work. In \cite[Lem. 2.11]{BC89},
Bismut and Cheeger prove the corresponding result for operators
of the form we are considering here. Their proof uses Getzler's
local index theory techniques for twisted Dirac operators,
adapted to odd dimensional base spaces, in a similar way as we
have described in Section \ref{VarLocalInd}.
\end{remark*}

\section{Elements of Bismut's Local Index Theory for
Families}\label{LocalFamIndex}

To discuss Dai's adiabatic limit formula, we need to recall some
aspects of Bismut's local index theory for families. We will be
rather sketchy and refer to the original article \cite{B86} and
the treatment in \cite[Ch.'s 9 \& 10]{BGV} for more details. The
survey article \cite{B98} is also recommended. For convenience
and since we will not need a greater generality, we restrict to
the case of the signature operator.

\subsection{The Index Theorem for Families}

The predecessor of local index theorem for families is the
$K$-theoretic version by Atiyah and Singer \cite{AS4}, which we
briefly recall. As announced we consider only the case of the
signature operator. Let $F\hookrightarrow M\xrightarrow{\pi}B$ be
an oriented fiber bundle of closed manifolds, where $F$ is assumed
to be even dimensional. We choose a vertical projection and a
vertical metric $g_v$. Let $\nabla^v$ be the associated
connection on $T^vM$, and let $D_v^+$ be the vertical signature
operator defined by the vertical chirality operator. As in
Proposition \ref{VertDeRhamKer} we can view $\ker D^+_v$ and
$\coker D^+_v$ as (spaces of sections of) finite dimensional
vector bundles over $B$.

\begin{dfn}\label{IndexBundle}
The \emph{index bundle} associated to $D_v^+$ is defined by
\[
\Ind D^+_v := [\ker D^+_v] - [\coker D^+_v] \in K^0(B).
\]
\end{dfn}
Note that as we are considering only families of signature
operators we do not need the beautiful construction for the case
of varying dimensions as in \cite{AS4}. The Chern character
defines a map in $K$-theory,
\[
\ch: K^0(B)\otimes \C \to H^{\ev}(B).
\]
Then cohomological version of the families index theorem as in
\cite[Thm. 5.1]{AS4} is

\begin{theorem}[Atiyah-Singer]\label{ASFam}
The Chern character of the index bundle associated to the
signature operator is given by
\[
\ch(\Ind D^+_v) = \Big[\int_{M/B}L(T^vM,\nabla^v)\Big] \in
H^{\ev}(B),
\]
where $L(T^vM,\nabla^v)$ is the Hirzebruch $L$-form of $T^vM$
defined via Chern-Weil theory as in \eqref{L-L-Hat} in terms of
$\nabla^v$.
\end{theorem}

\subsection{Superconnections and Associated Dirac Operators}

Quillen \cite{Qui85b} introduced superconnections to study
Chern-Weil theory for the Chern character of a difference bundle.
We briefly recall the basic definitions, and refer to \cite[Sec.
1.4]{BGV} and \cite{Qui85b} for details. Let $B$ be a closed,
oriented manifold, and let $E\to B$ be a complex vector bundle.

\begin{dfn}\label{SuperConnection}
A differential operator $\A$ on $\gO^\bullet(B,E)$ is called a
\emph{generalized connection} on $E$ if it satisfies the Leibniz
rule
\[
\A(\ga\wedge \gb)= d\ga \wedge \gb + (-1)^{|\ga|} \ga\wedge\A\gb,
\]
where $\ga\in\gO^\bullet(B)$ and $\gb\in\gO^\bullet(B,E)$. The
\emph{curvature} of $\A$ is defined as
\[
\A^2\in \gO^\bullet\big(B,\End(E)\big).
\]
If $E$ is $\Z_2$-graded and $\A$ is of odd parity with respect to
the total grading on $\gO^\bullet(B,E)$, then $\A$ is called a
\emph{superconnection}.
\end{dfn}

\begin{remark*}\quad\nopagebreak
\begin{enumerate}
\item If $E=E^+\oplus E^-$ is $\Z_2$-graded, the total grading of
$\gO^\bullet(B,E)$ referred to above is defined by
\[
\gO(B,E)^\pm:= \gO^{\ev}(B,E^\pm) + \gO^{\odd}(B,E^\mp).
\]
This should not be confused with the grading $\gO^\pm(B,E)$,
induced by the chirality operator $\tau_B$.
\item The fact that the curvature $\A^2$ is indeed given by the
action of an element in $\gO^\bullet\big(B,\End(E)\big)$ works as
in the case of a usual connection, see \cite[Prop. 1.38]{BGV}.
\item A generalized connection $\A$ is determined by its homogeneous
components
\[
\A = \A_{[0]} + \A_{[1]} + \A_{[2]} + \ldots,
\]
where $\A_{[p]}\in \gO^p\big(B,\End(E)\big)$ for $p\neq 1$ and
$\A_{[1]}$ is a connection on $E$. In the case that $\A$ is a
superconnection, the connection part $\A_{[1]}$ preserves the
splitting $E=E^+\oplus E^-$.
\end{enumerate}
\end{remark*}

\begin{dfn}\label{ChernSuperConn}
Let $E$ be $\Z_2$-graded, and let $\A$ be a superconnection on
$E$. Then we define the \emph{Chern character form} of $\A$ as
\[
\ch_s(E,\A):= \str_E\big(\gamma\exp(-\A^2)\big).
\]
Here, $\gamma:\gO^\bullet(B)\to \gO^\bullet(B)$ is a normalization
function, defined for forms of homogeneous degree by
\begin{equation}\label{NormSuperConnDef}
\gamma(\ga):=\big(\lfrac1{\sqrt{2\pi i}}\big)^{|\ga|}
\ga,\quad\text{where } \sqrt i= e^{\frac{i\pi}4}.
\end{equation}
\end{dfn}

The discussion in Appendix \ref{CharClass} generalizes to the
case of superconnections. In particular, $\ch_s(E,\A)$ is a
closed differential from on $B$ whose cohomology class is
independent of the superconnection $\A$. Since the supertrace
vanishes on endomorphisms of odd parity, one can also check that
$\ch_s(E,\A) \in \gO^{\ev}(B,E)$. We also note that if
$\A=\nabla$ is a connection in the usual sense, which decomposes
with respect to the splitting $E=E^+\oplus E^-$ into
$\nabla=\nabla^+\oplus \nabla^-$, then
\begin{equation}\label{ChernSuperConnRem}
\ch_s(E,\nabla) = \ch(E^+,\nabla^+) - \ch(E^-,\nabla^-),
\end{equation}
where the right hand side is as in Definition \ref{ChernCharDef}.
Then the main idea of \cite{Qui85b} can be summarized as

\begin{theorem}[Quillen]
Let $E\to B$ be a $\Z_2$-graded Hermitian bundle over $B$. Let
$[E^+]-[E^-]$ be the induced element in $K^0(B)$, and let $\A$ be
a superconnection on $E$. Then
\[
\ch\big([E^+]-[E^-]\big) = \big[\ch_s(E,\A)\big]\in H^{\ev}(B).
\]
\end{theorem}

\noindent\textbf{Generalized Clifford Connections and Dirac
Operators.} Whenever $E$ is endowed with the structure of a
Clifford module, one can associate a Dirac operator to a
generalized connection $\A$. Let
\[
\boldsymbol{\gs}^{-1}: \gL^\bullet T^*B\to \cl(T^*B)
\]
be the quantization map. If $c:\cl(T^*B)\to \End(E)$ denotes
Clifford multiplication, we get a natural Clifford contraction
\[
\gO^\bullet(B,E)\xrightarrow{\boldsymbol{\gs}^{-1}}
C^\infty\big(B,\cl(T^*B)\otimes E\big)\xrightarrow{c}
C^\infty(B,E).
\]
Since $\A$ maps $C^\infty(B,E)$ to $\gO^\bullet(B,E)$, we can
define
\begin{equation}\label{DiracSuperConn}
D_\A:=c\circ\boldsymbol{\gs}^{-1}\circ \A: C^\infty(B,E)\to
C^\infty(B,E).
\end{equation}
Clearly, this defines an elliptic operator of first order. In
order for $D_\A$ to be formally self-adjoint, we certainly have
to require that the connection part $\A_{[1]}$ of $\A$ is a
Clifford connection in sense of Definition \ref{GeomDiracOp}.
Furthermore, some condition has to be imposed on the other
homogeneous components of $\A$, which we derive now, see also
\cite[p. 117]{BGV}. For $p\neq 1$ and a local orthonormal frame
$\{f_a\}$ for $TB$ we can write locally
\[
A_{[p]} = \lfrac 1{p!}f^{a_1}\wedge\ldots\wedge f^{a_p}\wedge
T_{a_1\ldots a_p},\quad  \text{with}\quad T_{a_1\ldots a_p}\in
C^\infty\big(B,\End(E)\big).
\]
The contribution to $D_\A$ is then given by
\[
c\circ\boldsymbol{\gs}^{-1}\circ A_{[p]} = \lfrac
1{p!}c^{a_1}\ldots c^{a_p}T_{a_1\ldots a_p},
\]
where $c^{a_j}$ is short for Clifford multiplication with
$f^{a_j}$. For $e,\widetilde e \in C^\infty(B,E)$ one computes
that
\[
\begin{split}
\Scalar{c^{a_1}\ldots c^{a_p}T_{a_1\ldots a_p} e}{\widetilde e} &=
(-1)^p \Scalar{e}{ T_{a_1\ldots a_p}^* c^{a_p}\ldots
c^{a_1}\widetilde e }\\ &= (-1)^{\lfrac {p(p+1)}{2}}
\Scalar{e}{T_{a_1\ldots a_p}^* c^{a_1}\ldots c^{a_p}\widetilde e}.
\end{split}
\]
This motivates the following

\begin{dfn}\label{CliffordSuperConn}
Let $E$ be a Hermitian vector bundle. A generalized connection
$\A$ is called \emph{unitary} if its connection part is a unitary
connection and if for $p\neq 1$
\[
\A_{[p]}^* = (-1)^{\lfrac {p(p+1)}{2}} \A_{[p]}.
\]
Here, taking the adjoint is meant with respect to the endomorphism
part only. $\A$ is is called a \emph{generalized Clifford
connection}, if in addition, its connection part $\A_{[1]}$ is a
Clifford connection, and if for $p\neq 1$ and $\xi\in \gO^1(B)$
\begin{equation*}
\A_{[p]} c(\xi) = (-1)^p c(\xi)\A_{[p]},
\end{equation*}
where again the product is to be understood in the endomorphism
part.
\end{dfn}

The essential part of the following result is a consequence of the
discussion preceding Definition \ref{CliffordSuperConn}. The
other claims are immediate.

\begin{prop}
Let $E\to B$ be a Hermitian vector bundle endowed with a Clifford
structure, and let $\A$ be a generalized Clifford connection. Then
$D_\A$ is formally self-adjoint and $D_\A^2$ is a generalized
Laplacian. The symbol of $\A$ is given by Clifford multiplication
$c:T^*B\otimes E\to E$. If in addition $E$ is $\Z_2$-graded and
$\A$ is a superconnection, then $D_\A$ is $\Z_2$-graded.
\end{prop}

\begin{remark*}
As we have pointed out in Remark \ref{GeomDiracOpRem}, not every
Dirac operator arises as a geometric Dirac operator associated to
a Clifford connection. However, it is shown in \cite[Prop.
3.42]{BGV} that there is a 1-1 correspondence between Clifford
superconnections and $\Z_2$-graded Dirac operators with symbol
being the given Clifford structure. Going through the proof of
loc.cit. one sees that the same statement is true for ungraded
Dirac operators and generalized Clifford connections as defined
above. Note, however, that in contrast to \cite{BGV}, we require
Dirac operators to be formally self-adjoint.
\end{remark*}

\subsection{The Families Index Theorem for the Signature
Operator}\label{LocFamIndSign}

Generalizing Quillen's construction to infinite dimensional
bundles, Bismut \cite{B86} found a heat equation formula for the
Chern character form of the index bundle. We describe the setup
briefly, again restricting to the case of the untwisted signature
operator.\\

\noindent\textbf{Bismut's Superconnection.} Let $F\hookrightarrow
M\xrightarrow{\pi}B$ be an oriented fiber bundle of closed
manifolds. We choose a vertical projection, and let
$\gO_v^\bullet(M)$ be the $C^\infty(B)$-module of vertical
differential forms. We formally interpret this as the space of
sections of an infinite dimensional bundle $\sE$ over $B$, where
the fiber $\sE_y$ over $y\in B$ is given by the space of
differential forms over $\pi^{-1}(y)$. Since this picture has
only motivational character, we do not give any details of how
this bundle of Fr\'{e}ch\'{e}t spaces is defined rigorously. We
can then view the space of all differential forms $\gO^\bullet(M)$
as
\[
\gO^\bullet(M) \cong \gO^\bullet(B,\sE),
\]
compare with \eqref{FormsSplit}. Proposition \ref{dSplit} shows
that the total exterior differential $d_M$ on $\gO^\bullet(M)$
splits as
\begin{equation}\label{dSplitSuperConn}
d_M  = d_v + d_h + \imu(\gO),
\end{equation}
where with respect to any choice of metric $g_B$ in a local
orthonormal frame $\{f_a\}$ for $TB$
\begin{equation}\label{dSplitSuperConn:2}
d_h = f^a\wedge \widetilde \nabla^\oplus_a,\quad \imu(\gO) =
\lfrac 12 f^a\wedge f^b\wedge \imu(\gO_{ab}).
\end{equation}
Recall that when restricted to $\gO^\bullet_v(M)$, the connection
$\widetilde\nabla^\oplus_a$ is defined as $\widetilde \nabla^v_a$,
see Definition \ref{CanVerticalConn2}. We view the latter as a
natural connection $\nabla^\sE$ on the infinite dimensional
bundle $\sE$. Then \eqref{dSplitSuperConn} and
\eqref{dSplitSuperConn:2} express $d_M$ as a generalized
connection on $\sE$ with connection part $\nabla^{\sE}$. It is a
superconnection with respect to the even/odd grading on
$\gO^\bullet_v(M)$. In this interpretation, the property
$d_M^2=0$ states that $d_M$ is a \emph{flat} superconnection on
the bundle of vertical differential forms, an observation which
is due to \cite[Sec. III (b)]{BL95}. For this reason, $d_M$
together with its interpretation as a superconnection is
sometimes called the \emph{Bismut-Lott superconnection}.

\begin{dfn}\label{BismutSuperconnDef}
Assume that $F\hookrightarrow M\xrightarrow{\pi} B$ is endowed
with a vertical metric $g_v$ and a vertical projection. Let
$\nabla^v$ be the associated canonical connection, and let $D_v$
be the vertical de Rham operator. With respect to any choice of
$g_B$ and a local orthonormal frame $\{f_a\}$ on $B$ define
\begin{equation}\label{BismutSuperconnConn}
\nabla^{\sE,u} := f^a\wedge \big(\nabla^v_a + \lfrac 12
k_v(f_a)\big): \gO^\bullet_v(M)\to \gO^{1,\bullet}(M),
\end{equation}
where $k_v$ is the mean curvature form. Then, the \emph{Bismut
superconnection} is defined as
\[
\B:= \lfrac 12 D_v + \nabla^{\sE,u} - \lfrac12
c_v(\gO):\gO^\bullet_v(M)\to \gO^\bullet(M),
\]
where $c_v$ denotes the vertical Clifford multiplication on
$\gO^\bullet(M)$, and locally
\[
c_v(\gO) = \lfrac 12 f^a\wedge f^b\wedge c_v(\gO_{ab}).
\]
\end{dfn}

The Bismut superconnection naturally extends to an operator
$\gO^\bullet(M)\to \gO^\bullet(M)$, if we replace $\nabla^v$ with
$\nabla^\oplus$ in \eqref{BismutSuperconnConn}. For this note that
for $\ga\in \gO^p(B)$ and $\go\in \gO^\bullet_v(M)$,
\[
f^a\wedge \nabla^\oplus_a\big((\pi^*\ga)\wedge \go)\big) =
\pi^*(d_B\ga)\wedge \go + (-1)^p \pi^*\ga\wedge f^a\wedge
\nabla^v_a \go.
\]
This relation also shows that replacing $\nabla^v$ with
$\nabla^\oplus$ is the same as extending $\B$ from
$\gO^\bullet_v(M)=C^\infty(B,\sE)$ to $\gO^\bullet(M) =
\gO^\bullet(B,\sE)$ by requiring the Leibniz rule. Moreover, we
find that the term $f^a\wedge \nabla^\oplus_a$ is independent of
the chosen metric $g_B$, since this is true for the connection
$\nabla^v$, see Proposition \ref{CanVerticalConn}.

\begin{remark}\label{BismutUnitaryRem}
We want to point out that the definition of $\B$ can be motivated
by an infinite dimensional version of Lemma \ref{NablaCompPart}.
If we choose a vertical metric $g_v$ on $T^vM$, we can view the
pairing $(.,.)_{M/B}$ in \eqref{InfiniteMetric} as a metric on
the bundle $\sE$. As in Proposition \ref{FlatSignBundle}, the
connection part $\nabla^\sE=f^a\wedge \widetilde \nabla^v_a$ of
$d_M$ is compatible with the vertical intersection pairing, but
not necessarily with $(.,.)_{M/B}$. If we proceed as in Lemma
\ref{NablaCompPart}---using in particular \eqref{NablaNat} and
\eqref{LieTauComm}---we see that the unitary connection
associated to $\nabla^\sE$ is given by
\[
\begin{split}
f^a\wedge \widetilde \nabla^v_a + \lfrac 12 f^a\wedge \tau_v
\big[\widetilde \nabla^v_a,\tau_v\big] &= f^a\wedge
\big(\widetilde \nabla^v_a + B(f_a) +   \lfrac 12 k_v(f_a)\big)\\
&= f^a\wedge \big(\nabla^v_a + \lfrac 12 k_v(f_a)\big).
\end{split}
\]
This is precisely the connection part $\nabla^{\sE,u}$ of the
Bismut superconnection. In order to get the unitary
superconnection associated to $d_M$ we proceed as in Definition
\ref{CliffordSuperConn} and replace the other homogeneous
components $d_v$ and $\imu(\gO)$ of $d_M$ with
\[
\lfrac 12 (d_v+ d_v^t) = \lfrac 12 D_v\quad\text{and}\quad \lfrac
12 (\imu(\gO) - \imu(\gO)^t) = - \lfrac 12 c_v(\gO),
\]
which explain the remaining terms in the definition of $\B$.
\end{remark}

Since Remark \ref{BismutUnitaryRem} is the underlying motivation
for large parts of the treatment in this section, we extract the
following result, adding some observations which are immediate.

\begin{prop}\label{BismutUnitary}
The Bismut superconnection is the generalized unitary connection
associated to the flat superconnection $d_M$. It is a
superconnection with respect to the even/odd grading on
$\gO_v^\bullet(M)$. If the fiber is even dimensional, it is also
a superconnection with respect to the grading induced by the
vertical chirality operator $\tau_v$.
\end{prop}

\begin{remark}\label{AuxilliaryGrassmann}
In the context of the signature operator we are interested in the
grading given by the vertical chirality operator. However, we get
a superconnection with respect to this grading only if the fiber
is even dimensional. In the case that the fiber is odd
dimensional, one can turn $\B$ into a superconnection by adding
an auxiliary Grassmann variable, see \cite[Sec.'s II (b) \&
(f)]{BF2}. This is based on Quillen's ideas in the finite
dimensional case as in \cite[Sec. 5]{Qui85b}.
\end{remark}

\noindent\textbf{The Chern Character of the Bismut
Superconnection.} As in Definition \ref{SuperConnection}, the
curvature of the Bismut superconnection is defined as the
differential operator
\[
\B^2: \gO^\bullet(M)\to \gO^\bullet(M).
\]
In analogy with the finite dimensional situation in Section
\ref{IndefiniteMetric}, we cannot expect that the Bismut
superconnection is flat.

\begin{prop}
The curvature of $\B$ is a fiberwise elliptic operator. It is of
second order with leading term given by the vertical Laplacian
\[
D_v^2: \gO^\bullet(M,E)\to \gO^\bullet(M,E).
\]
\end{prop}

\begin{proof}
According to Definition \ref{FiberEllOpDef}, we have to check
first that $\B^2$ is $C^\infty(B)$ linear. For all $\gf\in
C^\infty(B)$,
\[
[\B,\pi^*\gf] = \emu(\pi^*d_B\gf),
\]
and thus,
\[
\big[\B^2,\pi^*\gf\big] = \B\circ \emu(\pi^*d_B\gf) +
\emu(\pi^*d_B\gf)\circ \B.
\]
Now, the operators $D_v$, $c_v(\gO)$ and $f^a\wedge k_v(f_a)$ all
anti-commute with exterior multiplication with a horizontal
1-form. Hence,
\[
\big[\B^2,\pi^*\gf\big] = f^a\wedge\nabla^\oplus_a(\pi^*d_B\gf) =
\pi^*d_B^2\gf =0,
\]
so that $\B^2$ is indeed a fiberwise differential operator. Now,
$\B^2$ contains $D_v^2$ as term of second order but a priori there
might be other contributions coming from $(\nabla^{\sE,u})^2$ and
the anti-commutator of $D_v$ and $\nabla^{\sE,u}$. To see that
this is not the case we note that $\nabla^{\sE,u}$ agrees with
$d_h$ up to a term of order 0. Moreover, we know from Corollary
\ref{dSquare} that
\[
d_h^2 = - \big\{d_v,\imu(\gO)\big\}\quad\text{and}\quad
\big\{D_v,d_h\big\} = \big\{d_v^t,d_h\big\},
\]
where both terms are of order $\le 1$, see Proposition
\ref{dvdhCommutator} for the second term. This implies that
$(\nabla^{\sE,u})^2$ and $\{D_v, \nabla^{\sE,u}\}$ are also of
order $\le1$.
\end{proof}


Since fiberwise elliptic operators are intimately related to
families of elliptic operators, the following version of Theorem
\ref{DuhamelThm} should be plausible. We do not give a proof,
referring to \cite[Thm.'s 9.48 \& 9.51]{BGV} for more details.

\begin{theorem}
The operator $e^{-\B^2}: \gO^\bullet(M)\to \gO^\bullet(M)$ is a
well-defined, fiberwise smoothing operators with coefficients in
$\gO^\bullet(B)$.
\end{theorem}

Using the fiberwise supertrace, this result allows us to study the
Chern character of the superconnection $\B$. To keep the
motivational part of this section reasonably short we skip the
discussion of the odd dimensional case and assume henceforth that
the fiber is even dimensional.

\begin{dfn}
Let $F\hookrightarrow M\xrightarrow{\pi} B$ have even dimensional
fiber $F$, and let $g_v$ be a vertical metric. We define the
\emph{Chern character of the Bismut superconnection} as
\[
\ch_s(\sE,\B):=\Str_v\big(\gamma e^{-\B^2}\big)\in \gO^\bullet(B),
\]
where $\gamma$ is the normalization function as in
\eqref{NormSuperConnDef}.
\end{dfn}

\noindent\textbf{Transgression and the Rescaled Superconnection.}
One of Bismut's main observations in \cite{B86} is that
$\ch_s(\sE,\B)$ is the right candidate for the Chern character
form for the index bundle as in Definition \ref{IndexBundle}.
This becomes apparent when rescaling the vertical metric, see
\cite[Sec. III (c)]{B86}.

For $t\in (0,\infty)$ we rescale the vertical metric with a
factor of $t^{-1}$. As in \eqref{RescaledClifford} this means
that Clifford multiplication with a vertical 1-form $\xi$ has to
be replaced with
\[
c_{v,t}(\xi) = \sqrt t \big(\emu(\xi)-\imu(\xi)\big),
\]
where inner multiplication is with respect to the fixed metric
$g_v$. Also, if $\{e_i\}$ is a local orthonormal frame for $T^vM$
with respect to $g_v$, then $\{\sqrt t e_i\}$ is orthonormal with
respect to the rescaled metric. Since $\nabla^v$ and $k_v$ are
independent of $t$, this motivates the following

\begin{dfn}\label{RescaledBismut}
For $t\in (0,\infty)$ define the \emph{rescaled Bismut
superconnection} by
\[
\B_t := \lfrac{\sqrt t}2 D_v + \nabla^{\sE,u} - \lfrac1{2\sqrt t}
c_v(\gO).
\]
\end{dfn}

Now, the infinite dimensional version of the transgression formula
for the Chern character is the following, see \cite[Thm.
9.17]{BGV}.

\begin{theorem}\label{SuperConnTrans}
For all $t\in (0,\infty)$, the differential form $\ch_s(\sE,\B_t)$
is closed and satisfies the transgression formula
\[
\frac{d}{dt} \ch_s(\sE,\B_t) = -  \gamma\, d_B
\Str_v\Big(\frac{d\B_t}{dt}e^{-\B_t^2}\Big).
\]
\end{theorem}

Since this transgression form will be of importance in the next
section, we make the following abbreviation.

\begin{dfn}\label{TransgressionForm}
The \emph{transgression form} associated to the rescaled Bismut
superconnection is given by
\begin{equation*}
\ga(\B_t):=\frac {\gamma}{\sqrt{2\pi i}}
\Str_v\Big(\frac{d\B_t}{dt}e^{-\B_t^2}\Big)\in \gO^{\odd}(B).
\end{equation*}
\end{dfn}

\noindent\textbf{Bismut's Local Index Theorem for Families.} With
the above ingredients, we can now summarize Bismut's main results
\cite[Thm.'s  3.4, 4.12 \& 4.16]{B86} in the case of the
signature operator, see also \cite[Thm.'s 10.21, 10.23 \&
10.32]{BGV}.

\begin{theorem}[Bismut]\label{Fam:IndexTheorem}
Let $F\hookrightarrow M\xrightarrow{\pi} B$ be endowed with a
vertical metric $g_v$ and a vertical projection, and assume that
$F$ is even dimensional. Then the rescaled superconnection $\B_t$
satisfies
\begin{equation}\label{Fam:McKeanSinger}
\big[\ch_s(\sE,\B_t)\big] = \ch(\Ind D^+_v) \in H^{ev}(B).
\end{equation}
Moreover, with respect to the $C^\infty$-topology on
$\gO^\bullet(B)$ and for $t\to 0$
\begin{equation}\label{Fam:LocalIndex}
\lim_{t\to 0}\ch_s(\sE,\B_t) = \int_{M/B}L(T^vM,\nabla^v)
\quad\text{and}\quad \ga(\B_t)= O(1),
\end{equation}
where $L(T^vM,\nabla^v)$ is the Hirzebruch $L$-form of $T^vM$
with respect to the connection $\nabla^v$. In particular,
\begin{equation}\label{BismutChernTrans:1}
\ch_s(\sE,\B_t) = \int_{M/B}L(T^vM,\nabla^v) -
d\int_0^t\ga(\B_s)ds.
\end{equation}
\end{theorem}

\begin{remark*}\quad\nopagebreak
\begin{enumerate}
\item An immediate consequence of Bismut's local index theorem for
families is Atiyah and Singer's earlier result as stated in
Theorem \ref{ASFam}. This is because \eqref{Fam:McKeanSinger} and
\eqref{BismutChernTrans:1} imply that
\[
\ch(\Ind D^+_v) = \Big[\int_{M/B}L(T^vM,\nabla^v)\Big] \in
H^{\ev}(B).
\]
\item We also want to remark that \eqref{Fam:McKeanSinger} is
a generalization of the McKean-Singer formula in Theorem
\ref{McKeanSinger}. In the case that $B$ is a point, the term
$\ch(\Ind D^+_v)$ coincides with the numerical index. On the
other hand, $e^{-\B_t^2}= e^{-tD_v^2}$, since the higher degree
terms in the Bismut superconnection vanish. Hence, if $B$ is a
point, the equation in \eqref{Fam:McKeanSinger} is equivalent to
the usual McKean-Singer formula.
\end{enumerate}
\end{remark*}

To understand why \eqref{Fam:McKeanSinger} is related to the limit
$t\to \infty$ also in higher dimensions, we summarize some
results from \cite[Sec. 9.3]{BGV}. Let $P_{\ker(D_v)}$ be the
projection $\sE\to \ker D_v$, and define
\begin{equation}\label{UnitaryConnHodge}
\nabla^{\sH_v,u}:=P_{\ker(D_v)} \circ \nabla^{\sE,u} \circ
P_{\ker(D_v)}.
\end{equation}
According to Lemma \ref{NablaCompPart} and Proposition
\ref{BismutUnitary}, the connection $\nabla^{\sH_v,u}$ is a
unitary connection on the bundle $\sH^\bullet_v(M)\to B$,
compatible with the grading given by $\tau_v$. We can thus define
the associated Chern character form
$\ch_s\big(\sH^\bullet_v(M),\nabla^{\sH_v,u}\big)$ as in
\eqref{ChernSuperConnRem}. Then there is the following result, see
\cite[Thm.'s 9.19 \& 9.23]{BGV}.

\begin{theorem}[Berline-Getzler-Vergne]\label{Fam:KerConst}
With respect to the $C^\infty$-topology on $\gO^\bullet(B)$ and
for $t\to \infty$
\begin{equation*}
\lim_{t\to \infty} \ch_s(\sE,\B_t) =
\ch_s\big(\sH^\bullet_v(M),\nabla^{\sH_v,u}\big),\quad\text{and}\quad
\ga(\B_t) = O(t^{-3/2}).
\end{equation*}
In particular,
\begin{equation}\label{BismutChernTrans:2}
\ch_s\big(\sH^\bullet_v(M),\nabla^{\sH_v,u}\big) = \ch_s(\sE,\B_t)
- d\int_t^\infty \ga(\B_s) ds.
\end{equation}
\end{theorem}

It is immediate from Definition \ref{IndexBundle}, the definition
of $\nabla^{\sH_v,u}$ and \eqref{ChernSuperConnRem} that
\[
\big[\ch_s\big(\sH^\bullet_v(M),\nabla^{\sH_v,u}\big)\big] =
\ch(\Ind D^+).
\]
so that \eqref{BismutChernTrans:2} is an extension to
differential forms of \eqref{Fam:McKeanSinger}, and explains why
the latter is related to taking the limit $t\to \infty$.

\begin{remark*}
We want to point out that Theorem \ref{Fam:KerConst} holds for
more general superconnections but relies on the fact that the
dimensions of the kernels of its homogeneous component of degree 0
do not jump. In contrast, Theorem \ref{Fam:IndexTheorem}
continues to hold without any assumption on the kernels but only
for the Bismut superconnection associated to a family of Dirac
operators.
\end{remark*}


\section{A General Formula for Rho Invariants}\label{GenAdiabaticLimit}

\subsection{Transgression and Adiabatic
Limits}\label{GenAdiabaticHeruistic}

We now want to relate the discussion in the last section to
adiabatic limits of Eta invariants and give a motivation for
Dai's general adiabatic limit formula. A large part of this
subsection will be heuristic without rigorous arguments. We hope
that nevertheless, our discussion helps to give some intuition
underlying the sophisticated theory.\\

\noindent\textbf{Statement of the Reduced Adiabatic Limit
Formula.} Let $F\hookrightarrow M\xrightarrow{\pi} B$ be an
oriented fiber bundle of closed manifolds, where $M$ is odd
dimensional. Since this section has only motivational character,
we assume for simplicity that $F$ is even dimensional and that
$B$ is odd dimensional. We endow the fiber bundle with a vertical
projection and a submersion metric $g=g_B\oplus g_v$. Let
$\ga(\B_t)$ be the transgression form associated to the rescaled
Bismut superconnection, see Definition \ref{TransgressionForm}.
Theorem \ref{Fam:IndexTheorem} and Theorem \ref{Fam:KerConst}
show that $\ga(\B_t)=O(1)$ as $t\to 0$ and
$\ga(\B_t)=O(t^{-3/2})$ as $t\to \infty$. This implies that we
can make the following definition, which goes back to \cite{BC89}.

\begin{dfn}\label{BCForm:Even}
If $F$ is even dimensional and $B$ is odd dimensional, we define
the \emph{Bismut-Cheeger Eta form}
\[
\widehat \eta := \int_0^\infty \ga(\B_t)dt.
\]
\end{dfn}

Now let $g_\eps$ be the adiabatic metric \eqref{AdiabaticMetric}
on $M$, and consider the associated adiabatic family of de Rham
operators,
\[
D_{M,\eps} = D_v+ \eps\cdot D_h+ \eps^2\cdot T:\gO^\bullet(M)\to
\gO^\bullet(M).
\]
Moreover, recall that in Definition \ref{CohomTwistSignOp} we
have introduced the odd signature operator on
$B$ with values in the bundle of vertical cohomology groups,
\[
D_B\otimes \nabla^{\sH_v}:\gO^\bullet\big(B,\sH_v^\bullet(M)\big)
\to \gO^\bullet\big(B,\sH_v^\bullet(M)\big).
\]
Then we have the following version of \cite[Thm.'s 0.1 \&
4.4]{Dai91}.

\begin{theorem}[Dai]\label{DaiRed}
Assume that $b:=\dim B$ is odd and that $\dim F$ is even. Then the
adiabatic limit of the reduced $\xi$-invariant is given by
\[
\lim_{\eps\to 0}\xi\big(\tau_M D_{M,\eps}\big) =
\xi\big(D_B\otimes \nabla^{\sH_v}\big) +
2^{\frac{b+1}2}\int_B\widehat L(TB,\nabla^B)\wedge\widehat
\eta\mod \Z.
\]
\end{theorem}

This formula has a long history which started with Witten's famous
holonomy formula in \cite{Wit85}. The first mathematically
rigorous treatments were given by Bismut and Freed in \cite{BF2}
as well as Cheeger in \cite{Che87}, both for the case that
$B=S^1$. Bismut and Freed emphasize the relation to Bismut's
local index theory, whereas Cheeger gives an independent proof
based on Duhamel's formula and finite propagation speed methods.
Later in \cite{BC89}, Bismut and Cheeger generalized the
adiabatic limit formula to higher dimensional base spaces under
the assumption that the vertical Dirac operator is invertible.
Note that in this case the twisted Eta invariant on the base does
not appear. Dai's remarkable generalization in \cite{Dai91}
allows the vertical Dirac operators to have non-trivial
kernels---which, however, need to form a vector bundle over $B$.
This is why his result applies in particular to the odd signature
operator. We will state Dai's result in its full generality
below, see Theorem \ref{DaiMain}.\\

\noindent\textbf{The Total Dirac Operator.} For the rest of this
subsection we want to give a heuristic derivation of Theorem
\ref{DaiRed}. Recall that we have seen in Theorem
\ref{SuperConnTrans} that the form $\ga(\B_t)$ plays the role of
the transgression form associated to the Chern character of the
superconnection $\B_t$. The rough idea is now that the adiabatic
limit formula of Theorem \ref{DaiRed} is an infinite dimensional
analog of the variation formula for the Eta invariant as in
Proposition \ref{EtaConnVar}, respectively Corollary
\ref{RhoConnVar}.

This idea does not apply directly to the adiabatic family of odd
signature operators, and as a tool we have to introduce another
natural Dirac operator acting on $\gO^\bullet(M)$. Let $c:T^*M\to
\End\big(\gL^\bullet T^*M\big)$ denote the natural Clifford
structure, and define a connection on $\gL^\bullet T^*M$ via
\begin{equation*}
\nabla^S:=\nabla^\oplus +\lfrac 12 c(\gt),
\end{equation*}
where $\gt$ is the tensor as in Definition \ref{SThetaDef}. One
easily checks that $\nabla^S$ is a Clifford connection as in
Definition \ref{GeomDiracOp}, see e.g. \cite[Prop. 10.10]{BGV}.

\begin{dfn}
The \emph{total Dirac operator} $D_S$ on $M$ is the geometric
Dirac operator associated to $\nabla^S$,
\[
D_S= c\circ \nabla^S: \gO^\bullet(M)\to \gO^\bullet(M).
\]
\end{dfn}

\begin{remark*}
We recall from Lemma \ref{LCDiff} that on $\gL^\bullet T^*M$, the
difference of the Levi-Civita connection $\nabla^g$ and the
connection $\nabla^\oplus$ is given by
\[
\nabla^g = \nabla^\oplus +\lfrac 12 \big(c(\gt)-\widehat
c(\gt)\big),
\]
where $\widehat c$ is the transposed Clifford structure.
Therefore, the total Dirac operator $D_S$ does in general not
coincide with the de Rham operator $D_M$.
\end{remark*}

\noindent\textbf{Relation to the Bismut Superconnection.} As in
Section \ref{LocFamIndSign}, we interpret $\gO^\bullet(M)$ as the
space of differential forms on $B$ with values in the infinite
dimensional vector bundle $\sE\to B$ of vertical differential
forms. Then the total Dirac operator is a differential operator
\[
D_S : \gO^\bullet(B,\sE)\to \gO^\bullet(B,\sE).
\]
Using the local description \eqref{ThetaCoord} for the tensor
$\gt$, one then checks that with respect to a local orthonormal
frame $\{f_a\}$ for $TB$,
\begin{equation}\label{TotalDiracLoc}
D_S = D_v + c^a\nabla_a^{\sE,u} - \lfrac 18 c^ac^b c_v(\gO_{ab}),
\end{equation}
where $\nabla^{\sE,u}$ is defined as in
\eqref{BismutSuperconnConn}. The fact that the terms appearing in
this formula are reminiscent of the terms appearing in the
definition of the Bismut superconnection seems to be one of the
main ideas which lead Bismut to the definition of $\B$, compare
with \cite[Thm. 2.24]{B98} and \cite[Sec. V]{B86}. We want to make
this relation more precise now.

In the interpretation of Proposition \ref{BismutUnitary}, the
Bismut superconnection is a superconnection on $\gL^\bullet
T^*B\otimes \sE$. Comparing with Definition
\ref{CliffordSuperConn}, one finds that it is a Clifford
connection with respect to the Clifford structure induced by
$c:T^*B\to \gL^\bullet T^*B$. Now formally using
\eqref{DiracSuperConn}, we can associate a Dirac operator to
this, which is locally given by
\[
\lfrac 12 D_v + c^a\nabla_a^{\sE,u} - \lfrac 14 c^ac^b
c_v(\gO_{ab}).
\]
Up to factors of $\frac 12$, this is coincides with
\eqref{TotalDiracLoc}. Hence, the total Dirac operator is
essentially the Dirac operator on $\gO^\bullet(B,\sE)$ associated
to the Bismut superconnection. To state the relation more
precisely, consider a metric of the form $g_B\oplus \lfrac 1t
g_v$, where $t\in (0,\infty)$, and let $D_{S,t}$ be the
associated total Dirac operator. As in the discussion preceding
Definition \ref{RescaledBismut} one infers that
\[
D_{S,t} = \sqrt{t} D_v + c^a\nabla^{\sE,u}_a - \lfrac {1}
{8\sqrt{t}}\, c^ac^b c_v(\gO_{ab}).
\]
Then we have as in \cite[Thm. 10.19]{BGV}

\begin{prop}\label{TotalDiracBismut}
Consider the Dirac operator on $\gO^\bullet(B,\sE)$ associated to
the rescaled Bismut superconnection $\B_t$ on $\gL^\bullet
T^*B\otimes \sE$, i.e.,
\[
D_{\B_t}:= c \circ\boldsymbol{\gs}^{-1}\circ \B_t:
\gO^\bullet(B,\sE)\to\gO^\bullet(B,\sE)
\]
Then, under the identification
$\gO^\bullet(M)=\gO^\bullet(B,\sE)$, we have $D_{\B_{4t}} =
D_{S,t}$.
\end{prop}

\begin{remark*}\quad\nopagebreak
\begin{enumerate}
\item The factor 4 occurs because we have defined the Bismut
superconnection in a slightly different way compared to
\cite{B86}. This is because we wanted Proposition
\ref{BismutUnitary} to hold without any constants appearing in
$\B$, see also \cite[Rem. 3.10]{BL95}.
\item We also want to stress that Proposition
\ref{TotalDiracBismut}---in the same way as Proposition
\ref{BismutUnitary} earlier---is more a formal interpretation
than a mathematical statement. Nevertheless, this helps in
understanding some of the ideas underlying the technicalities of
local families index theory.
\end{enumerate}
\end{remark*}

\noindent\textbf{Variation of the Eta Invariant of the Total Dirac
Operator.} After having introduced the total Dirac operator and
its interpretation in term of the rescaled Bismut superconnection
our aim is now to give a heuristic explanation of Theorem
\ref{DaiRed}, in the case that the adiabatic family of de Rham
operators is replaced with a corresponding family of total Dirac
operators.

We first note that since $D_{S,t}$ is a path of formally
self-adjoint elliptic operators of first order on
$\gO^\bullet(M)$, the general variation formula for the
$\xi$-invariant in Proposition \ref{RedEtaDerApp} shows that
\begin{equation}\label{EtaVarInfinite}
\xi\big(\tau_M D_{S,t_0}\big) - \xi\big(\tau_M D_{S,t_1}\big) =
\int_{t_0}^{t_1} \frac 1{\sqrt\pi}\LIM_{u\to 0} \sqrt
u\,\Tr\Big[\tau_M\lfrac{dD_{S,t}}{dt}e^{-uD_{S,t}^2}\Big]
dt\mod\Z,
\end{equation}
where $\LIM_{u\to 0}$ means taking the constant term in the
asymptotic expansion as $u\to 0$. Given the interpretation of
Proposition \ref{TotalDiracBismut}, we now formally apply
Proposition \ref{EtaConnVarLoc} to the case at hand and get
\begin{equation}\label{LocalVarInfinite}
\begin{split}
\frac 1{\sqrt\pi}\LIM_{u\to 0} &\sqrt
u\,\Tr\Big[\tau_M\lfrac{dD_{S,t}}{dt}e^{-uD_{S,t}^2}\Big]\\ &= -
2^{\frac{b+1}2}\int_B \widehat L(TB,\nabla^B)\wedge
\Tr_v\Big[\tau_v \lfrac{\gamma}{\sqrt{2\pi i}}\big(\lfrac
{d\B_{4t}}{dt}e^{-\B_{4t}^2}\big) \Big]\\
&= 2^{\frac{b+1}2}\int_B \widehat L(TB,\nabla^B)\wedge 4
\ga(\B_{4t}).
\end{split}
\end{equation}
\begin{remark*}
Note that in contrast to Proposition \ref{EtaConnVarLoc} we are
not only in the situation that the twisting bundle $\sE$ is
infinite dimensional, but are also dealing with a path of
superconnections rather than usual connections. This requires
special considerations already in the case of a finite dimensional
twisting bundle, see \cite[Sec. 2]{Get94}.
\end{remark*}

Despite this apparent lack of mathematical rigor, we assume for
the rest of this motivational part that the formula in
\eqref{LocalVarInfinite} is valid.\\

\noindent\textbf{The Limit as $\boldsymbol{t\to \infty}$.} Note
that for the total Dirac operator, the component
$c^a\nabla^{\sE,u}_a$ is the analog of the horizontal de Rham
operator $D_h$ in Definition \ref{PartDeRhamOpDef}. In analogy to
Definition \ref{CohomTwistSignOp} we thus define an operator
\begin{equation}\label{CohomUnitaryTwistSignOp}
\begin{split}
&D_B\otimes \nabla^{\sH_v,u}:
\gO^\bullet\big(B,\sH_v^\bullet(M)\big) \to
\gO^\bullet\big(B,\sH_v^\bullet(M)\big),\\
&D_B\otimes \nabla^{\sH_v,u} := P_{\ker D_v}\circ
\tau_M\big(c^a\nabla_a^{\sE,u}\big) \circ P_{\ker D_v},
\end{split}
\end{equation}
compare with \eqref{UnitaryConnHodge}. Then the remarkable result
\cite[Thm. 1.6]{Dai91} guarantees the following

\begin{theorem}[Dai]\label{KernelCollaps}
Abbreviate $D_{S,\infty}:= D_B\otimes \nabla^{\sH_v,u}$. Then, for
$N$ large enough, there exist positive constants $C_1$ and $C_2$
such for all $s>0$ as $t\to \infty$
\[
\Big| \Tr'\big(\tau_M D_{S,t}e^{-sD_{S,t}^2} \big)
-\Tr\big(D_{S,\infty}e^{-sD_{S,\infty}^2} \big) \Big| \le \frac
{C_1}{\sqrt t s^N} e^{-C_2 s}.
\]
Here, $\Tr'$ indicates taking the trace over those eigenvalues of
$D_{S,t}$ which are bounded away from 0 as $t\to \infty$.
\end{theorem}
It is shown in \cite[Thm. 1.5]{Dai91} that there are only
finitely many eigenvalues of $D_{S,t}$ which converge to zero as
$t\to \infty$. We also note that the formulation in \cite{Dai91}
is slightly different from what is stated here. This is due to
the fact that we are using a different scaling and a different
parameter, see Remark \ref{RescaleLargeTime} (iii) below. For the
time being, we only use that---very roughly---Theorem
\ref{KernelCollaps} means that
\[
\lim_{t\to\infty} \xi\big(\tau_M D_{S,t}\big) = \xi\big(D_B\otimes
\nabla^{\sH_v,u}\big)\mod\Z.
\]
Together with \eqref{EtaVarInfinite} and \eqref{LocalVarInfinite}
we arrive at the following heuristic formula
\begin{equation}\label{InfTransg:1}
\xi\big(\tau_M D_S\big) = \xi\big(D_B\otimes
\nabla^{\sH_v,u}\big) + 2^{\frac{b+1}2}\int_B \widehat
L(TB,\nabla^B)\wedge \int_1^{\infty}4\ga(\B_{4t})dt\mod\Z.
\end{equation}

\noindent\textbf{The Adiabatic Limit Formula for
$\boldsymbol{D_S}$.} To understand how an adiabatic limit comes
into play, we consider a metric of the form
\begin{equation*}
g_{\eps,t} = \lfrac 1{\eps^2}g_B \oplus \lfrac 1t g_v.
\end{equation*}
Let $D_{S,\eps,t}$ be the associated total Dirac operator. Since
multiplying $g_{\eps,t}$ by $\eps^{-1}$ does not change the
$\xi$-invariant, we have
\[
\xi\big(D_{S,\eps,t}\big) = \xi\big(\sqrt{\eps} D_{S,\eps,t}\big).
\]
We note that explicitly,
\[
\sqrt{\eps} D_{S,\eps,t} = \sqrt{\eps t} D_v + \sqrt{\eps}^3
c^a\nabla^{\sE,u}_a - \lfrac {\eps^3} {8\sqrt{\eps t}}\, c^ac^b
c_v(\gO_{ab}).
\]
As we have seen in Section \ref{RhoLocalInd}, the proof of
Proposition \ref{EtaConnVarLoc} uses Getzler's rescaling by
$\sqrt{u}$. Since \eqref{LocalVarInfinite} is a formal transition
to the case at hand, it is reasonable to expect that the same
rescaling is involved there. Without going into detail, we note
that using the rescaling $\sqrt{u\eps^3}$ instead, one formally
obtains
\[
\frac 1{\sqrt\pi}\LIM_{u\to 0} \sqrt
u\,\Tr\Big[\tau_M\lfrac{dD_{S,\eps,t}}{dt}e^{-uD_{S,\eps,t}^2}\Big]=
2^{\frac{b+1}2}\int_B \widehat L(TB,\nabla^B)\wedge 4\eps
\ga(\B_{4\eps t}).
\]
Inserting this in \eqref{InfTransg:1} and substituting $s=4\eps
t$, we arrive at
\begin{equation}\label{InfTransg:2}
\xi\big(\tau_M D_{S,\eps}\big) = \xi\big(D_B\otimes
\nabla^{\sH_v,u}\big) + 2^{\frac{b+1}2}\int_B \widehat
L(TB,\nabla^B)\wedge \int_{4\eps}^{\infty}\ga(\B_s)ds\mod\Z,
\end{equation}
where we have used that the term $\xi\big(D_B\otimes
\nabla^{\sH_v,u}\big)$ is independent of $\eps$. This is because
$D_B\otimes \nabla^{\sH_v,u}$ is an operator over $B$ and
rescaling the full metric does not change the $\xi$-invariant.
Now letting $\eps\to 0$ yields the analog of Theorem \ref{DaiRed}
for the total Dirac operator $D_{S,\eps}$.

\begin{remark*}
We want to point out again, that our considerations leading to
\eqref{InfTransg:2} are heuristic and have no rigorous
mathematical foundation. In fact, we do not even expect
\eqref{EtaVarInfinite}, \eqref{InfTransg:1} and
\eqref{InfTransg:2} to be correct without any changes. However, it
would be interesting to find estimates on the correction term in
\eqref{InfTransg:2} along the line of thoughts we have presented
to give an alternative proof of Theorem \ref{DaiRed} for the
total Dirac operator.
\end{remark*}

\noindent\textbf{Relation between $\boldsymbol{D_{S,\eps}}$ and
$\boldsymbol{D_{M,\eps}}$.} What is still missing in our
heuristic explanation of Theorem \ref{DaiRed} is why we can
replace $D_{S,\eps}$ with $D_{M,\eps}$ and $D_B\otimes
\nabla^{\sH_v,u}$ with $D_B\otimes \nabla^{\sH_v}$. First of all,
the relation between the latter two operators fits exactly into
the framework of Section \ref{IndefiniteMetric}. In particular,
Theorem \ref{CohomTwistEtaDiff} yields that
\[
\xi\big(D_B\otimes \nabla^{\sH_v,u}\big) = \xi\big(D_B\otimes
\nabla^{\sH_v}\big)\mod\Z.
\]
Formally, the relation between $D_{S,\eps}$ and $D_{M,\eps}$ is an
infinite dimensional analog of the same situation, where the odd
signature operator $\tau_M D_{M,\eps}$ plays the role of
$D_B\otimes \nabla$, whereas $\tau_M D_{S,\eps}$ is an analog of
the operator $D_B\otimes \nabla^u$, compare with Lemma
\ref{SignCompPart} and Proposition \ref{BismutUnitary}. Then the
heuristic analog of Theorem \ref{CohomTwistEtaDiff} is that we
can make the required substitution in Theorem \ref{DaiRed}. For a
technically precise explanation, see \cite[p. 304]{Dai91} and
\cite[p. 374]{BC92}.

\subsection{Dai's Adiabatic Limit Formula}\label{Dai:Gen}

After this heuristic digression, we now want to state Dai's
general adiabatic limit formula for the Eta invariant of the odd
signature operator in its full generality. From now on we will
also include the case that $\dim B$ is even and $\dim F$ is odd.
Hence, we first need the analog of the Bismut-Cheeger Eta form in
Definition \ref{BCForm:Even} for the case of odd dimensional
fibers, see \cite[Def. 4.93 \& Rem. 4.100]{BC89}.

\begin{dfn}\label{BCForm:Odd}
Assume that $\dim F$ is odd and $\dim B$ is even, and let $\B_t$
be the rescaled Bismut superconnection as in Definition
\ref{RescaledBismut}. Then we define
\[
\widehat \eta := \frac 1{\sqrt\pi} \int_0^\infty \gamma
\Tr_v^{\ev}\big[\tau_v\lfrac {d\B_t}{dt}e^{-\B_t^2}\big]dt\in
\gO^{\ev}(B),
\]
where $\gamma$ is the normalization function as in
\eqref{NormSuperConnDef}, and $\Tr_v^{\ev}$ indicates taking the
even form part of $\Tr_v$.
\end{dfn}

\begin{remark*}\quad\nopagebreak
\begin{enumerate}
\item The vertical chirality operator enters in Definition
\ref{BCForm:Odd}, since in contrast to \cite{BC89} we are dealing
with the operator $\tau_M D_M$ rather than the spin Dirac
operator. Nevertheless, $\Tr_v\circ \tau_v$ should not be viewed
as a supertrace, since in the case that $\dim F$ is odd, $\tau_v$
commutes with vertical Clifford multiplication.
\item The fact that the integral defining $\widehat \eta$ is indeed
convergent relies on the odd dimensional analogs of the small and
large time behaviour of $\ga(\B_t)$ in Theorem
\ref{Fam:IndexTheorem} and Theorem \ref{Fam:KerConst}, i.e.,
\[
\Tr_v^{\ev}\big[\tau_v\lfrac {d\B_t}{dt}e^{-\B_t^2}\big] =
O(1)\quad\text{as $t\to 0$},\quad\Tr_v^{\ev}\big[\tau_v\lfrac
{d\B_t}{dt}e^{-\B_t^2}\big] = O(t^{-3/2})\quad\text{as $t\to
\infty$}.
\]
It is difficult to find an explicit proof in the literature,
since the treatment is usually in the superconnection formalism
which does not apply to the operator in question. However, as
pointed out in Remark \ref{AuxilliaryGrassmann} one can overcome
this difficulty by introducing an extra Grassmann variable to
turn $\B_t$ into a superconnection of the required form, see
again \cite[Sec.'s II (b) \& (f)]{BF2},  and also \cite[Sec.
5.2.2]{BunMa} for additional remarks and references. The reader
who feels uncomfortable with this can equally well consider
\[
\widehat \eta(s) := \frac 1{\sqrt\pi} \int_0^\infty \gamma
\Tr_v^{\ev}\big[\tau_v\lfrac {d\B_t}{dt}e^{-\B_t^2}\big]t^s
dt,\quad \Re(s)\text{ large},
\]
and define the Eta form as the constant term in the Laurent series
around $s=0$ of the meromorphic continuation. Then the whole
discussion to follow goes through with only minor changes---the
sole problem being a more awkward notation.
\item In contrast to the case that $F$ is even dimensional, the
Bismut-Cheeger Eta form is now a differential form on $B$ of
\emph{even} degree. One easily checks that the degree 0 term is
given by
\[
\widehat \eta_{[0]} = \frac 1{2 \sqrt\pi}\int_0^\infty
u^{-1/2}\Tr_v\big[\tau_v D_v e^{-u D_v^2}\big] ds,
\]
where we have substituted $4u = t$. Hence, if $B_v^{\ev}= \tau_v
D_v|_{\gO^{\ev}_v(M)}$ is the family of vertical odd signature
operators, we see that
\begin{equation}\label{RhoFormDeg0}
\widehat \eta_{[0]} = \eta(B_v^{\ev})\in C^\infty(B),
\end{equation}
which is the function that associates to each point $y\in B$ the
Eta invariant of the fiber $\pi^{-1}(y)$. Note that we have used
\eqref{EtaDirect} which is possible since $B_v^{\ev}$ is a family
of geometric Dirac operators so that the Eta function can be
defined without making use of a meromorphic extension.
\end{enumerate}
\end{remark*}

\noindent\textbf{Preliminary Adiabatic Limit Formula.} The
starting point for the rigorous treatment of the adiabatic limit
formula is \cite[Prop. 1.4]{Dai91}, which in the case of the odd
signature operator reads

\begin{prop}[Dai]\label{DaiPrep}
Let $D_{\eps}$ be the family of de Rham operator associated to the
adiabatic metric \eqref{AdiabaticMetric}. Then there exists a
small positive number $\ga$ such that
\[
\lim_{\eps\to 0}\eta(\tau_M D_{\eps}) = 2^{[\frac{b+1}2]+1}\int_B
\widehat L(TB,\nabla^B)\wedge \widehat \eta + \lim_{\eps\to 0}
\frac 1{\sqrt\pi}\int_{\eps^{-\ga}}^\infty
u^{-1/2}\Tr\big[\tau_MD_\eps e^{-uD_{\eps}^2}\big]du,
\]
provided either one of the limits exists.
\end{prop}

\begin{remark}\quad\nopagebreak\label{RescaleLargeTime}
\begin{enumerate}
\item Dai deduces Proposition \ref{DaiPrep} from the main result in
\cite{BC89} which---translated to the situation at hand---gives
an explicit formula for
\[
\lim_{\eps\to 0}\Tr\big[\tau_M D_{\eps}e^{-uD_{\eps}^2}\big],
\]
together with remainder estimates, which are uniform in $\eps$ for
compact $u$-intervals, see \cite[(4.40)]{BC89}. We want to point
out, that Bismut's and Cheeger's proof is rather involved, one
difficulty being the presence of the term $\Tr\big[\tau_M \eps
D_h e^{-uD_{\eps}^2}\big]$, which behaves ``singular'' with
respect to Getzler's rescaling, see \cite[Rem. 3.4]{BC89}. This
part should simplify if one could find a proof along the lines of
the heuristic discussion in Section \ref{GenAdiabaticHeruistic}.
\item Note that for the odd signature operator, we already know from
Proposition \ref{AdiabaticLimitExist} that the limit on the left
hand side in Proposition \ref{DaiPrep} exists. By comparison with
Theorem \ref{DaiRed} we see that the second term on the right
hand side will produce the twisted Eta invariant of the base as
well as the integer contribution which we omitted so far.

\item To relate this term to the discussion in Section
\ref{GenAdiabaticHeruistic}, we substitute $s=\eps^2 u$, which
corresponds to rescaling the metric $g_\eps$ to $\eps^2 g_\eps$.
Then
\[
\int_{\eps^{-\ga}}^\infty u^{-1/2}\Tr\big[\tau_MD_\eps
e^{-uD_{\eps}^2}\big]du = \int_{\eps^{2-\ga}}^\infty
s^{-1/2}\Tr\big[\tau_M(\lfrac 1{\eps}D_\eps) e^{-s(\frac 1{\eps}
D_{\eps})^2}\big]ds.
\]
We now note that if we rename $\sqrt{t} = \eps^{-1}$, then $\lfrac
1{\eps} D_{\eps}$ becomes
\[
\sqrt t D_{(\sqrt t)^{-1}} = \sqrt t D_v + D_h + \lfrac1{\sqrt t}
T.
\]
where the individual terms are defined as in Definition
\ref{PartDeRhamOpDef}. This means that the operator $\lfrac
1{\eps} D_{\eps}$ plays essentially the role of the operator
$D_{S,t}$ in the discussion of Section
\ref{GenAdiabaticHeruistic}, and that taking the limit $\eps\to
0$ on the right hand side in Proposition \ref{DaiPrep} is related
to the limit $t\to \infty$ in our heuristic explanation of Theorem
\ref{DaiRed}. We also note without giving the details that this
also explains why Theorem \ref{KernelCollaps} as we have stated
it is indeed a reformulation of \cite[Thm. 1.6]{Dai91}.
\end{enumerate}
\end{remark}

\noindent\textbf{Behaviour of Small Eigenvalues.} The limit on the
right hand side of the formula in Proposition \ref{DaiPrep} is
closely related to eigenvalues of $\tau_MD_{\eps}$ which approach
zero as $\eps\to 0$. This is roughly because of their presence,
there is no uniform bound of the form $Ce^{-cu}$ for the term
$\Tr\big[\tau_MD_\eps e^{-uD_{\eps}^2}\big]$ as $\eps\to 0$ for
large $u$, compare with Lemma \ref{Duhamel:2:Lem}. The following
result provides the essential analysis of the spectrum of $D_\eps$
as $\eps\to 0$, see \cite[Thm. 1.5]{Dai91}.

\begin{theorem}[Dai]\label{AdiabaticSpectrum}
For $\eps>0$ the eigenvalues of $\tau_MD_\eps$ depend analytically
on $\eps$. There exist analytic functions
$\setdef{\gl_\eps^i}{i\in\Z}$ such that $\spec(\tau_M
D_\eps)=\bigcup_{i\in\Z}\gl_\eps^i$ for all $\eps>0$. Moreover,
the functions $\gl_\eps^i$ have the following properties.
\begin{enumerate}
\item There exists a positive constant $\gl_0$ such that either
for some $\eps_0$
\[
|\gl_\eps^i|\ge \gl_0>0,\quad\text{for $\eps\le \eps_0$},
\]
or $\gl_\eps^i$ has a complete asymptotic expansion
\[
\gl_\eps^i \sim \sum_{k\ge 1}\mu_k^i\eps^k\quad\text{as $\eps\to
0$},
\]
where $\mu_1^i\in \spec\big(D_B\otimes \nabla^{\sH_v}\big)$, see
Definition \ref{CohomTwistSignOp}. The latter gives a bijective
correspondence
\[
\Bigsetdef{\gl_\eps\in \spec(\tau_MD_\eps)}{\gl_\eps=O(\eps)\text{
as }\eps\to
0}\overset{1:1}{\longleftrightarrow}\spec\big(D_B\otimes
\nabla^{\sH_v}\big).
\]
\item Assume that $\gl_\eps^i=O(\eps)$ as $\eps\to 0$, and that $\mu^i_1\neq
0$. Then there is a uniform remainder estimate of the form
\[
\gl_\eps^i = \mu_1^i\eps + \eps^2 C(\eps)(\mu_1^i)^2,\quad
|C(\eps)|\le \const.
\]
\item For every $K>0$,
\[
\#\Bigsetdef{i\in \Z}{\gl_\eps^i=O(\eps)\text{ as }\eps\to
0,\text{ and }|\mu_1^i|\le K}<\infty.
\]
\end{enumerate}
\end{theorem}

\begin{remark*}
As one might expect, the proof of Theorem \ref{AdiabaticSpectrum}
in \cite{Dai91} relies on standard perturbation theory as in
\cite{K}, and this part is in fact not very difficult. However, to
prove that the eigenvalues have a complete asymptotic expansion,
Dai makes use of Melrose's theory of degenerate elliptic problems
as used in \cite{MazMel}. We refer to \cite[Sec. 2]{Dai91} for
more details and references.
\end{remark*}

\noindent\textbf{The General Adiabatic Limit Formula.} For $r\in
\N$ define
\begin{equation}\label{EigenvalueOrders}
\gL(\eps^r):=\bigsetdef{\gl_\eps\in\spec(\tau_MD_\eps)}{\gl_\eps=
O(\eps^r) ,\text{ as }\eps\to 0}.
\end{equation}
Note that part (iii) of Theorem \ref{AdiabaticSpectrum} shows
that $\#\gL(\eps^r)<\infty$ for all $r\ge 2$. Then Dai's main
results can be stated in the following way, see \cite[Cor. 1.6 \&
Prop. 1.8]{Dai91}. The formal but illuminating outline of the
proof in \cite[p. 275]{Dai91} is also recommended.

\begin{theorem}[Dai]\label{DaiMain}
With the notation of Proposition \ref{DaiPrep},
\[
\lim_{\eps\to 0} \frac 1{\sqrt\pi}\int_{\eps^{-\ga}}^\infty
u^{-1/2}\Tr\big[\tau_MD_\eps e^{-uD_{\eps}^2}\big]du =
\eta\big(D_B\otimes \nabla^{\sH_v}\big) + \lim_{\eps\to 0}
\sum_{\gl_\eps\in \gL(\eps^2)} \sgn(\gl_\eps).
\]
In particular,
\[
\begin{split}
\lim_{\eps\to 0}\eta(\tau_M D_{\eps}) = 2^{[\frac{b+1}2]+1}\int_B
\widehat L(TB&,\nabla^B)\wedge \widehat \eta\\ &+
\eta\big(D_B\otimes \nabla^{\sH_v}\big) + \lim_{\eps\to 0}
\sum_{\gl_\eps\in \gL(\eps^2)} \sgn (\gl_\eps).
\end{split}
\]
\end{theorem}

This separates the computation of the adiabatic limit of the Eta
invariant on the total space of the fiber bundle into three terms,
which are all of a very different nature. Intuitively, the first
term contains global information about the fiber, but is local on
$B$. The second term is global on the base and contains
cohomological information about the fiber. The third term is
global on both, the base and the fiber. It fits into the
heuristic discussion of Section \ref{GenAdiabaticHeruistic} as an
analog of the spectral flow term in Corollary \ref{RhoConnVar}.
Following again \cite{Dai91}, we will see in the next subsection
that for the odd signature operator, this term is expressible in
completely topological terms.

\subsection{Small Eigenvalues and the Leray Spectral
Sequence}\label{SF:SpecSeq}

Let $\gD_\eps$ be the Laplace operator associated to an adiabatic
metric $g_\eps$ on an oriented fiber bundle $F\hookrightarrow
M\xrightarrow{\pi} B$ of closed manifolds. In \cite{MazMel},
Mazzeo and Melrose analyze the space of harmonic forms
$\sH_\eps^\bullet(M)=\ker(\gD_\eps)$ as $\eps \to 0$, and show
that it has a basis which extends smoothly to $\eps=0$. Using a
Taylor series analysis to determine which forms lie in this
limiting space, they find a Hodge theoretic version of the
Leray-Serre spectral sequence. We review the result briefly, and
refer to \cite{MazMel} as well as \cite{F95} and \cite[Sec.
4.2]{Dai91} for details. As in the latter reference, we restrict
to the odd dimensional case and use the formulation in terms of
the odd signature operator.\\

\noindent\textbf{The Hodge Theoretic Spectral Sequence.} Assume
that $\dim M$ is odd. For $\gL(\eps^r)$ as in
\eqref{EigenvalueOrders} we define
\[
G_{\gL(\eps^r)}:= \sum_{\gl_\eps\in
\gL(\eps^r)}\ker\big(\tau_MD_\eps - \gl_\eps \big).
\]
We view this as a family of subspaces of $\gO^\bullet(M)$,
parametrized by $\eps\in (0,\infty)$. Note that Theorem
\ref{AdiabaticSpectrum} implies that for $r\ge 2$ each
$G_{\gL(\eps^r)}$ is finite dimensional. Now the analysis of
\cite{MazMel} adapted to the case at hand yields the following,
see \cite[Thm. 0.2 \& Prop. 4.2]{Dai91}.

\begin{theorem}[Mazzeo-Melrose, Dai]\label{DaiSpecSeq}
For $r\ge 2$ the family $G_{\gL(\eps^r)}$ depends smoothly on
$\eps$ down to $\eps=0$, i.e., there exists a smooth family of
orthonormal bases for $G_{\gL(\eps^r)}$, parametrized by
$\eps\in(0,\infty)$, that extends smoothly to $\eps=0$. One can
then define
\begin{equation}\label{DaiSpecSeqDef}
G_r:= \lim_{\eps\to 0} G_{\gL(\eps^r)},\quad d_r:=\lim_{\eps\to
0} \eps^{-r}d_\eps:G_r\to G_r,
\end{equation}
where $d_\eps= d_v+\eps d_h +\eps^2\imu(\gO)$. This defines a
spectral sequence that is isomorphic to the Leray-Serre spectral
sequence, i.e.,
\[
(G_r,d_r) \cong (E_r^{\bullet,\bullet},d_r),\quad r\ge 2.
\]
\end{theorem}

\begin{remark*}\quad\nopagebreak
\begin{enumerate}
\item Note that the spectral sequence $(G_r,d_r)$ is not defined
as a spectral sequence associated to a filtered complex as in
\eqref{SpecSeqDef}. Nevertheless, as subspaces of
$\gO^\bullet(M)=\bigoplus_{p,q} \gO^{p,q}(M)$, the spaces $G_r$
are naturally bigraded.
\item Since we are working on $\gO^\bullet(M)$ with the
fixed reference metric $g$ we are using the modified de Rham
differential $d_\eps$, compare with Remark \ref{DeRhamRescaleRem}.
\item It is immediate that each $d_\eps$ maps $G_{\gL(\eps^r)}$ to
itself, since $d_\eps$ commutes with $\tau_M D_\eps$. However,
note that for an element $\go_\eps$ of a basis of
$G_{\gL(\eps^r)}$ as in Theorem \ref{DaiSpecSeq}, one has
\[
(\tau_M D_\eps)\go_\eps = \gl_\eps\go_\eps = O(\eps^r),\quad
\text{as }\eps\to0.
\]
Since $\tau_MD_\eps = \tau_M d_\eps +d_\eps\tau_M$, this implies
that $d_\eps\go_\eps=O(\eps^r)$ as well. This should serve as a
motivation for the factor $\eps^{-r}$ in \eqref{DaiSpecSeqDef}.
\end{enumerate}
\end{remark*}

\noindent\textbf{The case $\boldsymbol{r=1}$.} Because of the
extra difficulty which arises from the fact that $G_{\gL(\eps)}$
is not finite dimensional the case $r=1$ is excluded from Theorem
\ref{DaiSpecSeq}. Nevertheless, for motivational purposes, we now
want to make a formal digression on this case, since it might
give an idea of the mechanism lying behind Theorem
\ref{DaiSpecSeq}. For a rigorous treatment, we refer to
\cite[Sec.'s 3 \& 4]{MazMel}.

Let $\gl_\eps$ be a family of eigenvalues of $\tau_MD_\eps$,
which is of order $\eps$ as $\eps\to 0$. Since for $\eps\in
(0,\infty)$, the family $\tau_MD_\eps$ depends analytically on
$\eps$, standard perturbation theory ensures that $\gl_\eps$
depends analytically on $\eps$ and that there exists an analytic
family of eigenforms $\go_\eps\in \gO^\bullet(M)$ with eigenvalue
$\gl_\eps$. We now assume without justification that $\gl_\eps$
and $\go_\eps$ extend analytically to $[0,\infty)$ so that we can
write
\[
\gl_\eps = \sum_{k\ge 1} \gl_k \eps^k,\quad\text{and}\quad
\go_\eps = \sum_{k\ge 0} \go_k \eps^k,\quad \text{as }\eps\to 0.
\]
From Definition \ref{AdiabaticFamily} we know that
\[
\tau_M D_\eps = \tau_M D_v + \eps \tau_M D_h +\eps^2 \tau_M
T:\gO^\bullet(M)\to\gO^\bullet(M).
\]
Then comparing $\eps$-powers in the identity $(\tau_M D_\eps)
\go_\eps= \gl_\eps \go_\eps$ one finds that
\begin{equation*}\label{MazMelMot:1}
D_v \go_0 =0,\quad \tau_M\big(D_v\go_1 + D_h \go_0\big) = \gl_1
\go_0.
\end{equation*}
In particular, $\go_0\in \ker
(D_v)=\gO^\bullet\big(B,\sH_v^\bullet(M)\big)$ and thus also
$d_v\go_0 =0$. Then the second equation yields
\begin{equation}\label{MazMelMot:2}
P_{\ker(D_v)^\perp}\big(\tau_M\big(D_v\go_1 + D_h
\go_0\big)\big)=0.
\end{equation}
On the other hand, when we compare $\eps$-powers in the identity
$\LScalar{\tau_M d_\eps \go_\eps }{d_\eps\tau_M\go_\eps} =0$, we
can deduce that
\[
\LScalar{\tau_Md_v\go_1 +\tau_Md_h\go_0}{d_v\tau_M\go_1
+d_h\tau_M\go_0} = 0.
\]
Using this and \eqref{MazMelMot:2} one infers that
\[
P_{\ker(D_v)^\perp}\big(d_v\go_1 +d_h\go_0\big) = 0.
\]
Hence,
\[
\begin{split}
\lim_{\eps\to 0} \eps^{-1}d_\eps \go_\eps &= \lim_{\eps\to 0}
\eps^{-1}\big(d_v\go_0 +
\eps(d_v\go_1+d_h\go_0)+\eps^2(\ldots)\big) = d_v\go_1 +d_h\go_0\\
&= P_{\ker(D_v)}(d_h\go_0).
\end{split}
\]
This shows at least formally that the construction of Theorem
\ref{DaiSpecSeq} extends to the case $r=1$ and gives
\[
G_1 = \gO^\bullet\big(B,\sH_v^\bullet(M)\big),\quad d_1 =
P_{\ker(D_v)}\circ d_h.
\]
Now Proposition \ref{FlatConnHodge} shows that via the
Hodge-de-Rham isomorphism,
\[
(G_1,d_1) \cong \big(\gO^p\big(B,H^\bullet_v(M)\big),\Bar
d_h\big),
\]
with the differential $\Bar d_h$ on
$\gO^p\big(B,H^\bullet_v(M)\big)$ associated to the flat
connection $\nabla^{H_v}$, see Definition \ref{FlatVertCohom}.
According to Lemma \ref{SpecSeqLem} this means that $(G_1,d_1)$ is
isomorphic to the $E_1$-term of the Leray-Serre spectral sequence
of the fiber bundle.

\begin{remark*}
In principle, one could now continue and give a formal derivation
of Theorem \ref{DaiSpecSeq} along the lines just presented.
However, already in the case $r=2$, the Hodge theoretical
expression of the differential becomes very unpleasant. Therefore,
we end the digression and refer to \cite{MazMel, Dai91, F95} for
more details.
\end{remark*}

We also want to point out that Theorem \ref{DaiSpecSeq} will not
be explicitly used later on. Yet, we have included the above
discussion to motivate why it is reasonable to expect that the
term
\[
\lim_{\eps\to 0} \sum_{\gl_\eps\in \gL(\eps^2)} \sgn (\gl_\eps)
\]
appearing in Dai's general adiabatic limit formula has a
topological interpretation in terms of the Leray-Serre spectral
sequence. We shall make this more precise now.\\

\noindent\textbf{Multiplicative Structure on the Leray-Serre
Spectral Sequence.} From now on we also incorporate a flat
twisting bundle $E$ over $M$ with connection $A$. The entire
treatment of Section \ref{GenAdiabaticLimit} carries over
verbatim; we have omitted it so far only for notational
convenience. However, in the discussion to follow, there are some
small distinctions to be made. We use the notation
$(E_{A,r}^{\bullet,\bullet},d_{A,r})$ for the spectral sequence
associated to the flat bundle, and reserve the notation
$(E_r^{\bullet,\bullet},d_r)$ from \eqref{SpecSeqDef} for the
spectral sequence associated to the trivial connection.

Recall that there is a multiplicative structure of the form
\[
\cdot:E^{p,q}_{A,r}\times E_{A,r}^{s,t}\to E_r^{p+s,q+t},
\]
which is canonically induced by the wedge-product and the metric
$h:E\otimes E\to \C$ on the bundle $E$, see \cite[pp.
174--177]{BT}, \cite[Thm. 5.2]{McC} or \cite[Thm. 9.24]{DavKir}.
For brevity, we introduce the notation
\[
E^k_{A,r}:= \bigoplus_{p+q=k} E^{p,q}_{A,r},\quad \cdot:
E^k_{A,r}\times E^l_{A,r}\to E^{k+l}_r.
\]
Then the differential and the multiplicative structure satisfy
the relation
\begin{equation}\label{SpecSeqMult}
d_r(\go\cdot \eta) = d_{A,r}(\go)\cdot \eta + (-1)^k \go\cdot
d_{A,r}(\eta),\quad \go\in E_{A,r}^k,\quad \eta\in E_{A,r}^l.
\end{equation}
Denote $m=\dim M$ and $b=\dim B$. Since we are assuming that the
fiber bundle is oriented, the bundle $H_v^{m-b}(M)\to B$ is
trivializable, where each vertical volume form gives a canonical
trivialization. We then have
\[
E_2^m = E_2^{b,m-b}= H^b\big(B,H_v^{m-b}(M)\big)\cong \C,\quad ,
\]
where a natural basis $\xi_2\in E_2^m$ is induced by any volume
form on $M$. Since $M$ is closed $H^m(M)\cong \C$, which implies
that $E_r^m$ for all $r\ge 2$. Moreover, as the isomporphism
$H^m(M)\cong \C$ is also canonically induced by any volume form
on $M$, we obtain natural bases $\xi_r$ for all $E_r^m$ with $r\ge
2$. Using the multiplicative structure on
$E_{A,r}^{\bullet,\bullet}$, we can thus define a natural pairing
\begin{equation}\label{SpecSeqPair}
Q_{A,r}: E_{A,r}^k\times E_{A,r}^{m-k} \to \C,
\end{equation}
by requiring that
\[
Q_{A,r}(\go,\eta)\,\xi_r = \go\cdot\eta, \quad \text{for all
}\quad (\go,\eta)\in E_{A,r}^k\times E_{A,r}^{m-k}.
\]

\begin{remark*}
Due to the presence of the pairing $E\times E\to \C$ in its
definition, the pairing \eqref{SpecSeqPair} is complex anti-linear
in the first variable. For even dimensional manifolds, the above
pairing has been analyzed already by \cite{CHS} in the context of
the signature of fiber bundles.
\end{remark*}

As in \cite[Sec. 4.3]{Dai91} we now define the analog for the odd
dimensional case,
\begin{equation}\label{SpecSeqPairOdd}
P_{A,r}: E_{A,r}^k\times E_{A,r}^{m-k-1}\to \C,\quad
P_{A,r}(\go,\eta):=Q_{A,r}\big(\go,d_{A,r}\eta\big).
\end{equation}

\begin{lemma}\label{SpecSeqPairOddLem}
Let $\go\in E_{A,r}^k$ and $\eta\in E_{A,r}^{m-k-1}$. Then
\[
P_{A,r}(\go,\eta) = (-1)^{(k+1)(m-k)}
\overline{P_{A,r}(\eta,\go)}.
\]
In particular, if $k=\frac{m-1}{2}$, then $P_{A,r}$ is Hermitian
if $m=4n-1$ and skew Hermitian if $m=4n-3$.
\end{lemma}

\begin{proof}
Since $E_r^m\cong \C$ for each $r\ge 2$ it follows that
$d_r|_{E_r^{m-1}}\equiv 0$. Then \eqref{SpecSeqMult} implies that
\[
\begin{split}
Q_{A,r}\big(\go,d_{A,r}\eta\big)\xi_r &= \go\cdot d_{A,r}\eta =
(-1)^{k+1} (d_{A,r}\go)\cdot \eta\\ &=
(-1)^{(k+1)(m-k)}\overline{\eta\cdot d_{A,r}\go}\\ &=
(-1)^{(k+1)(m-k)}
\overline{Q_{A,r}\big(\eta,d_{A,r}\go\big)}\xi_r.
\end{split}
\]
This implies the first assertion. Setting $k=\frac{m-1}2$, the
exponent becomes $(\frac{m+1}2)^2$, which yields the second
assertion.
\end{proof}

\noindent\textbf{Dai's Correction Term.} Using the pairing
$P_{A,r}$ introduced in \eqref{SpecSeqPairOdd}, we now come to
the topological interpretation of the third term in the general
adiabatic limit formula of Theorem \ref{DaiMain}.

\begin{dfn}\label{DaiCorrect}
For each $r\ge 2$ we define
\[
\gs_{A,r}:= \Sign\big(P_{A,r}:E_{A,r}^k\times E_{A,r}^k\to
\C\big),\quad \text{with}\quad k:=\lfrac{m-1}2,
\]
where as before the signature of a skew form is defined as the
number of positive imaginary eigenvalues minus the number of
negative imaginary ones. Moreover, we write
\[
\gs_A:=\sum_{r\ge 2} \gs_{A,r}.
\]
\end{dfn}

Note that the sum defining $\gs_A$ is finite since the spectral
sequence collapses after finitely many steps so that for large
$r$ the number $\gs_{A,r}$ is always 0. This is due to the
presence of the differential $d_{A,r}$ in the definition of
$Q_{A,r}$. Now we can state \cite[Thm. 4.4]{Dai91}.

\begin{theorem}[Dai]\label{DaiGen:2}
Let $E$ be a flat Hermitian bundle with connection $A$ over the
odd dimensional total space of an oriented fiber bundle
$F\hookrightarrow M\xrightarrow{\pi}B$ of closed manifolds. For
any adiabatic metric on $M$, let $\gL_A(\eps^2)$ be defined in
analogy to \eqref{EigenvalueOrders} with respect to the operator
$\tau_M D_{A,\eps}$. Then
\[
\lim_{\eps\to 0} \sum_{\gl_\eps\in \gL_A(\eps^2)} \sgn (\gl_\eps)
= 2\gs_A.
\]
Consequently, the following adiabatic limit formula holds
\[
\lim_{\eps\to 0}\eta(B_{A,\eps}^{\ev}) = 2^{[\frac{b+1}2]}\int_B
\widehat L(TB,\nabla^B)\wedge \widehat \eta_A + \lfrac 12
\eta\big(D_B\otimes \nabla^{\sH_{A,v}}\big) + \gs_A.
\]
Here, $\widehat \eta_A\in \gO^\bullet(B)$ is defined in analogy to
the untwisted case in Definition \ref{BCForm:Even}, respectively
Definition \ref{BCForm:Odd}.
\end{theorem}

Note that we have divided the adiabatic limit formula in Theorem
\ref{DaiMain} by a factor of 2 to get the corresponding formula
for the odd signature operator, see Remark \ref{OddSignRem}.\\

\noindent\textbf{An Adiabatic Limit Formula for the Rho
Invariant.} An immediate consequence of Theorem \ref{DaiGen:2} is
the result we were aiming at in this section. In analogy to
Definition \ref{RhoDef} we first make the following

\begin{dfn}\label{RhoFormDef}
Let $E$ be a flat bundle of rank $k$ with connection $A$ over the
odd dimensional total space of an oriented fiber bundle of closed
manifolds. For every submersion metric, we define the
\emph{Bismut-Cheeger Rho form} as
\[
\widehat\rho_A := \widehat \eta_A - k\cdot \widehat \eta\in
\gO^\bullet(B),
\]
and the \emph{Rho invariant of the bundle of vertical cohomology
groups} as
\[
\rho_{\sH_{A,v}}(B):=\lfrac12 \eta\big(D_B\otimes
\nabla^{\sH_{A,v}}\big) - k\cdot \lfrac12\eta\big(D_B\otimes
\nabla^{\sH_v}\big) \in \R.
\]
\end{dfn}

We include the factor $\lfrac12$ in the definition of
$\rho_{\sH_{A,v}}(B)$ since the operator $D_B\otimes
\nabla^{\sH_{A,v}}$ is essentially two copies of a usual odd
signature operator, see Remark \ref{CohomTwistSignOpRem} (ii). Now
an immediate---yet interesting---consequence of Theorem
\ref{DaiGen:2} and Corollary \ref{RhoAdiabatic} is

\begin{theorem}\label{RhoGen}
Let $E$ be a flat Hermitian bundle with connection $A$ over the
odd dimensional total space of an oriented fiber bundle
$F\hookrightarrow M\xrightarrow{\pi}B$ of closed manifolds. Then
with respect to every submersion metric
\[
\rho_A(M) = 2^{[\frac{b+1}2]}\int_B \widehat L(TB,\nabla^B)\wedge
\widehat \rho_A + \rho_{\sH_{A,v}}(B) + \gs_A - k\cdot \gs,
\]
where $\gs_A$ and $\gs$ refer to Dai's correction term as in
Definition \ref{DaiCorrect}.
\end{theorem}

\section{Circle Bundles Revisited}

We now want to use Theorem \ref{RhoGen}, to compute the
$\U(1)$-Rho invariant for principal circle bundles over Riemann
surfaces thus giving a different proof of Theorem
\ref{RhoCircBund}. Related results are due to Zhang \cite{Zha94}
and Dai-Zhang \cite{DaiZha95}. In the first reference, the case
of the untwisted spin Dirac operator for circle bundles over
general even dimensional spin manifolds is studied, see also
\cite[Sec. 3 a)]{Goet}. Dai and Zhang use the main result of
\cite{Zha94} and a computation of the Eta invariant for the
untwisted odd signature operator to determine the Kreck-Stolz
invariant for circle bundles. We want to point out that the
strategy in \cite{DaiZha95} for the untwisted odd signature
operator is to apply the signature theorem for manifolds with
boundary to the disk bundle associated to the given circle
bundle. Related discussions can be found in \cite{BeSa, Kom}.
However, this approach does not carry over to non-trivial flat
twisting bundles, since one would need an extension to a
\emph{flat} bundle over the disk bundle, see also Remark
\ref{RhoRem} (iv). In this
respect our strategy is of a purely intrinsic nature.\\

\noindent\textbf{The Bismut-Cheeger Rho Form.} As in Section
\ref{S1Bundles} let $\gS$ be a closed, oriented Riemann surface
of unit volume, and let $S^1\hookrightarrow M\xrightarrow{\pi}
\gS$ a principal $S^1$-bundle over $\gS$ of degree $0\neq l\in
\Z$. Choose a connection $i\go\in\gO^1(M,i\R)$ on $M$, and use
this and the metric on $\gS$ to equip $M$ with a submersion
metric. Let $L_A\to \gS$ be a holomorphic line bundle of degree
$k$, the holomorphic structure being induced by a unitary
connection $A$. As before, we assume that
\[
\lfrac i{2\pi}F_\go =  l\cdot \vol_\gS,\quad\text{and}\quad \lfrac
i{2\pi}F_A = k\vol_\gS.
\]
We endow $L:=\pi^*L_A\to M$ with the flat connection of Lemma
\ref{FlatPullBack}, i.e.,
\[
A_q= \pi^*A - iq\, \go, \quad q:=k/l.
\]

\begin{prop}\label{RhoFormS1Bundle}
The Bismut-Cheeger Rho form associated of the flat line bundle
$L=\pi^*L_A$ is given by
\[
\widehat \rho_{A_q} = 2P_1(q) + l\big(P_2(q)-\lfrac
16\big)\vol_\gS\in \gO^{\ev}(\gS).
\]
Here, $P_1$ and $P_2$ are the first and second periodic Bernoulli
functions of Definition \ref{PeriodicBernoulli}.
\end{prop}

\begin{proof}
First, we infer from \eqref{RhoFormDeg0} that
\[
\big(\widehat \rho_{A_q}\big)_{[0]} = \eta(B_{q,v})-
\eta(B_{0,v}),
\]
where $B_{q,v}$ is the vertical odd signature operator
\[
B_{q,v} = -i(\sL_e -iq) : C^\infty(M,L)\to C^\infty(M,L),
\]
see Proposition \ref{OddSignCirc}. Hence, we can use Remark
\ref{RhoRem} (iii) to deduce that
\[
\big(\widehat \rho_{A_q}\big)_{[0]} = 2P_1(q).
\]
To identify the 2-form part of the Rho form, we recall from
Definition \ref{RescaledBismut} that the rescaled Bismut
superconnection is given by
\[
\B_t = \lfrac{\sqrt t}2 D_{q,v} + \nabla^{\sE,u}  -
\lfrac1{2\sqrt t} c_v(\gO).
\]
First, we need to identify the terms appearing here. In the same
way as in Proposition \ref{OddSignCirc} one finds that the
vertical de Rham operator associated to the connection $A_q$ is
given by
\[
D_{q,v} = -i \tau_v(\sL_e - iq):\gO^\bullet_v(M,L)\to
\gO^\bullet_v(M,L).
\]
Moreover, Lemma \ref{CircleBundleCalc} shows that the mean
curvature of the fiber bundle vanishes, and that the connections
$\nabla^v$ and $\widetilde \nabla^v$ coincide. This implies that
the connection part of $\B$ coincides with the horizontal part of
the exterior differential,
\[
\nabla^{\sE,u}=d_{A_q,h} : \gO^\bullet_v(M,L) \to
\gO^{1,\bullet}(M,L).
\]
Lastly, one easily checks---again as in Proposition
\ref{OddSignCirc}---that
\begin{equation}\label{RhoFormS1Bundle:0}
c_v(\gO) = 2\pi i l  \vol_\gS\wedge \tau_v: \gO^\bullet_v(M,L) \to
\gO^{2,\bullet}(M,L).
\end{equation}
We now claim that in the case at hand
\begin{equation}\label{RhoFormS1Bundle:1}
\B_t^2 = \lfrac{1}4t D_{q,v}^2 + \lfrac12 c_v(\gO)D_{q,v}.
\end{equation}
\begin{proof}[Proof of \eqref{RhoFormS1Bundle:1}]
From \eqref{DvDh:anticommute:1} and Corollary \ref{dSquare} we
have the anti-commutator relations
\[
\big\{\nabla^{\sE,u}, D_{q,v}\big\} = \big\{\nabla^{\sE,u},
c_v(\gO) \big\}= 0.
\]
Moreover, the explicit formul{\ae} above easily yield
\[
c_v(\gO)D_{q,v} = D_{q,v}c_v(\gO) = 2\pi l \vol_{\gS}\wedge
(\sL_e-iq).
\]
Since $\nabla^{\sE,u}$ agrees with $d_{A_q,h}$ we infer from
Corollary \ref{dSquare} that
\[
(\nabla^{\sE,u})^2 = - \big\{d_{q,v},\imu(\gO)\big\},
\]
and one verifies that in the case at hand,
\[
\big\{d_{q,v},\imu(\gO)\big\} = - 2\pi l \vol_{\gS}\wedge
(\sL_e-iq) = - c_v(\gO)D_{q,v}.
\]
Lastly as $\gS$ is 2-dimensional, we have $c_v(\gO)^2=0$. Putting
all pieces together, we obtain
\[
\B_t^2 = \lfrac{1}4t D_{q,v}^2 - \lfrac12 c_v(\gO)D_{q,v} -
\big\{d_{q,v},\imu(\gO)\big\} = \lfrac{1}4t D_{q,v}^2 + \lfrac12
c_v(\gO)D_{q,v}.\qedhere
\]
\end{proof}

Having established the formula in \eqref{RhoFormS1Bundle:1}, we
continue with the proof of Proposition \ref{RhoFormS1Bundle}.
Since $D_{q,v}^2$ and $c_v(\gO)D_{q,v}$ commute, one can use
\eqref{RhoFormS1Bundle:1} and Duhamel's formula to show that
\[
\Tr_v\big[\tau_v\lfrac{d\B_t}{dt}e^{-\B_t^2}\big] =
\Tr_v\Big[\lfrac{1}{4\sqrt
t}\tau_v\big(D_{q,v}+\lfrac{c_v(\gO)}t\big)e^{-\frac{1}4t
D_{q,v}^2}\big(1- \lfrac12 c_v(\gO)D_{q,v}\big)\Big].
\]
According to Definition \ref{BCForm:Odd}, one thus finds that the
2-form part of $\widehat \eta_{A_q}$ is given by
\begin{equation}\label{RhoFormS1Bundle:2}
\begin{split}
(\widehat \eta_{A_q})_{[2]} &= \frac 1{2\pi i} \frac 1{\sqrt\pi}
\int_0^\infty \Tr_v\Big[\lfrac{1}{4\sqrt t}  \tau_vc_v(\gO)
\big(-\lfrac12 D_{q,v}^2+ \lfrac1t \big)e^{-\frac{1}4t D_{q,v}^2}
\Big]dt\\
&= \frac{l\vol_\gS}{4\sqrt{\pi}}\int_0^\infty\left( -
u^{-1/2}\Tr_v\big[D_{q,v}^2e^{-uD_{q,v}^2}\big] + \lfrac 12
u^{-3/2}\Tr_v\big[e^{-uD_{q,v}^2}\big]\right)du.
\end{split}
\end{equation}
Note that in the second equality we have used
\eqref{RhoFormS1Bundle:0} to replace $\tau_vc_v(\gO)$ with $2\pi
i l\vol_\gS$, and then made the substitution $t=4u$. We now
introduce a complex parameter $s$ with $\Re(s)\ge 0$ and define
\[
(\widehat \eta_{A_q})_{[2]}(s) = \frac{l\vol_\gS}{4
\gG\big(\lfrac{s+1}{2}\big)} \int_0^\infty\left( -
u^{\frac{s-1}2}\Tr_v\big[D_{q,v}^2e^{-uD_{q,v}^2}\big] +\lfrac 12
u^{\frac{s-3}2}\Tr_v\big[e^{-uD_{q,v}^2}\big]\right)du,
\]
see Remark \ref{RhoFormS1BundleRem} (i) below. Now for $\Re(s)$
large enough, we can split up the integral and compute using the
Mellin transform that
\[
\begin{split}
\frac{1}{\gG\big(\lfrac{s+1}{2}\big)}\int_0^\infty
u^{\frac{s-1}2}\Tr_v\big[D_{q,v}^2e^{-uD_{q,v}^2}\big]du =
\sum_{\gl\in\spec(D_{q,v}^2)}  \gl^{-\frac{s-1}2},
\end{split}
\]
and
\[
\begin{split}
\frac{1}{\gG\big(\lfrac{s+1}{2}\big)} \int_0^\infty
u^{\frac{s-3}2}\Tr_v\big[e^{-uD_{q,v}^2}\big]du &=
\sum_{\gl\in\spec(D_{q,v}^2)}
\frac{1}{\gG\big(\lfrac{s+1}{2}\big)} \int_0^\infty
u^{\frac{s-3}2}e^{-u\gl}du \\
&= \sum_{\gl\in\spec(D_{q,v}^2)}
\frac{\gG\big(\lfrac{s-1}{2}\big)}{\gG\big(\lfrac{s+1}{2}\big)}
\gl^{-\frac{s-1}2}  = \sum_{\gl\in\spec(D_{q,v}^2)} \lfrac2{s-1}
\gl^{-\frac{s-1}2}.
\end{split}
\]
Now we can use the computation of the spectrum of the vertical
odd signature operator in Lemma \ref{VerticalEigenspaces}. Note,
however, that here we are using $\tau_v D_{q,v}$ which corresponds
to two copies of the operator considered there. We find, again for
$\Re(s)$ large enough, that
\[
\begin{split}
\sum_{\gl\in\spec(D_{q,v}^2)}  \gl^{-\frac{s-1}2} =
\sum_{\gl\in\spec(\tau_v D_{q,v})} |\gl|^{1-s}
= 2 \sum_{\begin{smallmatrix} n \in\Z \\
n\neq q
\end{smallmatrix}}|n -q|^{1-s}  = 2 \widetilde \gz_q(s-1),
\end{split}
\]
where $\widetilde \gz_q$ is the Zeta function in Lemma
\ref{ZEtaCalc}. Hence,
\begin{equation}\label{RhoFormS1Bundle:3}
(\widehat \eta_{A_q})_{[2]}(s) = \frac{l\vol_\gS}2
\frac{2-s}{s-1} \widetilde \gz_q(s-1).
\end{equation}
We know the value of the meromorphic continuation $\widetilde
\gz_q(s)$ to $s=-1$ from Lemma \ref{ZEtaCalc}. Using this we
arrive at
\[
(\widehat \eta_{A_q})_{[2]} =  l P_2(q) \vol_\gS,
\]
so that indeed
\[
(\widehat \rho_{A_q})_{[2]} = (\widehat \eta_{A_q})_{[2]} -
\widehat \eta_{[2]} = l\big(P_2(q)-\lfrac 16\big)\vol_\gS.\qedhere
\]
\end{proof}

\begin{remark}\label{RhoFormS1BundleRem}\quad\nopagebreak
\begin{enumerate}
\item We want to point out that introducing a complex parameter
to split up the sum in \eqref{RhoFormS1Bundle:2} is necessary.
For this note that the individual terms do not give functions
which are holomorphic for $\Re(s)\ge 0$. However---at least in
the case that $q\notin \Z$---their sum is holomorphic for
$\Re(s)\ge 0$, since Lemma \ref{ZEtaCalc} implies that the poles
and zeros of $\frac{2-s}{s-1}$ and $\widetilde \gz_q(s-1)$ in
\eqref{RhoFormS1Bundle:3} cancel each other out.
\item Even though the above computations are very similar to the
ones in the proof of Proposition \ref{EtaTruncSign}, we want to
point out that there is a conceptual difference. Before, we had to
incorporate the operator $D_{A_q,h}$ and work on
$\gO^\bullet(M,L)$, whereas now, we work only on
$\gO_v^\bullet(M,L)$.
\end{enumerate}
\end{remark}

Using Proposition \ref{RhoFormS1Bundle} we can identify the first
term in the general formula of Theorem \ref{RhoGen}. Since $\gS$
is 2-dimensional, we have $\widehat L(T\gS,\nabla^\gS)=1$ so that
\begin{equation}\label{RhoFormS1BundleCor}
2\int_\gS \widehat L(T\gS,\nabla^\gS)\wedge \widehat \rho_{A_q} =
2l\big(P_2(q)-\lfrac 16\big).
\end{equation}

\noindent\textbf{Dai's Correction Term.} We now want to
understand the remaining terms appearing in Theorem \ref{RhoGen}
for the example at hand. Using \eqref{RhoFormS1BundleCor} it is
then immediate that our second proof of Theorem \ref{RhoCircBund}
will follow from the following result.

\begin{prop}\label{TopCorrTermCirc}
Let $l\neq 0$ be the degree of the principal circle bundle, and
let $A_q$ be the flat connection on $L=\pi^*L_A$ as before. Then
$\rho_{\sH_{A_q,v}}(\gS)=0$, and
\[
\gs_{A_q} = \begin{cases} - \sgn (l), &\text{if $A_q$ is the
trivial
connection},\\
\quad\; 0, &\text{if $A_q$ is non-trivial.}
\end{cases}
\]
\end{prop}

\begin{proof}
Recall that $L_\go\to\gS$ denotes the line bundle associated to
the principal bundle $\pi:M\to \gS$, endowed with the connection
$A_\go$ naturally induced by $\go$. As in
\eqref{TruncSignOpKern:1} we can identify
\begin{equation}\label{TopCorrTermCirc:1}
\ker(D_{q,v}) \cong \begin{cases}
\hphantom{\gO^\bullet(\gS,L_B)\otimes\big(}\{0\}, &\text{\rm if
$q\notin \Z$}, \\ \gO^\bullet(\gS,L_B)\oplus\big(
\gO^\bullet(\gS,L_B)\otimes \C\,[\go]\big), &\text{\rm if $q\in
\Z$}.
\end{cases}
\end{equation}
Here, $L_B = L_A\otimes L_\go^{-q}$, which is endowed with the
connection $B=A\otimes 1 +1\otimes qA_\go$. Note that in addition
to \eqref{TruncSignOpKern:1}, the term
$\gO^\bullet(\gS,L_B)\otimes \C\,[\go]$ appears since we do not
restrict $D_{q,v}$ to $\gO^{\bullet,0}(M,L)$.

If $q\notin \Z$, we deduce from \eqref{TopCorrTermCirc:1} that
the bundle of vertical cohomology groups vanishes, which
certainly implies that $\eta\big(D_\gS\otimes
\nabla^{\sH_{A,v}}\big)$ and $\gs_{A_q}$ are both zero. If
$q\in\Z$, we conclude as in the proof of Proposition
\ref{TruncSignOpKern} that $L_B$ is isomorphic to the trivial
line bundle and that the connection $B$ is flat. This implies that
the operator $D_\gS\otimes \nabla^{\sH_{A,v}}$ is unitarily
equivalent to several copies of a twisted odd signature operator
on $\gS$, see Remark \ref{CohomTwistSignOpRem} (ii) and Lemma
\ref{SignCompPart}. Since $\gS$ is even dimensional, we thus infer
that $\eta\big(D_\gS\otimes \nabla^{\sH_{A,v}}\big)=0$. To finish
the proof, it remains to compute $\gs_{A_q}$ for $q\in \Z$.

First of all, we explicitly describe the $E_2$-term of the
Leray-Serre spectral sequence. Using Lemma \ref{SpecSeqLem} we
proceed as in Proposition \ref{TruncSignOpKern} to conclude from
\eqref{TopCorrTermCirc:1} that for $q\in\Z$,
\begin{equation}\label{TopCorrTermCirc:2}
E_{A_q,2}^{\bullet,0} \cong H^\bullet(\gS,L_B),\quad
E_{A_q,2}^{\bullet,1} \cong H^\bullet(\gS,L_B)\otimes \C[\go].
\end{equation}
Now, according to Definition \ref{DaiCorrect}, we have to compute
the signature of the pairing
\begin{equation}\label{TopCorrTermCirc:3}
P_{A_q,2}= Q_{A,2}\big(\,.\,,d_{A_q,2}(\,.\,)\big):
E_{A_q,2}^{0,1}\times E_{A_q,2}^{0,1}\to\C,
\end{equation}
see \eqref{SpecSeqPairOdd}. Here, we have used that concerning the
signature of $P_{A_q,2}$ we can neglect the spaces
$E_{A_q,2}^{1,0}$, since the differential $d_{A_q,2}$ is of
bidegree $(2,-1)$ and thus zero on $E_{A_q,2}^{1,0}$, see Figure
\ref{Fig:SpecSeqCirc}.

\begin{figure}[htbp]
\centering
\includegraphics[width=0.43\linewidth]{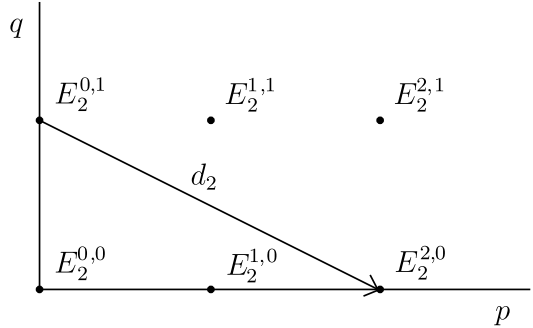}
\caption{The $E_2$-term of the spectral
sequence}\label{Fig:SpecSeqCirc}
\end{figure}

Now, if $A_q$ is a non-trivial connection, the line bundle $L_B =
L_A\otimes L_\go^{-q}$ with its natural connection is not
isomorphic to the trivial flat line bundle. According to
\eqref{TopCorrTermCirc:2} this implies that
$E_{A_q,2}^{0,1}=\{0\}$, see also Proposition
\ref{TruncSignOpKern}. Hence, in this case $\gs_{A_q}=0$. Now we
assume that $A_q$ is the trivial connection, and drop the
subscripts $A_q$ from the notation. Using
\eqref{TopCorrTermCirc:1}, one checks directly that
\[
E_2^{0,1} = \C\,[\go]\quad\text{and}\quad E_2^{2,0} =
\C\,[\vol_\gS].
\]
According to its definition in \eqref{SpecSeqPair}, the pairing
$Q_2$ in \eqref{TopCorrTermCirc:3} is induced by the wedge
product, followed by evaluation on the fundamental class. Using
this, one easily verifies that $Q_2$ satisfies
\[
Q_2: E_2^{0,1}\times E_2^{2,0}\mapsto \C,\quad
Q_2\big([\go],[\vol_\gS]\big) =1.
\]
Moreover, since the differential $d_2$ is naturally induced by
exterior differentiation, one obtains from Proposition
\ref{dSplit} that
\[
d_2:E_2^{0,1}\to E_2^{2,0},  \quad d_2[\go] = [\imu(\gO)\go] =
-2\pi l[\vol_\gS].
\]
According to \eqref{TopCorrTermCirc:3}, this means
\[
P_2\big([\go],[\go]\big) = Q_2\big([\go],d_2[\go]\big) = -2\pi l.
\]
This yields that $\Sign(P_2) = -\sgn (l)$. Since the spectral
sequence collapses at the $E_3$-stage, there are no higher
signatures. This finishes the computation of $\gs_{A_q}$ in the
case that $A_q$ is the trivial connection.
\end{proof}
\cleardoublepage
\chapter{3-dimensional Mapping Tori}\label{3dimMapTor}

After having described the abstract theory leading to a general
formula for Rho invariants of a fiber bundle, we now want to use
this to discuss a second class of examples in detail. In contrast
to the previous class considered, we now reverse the role of fiber
and base, and consider fiber bundles over $S^1$ with fiber a
closed, oriented surface.

We start this chapter to describe a particular way of
constructing submersion metrics. Here, it is convenient to use a
formulation in terms of symplectic forms and almost complex
structures. We then identify the geometric objects in this
setting, most notably the bundle of vertical cohomology groups
and the transgression form of the Bismut superconnection.

Under the assumption that the mapping torus is of finite order,
we derive a formula which expresses the Rho invariant in terms of
Hodge-de-Rham cohomology. Here, we can treat the case of a higher
dimensional gauge group and arbitrary genus of the surface fiber
without effort.

From then on we restrict to the case that the fiber is a
2-dimensional torus. We describe in some detail the geometric
setup, including a discussion of the spectrum of the Laplace
operator on a torus twisted by a flat $\U(1)$-connection. Using
this, we obtain a formula for the Rho form. We then employ ideas
related to the classical Kronecker limit formula, to cast this
expression into a different form which is more accessible for
direct computations.

To obtain explicit formul{\ae} for the Rho invariant, we now have
to distinguish between the cases that the monodromy of the
mapping torus is elliptic, parabolic or hyperbolic. In the
elliptic case, we can specialize the previous result about finite
order mapping tori to obtain a simple formula for the Rho
invariant. The parabolic case turns out to be more involved. Most
notably, the Eta invariant of the odd signature operator with
values in the bundle of vertical cohomology groups has to be
analyzed carefully. Yet, we shall arrive again at a very explicit
formula for $\U(1)$-Rho invariants.

The last part of our discussion will be concerned with the case of
a hyperbolic mapping torus. Here, the difficulty lies in
identifying the Rho form. Generalizing considerations by Atiyah
\cite{Ati87}, we relate this term to the logarithm of a
generalized Dedekind Eta function. We will find that the
transformation property of this logarithm under the action of the
modular group determines the value of the Rho invariant of a
hyperbolic mapping torus. Using a result due to Dieter
\cite{Die59}, we can then express the Rho invariant as the
difference of certain Dedekind sums. Simplifying this expression
we finally arrive at a general formula for $\U(1)$-Rho invariants
in the hyperbolic case as well.

\clearpage

\section{Geometric Preliminaries}\label{3dimMapTorPrelim}

Let $\gS$ be a closed, oriented surface, and let $f\in
\Diff^+(\gS)$ be an orientation preserving diffeomorphism of
$\gS$. Consider the mapping torus
\begin{equation}\label{MapTorDef}
\gS_f := \big(\gS\times \R\big)/\Z,
\end{equation}
where $\Z$ acts on $\gS\times \R$ via
\begin{equation}\label{MapTorDefConvSurf}
k\cdot(x,t) = \big(f^{-k}(x),t+k\big),\quad (x,t)\in \gS\times
\R,\quad k\in\Z.
\end{equation}
Then $\gS_f$ is naturally the total space of a fiber bundle over
$S^1$, see Appendix \ref{MapTorFlatConn}. Moreover, according to
Lemma \ref{MapTorIsotop}, the diffeomorphism type of the mapping
torus depends only on the isotopy class of $f$, i.e., on the
image of $f$ in the mapping class group
\[
\Diff^+(\gS) / \Diff_0(\gS).
\]
We now describe the geometric structure we need on a mapping
torus in some detail. Related material can be found in the
context of Seiberg-Witten equations in \cite[Sec. 8]{Sal}.\\

\noindent\textbf{Moser's Trick.} The freedom of varying $f$ in
its isotopy class allows us to fix particularly convenient choices
for the monodromy map. To exhibit such a choice we need the
following result of \cite{Mos} which is sometimes called
\emph{Moser's trick}, see also \cite[Sec. 3.2]{McDSal}.

\begin{prop}\label{MoserTrick}
Let $\go \in \gO^2(\gS)$ be a symplectic form on the closed,
oriented surface $\gS$, and let $f\in \Diff^+(\gS)$. Then there
exists $\widetilde f\in \Diff^+(\gS)$, isotopic to $f$, with
$\widetilde f^*\go =\go$.
\end{prop}

\begin{proof}[Sketch of proof]
We first note that since $\gS$ is of dimension 2, a symplectic
form $\go$ on $\gS$ is simply a 2-form which is non-degenerate in
the sense that $\int_\gS \go\neq 0$. Moreover, as $f$ is
orientation preserving, we have
\[
[f^*\go]=[\go]\in H^2(\gS,\R).
\]
Thus, there exists $\ga\in \gO^1(\gS)$ with
\[
f^*\go = \go +d\ga.
\]
We now define a time dependent vector field $X:\R\to
C^\infty(\gS,T\gS)$ by requiring that
\[
\imu(X_t)(\go + td\ga) = -\ga,\quad t\in \R.
\]
This is well-defined because $\go$ is non-degenerate. Let
$\gF:\R\to \Diff^+(\gS)$ be the flow uniquely defined by the
initial value problem
\[
\frac{d}{dt} \gF_t = X_t\circ \gF_t,\quad \gF_0=\id_\gS.
\]
Then one checks using Cartan's formula that
\[
\begin{split}
\frac d{dt} \gF^*_t(\go + td\ga) &= \gF_t^* d\ga +
\gF_t^*\big(\sL_{X_t}(\go + td\ga)\big)\\ &= \gF_t^* d\ga+
\gF_t^*\big(d\circ \imu(X_t)(\go +
td\ga)\big)\\
&= \gF_t^* d\ga - \gF_t^* d\ga=0,
\end{split}
\]
where we have used the definition of $X_t$ in the last line.
Hence, $\gF^*_t(\go + td\ga)$ is independent of $t$ and so
\[
\go=\gF_0^*\go = \gF^*_1(\go +d\ga) = \gF^*_1(f^*\go).
\]
Now the claim follows with $\widetilde f:=f\circ \gF_1$.
\end{proof}

\noindent\textbf{Metrics on Mapping Tori.} We can use Proposition
\ref{MoserTrick} to define particular Riemannian metrics on a
mapping torus $\gS_f$. As remarked before, a symplectic form
$\go\in \gO^2(\gS)$ is the same as a volume form on $\gS$. The
space of metrics on $\gS$ with $\go$ as a volume form has the
following description, see e.g. \cite[Sec. 4.1]{McDSal}.

Recall that an almost complex structure is an endomorphism $J\in
C^\infty\big(\gS,\End(T\gS)\big)$ with $J^2=-\Id$. It is called
\emph{$\go$-compatible} if for all $v,w\in T\gS$ with $v\neq 0$,
\[
\go(Jv,Jw)= \go(v,w),\quad \text{and}\quad \go(v,Jv)>0.
\]
Then we get a Riemannian metric with volume form $\go$ by letting
\begin{equation}\label{AlmCompMetDef}
g_J(v,w):= \go(v,Jw), \quad v,w\in T\gS.
\end{equation}
Moreover, the space of metrics with volume form $\go$ is naturally
isomorphic to
\[
\sJ_\go:=\bigsetdef{J\in
C^\infty\big(\gS,\End(T\gS)\big)}{\text{$J$ is an
$\go$-compatible almost complex structure}}.
\]

Given a mapping torus $\gS_f$, we now fix a symplectic form
$\go\in \gO^2(\gS)$ of unit volume. Invoking Proposition
\ref{MoserTrick} and Lemma \ref{MapTorIsotop}, we may assume that
$f^*\go=\go$. Note that this implies that $\sJ_\go$ is invariant
under conjugation with $f_*$. Since $\sJ_\go$ is easily seen to
be path connected---in fact, even contractible---we can choose a
path $J_t:\R\to \sJ_\go$ of $\go$-compatible almost complex
structures on $\gS$ satisfying
\[
J_{t+1} = f_*^{-1}\circ J_t\circ f_*.
\]
Let $g_t$ be the path of Riemannian metrics defined by $J_t$ and
$\go$ as in \eqref{AlmCompMetDef}, and define a metric on
$\gS\times \R$ by
\begin{equation}\label{MapTorMet}
g:= dt\otimes dt + g_t.
\end{equation}
It is immediate from the fact that $f^*\go=\go$ and the
convention \eqref{MapTorDefConvSurf} of how to define the mapping
torus as a quotient of $\gS\times \R$ that $g$ descends to a
Riemannian metric on the mapping torus $\gS_f$.\\

\noindent\textbf{Calculus on Mapping Tori.} We now want to
identify the various quantities introduced in Section
\ref{FiberedCalc}. Since the base of the fiber bundle is
1-dimensional, the curvature form $\gO$ vanishes. Clearly, each
vertical vector field on $\gS_f$ is induced by a path
\begin{equation}\label{VertVectMapTor}
V:\R\to C^\infty(\gS,T\gS),\quad V_{t+1}= f_*V_t,
\end{equation}
and each horizontal vector field $X$ can be identified with
\[
X=\gf_t\pd_t,\quad \text{with}\quad \gf:\R\to C^\infty(\gS),\quad
\gf_{t+1}= \gf_t\circ f.
\]

\begin{lemma}\label{MapTorCalc}
Let $U$ and $V$ be vertical vector fields on $\gS$. With respect
to a metric $g$ as in \eqref{MapTorMet}, the Levi-Civita
connection $\nabla^g$ on $\gS_f$ is given by
\begin{equation}\label{MapTorCalc:1}
\nabla^g_{\pd_t} V\big|_t = \pd_t V_t +\lfrac 12 \dot J_t J_t
V_t,\quad \nabla^g_UV\big|_t = \nabla^{g_t}_{U_t}V_t -\lfrac 12
\go(U_t, \dot J_t V_t)\pd_t,\quad \nabla^g_{\pd_t}\pd_t =0,
\end{equation}
where $\nabla^{g_t}$ is the Levi-Civita connection on $\gS$ with
respect to the metric $g_t$ and $\dot J_t =[\pd_t,J_t]$.
Moreover, the natural vertical connection $\nabla^v$ is given by
\begin{equation}\label{MapTorCalc:2}
\nabla^v_UV\big|_t =  \nabla^{g_t}_{U_t}V_t,\quad \nabla^v_{\pd_t}
V\big|_t = \pd_t V_t +\lfrac 12 \dot J_t J_t V_t.
\end{equation}
In particular, the difference tensor $S$ of Definition
\eqref{SThetaDef} is given by
\[
S(U,V)\big|_t = -\lfrac 12 \go(U_t,\dot J_t V_t)\pd_t,
\]
and the mean curvature form $k_v$ vanishes.
\end{lemma}

\begin{proof}[Sketch of proof]
The description of the Levi-Civita connection easily follows from
the explicit formula, see e.g. \cite[Sec. 1.2]{BGV},
\[
\begin{split}
2g\big(\nabla^g_XY,Z\big) &= g\big([X,Y],Z\big) -
g\big([Y,Z],X\big) + g\big([Z,X],Y\big)\\
&\qquad + Xg(Y,Z) + Yg(Z,X) - Zg(X,Y).
\end{split}
\]
For example,
\[
\begin{split}
2g\big(\nabla^g_{\pd_t} V,U\big) &= g\big([\pd_t,V],U\big) +
g\big([U,\pd_t],V\big) + \pd_t g(V,U)\\
&= g\big(\pd_t V ,U\big) -
g\big(V,\pd_t U\big) + \pd_t \go(V,JU) \\
&= 2 g\big(\pd_t V, U\big) + \go(V,\dot J U)\\
&= 2 g\big(\pd_t V, U\big) + g(J V,\dot J U) = 2 g\big(\pd_t V
+\lfrac 12 \dot J J V,\, U\big),
\end{split}
\]
where we have used that $\dot J_t$ is self-adjoint with respect to
$g_t$ for each $t$. The second equation in \eqref{MapTorCalc:1} is
proven similarly, while the third is clear. Then
\eqref{MapTorCalc:2} follows by taking the vertical projection of
$\nabla^g$. Taking differences yields the formula for $S$. If
$\{e_1,e_2\}$ is a vertical local orthonormal frame, Lemma
\ref{MeanCurvFiber:alt} implies that
\[
k_v(\pd_t)=\lfrac 12 g\big(\dot J J e_1,e_1\big) + \lfrac 12
g\big(\dot J J e_2,e_2\big),
\]
which is zero since $\dot J_t J_t$ is easily seen to be
skew-adjoint with respect to $g_t$.
\end{proof}

\noindent\textbf{Flat Connections and the Bundle of Vertical
Cohomology Groups.} In Appendix \ref{MapTorFlatConn} we have
included a detailed description of the moduli space of flat
$\U(k)$-connections over mapping tori. According to Proposition
\ref{MapTorModuliGeom}, a flat Hermitian vector bundle over
$\gS_f$ is given---up to isomorphism---by a pair $(a,u)$, where
$u\in C^\infty\big(\gS,\U(k)\big)$ is a gauge transformation, and
$a$ is a flat $\U(k)$-connection over $\gS$ satisfying
\begin{equation}\label{MapTorConnTransSurf}
a = u^{-1} (f^*a) u + u^{-1}du.
\end{equation}
We briefly recall from Section \ref{MapTorFlatConn} how this data
defines a flat vector bundle over $\gS_f$. First, let
\[
\widehat f_u:\gS\times \C^k\to \gS\times \C^k,\quad \widehat
f_u(x,z) = \big(f(x),u(x)z\big),
\]
be the automorphism of the trivial bundle over $\gS$ defined by
$u$, see Remark \ref{MapTorBundleTriv}. Then we define a vector
bundle $E_u\to \gS_f$ as the mapping torus
\[
\big((\gS\times \C^k)\times \R\big)/\sim, \quad
\big((x,z),t+1\big)\sim \big(\widehat f_u(x,z),t\big).
\]
Viewing $a$ as a constant path of Lie algebra valued 1-form, $a
\in C^\infty\big(\R,\gO^1(\gS,\cu(k)\big)$, condition
\eqref{MapTorConnTransSurf} ensures that we can define a
connection $A$ on $E_u$. Since $a$ is flat, the same is true for
$A$, see \eqref{MapTorConnDef} and \eqref{MapTorFlatDef}.

\begin{lemma}\label{MapTorCohomBundle}
Let $u\in C^\infty\big(\gS,\U(k)\big)$ be a gauge transformation
defining a bundle $E_u\to \gS_f$, and let $A$ be a flat
$\U(k)$-connection over $T^2_M$ defined by a pair $(a,u)$
satisfying \eqref{MapTorConnTransSurf}.
\begin{enumerate}
\item The bundle automorphism $\widehat f_u :\gS\times \C^k\to
\gS\times \C^k$ induces an isomorphism
\[
\widehat f_u^*: H^\bullet(\gS,E_a) \to H^\bullet(\gS,E_a).
\]
\item Let $\nabla^{H_{A,v}}$ be the natural flat connection on
the bundle $H^\bullet_{A,v}(\gS_f)\to S^1$ of vertical cohomology
groups, see Definition \ref{FlatVertCohom}. Then its holonomy
representation is given by
\[
\hol_{\nabla^{H_{A,v}}}:\pi_1(S^1)\to
\GL\big(H^\bullet(\gS,E_a)\big),\quad
\hol_{\nabla^{H_{A,v}}}(\gamma) = \widehat f_u^*,
\]
where $\gamma\in \pi_1(S^1)$ is the canonical generator.
\end{enumerate}
\end{lemma}

\begin{proof}
The map $\widehat f_u$ acts on $\gO^\bullet(\gS,\C^k)$ via
\[
\widehat f_u^*\ga = u^{-1} f^*\ga,\quad \ga \in
\gO^\bullet(\gS,\C^k).
\]
Condition \eqref{MapTorConnTransSurf} is easily seen to be
equivalent to $\widehat f_u^* \circ d_a = d_a\circ \widehat
f_u^*$. This implies that $\widehat f_u^*$ descends to an
isomorphism on cohomology, which proves part (i). As in
\eqref{FormsMapTor}, we have the following identification of the
space of vertical, $E_u$-valued differential forms
\begin{equation}\label{VerticalFormsMapTor}
\gO^\bullet_v(\gS_f,E_u) = \bigsetdef{\ga_t:\R\to
\gO^\bullet(\gS,\C^k)}{\ga_{t+1} = \widehat f^*_u\ga_t}.
\end{equation}
With respect to this identification, the vertical differential
$d_{A,v}$ coincides with $d_a$ applied pointwise for each $t$.
Moreover, since the connection $A$ has no $dt$ component, the
horizontal differential
\[
d_{A,h}:\gO^\bullet_v(\gS_f,E_u)\to \gO^{1,\bullet}(\gS_f,E_u)
\]
is given with respect to \eqref{VerticalFormsMapTor} by
\[
d_{A,h}\ga_t = dt\wedge \pd_t \ga_t.
\]
From these observations one finds that the space of sections of
$H^\bullet_{A,v}(\gS_f)\to S^1$ can be described as
\[
C^\infty\big(S^1,H_{A,v}(\gS_f)\big) = \bigsetdef{[\ga]_t:\R\to
H^\bullet(\gS,E_a)}{[\ga]_{t+1} = f^*_u[\ga]_t}.
\]
Now it is immediate from the definition of the connection
$\nabla^{H_{A,v}}$ that the holonomy representation has the
claimed form.
\end{proof}

\noindent\textbf{The Bismut Superconnection.} To identify the
terms appearing in Dai's adiabatic limit formula, we first need to
understand the connection $\nabla^{A,v}$ acting on forms, and
then the Bismut superconnection in the setting at hand.

\begin{lemma}\label{VerticalConnMapTor}
With respect to the identification \eqref{VerticalFormsMapTor},
the connection $\nabla^{A,v}$ acting on
$\gO^\bullet_v(\gS_f,E_u)$ is given by
\begin{equation*}
\nabla^{A,v}_{\pd_t} = \pd_t - \lfrac 12 \dot \tau_v \tau_v,\quad
\nabla^{A,v}_V\big|_t  = \nabla^{a,g_t}_{V_t},
\end{equation*}
where $V$ is a vertical vector field, and $\nabla^{a,g_t}$
denotes the connection on $\gO^\bullet(\gS,\C^k)$ induced by $g_t$
and $a$, acting pointwise for each $t$.
\end{lemma}

\begin{proof}
The assertion about $\nabla^{A,v}_V$ is true by definition of
$\nabla^{A,v}$ and the fact that $A$ is independent of $t$.
Concerning $\nabla^{A,v}_{\pd_t}$, we observe that $d_{A,h}=
dt\wedge \pd_t$, which implies that
\begin{equation}\label{CanVerticalConnMapTor}
\widetilde \nabla^{A,v}_{\pd_t} = \pd_t.
\end{equation}
Since the mean curvature form vanishes, we can use
\eqref{NablaNat} and \eqref{LieTauComm} to deduce that
\[
\nabla^{A,v}_{\pd_t} = \pd_t + \lfrac 12 \tau_v[\pd_t,\tau_v] =
\pd_t - \lfrac 12 \dot \tau_v \tau_v,
\]
where we have used that $\dot \tau_v = [\pd_t,\tau_v]$ and
$\dot\tau_v\tau_v = -\tau_v \dot\tau_v$.
\end{proof}

Using Lemma \ref{VerticalConnMapTor} we can now find a more
explicit expression for the Bismut superconnection and its
transgression form.

\begin{prop}\label{MapTorBCForm}
The rescaled Bismut superconnection associated to $A$ is given
for $s\in (0,\infty)$ by
\[
\B_s = \lfrac {\sqrt s}2 D_{A,v} + dt\wedge \big(\pd_t -\lfrac 12
\dot \tau_v \tau_v\big): \gO^\bullet_v(\gS_f,E_u) \to
\gO^\bullet(\gS_f,E_u),
\]
where $D_{A,v}$ is the vertical de Rham operator. Moreover,
\begin{equation*}
\begin{split}
\ga(\B_s) = \frac i{16\pi} dt\wedge
\Tr_v\Big(\dot\tau_v\big(d_{A,v}^td_{A,v}-d_{A,v}d_{A,v}^t\big)e^{-\frac
s4 D_{A,v}^2}\Big) \in \gO^1(S^1).
\end{split}
\end{equation*}
\end{prop}

\begin{proof}
The formula for $\B_s$ follows immediately from Lemma
\ref{VerticalConnMapTor} and the fact that $k_v$ and $\gO$ are
zero, see \eqref{BismutSuperconnConn} and Definition
\ref{RescaledBismut}. For the second assertion, we drop the
connection $A$ from the notation. Since the base is
1-dimensional, and $dt$ anti-commutes with $D_v$, one finds that
\[
\B_s^2\big|_{s=4r} = r D_v^2 + \sqrt r dt\wedge
\big[\nabla^v_{\pd_t},D_v\big].
\]
As in \cite[Lem. 9.42]{BGV} and \cite[Thm. 3.3]{BF2} an
application of Duhamel's formula then yields that
\[
e^{-r D_v^2 - \sqrt r dt\wedge[\nabla^v_{\pd_t},D_v]} = e^{-r
D_v^2} - \sqrt r\, dt\wedge  \int_0^1 e^{-r' r
D_v^2}\big[\nabla^v_{\pd_t},D_v\big] e^{-(1-r')r D_v^2}dr'.
\]
where the higher correction terms vanish, again since the base is
1-dimensional. Therefore,
\[
\begin{split}
\Str_v\Big(\lfrac{d\B_s}{ds}e^{-\B_s^2}\Big)\Big|_{s=4r}&= \lfrac
1{8\sqrt r} \Str_v\Big(D_ve^{-r D_v^2 - \sqrt r dt\wedge
[\nabla^v_{\pd_t},D_v]}\Big)\\ &= \lfrac 1{8\sqrt r}
\Str_v\big(D_v e^{-rD_v^2}\big) + \lfrac 18 dt\wedge
\Str_v\Big(D_v[\nabla^v_{\pd_t},D_v] e^{-r D_v^2}\Big).
\end{split}
\]
Note that a factor of $-1$ enters in front of the second term
since we have interchanged $dt$ and $D_v$. Since $\tau_v$
anti-commutes with $D_v$, the first term vanishes. To simplify
the second term, we also drop the sub- and subscripts $v$ from
the notation. Then one verifies using Lemma
\ref{VerticalConnMapTor} and the relations $D=d - \tau d \tau $,
$[\pd_t,\tau]=\dot \tau$ as well as $\dot\tau \tau=-\tau\dot\tau$
that
\[
\tau [\nabla_{\pd_t},D] = -\lfrac 12\big(\tau [\dot\tau\tau,d] -
[\dot\tau\tau,d]\tau \big) = \lfrac 12 \big(\dot \tau d + \tau
d\dot \tau\tau +\dot\tau \tau d\tau -d\dot\tau\big)
\]
Now, since $\tau D = -D \tau$, we can use the trace property and
the fact that $e^{-rD^2}$ is a semi-group of smoothing operators
to find that
\[
\Tr\big(\tau d\dot \tau\tau D e^{-rD^2}\big) = - \Tr\big( d\dot
\tau D e^{-rD^2}\big),\quad \Tr\big(\dot \tau \tau d \tau D
e^{-rD^2}\big) = \Tr\big(\dot\tau d  D e^{-rD^2}\big).
\]
Thus, by repeatedly making use of the trace property, one obtains
\[
\begin{split}
\Str\big(D[\nabla_{\pd_t},D]e^{-rD^2}\big) & =
\Str\big([\nabla_{\pd_t},D]De^{-rD^2}\big)\\
&= -\Tr\big((d\dot\tau -\dot\tau d)De^{-rD^2}\big)
=-\Tr\big(\dot\tau(d^td-dd^t)e^{-rD^2}\big).
\end{split}
\]
Recalling the normalization factor in Definition
\ref{TransgressionForm}, we arrive at the claimed formula for the
transgression form.
\end{proof}

\section{Finite Order Mapping Tori}\label{FiniteMapTor}

Proposition \ref{MapTorBCForm} shows that if we can achieve that
$\dot \tau_v\equiv 0$, then the Bismut-Cheeger Eta form associated
to a flat connection over $\gS_f$ vanishes. From the discussion in
Section \ref{3dimMapTorPrelim} we know that $\dot \tau_v\equiv 0$
is equivalent to finding an $f$-invariant metric $g_\gS$ on $\gS$.
Clearly, such a metric will not exist for arbitrary $f\in
\Diff^+(\gS)$.

\begin{lemma}\label{f:InvMet}
If $f\in \Diff^+(\gS)$ is of finite order $n$, there exists a
metric $g_\gS$ of unit volume with $f^*g_\gS = g_\gS$.
\end{lemma}

\begin{proof}
Choose an arbitrary metric $g_\gS$ on $\gS$ of unit volume and
define
\[
\widetilde g_\gS := \lfrac 1n \sum_{j=0}^{n-1} (f^j)^* g_\gS.
\]
Then $\widetilde g_\gS$ is again a metric of unit volume.
Moreover, since $f^n=\id_\gS$, one finds that indeed
$f^*\widetilde g_\gS=\widetilde g_\gS$.
\end{proof}

\begin{remark}\label{HurwitzRem}\quad\nopagebreak
\begin{enumerate}
\item A metric $g_\gS$ defines an almost complex structure on
$\gS$ which is integrable, see Section \ref{HolomAspects}. If
$g_\gS$ is $f$-invariant for some $f\in \Diff^+(\gS)$, then $f$
is holomorphic with respect to the complex structure defined by
$g_\gS$.
\item It can be shown that if there exists an $f$-invariant metric
$g_\gS$, then the mapping class $[f]\in
\Diff^+(\gS)/\Diff_0(\gS)$ is necessarily of finite order. If the
genus of $\gS$ is 0 or 1, this can be checked directly, see
Proposition \ref{FixPointClass} below for $\gS=T^2$. For higher
genera, one can use for example the holomorphic description, and
invoke the Riemann-Hurwitz formula, see \cite[Ch. V]{FarKra}.
\end{enumerate}
\end{remark}

\noindent\textbf{Rho Invariants of a Finite Order Mapping Torus.}
Assume from now on that $f\in \Diff^+(\gS)$ is of finite order,
and that an $f$-invariant metric $g_\gS$ has been chosen. We also
fix a flat connection $A$ on a Hermitian vector bundle $E_u\to
\gS_f$ defined by a pair $(a,u)$ satisfying $\widehat f_u^*a=a$
as in \eqref{MapTorConnTransSurf}. In Lemma
\ref{MapTorCohomBundle} we have given a description of the bundle
of vertical cohomology groups in terms of de Rham cohomology.
However, the Rho invariant of the bundle of vertical cohomology
groups in Definition \ref{RhoFormDef}---which appears in the
formula for the Rho invariant in Theorem \ref{RhoGen}---is
defined using the Hodge theoretic description. In the case of a
finite order mapping torus, we have the following extension of
Lemma \ref{MapTorCohomBundle}.

\begin{lemma}\label{FiniteMapTorCohomBundle}\quad\nopagebreak
\begin{enumerate}
\item With respect to the induced metric on $\sH^\bullet(\gS,E_a)$,
the bundle map $\widehat f_u$ defines an isometry
\[
\widehat f_u^* : \sH^\bullet(\gS,E_a) \to
\sH^\bullet(\gS,E_a),\quad \sH^\bullet(\gS,E_a) =
\ker(d_a+d_a^t)\subset \gO^\bullet(\gS,\C^k).
\]
The splitting into $\pm 1$-eigenspaces of $\tau_\gS$,
\[
\sH^\bullet(\gS,E_a) = \sH^+(\gS,E_a) \oplus \sH^-(\gS,E_a),
\]
is invariant with respect to $\widehat f_u^*$.
\item The flat
connection $\nabla^{H_{A,v}}$ on the bundle of vertical
cohomology groups is compatible with the metric
$\scalar{.}{.}_{\sH_{A,v}}$ of Definition \ref{CohomBundleMet},
and its holonomy representation is given by
\[
\hol_{\nabla^{H_{A,v}}}(\gamma) = \widehat f_u^*\in
\GL\big(\sH^\bullet(\gS,E_a)\big),
\]
where $\gamma\in \pi_1(S^1)$ is the canonical generator.
\end{enumerate}
\end{lemma}

\begin{proof}
Since $f$ is an isometry with respect to the metric $g_\gS$, the
pullback
\[
f^*:\gO^\bullet(\gS)\to \gO^\bullet(\gS)
\]
commutes with the chirality operator $\tau_\gS$. Moreover, as
$\widehat f_u^*=u^{-1}f^*$, the same is true for $\widehat f_u^*$.
This, and the fact that $f_u^*\circ d_a = d_a\circ \widehat
f_u^*$ implies part (i). Since $\dot\tau\equiv 0$, we can deduce
from \eqref{LieTauComm} and Proposition \ref{FlatSignBundle} that
the connection $\nabla^{H_{A,v}}$ is indeed compatible with the
metric $\scalar{.}{.}_{\sH_{A,v}}$. Moreover, if
\[
\Psi: \sH^\bullet(\gS,E_a) \to H^\bullet(\gS,E_a)
\]
denotes the Hodge-de-Rham isomorphism, the diagram
\[
\begin{CD}
\sH^\bullet(\gS,E_a)  @>{\widehat f_u^*}>> \sH^\bullet(\gS,E_a) \\
@V{\Psi}VV @V{\Psi}VV \\
H^\bullet(\gS,E_a) @>{\widehat f_u^*}>> H^\bullet(\gS,E_a),
\end{CD}
\]
is commutative. This is because it is naturally induced by
$\widehat f_u^*: \gO^\bullet(\gS,\C^k)\to \gO^\bullet(\gS,\C^k)$.
Hence, the description of the holonomy representation follows
from Lemma \ref{MapTorCohomBundle}.
\end{proof}

%
%

\begin{theorem}\label{RhoFiniteMapTor}
Let $f\in \Diff^+(\gS)$ be of finite order. Let $A$ be a flat
$\U(k)$-connection over $\gS_f$, defined by a pair $(a,u)$ of flat
connection and gauge transformation over $\gS$ satisfying
$\widehat f_u^* a=a$. Then with respect to every $f$-invariant
metric $g_\gS$ on $\gS$
\[
\begin{split}
\rho_A(\gS_f) =& 2 \tr\log\big[ \widehat
f_u^*|_{\sH^+(\gS,E_a)\cap \gO^1}\big] - \rk\big[(\widehat
f_u^*-\Id)|_{\sH^+(\gS,E_a)\cap \gO^1}\big]\\ &\quad - 2
\tr\log\big[ \widehat f_u^*|_{\sH^-(\gS,E_a)\cap \gO^1}\big]  +
\rk\big[(\widehat f_u^*-\Id)|_{\sH^-(\gS,E_a)\cap \gO^1}\big] \\
&\quad -4k \tr\log\big[f^*|_{\sH^+(\gS)\cap \gO^1}\big] + 2k
\rk\big[(f^*-\Id)|_{\sH^+(\gS)\cap \gO^1}\big].
\end{split}
\]
Here, ``$\tr\log$'' is defined for a unitary map $T\in\U(n)$ as
\begin{equation*}
\tr\log T:= \sum_{j=1}^n \gt_j \in \R,
\end{equation*}
where $e^{2\pi i\gt_j}$ are the eigenvalues of $T$, and where we
require $\gt_j\in [0,1)$.
\end{theorem}

\begin{proof}
As we have already noted before, Proposition \ref{MapTorBCForm}
implies that the Bismut-Cheeger Rho form vanishes, because
$g_\gS$ is $f$-invariant. Moreover, the Leray-Serre spectral
sequence associated to $A$---respectively the trivial
connection---collapses at the $E_2$-stage, since the base is
1-dimensional, see Theorem \ref{SpecSeqThm}. This implies that
Dai's correction term in Definition \ref{DaiCorrect} vanishes.
Hence, Theorem \ref{RhoGen} yields
\begin{equation}\label{RhoFiniteMapTor:0}
\rho_A(\gS_f) = \rho_{\sH_{A,v}}(S^1) = \lfrac12
\eta\big(D_{S^1}\otimes \nabla^{\sH_{A,v}}\big) - k \cdot
\lfrac12\eta\big(D_{S^1}\otimes \nabla^{\sH_v}\big),
\end{equation}
where $D_{S^1}\otimes \nabla^{\sH_{A,v}}$ and $D_{S^1}\otimes
\nabla^{\sH_v}$ are as in Definition \ref{CohomTwistSignOp}. We
have seen in Lemma \ref{FiniteMapTorCohomBundle} that the flat
connection $\nabla^{\sH_v}$ is unitary with respect to the metric
$\scalar{.}{.}_{\sH_v}$. Hence, we deduce from Lemma
\ref{SignCompPart} that
\begin{equation}\label{RhoFiniteMapTor:1}
D_{S^1}\otimes \nabla^{\sH_v} = \begin{pmatrix} D_{S^1}\otimes
\nabla^{\sH_v,+} &0\\ 0 & -D_{S^1}\otimes \nabla^{\sH_v,-}
\end{pmatrix},
\end{equation}
where
\begin{equation}\label{RhoFiniteMapTor:2}
D_{S^1}\otimes \nabla^{\sH_v,\pm}= \tau_{S^1}
d_{\nabla^{\sH_v,\pm}} + d_{\nabla^{\sH_v,\pm}} \tau_{S^1},
\end{equation}
and $\nabla^{\sH_v,\pm}$ denotes the restriction of
$\nabla^{\sH_v}$ to $\sH^\pm_v(\gS_f)$. The same formula holds for
$\sH^\bullet_v(\gS_f)$ replaced with $\sH_{A,v}^\bullet(\gS_f)$.
Let us abbreviate
\begin{equation}\label{RhoFiniteMapTor:3}
B_{\pm} := \tau_{S^1} d_{\nabla^{\sH_v,\pm}}
:C^\infty\big(S^1,\sH^\pm_v(\gS_f)\big)\to
C^\infty\big(S^1,\sH^\pm_v(\gS_f)\big),
\end{equation}
and define $B_{A,\pm}$ correspondingly. Then
\eqref{RhoFiniteMapTor:0} and \eqref{RhoFiniteMapTor:1} show that
\begin{equation}\label{RhoFiniteMapTor:4}
\rho_A(\gS_f) = \eta(B_{A,+}) - \eta(B_{A,-})- k\cdot \eta(B_+) +
k\cdot \eta(B_-).
\end{equation}
Note that the factors $\frac 12$ disappear as the operators in
\eqref{RhoFiniteMapTor:2} are equivalent to two copies of the
operators in \eqref{RhoFiniteMapTor:3}, see Remark
\ref{CohomTwistSignOpRem} (ii). We deduce from Lemma
\ref{FiniteMapTorCohomBundle} (i) that $\widehat f_u^*$ restricts
to a unitary map on $\sH^\pm(\gS,E_a)$. Hence, we can find a
basis of $\sH^\pm(\gS,E_a)$ such that
\[
\widehat f_u^*|_{\sH^\pm(\gS,E_a)} = \diag\big(e^{2\pi i
\gt_1^\pm},\ldots, e^{2\pi i \gt_n^\pm}\big),
\]
where $\gt_j^\pm \in [0,1)$, and $n=\dim \sH^+(\gS,E_a)= \dim
\sH^-(\gS,E_a)$. Note that the equality of dimensions follows from
the fact that $\Sign_a(\gS)=0$, which in turn is a consequence of
the signature formula, see Theorem \ref{TwistSignIndThm}. Now
Lemma \ref{FiniteMapTorCohomBundle} (ii) yields that the
restriction $\widehat f_u^*|_{\sH^\pm(\gS,E_a)}$ defines the
holonomy representation of the connection
$\nabla^{\sH_{A,v},\pm}$. From this we deduce as in Remark
\ref{RhoRem} (iii) that
\[
\begin{split}
\eta(B_{A,+}) - \eta(B_{A,-}) =& \sum_{\gt_j^+\neq 0}
(2\gt_j^+ - 1) - \sum_{\gt_j^-\neq 0} (2\gt_j^- - 1)\\
=& 2 \tr\log\big[ \widehat f_u^*|_{\sH^+}\big] -
\rk\big[(\widehat f_u^*-\Id)|_{\sH^+}\big]\\ &\quad - 2
\tr\log\big[ \widehat f_u^*|_{\sH^-}\big]  +  \rk\big[(\widehat
f_u^*-\Id)|_{\sH^-}\big],
\end{split}
\]
where for convenience we have abbreviated $\sH^\pm =
\sH^\pm(\gS,E_a)$ in the last equality. Continuing with obvious
abbreviations, we now decompose
\[
\sH^\pm = \sH^\pm\cap (\gO^0\oplus \gO^2) \oplus \sH^\pm\cap
\gO^1.
\]
Each element of $\sH^\pm\cap \big(\gO^0\oplus \gO^2\big)$ is of
the form $\gf \pm \tau_\gS\gf $ with $\gf\in \sH^0$. Hence, we
have a natural isomorphism (compare also with the proof of
Proposition \ref{Sign=Ind}),
\[
\sH^+\cap \big(\gO^0\oplus \gO^2\big) \xrightarrow{\cong}
\sH^-\cap \big(\gO^0\oplus \gO^2\big),\quad \gf+\tau_\gS\gf
\mapsto \gf - \tau_\gS\gf.
\]
Since $\widehat f_u^*$ commutes with $\tau_\gS$, we can conclude
that
\[
\begin{split}
2 \tr\log\big[ &\widehat f_u^*|_{\sH^+\cap
(\gO^0\oplus\gO^2)}\big] - \rk\big[(\widehat
f_u^*-\Id)|_{\sH^+\cap (\gO^0\oplus\gO^2)}\big]\\
&\qquad = 2 \tr\log\big[ \widehat f_u^*|_{\sH^-\cap
(\gO^0\oplus\gO^2)}\big] - \rk\big[(\widehat
f_u^*-\Id)|_{\sH^-\cap (\gO^0\oplus\gO^2)}\big].
\end{split}
\]
Therefore,
\[
\begin{split}
\eta(B_{A,+}) - \eta(B_{A,-}) =& 2 \tr\log\big[ \widehat
f_u^*|_{\sH^+\cap\gO^1}\big] - \rk\big[(\widehat
f_u^*-\Id)|_{\sH^+\cap \gO^1}\big]\\ &\quad - 2 \tr\log\big[
\widehat f_u^*|_{\sH^-\cap\gO^1}\big]  +  \rk\big[(\widehat
f_u^*-\Id)|_{\sH^-\cap\gO^1}\big].
\end{split}
\]
This identifies the twisted terms appearing in the formula of
Theorem \ref{RhoFiniteMapTor}. In the case that $a$ is the trivial
connection and $u\equiv 1$, we can simplify this further. Since
we are considering complex valued forms, we have a conjugation
\[
\sH^1(\gS)\to \sH^1(\gS),\quad \ga\mapsto \Bar \ga.
\]
The chirality operator $\tau_\gS$ is readily seen to anti-commute
with conjugation. This yields an anti-linear isomorphism
\[
\sH^+(\gS)\cap \gO^1(\gS) \xrightarrow{\cong} \sH^-(\gS)\cap
\gO^1(\gS).
\]
Since $f^*$ is the complex linear extension of a real
automorphism, it commutes with conjugation. From this one readily
deduces that the eigenvalues of $f^*|_{\sH^-\cap \gO^1}$ are
complex conjugate to the eigenvalues of $f^*|_{\sH^+\cap \gO^1}$.
By checking the definition of ``$\tr\log$'' carefully, one
conlcudes
\[
\tr\log\big[ f^*|_{\sH^-\cap \gO^1}\big] =
\rk\big[(f^*-\Id)|_{\sH^+\cap \gO^1}\big] - \tr\log\big[
f^*|_{\sH^+\cap \gO^1}\big].
\]
Now, as $\rk\big[(f^*-\Id)|_{\sH^-\cap \gO^1}\big] =
\rk\big[(f^*-\Id)|_{\sH^+\cap \gO^1}\big]$, we finally get
\[
\begin{split}
\eta(B_+) - \eta(B_-) &= 4\tr\log\big[ f^*|_{\sH^+\cap
\gO^1}\big] -2 \rk\big[(f^*-\Id)|_{\sH^+\cap \gO^1}\big],
\end{split}
\]
which is precisely the untwisted term in the claimed formula.
\end{proof}

\begin{remark*}\quad\nopagebreak
\begin{enumerate}
\item The last step of the above proof is essentially equivalent
to observing that $f^*$ acting on $\sH^1(\gS)$ is the
complexification of a symplectic map. This explains why the
eigenvalues come in conjugate pairs. Clearly, this is no longer
true if we consider $u^{-1}f^*$ with $u\in \U(1)$, which also
defines a flat connection over $\gS_f$, see also Theorem
\ref{RhoFiniteTorusBundle} below.
\item In a similar direction, if the connection $a$ is
non-trivial, complex conjugation gives rise to an anti-linear
isomorphism
\[
\sH^+(\gS,E_a)\xrightarrow{\cong} \sH^-(\gS,E_{\Bar a}).
\]
From this, one can relate the eigenvalues of $\widehat f_u^*$
acting on $\sH^+(\gS,E_a)$ with the eigenvalues of $\widehat
f_{\Bar u}^*$ acting on $\sH^-(\gS,E_{\Bar a})$. However, this
only simplifies the formula of Theorem \ref{RhoFiniteMapTor} in
the case that $a$ and $u$ are real, in the sense that they arise
from an $\operatorname{O}(k)$-structure.
\item Although Theorem \ref{RhoFiniteMapTor} expresses the Rho
invariants of $\gS_f$ in terms of Hodge-de-Rham cohomology of
$\gS$, it is only an intermediate step to an expression in
completely topological terms. The next step would be to use the
ideas of the Atiyah-Bott fixed point formula---see \cite[Sec.
6.2]{BGV}---to relate the traces appearing in Theorem
\ref{RhoFiniteMapTor} to the fixed point data of $f$. This, in
turn, can be expressed in terms of the Seifert invariants of the
finite order mapping torus. We refer to \cite[Sec. 5]{And} and
\cite[Sec. 2.2]{Fuj} for a discussion of these ideas in the
context of the determinant line bundle over the moduli space of
flat connections associated to a finite order mapping torus.
\item In \cite{Mor} an interpretation of the untwisted Eta
invariant for finite order mapping tori is given in terms of
Meyer's cocycle for the mapping class group, see \cite{Mey73}. It
follows from the proof of Theorem \ref{RhoFiniteMapTor} that the
adiabatic limit of the untwisted Eta invariant is given by
\begin{equation}\label{UntwistedEtaFiniteMapTor}
4\tr\log\big[f^*|_{\sH^+(\gS)\cap\gO^1(\gS)}\big] - 2
\rk\big[(f^*-\Id)|_{\sH^+(\gS)\cap\gO^1(\gS)}\big].
\end{equation}
It would be interesting to relate this to the main result of
\cite{Mor}.
\end{enumerate}
\end{remark*}

%

\section[Torus Bundles over $S^1$, General Setup]{Torus Bundles
over $\boldsymbol{S^1}$, General Setup}\label{TorusBundlesGen}

We now consider the case that $\gS=T^2$ is the 2-dimensional
torus. In the same setting, Atiyah \cite{Ati87} studies a rich
interplay between the untwisted Eta invariant and other
topological invariants, the Dedekind Eta function and also number
theoretical $L$-series. As a tool Atiyah also makes intensive use
of the idea of adiabatic limits, and much of our discussion is
influenced by the treatment in \cite{Ati87}. We shall restrict to
the case of $\U(1)$-connections, which already contains many
important ideas. However, in view of the computations of
Chern-Simons invariants for torus bundles in \cite{Jef92, KK90}
the generalization to higher gauge groups would be extremely
interesting. \\

\subsection[Geometry of Torus Bundles over
$S^1$]{Geometry Torus Bundles over $\boldsymbol{S^1}$}

\noindent\textbf{Complex Structures on $\boldsymbol{T^2}$.} We
fix the standard torus $T^2=\R^2/\Z^2$, endowed with the volume
form induced by $\go = dx\wedge dy$. As in Section
\ref{3dimMapTorPrelim} we are interested in the space of all
metrics which have $\go$ as a volume form. Equivalently, we need
to understand the space $\sJ_\go$ of all $\go$-compatible almost
complex structures. It is well known that $\sJ_\go$ is the
Teichm\"{u}ller space of $T^2$, i.e., the upper half
plane\footnote{We are using the letter $\gs$ for elements in $\H$
rather than the more common letter $\tau$ to avoid confusion with
the chirality operator.}
\[
\H:=\bigsetdef{\gs=\gs_1+ i\gs_2\in \C}{ v>0},
\]
see \cite[Thm. 2.7.2]{JosRS}. For definiteness, we will use the
following explicit isomorphism. Note that each almost complex
structure $J\in \sJ_\go$ can be identified with a matrix in
$\operatorname{M}_2(\R)$ because the tangent space of $T^2$ is
canonically isomorphic to $\R^2$.

\begin{lemma}\label{TeichmuellerAlmComp}
The map
\[
\gF:\H\to \sJ_\go,\quad \gF(\gs)= \frac 1{\gs_2}
\begin{pmatrix} -\gs_1 &-|\gs|^2\\ 1 &\gs_1
\end{pmatrix},\quad \gs=\gs_1+i\gs_2,
\]
is a bijection. The metric on $T^2$ defined by $\gF(\gs)$ as in
\eqref{AlmCompMetDef} is given with respect to the standard
coordinate basis as
\[
g_\gs  = \frac 1{\gs_2} \Big(dx\otimes dx + \gs_1(dx\otimes
dy+dy\otimes dx) + |\gs|^2 dy\otimes dy\Big).
\]
\end{lemma}

\begin{proof}
If we identify $J \in \sJ_\go$ with an element in
$\operatorname{M}_2(\R)$, one easily checks that $J^2=-\Id$ is
equivalent to $\det(J)=1$ and $\tr(J)=0$. Hence there exist
$r,s,t\in \R$ such that
\begin{equation}\label{TeichmuellerAlmComp:1}
J= \begin{pmatrix} -r &t\\ s &r
\end{pmatrix}, \quad r^2+st =-1.
\end{equation}
Let $J_0$ be the almost complex structure which induces the
standard scalar product on $\R^2$, i.e.,
\[
J_0 := \begin{pmatrix} 0 &-1\\ 1 &0
\end{pmatrix},\quad g_0= \go(.,J_0.)=dx\otimes dx + dy\otimes dy.
\]
Then one verifies that $J$ is $\go$-compatible if and only if
$-J_0 J$ is positive definite and symmetric. Now
\begin{equation}\label{TeichmuellerAlmComp:2}
-J_0 J = \begin{pmatrix} s &r\\ r &-t
\end{pmatrix},
\end{equation}
and so \eqref{TeichmuellerAlmComp:1} implies that $-J_0 J$ is
positive definite if and only if $s>0$. Using this one finds that
$\gF$ is well-defined with inverse given by
\[
\gF^{-1}(J) = \frac 1s (r+i).
\]
Moreover, \eqref{TeichmuellerAlmComp:2} relates $\gF(\gs)$ and
the associated metric $g_\gs$ and easily yields the second claim.
\end{proof}

\begin{remark}\label{EllipticCurve}
Viewed from a complex analytic perspective, the above
identification might seem a bit cumbersome. As in \cite[Sec.
2.7]{JosRS}, any $\gs\in \H$ defines a lattice
\[
\gL(\gs):= \bigsetdef{m+n\gs}{(m,n)\in\Z^2}\subset \C,
\]
and the quotient torus $\C/\gL(\gs)$ is naturally a complex
manifold, with complex structure induced by the one of the
complex plane, and metric induced by $\frac 1{\gs_2}(dx^2+dy^2)$.
Note that one has to divide by $\gs_2$ to get a metric of unit
volume. It is easy to check that
\[
\psi_\gs: \R^2\to \C,\quad (x,y)\mapsto x+\gs y,
\]
descends to an isometry
\begin{equation}\label{EllitpicCurveIsom}
\psi_\gs:(T^2,g_\gs)\to \C/\gL(\gs),
\end{equation}
where the metric $g_\gs$ is defined as in Lemma
\ref{TeichmuellerAlmComp}. The complex analytic description is
better suited for explicit computations if $\gs$ is fixed.
However, if $\gs$ varies the underlying manifold varies as well.
This is sometimes inconvenient in the study of families. In the
following, we will use both descriptions. To avoid confusion we
will reserve $T^2$ for the standard 2-torus and use the notation
$\C/\gL(\gs)$ whenever we prefer to think in complex analytic
terms.
\end{remark}

Working on $\C/\gL(\gs)$ has the advantage that the metric is up
to a constant factor induced by the standard metric. In
particular, we will consider the 1-forms $dz = dx+idy$ and $d\Bar
z = dx - idy$. Note that
\[
\LScalar{dz}{d\Bar z}=0,\quad \|dz \|_{L^2} = \|d\Bar z\|_L^2=
\sqrt{2\gs_2},
\]
and, with the chirality operator $\tau$,
\[
\tau dz  = dz,\quad \tau d\Bar z = -d \Bar z, \quad\text{and}\quad
\tau(dz\wedge d\Bar z) = -2\gs_2.
\]
Using the isometry $\psi_\gs$ of \eqref{EllitpicCurveIsom}, one
translates this easily to $(T^2,g_\gs)$. More precisely, we define
\begin{equation}\label{Harmonic1FormsDef}
\go_\gs:= \psi_\gs^*(dz) = dx+\gs dy,\quad \text{and}\quad
\go_{\Bar \gs} := \psi_\gs^*(d\Bar z) = dx + \Bar \gs dy.
\end{equation}
Then we have a natural basis $(\go_\gs,\go_{\Bar \gs})$ for the
$C^\infty(T^2)$-module $\gO^1(T^2)$ satisfying
\begin{equation}\label{Harmonic1FormsProp}
\begin{split}
\Lscalar{\go_\gs}{\go_{\Bar\gs}}=0,\quad &\|\go_\gs\|_{L^2} =
\|\go_{\Bar\gs}\|_{L^2} = \sqrt{2\gs_2},\\ \tau \go_\gs =
\go_\gs,\quad&\text{and} \quad \tau \go_{\Bar\gs} = -
\go_{\Bar\gs}.
\end{split}
\end{equation}

\noindent\textbf{Flat Connections and Dolbeault Operators.} The
moduli space of flat $\U(1)$-connections over $T^2$ has a simple
structure. Since $T^2$ is the quotient of $\R^2$ by the standard
lattice $\Z^2$, we have
\[
\pi_1(T^2)\cong \Z e_1\oplus \Z e_2\subset \R^2,
\]
where $(e_1,e_2)$ is the standard basis of $\R^2$. Then, since
$\U(1)$ is abelian,
\[
\cM\big(T^2,\U(1)\big) \cong \Hom\big(\pi_1(T^2),\U(1)\big) \cong
\U(1)\times \U(1).
\]
As we are usually working with connections rather than
representations of the fundamental group, we summarize how the
above isomorphism works explicitly.

\begin{lemma}\label{TorusFlatConn}\quad \nopagebreak
\begin{enumerate}
\item Up to gauge equivalence, a flat $\U(1)$-connection over $T^2$
is induced by a $\Z^2$-invariant 1-form with constant
coefficients,
\begin{equation}
a_\nu = - 2\pi i (\nu_1 dx + \nu_2 dy) \in \gO^1(\R^2,i\R), \quad
\nu=( \nu_1, \nu_2)  \in \R^2.
\end{equation}
Two connections $a_\nu$ and $a_{\nu'}$ are gauge equivalent if
and only if $\nu-\nu'\in \Z^2$.
\item In terms of the generators
$e_1,e_2$ of $\pi_1(T^2)$, the holonomy representation of $a_\nu$
is given by
\[
\hol_{a_\nu}:\pi_1(T^2)\to \U(1),\quad \hol_{a_\nu}(e_j)=e^{2\pi i
\nu_j}.
\]
\end{enumerate}
\end{lemma}

\begin{proof}
A flat connection $a$ over $T^2$ is a $\Z^2$-invariant 1-form,
satisfying $da =0$. Therefore, it gives an element in de Rham
cohomology. Since
\[
H^1(T^2,\R) = \R^2,
\]
we can find a $\Z^2$-invariant function $f:\R^2\to \R$ and $\nu=(
\nu_1, \nu_2)  \in \R^2$ such that
\[
a - i df= - 2\pi i (\nu_1 dx + \nu_2 dy).
\]
Defining $u:= \exp(if)$, we get a gauge transformation on $T^2$
which brings $a$ into the claimed form. If $a_\nu$ and $a_{\nu'}$
are gauge equivalent, there exists a $\Z^2$-invariant function
$u:\R^2\to \U(1)$ satisfying
\[
2\pi i \big((\nu_1-\nu_1')dx + (\nu_2-\nu_2')dy)\big) = u^{-1}du.
\]
This easily implies that $u$ is of the form
\[
u = C\cdot \exp\big(2\pi i \Scalar{\nu-\nu'}{\left(\begin{smallmatrix} x\\
y
\end{smallmatrix}\right)}\big),\quad C\in \U(1),
\]
and this is $\Z^2$-invariant precisely if $\nu-\nu'\in \Z^2$.
This proves part (i). Concerning (ii), we compute using
Definition \ref{HolDef}
\[
\hol_{a_\nu}(e_j) = \exp\Big(- \int_0^1 a_\nu\big|_{s e_i}(e_i)
ds\Big) = \exp\big(2\pi i \nu_i\big).\qedhere
\]
\end{proof}

As pointed out in Remark \ref{EllipticCurve} it is often
convenient to work on $\C/\gL(\gs)$ rather than $T^2$. We collect
the following formul{\ae}, for definitions see Appendix
\ref{HolomAspects}, in particular \eqref{DolbDef}.

\begin{prop}\label{DolbeaultTorus}
Let $a_\nu$ be a flat $\U(1)$-connection over $T^2$ as in Lemma
\ref{TorusFlatConn}, and let $\gs=\gs_1+i\gs_2\in \H$.
\begin{enumerate}
\item The pullback of $a_\nu$ to $\C/\gL(\gs)$ is given by
\[
a= (\psi_\gs^{-1})^* a_\nu = - \Bar w_\nu dz + w_\nu d\Bar
z,\quad\text{where} \quad w_\nu:= \lfrac{\pi}{\gs_2}(\nu_2 -
\gs\nu_1).
\]
\item Let $\pd_a$ and $\Bar \pd_a$ be the Dolbeault operators
associated to $a$. Then the twisted Laplace operator on
$C^\infty\big(\C/\gL(\gs)\big)$ is given by
\[
\gD_a = 2\pd_a^t \pd_a = 2\Bar\pd_a^t\Bar\pd_a = -4\gs_2\Big(
\frac{\pd^2}{\pd z\pd \Bar z}-\Bar w_\nu\frac{\pd}{\pd\Bar z} +
w_\nu \frac{\pd}{\pd z} - |w_\nu|^2\Big).
\]
\item If $\gf\in C^\infty(\C)$ is $\gL(\gs)$-invariant, then
\[
\begin{split}
\pd_a\pd_a^t (\gf dz) &= \Bar \pd_a^t \Bar\pd_a(\gf dz) = (\lfrac
12 \gD_a\gf )dz, \\
\Bar \pd_a\Bar\pd_a^t (\gf d\Bar z) &= \pd_a^t \pd_a(\gf d \Bar z)
= (\lfrac 12 \gD_a\gf )d\Bar z,
\end{split}
\]
and
\[
\begin{split}
\Bar \pd_a\pd_a^t (\gf dz) & = - \pd_a^t \Bar \pd_a (\gf dz) =
-2\gs_2 \Big(\frac{\pd^2 \gf}{\pd\Bar z^2} + 2w_\nu
\frac{\pd\gf}{\pd\Bar z} + w_\nu^2 \gf \Big)d\Bar z, \\
\pd_a\Bar \pd_a^t (\gf d\Bar z) & = - \Bar \pd_a^t \pd_a (\gf
d\Bar z) = -2\gs_2 \Big(\frac{\pd^2 \gf}{\pd z^2} - 2\Bar w_\nu
\frac{\pd\gf}{\pd z} + \Bar w_\nu^2 \gf \Big)dz.
\end{split}
\]
In particular,
\[
\gD_a (\gf dz ) = (\gD_a\gf)dz\quad\text{and}\quad \gD_a (\gf
d\Bar z ) = (\gD_a\gf)d\Bar z.
\]
\end{enumerate}
\end{prop}

\begin{proof}[Sketch of proof]
Although Proposition \ref{DolbeaultTorus} is a standard exercise
in complex analysis, we want to give some remarks on the proof to
clarify the sign conventions we are using. For part (i), we use
the 1-forms $\go_\gs$ and $\go_{\Bar \gs}$ of
\eqref{Harmonic1FormsDef} to express
\[
dy = \lfrac 1{2i \gs_2} (\go_\gs - \go_{\Bar \gs}),\quad dx =
\lfrac 1{2i \gs_2}(\gs \go_{\Bar\gs} - \Bar \gs \go_\gs).
\]
Since $\go_\gs= \psi_\gs^*(dz)$ and $\go_{\Bar \gs}=
\psi_\gs^*(d\Bar z)$ we deduce that
\[
(\psi_\gs^{-1})^* a_\nu = - \lfrac{\pi}{\gs_2}\big( \nu_1(\gs
d\Bar z - \Bar\gs dz) + \nu_2(dz - d\Bar z)\big) = - \Bar w_\nu
dz + w_\nu d\Bar z,
\]
with $w_\nu= \frac{\pi}{\gs_2}(\nu_2 - \gs\nu_1)$ as claimed. For
part (ii) and (iii) we note that
\[
\pd_a = \pd -\emu(\Bar w_\nu dz),\quad \Bar \pd_a = \Bar\pd
+\emu(w_\nu d\Bar z),
\]
and, since we are using the complex linear chirality operator,
\[
\pd_a^t = -\tau \circ \Bar \pd_a\circ \tau,\quad \Bar \pd_a^t =
-\tau \circ \pd_a\circ \tau.
\]
Then one computes that for every $\gf\in C^\infty(\C)$,
\[
\pd_a^t \pd_a \gf = -\tau \Bar \pd_a \tau \Big( \frac{\pd\gf}{\pd
z} - \Bar w_\nu\gf\Big)dz = -\tau \Big(\frac{\pd^2\gf}{\pd z\pd
\Bar z}-\Bar w_\nu\frac{\pd\gf}{\pd\Bar z} + w_\nu
\frac{\pd\gf}{\pd z} - |w_\nu|^2\gf \Big)d\Bar z\wedge dz,
\]
where we have used that $\tau dz = dz$. Since $\tau(d\Bar z\wedge
dz) = 2\gs_2$, this yields the claimed formula for $\pd_a^t \pd_a
\gf$. In a similar way one computes $\Bar\pd_a^t \Bar \pd_a \gf$.
Then part (ii) follows since
\[
\gD_a \gf = (\pd_a+\Bar\pd_a)^t (\pd_a+\Bar\pd_a) \gf = \pd_a^t
\pd_a \gf + \Bar\pd_a^t \Bar \pd_a \gf.
\]
Using the same ideas one easily verifies part (iii).
\end{proof}

Using Proposition \ref{DolbeaultTorus} it is straightforward to
determine the spectrum of $\gD_a$.

\begin{prop}\label{TorusSpec}
Let $\gs\in \H$. For $n=(n_1,n_2)\in \Z^2$ define
\[
\gf_n:= e^{\Bar w_n z - w_n\Bar z}:\C\to \U(1),\quad w_n =
\lfrac\pi{\gs_2} (n_2- \gs n_1).
\]
Then $\gf_n$ is $\gL(\gs)$-invariant, and $\setdef{\gf_n}{n\in
\Z^2}$ is a orthonormal basis for $L^2\big(\C/\gL(\gs)\big)$. If
$a = - \Bar w_\nu dz + w_\nu d\Bar z$ is a flat
$\U(1)$-connection as before, then
\[
\gD_a \gf_n = \gl_{n,\nu} \gf_n,\quad\text{where}\quad
\gl_{n,\nu}= 4\gs_2 |w_{n-\nu}|^2 \gf_n.
\]
Moreover,
\[
\Bigsetdef{\frac{\gf_ndz}{\sqrt{2\gs_2}}, \frac {\gf_nd \Bar
z}{\sqrt{2\gs_2}}}{n\in \Z^2}
\]
gives an orthonormal basis for the space of 1-forms consisting of
eigenforms for $\gD_a$ with respect to the same eigenvalues
$\gl_{n,\nu}$.
\end{prop}

Since Proposition \ref{TorusSpec} is well known, we shall proceed
without further comments on the proof. However, we want to note
that we can use Proposition \ref{TorusSpec} and the Hodge-de-Rham
isomorphism to determine the twisted cohomology groups of $T^2$.
Recall that in the proof of Proposition \ref{TruncSignOpKern} we
have already done this using topological methods.

\begin{cor}\label{HarmonicFormsTorus}
Let $a_\nu$ be a flat $\U(1)$-connection over $T^2$, and let
$\gs\in \H$ determine the metric $g_\gs$. Then the cohomology of
$T^2$ with values in the line bundle $L_{a_\nu}$ is given in terms
of harmonic forms by
\[
\sH^\bullet(T^2,L_{a_\nu}) = \begin{cases} \C\oplus (\C \go_\gs
\oplus \C
\go_{\Bar\gs}) \oplus \C dx\wedge dy, &\text{\rm if $\nu\in \Z^2$},\\
\hphantom{ \C\oplus (\C \go_\gs \oplus \C }\{0\},&\text{\rm if
$\nu\notin \Z^2$},
\end{cases}
\]
where $\go_\gs$ and $\go_{\Bar\gs}$ are as in
\eqref{Harmonic1FormsDef}.
\end{cor}

\noindent\textbf{Mapping Tori with Fiber $\boldsymbol{T^2}$.} It
is well known that the mapping class group of $T^2$ is isomorphic
$\SL_2(\Z)$, see \cite[Sec. 2.9]{Iva}. Here, the action of an
element $M\in \SL_2(\Z)$ on $T^2$ is the one induced by matrix
multiplication
\[
\R^2\to \R^2,\quad  \big(\begin{smallmatrix} x\\ y
\end{smallmatrix}\big)\mapsto M \big(\begin{smallmatrix} x\\ y
\end{smallmatrix}\big).
\]
As every $M\in \SL_2(\Z)$ has determinant 1, it preserves the
volume form $dx\wedge dy$ and we do not have to invoke Proposition
\ref{MoserTrick}. On the other hand, $\SL_2(\Z)$ acts on the
Riemann sphere $\widehat \C$ by fractional linear
transformations, and this restricts to an action on $\H$,
\begin{equation}\label{FractLin}
M:\H\to \H,\quad M\gs := \frac {a\gs + b}{c\gs+d},\quad
M=\begin{pmatrix} a &b\\c&d
\end{pmatrix}\in \SL_2(\Z),
\end{equation}
see for example \cite[p. 6]{Shi}. Unfortunately, the isometry
\eqref{EllitpicCurveIsom} does not behave equivariantly with
respect to these two $\SL_2(\Z)$-actions. We can remedy this,
using the following involution on $\SL_2(\Z)$,
\[
\SL_2(\Z)\to \SL_2(\Z),\quad M=\begin{pmatrix} a &b\\ c&d
\end{pmatrix}\longmapsto M^{\op}:=\begin{pmatrix} d &b\\ c&a
\end{pmatrix}.
\]
Note that $(M_1M_2)^{\op} = M_2^{\op}M_1^{\op}$, so that we can
use the involution to turn a left action of $\SL_2(\Z)$ into a
right action.

\begin{lemma}\label{TeichmAlmComp:Equiv}
Let $\gF: \H\to \sJ_\go$ be the map of Lemma
\ref{TeichmuellerAlmComp}. Then for $M\in \SL_2(\Z)$
\begin{equation*}
\gF(M\gs) = (M^{\op})^{-1}\gF(\gs) M^{\op}.
\end{equation*}
\end{lemma}

\begin{proof}[Sketch of proof]
The group $\SL_2(\Z)$ is generated by the elements
\[
S:=\begin{pmatrix} 0 &-1\\ 1&0\end{pmatrix},\quad
T:=\begin{pmatrix} 1 & 1\\ 0&1\end{pmatrix},\quad
S^2=(ST)^3,\quad S^4=\Id,
\]
see \cite[pp. 16--17]{Shi}. As fractional linear transformations
they act as
\[
S(\gs)= -\frac 1{\gs},\quad T(\gs) = \gs+1.
\]
One then computes that for $\gs=\gs_1+i\gs_2\in \H$,
\[
\begin{split}
\gF(S\gs) = \gF\Big(-\frac {\bar \gs}{|\gs|^2}\Big) &= \frac
1{\gs_2}
\begin{pmatrix} \gs_1 &-1\\ |\gs|^2 &-\gs_1\end{pmatrix}\\ &=
\frac 1{\gs_2}\begin{pmatrix} 0 &1\\ -1&0\end{pmatrix}
\begin{pmatrix}
-\gs_1 &-|\gs|^2\\ 1& \gs_1\end{pmatrix}\begin{pmatrix} 0 &-1\\
1&0\end{pmatrix} = S^{-1} \gF(\gs) S,
\end{split}
\]
and
\[
\gF(T\gs)=\gF(\gs+1) = \frac 1{\gs_2} \begin{pmatrix} -\gs_1-1 &-|\gs_1+1|^2\\
1 &\gs_1+1\end{pmatrix}=\ldots =T^{-1} \gF(\gs) T.
\]
Now one has to verify that the formula of Lemma
\ref{TeichmAlmComp:Equiv} holds for all words in $S$ and $T$.
This is almost tautologically true. For example, if we consider
$M=ST$, then
\[
\gF(M\gs) = S^{-1} \gF(T\gs) S = (TS)^{-1} \gF(\gs) TS =
(M^{\op})^{-1} \gF(\gs) M^{\op}.\qedhere
\]
\end{proof}

The above lemma suggests that using the involution $M\mapsto
M^{\op}$ we should either redefine the action of $\SL_2(\Z)$ as
the mapping class group or turn the natural left action
\eqref{FractLin} of $\SL_2(\Z)$ into a right action. We opt for
the latter, although this leads to an unfortunate difference in
notation compared to the literature. However, redefining the
action of $\SL_2(\Z)$ as the mapping class group seems more
unnatural.

\begin{dfn}\label{TorusFiberMet}
Let $M\in SL_2(\Z)$, and let $\gs(t):\R\to \H$ be $M$-invariant in
the sense that $M^{\op}\gs(t) = \gs(t+1)$. Then we define
$(T^2_M,g_\gs)$ to be the mapping torus
\[
\big(T^2 \times \R\big)/\sim, \quad \big(M\big(\begin{smallmatrix}
x\\ y
\end{smallmatrix}\big),t\big) \sim \big( \big(\begin{smallmatrix} x\\ y
\end{smallmatrix}\big),t+1\big),\quad \big(\big(\begin{smallmatrix} x\\ y
\end{smallmatrix}\big),t\big)\in T^2 \times
\R,
\]
endowed with the metric induced by
\[
g_\gs:= dt\otimes dt + g_{\gs(t)}.
\]
Here, for each $t\in \R$ the metric $g_{\gs(t)}$ on $T^2$ is
defined as in Lemma \ref{TeichmuellerAlmComp}.
\end{dfn}

\noindent\textbf{Flat $\boldsymbol{\U(1)}$-connections over
$\boldsymbol{T^2_M}$.} We now give an explicit description of
flat $\U(1)$-connections over $T^2_M$ up to gauge equivalence.

\begin{prop}\label{TorusBundleFlatConn}
Let $M\in SL_2(\Z)$. Every flat $\U(1)$-connection $A$ over the
mapping torus  $T^2_M$ is equivalent to one induced by a flat
connection $a_\nu$ over $T^2$ as in Lemma \ref{TorusFlatConn} and
a gauge transformation $u\in C^\infty\big(T^2,\U(1)\big)$
satisfying
\begin{equation}\label{TorusBundleFlatConn:1}
\begin{pmatrix} m_1\\ m_2\end{pmatrix}:= (\Id - M^t)\begin{pmatrix}
\nu_1\\ \nu_2\end{pmatrix} \in \Z^2,
\end{equation}
and
\begin{equation*}\label{TorusBundleFlatConn:2}
u=\exp \big[- 2\pi i\big(m_1 x + m_2 y +\gl \big)\big]
\end{equation*}
for some $\gl\in [0,1)$.
\end{prop}

\begin{proof}
Let $a_\nu= -2\pi i (\nu_1dx +\nu_2dy)$ be a connection over $T^2$
as in Lemma \ref{TorusFlatConn}. Since $M$ acts by matrix
multiplication on $\R^2$, the pullback of $a_\nu$ by $M$ is given
by
\[
M^*a_\nu = -2\pi i(\mu_1 dx + \mu_2 dy),\quad \begin{pmatrix} \mu_1\\
\mu_2
\end{pmatrix} = M^t \begin{pmatrix} \nu_1\\ \nu_2
\end{pmatrix}.
\]
Now the condition $\widehat M_u^* a_\nu= a_\nu$ of
\eqref{MapTorConnTransSurf} means that connection $a_\nu$ is the
restriction of a connection over $T^2_M$ if and only if there
exists a $\Z^2$-invariant function $u:\R^2\to \U(1)$ such that
\begin{equation}\label{TorusBundleFlatConn:0}
\begin{pmatrix} m_1\\ m_2\end{pmatrix}= \big(\Id - M^t\big)
\begin{pmatrix} \nu _1\\ \nu_2
\end{pmatrix} = \frac {i}{2\pi}
\begin{pmatrix} u^{-1}\pd_x u\\ u^{-1}\pd_y u
\end{pmatrix}.
\end{equation}
Clearly, a function $u:\R^2\to \U(1)$ satisfying
\eqref{TorusBundleFlatConn:0} is necessarily of the form
\[
u = \exp \big[-2\pi i\big(m_1 x + m_2 y +\gl \big)\big], \quad
\gl\in [0,1),
\]
and this is $\Z^2$-invariant precisely if $(m_1,m_2)\in \Z^2$.
\end{proof}

\begin{remark}\label{TorusBundleFlatConnRem}\quad \nopagebreak
\begin{enumerate}
\item Note that the gauge transformation $u$ is---up to the number
$\gl\in [0,1)$---determined by \eqref{TorusBundleFlatConn:1}. For
simplicity we will sometimes consider only the case $\gl=0$ and
neglect the gauge transformation from the notation.
\item In Remark \ref{MapTorTrivFlat} (ii) we have pointed out that
the bundle $L\to T^2_M$ on which the flat connection $A$ is
defined is not necessarily trivializable. In the case at hand this
topological data is encoded in \eqref{TorusBundleFlatConn:1}. One
verifies---using for example Proposition
\ref{MapTorIsomCond}---that $L$ is trivializable if and only if
\[
\begin{pmatrix} m_1\\ m_2\end{pmatrix} \in
\im\big(\Id-M^t:\Z^2\to \Z^2\big).
\]
For example, if $M=\left(\begin{smallmatrix}  1 &2\\ 0 &1
\end{smallmatrix}\right)$, then $\Id-M^t = \left(
\begin{smallmatrix}  0 &0 \\ 2 &0 \end{smallmatrix}\right)$ and so
\[
\begin{pmatrix} 0\\ 1\end{pmatrix} \notin
\im\big(\Id-M^t :\Z^2\to \Z^2\big).
\]
This implies that in this case, the gauge transformation $u =
\exp(-2\pi i y)$ defines a flat bundle $L\to T^2_M$ which is not
trivializable.
\end{enumerate}
\end{remark}

\subsection{The Bismut-Cheeger Eta Form}

We now want to use Proposition \ref{MapTorBCForm} to express the
Eta form in terms of the data introduced in the last paragraphs.
First of all, we need to understand the vertical chirality
operator and its variation.

Let $\gs(t)=\gs_1(t)+i\gs_2(t)$ be an $M$-invariant path in $\H$
in the sense of Definition \ref{TorusFiberMet}, and let
$\go_{\gs(t)}$ and $\go_{\Bar \gs(t)}$ be the associated path of
1-forms as defined in \eqref{Harmonic1FormsDef}. Differentiating
with respect to $t$ easily yields that
\[
\dot \go_{\gs(t)} = \lfrac{\dot\gs(t)}{2i\gs_2(t)}(\go_{\gs(t)} -
\go_{\Bar\gs(t)})\quad\text{and}\quad \dot \go_{\Bar \gs(t)} =
\lfrac{\dot{\Bar\gs}(t)}{2i\gs_2(t)}(\go_{\gs(t)} -
\go_{\Bar\gs(t)}).
\]
Let $\tau_t$ be the associated path of chirality operators on
$\gO^\bullet(T^2)$. Since the volume form on $T^2$ does not vary
with $t$ the action of $\tau_t$ on $\gO^0(T^2)$ and $\gO^2(T^2)$
is independent of $t$. On $\gO^1(T^2)$ it is determined by $\tau_t
\go_{\gs(t)} = \go_{\gs(t)}$ and $\tau_t \go_{\Bar \gs(t)} = -
\go_{\Bar \gs(t)}$, see \eqref{Harmonic1FormsProp}. This readily
implies that the derivative of $\tau_t$ with respect to $t$ is
given by
\begin{equation}\label{DotTauTorus}
\dot\tau_t|_{\gO^0(T^2)}=0,\quad \dot \tau_t|_{\gO^2(T^2)}=0,\quad
\dot \tau_t \go_{\gs(t)} = \lfrac{i\dot \gs(t)}{\gs_2(t)}
\go_{\Bar\gs(t)},\quad \dot \tau_t \go_{\Bar \gs(t)} =
\lfrac{i\dot{\Bar \gs}(t)}{\gs_2(t)} \go_{\gs(t)}.
\end{equation}

\begin{prop}\label{TorBundleBCForm}
Let $a_\nu$ be a flat connection over $T^2$ as in Proposition
\ref{TorusBundleFlatConn}. Denote by $A$ be the associated flat
connection on the line bundle $L\to T^2_M$, and let $\gs(t)$ be an
$M$-invariant path in $\H$. Then the Bismut-Cheeger Eta form is
given by
\[
\widehat \eta_A = \frac 1{2\pi} \Re \Big(\dot \gs(t)
\int_0^\infty F_\nu\big(\gs(t),u\big) du\Big)dt,
\]
where for every $\gs=\gs_1+i\gs_2 \in \H$ and $u\in \R^+$
\begin{equation}\label{F:Def}
F_\nu(\gs,u) = \sum_{n\in \Z^2}\frac {\pi^2}{\gs_2^2}
\big(n_2-\nu_2 -\Bar\gs(n_1-\nu_1)\big)^2 e^{-u
\frac{\pi^2}{\gs_2}|n_2-\nu_2 -\gs(n_1-\nu_1)|^2 }.
\end{equation}
The sum in \eqref{F:Def} converges absolutely and there are
estimates, locally uniform in $\gs$,
\[
|F_\nu(\gs,u)| \le Ce^{-cu}\quad\text{as $u\to
\infty$},\qquad|F_\nu(\gs,u)| \le Ce^{-\frac cu}\quad\text{as
$u\to 0$}.
\]
\end{prop}

\begin{proof}
According to Proposition \ref{MapTorBCForm}, the Eta form
associated to the connection $A$ and the path $\gs(t)$ is given by
\[
\widehat \eta_A = \frac i{16\pi} dt\wedge \int_0^\infty
\Tr_v\Big(\dot\tau_v
\big(d_{A,v}^td_{A,v}-d_{A,v}d_{A,v}^t\big)e^{-\frac u4
D_{A,v}^2}\Big) du.
\]
It follows from Lemma \ref{VerticalConnMapTor} that under the
identification
\[
\gO^\bullet_v(T^2_M,L) = \bigsetdef{\ga_t:\R\to
\gO^\bullet(T^2)}{\ga_{t+1} = \widehat M^*\ga_t},
\]
the operator $d_{A,v}$ coincides with $d_{a_\nu}$ applied
pointwise for each $t$. The same is true for $d_{A,v}^t$ and
$d_{a_\nu}^t$, where the transpose has to be taken pointwise for
each $t$ with respect to the metric induced by $\gs(t)$. Now, the
operators $d_{a_\nu}d_{a_\nu}$, $d_{a_\nu}d_{a_\nu}^t$ and
$e^{-\frac u4 D_{a_\nu}^2}$ all preserve the decomposition
\[
\gO^\bullet(T^2) = \gO^0(T^2)\oplus \gO^1(T^2)\oplus \gO^2(T^2).
\]
Moreover, we know from \eqref{DotTauTorus} that the operator
$\dot \tau_v$ acts trivially on $\gO^0\oplus \gO^2$. Therefore,
\[
\Tr_v\Big(\dot\tau_v
\big(d_{A,v}^td_{A,v}-d_{A,v}d_{A,v}^t\big)e^{-\frac u4
D_{A,v}^2}\Big) = \Tr_t\Big(\dot\tau_t
\big(d_{a_\nu}^td_{a_\nu}-d_ad_{a_\nu}^t\big)e^{-\frac u4
D_{a_\nu}^2}\big|_{\gO^1(T^2)}\Big),
\]
where the subscripts $t$ indicate that we consider the right hand
side as a function of $t$. Instead of working over $(T^2,g_\gs)$
we now switch to $\C/\gL(\gs)$ to be able to use Proposition
\ref{DolbeaultTorus} and Proposition \ref{TorusSpec}. Using the
notation introduced used there, we write
\[
a= (\psi_\gs^{-1})^* a_\nu,\quad w_n = \lfrac\pi{\gs_2} (n_2- \gs
n_1),\quad \gf_n= e^{\Bar w_n z - w_n\Bar z}:\C\to \U(1),
\]
where $n=(n_1,n_2)\in \Z^2$. Note that for notational
convenience, we have now dropped the reference to the $t$
dependence. Then Proposition \ref{DolbeaultTorus} yields that
\[
\begin{split}
\big(d_a^td_a-d_ad_a^t\big) (\gf_n dz) &= 2\Bar\pd_a \pd_a^t
(\gf_n dz) = -4\gs_2 ( w_n^2 -2w_\nu w_n +w_\nu^2 )\gf_n d\Bar z\\
&= -4\gs_2 ( w_{n-\nu})^2\gf_n d\Bar z.
\end{split}
\]
Under the isometry $\psi_\gs:(T^2,\gs)\to \C/\gL(\gs)$, the pair
$(dz,d\Bar z)$ pulls back to $(\go_\gs,\go_{\Bar\gs})$. Hence, we
deduce from \eqref{DotTauTorus} that
\[
\dot\tau \big(d_a^td_a-d_ad_a^t\big) (\gf_n dz) =
-4i\dot{\bar\gs} (w_{n-\nu})^2 \gf_n dz.
\]
Similarly, one finds that
\[
\dot\tau \big(d_a^td_a-d_ad_a^t\big) (\gf_n d\Bar z) = -4i\dot\gs
(\Bar w_{n-\nu})^2 \gf_n d\Bar z.
\]
The Laplace operators $\gD_a$ and $D_{a_\nu}^2$ coincide via
$\psi_\gs:(T^2,\gs)\to \C/\gL(\gs)$. Using Proposition
\ref{TorusSpec} and the fact that $e^{-\frac u4 \gD_a^2}$ is an
operator with smooth kernel, one then concludes that
\begin{equation*}
\begin{split}
\frac i{16\pi}\Tr\Big(\dot\tau \big(d_a^td_a
-d_ad_a^t\big)e^{-\frac u4 \gD_a^2}\Big)&= \frac 1{2\pi}
\sum_{n\in \Z^2} \Re\big( \dot \gs  (\Bar w_{n-\nu})^2\big) e^{ -
u\gs_2 |w_{n-\nu}|^2}\\ &= \frac 1{2\pi} \Re\big(\dot\gs
F_\nu(\gs,u)\big),
\end{split}
\end{equation*}
where in the last step we have simply used the definition of
$F_\nu(\gs,u)$. Concerning the absolute convergence of
$F_\nu(\gs,u)$ and the estimate as $u\to \infty$ we assume for
simplicity that $\nu=0$. The general case requires only minor
changes. Define
\[
\begin{split}
r_\gs:=\min\bigsetdef{|x_2 - \gs x_1|}{(x_1,x_2)\in \R^2,\
|x_1|+|x_2|=1},\\
R_\gs:=\max\bigsetdef{|x_2 - \gs x_1|}{(x_1,x_2)\in \R^2,\
|x_1|+|x_2|=1}.
\end{split}
\]
Clearly, $r_\gs$ and $R_\gs$ depend continuously on $\gs$. For
$n\in \Z^2$ and some constants $c$ and $C$, not depending on $n$
and $\gs$, we have
\[
|w_n|^2 \ge c r_\gs|n|^2 ,\quad  |w_n|^2 \le C R_\gs |n|^2,
\]
so that
\[
\big|\Bar w_n^2 e^{ - u\gs_2 |w_n|^2}\big| \le C  |n|^2 e^{-u c
 |n|^2 },
\]
where the constants $c$ and $C$ now depend continuously on $\gs$.
This implies absolute convergence of the series in \eqref{F:Def}.
Concerning the estimate as $u\to \infty$ one now proceeds exactly
as in the proof of Lemma \ref{BasicKernelEst} and we will skip
the details. However, the estimate for $u\to 0$ cannot be read off
in the same way, and the general theory only yields that
$F_\nu(\gs,u)$ is bounded as $u\to 0$, see Theorem
\ref{Fam:IndexTheorem}. To obtain the required estimate, we
proceed as in \cite[Thm. 2.15]{BC92}. Recall that the Poisson
summation formula states that for any rapidly decreasing function
$f$ on $\R^d$
\begin{equation}\label{Poisson}
\sum_{n\in \Z^d} f(n) = \sum_{n\in \Z^d} \widehat f(n),\quad
\widehat f(n) = \int_{\R^d} f(x) e^{-2\pi i \scalar{n}{x}} dx,
\end{equation}
see e.g., \cite[Sec. 20.1]{Lang}. We now let $\nu\in \R^2$ be
arbitrary again and use the Poisson summation formula to bring
$F_\nu(\gs,u)$ into a different form. Since we will need the
formula only to obtain the estimate as $u\to 0$, we give only some
intermediate steps
\[
\begin{split}
F_\nu(\gs,u) &=  \sum_{n\in \Z^2} (\Bar w_{n-\nu})^2 e^{ - u\gs_2
|w_{n-\nu}|^2} = \sum_{n\in \Z^2} \int_{\R^2} (\Bar w_{x-\nu})^2
e^{ - u\gs_2 |w_{x-\nu}|^2} e^{-2\pi i
\scalar{x}{n}}dx\\
&= \sum_{n\in \Z^2} e^{-2\pi i \scalar{\nu}{n}} \frac 1{\gs_2 u^2}
\int_{\R^2} (x_1 + ix_2)^2 e^{-\pi |x|^2}e^{-2\pi i
\scalar{x}{\xi_n}}dx\Big|_{\xi_n = \frac 1{\sqrt{\pi u \gs_2}}
\left(\begin{smallmatrix}\gs_2 n_2\\
n_1+\gs_1n_2\end{smallmatrix}\right)},
\end{split}
\]
where the last line follows from a suitable substitution. For
arbitrary $\xi=(\xi_1,\xi_2) \in \R^2$ one computes that
\[
\int_{\R^2} (x_1 + ix_2)^2 e^{-\pi |x|^2}e^{-2\pi i
\scalar{x}{\xi_n}}dx = \big(-i\xi_1 + \xi_2\big)^2 e^{-\pi
|\xi|^2}.
\]
Moreover, with $\xi= \frac 1{\sqrt{\pi u \gs_2}}
\left(\begin{smallmatrix}\gs_2 n_2\\
n_1+\gs_1n_2\end{smallmatrix}\right)$ as above,
\[
-i\xi_1 + \xi_2 = \frac 1{\sqrt{\pi u \gs_2}} (n_1 + \Bar\gs n_2),
\]
and thus,
\begin{equation}\label{F:Alt}
F_\nu(\gs,u) = \frac 1{\pi \gs_2^2} \sum_{n\in \Z^2} e^{-2\pi i
\scalar{\nu}{n}} (\Bar w_n^*)^2 u^{-3}
e^{-\frac{|w_n^*|^2}{u\gs_2}},\quad w_n^* := n_1 + \gs n_2.
\end{equation}
Since $\Bar w_n=0$ for $n=(0,0)$, we see that the only possible
term which might fail to decrease exponentially as $u\to 0$ drops
out. Hence, we can proceed again as in Lemma \ref{BasicKernelEst}
to deduce that $|F_\nu(\gs,u)| \le Ce^{-\frac cu}$ as $u\to 0$,
with constants $c$ and $C$ depending continuously on $\gs$.
\end{proof}

We now want to simplify the expression for $\widehat \eta_A$
further and get an expression for the integral of the Eta form
$\widehat \eta_A$ over the base. We note that the function $F_\nu$
in \eqref{F:Def} depends on $\nu\in \R^2$ only modulo $\Z^2$, see
also \eqref{F:Alt}. Hence, we will often assume in the following
that $0\le \nu_1<1$ or even that $\nu \in [0,1)^2$.

\begin{theorem}\label{F:Thm}
For $\nu=(\nu_1,\nu_2)\in \R^2$ with $0\le\nu_1<1$ and $\gs\in \H$
write $z= \nu_1\gs - \nu_2$.
\begin{enumerate}
\item Employing the notation
\[
q_\gs= e^{2\pi i \gs},\quad \text{and}\quad q_z= e^{2\pi i z},
\]
we define
\[
E_\nu(\gs) := \sum_{n_1>0}\sum_{n_2>0}\lfrac 1{n_2}
(q_z+q_z^{-1})^{n_2}q_\gs^{n_1n_2},\quad \text{if $\nu_1=0$},
\]
and
\[
E_\nu(\gs) := \sum_{n_2>0}\lfrac 1{n_2} q_z^{n_2}+
\sum_{n_1>0}\sum_{n_2>0}\lfrac 1{n_2}
(q_z+q_z^{-1})^{n_2}q_\gs^{n_1n_2},\quad \text{if $\nu_1\neq 0$}.
\]
Then the sum defining $E_\nu(\gs)$ converges absolutely to a
function which is holomorphic on $\H$.

\item Let $F_\nu(\gs,u)$ be as in \eqref{F:Def}. Then
\[
\frac 1{2\pi}  \int_0^\infty F_0\big(\gs,u\big) du = \frac 16
-\frac 1{2\pi \gs_2} + \frac i{\pi}\frac{\pd}{\pd \gs} E_0(\gs),
\]
and for $\nu\notin\Z$
\[
\frac 1{2\pi}  \int_0^\infty F_\nu\big(\gs,u\big) du = P_2(\nu_1)
+ \frac i{\pi} \frac{\pd}{\pd \gs} E_\nu(\gs),
\]
where $P_2$ is the second periodic Bernoulli function, see
Definition \ref{PeriodicBernoulli}.
\end{enumerate}
\end{theorem}

\begin{remark}\label{E:Rem}\quad\nopagebreak
\begin{enumerate}
\item The function $E_\nu(\gs)$ is related to the logarithm of a
generalized Dedekind Eta function, see \eqref{E:DedRel:1} and
Lemma \ref{E:GenDedRel} below. As such it appears in the constant
term of the Laurent series at $s=1$ of certain Eisenstein series.
Without going into details about the exact relation, we recall
that the determination of this constant term is classically
referred to as a \emph{Kronecker limit formula}, see \cite[Ch.
20]{Lang} and \cite[Sec. 4]{RS73}. In the proof of Theorem
\ref{F:Thm} below, we mimic a combined proof of the \emph{first}
and the \emph{second} Kronecker limit formula as in \cite{Lang}.
\item Since the sum defining $E_\nu(\gs)$ converges absolutely, we
can interchange the summation over $n_1$ and $n_2$. Since
$|q_\gs^{n_2}|<1$ for $n_2>0$ we find
\[
\sum_{n_1>0} q_\gs^{n_1n_2} = \frac{q_\gs^{n_2}}{1- q_\gs^{n_2}}
= \frac{q_\gs^{n_2/2}}{q_\gs^{-n_2/2}- q_\gs^{n_2/2}} = \frac
i2\big(\cot(\pi n_2 \gs ) + i\big).
\]
Hence, for $\nu_1=0$
\[
E_\nu(\gs) = i \sum_{n_2>0} \lfrac 1{n_2} \cos(2\pi z n_2)
\big(\cot(\pi n_2 \gs ) + i\big).
\]
In the case that $\nu_1\neq 0$ one can combine the two sums over
$n_2$. Then
\[
E_\nu(\gs) = i \sum_{n_2\ge 1} \lfrac 1{n_2}\big(\cos(2\pi n_2
z)\cot(\pi n_2\gs)  + \sin(2\pi n_2 z)\big).
\]
\end{enumerate}
\end{remark}

\begin{proof}[Proof of Theorem \ref{F:Thm}]
We can assume for simplicity that $0\le \nu_2<1$ as well, since
the claims of Theorem \ref{F:Thm} depend on $\nu_2$ only modulo
$\Z$. As an auxiliary tool we define
\begin{equation}\label{G:Def}
G_\nu(\gs,s):= \frac 1{2\pi\gG(s)}\int_0^\infty
u^{s-1}F_\nu(\gs,u)du.
\end{equation}
The estimates in Proposition \ref{TorBundleBCForm} ensure that
$G_\nu(\gs,s)$ is a holomorphic function for all $s\in \C$.
Clearly,
\[
G_\nu(\gs,1) = \frac 1{2\pi}\int_0^\infty F_\nu(\gs,u)du.
\]
Since the sum over $n=(n_1,n_2)\in \Z^2$ defining $F_\nu(\gs,u)$
converges absolutely, we can first extract possible terms with
$n_1=\nu_1=0$, and then sum up the remaining terms. More
precisely, define
\[
F_\nu^0(\gs,u):= \begin{cases}\frac {\pi^2}{\gs_2^2} \sum_{n_2\in
\Z} |n_2-\nu_2|^2 e^{-u \frac{\pi^2}{\gs_2}|n_2-\nu_2|^2 },
&\text{if $\nu_1=0$},\\ \hphantom{\frac {\pi^2}{\gs_2^2}
\sum_{n_2\in \Z} \big(n_2}0, &\text{if $\nu_1\neq 0$},
\end{cases}
\]
and
\[
F_\nu^1(\gs,u) := \frac {\pi^2}{\gs_2^2} \sum_{n_1\neq
\nu_1}\sum_{n_2\in \Z} \big(n_2 -\Bar\gs n_1 + \Bar z \big)^2
e^{-u \frac{\pi^2}{\gs_2}|n_2-\gs n_1 +z|^2 },
\]
where $z = \nu_1\gs - \nu_2$. Accordingly, we can split
$G_\nu(\gs,s)$ for $\Re(s)$ large enough as
\[
\begin{split}
G_\nu(\gs,s) &= \frac 1{2\pi\gG(s)}\int_0^\infty
u^{s-1}F_\nu^0(\gs,u)du + \frac 1{2\pi\gG(s)}\int_0^\infty
u^{s-1}F_\nu^1(\gs,u)du \\
&=:G_\nu^0(\gs,s) + G_\nu^1(\gs,s).
\end{split}
\]
Now, again for $\Re(s)$ large enough, we can interchange summation
and integration, so that for $\nu_1=0$ the substitution $u\mapsto
\frac{\pi^2}{\gs_2}|n_2-\nu_2|^2 u$ yields
\[
\begin{split}
G_\nu^0(\gs,s) &= \frac \pi {2\gs_2^2} \sum_{n_2\in \Z}
|n_2-\nu_2|^2 \frac 1{\gG(s)}\int_0^\infty u^{s-1}e^{-u
\frac{\pi^2}{\gs_2}|n_2-\nu_2|^2 }du\\
&= \frac {\gs_2^{s-2}}{2\pi^{2s-1}} \sum_{n_2\in \Z}
|n_2-\nu_2|^{2-2s} = \frac {\gs_2^{s-2}}{2\pi^{2s-1}} \widetilde
\gz_{\nu_2}(s-1),
\end{split}
\]
where $\widetilde \gz_{\nu_2}(s-1)$ is the periodic Zeta function
in Proposition \ref{ZEtaCalc}. This implies that $G_\nu^0(\gs,s)$
admits a meromorphic continuation to the whole $s$-plane, and
\begin{equation}\label{F:Thm:1}
G_\nu^0(\gs,1) = \begin{cases}\hphantom{-(2\pi} 0, &\text{if
$\nu_2\neq 0$},\\ -(2\pi \gs_2)^{-1}, &\text{if $\nu_1=\nu_2= 0$}.
\end{cases}
\end{equation}
To identify $G_\nu^1(\gs,s)$ we assume again that $\Re(s)$ is
large enough, so that we can freely interchange summation and
integration. Then, with the substitution $u\mapsto
\frac{\pi^2}{\gs_2} u$, we get
\[
G_\nu^1(\gs,s) = \frac {\gs_2^{s-2}}{2\pi^{2s-1}} \sum_{n_1\neq
\nu_1} \frac 1{\gG(s)}\int_0^\infty u^{s-1} \sum_{n_2\in\Z}
\big(n_2 -\Bar\gs n_1 + \Bar z \big)^2 e^{-u |n_2-\gs n_1 + z|^2
} du.
\]
We now apply the Poisson summation formula \eqref{Poisson} to the
inner sum and compute that
\[
\begin{split}
\sum_{n_2\in\Z}& \big(n_2 -  \Bar\gs n_1 + \Bar z \big)^2 e^{-u
|n_2-\gs n_1 +z|^2 }\\
&= \sum_{n_2\in\Z} \int_\R \big(x - \Bar\gs n_1 + \Bar z\big)^2
e^{-u |x-\gs n_1 +z|^2 } e^{-2\pi i x n_2} dx\\
&= \sum_{n_2\in\Z} e^{-2\pi i\Re(\gs n_1 - z)n_2} e^{-u(\Im(\gs
n_1 - z))^2}\int_\R \big(x + i\Im(n_1\gs -z)\big)^2
e^{-u x^2} e^{-2\pi i x n_2} dx,\\
\end{split}
\]
where we have separated the real and imaginary parts and then made
the substitution $x\mapsto x - \Re\big(\gs n_1 - z\big)$.
Clearly, the integral in the last expression decays exponentially
as $|n_2|\to \infty$. Hence, the sum converges absolutely, and we
can rearrange the order of summation again. Write
\begin{equation}\label{F:Thm:2}
G_\nu^1(\gs,s) = G_\nu^{10}(\gs,s) + G_\nu^{11}(\gs,s),
\end{equation}
where $G_\nu^{10}(\gs,s)$ is the contribution coming from
$n_2=0$, i.e.,
\[
G_\nu^{10}(\gs,s) = \frac {\gs_2^{s-2}}{2\pi^{2s-1}} \frac
1{\gG(s)} \sum_{n_1\neq \nu_1} \int_0^\infty u^{s-1}
e^{-u(\Im(\gs n_1 - z))^2} \int_\R \big(x + i\Im(n_1\gs
-z)\big)^2 e^{-u x^2} dx.
\]
Setting $a:= \Im(n_1\gs -z)$ we have
\[
\int_\R (x + ia)^2 e^{-u x^2} dx = \sqrt\pi \big(\lfrac 12
u^{-3/2} - a^2 u^{-1/2} \big).
\]
Therefore, standard manipulations involving the Gamma function
yield
\[
\begin{split}
\int_0^\infty u^{s-1} e^{-ua^2} \int_\R (x + ia)^2 e^{-u x^2} dx
& = \sqrt\pi |a|^{3-2s}\big(\lfrac12 \gG(s-\lfrac32) -
\gG(s-\lfrac12)\big) \\
&= \sqrt\pi |a|^{3-2s} \gG(s-\lfrac12)\frac{4-2s}{2s-3}.
\end{split}
\]
Recalling that $a= \Im(n_1\gs -z)$ we find that for $\Re(s)$
large enough,
\[
\begin{split}
G_\nu^{10}(\gs,s) &= \sqrt\pi \frac {\gs_2^{s-2}}{2\pi^{2s-1}}
\frac{\gG(s-\lfrac12)}{\gG(s)}\frac{4-2s}{2s-3} \sum_{n_1\neq
\nu_1} \big|\Im(n_1\gs -z)\big|^{3-2s}\\
&= \sqrt\pi \frac {\gs_2^{1-s}}{2\pi^{2s-1}}
\frac{\gG(s-\lfrac12)}{\gG(s)}\frac{4-2s}{2s-3} \widetilde
\gz_{\nu_1}(2s-3),
\end{split}
\]
where we have used that $\Im(n_1\gs -z) = \gs_2(n_1-\nu_1)$.
Hence we have found an expression for $G_\nu^{10}(\gs,s)$ which
can be extended to a meromorphic function on the whole $s$-plane.
It follows from Proposition \ref{ZEtaCalc} that $s=1$ is not a
pole, and that
\begin{equation}\label{F:Thm:3}
G_\nu^{10}(\gs,1) = P_2(\nu_1) = \nu_1^2 -\nu_1 +\lfrac 16.
\end{equation}
Now we have to consider the general term $G^{11}_\nu(\gs,s)$ in
\eqref{F:Thm:2}. Write $a_n= |\Im(n_1\gs -z)|$ and $b_n=
\pi|n_2|$. Note that $a_n,b_n >0$ if $n_1\neq \nu_1$ and $n_2\neq
0$. Then, for $\Re(s)$ large,
\begin{equation}\label{F:Thm:4}
\begin{split}
G^{11}_\nu(\gs,s) = \frac {\gs_2^{s-2}}{2\pi^{2s-1}} \frac
1{\gG(s)} \sum_{n_1\neq \nu_1}&\sum_{n_2\neq 0} e^{-2\pi i\Re(\gs
n_1 - z)n_2} \int_0^\infty u^{s-1}e^{-ua_n^2}  \\ &\int_\R \Big(x
+
 \frac{i a_n}{\scriptstyle\sgn(n_1-\nu_1)}\Big)^2 e^{-u x^2} e^{-2 i
x \sgn(n_2) b_n} dx du.
\end{split}
\end{equation}
To compute the integrals in the sum, we replace $a_n$ and $b_n$
by real parameters $a, b>0$. Then
\[
\begin{split}
\int_\R \Big(x + \frac{i a}{\scriptstyle\sgn(n_1-\nu_1)}\Big)^2 &
e^{-u x^2} e^{-2 i x \sgn(n_2) b} dx \\& = - \Big(\frac{
\pd_b}{2\scriptstyle\sgn(n_2)}  + \frac{
a}{\scriptstyle\sgn(n_1-\nu_1)}\Big)^2 \int_\R e^{-u x^2} e^{-2 i
x \sgn(n_2) b} dx\\
&= - \Big(\frac{ \pd_b}{2\scriptstyle\sgn(n_2)}  + \frac{
a}{\scriptstyle\sgn(n_1-\nu_1)}\Big)^2 \sqrt{\frac\pi
u}e^{-b^2/u}.
\end{split}
\]
Therefore, the $u$-integral in \eqref{F:Thm:4} in terms of the
parameters $a$ and $b$ is given by
\begin{equation}\label{F:Thm:5}
- \sqrt\pi \Big(\frac{ \pd_b}{2\scriptstyle\sgn(n_2)}  + \frac{
a}{\scriptstyle\sgn(n_1-\nu_1)}\Big)^2 K_{s-\frac12}(a,b),
\end{equation}
where $K_s(a,b)$ is the \emph{Bessel K-function} \cite[Sec.
20.3]{Lang}
\[
K_s(a,b) = \int_0^\infty u^{s-1} e^{-(a^2 u + b^2/u)}du.
\]
Moreover, for fixed $s$, one has
\[
\pd_b K_s(a,b) = -2b K_{s-1}(a,b),
\]
so that \eqref{F:Thm:5} is actually a sum of Bessel $K$-functions
for different $s$-parameters. We also collect from \cite[Sec.
20.3]{Lang} that $K_s(a,b)$ is holomorphic on the whole $s$-plane
and satisfies estimates, locally uniform in $s$, of the form
\[
|K_s(a,b)| \le C \Big(\frac b a\Big)^s e^{-2ab},\quad ab\to
\infty.
\]
This implies that the summand in \eqref{F:Thm:4} decays
exponentially as $|(n_1,n_2)|\to \infty$, locally uniform in $s$.
From this one deduces that $G^{11}_\nu(\gs,s)$ can be extended
holomorphically to the whole $s$-plane, and that we can simply put
$s=1$ to find the value we are interested in. Now, $K_{\frac
12}(a,b) = \frac {\sqrt\pi} a e^{-2ab}$, see \cite[p. 271]{Lang}.
Using this, one verifies without effort that the value of
\eqref{F:Thm:5} at $s=1$ is equal to
\[
-2\pi a \big(1- \sgn(n_2(n_1-\nu_1))\big) e^{-2ab}.
\]
Using this one finds that
\[
G^{11}_\nu(\gs,1) = - \sum_{n_1\neq \nu_1}\sum_{n_2\neq 0}
|n_1-\nu_1| \big(1- \sgn(n_2(n_1-\nu_1))\big) e^{-2\pi
|n_2||\Im(\gs n_1-z)|} e^{-2\pi i\Re(\gs n_1 - z)n_2}.
\]
Since $n_2\neq 0$ and $n_1\neq \nu_1$, the above sum converges
absolutely. Moreover, all the terms with $\sign(n_2) =
\sign(n_1-\nu_1)$ drop out. Using the notation $q_\gs= e^{2\pi i
\gs}$ and $q_z= e^{2\pi i z}$, one obtains
\begin{equation}\label{F:Thm:6}
G^{11}_\nu(\gs,1) =  -2 \sum_{n_2>0} \nu_1 q_z^{n_2}  - 2
\sum_{n_1>0}\sum_{n_2>0}\Big[ (n_1+\nu_1) (q_zq_\gs^{n_1})^{n_2}
+ (n_1-\nu_1)(q_z^{-1}q_\gs^{n_1})^{n_2}\Big].
\end{equation}
Note that for $\nu_1=0$ the first term is equal to zero. Now,
$q_z$ and $q_\gs$ are holomorphic as functions of $\gs$, and
\[
\nu_1 q_z^{n_2} = \frac1{2\pi i n_2}\frac{\pd}{\pd \gs}
q_z^{n_2},\quad (n_1\pm \nu_1) (q_zq_\gs^{n_1})^{n_2}= \frac1{2\pi
i n_2}\frac{\pd}{\pd \gs} (q_z^{\pm 1}q_\gs^{n_1})^{n_2},
\]
In the case that $\nu_1\neq 0$ we have $z\in \H$. This yields
that $\frac1{n_2} q_z^{n_2}$ decays exponentially as $n_2\to
\infty$. For arbitrary $\nu_1$ and $n_1>0$ the term $\frac1{n_2}
(q_z^{\pm 1}q_\gs^{n_1})^{n_2}$ decays exponentially in both,
$n_1$ and $n_2$. Moreover, this decay is certainly locally uniform
in $\gs$. This implies that in \eqref{F:Thm:6} we can interchange
summation and differentiation to find that
\begin{equation}\label{F:Thm:7}
G^{11}_\nu(\gs,1) =  \frac i{\pi} \frac{\pd}{\pd \gs} E_\nu(\gs),
\end{equation}
where $E_\nu(\gs)$ is defined as in part (i) of the theorem. Since
the sums converge absolutely and locally uniform in $\gs$, we
conclude that $E_\nu(\gs)$ defines a holomorphic function on
$\H$, which proves part (i). Moreover, we have split the
auxiliary function in \eqref{G:Def} for $\Re(s)$ large as
\[
G_\nu(\gs,s) = G^0_\nu(\gs,s) + G^{10}_\nu(\gs,s)+
G^{11}_\nu(\gs,s).
\]
As we have seen, the terms on the right hand side extend to
meromorphic function on the $s$-plane, and thus, the above
equality continues to hold for all $s$.\footnote{Here, we
encounter a similar situation as in Remark
\ref{RhoFormS1BundleRem} (i). The terms $G^0_\nu(\gs,s)$ and
$G^{10}_\nu(\gs,s)$ are not necessarily holomorphic, whereas
$G_\nu(\gs,s)$ and $G^{11}_\nu(\gs,s)$ are. This implies that the
poles of $G^0_\nu(\gs,s)$ and $G^{10}_\nu(\gs,s)$ have to cancel
each other out. With some effort one can check this directly. We
will not go into further details, since the value we are
interested in is $s=1$, which is no pole for any of the
summands.} Therefore, we can insert the values at $s=1$, which we
have computed in \eqref{F:Thm:1}, \eqref{F:Thm:3} and
\eqref{F:Thm:7}, and deduce that
\[
G_\nu(\gs,1) = P_2(\nu_1) + \frac i{\pi} \frac{\pd}{\pd \gs}
E_\nu(\gs) + \begin{cases}\hphantom{-(2\pi} 0, &\text{if $\nu_2\neq 0$},\\
-(2\pi \gs_2)^{-1}, &\text{if $\nu_1=\nu_2= 0$},\end{cases}
\]
which proves part (ii) of Theorem \ref{F:Thm}.
\end{proof}

As a consequence of Theorem \ref{F:Thm}, we obtain the expression
for the integral of Bismut-Cheeger Eta form we were aiming at. We
know from Proposition \ref{TorBundleBCForm} that
\[
\widehat \eta_A = \frac 1{2\pi} \Re \Big(\dot \gs(t)
\int_0^\infty F_\nu\big(\gs(t),u\big)du\Big)dt
\]
Hence, integrating the formula in Theorem \ref{F:Thm} (ii) with
respect to $t$ one easily arrives at

\begin{theorem}\label{TorBundleRhoForm}
Let $M\in \SL_2(\Z)$, let $\gs(t)$ be an $M$-invariant path in
$\H$, and use $\gs(t)$ to endow the mapping torus $T^2_M$ with a
metric.
\begin{enumerate}
\item The untwisted Eta form $\widehat \eta$ satisfies
\[
\int_0^1\widehat \eta = \frac 1\pi \Re\Big[\pi \gs P_2(0) + i
E_\nu(\gs) \Big]_{\gs(0)}^{\gs(1)} - \frac 1{2\pi}\int_0^1
\frac{\dot \gs_1(t)}{\gs_2(t)}dt,
\]
where we use the abbreviation $[f(\gs)]_{\gs(0)}^{\gs(1)} =
f(\gs(1))-f(\gs(0))$.
\item Let $\nu\in \R^2\setminus \Z^2$ with $0\le \nu_1<1$ satisfy
$(\Id - M^t)\nu\in \Z^2$, and let $A$ be the corresponding flat
$\U(1)$-connection over the mapping torus $T^2_M$. Then
\[
\int_0^1\widehat \eta_A = \frac 1\pi \Re\Big[\pi \gs P_2(\nu_1) +
i E_\nu(\gs) \Big]_{\gs(0)}^{\gs(1)}.
\]
In particular, the Rho form $\widehat \rho_A$ satisfies
\[
\int_0^1\widehat \rho_A = \frac 1\pi \Re\Big[\pi \gs
\big(P_2(\nu_1)-\lfrac 16\big) + i \big(E_\nu(\gs)-E_0(\gs)\big)
\Big]_{\gs(0)}^{\gs(1)} + \frac 1{2\pi}\int_0^1 \frac{\dot
\gs_1(t)}{\gs_2(t)}dt.
\]
\end{enumerate}
\end{theorem}

\begin{remark}\label{E:Ext}
The forms $\widehat \eta_A$ and $\widehat \rho_A$ depend only on
$\nu$ modulo $\Z^2$. Also, by its very definition, $P_2(\nu_1)$
depends on $\nu_1$ only modulo $\Z$. Therefore, it is
reasonable---and convenient---to extend the definition of
$E_\nu(\gs)$ to arbitrary $(\nu_1,\nu_2)\in \R^2$ by letting
\[
E_{(\nu_1,\nu_2)}(\gs) := E_{(\nu_1-[\nu_1],\nu_2)}(\gs),
\]
where $[\nu_1]$ is the largest integer less or equal than
$\nu_1$. Then Theorem \ref{TorBundleRhoForm} (ii) continues to
hold without the assumption on $\nu_1$.
\end{remark}

\section[Torus Bundles over $S^1$, Explicit Computations]{Torus
Bundles over $\boldsymbol{S^1}$, Explicit
Computations}\label{TorusBundlesExp}

In this section, we want to give a more explicit formula for the
Rho invariants of a mapping torus $T^2_M$ with $M\in \SL_2(\Z)$.
The result depends considerably on whether $M$ is elliptic,
parabolic or hyperbolic---see Definition \ref{TorusMapClasses}
below---and we have to treat all three cases separately. Explicit
formul{\ae} for the untwisted Eta invariant have been obtained in
\cite{Ati87} and \cite[App. 3]{Che87}. Both references make use
of adiabatic limits, and much of our treatment parallels their
discussion. In \cite{Che87} the focus is on the hyperbolic case,
and the Eta invariant is identified with the value of certain
number theoretical $L$-series, see also \cite{ADS, BC92, Mue84}.
This has its origin in Hirzebruch's work \cite{Hir73}, where a
topological interpretation of the aforementioned $L$-series was
conjectured. A similar relation can also be found for twisted Eta
invariants. However, our aim is to get a simple formula for the
Rho invariant, and values of $L$-series are certainly not easy to
compute. Fortunately, Atiyah \cite{Ati87} found a number of very
different ways to express the untwisted Eta invariant of $T^2_M$,
and we shall derive a formula for the Rho invariant along those
lines.\\

\noindent\textbf{Rough Classification of Elements in
$\boldsymbol{\SL_2(\Z)}$.} For explicit computations we now have
to find $M$-invariant paths in $\H$. For this we will use that
elements in $\SL_2(\Z)$ split into three natural classes.

\begin{dfn}\label{TorusMapClasses}
Let $M\in \SL_2(\Z)$, and let $\gD:=(\tr M)^2-4$. Then $M$ is
called
\begin{enumerate}
\item \emph{elliptic}, if $\gD<0$,
\item \emph{parabolic}, if $\gD=0$, and
\item \emph{hyperbolic}, if $\gD>0$.
\end{enumerate}
\end{dfn}

\begin{remark*}\label{TorBundleDiff}
Recall that according to Lemma \ref{MapTorIsotop}, the
diffeomorphism type of $T^2_M$ depends only on the conjugacy
class of $M$ in $\SL_2(\Z)$. Moreover, $T^2_{M^{-1}}$ and $T^2_M$
are related by an orientation reversing diffeomorphism. In
addition, one verifies that for all $M\in \SL_2(\Z)$
\[
S^{-1} M^t S = M^{-1},\quad S=\begin{pmatrix} 0 &-1\\
1&0\end{pmatrix}.
\]
This implies that in the case at hand, there also exists an
orientation reversing diffeomorphism $T^2_M\cong T^2_{M^t}$.
Since Rho invariants depends only on the oriented diffeomorphism
type of the mapping torus $T^2_M$, and the relation among Rho
invariants for different orientations is determined by Lemma
\ref{EtaProp} (ii), we are interested in elements of $\SL_2(\Z)$
only up to conjugation, taking inverses and transposes. Note that
$\gD$ in Definition \ref{TorusMapClasses} is invariant under
these operations so that $M$, $M^{-1}$ and $M^t$ all belong to
the same class.
\end{remark*}

We first collect some well-known facts, see for example
\cite[Sec. 1.4]{Shi}.

\begin{prop}\label{ConjugacyClasses}
Let $M\in \SL_2(\Z)$.
\begin{enumerate}
\item $M$ is parabolic if and only if $M$ is conjugate in
$\SL_2(\Z)$ to $\pm \left(\begin{smallmatrix} 1 &l\\
0&1\end{smallmatrix}\right)$ with $l\in \Z$.
\item $M$ is elliptic if and only if $M$ it is of finite order with
$M\neq \pm \Id$. In this case, $M$ is of order 3,4 or 6, and
conjugate in $\SL_2(\Z)$ to an element of the form
\[
\pm \begin{pmatrix}0 &-1\\
1&0 \end{pmatrix},\quad
\pm \begin{pmatrix} 0 &1\\
-1&1\end{pmatrix},\quad  \pm \begin{pmatrix} 1 &-1\\
1&0\end{pmatrix}.
\]
\end{enumerate}
\end{prop}

\begin{proof}[Sketch of proof]
The eigenvalues of $M$ are easily seen to be
\begin{equation}\label{SL2Z:Eigenvalues}
\gk = \lfrac 12 \big(\tr M + \sqrt \gD\big),\quad \gk^{-1} =
\lfrac 12 \big(\tr M - \sqrt \gD\big),
\end{equation}
where we fix the complex square root with $\sqrt{-1}=i$. By
definition, $M$ is parabolic if and only if $\gk=\gk^{-1}=\pm 1$,
so that the ``if'' part of (i) is clear. To prove the
``only if'' part, write $M=\left(\begin{smallmatrix} a &b\\
c&d\end{smallmatrix}\right)$, and assume that $M$ is parabolic
with $M\neq \pm \Id$. We can then assume---modulo conjugation by
$S$---that $c\neq 0$. Replacing $M$ with $-M$ if necessary we can
also achieve that $\tr M= a+d =2$. Define
\[
g:=\gcd\big(a-d,2c\big),\quad p:=\frac{a-d}g,\quad q:= \frac
{2c}g.
\]
It follows from $a+d=2$ and $ad-bc =1$ that
\[
ap+bq = p,\quad cp+dq =q.
\]
Moreover, $\gcd(p,q)=1$ so that we can find $r,s\in\Z$ with $pr-qs
=1$. Then $\left(\begin{smallmatrix} p &s\\
q&r\end{smallmatrix}\right)\in \SL_2(\Z)$, and one verifies that
\[
\begin{pmatrix} p &s\\ q&r \end{pmatrix}^{-1}
\begin{pmatrix} a &b\\ c&d \end{pmatrix}
\begin{pmatrix} p &s\\ q&r \end{pmatrix} =
\begin{pmatrix} 1 &l\\ 0&1 \end{pmatrix},
\]
for some $l\in \Z$. This proves part (i).

Concerning part (ii), we first note that part (i) implies that the
only parabolic elements of finite order are $\pm \Id$. Hence, we
can assume that $\gk\neq\gk^{-1}$. Then $M$ is conjugate in
$\GL_2(\C)$ to $\left(\begin{smallmatrix} \gk &0\\ 0 &\gk^{-1}
\end{smallmatrix}\right)$. Hence, $M$ is of finite order if and
only if $\gk$ is a root of unity. Then $\gk^{-1}= \Bar \gk$, and
$M$ is elliptic because
\[
|\tr M| = 2 |\Re(\gk)| <2,\quad \text{since}\quad \gk\neq\Bar \gk.
\]
For the reverse direction, we only note that if $M$ is elliptic,
then $\tr M\in \{-1,0,1\}$ and one easily checks by hand that
$\gk$ is a root of unity---in fact, $\gk= e^{i\frac {2\pi}3},
e^{i\frac \pi 2}$ or $e^{i\frac\pi 3}$. From this it is not
difficult to determine explicitly all conjugacy classes of
elliptic elements. We refer to \cite[pp. 14--15]{Shi}.
\end{proof}

We now analyze the action of $\SL_2(\Z)$ on $\H$ in some more
detail. Recall from Definition \ref{TorusFiberMet} that we let
$M=\left(\begin{smallmatrix} a &b\\ c&d\end{smallmatrix}\right)
\in \SL_2(\Z)$ act on $\H$ by the restriction of the fractional
linear transformation
\[
M^{\op}:\widehat \C\to \widehat \C,\quad M^{\op}z = \frac
{dz+b}{cz+a}.
\]
The following results are well-known but to fix notation we
sketch the proof.

\newpage

\begin{prop}\label{FixPointClass}
Let $M\in \SL_2(\Z)$ with $M\neq \pm \Id$.
\begin{enumerate}
\item If $M$ is parabolic of the form $M=\pm\left(\begin{smallmatrix} 1
&l\\ 0&1\end{smallmatrix}\right)$ with $l\neq 0$, then $M^{\op}$
has no fixed points in $\C$, and horizontal lines $\setdef{\gs\in
\H}{\Im(\gs)=\const}$ are invariant under the action of $M^{\op}$.
\item If $M$ is elliptic, the fractional linear transformation given by
$M^{\op}$ has exactly one fixed point in $\H$.
\item $M$ is hyperbolic, if and only if the fractional
transformation given by $M^{\op}$ has two distinct fixed points
$\ga,\gb \in \R \subset \C$, and the circle
\[
\Bigsetdef{\gs \in \H}{ \big|\gs-\lfrac{\ga+\gb}2 \big| =
\big|\lfrac{\ga-\gb}2\big|}
\]
is invariant under the action of $M^{\op}$.
\end{enumerate}
\end{prop}

\begin{proof}[Sketch of proof]
Part (i) is immediate since for all $z\in \C$ we have
$M^{\op}z=z+l$. If $M=\left(\begin{smallmatrix} a &b\\
c&d\end{smallmatrix}\right)$ is not parabolic, one easily verfies
that the eigenvalues $\gk$ and $\gk^{-1}$ as in
\eqref{SL2Z:Eigenvalues} cannot be integers. Hence, $M$ is not in
diagonal or triangular form, which implies that $b,c\neq 0$. Then
the fixed points of $M^{\op}$ acting on $\widehat \C$ are easily
seen to be
\begin{equation*}
\ga = \frac {\gk - a}c\quad\text{and}\quad \gb =
\frac{\gk^{-1}-a}c.
\end{equation*}
If $M$ is elliptic, then $\Im(\gk)>0$ and $\Bar \ga = \gb$. Thus,
the unique fixed point of $M^{\op}$ as claimed in part (ii) is
given by $\ga\in \H$ if $c>0$, and by $\gb\in \H$ if $c<0$. Let
us now assume that $M$ is hyperbolic. Then the eigenvalues are
real and $\gk>\gk^{-1}$. If we also assume for simplicity that
$c>0$, we get $\gb<\ga$. Then one verifies using elementary
linear algebra that for all $\gs\in \H$
\begin{equation}\label{FixPointClass:1}
\Big|\gs-\frac{\ga+\gb}2 \Big| =
\frac{\ga-\gb}2\quad\Longleftrightarrow\quad
\Re\Big(\frac{\gs-\ga}{\gs-\gb}\Big) = 0.
\end{equation}
On the other hand, $M^{\op}\gs$ is uniquely defined by the normal
form of the fractional linear transformation
\begin{equation}\label{FixPointClass:2}
\frac{M^{\op}\gs -\ga}{M^{\op}\gs -\gb} =
\gk^{-2}\frac{\gs-\ga}{\gs-\gb}.
\end{equation}
Since $M$ is hyperbolic we have $\gk\in \R$. Hence, one finds
from \eqref{FixPointClass:1} and \eqref{FixPointClass:2} that the
circle
\[
\Bigsetdef{\gs \in \H}{ \big|\gs-\lfrac{\ga+\gb}2 \big| =
\lfrac{\ga-\gb}2}
\]
is indeed invariant under the action of $M^{\op}$. This proves
part (iii).
\end{proof}

\subsection{The Elliptic Case}

Assume that $M=\left(\begin{smallmatrix} a&b\\ c&d
\end{smallmatrix}\right)\in \SL_2(\Z)$ is elliptic. Then,
according to Proposition \ref{ConjugacyClasses}, $M$ is of finite
order, so that we are in the situation considered in Section
\ref{FiniteMapTor}. As in Proposition \ref{FixPointClass} there
are precisely two fixed points for $M^{\op}$ acting on $\C$, one
of which lies in $\H$, explicitly given by
\[
\gs = \frac{\gk - a}c,\quad \Bar\gs = \frac{\Bar\gk - a}c,
\]
where,
\begin{equation}\label{EigenvalueElliptic}
\gk = \lfrac 12\big(\tr M + i\sqrt{4-(\tr M)^2} \big) = e^{2\pi i
\gt},\quad \gt\in (0,\lfrac 12).
\end{equation}
Actually, one easily checks that $\gt\in \{\lfrac 16, \lfrac
14,\lfrac 13\}$.

\begin{theorem}\label{RhoFiniteTorusBundle}
Let $M=\left(\begin{smallmatrix} a&b\\ c&d
\end{smallmatrix}\right)\in \SL_2(\Z)$ be elliptic, and let $A$
be a flat connection over $T^2_M$ defined by a pair $(a_\nu,u)$ as
in Proposition \ref{TorusBundleFlatConn}.
\begin{enumerate}
\item If $\nu\notin \Z^2$, and $\gt$ is as in
\eqref{EigenvalueElliptic}, then
\[
\rho_{A}(T^2_M) = (2-4\gt)\sgn(c).
\]
\item If $a_\nu=0$ is the trivial connection, so that
$u\equiv e^{-2\pi i \gl}\in \U(1)$, then
\[
\rho_{A}(T^2_M) = \begin{cases}\hphantom{\sgn} 0, &\text{\rm if
$\Re(u)<\Re(\gk)$},\\
\,\,\sgn(c) , &\text{\rm if
$\Re(u)=\Re(\gk)$},\\
2\sgn(c), &\text{\rm if $\Re(u)>\Re(\gk)$}.
\end{cases}
\]
\end{enumerate}
\end{theorem}

\begin{proof}
If $\nu\notin \Z^2$, it follows from Lemma \ref{TorusFlatConn}
(iii) that $H^\bullet(T^2,L_{a_\nu}) = \{0\}$. Hence, Theorem
\ref{RhoFiniteMapTor} yields that in this case the only
contribution to the Rho invariant comes from the Eta invariant of
the trivial connection. More precisely,
\begin{equation*}
\rho_A(T^2_M)=-4 \tr\log\big[M^*|_{\sH^+(T^2)\cap \gO^1}\big] + 2
\rk\big[(M^*-\Id)|_{\sH^+(T^2)\cap \gO^1}\big].
\end{equation*}
Let $\gs$ respectively $\Bar\gs$ be the fixed point of $M^{\op}$
in $\H$ as above, depending on whether $c>0$ or $c<0$. Use this to
define an $M$-invariant metric on $T^2$ as in Lemma
\ref{TeichmuellerAlmComp}. With the notation of
\eqref{Harmonic1FormsDef}, it follows from Corollary
\ref{HarmonicFormsTorus} that if $c>0$, then
\[
\sH^+(T^2)\cap \gO^1 = \C\go_\gs,\quad \sH^-(T^2)\cap \gO^1 =
\C\go_{\Bar\gs},
\]
and, if $c<0$, then
\[
\sH^+(T^2)\cap \gO^1 = \C\go_{\Bar\gs},\quad \sH^-(T^2)\cap \gO^1
= \C\go_\gs.
\]
Moreover, it is immediate that
\[
M^*\go_\gs = \gk\cdot \go_\gs,\quad\text{and}\quad
M^*\go_{\Bar\gs} = \Bar \gk \cdot\go_{\Bar\gs}.
\]
This implies that
\begin{equation}\label{RhoFiniteTorusBundle:1}
M^*|_{\sH^\pm (T^2)\cap \gO^1}= \gk,\quad M^*|_{\sH^\mp (T^2)\cap
\gO^1}=\Bar\gk,\quad \text{if $\pm c>0$}.
\end{equation}
Hence, with $\gt$ as in \eqref{EigenvalueElliptic} and the
definition of ``$\tr\log$'' in Theorem \ref{RhoFiniteMapTor}, one
finds that if $\nu\notin \Z^2$,
\[
\rho_A(T^2_M) = \begin{cases}\quad -4\gt +2, &\text{if $c>0$},\\
-4(1-\gt) +2,&\text{if $c<0$}.
\end{cases}
\]
This proves part (i) of Theorem \ref{RhoFiniteTorusBundle}. Now
assume that $a_\nu$ is the trivial connection. Lemma
\ref{TorusFlatConn} and Proposition \ref{TorusBundleFlatConn} then
imply that we can choose $\nu = 0$ and $u$ to be the constant
gauge transformation $e^{-2\pi i \gl}$, with $\gl\in [0,1)$. Then
part (i) and Theorem \ref{RhoFiniteMapTor} show that the Rho
invariant of $A$ is given by
\[
\begin{split}
2 \tr&\log\big[ u^{-1}M^*|_{\sH^+(T^2)\cap
\gO^1}\big] - \rk\big[(u^{-1}M^*-\Id)|_{\sH^+(T^2)\cap \gO^1}\big]\\
& - 2 \tr\log\big[ u^{-1}M^*|_{\sH^-(T^2)\cap \gO^1}\big]  +
\rk\big[(u^{-1}M^*-\Id)|_{\sH^-(T^2)\cap \gO^1}\big]
+(2-4\gt)\sgn(c).
\end{split}
\]
To compute the above quantities, we have to replace $M^*$ in
\eqref{RhoFiniteTorusBundle:1} with $u^{-1}M^*$. We assume for
simplicity that $c>0$; the other case works analogously. Then
\begin{equation*}
u^{-1}M^*|_{\sH^+ (T^2)\cap \gO^1}= u^{-1}\gk,\quad
u^{-1}M^*|_{\sH^- (T^2)\cap \gO^1}=u^{-1}\Bar\gk.
\end{equation*}
Now if $\Re(u)<\Re(\gk)$, then $\gl\in [0,\gt)$ or $\gl\in
(1-\gt,1)$. One readily verifies that in this case,
\[
\rk\big[(u^{-1}M^*-\Id)|_{\sH^+(T^2)\cap \gO^1}\big]=
\rk\big[(u^{-1}M^*-\Id)|_{\sH^-(T^2)\cap \gO^1}\big],
\]
and
\[
2 \tr\log\big[ u^{-1}M^*|_{\sH^+(T^2)\cap \gO^1}\big] - 2
\tr\log\big[ u^{-1}M^*|_{\sH^-(T^2)\cap \gO^1}\big] = 4\gt - 2.
\]
This implies that if $\Re(u)<\Re(\gk)$, then $\rho_A(T^2_M)=0$.
Similarly, if $\Re(u)>\Re(\gk)$, one computes that
\[
\rk\big[(u^{-1}M^*-\Id)|_{\sH^+(T^2)\cap \gO^1}\big]=
\rk\big[(u^{-1}M^*-\Id)|_{\sH^-(T^2)\cap \gO^1}\big],
\]
and
\[
2 \tr\log\big[ u^{-1}M^*|_{\sH^+(T^2)\cap \gO^1}\big] - 2
\tr\log\big[ u^{-1}M^*|_{\sH^-(T^2)\cap \gO^1}\big] = 4\gt,
\]
so that $\rho_A(T^2_M)=2$. Lastly, if $\Re(u)=\Re(\gk)$, then
either $\gl=\gt$ or $\gl=1-\gt$. In the first case,
\[
\rk\big[(u^{-1}M^*-\Id)|_{\sH^+(T^2)\cap \gO^1}\big]=1,\quad
\rk\big[(u^{-1}M^*-\Id)|_{\sH^-(T^2)\cap \gO^1}\big]=0,
\]
and
\[
2 \tr\log\big[ u^{-1}M^*|_{\sH^+(T^2)\cap \gO^1}\big]= 4\gt,\quad
2\tr\log\big[ u^{-1}M^*|_{\sH^-(T^2)\cap \gO^1}\big] = 0.
\]
This yields $\rho_A(T^2_M) =1$. In a similar way one deals with
the case $\gl=1-\gt$. If $c<0$, one has to replace $\gk$ with
$\Bar\gk$ in the above computations, and one easily verifies that
the result is the negative of what we computed in the case $c>0$.
\end{proof}

In part (i) of Theorem \ref{RhoFiniteTorusBundle} the twisting
connection does not contribute. Also recall that through Theorem
\ref{RhoFiniteMapTor} we have used Theorem \ref{RhoGen}, which
expresses the Rho invariant as the difference of adiabatic limits
of Eta invariants. Hence, part (i) of Theorem
\ref{RhoFiniteTorusBundle} can be rephrased as
\[
\lim_{\eps\to 0} \eta(B_\eps^{\ev}) = (4\gt-2)\sgn(c).
\]
Here, $B_\eps^{\ev}$ is the adiabatic family of untwisted odd
signature operators associated to the $M$-invariant metric on
$T^2$ induced by $\gs$ respectively $\Bar \gs$. Now in the case
at hand, the family $\eta(B_\eps^{\ev})$ is independent of
$\eps$, see \cite[p. 360]{Ati87}. Hence, we arrive at the
following

\begin{cor}\label{EtaElliptic}
Let $M=\left(\begin{smallmatrix} a &b\\
c&d\end{smallmatrix}\right)\in \SL_2(\Z)$ be elliptic, and endow
$T^2_M$ with the metric induced by an $M$-invariant metric on
$T^2$. Then, with $\gt$ is as in \eqref{EigenvalueElliptic},
\begin{equation*}
\eta(B^{\ev}) = (4\gt-2)\sgn(c).
\end{equation*}
\end{cor}

To check consistency with previous results, we now use Corollary
\ref{EtaElliptic} to compute the Eta invariant for the examples
considered in \cite[p. 372]{Ati87}, respectively \cite[p.
48]{Mor}.
\begin{enumerate}
\item $M=\begin{pmatrix}  0 &-1\\
1&1 \end{pmatrix}$. Then $\gk =e^{\frac{\pi i} 3}$, so that $\gt =
\lfrac 16$. Hence, $\eta(B^{\ev}) = \frac 46 -2 = -\frac 43$.
\item $M=\begin{pmatrix}  -1 &-1\\
1&0 \end{pmatrix}$. Then $\gk= e^{\frac {2\pi i}3}$ and
$\eta(B^{\ev}) = \frac 43 -2 = -\frac 23$.
\item $M=\begin{pmatrix}  0 &-1\\
1&0 \end{pmatrix}$. Then $\gk= e^{\frac {\pi i} 2}$, so that
$\eta(B^{\ev}) = \frac 44 -2 = -1$.
\end{enumerate}
Therefore, we obtain the same values as in \cite{Ati87} and
\cite{Mor}. Yet, the underlying abstract ideas we have used are
very different from what is considered there.\\

\subsection{The Parabolic Case}

If $M$ is parabolic, we know from Proposition
\ref{ConjugacyClasses} that $M$ is conjugate
to an element of the form $\pm \big(\begin{smallmatrix} 1 &l \\
0&1\end{smallmatrix}\big)$. Therefore, we can always choose
$\gs(t) := tl + i$ as an $M$-invariant path in $\H$. We first
compute the integral of the Rho form using Theorem
\ref{TorBundleRhoForm}.

\begin{prop}\label{RhoFormParab}
Let $\nu=(\nu_1,\nu_2)\in \R^2$ with $\nu\notin \Z^2$ satisfy
$(M^t-\Id)\nu \in \Z^2$. Let $A$ be the corresponding flat
connection over the mapping torus $T^2_M$. Then the integral of
the Rho form with respect to the metric induced by $\gs(t) := tl
+ i$ is given by
\[
\int_0^1 \widehat \rho_A = l\big(P_2(\nu_1)-\lfrac 16\big) + \frac
l{2\pi},
\]
where $P_2$ is the second periodic Bernoulli function.
\end{prop}

\begin{proof}
Since both sides of the equation in Proposition
\ref{RhoFormParab} depend on $\nu$ only modulo $\Z^2$, we can
assume that $\nu\in [0,1)^2$. Since $\gs(0)=i$ and $\gs(1)= l +
i$, we have
\[
\cot(\pi \gs(1) n) = \cot(\pi \gs(0)n)\quad \text{for all $n\in
\N$}.
\]
This implies---using the notation of Theorem \ref{F:Thm} and
Remark \ref{E:Rem} (ii)---that
\[
E_0\big(\gs(1)\big) = E_0\big(\gs(0)\big).
\]
We now claim that also
\begin{equation}\label{RhoFormParab:1}
E_\nu\big(\gs(1)\big) = E_\nu\big(\gs(0)\big).
\end{equation}
Indeed, if $M= \big(\begin{smallmatrix} 1 &l \\
0&1\end{smallmatrix}\big)$, the condition $(M^t-\Id)\nu \in \Z^2$
guarantees that $l\nu_1 \in \Z$, whereas $\nu_2$ is arbitrary.
Thus, with $z(t) = \nu_1\gs(t) - \nu_2$ as in Theorem \ref{F:Thm},
we get
\[
z(1) = z(0) + l\nu_1 \in z(0)+ \Z.
\]
This easily yields \eqref{RhoFormParab:1} in the case at hand.
If $M= - \big(\begin{smallmatrix} 1 &l \\
0&1\end{smallmatrix}\big)$, then $(M^t-\Id)\nu \in \Z^2$ means
that
\[
2\nu_1\in \Z,\quad l\nu_1 + 2\nu_2 \in \Z.
\]
Since we are assuming that $\nu\in [0,1)^2$, there are only a few
possible values for $\nu$. First of all, if $\nu_1 =0$, then $z(1)
= z(0)$, so that \eqref{RhoFormParab:1} holds again. If $\nu_1=
\lfrac 12$, then
\[
\nu_2\in \begin{cases} \{0, \lfrac 12\}, &\text{if $l$ is even},\\
\{\lfrac14 , \lfrac 34\}, &\text{if $l$ is odd}.\end{cases}
\]
Moreover, we have $z(1)=z(0) + \frac l2$ and so
\begin{equation}\label{RhoFormParab:2}
\begin{split}
\cos\big(2\pi z(1)n\big) &= (-1)^{nl}\cos\big(2\pi z(0 )n\big),\\
\sin\big(2\pi z(1)n\big) &= (-1)^{nl}\sin\big(2\pi z(0 )n\big).
\end{split}
\end{equation}
This implies that only summands such $nl$ is odd contribute to
$E_\nu\big(\gs(1)\big) - E_\nu\big(\gs(0)\big)$. Hence, if $l$ is
even we are done. Let us thus assume that $l$ is odd and---for
definiteness---that $\nu_2 = \lfrac 14$. The case $\nu_2 = \lfrac
34$ is analogous. Then, for odd $n\in \N$,
\[
\cos\big(2\pi z(0)n\big) = \cos\big(\pi i n - \lfrac \pi 2 n\big)
= i^{n-1}\sin(\pi i n),\quad \sin\big(2\pi z(0)n\big) = -
i^{n-1}\cos(\pi i n),
\]
and so
\[
\cos\big(2\pi n z(0)\big)\cot(\pi n i ) + \sin\big(2\pi n
z(0)\big) = 0.
\]
This, together with \eqref{RhoFormParab:2}, implies
\eqref{RhoFormParab:1} in this last case as well. We can now use
\eqref{RhoFormParab:1} and Theorem \ref{TorBundleRhoForm} to
deduce that in all cases
\[
\int_0^1 \widehat \rho_A = l\big(P_2(\nu_1)- \lfrac 16) + \frac
1{2\pi}\int_0^1 ldt = l\big(P_2(\nu_1)- \lfrac 16) + \frac
l{2\pi}.\qedhere
\]
\end{proof}

To conclude the computation of the Rho invariants of $T^2_M$ for
parabolic $M\in \SL_2(\Z)$, we still have to determine the Rho
invariant of the bundle of vertical cohomology groups
$\rho_{\sH_{A,v}}(S^1)$ appearing in Theorem \ref{RhoGen}. Recall
that Dai's correction term vanishes because the base is
1-dimensional, see the proof of Theorem \ref{RhoFiniteMapTor}.

\begin{prop}\label{RhoCohomParab}
For $\eps=\pm 1$ and $l\in \Z$ let $T^2_M$ be the mapping
torus of $M=\eps \big(\begin{smallmatrix} 1 &l \\
0&1\end{smallmatrix}\big)$. Endow $T^2_M$ with the metric given
by $\gs(t) = tl+ i$. Then, for all connections $A$ as in
Proposition \ref{RhoFormParab},
\[
\rho_{\sH_{A,v}}(S^1) = \begin{cases} \hphantom{-\,-\,}0,
&\text{\rm if
$l=0$,}\\ - \lfrac l\pi + \sgn (l), &\text{\rm if $\eps=1$, $l\neq 0$,}\\
\hphantom{-\,}- \lfrac l\pi, &\text{\rm if $\eps=-1$.}
\end{cases}
\]
\end{prop}

The proof turns out to be somewhat involved, and we sketch the
strategy first. We know from Corollary \ref{HarmonicFormsTorus}
that the twisted cohomology groups of $T^2$ vanish except for the
trivial connection. Hence, for $A$ as in Proposition
\ref{RhoFormParab} we can argue as in the proof of part (i) of
Theorem \ref{RhoFiniteTorusBundle} that
\begin{equation}\label{RhoCohomParabRem}
\rho_{\sH_{A,v}}(S^1) = - \lfrac12\eta\big(D_{S^1}\otimes
\nabla^{\sH_v}\big),
\end{equation}
where $D_{S^1}\otimes \nabla^{\sH_v}$ is as in Definition
\ref{CohomTwistSignOp}. However, unlike in the case of elliptic
elements, the connection $\nabla^{\sH_v}$ on the bundle of
vertical cohomology groups is not unitary, so that it is
difficult to compute the above Eta invariant directly. The idea
of our proof is to study the difference between $D_{S^1}\otimes
\nabla^{\sH_v}$ and the odd signature operator associated to the
unitary connection $\nabla^{\sH_v,u}$, see
\eqref{UnitaryConnHodge} and \eqref{CohomUnitaryTwistSignOp}.
More precisely, we will compute $\eta\big(D_{S^1}\otimes
\nabla^{\sH_v,u}\big)$ and then use the variation formula of
Proposition \ref{EtaDiffSF} to obtain $\eta\big(D_{S^1}\otimes
\nabla^{\sH_v}\big)$. Here, the considerations of Section
\ref{IndefiniteMetric} will play a role.

\begin{proof}[Proof of Proposition \ref{RhoCohomParab}]
We split the bundle of vertical cohomology groups as
\[
\sH_v^\bullet(T^2_M) = \sH_v^0(T^2_M) \oplus \sH_v^1(T^2_M)
\oplus \sH_v^2(T^2_M).
\]
It follows from Corollary \ref{HarmonicFormsTorus}, that
$\sH_v^0(T^2_M)$ and $\sH_v^2(T^2_M)$ can be trivialized by the
constant sections $1$ respectively $dx\wedge dy$. With respect to
this trivialization, the connection $\nabla^{\sH_v}$ is the
trivial connection, see \eqref{CanVerticalConnMapTor}. According
to Remark \ref{RhoRem} (iii) the Eta invariant of the untwisted
odd signature operator over $S^1$ vanishes, so that we only have
to compute the contribution to $\eta\big(D_{S^1}\otimes
\nabla^{\sH_v}\big)$ coming from $\sH_v^1(T^2_M)$.

Let $\go_{\gs(t)}$ and $\go_{\Bar \gs(t)}$ be as in
\eqref{Harmonic1FormsDef} with respect to $\gs(t) = tl +i$, and
define
\[
\ga_t := \go_{\gs(t)} + \go_{\Bar \gs(t)},\quad \gb_t :=
\go_{\gs(t)} - \go_{\Bar \gs(t)}.
\]
It is immediate from Corollary \ref{HarmonicFormsTorus} that for
each $t$ the pair $(\ga_t,\gb_t)$ forms an orthogonal basis of
$\sH^1(T^2_M,g_{\gs(t)})$. However, it is not necessarily a
trivialization of the bundle of vertical cohomology groups. For
this note that, with $\eps=\pm 1$ as in the statement of the
proposition,
\[
M^* \go_{\gs(t)} = \eps\go_{\gs(t+1)},\quad M^* \go_{\Bar\gs(t)}
= \eps\go_{\Bar\gs(t+1)},
\]
so that also
\[
M^*\ga_t = \eps \ga_{t+1},\quad M^*\gb_t = \eps \gb_{t+1}.
\]
Nevertheless, this means that we can write every section of
$\sH_v^1(T^2_M)\to S^1$ as
\[
\gf_\ga(t) \ga_t + \gf_\gb(t) \gb_t,
\]
with functions $\gf_\ga$ and $\gf_\gb$ on $\R$ satisfying the
condition
\begin{equation}\label{ParabBoundCond}
\gf_\ga(t+1) = \eps\gf_\ga(t),\quad\gf_\gb(t+1) = \eps\gf_\gb(t).
\end{equation}
Now, note that
\[
\pd_t \go_{\gs(t)}= \pd_t \go_{\Bar\gs(t)} = l dy = \frac
{l}{2i}\big(\go_{\gs(t)} - \go_{\Bar\gs(t)}\big).
\]
Using this we see that the flat connection $\nabla^{\sH_v}$ on
$\sH_v^1(T^2_M)$ is given by
\begin{equation}\label{ParabCanConn}
\nabla^{\sH_v}_{\pd_t} \ga_t = \pd_t \ga_t = -il \gb_t,\quad
\nabla^{\sH_v}_{\pd_t} \gb_t = \pd_t \gb_t =0.
\end{equation}
Moreover, one verifies using \eqref{Harmonic1FormsProp} and
\eqref{DotTauTorus} that
\begin{equation}\label{ParabTau}
\tau_t \ga_t =\gb_t,\quad \tau_t \gb_t =\ga_t,\quad \dot \tau_t
\ga_t = il \ga_t,\quad \dot \tau_t \gb_t = - il \gb_t.
\end{equation}
According to Lemma \ref{VerticalConnMapTor}, this means that the
unitary connection $\nabla^{\sH_v,u}$ of \eqref{UnitaryConnHodge}
on $\sH_v^1(T^2_M)$ is given by
\begin{equation}\label{ParabUnitaryConn}
\nabla^{\sH_v,u}_{\pd_t} \ga_t = - \frac{il}2 \gb_t,\quad
\nabla^{\sH_v,u}_{\pd_t} \gb_t = - \frac{il}2 \ga_t.
\end{equation}
With the operators of Definition \ref{CohomTwistSignOp} and
\eqref{CohomUnitaryTwistSignOp} we introduce the abbreviations
\[
D:= D_{S^1}\otimes \nabla^{\sH_v}|_{C^\infty(S^1,
\sH^1_v(T^2_M))},\quad D^u:=D_{S^1}\otimes
\nabla^{\sH_v,u}|_{C^\infty(S^1, \sH^1_v(T^2_M))}.
\]
Then, using the splitting of the total chirality operator as in
Lemma \ref{TauSplit} one verifies that
\[
D= -i\tau_t \nabla^{\sH_v}_{\pd_t},\quad D^u= -i\tau_t
\nabla^{\sH_v,u}_{\pd_t}.
\]
From the considerations at the beginning of this proof we deduce
that
\begin{equation}\label{ParabEta}
\eta(D) = \lfrac12\eta\big(D_{S^1}\otimes \nabla^{\sH_v}\big),
\end{equation}
where the factor enters for the same reason as in Remark
\ref{CohomTwistSignOpRem} (ii) and Remark \ref{OddSignRem} (i).
Similarly,
\[
\eta(D^u) = \lfrac12\eta\big(D_{S^1}\otimes \nabla^{\sH_v,u}\big).
\]
Now let
\[
D_s := D^u + s(D-D^u), \quad s\in [0,1].
\]
Then $D_s$ is a family of self-adjoint operators, which is
precisely of the form considered in Theorem
\ref{CohomTwistEtaDiff}. Hence, the ``local variation'' of the
Eta invariant vanishes , so that the general variation formula of
Proposition \ref{EtaDiffSF} reduces to
\begin{equation}\label{ParabEtaDiff}
\eta(D) = \eta(D^u) + \dim(\ker D^u) -\dim(\ker D)  +  2
\SF(D_s)_{s\in[0,1]}.
\end{equation}
To determine the terms appearing in \eqref{ParabEtaDiff}, we now
explicitly compute $\spec(D_s)$. If $\gf_\ga \ga_t + \gf_\gb
\gb_t$ is a section of $\sH^1_v(T^2_M)$, then
\eqref{ParabCanConn}, \eqref{ParabTau} and
\eqref{ParabUnitaryConn} imply that $D_s$ acts in terms of the
coordinate functions $(\gf_\ga,\gf_\gb)$ as
\[
D_s\begin{pmatrix}\gf_\ga \\ \gf_\gb \end{pmatrix} = -i \left[
\begin{pmatrix}0 &\pd_t\\ \pd_t &0\end{pmatrix} - \frac
{il}2\begin{pmatrix}1+s &0\\ 0 & 1-s\end{pmatrix} \right]
\begin{pmatrix}\gf_\ga \\ \gf_\gb \end{pmatrix}.
\]
Hence, if $\gf_\ga \ga_t + \gf_\gb\gb_t$ is an eigenvector with
eigenvalue $\gl(s)\in \R$, then
\begin{equation}\label{DsCoord}
\pd_t \begin{pmatrix}\gf_\ga \\ \gf_\gb \end{pmatrix} =
i\begin{pmatrix} 0 &\gl(s)+(1-s)\lfrac l2\\
\gl(s)+(1+s)\lfrac l2 &0\end{pmatrix}
\begin{pmatrix}\gf_\ga \\ \gf_\gb \end{pmatrix} =: iT_{\gl(s)}
\begin{pmatrix}\gf_\ga \\ \gf_\gb
\end{pmatrix},
\end{equation}
which is an ordinary linear differential equation with constant
coefficients. Let us assume from now on that $l\neq 0$. The case
$l=0$ will be dealt with separately at the end. The characteristic
equation for the eigenvalues $\gk$ of $T_{\gl(s)}$ is
\[
\gk^2 = \big(\gl(s)+\lfrac l2\big)^2 -  s^2\big(\lfrac l2\big)^2.
\]
Therefore, unless $\big(\gl(s)+\lfrac l2\big)^2 =  s^2\big(\lfrac
l2\big)^2$, the matrix $T_{\gl(s)}$ has two distinct eigenvalues
$\gk$ and $-\gk$, and can be brought into diagonal form. Now a
solution to \eqref{DsCoord} satisfies the condition
\eqref{ParabBoundCond} if and only if $e^{i\gk} = \eps$. Write
$\eps = e^{i\gt}$ with $\gt\in\{0,\pi\}$. Then the condition
$e^{i\gk} = \eps$ is equivalent to $\gk = 2\pi n +\gt$ with $n\in
\N$. Hence, $\gl(s)$ is an eigenvalue of $D_s$ if and only if
\[
\gl(s) = \gl^\pm_n(s) := \pm \sqrt{\gk_n^2 + s^2\big(\lfrac
l2\big)^2} - \lfrac l2,\quad \gk_n = 2\pi n +\gt\,\text{ for some
$n\in \N$.}
\]
Moreover, unless $n=0$ and $\gt=0$, the standard procedure for
solving \eqref{DsCoord} gives us two linearly independent
solutions for both eigenvalues $\gl^+_n(s)$ and $\gl^-_n(s)$. In
the special case $n=0$ and $\gt=0$ we have
\[
\gl^\pm_0(s) = \lfrac 12(\pm s |l| - l),\quad \text{and so}\quad
T_{\gl^{\pm}_0(s)} = \frac s2 \begin{pmatrix} 0 & \pm |l| - l  \\
 \pm |l| + l  &0\end{pmatrix}.
\]
If $s=0$, this matrix vanishes so that the eigenvalues
$\gl^{\pm}_0(0)$ have multiplicity 2. If $s\neq 0$, this matrix
is in triangular form. This implies that the eigenvalues
$\gl^{\pm}_0(s)$ have multiplicity 1. To see this, consider for
example the case $l>0$, and $T= T_{\gl^+_0(s)}$. Then
\[
e^{itT} = \begin{pmatrix} 1 & 0  \\
0  & 1\end{pmatrix} + i t\begin{pmatrix} 0 & 0  \\
sl  & 0\end{pmatrix},
\]
so that the only solutions of \eqref{DsCoord} satisfying the
condition \eqref{ParabBoundCond} are constant multiples of
$(\gf_\ga,\gf_\gb)= (0,1)$.

Since the operator $D^u$ coincides with $D_0$ we see in
particular that
\[
\spec(D^u) = \bigsetdef{\pm (2\pi n+ \gt)-\lfrac l2}{n\in \N},
\]
where all eigenvalues have multiplicity 2. Hence, for $\Re(z)>1$
\[
\eta(D^u,z) = 2 \sum_{n \in\Z}\frac {\sgn(2\pi n +\gt -\lfrac
l2)}{|2\pi n +\gt -\lfrac l2|^z} =  2 (2\pi)^z \sum_{n \in\Z}\frac
{\sgn(n -\lfrac {l-2\gt}{4\pi})}{|n -\lfrac {l-2\gt}{4\pi}|^z},
\]
which up to a factor is the Eta function considered in
Proposition \ref{ZEtaCalc}. Hence, the value at $z=0$ of the
meromorphic continuation of $\eta(D^u,z)$ is given as follows:
Let $m\in \Z$ be such that
\begin{equation}\label{ParabUnitaryEta}
\frac l{4\pi}- \frac \gt{2\pi} - m \in
(0,1).\quad\text{Then}\quad  \eta(D^u) = \frac l{\pi}
-2\frac{\gt}\pi - 4m -2.
\end{equation}
This identifies the first term in \eqref{ParabEtaDiff}.

To compute the spectral flow term, we first assume that $l>0$.
Then the zero eigenvalues of $D_s$ for $s\in(0,1)$ are given by
those $\gl^+_n(s)$ for which $s$ and $n$ are related by
\[
\gk_n^2 = (1-s^2)(\lfrac l2)^2.
\]
The family $\gl^+_n(s)$ is strictly increasing with $s$ and all
eigenvalues have multiplicity 2. For the latter note that since we
are assuming that $s\neq 1$, the eigenvalues of multiplicity 1,
which we have found in the case $\gt=0$, are never zero.
Therefore, each zero eigenvalue will contribute $+2$ to
$\SF(D_s)_{s\in[0,1]}$. Since $1-s^2$ maps $(0,1)$ bijectively
onto itself, we have to count the number of $n\in \N$ for which
$0<\gk_n < \frac l2$, or---equivalently---for which $-\frac
{\gt}{2\pi} < n  < \frac l{4\pi}-\frac {\gt}{2\pi}$. Now, with
$m$ as in \eqref{ParabUnitaryEta}, it is immediate to check that
\[
\#\bigsetdef{n\in \N}{-\lfrac {\gt}{2\pi} < n  < \lfrac
l{4\pi}-\lfrac {\gt}{2\pi}} = \begin{cases} \hphantom{+}m,
&\text{if
$\gt=0$,}\\
m+1 , &\text{if $\gt =\pi$.}
\end{cases}
\]
Note that since we are assuming that $l>0$, we certainly have
$m\ge 0$ if $\gt=0$ and $m\ge -1$ if $\gt=\pi$. Concerning the
endpoints of the path, there are no zero eigenvalues for $s=0$.
For $s=1$ we only have one if $\gt = 0$, and this is the
eigenvalue $\gl_0^+(1)$ of multiplicity 1. Putting all
information together, we find that for $l>0$,
\[
\dim(\ker D^u) -\dim(\ker D)  +  2 \SF(D_s)_{s\in[0,1]} =
\begin{cases} -1 + 2(2m+1),&\text{if $\gt=0$,}\\
\hphantom{-1\,}4(m+1), &\text{if $\gt=\pi$.}\end{cases}
\]
Together with \eqref{ParabUnitaryEta} and \eqref{ParabEtaDiff},
we conclude that in the case that $l>0$,
\[
\eta(D) = \begin{cases} \lfrac l\pi - 1, &\text{if $\gt=0$,}\\
\hphantom{-\,}\lfrac l\pi, &\text{if $\gt=\pi$.}
\end{cases}
\]
Let us now assume that $l<0$. Then the role of $\gl^+_n(s)$ in
the preceding discussion is replaced by $\gl^-_n(s)$, which
strictly decreases with $s$. Hence, the contribution to the
spectral flow is $-2$ for each zero eigenvalue. With $m$ as in
\eqref{ParabUnitaryEta} we now have $m\le -1$ and for both values
of $\gt$
\[
\#\bigsetdef{n\in \N}{-\lfrac {\gt}{2\pi} < n  < - \lfrac
l{4\pi}-\lfrac {\gt}{2\pi}} = -m-1.
\]
For $l<0$ and $\gt=0$, the zero eigenvalue $\gl^-_0(1)$ of
multiplicity 1 does not contribute to the spectral flow.
Therefore, we arrive at
\[
\dim(\ker D^u) -\dim(\ker D)  +  2 \SF(D_s)_{s\in[0,1]} =
\begin{cases} -1 + 4(m+1),&\text{if $\gt=0$,}\\
\hphantom{-1\,}4(m+1), &\text{if $\gt=\pi$.}\end{cases}
\]
Hence, we conclude that for $l<0$
\[
\eta(D) = \begin{cases} \lfrac l\pi + 1, &\text{if $\gt=0$,}\\
\hphantom{-\,}\lfrac l\pi, &\text{if $\gt=\pi$.}
\end{cases}
\]
Hence, using \eqref{RhoCohomParabRem} and \eqref{ParabEta} we
have proved Proposition \ref{RhoCohomParab} in the case that
$l\neq 0$. If $l=0$, one easily checks that
\[
\spec(D) = \bigsetdef{2\pi n + \gt}{n\in \Z},
\]
where all eigenvalues have multiplicity 2. This implies that
$\spec(D)$ is symmetric, so that $\eta(D)=0$.
\end{proof}

We can now combine Propositions \ref{RhoFormParab} and
\ref{RhoCohomParab} to obtain the formula for $U(1)$-Rho
invariants for mapping tori with parabolic monodromy. According
to Theorem \ref{RhoGen}, we have
\[
\rho_A(T^2_M) = 2 \int_0^1 \widehat \rho_A +
\rho_{\sH_{A,v}}(S^1).
\]
Therefore, we arrive at the following

\begin{theorem}\label{RhoParab}
Let $\eps=\pm 1$ and $l\in \Z$, and let $T^2_M$ be the mapping
torus of the parabolic element $M=\eps \big(\begin{smallmatrix} 1 &l \\
0&1\end{smallmatrix}\big)$. Let $A$ be a flat connection over the
mapping torus $T^2_M$, defined by $\nu\in \R^2$ with $\nu\notin
\Z^2$, satisfying $(M^t-\Id)\nu \in \Z^2$. If $l=0$, the Rho
invariant $\rho_A(T^2_M)$ vanishes. For $l\neq 0$ we have
\[
\rho_A(T^2_M) = 2l\big(P_2(\nu_1)-\lfrac 16\big) + \begin{cases}
\sgn (l), &\text{\rm if $\eps=1$}\\ \hphantom{sg} 0, &\text{\rm if
$\eps=-1$.}
\end{cases}
\]
\end{theorem}

\newpage

\begin{remark}\label{ParabRhoRem}\quad
\begin{enumerate}
\item We want to point out that the assumption that $\nu\notin
\Z^2$ excludes possibly non-trivial flat connections on $T^2_M$
which restrict to the trivial connection over $T^2$. Note that for
elliptic elements in Theorem \ref{RhoFiniteTorusBundle} (ii) we
included a discussion. However, in the case of parabolic
elements---and also in the hyperbolic case below---the case
$\nu\notin \Z^2$ is much more interesting and a parallel treatment
of the remaining case would lead to more notational inconvenience
and a tedious distinction between all cases. Since the insight
gained seemed not worth the effort, we opted to work under the
assumption that $\nu\notin \Z^2$ only.
\item Note that if $\eps =1$ in Theorem \ref{RhoParab}, then $\nu_1
l\in \Z$, so that $\nu_1 = k/l$ for some $k\in \Z$. Hence, the
formula for the Rho invariant is the same as the formula for the
Rho invariant for a principal circle bundle of degree $l$ over
$T^2$ in Theorem \ref{RhoCircBund}. The underlying reason is that
for $\eps=1$, the mapping torus $T^2_M$ is at the same time a
principal $S^1$-bundle of degree $l$ over $T^2$, see \cite[p.
470]{Sco83}.
\end{enumerate}
\end{remark}

\subsection{The Hyperbolic Case}

Now we turn to the generic---and most interesting---case that $M$
is hyperbolic. This section is less self-contained than the
previous sections, since we will deduce the main result from a
well-known transformation formula for certain generalized Dedekind
Eta functions. Since this would lead to far afield, we shall not
attempt to give a detailed treatment but refer to the literature
for proofs.\\

\noindent\textbf{$\boldsymbol{M}$-invariant Paths in
the upper half plane.} Assume that $M=\big(\begin{smallmatrix} a &b \\
c&d\end{smallmatrix}\big)\in \SL_2(\Z)$ is hyperbolic. As in the
proof of Proposition \ref{FixPointClass}, we know that $b,c\neq
0$, and that the fixed points of $M^{\op}$ acting on $\widehat \C$
are given by
\begin{equation}\label{AlphaBetaDef}
\ga = \frac {\gk - a}c,\quad \gb =
\frac{\gk^{-1}-a}c,\quad\text{where}\quad \gk = \lfrac 12\big(a+d
+ \sqrt\gD\big).
\end{equation}
For $t\in \R$ define
\begin{equation}\label{SigmaHypDef}
\gs(t):=\frac 1{|\gk|^{2t}+|\gk|^{-2t}}\big(\ga |\gk|^{2t} +\gb
|\gk|^{-2t} + i|\ga-\gb| \big).
\end{equation}

\begin{lemma}\label{SigmaHypLem}
The path $\gs(t)$ in $\H$ lies on the circle
\[
\Bigsetdef{\gs \in \H}{ \big|\gs-\lfrac{\ga+\gb}2 \big| =
\big|\lfrac{\ga-\gb}2\big|}
\]
and satisfies $M^{\op}\gs(t) = \gs(t+1)$.
\end{lemma}

\begin{proof}
We proof the second assertion first. By comparison with
\eqref{FixPointClass:2}, we thus have to show that
\begin{equation}\label{SigmaHypLem:1}
\frac{\gs(t+1)-\ga}{\gs(t+1)-\gb} = \gk^{-2}
\frac{\gs(t)-\ga}{\gs(t)-\gb}
\end{equation}
Now,
\[
\begin{split}
\gs(t) - \ga &= \frac 1{|\gk|^{2t}+|\gk|^{-2t}}\big(\ga |\gk|^{2t}
+\gb |\gk|^{-2t} - \ga |\gk|^{2t} - \ga |\gk|^{-2t} + i|\ga-\gb|
\big)\\
&= \frac 1{|\gk|^{2t}+|\gk|^{-2t}}\big((\gb-\ga)|\gk|^{-2t} +
i|\ga-\gb| \big),
\end{split}
\]
and, similarly,
\[
\gs(t) - \gb = \frac
1{|\gk|^{2t}+|\gk|^{2t}}\big((\ga-\gb)|\gk|^{2t} + i|\ga-\gb|
\big).
\]
Let $t,t'\in \R$, and abbreviate $\eps:=\sgn(\ga-\gb)=\sgn(c)$.
Then
\[
\begin{split}
\frac{\big(\gs(t')-\ga\big)\big(\gs(t)-\gb\big)}{\big(\gs(t')-
\gb\big)\big(\gs(t)-\ga\big)} &= \frac{\big(-|\gk|^{-2t'}
+i\eps\big)\big(|\gk|^{2t}+
i\eps\big)}{\big(|\gk|^{2t'}+i\eps\big)\big(-|\gk|^{-2t}+
i\eps\big)}\\ &= \frac{ -|\gk|^{-2(t'-t)}-
1+i\eps\big(|\gk|^{2t}-|\gk|^{-2t'}\big)}{ -|\gk|^{2(t'-t)}-
1+i\eps\big(|\gk|^{2t'}-|\gk|^{-2t}\big)} = |\gk|^{-2(t'-t)}.
\end{split}
\]
In particular, for $t'=t+1$ we obtain the formula in
\eqref{SigmaHypLem:1}, which in turn shows that $\gs(t)$ is
$M$-invariant. Moreover,
\[
\frac{\gs(t)-\ga}{\gs(t)-\gb} = \gk^{-2t}
\frac{\gs(0)-\ga}{\gs(0)-\gb}.
\]
Now, as
\[
\gs(0) = \frac{\ga+\gb}2 + i \frac{|\ga-\gb|}{2},
\]
this implies as in the proof of Proposition \ref{FixPointClass}
that all points $\gs(t)$ lie on the circle
\[
\Bigsetdef{\gs \in \H}{ \big|\gs-\lfrac{\ga+\gb}2 \big| =
\big|\lfrac{\ga-\gb}2\big|}.\qedhere
\]
\end{proof}

\noindent\textbf{The Rho invariant of the bundle of vertical
cohomology groups.} Having found an $M$-invariant path in $\H$ we
now need to compute the integral over the Rho form and the Rho
invariant of the bundle of vertical cohomology groups. We start
with the latter, which is more straightforward than in the
parabolic case.

\begin{prop}\label{RhoCohomHyp}
Let $T^2_M$ be the mapping torus of a hyperbolic element
$M\in\SL_2(\Z)$. Endow $T^2_M$ with the metric given by
\eqref{SigmaHypDef}, and let $A$ be a flat connection determined
by $\nu=(\nu_1,\nu_2)\in \R^2$ with $\nu\notin \Z^2$ and
$(M^t-\Id)\nu \in \Z^2$. Then
\[
\rho_{\sH_{A,v}}(S^1) = 0.
\]
\end{prop}

\begin{proof}
Again we know from Corollary \ref{HarmonicFormsTorus} that the
twisted cohomology groups of $T^2$ vanish except for the case
that the underlying connection is the trivial one. Thus, as in
\eqref{RhoCohomParabRem}
\begin{equation*}
\rho_{\sH_{A,v}}(S^1) = - \lfrac12\eta\big(D_{S^1}\otimes
\nabla^{\sH_v}\big).
\end{equation*}
Moreover, as explained in the proof of Proposition
\ref{RhoCohomParab}, we only need to study the restriction of
$D_{S^1}\otimes \nabla^{\sH_v}$ to
$C^\infty\big(S^1,\sH^1_v(T^2_M)\big)$.

In view of the rather complicated formula for $\gs(t)$ it is
inconvenient to work directly with the basis
$(\go_{\gs(t)},\go_{\Bar \gs(t)})$ of $\sH^1(T^2_M,g_{\gs(t)})$,
given by Corollary \ref{HarmonicFormsTorus}. Instead we define
\[
\go_\ga(t) := \frac {|\gk|^t - i \eps|\gk|^{-t}}2\, \go_{\gs(t)}
+ \frac {|\gk|^t + i \eps |\gk|^{-t}}2\, \go_{\Bar \gs(t)},
\]
and
\[
\go_\gb(t) := \frac {|\gk|^{-t} + i \eps|\gk|^t}2\, \go_{\gs(t)}
+ \frac {|\gk|^{-t} - i \eps|\gk|^t}2\, \go_{\Bar \gs(t)},
\]
where as before $\eps = \sgn(c)$. Then
$\big(\go_\ga(t),\go_\gb(t)\big)$ is a clearly basis of
$\sH^1(T^2_M,g_{\gs(t)})$ for each $t$. Moreover,
\eqref{Harmonic1FormsProp} implies that
\begin{equation}\label{RhoCohomHyp:1}
\tau_t \go_\ga(t) = -i \go_\gb(t),\quad \tau_t\go_\gb(t) =
i\go_\ga(t),
\end{equation}
where $\tau_t$ is the chirality operator defined by $\gs(t)$.
Now, a straightforward calculation---which we skip---shows that
\begin{equation}\label{RhoCohomHyp:2}
\go_\ga(t) = |\gk|^t (dx + \ga dy),\quad \go_\gb(t) = |\gk|^{-t}
(dx + \gb dy).
\end{equation}
Thanks to this identity, we obtain---without having to compute
the derivatives of $\go_{\gs(t)}$ and $\go_{\Bar \gs(t)}$
explicitly---that
\begin{equation}\label{RhoCohomHyp:3}
\pd_t \go_\ga(t) = \log|\gk|\cdot \go_\ga(t),\quad
\pd_t\go_\gb(t)= -\log|\gk|\cdot \go_\gb(t).
\end{equation}
From the definition of $\ga$ and $\gb$ in \eqref{AlphaBetaDef} we
immediately see that $c\ga + a=\gk$ and $c\gb + a = \gk^{-1}$.
Moreover, using that $ad-bc=1$ and that $\gk+\gk^{-1}= a+d$, one
computes
\[
d\ga + b = \frac {d\gk - ad}{c} + b = \gk \frac {d-\gk^{-1}}{c} =
\gk \ga,\quad d\gb + b=\ldots= \gk^{-1}\gb.
\]
This means that $(1,\ga)$ and $(1,\gb)$ are eigenvectors of
$M^t=\left(\begin{smallmatrix} a & c\\ b&d
\end{smallmatrix}\right)$ with eigenvalues $\gk$ and
$\gk^{-1}$, respectively. Therefore, it follows from
\eqref{RhoCohomHyp:2} that
\[
M^* \go_\ga(t) = |\gk|^t\big((a+c\ga) dx + (b +d\ga) dy\big) =
|\gk|^t \big(\gk dx + \ga\gk dy\big) = \sgn(\gk) \go_\ga(t+1),
\]
and similarly,
\[
M^* \go_\gb(t) =\sgn(\gk) \go_\gb(t+1).
\]
Hence, any section of $\sH^1_v(T^2_M)\to S^1$ can be written as
\begin{equation*}
\gf_\ga(t) \go_\ga(t) + \gf_\gb(t) \go_\gb(t),
\end{equation*}
where $\gf_\ga, \gf_\gb\in C^\infty(\R)$ satisfy
\[
\gf_\ga(t+1)=\sgn(\gk) \gf_\ga(t),\quad \gf_\gb(t+1)=\sgn(\gk)
\gf_\gb(t)
\]
We deduce from \eqref{RhoCohomHyp:1} and \eqref{RhoCohomHyp:3}
that the operator $D:=-i\tau_t \pd_t$ on
$C^\infty\big(S^1,\sH^1_v(T^2_M)\big)$ acts in terms of the
coordinate functions $(\gf_\ga,\gf_\gb)$ as
\[
D \begin{pmatrix} \gf_\ga \\ \gf_\gb \end{pmatrix} =
\begin{pmatrix} 0 &\pd_t -\log|\gk| \\ -\pd_t -\log|\gk| &0
\end{pmatrix} \begin{pmatrix} \gf_\ga \\ \gf_\gb \end{pmatrix}
\]
Hence, it becomes clear that if $(\gf_\ga,\gf_\gb)$ defines an
eigenvector of $D$ with eigenvalue $\gl$, then
$(\gf_\ga,-\gf_\gb)$ gives rise to an eigenvector with eigenvalue
$-\gl$. This means that $\spec(D)$ is symmetric, so that
$\eta(D)=0$. Since the operator $D$ is precisely the restriction
of $D_{S^1}\otimes \nabla^{\sH_v}$ to
$C^\infty\big(S^1,\sH^1_v(T^2_M)\big)$, we obtain the desired
result.
\end{proof}

\noindent\textbf{The Logarithm of the Dedekind Eta Function.} The
discussion of the Rho form is more transparent, if we consider the
twisted and the untwisted Eta forms separately. We start with the
untwisted case. This case has already received a far-reaching
treatment in the beautiful article \cite{Ati87}, from which we
borrow the main ideas.

Recall that the classical \emph{Dedekind Eta function} is defined
as
\[
\boldsymbol{\eta}(\gs):= q_\gs^{\frac 1{24}}\prod_{n=1}^\infty
(1-q_\gs^n),\quad \gs\in \H,\quad q_\gs:= e^{2\pi i \gs},
\]
see \cite[Sec. 18.5]{Lang}. As in \cite{Ati87} we use the bold
symbol $\boldsymbol \eta$ to avoid confusion with an Eta function
in the sense of Definition \ref{EtaFctnDef}. Using the power
series expansion
\[
\log(1-z) = -\sum_{m=1}^\infty \frac {z^m}{m},\quad |z|<1,
\]
one can define a logarithm of $\boldsymbol{\eta}(\gs)$ by
\begin{equation}\label{DedLog}
\log \boldsymbol{\eta}(\gs) := \frac{\pi i \gs}{12} -
\sum_{n>0}\sum_{m>0} \frac {q_\gs^{mn}}{n}.
\end{equation}
The sum in \eqref{DedLog} is of the same form as the one in
defining $E_0(\gs)$ in Theorem \ref{F:Thm} (i). Hence, we can
make Remark \ref{E:Rem} (i) more precise and note that
\begin{equation}\label{E:DedRel:1}
\log \boldsymbol{\eta}(\gs) = \frac{\pi i \gs}{12} -\frac{
E_0(\gs)}2= \frac 12\big(\pi i \gs P_2(0) - E_0(\gs)\big),
\end{equation}
where as always, $P_2$ is the second periodic Bernoulli function.
Therefore, we can reformulate Theorem \ref{TorBundleRhoForm} (i)
as
\begin{equation}\label{TorBundleEtaForm:Alt}
\int_0^1\widehat \eta = \frac 2\pi \Im \big[\log
\boldsymbol{\eta}\big(M^{\op}\gs(0)\big) - \log
\boldsymbol{\eta}\big(\gs(0)\big)\big] - \frac 1{2\pi}\int_0^1
\frac{\dot \gs_1(t)}{\gs_2(t)}dt,
\end{equation}
where $\gs(t)=\gs_1(t)+i\gs_2(t)$ is an $M$-invariant path in
$\H$. Hence, the Eta invariant of $T^2_M$ is related to the
transformation property of $\log \boldsymbol{\eta}$ under modular
transformations.

The study of this has a long history, starting with Dedekind's
work \cite{Ded}. There are several different proofs of the
following theorem, see for example \cite{Scz78} for references
and a beautifully simple proof. A short discussion of the
Dedekind sums appearing below is included in Appendix
\ref{CompDedekind}.

\begin{theorem}[Dedekind]\label{LogDedTrans}
Let $\gs\in \H$, and let $M=\left(\begin{smallmatrix} a &b\\ c&d
\end{smallmatrix}\right)\in \SL_2(\Z)$ with $c\neq 0$. Then
\[
\log \boldsymbol{\eta}(M^{\op}\gs) - \log \boldsymbol{\eta}(\gs)
= \frac 12 \log\Big(\frac{c\gs +a}{\sgn(c)i}\Big) + \pi
i\Big(\frac{a+d}{12c} - \sgn(c)s(a,c)\Big),
\]
where the logarithm on the right hand side is the standard branch
on $\C\setminus\R^-$, and $s(a,c)$ is the classical Dedekind sum,
see \eqref{DedDef},
\[
s(a,c) = \sum_{k=1}^{|c|-1} P_1\big(\lfrac
{ak}{c}\big)P_1\big(\lfrac k{c}\big).
\]
\end{theorem}

\begin{remark*}
Note that since we have defined the action of $\SL_2(\Z)$ on $\H$
using the involution $M\mapsto M^{\op}$ as in Lemma
\ref{TeichmAlmComp:Equiv}, we have to interchange $a$ and $d$ in
the classical formula. However, $s(a,c)$ is not affected by this,
see \eqref{DedSymm}.
\end{remark*}

\noindent\textbf{The Untwisted Eta Invariant.} From Theorem
\ref{LogDedTrans} we can deduce the formula for the Eta invariant
of $T^2_M$ for hyperbolic $M$. The formula we shall obtain appears
as a signature cocycle for the mapping class group the formula
already in \cite{Mey73}, and as a signature defect in
\cite{Hir73}. However, its derivation using Theorem
\ref{LogDedTrans} and the adiabatic limit formula as well as an
explanation of the relation among these different invariants are
due to Atiyah \cite{Ati87}.

\begin{theorem}[Atiyah, Hirzebruch, Meyer]\label{EtaHyp}
Let $M=\left(\begin{smallmatrix} a &b\\ c&d
\end{smallmatrix}\right)\in \SL_2(\Z)$ by hyperbolic. Let $g$ be
the metric on $T^2_M$ defined by by $\gs(t)$ as in
\eqref{SigmaHypDef}, and let $B^{\ev}$ be the associated odd
signature operator on $T^2_M$. Then
\[
\eta(B^{\ev}) = \frac {a+d}{3c} - 4\sgn(c) s(a,c) -
\sgn\big(c(a+d)\big).
\]
\end{theorem}

\begin{proof}
Let $g_\eps$ be the adiabatic metric associated to $g$, and denote
by $B_{\eps}^{\ev}$ the corresponding adiabatic family of odd
signature operators. It follows from Proposition
\ref{RhoCohomHyp} and its proof that the Eta invariant of the
bundle of vertical cohomology groups vanishes. We thus deduce from
Theorem \ref{DaiGen:2} that
\[
\lim_{\eps\to 0} \eta(B_{\eps}^{\ev}) = 2\int_0^1 \widehat \eta =
\frac 4\pi \Im \big[\log \boldsymbol{\eta}\big(M^{\op}\gs(0)\big)
- \log \boldsymbol{\eta}\big(\gs(0)\big)\big] - \frac
1{\pi}\int_0^1 \frac{\dot \gs_1(t)}{\gs_2(t)}dt,
\]
where we have used \eqref{TorBundleEtaForm:Alt} for the last
equality. Hence, Theorem \ref{LogDedTrans} implies that
\begin{equation}\label{EtaHyp:0}
\lim_{\eps\to 0} \eta(B_{\eps}^{\ev}) = \frac {a+d}{3c} -
4\sgn(c) s(a,c) +\frac 1\pi \Big[2\Im \log\Big(\frac{c\gs(0)
+a}{\sgn(c)i}\Big)- \int_0^1 \frac{\dot
\gs_1(t)}{\gs_2(t)}dt\Big].
\end{equation}
We now note that it follows from Lemma \ref{SigmaHypLem} that
\[
\gs(t) =\lfrac {\ga+\gb}2 + \lfrac{|\ga-\gb|}2 e^{i\gf(t)},\quad
\gf(t) = \arg\big(\gs(t)- \lfrac {\ga+\gb}2\big).
\]
Here, the argument function is such that for $z\in
\C\setminus\R^-$, one has $\arg(z)\in (-\pi,\pi)$. We obtain
\[
\int_0^1 \frac{\dot \gs_1(t)}{\gs_2(t)}dt = - \arg\big(\gs(1)-
\lfrac {\ga+\gb}2\big) + \arg\big(\gs(0)- \lfrac {\ga+\gb}2\big).
\]
Using the explicit formula \eqref{SigmaHypDef} for $\gs(t)$, one
finds that $\arg\big(\gs(0)- \frac {\ga+\gb}2\big)=\lfrac \pi 2$,
and
\[
\gs(1)- \frac {\ga+\gb}2 = \frac 12
\frac{(\ga-\gb)(\gk^2-\gk^{-2})}{\gk^2+\gk^{-2}} + i\frac
{|\ga-\gb|}{\gk^2+\gk^{-2}}.
\]
Hence, using the abbreviations
\begin{equation}\label{xyDef}
x:= \lfrac 12\sgn(c) (\gk+\gk^{-1}),\quad y:= \lfrac
12(\gk-\gk^{-1}),
\end{equation}
we find that
\[
\arg\big(\gs(1)- \lfrac {\ga+\gb}2\big) = \arg\big(2xy + i\big).
\]
Note that $y^2 = x^2- 1$, so that $2xy + i = i(x-iy)^2$.
Moreover, $y>0$ and $y<|x|$ so that
\[
\arg(x+iy)\in \begin{cases}\,\,\, \big(0,\frac \pi 4\big),
&\text{if $x>0$},\\ \big(\frac {3\pi}4, 2\pi\big), &\text{if
$x<0$}.
\end{cases}
\]
Thus, we obtain
\begin{equation}\label{EtaHyp:1}
\int_0^1 \frac{\dot \gs_1(t)}{\gs_2(t)}dt =
\arg\big(-i(x+iy)^2\big) + \lfrac \pi 2 = 2\arg(x+iy) -
\begin{cases}\,\, 0 , &\text{if $x>0$},\\
2\pi , &\text{if $x<0$}.\end{cases}
\end{equation}
On the other hand, it follows from the definition of $\ga$ and
$\gb$ that
\[
c\gs(0) +a = c\big(\lfrac{\ga+\gb}2 + i \lfrac{|\ga-\gb|}2\big) +
a = \lfrac 12 \big(\gk+\gk^{-1} +i\sgn(c)(\gk-\gk^{-1})\big) =
\sgn(c)(x+iy),
\]
with $x$ and $y$ as in \eqref{xyDef}. Since we are using the
standard branch of the logarithm, we have
\begin{equation*}
\Im \log\Big(\frac{c\gs(0) +a}{\sgn(c)i}\Big) =
\arg\big(-i(x+iy)\big) = \arg(x+iy) - \lfrac \pi 2.
\end{equation*}
Combining this with \eqref{EtaHyp:1} we find that
\[
2 \Im \log\Big(\frac{c\gs(0) +a}{\sgn(c)i}\Big) - \int_0^1
\frac{\dot \gs_1(t)}{\gs_2(t)}dt  = -\sgn(x) \pi.
\]
As $\gk+\gk^{-1} = a+d$ we have $\sgn(x) =\sgn\big(c(a+d)\big)$
so that using \eqref{EtaHyp:0}, we finally arrive at
\[
\lim_{\eps\to 0} \eta(B_{\eps}^{\ev}) = \frac {a+d}{3c} -
4\sgn(c) s(a,c) - \sgn\big(c(a+d)\big).
\]
Hence, it remains to argue that in the case at hand,
\[
\eta(B^{\ev}) = \lim_{\eps\to 0} \eta(B_{\eps}^{\ev}),
\]
This is precisely \cite[Lem. 5.56]{Ati87} and we will not repeat
the argument here.
\end{proof}

\begin{remark*}
The proof of Theorem \ref{EtaHyp} in \cite{Ati87} is along
different lines than our discussion. In \cite[Thm. 5.60]{Ati87},
the Eta invariant of $T^2_M$ is seen to be equal to a large
number of quantities, including a signature defect. Then in
\cite[Sec. 6]{Ati87}, the formula for $\eta(B^{\ev})$ is obtained
by explicitly constructing a bounding manifold and a computation
of the signature defect. In particular, the transformation formula
of the Dedekind Eta function in Theorem \ref{LogDedTrans} is not
used. However, since our focus is the application of the
adiabatic limit formula to compute Eta respectively Rho
invariants, we have to use Theorem \ref{LogDedTrans} in some form.
\end{remark*}

\noindent\textbf{The Generalized Dedekind Eta Function.} To obtain
the formula for $\U(1)$-Rho invariants of $T^2_M$ in the spirit
of the discussion of the untwisted case, we now need a twisted
version of $\boldsymbol{\eta}(\gs)$ and a transformation formula
for its logarithm. Fortunately, a corresponding treatment can be
found in \cite{Die59}. As in \cite[p. 38]{Die59} we make the
following

\begin{dfn}\label{GenDedEta}
For $g,h\in \R$ and $\gs\in \H$ let $z:= g\gs + h$, $q_z =
e^{2\pi i z}$ and $q_\gs = e^{2\pi i\gs}$. Define
\[
\boldsymbol{\eta}_{g,h}(\gs):= \xi(g,h)\, q_\gs^{\frac{P_2(g)}2}
(1-q_z) \prod_{n=1}^{\infty} (1-q_zq_\gs^m)(1-q_z^{-1}q_\gs^m),
\]
where
\[
\xi(g,h) := \begin{cases} e^{2\pi i (g-\frac 12)P_1(h)}, &\text{if
$g\in \Z$}\\  \,\,\,e^{2\pi i [g]P_1(h)}, &\text{if $g\notin \Z$}.
\end{cases}
\]
As always, $P_1$ is the first periodic Bernoulli function, and
$[g]$ is the largest integer less or equal than $g$.
\end{dfn}

\begin{remark}\label{GenDedEtaRem}\quad \nopagebreak
\begin{enumerate}
\item Since $\gs\in \H$, the term $q_\gs^m$ decays exponentially
with $m$. This implies that $\boldsymbol{\eta}_{g,h}(\gs)$ is
well-defined.
\item Definition \ref{GenDedEta} might look slightly different than
the formula in \cite{Die59}. Yet, writing $g$ and $h$ as $\Tilde
g/f$ and $\Tilde h/f$ with integers $\Tilde g,$ $\Tilde h$ and
$f$, the function $\boldsymbol{\eta}_{g,h}(\gs)$ is easily seen
to be equal to what is denoted $\eta_{\Tilde g,\Tilde h}(\gs)$ in
loc.cit.
\item The reason for the factor $\xi(g,h)$ is to achieve that
$\boldsymbol{\eta}_{g,h}(\gs)$ depends on $g$ and $h$ only modulo
$\Z$, see \cite[p. 39]{Die59}.
\item The Dedekind Eta function $\boldsymbol{\eta}(\gs)$ is not
equal to $\boldsymbol{\eta}_{0,0}(\gs)$, since the latter
obviously vanishes. Dropping the factor $1-q_z$ from the
definition of $\boldsymbol{\eta}_{g,h}(\gs)$ one would get a
direct generalization of $\boldsymbol{\eta}(\gs)^2$. However,
Definition \ref{GenDedEta} allows us to use the results of
\cite{Die59} without too many changes.
\end{enumerate}
\end{remark}

One defines $\log \boldsymbol{\eta}_{g,h}(\gs)$ in analogy to
\eqref{DedLog}, see \cite[p. 40]{Die59}.

\begin{dfn}\label{GenDedLog}
Let $(g,h)\in \Q^2\setminus \Z^2$. If $0\le g<1$, we define
\begin{equation*}
\log \boldsymbol{\eta}_{g,h}(\gs):= \pi
i\big(\gf(g,h)+P_2(g)\big) - \sum_{n>0}\lfrac 1{n} q_z^{n} -
\sum_{m>0}\sum_{n>0}\lfrac 1{n} (q_z+q_z^{-1})^{n}q_\gs^{mn},
\end{equation*}
where
\[
\gf(g,h) = \begin{cases} -P_1(h), &\text{if $g=0$},\\
\hphantom{-P}0, &\text{if $g\neq 0$.}
\end{cases}
\]
For general $g$ we define
\[
\log \boldsymbol{\eta}_{g,h}(\gs):= \log
\boldsymbol{\eta}_{g-[g],h}(\gs).
\]
\end{dfn}

The transformation formula of $\log \boldsymbol{\eta}_{g,h}(\gs)$
is then given by \cite[Thm. 1]{Die59},

\begin{theorem}[Dieter]\label{GenDedTrans}
Let $M=\left(\begin{smallmatrix} a& b\\ c&
d\end{smallmatrix}\right)\in \SL_2(\Z)$ with $c\neq 0$, let
$(g,h)\in \Q^2\setminus \Z^2$, and define
\[
\begin{pmatrix} g'\\ h'\end{pmatrix}:= \begin{pmatrix} a &-c \\ -b
& d\end{pmatrix} \begin{pmatrix} g\\ h\end{pmatrix}.
\]
Then for all $\gs\in \H$
\[
\log \boldsymbol{\eta}_{g',h'}(M^{\op}\gs) - \log
\boldsymbol{\eta}_{g,h}(\gs) = \pi i \Big(\frac ac P_2(g) + \frac
dc P_2(g') - 2\sgn(c) s_{g',h'}(d,c) \Big),
\]
where $s_{g',h'}(d,c)$ is the generalized Dedekind sum, see
Definition \ref{GenDedDef},
\[
s_{g',h'}(d,c) = \sum_{k=0}^{|c|-1} P_1\big(d\lfrac
{k+g'}{c}+h'\big)P_1\big(\lfrac {k+g'}{c}\big).
\]
\end{theorem}

\begin{remark*}\quad\nopagebreak
\begin{enumerate}
\item As for the transformation formula for
$\log\boldsymbol{\eta}(\gs)$, we have formulated Theorem
\ref{GenDedTrans} in terms of $M^{\op}$ acting on $\H$, which
means that $a$ and $d$ have been interchanged in comparison to
\cite[Thm. 1]{Die59}.
\item The generalized Dedekind sums appeared first in
\cite{Mey57}. A brief discussion of the aspects we need is
contained in \ref{CompDedekind}.
\item We want to point out that the proof of Theorem
\ref{GenDedTrans} in \cite{Die59} is rather involved. The simple
proof in \cite{Scz78} of the transformation formula for
$\log\boldsymbol{\eta}(\gs)$ carries over with minor changes in
the case that $g\in \Z$. It would be interesting to know if there
is a proof for the general case of Theorem \ref{GenDedTrans}
along the lines of \cite{Scz78}.
\end{enumerate}
\end{remark*}

\noindent\textbf{Application to the Rho Form.} As in the
untwisted case, the structure of the formula in Definition
\ref{GenDedLog} resembles what we have encountered in Theorem
\ref{TorBundleRhoForm} (ii). In fact,

\begin{lemma}\label{E:GenDedRel}
With the notation of Theorem \ref{F:Thm} and Remark \ref{E:Ext},
we have for all $\nu \in \Q^2\setminus\Z$.
\[
\Im \big(\log \boldsymbol{\eta}_{\nu_1,-\nu_2}(\gs)\big) = \Im
\big(\pi i \gs P_2(\nu_1) - E_\nu(\gs)\big).
\]
\end{lemma}

\begin{proof}
Both sides of the equation are defined in terms of
$\nu_1-[\nu_1]$ and are $\Z$-periodic in $\nu_2$ Hence, we can
assume that $\nu\in [0,1)^2$. Then, if $\nu_1\neq 0$, the
relation is immediate---and clearly holds for the real parts as
well. If $\nu_1=0$, one observes that
\begin{equation}\label{E:DedRel:2:Rem}
\Im \Big(\pi iP_1(\nu_2)  - \sum_{n>0}\lfrac 1{n} e^{-2\pi i
\nu_2}\Big) = - \Im \Big( \sum_{n>0}\lfrac 1{n} \cos(2\pi
\nu_2)\Big) = 0,
\end{equation}
where the first equality follows from the Fourier series expansion
\[
P_1(\nu_2)= \nu_2-\lfrac12 = -\frac 1{2\pi i} \sum_{n>0}\lfrac
1{n} \big(e^{2\pi i \nu_2} - e^{-2\pi i \nu_2}\big),\quad
\nu_2\notin \Z,
\]
whose proof is a standard exercise. Then \eqref{E:DedRel:2:Rem}
implies that the imaginary parts of the extra terms in Definition
\ref{GenDedLog} cancel each other out so that the result is
indeed the right hand side of the formula in Lemma
\ref{E:GenDedRel}.
\end{proof}

\begin{prop}\label{RhoFormHyp}
Let $M=\left(\begin{smallmatrix} a& b\\ c&
d\end{smallmatrix}\right)\in \SL_2(\Z)$ be hyperbolic, and
$\nu\in \R^2\setminus \Z^2$ with $(\Id -M^t)\nu \in \Z^2$. Let $A$
be the corresponding flat $\U(1)$-connection over the mapping
torus $T^2_M$, and use \eqref{SigmaHypDef} to define a metric.
Then
\[
\int_0^1\widehat \eta_A = \lfrac {a+d}c P_2(\nu_1) - 2\sgn(c)
s_{\nu_1,\nu_2}(a,c).
\]
\end{prop}

\begin{proof}
Certainly, we want to use Lemma \ref{E:GenDedRel} to apply Theorem
\ref{GenDedTrans} to part (ii) of Theorem \ref{TorBundleRhoForm}.
We first note that $(\Id -M^t)\nu \in \Z^2$ implies that $\nu\in
\Q^2$, so that we are precisely in the situation of Lemma
\ref{E:GenDedRel}. Abbreviate $\nu':= M^t\nu$, so that by
assumption $\nu-\nu'\in \Z$. Thus, according to Definition
\ref{GenDedLog},
\[
\log\boldsymbol{\eta}_{\nu_1',-\nu_2'}(\gs) =
\log\boldsymbol{\eta}_{\nu_1,-\nu_2}(\gs),\quad \text{for all $
\gs\in \H$}.
\]
Moreover, $P_1(\nu_1') = P_1(\nu_1)$, and
\[
\begin{pmatrix} \nu_1'\\ -\nu_2'\end{pmatrix}= \begin{pmatrix} a &-c \\ -b
& d\end{pmatrix} \begin{pmatrix} \nu_1\\ -\nu_2\end{pmatrix}.
\]
Hence, Lemma \ref{E:GenDedRel}, Theorem \ref{GenDedTrans} and
Theorem \ref{TorBundleRhoForm} (ii) imply that
\[
\int_0^1\widehat \eta_A = \lfrac {a+d}c P_2(\nu_1) - 2\sgn(c)
s_{\nu_1',-\nu_2'}(d,c).
\]
Now,
\[
\begin{split}
s_{\nu_1',-\nu_2'}(d,c) &= \sum_{k=0}^{|c|-1} P_1\big(d \lfrac
{k+\nu_1'}{c}-\nu_2'\big)P_1\big(\lfrac {k+ \nu_1'}{c}\big)  =
\sum_{k=0}^{|c|-1} P_1\big(\lfrac
{dk+d\nu_1'-c\nu_2'}{c}\big)P_1\big(\lfrac {k+
a\nu_1+c\nu_2}{c}\big)\\
&= \sum_{k=0}^{|c|-1} P_1\big(\lfrac {k+\nu_1}{c}\big)P_1\big(a
\lfrac {k+ \nu_1}{c} +\nu_2\big) = s_{\nu_1,\nu_2}(a,c),
\end{split}
\]
where we have rewritten $\nu'$ in terms of $\nu$, and then used
that $\setdef{ak}{k=0,\ldots,|c|-1}$ is a representation system
of $\Z$ modulo $c$, see Appendix \ref{CompDedekind} for more
details. This implies the desired result.
\end{proof}

\noindent\textbf{Rho Invariants of Hyperbolic Mapping Tori.}
After this preparation, we finally arrive at the main result of
this section.

\begin{theorem}\label{RhoHyp}
Let $M=\left(\begin{smallmatrix} a& b\\ c&
d\end{smallmatrix}\right)\in \SL_2(\Z)$ be hyperbolic, and let
$\nu\in \R^2\setminus \Z^2$ satisfy
\[
\begin{pmatrix} m_1 \\ m_2 \end{pmatrix} = (\Id -M^t)
\begin{pmatrix}
\nu_1\\ \nu_2 \end{pmatrix}\in \Z^2.
\]
Let $A$ be the corresponding flat $\U(1)$-connection over the
mapping torus $T^2_M$, and define $r\in \{0,\ldots |c|-1\}$ by
requiring that $m_1\equiv r\, (c)$. Then
\[
\begin{split}
\rho_A(T^2_M) & = \lfrac {2(a+d)- 4}c \big( P_2(\nu_1) -\lfrac
16\big) - 4\sum_{k=1}^{|c|-r} P_1\big(\lfrac{dk}{c}\big)  +
\sgn\big(c(a+d)\big) -
\sgn(c)\gd(\nu_1)\big(1-\gd(\lfrac{m_1}c)\big)
\\  & \qquad\quad  -2 P_1\big(\lfrac
{dm_1}{c}\big) -2\gd(\nu_1)\Big( P_1\big(\lfrac{m_1}{c}\big) -
P_1\big(\lfrac{dm_1}{c}\big)\Big),
\end{split}
\]
where $\gd$ is the characteristic function of $\R\setminus \Z$.
\end{theorem}

\begin{proof}
According to Proposition \ref{RhoCohomHyp}, the Rho invariant of
the bundle of vertical cohomology groups vanishes. Also, since
the base is 1-dimensional, Dai's correction term is zero. Hence,
we can use the general formula for Rho invariants in Theorem
\ref{RhoGen}, together with the formul{\ae} of Theorem
\ref{EtaHyp} and Proposition \ref{RhoFormHyp}, to deduce that
\begin{equation}\label{RhoHypPrep}
\rho_A(T^2_M) =  \lfrac {2(a+d)}c \big( P_2(\nu_1) -\lfrac
16\big) - 4\sgn(c)\big(s_{\nu_1,\nu_2}(a,c) - s(a,c)\big)   +
\sgn\big(c(a+d)\big).
\end{equation}
A formula for the difference of $s_{\nu_1,\nu_2}(a,c)$ and
$s(a,c)$ is given in Proposition \ref{GenDedRel}. With $r\in
\{0,\ldots |c|-1\}$ such that $m_1\equiv r\, (c)$, we have
\[
\begin{split}
s_{\nu_1,\nu_2}(a,c) - s(a,c) = \lfrac 1{|c|} \big(P_2(\nu_1)
&-\lfrac 16 \big) + \sum_{k=1}^{|c|-r}
P_1\big(\lfrac{dk}{|c|}\big) + \lfrac 12 P_1\big(\lfrac
{dm_1}{|c|}\big)
\\  & +\lfrac 12 \gd(\nu_1)\Big(
P_1\big(\lfrac{m_1}{|c|}\big) -
P_1\big(\lfrac{dm_1}{|c|}\big)\Big) + \lfrac
1{4}\gd(\nu_1)\big(1-\gd(\lfrac{m_1}c)\big).
\end{split}
\]
We now insert this into \eqref{RhoHypPrep}. Since $P_1$ is odd,
the factor $\sgn(c)$ in front of $s_{\nu_1,\nu_2}(a,c) - s(a,c)$
cancels the norms in the denominators. Then we arrive at the
formula of Theorem \ref{RhoHyp}.
\end{proof}

\noindent\textbf{Immediate Applications.} The main formula in
Theorem \ref{RhoHyp} might look more complicated than the
intermediate formula \eqref{RhoHypPrep}. Yet, it is more
satisfactory from a computational point of view, since the sum
$\sum_{k=1}^{|c|-r} P_1\big(\lfrac{dk}{c}\big)$ is much easier to
compute than the individual Dedekind sums. For concreteness, let
us use Theorem \ref{RhoHyp} for some explicit computations.

\begin{example*}\quad\nopagebreak
\begin{enumerate}
\item Consider
\[
M=\begin{pmatrix} -2 &1 \\ 1&-1 \end{pmatrix}, \quad \text{so
that}\quad \Id - M^t = \begin{pmatrix} 3 &-1 \\ -1&2
\end{pmatrix}.
\]
Since $\det(\Id - M^t) = 5$, a pair $\nu=(\nu_1,\nu_2)\in \R^2$
with $m=(\Id - M^t) \nu\in \Z^2$ has to consist of rational
numbers with denominator 5. Recall that we exclude the case
$\nu\in \Z^2$ and may restrict to $\nu\in [0,1)^2$. One then
verifies that to obtain a full set of representatives for the flat
connections on $T^2_M$ we are interested in, we need to consider
pairs $\nu$ and $m$ with
\[
\begin{matrix}
\nu&= &(\lfrac 15,\lfrac 35)  &(\lfrac 25,\lfrac 15)
 &(\lfrac 35,\lfrac 45) &(\lfrac 45,\lfrac 25) \\
m&= &(0,1) &(1,0) &(1,1) &(2,0).
\end{matrix}
\]
As $c=1$ and $\nu_1\neq 0$, the formula of Theorem \ref{RhoHyp}
reduces to
\[
\rho_A(T^2_M) = 2\big((a+d)- 2\big) ( \nu_1^2 -\nu_1 ) -
\sgn\big(c(a+d)\big) - \sgn(c) = -10(\nu_1^2 -\nu_1 ).
\]
Hence, one computes
\[
\begin{matrix}
\nu_1 &= &\lfrac 15   &\lfrac 25 &\lfrac 35 &\lfrac 45, \\ \\
\rho_A(T^2_M) &= &\lfrac 85 &\lfrac {12}5 &\lfrac {12}5 &\lfrac
85.
\end{matrix}
\]
\item As a further example, let us consider
\[
M=\begin{pmatrix} 3 &2 \\ 4&3 \end{pmatrix}, \quad \text{so
that}\quad \Id - M^t = \begin{pmatrix} -2 &-4 \\ -2&-2
\end{pmatrix}.
\]
Now, one easily verifies that we can represent the conjugacy
classes of flat connections of interest by
\[
\begin{matrix}
\nu&= &(0,\lfrac 12)  &(\lfrac 12,0)
 &(\lfrac 12,\lfrac 12), \\
m&= &(-2,-1) &(-1,-1) &(-3,-2),\\
r&= &2 &3 &1.
\end{matrix}
\]
Then
\[
\rho_A(T^2_M) = 2( \nu_1^2 -\nu_1 ) - 4\sum_{k=1}^{4-r}
P_1\big(\lfrac{3k}{4}\big)  + 1 -2 P_1\big(\lfrac {3r}{4}\big)
-2\gd(\nu_1)\Big( P_1\big(\lfrac{r}{4}\big) -
P_1\big(\lfrac{3r}{4}\big)\Big).
\]
For $\nu=(0,\frac 12)$, we have $r=2$, and so
\[
\rho_A(T^2_M) = - 4\big(P_1(\lfrac 34)  + P_1(\lfrac 12)\big) + 1
-2P_1(\lfrac 12) = 0,
\]
for $\nu=(\frac 12,0)$ with $r=3$,
\[
\rho_A(T^2_M) = 2\big(\lfrac 14 -\lfrac 12\big) - 4P_1(\lfrac34)
+1 - 2P_1(\lfrac34) = 0,
\]
and lastly, for $\nu=(\frac 12,\lfrac 12)$ with $r=1$,
\[
\rho_A(T^2_M) = 2\big(\lfrac 14 -\lfrac 12\big) + 1 - 2P_1(\lfrac
14) = 1,
\]
where we have used that $\sum_{k=1}^{3}
P_1\big(\lfrac{3k}{4}\big)=0$, see \eqref{SumFullRepSys}.
\end{enumerate}
\end{example*}

Recall from Corollary \ref{RhoConnVar} (ii) that the non-integer
part of the Rho invariant on a 3-dimensional manifold is
essentially the Chern-Simons invariant associated to the Chern
character. More precisely,
\begin{equation}\label{CSRhoRelRecall}
\rho_A(T^2_M) \equiv 4\CS(A) \mod\Z.
\end{equation}
In the case of torus bundles over surfaces, computations for
Chern-Simons invariants are contained in \cite{FreVaf, Jef92,
KK90}. For the case of $\U(1)$-connections, see for example
\cite[Thm. 7.22]{FreVaf}.

\begin{cor}\label{CSHyp}
Under the assumptions of Theorem \ref{RhoHyp} we have
\[
\rho_A(T^2_M) \equiv 2\big(\nu_2m_1 - \nu_1m_2\big)\mod \Z.
\]
\end{cor}

\begin{proof}
First note that for all $k\in \Z$
\[
2P_1\big(\lfrac kc \big) \equiv 2\lfrac kc \mod\Z.
\]
In particular,
\[
4\sum_{k=1}^{|c|-r} P_1\big(\lfrac{dk}{c}\big) \equiv 4\lfrac dc
\lfrac{(|c|-r)(|c|-r+1)}2 \equiv 2\big(\lfrac {dm_1^2}c - \lfrac
{dm_1}c \big)\mod \Z.
\]
Here, we have used that by definition $r\equiv m_1\,(c)$. We also
note that $(\Id -M^t)\nu = m$ means explicitly that
\begin{equation}\label{nu:m:expl}
\begin{pmatrix} m_1\\ m_2 \end{pmatrix} = \begin{pmatrix}
1-a & -c \\ -b &1-d \end{pmatrix} \begin{pmatrix} \nu_1\\ \nu_2
\end{pmatrix},
\quad \begin{pmatrix} \nu_1\\ \nu_2 \end{pmatrix} = \lfrac 1{a+d
-2} \begin{pmatrix} d-1  & -c \\ -b &a -1 \end{pmatrix}
\begin{pmatrix} m_1\\ m_2
\end{pmatrix}
\end{equation}
Let us assume now that $\nu_1\in \Z$. Then
\[
\rho_A(T^2_M) \equiv   - 2\big(\lfrac {dm_1^2}c - \lfrac {dm_1}c
\big) - 2 \lfrac {dm_1}{c} \equiv - 2 \lfrac {dm_1^2}c \mod\Z.
\]
It follows from \eqref{nu:m:expl} that $\frac{m_1}c \equiv
-\nu_2$ modulo $\Z$, and $dm_1\equiv m_1$ modulo $\Z$. Therefore,
\[
\rho_A(T^2_M) \equiv 2 \nu_2m_1 \mod\Z,
\]
which is the claim of Corollary \ref{CSHyp} in the case that
$\nu_1\in \Z$. If $\nu\notin\Z$, we have
\begin{equation}\label{CSHyp:1}
\rho_A(T^2_M) \equiv \lfrac {2(a+d)- 4}c ( \nu_1^2 -\nu_1) -
2\big(\lfrac {dm_1^2}c - \lfrac {dm_1}c \big) -2 \lfrac{m_1}{c}
\mod \Z.
\end{equation}
From \eqref{nu:m:expl} we know that
\[
\lfrac {(a+d)- 2}c \nu_1 = \lfrac{(d-1)m_1}c - m_2.
\]
Inserting this into \eqref{CSHyp:1}, one finds that
\[
\begin{split}
\rho_A(T^2_M) & \equiv 2\Big( \big(\lfrac{(d-1)m_1}c -
m_2\big)(\nu_1-1)
+ \lfrac{(d-1)m_1}c - \lfrac {dm_1^2}c \Big) \mod \Z \\
& \equiv 2\big( -\nu_1 m_2 +\lfrac{dm_1}c (\nu_1 -m_1) - \lfrac
{m_1}c \nu_1\big) \mod \Z\\
& \equiv 2\big( -\nu_1 m_2 + dm_1\nu_2 +\lfrac{(ad-1)m_1}c
\nu_1\big) \mod \Z,
\end{split}
\]
where in the last line we have used that $\nu_1 -m_1 = a\nu_1 +
c\nu_1$, see \eqref{nu:m:expl}. Now using that $ad-1=bc$ and
observing that $b\nu_1 +d\nu_2\equiv \nu_2$ modulo $\Z$, we
arrive at
\[
\rho_A(T^2_M) \equiv 2\big(-\nu_1 m_2 + m_1\nu_2\big)\mod\Z.
\qedhere
\]
\end{proof}

\begin{remark*}
The formula of Corollary \ref{CSHyp} also holds in the
parabolic case: Let $\eps \left(\begin{smallmatrix} 1& l\\ 0&
1\end{smallmatrix}\right)$ with $l\in \Z$ and $\eps=\pm 1$, and
let $\nu=(\nu_1,\nu_2)\in \R^2\setminus \Z^2$ satisfy
$m=(\Id-M^t)\nu \in \Z^2$. Then, if $\eps=1$,
\[
-l\nu_1 =m_2\in \Z,\quad m_1=0,\quad\text{so that}\quad
2\big(-\nu_1 m_2 + m_1\nu_2\big) = 2l\nu_1^2.
\]
According to Theorem \ref{RhoParab}, this is congruent to
$\rho_A(T_M^2)$ modulo $\Z$. If $\eps = -1$, then
\[
2\nu_1 = m_1\in \Z,\quad -l\nu_1 +2\nu_2 = m_2\in \Z,
\]
so that again
\[
2\big(-\nu_1 m_2 + m_1\nu_2\big) = 2l\nu_1^2 -2 \nu_1\nu_2 + 2
\nu_1\nu_2 = 2l\nu_1^2.
\]
\end{remark*}

\noindent\textbf{Jeffrey's Conjecture.} We end the main
discussion of this thesis with a remark concerning a possible
perspective for further research. According to
\eqref{CSRhoRelRecall}, Corollary \ref{CSHyp} identifies the
Chern-Simons invariant only modulo $\lfrac 14 \Z$, which might
seem a bit disappointing. Moreover, the methods of \cite{FreVaf,
Jef92, KK90} to obtain the formula for the Chern-Simons invariant
are much less involved than what we have presented. However, this
is precisely the real strength of Theorem \ref{RhoHyp}. It can be
used to compute the difference $\rho_A(T^2_M) - 4\CS(A) \in \Z$.

Recall from Corollary \ref{RhoConnVar} (i) that this is
essentially the spectral flow of the odd signature operator
between the trivial connection and $A$. For this reason it is
promising that a generalization of Theorem \ref{RhoHyp} to higher
gauge groups might be a way to prove Jeffrey's conjecture about
the mod 4 reduction of this spectral flow term, see \cite[Conj.
5.8]{Jef92}.

\cleardoublepage

\appendix

\addtocontents{toc}{\contentsline{chapter}{Appendix}{}{}}
\renewcommand{\chaptermark}[1]%
   {\markboth{\textsc{Appendix \thechapter.\ #1}}{}}
\renewcommand{\sectionmark}[1]%
   {\markright{\textsc{\thesection.\ #1}}}

\renewcommand{\thesection}{\Alph{chapter}.\arabic{section}}
\renewcommand{\thesubsection}{\Alph{chapter}.\arabic{section}.\arabic{subsection}}
\renewcommand{\thetheorem}{\Alph{chapter}.\arabic{section}.\arabic{theorem}}
\renewcommand{\theequation}{\Alph{chapter}.\arabic{equation}}

\chapter{Characteristic Classes and Chern-Simons Forms}\label{CharClass}

Although we assume that the reader is familiar with the theory of
characteristic classes, we include a short survey of Chern-Weil
theory in the way we will use it. We closely follow \cite[Sec.
1.5]{BGV} and \cite[Ch. 1]{Z} to which we also refer for more
details. We place some emphasis on transgression forms and
formulate the results about Chern-Simons invariants, which we use
in Section \ref{RhoLocalInd}.

\section{Chern-Weil Theory}

\subsection{Connections and Characteristic Forms}

We start with a short algebraic preliminary. Let $V$ be a complex
vector space. For $m\in\N$ consider $(\gL^{\ev} \C^m)\otimes V$
as module over the commutative algebra $\gL^{\ev} \C^m$. Then any
element $T\in \gL^{\ev} \C^m\otimes \End(V)$ may be viewed as a
module endomorphism. Upon choosing a basis for $V$, this is a
matrix with entries in $\gL^{\ev} \C^m$. In this way we can
define expressions like $T^n$ and $\det T$. We extend the trace
$\tr_V:\End(V)\to\C$ on $V$ in the natural way to a trace
\[
\tr_V: (\gL^{\ev}\C^m)\otimes \End(V)\to \gL^{\ev}\C^m.
\]
Let $f(z) = \sum_{n\ge 0}a_n z^n$ be a formal power series with
coefficients $a_n$ in $\C$, and assume that $T\in
(\gL^{2\bullet+2} \C^m)\otimes \End(V)$. An endomorphism $T$ of
this form is nilpotent, so that we can define
\[
f(T) = \sum_{n\ge 0}a_n T^n \in (\gL^{\ev} \C^m)\otimes \End(V).
\]
The following algebraic result is the main tool we use for
defining the characteristic forms we need. We skip the easy proof.

\begin{lemma}\label{DetTrLog}
For every $T\in (\gL^{2\bullet+2} \C^m)\otimes \End(V)$,
\begin{equation*}
\det\big(1 + T\big) = \exp\big(\tr_V\big[\log(1+T)\big]\big),
\end{equation*}
where the $\exp(.)$ is taken in the algebra $\gL^{\ev}\C^m$ and
the logarithm is defined using the formal power series
\[
\log(1+z) = \sum_{n\ge 0} \frac{(-1)^n}{n+1} z^{n+1}.
\]
\end{lemma}

It follows from Lemma \ref{DetTrLog} that if $f(z)=1+\sum_{n\ge 1}
a_n z^n$ is a normalized formal power series, then for all $T\in
(\gL^{2\bullet+2} \C^m)\otimes \End(V)$
\[
\det\big(f(T)\big) = \exp\big(\tr_V\big[\log f(T)\big]\big).
\]
Motivated by this we also define
\begin{equation}\label{SqrtDet}
{\det}^{1/2}\big(f(T)\big) := \exp\big(\lfrac 12\tr_V\big[\log
f(T) \big]\big).
\end{equation}

\noindent\textbf{Characteristic Forms of Complex Vector Bundles.}
Now let $E\to M$ be a complex vector bundle over an
$m$-dimensional manifold $M$. Let $\nabla$ be a connection on $E$
with curvature $F_\nabla\in \gO^2\big(M,\End(E)\big)$. In this
context $\gO^{\ev}\big(M,\End(E)\big)$ plays the role of
$\gL^{\ev}\C^m\otimes\End(V)$ in the above considerations.

If $T\in C^\infty\big(M,\End(E)\big)$, then the commutator
$[\nabla,T]$ is an element of $\gO^1\big(M,\End(E)\big)$. We can
extend this to a derivation
\begin{equation*}\label{GradComm}
[\nabla,\cdot]: \gO^\bullet\big(M,\End(E)\big)\to
\gO^{\bullet+1}\big(M,\End(E)\big),
\end{equation*}
by requiring that for $\ga\in \gO^\bullet(M)$ of pure degree
$|\ga|$, and $T\in \gO^\bullet\big(M,\End(E)\big)$,
\[
[\nabla,\ga\wedge T] = (d\ga)\wedge T + (-1)^{|\ga|}\ga\wedge
[\nabla,T].
\]
Then it is easy to check that for every such $T$,
\begin{equation}\label{dTr}
\tr_E [\nabla,T] = d\big(\tr_E T\big).
\end{equation}
Let $f(z) = \sum_{n\ge 0}a_n z^n$ be a formal power series. We
define
\begin{equation*}\label{CurvPowerSeries}
f(\nabla):= \sum_{n\ge 0}a_n \big(\lfrac
{i}{2\pi}F_\nabla\big)^n\in \gO^{\ev}\big(M,\End(E)\big).
\end{equation*}
From \eqref{dTr} and the fact that $[\nabla,F_\nabla]=0$, we
obtain
\begin{equation}\label{CharFormClosed}
d\tr_E\big(f(\nabla) \big) = \tr_E \big[\nabla,f(\nabla)\big]=0.
\end{equation}

\begin{dfn}\label{CharFormDef}
Let $E$ be a complex vector bundle over $M$ with connection
$\nabla$, and let $f(z) = \sum_{n\ge 0}a_n z^n$ be a formal power
series. Then we define the \emph{characteristic form} of $\nabla$
associated to $f$ by
\begin{equation*}
\tr_E\big[f(\nabla)\big]=\tr_E \big[\sum_{n\ge 0}a_n \big(\lfrac
{i}{2\pi}F_\nabla\big)^n\big] \in \gO^{\ev}(M).
\end{equation*}
\end{dfn}

\begin{dfn}\label{ChernCharDef}
Let $E$ be a complex vector bundle over $M$ with connection
$\nabla$.
\begin{enumerate}
\item The characteristic form associated to $\exp(z)$ is called
the \emph{Chern character form}
\begin{equation*}
\ch(E,\nabla):= \tr_E\big[ \exp\big(\lfrac
{i}{2\pi}F_\nabla\big)\big] \in \gO^{\ev}(M).
\end{equation*}
\item The characteristic form
\[
c(E,\nabla):= \det\big(1+\lfrac {i}{2\pi}F_\nabla \big)=
\exp\big(\tr_E\big[\log(1+\lfrac {i}{2\pi}F_\nabla)\big]\big) \in
\gO^{\ev}(M)
\]
is called the {\em total Chern form}.
\item The \emph{$j$-th Chern form}
\[
c_j(E,\nabla)\in \gO^{2j}(M)
\]
is defined as the component of degree $2j$ of the total Chern
form, i.e.,
\begin{equation*}\label{ChernFormDef}
c(E,\nabla)= \sum_{j=0}^{[m/2]} c_j(E,\nabla) = 1+
c_1(E,\nabla)+c_2(E,\nabla)+\ldots.
\end{equation*}
\end{enumerate}
\end{dfn}

\begin{remark*}\quad\nopagebreak
\begin{enumerate}
\item Note that it follows form Lemma \ref{DetTrLog} that the total
Chern form fits into the framework of Definition
\ref{CharFormDef} if we take $f(z) = \log(1+z)$ and exponentiate
in $\gO^{\ev}(M)$ after taking $\tr_E\big[f(\nabla)\big]$. Since
\[
d\big(\exp\circ \tr_E\big[f(\nabla)\big]\big)=
d\big(\tr_E\big[f(\nabla)\big]\big)\wedge \big(\exp\circ
\tr_E\big[f(\nabla)\big]\big),
\]
it follows from \eqref{CharFormClosed} that this construction also
gives closed forms.
\item When we decompose the Chern character form into its homogeneous
components
\[
\ch(E,\nabla) = \sum_{j}^{[m/2]}\ch_j(E,\nabla),
\]
then one easily finds relations between $\ch_j$ and the Chern
forms $c_j$ for small $j$. Here, we are dropping the reference to
$E$ and $\nabla$ for the moment. For example,
\[
\ch_0=\rk E,\quad \ch_1 = c_1,\quad \ch_2 = \lfrac 12 c_1^2
-c_2,\quad \ch_3 = \lfrac 16(3c_3 -3c_2c_1+c_1^3),\quad \ldots.
\]
\item If $E_1$ and $E_2$ are two complex vector bundles over $M$
endowed with connections $\nabla_1$ and $\nabla_2$, the following
relations are immediate.
\begin{enumerate}
\item The Chern character form satisfies
\[
\ch(E_1\oplus E_2,\nabla_1\oplus\nabla_2) = \ch(E_1,\nabla_1) +
\ch(E_2,\nabla_2)
\]
and
\[
\ch(E_1\otimes E_2,\nabla_1\otimes 1+ 1\otimes\nabla_2) =
\ch(E_1,\nabla_1) \wedge \ch(E_2,\nabla_2)
\]
\item The total Chern form satisfies
\[
c(E_1\oplus E_2,\nabla_1\oplus\nabla_2) = c(E_1,\nabla_1) \wedge
c(E_2,\nabla_2)
\]
\end{enumerate}
\item If $E$ is equipped with a metric and $\nabla$ is a
compatible connection, then the associated Chern forms and the
Chern character form are $\R$ valued forms. Moreover, assume that
$E$ is of rank $k$ and admits an $\SU(k)$-structure. The latter
means that the determinant line $\det(E)=\gL^kE$ is trivial. Then
for every compatible connection $\nabla$
\[
\ch_{2j+1}(E,\nabla) = 0,
\]
which is due to the fact that the trace of elements in the Lie
algebra $\mathfrak{su}(k)$ vanishes. In particular, the first
Chern form $c_1(E,\nabla)$ is trivial for $\SU(k)$-bundles.
\end{enumerate}
\end{remark*}

\noindent\textbf{Characteristic Forms of Real Vector Bundles.}
For our purposes it is enough to define the characteristic forms
which are obtained by complexifying the bundle and the
connection. In particular, we need not restrict to orthogonal
connections as one would need to define the Euler class.

\begin{dfn}\label{PontLAHatDef}
Let $M$ be an $m$-dimensional manifold, and let $\nabla$ be a
connection on a real vector bundle $E\to M$. Let $E^\C:=E\otimes
\C$ be endowed with the induced connection $\nabla^\C$.
\begin{enumerate}
\item We call the characteristic form
\[
p(E,\nabla):={\det}^{1/2}\big(1+
(\lfrac{i}{2\pi}F_{\nabla^\C})^2\big) = \exp\Big(\lfrac
12\tr_V\Big[\log \big(1+(\lfrac{i}{2\pi}F_{\nabla^\C})^2\big)
\Big]\Big)
\]
the {\em total Pontrjagin form} of $\nabla$.
\item
The \emph{$j$-th Pontrjagin form}
\[
p_j(E,\nabla)\in \gO^{4j}(M)
\]
is defined as the component of degree $4j$ of the total
Pontrjagin form, i.e.,
\begin{equation*}
p(E,\nabla)= \sum_{j=0}^{[m/4]} p_j(E,\nabla) = 1+
p_1(E,\nabla)+p_2(E,\nabla)+\ldots.
\end{equation*}
\item We define the \emph{Hirzebruch $\widehat L$-form} as
\begin{equation*}
\widehat L(E,\nabla):=
{\det}^{1/2}\left(\frac{\lfrac{i}{4\pi}F_{\nabla^\C}}{
\tanh\big(\lfrac{i}{4\pi}F_{\nabla^\C}\big)}\right) \in
\gO^{4\bullet}(M).
\end{equation*}
\item Moreover, the \emph{$\widehat A$-form} is defined as
\begin{equation*}
\widehat A(E,\nabla):=
{\det}^{1/2}\left(\frac{\lfrac{i}{4\pi}F_{\nabla^\C}}{
\sinh\big(\lfrac{i}{4\pi}F_{\nabla^\C}\big)}\right) \in
\gO^{4\bullet}(M).
\end{equation*}
\end{enumerate}
\end{dfn}

\begin{remark*}\quad\nopagebreak
\begin{enumerate}
\item The definition of the characteristic forms above varies in
the literature. First of all, some authors, e.g. \cite{BGV}, drop
the normalizing constants $\frac{i}{2\pi}$ from the definition.
We include them to get integer valued characteristic classes.
Moreover, the $\widehat L$-form is related to the classical
Hirzebruch $L$-form via
\begin{equation}\label{L-L-Hat}
2^{2n}\cdot \widehat L(E,\nabla)_{[4n]} = L(E,\nabla)_{[4n]},
\end{equation}
where $(\ldots)_{[n]}$ means taking the $n$-form component of a
differential form.

\item We have not yet remarked, why the $\widehat L$-form and the
$\widehat A$-form are well-defined. We give some brief remarks
and refer to \cite[App. B]{MS} for more details. Recall that the
Bernoulli numbers $B_n$ can be defined by the following
generating function:
\begin{equation}\label{Bernoulli}
\frac{z}{e^z-1} = \sum_{n= 0}^\infty B_n \frac{z^n}{n!}, \quad
|z|<2\pi,
\end{equation}
see \cite[Sec. 9.1]{Cohen}. With respect to this sign convention,
the first non-trivial $B_n$ are given by
\[
B_0 =1,\quad B_1=-\lfrac 12,\quad B_2=\lfrac 16,\quad B_4 =
-\lfrac 1{30},\quad B_6 = \lfrac 1{42},\quad \ldots
\]
Using \eqref{Bernoulli}, one finds that for $|z|<\pi$
\[
\frac{z/2}{\tanh(z/2)} = 1+\sum_{n\ge 1}\frac{ 1
}{(2n)!}B_{2n}z^{2n} = 1+\lfrac1{12} z^2 -\lfrac 1{720} z^4+\ldots
\]
and
\[
\frac{z/2}{\sinh(z/2)} = 1 +\sum_{n\ge
1}\frac{2^{2n-1}-1}{2^{2n-1}(2n)!}B_{2n}z^{2n}= 1-\lfrac 1{24}z^2
+ \lfrac 1{5760}z^4 +\ldots.
\]
This shows that both are normalized power series, which implies
that the Hirzebruch $\widehat L$-form and the $\widehat A$-form
are well-defined. Moreover, if $\dim M =4$, then
\begin{equation}\label{LA4Man}
\widehat L(E,\nabla) = 1+ \lfrac1{12} p_1(E,\nabla),\quad \widehat
A(E,\nabla) = 1 -\lfrac1{24} p_1(E,\nabla).
\end{equation}

\item Note that if $E$ is endowed with a bundle metric, and
compatible connection $\nabla$, then $\nabla^\C$ on $E^\C$
satisfies $F_{\nabla^\C}^t = -F_{\nabla^\C}$. Thus,
\[
\big(1+\lfrac{i}{2\pi}F_{\nabla^\C}\big)^t = 1-
\lfrac{i}{2\pi}F_{\nabla^\C}.
\]
Hence also,
\[
{\det}^{1/2}\big(1+ \lfrac{i}{2\pi}F_{\nabla^\C} \big) =
{\det}^{1/2}\big(1- \lfrac{i}{2\pi}F_{\nabla^\C} \big),
\]
and so
\[
\begin{split}
\sum_{l=0}^{[m/2]}c_{l}(E^\C,\nabla^\C) &=
\det\big(1+\lfrac{i}{2\pi}F_{\nabla^\C}\big)= {\det}^{1/2}\big(1+
\lfrac{i}{2\pi}F_{\nabla^\C}\big){\det}^{1/2}\big(1+
\lfrac{i}{2\pi}F_{\nabla^\C}\big)\\
&= {\det}^{1/2}\big(1 - ( \lfrac{i}{2\pi}F_{\nabla^\C})^2\big) =
\sum_{j=0}^{[m/4]} (-1)^jp_j(E,\nabla).
\end{split}
\]
From this we deduce that
\[
c_{2j} (E^\C,\nabla^\C) = (-1)^jp_j(E,\nabla),\quad\text{and}\quad
c_{2j-1} (E^\C,\nabla^\C) =0.
\]
\item If $E$ can be written as $E=E_1\oplus E_2$, and $\nabla$
decomposes as $\nabla=\nabla_1\oplus \nabla_2$, the total
Pontrjagin form $p$ and the forms introduced in Definition
\ref{PontLAHatDef} satisfy
\[
p(E) = p(E_1) \wedge p(E_2),\quad \widehat L(E) = \widehat
L(E_1)\wedge \widehat L(E_2),\quad \widehat A(E) = \widehat A(E_1)
\wedge \widehat A(E_2),
\]
where we are dropping the references to the connections.
\end{enumerate}
\end{remark*}

\subsection{Transgression and Characteristic Classes}

As we have seen in \eqref{CharFormClosed}, characteristic forms
associated to a formal power series $f$ as in Definition
\ref{CharFormDef} are closed. Therefore, they define de Rham
cohomology classes. The famous Chern-Weil theorem states that the
difference
\[
\tr_E\big[f(\nabla^1)\big] -\tr_E\big[f(\nabla^0)\big]
\]
for two connections $\nabla^0$ and $\nabla^1$ on $E$ is an exact
form. Therefore, the cohomology class is independent of the
connection. We refer to \cite[Prop. 1.41]{BGV} for a proof of the
following result.

\begin{theorem}\label{ChernWeil}
Let $E\to M$ be a complex vector bundle over a manifold $M$, and
let $f(z)=\sum_{n\ge 0}a_n z^n$ be a formal power series. If
$\nabla^t$ is a smooth path of connections on $E$, then
\[
\lfrac {d}{dt} \tr_E\big[f(\nabla^t)\big] = d \tr_E\big[\lfrac
i{2\pi}(\lfrac d{dt}\nabla^t) \wedge f'(\nabla^t)\big].
\]
In particular, if $a:= \nabla^1 - \nabla^0\in
\gO^1\big(M,\End(E)\big)$ is the difference of two connections,
then
\[
\tr_E\big[f(\nabla^1)\big] -\tr_E\big[f(\nabla^0)\big] = d
\int_0^1 \lfrac{i}{2\pi} \tr_E\big[a \wedge f'(\nabla^0+ta)\big]dt
\]
Therefore, we have an equality of cohomology classes
\[
\Big[\tr_E\big[f(\nabla^1)\big]\Big] =
\Big[\tr_E\big[f(\nabla^0)\big]\Big]\in H^{\ev}(M).
\]
\end{theorem}

\begin{dfn}\label{TransgressDef}
Let $E\to M$ be a complex vector bundle over a manifold $M$, and
let $f(z)=\sum_{n\ge 0}a_n z^n$ be a formal power series.
\begin{enumerate}
\item Let $\nabla$ be an arbitrary connection on $E$.
Then the cohomology class
\[
c_f(E):= \Big[\tr_E\big[f(\nabla)\big]\Big] \in H^{\ev}(M)
\]
is called the \emph{$f$-class} of $E$ or the \emph{characteristic
class} of $E$ associated to $f$.
\item If $M$ is closed and oriented, the number
\[
\big\langle c_f(E),[M] \big\rangle = \int_M c_f(E) \in \R
\]
is called the \emph{characteristic number} of $E$ associated to
$f$. If all characteristic numbers associated to $f$ are
integers, the $f$-class is called \emph{integer valued}.
\item If $\nabla^t$ is a path of connections, we call
\begin{equation*}
Tc_f(\nabla^t):=\int_0^1 \lfrac i{2\pi} \tr_E\big[(\lfrac
d{dt}\nabla^t) \wedge f'(\nabla^t)\big]dt \in \gO^{\odd}(M)
\end{equation*}
the \emph{transgression form of the $f$-class} associated to
$\nabla^t$. If $\nabla^t = \nabla^0+ ta$ we also use the notation
\[
Tc_f(\nabla^0,\nabla^1):= \int_0^1 \lfrac{i}{2\pi} \tr_E\big[a
\wedge f'(\nabla^0 + ta)\big]dt \in \gO^{\odd}(M).
\]
\item The transgression form of the Chern character is called the
\emph{Chern-Simons form} of $\nabla^1$ with respect to $\nabla^0$,
\begin{equation*}
\cs(\nabla^0,\nabla^1) = \int_0^1 \lfrac{i}{2\pi}  \tr_E\big[a
\wedge \exp(\nabla^0 + ta)\big]dt \in \gO^{\odd}(M).
\end{equation*}
\end{enumerate}
\end{dfn}

\begin{remark}\label{ChernWeilRem}\quad\nopagebreak
\begin{enumerate}
\item Theorem \ref{ChernWeil} also applies to characteristic forms
of the form $\exp\big(\tr_E\big[f(\nabla)\big]\big)$. For this
note that
\[
\begin{split}
\lfrac{d}{dt}\exp\big(\tr_E\big[f(\nabla^t)\big]\big)&=
\lfrac{d}{dt}\big(\tr_E\big[f(\nabla^t)\big]\big)\wedge
\exp\big(\tr_E\big[f(\nabla^t)\big]\big)\\
&= d\big( \lfrac i{2\pi}\tr_E\big[(\lfrac d{dt}\nabla^t) \wedge
f'(\nabla^t)\big]\big)\wedge \exp\big(
\tr_E\big[f(\nabla^t)\big]\big).
\end{split}
\]
This form is exact, since
$\exp\big(\tr_E\big[f(\nabla^t)\big]\big)$ is closed. Hence, the
transgression form in this case is
\[
\int_0^1 \lfrac i{2\pi} \tr_E\big[(\lfrac d{dt}\nabla^t) \wedge
f'(\nabla^t)\big] \wedge \exp\big(
\tr_E\big[f(\nabla^t)\big]\big)dt
\]
\item When considering the cohomology class of one of the
particular characteristic forms introduced in the last section,
we will call them \emph{Chern character}, \emph{Chern class},
\emph{$\widehat L$-class}, etc. The distinction between forms and
classes is done by incorporating the connection in the notation.
For example,
\[
\ch(E,\nabla)\in \gO^{\ev}(M),\quad\text{but}\quad \ch(E)\in
H^{\ev}(M).
\]
\item The Chern and Pontrjagin classes are integer valued due to
the normalization factor of $\lfrac{i}{2\pi}$, see \cite[App.
C]{MS}. The other characteristic classes we have defined are in
general only $\Q$ valued.
\item Often the term Chern-Simons form is reserved for the degree
3 part of what we have called the Chern-Simons form. Due to its
importance in 3-manifold topology, we want to derive an explicit
formula for it. We abbreviate $\nabla:=\nabla^0$ and let $F_t$
denote the curvature of $\nabla^t:=\nabla +ta$. Then
\[
F_t = F_\nabla + t(\nabla a) +t^2 a\wedge a.
\]
For the component of degree 4 of the Chern character form we have
$f(z)=z^2/2$, so that $f'(z)=z$. According to Definition
\ref{TransgressDef},
\[
\cs(\nabla^0,\nabla^1)_{[3]}= -\frac{1}{4\pi^2}\int_0^1
\tr_E\big[a \wedge \big(F_\nabla+t\nabla a+ t^2 a\wedge
a\big)\big]dt.
\]
Integrating this expression we get
\begin{equation}\label{ChernSimonsExpl}
\cs(\nabla^0,\nabla^1)_{[3]}= -\frac {1}{4\pi^2}\tr_E\big[a\wedge
F_\nabla + \lfrac 12 a\wedge\nabla a + \lfrac 13 a\wedge a\wedge
a \big].
\end{equation}
In particular, if $\nabla$ is a flat connection we get the
well-known expression
\begin{equation*}
\cs(\nabla^0,\nabla^1)_{[3]}  = -\frac
{1}{8\pi^2}\tr_E\big[a\wedge \nabla a + \lfrac 23 a\wedge a\wedge
a \big].
\end{equation*}
\end{enumerate}
\end{remark}

\section{Chern-Simons Invariants}\label{CSInvariants}

There is also a different description of transgression forms,
which we shall describe now. Let $E\to M$ be a complex vector
bundle over a manifold $M$, endowed with a path $\nabla^t$ of
connections. Over the cylinder $N:=[0,1]\times M$, we consider
the vector bundle $\pi^*E\to N$, where $\pi: N\to M$ is the
natural projection. The path $\nabla^t$ defines a connection on
$\pi^*E$ via
\begin{equation}\label{ConnCyl}
\widetilde \nabla := dt\wedge\lfrac d{dt} + \pi^*\nabla^t.
\end{equation}
Its curvature is easily seen to be given by
\[
F_{\widetilde \nabla} = dt\wedge\big(\lfrac
d{dt}\pi^*\nabla^t\big) + \pi^* F_{\nabla^t}.
\]
Since $dt\wedge dt=0$, one deduces from the trace property that
for all $n\ge 1$,
\[
\tr_{\pi^*E}\big[F_{\widetilde \nabla}^n \big] = dt\wedge \pi^*
\tr_E\big[n \big(\lfrac d{dt}\nabla^t\big)\wedge
F_{\nabla^t}^{n-1}\big] + \pi^*\tr_E\big[F_{\nabla^t}^n \big].
\]
This implies that for any formal power series $f(z)=\sum_n a_n
z^n$,
\begin{equation}\label{CharFormCyl}
\tr_{\pi^*E}\big[f(\widetilde \nabla)\big] = \lfrac
i{2\pi}dt\wedge \pi^* \tr_E\big[\big(\lfrac
d{dt}\nabla^t\big)\wedge f'(\nabla^t)\big] +
\pi^*\tr_E\big[f(\nabla^t) \big].
\end{equation}
Now, consider integration along the fiber as in Proposition
\ref{FiberInt},
\[
\int_{N/M}: \gO^\bullet(N)\to \gO^{\bullet-1}(M).
\]
Then, comparing \eqref{CharFormCyl} and the definition of the
transgression form in Definition \ref{TransgressDef}, one readily
obtains

\begin{lemma}\label{TransgressCyl}
If $\nabla^t$ is a path of connections over $M$, and $\widetilde
\nabla$ denotes the associated connection \eqref{ConnCyl} over the
cylinder $N:=[0,1]\times M$, then
\[
Tc_f(\nabla^t) = \int_{N/M} \tr_{\pi^*E}\big[f(\widetilde
\nabla)\big]\in\gO^{\odd}(M).
\]
\end{lemma}

Using this result we can now derive the following important
property of transgression forms, see \cite[Sec. 3] {CS}.

\begin{prop}\label{CSPathIndep}
Let $E\to M$ be a complex vector bundle over a manifold $M$. If
$\nabla^t$ is a closed path of connections on $E$, then
\[
Tc_f(\nabla^t) \in d\gO^{\ev}(M).
\]
\end{prop}

\begin{proof}
Since the space of connections on $E$ is contractible, we can
find a smooth two-parameter family $\nabla^{s,t}$ of connections
which gives a homotopy relative endpoints from $\nabla^t$ to the
constant path. On the cylinder $N$ we consider the one-parameter
family
\[
\widetilde \nabla^s := dt\wedge \lfrac d{dt} + \nabla^{s,t},
\]
where we are dropping the pullback with $\pi$ from the notation.
Using Theorem \ref{ChernWeil} and \eqref{CharFormCyl} one finds
that
\[
\begin{split}
\lfrac d{ds} \tr_E\big[f(\widetilde \nabla^s)\big] &= - \frac
1{4\pi^2}\, d_N \tr_{\pi^*E}\Big[\big(\lfrac
d{ds}\nabla^{s,t}\big)\wedge dt\wedge \big(\lfrac
d{dt}\nabla^{s,t}\big)\wedge f''(\nabla^{s,t})\Big]\\
&\qquad \qquad \qquad + d_N \tr_E\Big[\lfrac i{2\pi}\big(\lfrac
d{ds}\nabla^{s,t}\big)\wedge f'(\nabla^{s,t})\Big]\\
&=: d_N \big(dt\wedge \ga(s,t)\big) + d_N \gb(s,t),
\end{split}
\]
where $\ga(s,t)$ and $\gb(s,t)$ are two-parameter families of
differential forms on $M$. Then
\[
\begin{split}
\frac d{ds}\int_{N/M} \tr_E\big[f(\widetilde \nabla^s)\big] &=
\int_{N/M} d_N \big(dt\wedge \ga(s,t)\big) + \int_{N/M}d_N
\gb(s,t)\\
&= d_M \int_{N/M}dt\wedge \ga(s,t) +
\int_{N/M}dt\wedge\big(\lfrac d{dt} \gb(s,t)\big)\\
&= d_M \int_{N/M}dt\wedge \ga(s,t) + \gb(s,1) - \gb(s,0).
\end{split}
\]
By assumption, $\nabla^{s,1}=\nabla^{s,0}$ is constant for all
$s$. Checking the explicit formula for $\gb(s,t)$ one finds that
$\gb(s,1) = \gb(s,0)$. According to Lemma \ref{TransgressCyl},
this shows that
\[
\lfrac d{ds} Tc_f(\nabla^{s,t}) = \frac d{ds}\int_{N/M}
\tr_E\big[f(\widetilde \nabla^s)\big] \in d\gO^{\ev}(M),
\]
from which the result follows.
\end{proof}

\noindent\textbf{Chern-Simons Invariants.} The last result shows
that transgression forms can be used to define numerical
invariants associated to pairs of connections on odd dimensional
manifolds. In this respect they are odd analogues of
characteristic numbers.

\begin{dfn}\label{CSDef}
Let $M$ be a closed manifold, and let $f(z)=\sum_{n\ge 0}a_n z^n$
be a formal power series. If $\nabla^0$ and $\nabla^1$ are two
connections on a complex vector bundle $E\to M$ we define the
\emph{Chern-Simons invariant} of $\nabla^1$ with respect to
$\nabla^0$ associated to $f$ as
\[
\CS_f(\nabla^0,\nabla^1):=\int_M Tc_f(\nabla^0,\nabla^1).
\]
\end{dfn}

\begin{prop}\label{CSProp}
Let $M$ be a closed manifold, and let $f$ be a formal power
series. Consider two connections $\nabla^0$ and $\nabla^1$ on a
complex vector bundle $E\to M$.
\begin{enumerate}
\item If $\nabla^t$ is any path connecting $\nabla^0$ and
$\nabla^1$, then
\[
\CS_f(\nabla^0,\nabla^1) = \int_M Tc_f(\nabla^t).
\]
\item Let $\nabla^2$ be a third connection on $E$, then
\[
\CS_f(\nabla^0,\nabla^2) = \CS_f(\nabla^0,\nabla^1)
+\CS_f(\nabla^1,\nabla^2).
\]
\item Assume that $\nabla^N$ is a connection over the cylinder
$N=[0,1]\times M$ such that on a collar of the boundary it is of
the form \eqref{ConnCyl}. Then
\begin{equation}\label{CSCyl}
\CS_f(\nabla^0,\nabla^1) = \int_N \tr_E\big[f(\nabla^N)\big].
\end{equation}
\item If $W$ is a compact manifold with boundary
$M$, and $E$, $\nabla^0$ and $\nabla^1$ extend to $E_W$,
$\hat\nabla^0$ and $\hat\nabla^1$, then
\[
\CS_f(\nabla^0,\nabla^1) = \int_{W}
\tr_E\big[f(\hat\nabla^1)\big] - \int_{W}
\tr_E\big[f(\hat\nabla^0)\big].
\]
\item Assume that $f$ gives an integer valued characteristic
class, and let $\gF:E\to E$ be a bundle isomorphism. Then for
every connection $\nabla$ on $E$,
\[
\CS_f(\nabla,\gF^*\nabla) \in\Z.
\]
\end{enumerate}
\end{prop}

\begin{proof}[Sketch of proof]
Since $M$ is assumed to be closed, part (i) follows from
Proposition \ref{CSPathIndep}. Part (ii) is an immediate
consequence of (i). For (iii) let $\nabla^t$ be a path connecting
$\nabla^0$ and $\nabla^1$ such that on a collar of the boundary,
$\nabla^N$ and $\widetilde \nabla$ as in \eqref{ConnCyl} agree.
Theorem \ref{ChernWeil} implies that $\nabla^N - \widetilde
\nabla$ is the differential of a form on $N$ with compact support
away from the boundary. Then Stokes' Theorem readily yeilds (iii).
Part (iv) also follows from Theorem \ref{ChernWeil} and Stokes'
Theorem.\footnote{Alternatively, one could glue a cylinder to the
boundary of $W$ to interpolate between $\hat\nabla^0$ and
$\hat\nabla^1$, and then use Lemma \ref{TransgressCyl} as well as
the additivity under cutting and pasting of characteristic
numbers.} For part (v) denote by $\gf$ the map covered by $\gF$.
Then the mapping torus
\[
E_\gF:= \big([0,1]\times E\big)/ \sim\; ,\quad (1,x)\sim
(0,\gF(x))
\]
is a Hermitian vector bundle over the mapping torus $M_\gf$.
Endow $E_\gF$ with a connection $\nabla^\gF$, induced by
connecting $\nabla$ and $\gF^*\nabla$ over $M$. Then one easily
finds that
\[
\CS_f(\nabla,\gF^*\nabla) = \int_{M_\gf}
\tr_{E_\gF}\big[f(\nabla^\gF)\big].
\]
The right hand side is integer valued as $M_\gf$ is closed.
\end{proof}

\begin{remark}\label{CSMultSeq}
For a characteristic class of the form $\exp\big(
\tr_E\big[f(\nabla)\big]\big)$ we have seen in Remark
\ref{ChernWeilRem} that the transgression form is given by
\[
\int_0^1 \lfrac i{2\pi} \tr_E\big[(\lfrac d{dt}\nabla^t) \wedge
f'(\nabla^t)\big] \wedge \exp\big(
\tr_E\big[f(\nabla^t)\big]\big)dt.
\]
Lemma \ref{TransgressCyl} extends to this context: Let
$\widetilde \nabla$ be the connection over $N=[0,1]\times M$, and
write
\[
\tr_E\big[f(\widetilde \nabla)\big] = dt\wedge \ga(t) +\gb(t),
\]
where $\ga(t)$ and $\gb(t)$ contain no $dt$-factor. Then
\[
\exp\big(\tr_E\big[f(\widetilde \nabla)\big]\big) =
\exp\big(dt\wedge\ga(t)\big)\wedge \exp\big(\gb(t)\big) =
\big(1+dt\wedge\ga(t)\big) \wedge \exp\big(\gb(t)\big).
\]
Therefore,
\[
\begin{split}
\int_{N/M} \exp\big(\tr_E\big[f(\widetilde \nabla)\big]\big)&=
\int_{N/M} dt\wedge\ga(t)\wedge \exp\big(\gb(t)\big)\\
&= \int_{N/M} \tr_E\big[f(\widetilde \nabla)\big]\wedge
\exp\big(\tr_E\big[f(\nabla^t)\big]\big)\\
&= \int_0^1 \lfrac i{2\pi}\tr_E\big[(\lfrac d{dt}\nabla^t) \wedge
f'(\nabla^t)\big] \wedge \exp\big(
\tr_E\big[f(\nabla^t)\big]\big)dt
\end{split}
\]
where we have used \eqref{CharFormCyl} in the last line.
Similarly, one checks that Proposition \ref{CSPathIndep} and
Proposition \ref{CSProp} continue to hold in this context.
\end{remark}

\cleardoublepage
\chapter{Remarks on Moduli Spaces}\label{ModApp}

In this appendix we include some details concerning the moduli
space of flat connections and the moduli space of holomorphic
line bundles over a Riemann surface. Since the Rho invariant
depends only on the gauge equivalence class of the underlying
flat connection, understanding the moduli space of flat
connections is a prerequisite for the computation of Rho
invariants. Moreover, the interplay between flat connections and
representations of the fundamental group is often used in the
main body of this thesis. Therefore, we start with a detailed
discussion of these topics, in particular including some remarks
on the question of whether a given flat bundle is trivializable
or not.

We proceed with a discussion of the moduli space of flat
connections associated to a mapping torus. Here the objective is
to prove the facts we have used in Chapter \ref{3dimMapTor}.
After this, we add some remarks about the moduli space of
holomorphic line bundles over a Riemann surface and its relation
to the moduli space of (flat) connections. This will establish
some facts we have freely used in Section \ref{S1Bundles}.

\section{The Moduli Space of Flat Connections}\label{FlatConn}

\subsection{Flat Connections and Representations of the
Fundamental Group}\label{FlatConnRep}

Since many features become more transparent in a more general
setup, we start working with a principal $G$-bundles, where $G$ is
an arbitrary connected matrix Lie group. Ultimately we are
interested in flat Hermitian vector bundles and restrict to
$G=\U(k)$. A general reference for the contents of this section
are \cite[Sec. 2.1]{DK} and \cite[Ch. II]{KN}.

Denote by $\cg$ the Lie algebra of $G$. Since we are assuming that
$G$ is a matrix Lie group, $\cg$ is a matrix Lie algebra. We use
the notation ``$\Ad$'' for the adjoint action of $G$ on itself and
``$\ad$'' to denote the adjoint
action of $G$ on $\cg$.\\

\noindent\textbf{Connections and Curvature.} Let $M$ be a
connected manifold, and let $P\xrightarrow{\pi} M$ be a principal
$G$-bundle. Let $R_g$ denote the right-action of $g\in G$ on $P$.
Recall that a \emph{$G$-connection} on $P$ is a Lie algebra
valued 1-form $A\in \gO^1(P,\cg)$ satisfying
\begin{equation}\label{ConnDef}
R_g^* A = g^{-1}Ag,\quad\text{and}\quad A\big(\lfrac
d{dt}\big|_{t=0}p\cdot\exp(tX)\big) = X,\quad p\in P, \quad X\in
\mathfrak g.
\end{equation}
We denote the space of all $G$-connections on $P$ by $\cA(P)$.
The \emph{curvature} of $A$ is defined as
\[
F_A = d A + A\wedge A\in \gO^2(P,\cg),
\]
where $A\wedge A$ stands for taking the exterior product in the
form part and matrix multiplication in the Lie algebra part. The
curvature is easily seen to be $\ad$-equivariant and horizontal,
i.e.,
\[
R_g^* F_A = g^{-1}F_A g,\quad \imu\big(\lfrac
d{dt}\big|_{t=0}p\cdot\exp(tX)\big) F_A = 0.
\]
This implies that $F_A$ can also be viewed as a 2-form on $M$
with values in the bundle $\ad(\cg)= P\times_{\ad} \cg$. A
connection $A$ is called \emph{flat} if $F_A=0$, and we denote by
\[
\cF(P)=\setdef{A\in \cA(P)}{F_A=0}
\]
the space of flat $G$-connections on $P$.\\

\noindent\textbf{Gauge Transformations and the Moduli Space.} We
also recall that a \emph{gauge transformation} is a
$G$-equivariant bundle isomorphism,
\[
\gF:P\to P,\quad \gF(p\cdot g) = \gF(p)\cdot g.
\]
If one defines $u:P\to G$ by requiring that $\gF(p)= p\cdot u(p)$,
then $u$ is $\Ad$-equivariant,
\[
u : P\to G,\quad u(p\cdot g) =  g^{-1} u(p) g.
\]
Conversely, it is easy to see that every gauge transformation
arises this way. Hence, one of several equivalent ways to define
the group of gauge transformations is
\[
\cG(P) := C^\infty\big(M,\Ad(P)\big),\quad \text{where}\quad
\Ad(P)= P\times_{\Ad} G.
\]
The pullback of a connection by a gauge transformation gives a
natural action of $\cG(P)$ on $\cA(P)$. In terms of an
$\Ad$-equivariant map $u:P\to G$ this takes the form
\[
A\cdot u = u^{-1}A u + u^{-1}d u,\quad A\in \cA(P),\quad u\in
\cG(P).
\]
We point out that $u^{-1}d u$ is the pullback of the
Maurer-Cartan form on $G$ via $u$. The curvature behaves
equivariantly with respect to this action,
\[
F_{A\cdot u} = u^{-1} F_A u, \quad A\in \cA(P),\quad u\in \cG(P).
\]
In particular, the action of $\cG(P)$ on $\cA(P)$ leaves the space
$\cF(P)$ of flat connections invariant. One can thus define the
\emph{moduli space} of flat connections on $P$ as
\[
\cM(P):= \cF(P)/\cG(P).
\]

\begin{remark}\label{GaugeTrivBundle}
We will often encounter the situation that $P$ is trivializable.
If we fix a trivialization $P\cong M\times G$, we can identify
\[
\cA(M\times G)\cong \gO^1(M,\cg),\quad \cG(M\times G)\cong
C^\infty(M,G).
\]
\end{remark}

\noindent\textbf{The Holonomy Representation.} Fix a base point
$p_0\in P$, and let $x_0:=\pi(p_0)$. Consider a closed loop based
at $x_0$, i.e.,
\[
c:I \to M,\quad I=[0,1],\quad\text{with}\quad c(0)=c(1)=x_0.
\]
Since $I$ is contractible, the pullback $c^*P\to I$ is
trivializable. As we are assuming that $G$ is connected, we can
fix a lift $\widehat c:I\to P$ of $c$ such that $\widehat
c(0)=\widehat c(1)=p_0$. Now let $A$ be a $G$-connection---not
necessarily flat for the moment---and let
\begin{equation}\label{HolLocalConn}
A_t:=(\widehat c^{\,*} A)(\pd_t)\in C^\infty(I,\cg).
\end{equation}

\begin{dfn}\label{HolDef}
Let $g_t:I\to G$ be the unique solution of the ordinary
differential equation
\[
\pd_t g_t = - A_tg_t,\quad g_0=e,
\]
where $e\in G$ is the identity element. Then the \emph{holonomy}
of $A$ along $c$ with respect to the base point $p_0$ is defined
by
\[
\hol_A(c,p_0) := g_1\in G.
\]
\end{dfn}

Note that the definition gives no reference to the lift $\widehat
c$ we have fixed. The reason why we are allowed to do so is one of
the contents of the following result.

\begin{lemma}\label{HolLem}
Let $A$ be a connection on $P$, and let $c:I\to M$ be a closed
loop, based at $x_0$.
\begin{enumerate}
\item If $\gf:I\to I$ is an orientation preserving
re\-para\-metri\-zation, then
\[
\hol_{A}(c\circ \gf,p_0) = \hol_{A}(c,p_0).
\]
\item For every gauge transformation $u\in \cG(P)$,
\[
\hol_{A\cdot u}(c,p_0) = u(p_0)^{-1} \hol_A(c,p_0)\, u(p_0).
\]
In particular, $\hol_A(c,p_0)$ is independent of the lift
$\widehat c$ chosen in its definition.
\item Assume that $\widetilde c$ is another loop, based at $x_0$,
and denote by $c*\widetilde c$ the loop defined by first running
along $c$ and then along $\widetilde c$. Then
\[
\hol_A(c*\widetilde c,p_0) = \hol_A(\widetilde c,p_0)
\hol_A(c,p_0).
\]
\item Let $p_1\in P$ be a different base point, and $\widehat c_0:I\to P$
be a path connecting $p_0$ with $p_1$. Then there exists $g\in G$
such that
\[
\hol_A(c_0^{-1}*c*c_0,p_1) = g \hol_A(c,p_0)\, g^{-1},\quad
\text{where $c_0:=\pi\circ \widehat c_0$.}
\]
\end{enumerate}
\end{lemma}

\begin{proof}\footnote{Although the assertions in Lemma
\ref{HolLem} is standard, we include a proof as the discussion to
follow relies on similar arguments.} To prove part (i) let $g_t$
be as in Definition \ref{HolDef}. Then
\[
\pd_s g_{\gf(s)} = \gf'(s) (\pd_tg_t)|_{t=\gf(s)} = - \gf'(s)
A_{\gf(s)}g_{\gf(s)} = -(\widehat c\circ \gf)^* A(\pd_s)
g_{\gf(s)}.
\]
Moreover, we have $g_{\gf(0)}= g_0 = e$ and $g_{\gf(1)}= g_1$.
Then part (i) is true by definition. To prove (ii) define $u_t:=
u\circ \widehat c\in C^\infty(I,G)$, and $\widehat c_u:= \widehat
c\cdot u_t$. Then $\widehat c_u$ is a lift of $c$ with
\[
\widehat c_u(0)=\widehat c_u(1)=p_0\cdot u(p_0),
\quad\text{and}\quad \widehat c^{\,*} (A\cdot u)(\pd_t)=\widehat
c_u^* A(\pd_t) = A_t\cdot u_t.
\]
If $g_t$ is as in Definition \ref{HolDef}, then
\[
\begin{split}
\pd_t (u^{-1}_tg_t u_0) &=  (\pd_t u_t^{-1})g_t u_0 + u^{-1}(\pd_t
g_t)u_0\\ &= - u_t^{-1}(\pd_t u_t)(u^{-1}_tg_t u_0) - (u^{-1}_tA_t
u_t)(u^{-1}_tg_t u_0)\\ &= - (A_t\cdot u_t)(u^{-1}_tg_t u_0).
\end{split}
\]
Since $u^{-1}_0g_0 u_0 = e$, this implies that
\[
\hol_{A\cdot u_t}(c,p_0\cdot g) = u^{-1}_1g_1 u_0 =
u^{-1}(p_0)\hol_A(c,p_0) \, u(p_0).
\]
The second assertion of (ii) follows from the fact that every
lift of $c$ with base point $p_0$ is of the form $\widehat c\cdot
u$ with $u(p)=e$. Concerning part (iii), we define $\widetilde
A_t$ as in \eqref{HolLocalConn} with respect to a lift of
$\widetilde c$ and note without going into detail that
$\hol_A(c*\widetilde c,p_0)$ is given by $\widetilde g_1$, where
$\widetilde g_t$ is the unique solution to
\[
\pd_t\widetilde g_t= - \widetilde A_t \widetilde g_t,\quad
\widetilde g_0 = \hol_A(c,p_0).
\]
This readily yields $\widetilde g_1 = \hol_A(\widetilde c,p_0)
\hol_A(c,p_0)$. To prove (iv) we can solve the initial value
problem
\[
\pd_t g_t = - (\widehat c_0^{\,*} A )(\pd_t)g_t,\quad g_0 =e.
\]
Then one verifies without effort that the assertion holds with
$g:=g_1$.
\end{proof}

Let $\gO(M,x_0)$ be the \emph{based loop group} of $M$. This is
the set of all loops, based at $x_0$, modulo orientation
preserving reparametrization. For reasons of functoriality we
endow $\gO(M,x_0)$ with the product $c\cdot \widetilde c :=
\widetilde c*c$, where $\widetilde c* c$ is as in Lemma
\ref{HolLem}. Let $p_0\in P$ with $x_0=\pi(p_0)$. Using the above
results, one obtains a well-defined homomorphism
\begin{equation}\label{HolRepDef}
\hol_A: \gO(M,x_0)\to G,\quad c\mapsto \hol_A(c,p_0)
\end{equation}

\begin{dfn}
Let $A$ be a connection on $P$. Then the homomorphism
\eqref{HolRepDef} is called the \emph{holonomy representation} of
$A$ with respect to the base point $p_0$. We also define the
\emph{holonomy group} of $A$ with respect to $p_0$ as
\[
G_A(p_0):= \im\big(\hol_A:\gO(M,x_0)\to G\big)
\]
A connection $A$ is called \emph{irreducible}, if $G_A(p_0) = G$.
Otherwise, it is called \emph{reducible}. Moreover, the
\emph{isotropy group} of $A$ is defined as
\[
I(A):= \bigsetdef{u\in \cG(P)}{A\cdot u = A}.
\]
\end{dfn}

\begin{lemma}\label{CentralLem}
The conjugacy class of $G_A(p_0)$ is independent of $p_0$ and the
gauge equivalence class of $A$. Moreover, for fixed $p_0\in P$,
the map
\[
I(A) \to G,\quad u\mapsto u(p_0),
\]
maps $I(A)$ isomorphically to the centralizer of $G_A(p_0)$ in
$G$.
\end{lemma}

\begin{proof}
The first assertion is immediate from Lemma \ref{HolLem}. Now fix
$p_0$, and assume that $u\in I(A)$. Lemma \ref{HolLem} implies
that for every loop $c:I\to M$ , based at $x_0=\pi(p_0)$, we have
\[
\hol_A(c,p_0)=\hol_{A\cdot u}(c,p_0) = u(p_0)^{-1} \hol_A(c,p_0)\,
u(p_0).
\]
Thus, $u(p_0)$ lies in the centralizer of $G_A(p_0)$. Now let
$P_0$ be the set of all $p\in P$ such that there exists a
horizontal path $\widehat c:I\to P$ with $\widehat c(0)=p_0$ and
$\widehat c(1)=p$. Then $P_0$ intersects every fiber of $\pi:P\to
M$. This is because if $\widehat c:I\to P$ is an arbitrary path
with $\widehat c(0)=p_0$, we can solve
\[
\pd_t g_t = - A_tg_t,\quad g_0 =e,
\]
with $A_t$ as in \eqref{HolLocalConn}, to get a horizontal path
$\widehat c\cdot g_t:I\to P$ whose endpoint lies in the same
fiber as the endpoint of $\widehat c$.

To prove injectivity of the map $I(A)\to G$, $u\mapsto u(p_0)$,
let $u\in I(A)$ with $u(p_0)=e$. Since $u$ is $\Ad$-equivariant
and $P_0$ intersects every fiber of $\pi:P\to M$, it suffices to
show that $u|_{P_0}\equiv e$. Let $p\in P_0$, and let $\widehat c$
a horizontal path connecting $p_0$ with $p$. Then $A\cdot u =A$
implies that
\[
(u^{-1}\circ \widehat c) \pd_t (u\circ \widehat c )(t) =
u^{-1}du|_{\widehat c(t)}(\lfrac d{dt} \widehat c) = (A - u^{-1}A
u)|_{\widehat c(t)}(\lfrac d{dt} \widehat c )=0,
\]
where we have used that $\lfrac d{dt} \widehat c$ is horizontal.
Hence $u$ is constant along $\widehat c$ and thus, $u(p)=e$.
Next, assume that $g_0\in G$ lies in the centralizer of
$G_A(p_0)$. We need to define $u\in I(A)$ such that $u(p_0)=g_0$.
We first define $u|_{P_0}$ to be the constant map $g_0$. To see
that this defines a gauge transformation, we need to check that
$u(p\cdot g)=g^{-1} u(p) g$, whenever $p$ and $p\cdot g$ both lie
in $P_0$. Let $\widehat c$ and $\widehat c_g$ be horizontal paths
connecting $p_0$ with $p$, respectively with $p\cdot g$. Then
$\widehat c_g*(\widehat c^{-1}\cdot g)$ is a horizontal path
which connects $p_0$ with $p_0\cdot g$. This implies that $g\in
G_A(p_0)$. Since we have assumed that $g_0$ lies in the
centralizer of $G_A(p_0)$, we obtain
\[
u(p\cdot g)=g^{-1} u(p) g = g^{-1} g_0 g = g_0.
\]
Hence, $u|_{P_0}$ is $\Ad$-equivariant and can be extended to a
gauge transformation on $P$. To prove that $u\in I(A)$ first note
that the values of the 1-forms $A\cdot u$ and $A$ on vertical
vectors are both prescribed by \eqref{ConnDef}. To see that
$A\cdot u$ and $A$ also agree on horizontal vectors, it suffices
to consider $A|_p$ and $A\cdot u|_p$ for $p\in P_0$ since both,
$A\cdot u$ and $A$, are $\ad$-equivariant. Now if $v\in T_pP$ is
horizontal, there exists a horizontal path $\widehat
c:(-\eps,\eps)\to P$ with $\widehat c(0)=p$ and $\lfrac
d{dt}\widehat c(0) = v$. By definition of $P_0$ one easily checks
that $\im(\widehat c)\subset P_0$. Thus, $u$ is constant along
$\widehat c$ so that
\[
A\cdot u|_p(v) - A|_p(v) = u^{-1}du|_p(v) = u^{-1}(p)\lfrac
d{dt}\big|_{t=0} u\circ \widehat c = 0.\qedhere
\]
\end{proof}

\noindent\textbf{Flat Connections and the Fundamental Group.}
After this technical preparation, we turn our attention to flat
connections. The following result shows that they are of a
topological nature.

\begin{prop}\label{HolFlat}
If $A$ is flat, then the holonomy $\hol_A(c,p_0)$ depends only on
the homotopy class $[c]\in \pi_1(M)=\pi_1(M,x_0)$. In particular,
the holonomy representation defines a homomorphism
\[
\hol_A \in \Hom\big(\pi_1(M),G\big).
\]
Moreover, the assignment $\cF(P)\to \Hom\big(\pi_1(M),G\big)$,
$A\mapsto \hol_A$ gives well-defined map
\[
\cM(P)\to \Hom\big(\pi_1(M),G\big)/G.
\]
\end{prop}

\begin{proof}
Consider a homotopy
\[
c:I\times I\to M, \quad c(s,0)=c(s,1)=x_0.
\]
Since any fiber bundle has the homotopy lifting property, we can
choose a lift
\[
\widehat c:I\times I\to P,\quad \widehat c(s,0)=\widehat
c(s,1)=p_0
\]
Abusing notation we use the letter $A$ also to denote the 1-form
$\widehat c^{\,*}A\in \gO^1(I\times I,\cg)$. The flatness
condition $dA+A\wedge A=0$ written out with respect to the
coordinates $(s,t)\in I\times I$ is
\begin{equation}\label{HolProp:1}
\pd_s A(\pd_t) -\pd_tA(\pd_s) + A(\pd_s)A(\pd_t) -
A(\pd_t)A(\pd_s)=0.
\end{equation}
For fixed $s$ let $g_s=g_s(t):I\to  G$ denote the solution to
\begin{equation}\label{HolProp:2}
\pd_t g_s = - A(\pd_t) g_s,\quad g_s(0)=e.
\end{equation}
Since $A$ depends smoothly on $s$ and $t$, it follows from the
standard theory of ordinary differential equations that $g_s$
depends smoothly on $s$. We then compute
\[
\begin{split}
\pd_t\big(\pd_sg_s + A(\pd_s)g_s\big) &= \pd_s(\pd_tg_s) +
\big(\pd_t A(\pd_s)\big)g_s + A(\pd_s)(\pd_tg_s)\\
&= - \pd_s\big(A(\pd_t) g_s \big)+ \big(\pd_tA(\pd_s)\big)g_s -
A(\pd_s)A(\pd_t) g_s,
\end{split}
\]
where we have used \eqref{HolProp:2}. Then \eqref{HolProp:1}
implies that
\[
\begin{split}
\pd_t\big(\pd_sg_s + A(\pd_s)g_s\big) &= \big(A(\pd_s)A(\pd_t) -
A(\pd_t)A(\pd_s)\big) g_s -A(\pd_t)\pd_sg_s  - A(\pd_s)A(\pd_t)
g_s\\
&= - A(\pd_t)\big(\pd_sg_s + A(\pd_s)g_s\big).
\end{split}
\]
The initial condition in \eqref{HolProp:2} and the fact that
$\widehat c(s,0)$ is constant for $s\in I$ implies that
$\pd_sg_s|_{t=0}=0$ and $A(\pd_s)|_{t=0}=0$. Hence,
\[
\pd_sg_s = - A(\pd_s)g_s,\quad \text{for all $t\in I$.}
\]
Moreover, we have $A(\pd_s)|_{t=1}=0$ so that $g_s(1)$ is
independent of $s$. By definition of the holonomy and
\eqref{HolProp:2} this proves the first assertion of Proposition
\ref{HolFlat}. The other assertions are immediate from Lemma
\ref{HolLem}.
\end{proof}

\noindent\textbf{The Moduli Space of Representations.} Let
$\cM(M,G)$ be the moduli space of flat principal $G$-bundles,
i.e., the space of isomorphism classes of pairs $(P,A)$ where $P$
is a principal $G$-bundle and $A$ is a flat connection on $P$.
Our next goal is to show that the map in Proposition \ref{HolFlat}
induces an isomorphism
\[
\cM(M,G)\cong \Hom\big(\pi_1(M),G\big)/G.
\]

\begin{remark}\label{FlatTriv}
Note that in general $\cM(M,G)$ will be strictly larger than the
moduli space $\cM(P)$ for one fixed flat bundle $P$. This is
because there might exist flat bundles such that the underlying
principal $G$-bundles are not isomorphic. In the case that
$G=\U(k)$ we will say more about this in Section
\ref{FlatTrivSec} below.
\end{remark}

To associate a flat $G$-bundle to any representation
$\ga:\pi_1(M)\to G$, let $\widetilde M$ be the universal cover of
$M$. For definiteness we fix a base point $x_0\in M$ and identify
$\widetilde M$ with the space of homotopy classes of paths in $M$
starting at $x_0$. Then $\pi_1(M,x_0)$ naturally acts on
$\widetilde M$ from the right. For $\ga: \pi_1(M,x_0)\to G$ we
define the principal $G$-bundle
\[
P_\ga := \widetilde M\times_\ga G = (\widetilde M\times G)/\sim
\]
where
\begin{equation}\label{ActFundDef}
(\widetilde x, g)\sim (\widetilde x\cdot c, \ga(c)^{-1} g),\quad
(\widetilde x,g)\in \widetilde M\times G,\quad c\in \pi_1(M).
\end{equation}
Pulling back the Maurer-Cartan form $g^{-1}dg\in \gO^1(G,\cg)$ to
$\widetilde M\times G$ defines a natural flat connection on
$\widetilde M\times G$, which is invariant under the action
\eqref{ActFundDef} of $\pi_1(M)$. In this way we get an induced
flat connection $A_\ga$ on $P_\ga$. It is straightforward to
check that with respect to the base point $p_0:=[x_0,e]\in P_\ga$,
\[
\hol_{A_\ga}\big(c,p_0\big) = \ga(c),\quad c\in \pi_1(M,x_0).
\]
More generally, we have

\begin{prop}\label{FlatModuliRepVar}
Let $P$ be a principal $G$-bundle with flat connection $A$, and
let $\ga:\pi_1(M)\to G$ be a representation of the fundamental
group. Then $(P,A)$ is isomorphic to $(P_\ga,A_\ga)$ if and only
if there exists $g\in G$ with $\hol_A = g^{-1} \ga g$. In
particular, we have a bijection
\[
\cM(M,G)\xrightarrow{\cong} \Hom\big(\pi_1(M),G\big)/G,\quad
[P,A]\mapsto [\hol_A].
\]
\end{prop}

\begin{proof}[Sketch of proof]
The assertion that $(P,A)\cong (P_\ga,A_\ga)$ implies $\hol_A =
g^{-1} \ga g$ is an immediate generalization of Lemma
\ref{HolLem} (ii). For the reverse direction first consider a
representation $\ga$ and let $\widetilde \ga:= g^{-1}_0\ga g_0$
for some $g_0\in G$. Define
\[
\widetilde M\times G\to \widetilde M\times G,\quad (\widetilde
x,g) \mapsto (\widetilde x,g_0 g).
\]
This descend to a bundle map $P_\ga\to P_{\widetilde \ga}$ since
\[
\big(\widetilde x\cdot c, \widetilde \ga(c)^{-1} (g_0 g)\big) =
\big(\widetilde x\cdot c, g_0(\ga(c)^{-1}g)\big).
\]
One verifies that this gives an isomorphism of flat bundles. Now
assume that $\hol_A = \ga$, and fix a base point $p_0\in P$. For
$\widetilde x\in \widetilde M$, let $c_{\widetilde x}:I\to M$ be
a path in $M$ representing $\widetilde x$ and starting at
$x_0=\pi(p_0)$. Let $\widehat c_{\widetilde x}:I\to P$ be the
horizontal lift to $P$, starting at $p_0$. Using the same ideas
as in Proposition \ref{HolFlat} one finds that $\widehat
c_{\widetilde x}(1)$ depends only on the homotopy class of
$c_{\widetilde x}$. Hence, we get a well-defined map
\[
\gF: \widetilde M\times G\to P,\quad (\widetilde x,g)\mapsto
\widehat c_{\widetilde x}(1)\cdot g.
\]
The construction is in such a way that $\gF$ is $G$-equivariant
and surjective. Moreover, if $c\in \pi_1(M,x_0)$, then a
straightforward calculation shows that
\[
\widehat c_{(\widetilde x\cdot c)}(1) = \widehat c_{\widetilde
x}(1)\cdot\hol_A(c,p_0),
\]
and
\[
\gF^{-1}\big(\widehat c_{\widetilde x}(1)\big) =
\bigsetdef{\big(\widetilde x\cdot c, \hol_A(c,p_0)^{-1}\big)}{c\in
\pi_1(M,x_0)}.
\]
Since we are assuming that $\hol_A(c,p_0) = \ga(c)$, this implies
that $\gF$ descends to a bundle isomorphism $P_\ga\to P$.
Concerning the relation between the flat connections $A$ and
$A_\ga$ we remark without further comments that $\gF$ is defined
in such a way that it maps the horizontal distribution on the
trivial bundle $\widetilde M\times G$ to the horizontal
distribution on $P$ given by $A$. This implies that via the
isomorphism $P_\ga\cong P$ the connections $A_\ga$ and $A$ agree.
\end{proof}

\subsection[Flatness and Triviality of U$(k)$-Bundles]{Flatness and
Triviality of $\boldsymbol{\U(k)}$-Bundles}\label{FlatTrivSec}

As mentioned in Remark \ref{FlatTriv}, a flat principal
$G$-bundle is not necessarily trivializable. In this section we
give some details for the case $G=\U(k)$. Via the standard
representation of $\U(k)$ on $\C^k$, a principal $\U(k)$-bundle
$P$ defines a Hermitian vector bundle $E\to M$ and vice versa. We
will freely switch between the two equivalent notions.\\

\noindent\textbf{Flat Line Bundles.} We first recall that the
space of all Hermitian line bundles over a given manifold can be
described in terms of \v{C}ech cohomology, see \cite[Sec.
III.4]{W} for details. Let $M$ be a compact, connected manifold.
For a Lie group $G$ we denote by $\underline G$ the sheaf of
locally smooth functions on $M$ with values in $G$. Then
$H^1\big(M,\underline\U(1)\big)$ is isomorphic to the set of
Hermitian line bundles over $M$ up to isomorphism. Note that the
group structure on the latter is given by the tensor products of
line bundles. There is an exact sequence of sheaves.
\begin{equation}\label{U(1)CoeffSeq}
0\longrightarrow \Z\longrightarrow \underline
\R\xrightarrow{e^{2\pi ix}}\underline\U(1)\longrightarrow 0.
\end{equation}
Since the sheaf $\underline \R$ is fine (i.e., admits partitions
of unity), the cohomology $H^\bullet(M,\underline\R)$ vanishes
away from degree 0. The long exact sequence in cohomology then
produces natural isomorphisms
\[
H^p\big(M,\underline\U(1)\big) \cong H^{p+1}(M,\Z),\quad p\ge 1.
\]
For $p=1$, this isomorphism coincides with the integral first
Chern class
\begin{equation}\label{ChernCech}
c_1: H^1\big(M,\underline\U(1)\big) \xrightarrow{\cong} H^2(M,\Z).
\end{equation}

On the other hand, we have seen in Proposition
\ref{FlatModuliRepVar} that \emph{flat} Hermitian line bundles
are classified by representations of $\pi_1(M)$ in $\U(1)$. Here,
conjugation does not play a role here since $\U(1)$ is abelian.
Therefore, the moduli space of flat Hermitian line bundles has the
cohomological description
\[
\Hom\big(\pi_1(M),\U(1)\big) = H^1\big(M,\U(1)\big).
\]
Note that in terms of \v{C}ech cohomology, $H^p\big(M,\U(1)\big)$
refers to the sheaf of locally constant (rather than $C^\infty$)
functions with values in $\U(1)$. Similarly, we have to
distinguish between $H^\bullet(M,\R)$ and $H^\bullet(M,\underline
\R)$. We have a long exact coefficient sequence
\[
...\longrightarrow H^p(M,\Z)\longrightarrow H^p(M,\R)
\longrightarrow H^p\big(M,\U(1)\big)\longrightarrow
H^{p+1}(M,\Z)\longrightarrow...
\]
Here, the integral first Chern class appears again as the map
\[
c_1: H^1\big(M,\U(1)\big)\longrightarrow H^2(M,\Z).
\]
Moreover, the universal coefficient theorem shows that
\[
\ker\big(H^2(M,\Z)\to H^2(M,\R)\big) = \Tor\big(H^2(M,\Z)\big).
\]
Together with \eqref{ChernCech} this easily yields
\begin{lemma}\label{TorsionChern}
A line bundle $L\to M$ admits a flat connection if and only if
its integral first Chern class satisfies $c_1(L)\in
\Tor\big(H^2(M,\Z)\big)$.
\end{lemma}

\begin{remark*}
The above result also follows from Chern-Weil theory: The
representative of the first Chern class in de Rham cohomology is
given by $\big[\tr (\lfrac {i}{2\pi} F_A)\big]$, where $A$ is any
$\U(1)$-connection on $L$. Hence, $L$ admits a flat connection if
and only if $c_1(L)$ vanishes in $H^2(M,\R)$. This is precisely
the case if the integral first Chern class is a torsion class.
\end{remark*}

There is also a topological condition for triviality of a flat
Hermitian line bundle. It is straightforward to check that the
natural map of sheaves $\U(1)\to \underline \U(1)$ relates the
two exact coefficient sequences via
\[
\begin{CD}
0 @>>> H^1\big(M,\underline U(1)\big) @>{c_1}>> H^2(M,\Z)@>>>0\\
@AAA @AAA @| @AAA \\
H^1(M,\R) @>>> H^1\big(M,U(1)\big) @>{c_1}>> H^2(M,\Z)@>>>H^2(M,\R).\\
\end{CD}
\]

\begin{lemma}\label{TorsionChernTop}
Let $L_\ga\to M$ be a flat Hermitian line bundle on $M$ with
holonomy $\ga:H_1(M,\Z)\to \U(1)$. Then $L_\ga$ is trivializable
if and only if the restriction of $\ga$ to the torsion subgroup
$\Tor\big(H_1(M,\Z)\big)$ is trivial.
\end{lemma}

\begin{proof}
The line bundle $L_\ga$ is trivializable if and only if
$c_1(L_\ga)=0$ in $H^2(M,\Z)$. The above diagram shows that this
is precisely if
\[
L_\ga\in \im\big[H^1(M,\R)\to H^1\big(M,\U(1)\big)\big].
\]
The latter means that $\ga=\exp(2\pi i\widehat\ga)$ with
$\widehat\ga:H_1(M,\Z)\to \R$. Since $(\R,+)$ is a torsion free
abelian group, the restriction of $\widehat\ga$ to
$\Tor\big(H_1(M,\Z)\big)$ vanishes. Conversely, if $\ga$ is
trivial on $\Tor\big(H_1(M,\Z)\big)$, it can be lifted to a
homomorphism $\widehat\ga$ as above.
\end{proof}

\noindent\textbf{Flat $\boldsymbol{\U(k)}$-undles.} If we now
consider $\U(k)$ for $k>1$, Lemma \ref{TorsionChern} and Lemma
\ref{TorsionChernTop} do not generalize immediately, since we do
not have a simple cohomological description of the space of
(flat) $\U(k)$-bundles. This is mainly due to the fact that
$H^1\big(M,\underline \U(k)\big)$ is not a group since $\U(k)$ is
non-abelian. However, in the case that $\dim M\le 3$ the results
about flat $\U(1)$-bundles generalize to $\U(k)$. The underlying
reason is the following

\begin{prop}\label{SUBundleTriv}
Let $M$ be a manifold of dimension $\le 3$. Then every principal
$\SU(k)$-bundle $P\to M$ is trivializable.
\end{prop}

\begin{proof}[Idea of proof]
Let $E\SU(k)\to B\SU(k)$ be the universal $\SU(k)$-bundle, see
\cite[Sec. 8.6]{DavKir}. It has the property that $E\SU(k)$ is
contractible, and the set of isomorphism classes of
$\SU(k)$-bundles is isomorphic to the set of homotopy classes of
maps $M\to B\SU(k)$. Since the total space of the universal
$\SU(k)$-bundle is contractible, the long exact homotopy sequence
yields that
\begin{equation}\label{HomotopyGroups}
\pi_n\big(B\SU(k)\big) = \pi_{n-1}\big(\SU(k)\big),\quad n\ge 1.
\end{equation}
It is well known that $\SU(k)$ is 2-connected, see \cite[Sec.
6.14]{DavKir}. Hence, \eqref{HomotopyGroups} implies that
$B\SU(k)$ is 3-connected. Since we are assuming that $\dim M\le
3$ it follows that every map $M\to B\SU(k)$ is homotopic to a
constant map. Therefore, every $\SU(k)$-bundle over $M$ is
trivializable.
\end{proof}

We then have the following generalization of Lemma
\ref{TorsionChern} and Lemma \ref{TorsionChernTop}.

\begin{cor}
Let $M$ be a compact, connected manifold of dimension $\le 3$,
and let $E\to M$ be a Hermitian vector bundle over $M$ of rank
$k$.
\begin{enumerate}
\item The bundle $E$ admits a flat connection if and only if
its integral first Chern class satisfies $c_1(E)\in
\Tor\big(H^2(M,\Z)\big)$.
\item Assume that $E$ is flat with holonomy $\ga:\pi_1(M)\to
\U(k)$. Then $E$ is trivializable if and only if the restriction
of $\det(\ga):H^1(M,\Z)\to \U(1)$ to $\Tor\big(H_1(M,\Z)\big)$ is
trivial.
\end{enumerate}
\end{cor}

\begin{proof}
Let $\det(E):=\gL^k E\to M$ denote the determinant line bundle of
$E$. One concludes from the exact sequence
\[
0\longrightarrow \SU(k)\longrightarrow \U(k)\xrightarrow{\det}
\U(1)\longrightarrow 0,
\]
that $E\otimes \det(E)^{-1}$ admits an $\SU(k)$-structure. As we
are assuming $\dim M\le 3$, it follows from Proposition
\ref{SUBundleTriv} that $E\otimes \det(E)^{-1}$ is isomorphic to
$M\times \C^k$. This implies that $E$ is flat/trivializable if
and only if $\det(E)$ is flat/trivializable. Since
\[
c_1(E)=c_1\big(\det(E)\big)\in H^2(M,\Z),
\]
part (i) readily follows from Lemma \ref{TorsionChern}, whereas
(ii) is immediate from Lemma \ref{TorsionChernTop}.
\end{proof}

\begin{remark}\label{SurfaceBundleFlat}
If $\gS$ is a closed, oriented surface, then $H^2(\gS,\Z)\cong
\Z$. Therefore, there are no torsion elements and every flat
Hermitian vector bundle over $\gS$ is isomorphic to the trivial
bundle. For 3-manifolds there will in general exist non-trivial
flat Hermitian bundles, see for example Remark
\ref{MapTorTrivFlat} below or Section \ref{S1Bundles}.
\end{remark}

\section{Flat Connections over Mapping Tori}\label{MapTorFlatConn}

We let $M$ be a closed, oriented manifold, and let $f\in
\Diff^+(M)$ be an orientation preserving diffeomorphism. Then we
define the \emph{mapping torus} $M_f$ of $f$ as

\begin{figure}[htbp] \centering
\includegraphics[width=0.4\linewidth]{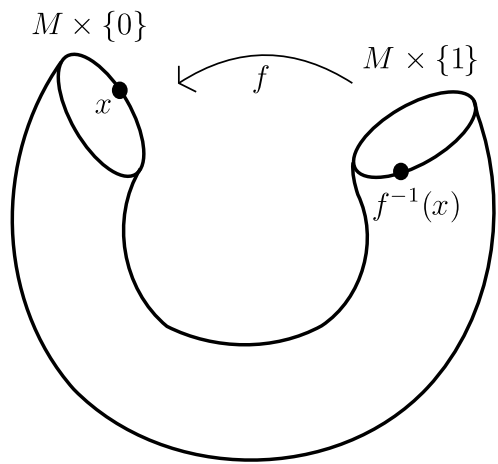}
\caption{The mapping torus of $f$}\label{Fig:MapTor}
\end{figure}
\[
M_f := \big(M\times \R\big)/\Z,
\]
where $\Z$ acts on $M\times \R$
via
\begin{equation}\label{MapTorDefConv} k\cdot(x,t) =
\big(f^{-k}(x),t+k\big),\quad (x,t)\in M\times \R,\quad k\in\Z.
\end{equation}
We use this definition rather than defining $M_f$ as a quotient
of $M\times [0,1]$, since it is convenient to work
$\Z$-equivariantly on $M\times \R$. Since we assume that $f$ is
orientation preserving, the product orientation of $M\times \R$
defines an orientation on $M_f$. The map
\[
M\times \R\to S^1,\quad (x,t)\mapsto \exp(2\pi i t)
\]
is invariant with respect to the action \eqref{MapTorDefConv},
and gives rise to a fiber bundle
\[
M\hookrightarrow M_f\xrightarrow{\pi} S^1.
\]

Let $\Diff_0(M)$ be the component of the identity in $\Diff^+(M)$
with respect to the $C^\infty$-topology. We collect the following
well-known material.

\begin{lemma}\label{MapTorIsotop}
The diffeomorphism class of $M_f$ depends only on the conjugacy
class of $f$ inside the mapping class group
\[
\Diff^+(M)/\Diff_0(M).
\]
Moreover, there exists a diffeomorphism $-M_f\cong M_{f^{-1}}$,
where $-M_f$ carries the reversed orientation.
\end{lemma}

\begin{proof}
To show that the diffeomorphism class of $M_f$ depends only on
the isotopy class of $f$, let $f_t:[0,1]\to \Diff^+(M)$ be an
isotopy. Possibly using a reparametrization of $[0,1]$ we may
assume that $f_t$ is constant near $0$ and $1$. Define $\gf_t:=
f_t^{-1}\circ f_0$ and extend $\gf_t$ to a path $\gf_t:\R\to
\Diff^+(M)$ by requiring that $\gf_{t+1} = f^{-1}_1\circ
\gf_t\circ f_0$. Then one easily checks that
\[
\gF:M\times \R\to M\times \R,\quad \gF(x,t):=
\big(\gf_t(x),t\big),
\]
is $\Z$-equivariant with respect to \eqref{MapTorDefConv}. Hence,
it descends to a diffeomorphism $\gF:M_{f_0}\rightarrow M_{f_1}$.
Similarly, if $g\in \Diff^+(M)$, then the map
\[
\Psi:M\times \R\to M\times \R,\quad \Psi(x,t):= \big(g(x),t\big)
\]
defines a diffeomorphism $\Psi:M_f\rightarrow M_{gfg^{-1}}$. Thus,
the diffeomorphism class of $M_f$ depends only on the conjugacy
class of $f$ in $\Diff^+(M)/\Diff_0(M)$. In a similar way one
verifies that an orientation reversing diffeomorphism $M_f\cong
M_{f^{-1}}$ is induced by
\[
M\times \R\to M\times \R,\quad (x,t)\mapsto (x,-t).\qedhere
\]
\end{proof}

\begin{remark}\label{MapTorFiberBundRem}
We also want to point out that every oriented fiber bundle
$M\hookrightarrow \widetilde M\xrightarrow{\pi} S^1$ arises in
the way just described: Identify $M$ with the fiber
$\pi^{-1}(1)$, and endow $M$ with a vertical projection
$P^v:TM\to T^vM$, see Section \ref{FiberedConn}. For $x\in M$
denote by $c_x:[0,1]\to M$ the horizontal lift of the path
\[
[0,1]\to S^1,\quad t\mapsto e^{2\pi i t}
\]
which starts at $x$. Then $c_x(1)\in M$, so that we can define
\[
f:M\to M,\quad f(x):= c_x(1).
\]
It follows from the standard theory of ordinary differential
equations that $f$ is a diffeomorphism. Moreover, since we have
assumed that $M\hookrightarrow \widetilde M\xrightarrow{\pi} S^1$
is an oriented fiber bundle, it follows that $f$ is orientation
preserving, i.e., $f\in \Diff^+(M)$. One then verifies that the
fiber bundle $\pi:\widetilde M\to S^1$ is indeed isomorphic to the
fiber bundle given by the mapping torus $M_f$. This
identification might seem to depend on the choice of $P^v$ and
the identification of $M$ as a fiber. However, since all vertical
projections are homotopic, one can check that the conjugacy class
of $f$ in the mapping class group $\Diff^+(M)/\Diff_0(M)$ does
not change, so that by Lemma \ref{MapTorIsotop} the isomorphism
class of the mapping torus is unambiguously defined.
\end{remark}

\subsection{Algebraic Description of the Moduli Space}

Let $G$ be a connected matrix Lie group. Ultimately $G$ will be
$\U(k)$. We use Proposition \ref{FlatModuliRepVar} to identify
\[
\cM(M_f,G)= \Hom\big(\pi_1(M_f),G\big)/G,
\]
where $G$ acts by conjugation. We fix a base point $x_0\in M$
assume for simplicity that $f(x_0)=x_0$. This is possible as we
can always find an element in the isotopy class of $f\in
\Diff^+(M)$ which fixes $x_0$. Then the path
\[
\gamma:\R\to M\times \R,\quad t\mapsto (x_0,t),
\]
descends to a closed path in $M_f$, whose homotopy class we also
denote by $\gamma\in \pi_1(M_f)$. Without including it in the
notation, we are using $[x_0,0]$ as a base point for $M_f$. On the
other hand, the inclusion $M=M\times \{0\} \subset M\times \R$
induces a map
\[
i_*:\pi_1(M)\to \pi_1(M_f).
\]
Then the following is easily verified.\footnote{Recall that for
fundamental groups we are using the product $[c]\cdot [\widetilde
c] := [\widetilde c * c]$, where $\widetilde c * c$ means first
$\widetilde c$ and then $c$.}

\begin{lemma}\label{MapTorFund}
The fundamental group of $M_f$ with respect to the base point
$[x_0,0]$ is given by
\[
\pi_1(M_f)= \big\langle \pi_1(M),\gamma\,\big|\,  \gamma^{-1}
c\gamma = f_*c,\;c\in \pi_1(M)\big\rangle,
\]
where $f_*:\pi_1(M)\to \pi_1(M)$ is the induced map on the
fundamental group.
\end{lemma}

The map $i_*:\pi_1(M)\to \pi_1(M_f)$ gives rise to a natural
homomorphism
\[
i^*: \Hom\big(\pi_1(M_f),G\big)\to \Hom\big(\pi_1(M),G\big),\quad
\ga\mapsto \ga\circ i_*.
\]
This map is $G$-equivariant and we get an induced map
\begin{equation}\label{FibrationAlg}
[i^*]: \cM(M_f,G)\to \cM(M,G),\quad [\ga]\mapsto [\ga\circ i_*].
\end{equation}
Similarly, $f^*:\Hom\big(\pi_1(M),G\big)\to
\Hom\big(\pi_1(M),G\big) $ is $G$-equivariant so that it descends
to a map
\[
[f^*]:\cM(M,G) \to \cM(M,G),\quad [\ga]\mapsto [\ga\circ f_*].
\]
For the following result see also \cite[Sec. 8]{And}.

\begin{prop}\label{MapTorModuliAlg}
The natural map \eqref{FibrationAlg} defines a surjection
\[
[i^*]: \cM(M_f,G)\to \Fix [f^*]\subset \cM(M,G).
\]
Moreover, if $\ga\in \Hom\big(\pi_1(M),G\big)$ is such that
$[\ga]\in \Fix [f^*]$, then
\begin{equation}\label{FiberAlg}
[i^*]^{-1}[\ga] \cong \bigsetdef{g \in G}{g^{-1}\ga g =
f^*\ga}/Z(\ga)
\end{equation}
where $Z(\ga):=\setdef{h\in G}{h^{-1}\ga h=\ga}$ is the
centralizer of $\ga$, and acts on $G$ by conjugation.
\end{prop}

\begin{proof}
With $\gamma$ as in Lemma \ref{MapTorFund} define
\[
\gF: \Hom\big(\pi_1(M_f),G\big) \to
\Hom\big(\pi_1(M),G\big)\times G,\quad \ga\mapsto \big(\ga\circ
i_*,\ga(\gamma) \big).
\]
Then Lemma \ref{MapTorFund} implies that $\gF$ induces an
isomorphism
\begin{equation}\label{MapTorModuliAlg:1}
\gF:\Hom\big(\pi_1(M_f),G\big)\xrightarrow{\cong}
\bigsetdef{(\ga,g)\in \Hom\big(\pi_1(M),G\big)\times G}{g^{-1}\ga
g = f^*\ga}.
\end{equation}
The action of $G$ by conjugation on $\Hom\big(\pi_1(M_f),G\big)$
translates to to the right hand side as $G$ acting diagonally by
conjugation, i.e., $h\in G$ acts on an element $(\ga,g)$ via
\[
(\ga,g)\cdot h := \big(h^{-1}\ga h, h^{-1}g h\big).
\]
Then
\[
\cM(M_f,G)\cong \bigsetdef{(\ga,g)\in
\Hom\big(\pi_1(M),G\big)\times G}{g^{-1}\ga g = f^*\ga}/G.
\]
Under this identification, the map $[i^*]$ in \eqref{FibrationAlg}
is given by the projection onto the first factor. Now $[\ga] =
[f^*\ga]$ if and only if there exists $g\in G$ with $g^{-1}\ga g =
f^*\ga$, which guarantees that
\[
[i^*]: \cM(M_f,G) \to \Fix[f^*]\subset \cM(M,G), \quad
[\ga,g]\mapsto [\ga],
\]
is well-defined and surjective. The inverse image of $[\ga]\in
\Fix[f^*]$ is given by
\[
[i^*]^{-1}[\ga]= \bigsetdef{(\widetilde \ga,g)\in
\Hom\big(\pi_1(M),G\big)\times G}{g^{-1}\widetilde \ga g =
f^*\widetilde \ga,\, [\widetilde \ga]=[\ga]}/G.
\]
Hence, if we fix a representative $\ga$, then
\[
[i^*]^{-1}[\ga]\cong  \bigsetdef{g \in G}{g^{-1}\ga g =
f^*\ga}/Z(\ga).\qedhere
\]
\end{proof}

\begin{remark*}\quad\nopagebreak
\begin{enumerate}
\item If we fix a different representative $\widetilde \ga = h^{-1}
\ga h$ in \eqref{FiberAlg}, then one easily checks that
\[
\bigsetdef{\widetilde g \in G}{\widetilde g^{-1}\widetilde \ga
\widetilde g = f^*\widetilde \ga} = \bigsetdef{h^{-1} g h \in
G}{g^{-1}\ga g = f^*\ga}.
\]
Moreover,
\[
Z(\widetilde \ga) = h^{-1} Z(\ga) h.
\]
This describes the relation among the isomorphisms in
\eqref{FiberAlg} for different choices of representatives for
$[\ga]\in \Fix[f^*]$.
\item We also want to point out that if we fix $g$ with $g^{-1} \ga g =
f^*\ga$, then
\[
\bigsetdef{g \in G}{g^{-1}\ga g = f^*\ga} = Z(\ga)g.
\]
Furthermore, note that if $\ga\in \Hom\big(\pi_1(M),G\big)$ is
irreducible, then $Z(\ga)$ coincides with the center of $G$,
\[
Z(\ga)=Z(G)=\bigsetdef{g\in G}{gh=hg\quad\text{for all $h\in G$}}.
\]
\end{enumerate}
\end{remark*}

\subsection{Geometric Description of the Moduli Space}

We now turn to a more geometric description of $\cM(M_f,G)$. A
related discussion, yet in a different context, can be found in
\cite[Sec. 5]{DosSal} and \cite[Sec. 7]{Sal}. Since we prefer to
work with Hermitian vector bundles over $M_f$, we assume from now
on that $G=\U(k)$.\\

\noindent\textbf{Vector Bundles over Mapping Tori.} Using the
ideas of Remark \ref{MapTorFiberBundRem} one finds that every
Hermitian vector bundle over $M_f$ is isomorphic to a mapping
torus
\begin{equation}\label{GaugeBundleDef}
E_{\widehat f}:=\big(E\times \R\big)/\Z \to \big(M\times
\R\big)/\Z= M_f.
\end{equation}
Here, $\widehat f$ is a bundle isomorphism covering $f$, and $\Z$
acts on $E\times \R$ via
\[
k\cdot(e,t)= \big(\widehat f^{-k}(e),t+k\big),\quad (e,t)\in
E\times \R,\quad k\in \Z.
\]
To keep the notation simple, we are identifying $E\times \R$ with
the pullback bundle $\pi_M^*E\to M\times\R$, where $\pi_M:M\times
\R\to M$. For a fixed Hermitian vector bundle $E\to M$, we denote
by $\cG(E)$ the group of gauge transformations and by $\cG_f(E)$
the set of isomorphism classes of bundle isomorphisms covering
$f$. Note that there is a free and transitive action of $\cG(E)$
on $\cG_f(E)$, given by
\[
\cG_f(E)\times \cG(E)\to \cG_f(E),\quad (\widehat f,u)\mapsto
\widehat f\circ u.
\]
Hence, upon fixing one particular $\widehat f\in \cG_f(E)$, the
space $\cG_f(E)$ is isomorphic to the group of gauge
transformations $\cG(E)$.

\begin{remark}\label{MapTorBundleTriv}
If $E=M\times \C^k$ is the trivial bundle, the group $\cG(E)$
coincides with $C^\infty\big(M,\U(k)\big)$, see Remark
\ref{GaugeTrivBundle}. Then for $u\in C^\infty\big(M,\U(k)\big)$,
we define
\[
\widehat f_u: M\times \C^k\to M\times \C^k,\quad \widehat
f_u(x,z):= \big(f(x),u(x)z\big).
\]
This gives a canonical identification
\[
C^\infty\big(M,\U(k)\big)\xrightarrow{\cong} \cG_f(M\times \C^k)
,\quad u\mapsto \widehat f_u.
\]
We then write $E_u$ for the bundle defined by $\widehat f_u$.
\end{remark}

\begin{lemma}\label{MapTorBundleIsom}
If $\widehat f_1,\widehat f_2 :E\to E$ are two bundle
isomorphisms covering $f$, then $E_{\widehat f_1}\cong
E_{\widehat f_2}$ if and only if there exists $\gf_t\in
C^\infty\big(\R,\cG(E)\big)$ such that
\begin{equation}\label{MapTorBundleIsomCond}
\gf_{t+1}= \widehat f_2^{-1}\circ \gf_t \circ \widehat f_1.
\end{equation}
\end{lemma}

\begin{proof}
If we define
\[
\gF:E\times \R \to E\times \R,\quad
\gF(e,t):=\big(\gf_t(e),t\big),
\]
then
\[
\gF\big(\widehat f_1^{-1}(e),t+1\big) = \big(\gf_{t+1}\circ
\widehat f_1^{-1}(e),t+1\big)= \big(\widehat f_2^{-1}\circ
\gf_t(e),t+1\big).
\]
This implies that $\gF$ is a $\Z$-equivariant bundle isomorphism
and so $E_{\widehat f_1}$ and $E_{\widehat f_2}$ are isomorphic.
Conversely, we can lift a bundle isomorphism $E_{\widehat
f_1}\cong E_{\widehat f_2}$ to a $\Z$-equivariant map $E\times \R
\to E\times \R$ and use this to define $\gf_t$ with the required
property.
\end{proof}

\begin{remark}\label{MapTorBundleIsomTriv}
Assume that $E=M\times \C^k$ is a trivial bundle, and let $u,v \in
C^\infty\big(M,\U(k)\big)$ define vector bundles $E_u$ and $E_v$
as in Remark \ref{MapTorBundleTriv}. Then
\eqref{MapTorBundleIsomCond} takes the form
\begin{equation*}
\gf_{t+1}= \widehat f_v^{-1}\circ \gf_t \circ \widehat f_u =
v^{-1}(\gf_t\circ f) u.
\end{equation*}
\end{remark}

\noindent\textbf{Additional Remarks for Line Bundles.} In the case
that $E=L$ is a line bundle, we can also give a more topological
interpretation of \eqref{MapTorBundleIsomCond}. Since $\U(1)$ is
abelian, the group of gauge transformations of $L$ is given by
$\cG(L)=C^\infty\big(M,\U(1)\big)$, irrespectively of whether $L$
is trivializable or not. From the long exact sequence associated
to the coefficient sequence \eqref{U(1)CoeffSeq}, we get an exact
sequence
\begin{equation}\label{ExactCoeffSeq}
0\longrightarrow \Z\longrightarrow C^\infty(M,\R)\longrightarrow
C^\infty\big(M,\U(1)\big)\longrightarrow H^1(M,\Z)\longrightarrow
0.
\end{equation}
Here, the map $C^\infty\big(M,\U(1)\big)\to H^1(M,\Z)$ is given by
\[
u\mapsto u_*\in \Hom\big(\pi_1(M),\Z\big) = H^1(M,\Z).
\]

\begin{remark}\label{IntPeriod}
The de Rham theorem defines a map
\[
\gO^p_{cl}(M,\R)\to H^p(M,\R),
\]
where the left hand side denotes the space of closed $p$-forms.
Since $H^1(M,\Z)$ has no torsion, we can use this to identify
\[
H^1(M,\Z) \cong \frac{\gO^1_{cl}(M,\Z)}{d C^\infty(M,\R)}.
\]
Here, $\gO^p_{cl}(M,\Z)$ denotes the space of closed $p$-forms
with \emph{integral periods}, i.e., the kernel of the projection
\[
\gO^p_{cl}(M,\R)\to H^p(M,\R)/H^p(M,\Z).
\]
\end{remark}

Using this, the last map in \eqref{ExactCoeffSeq} can be
expressed as
\begin{equation}\label{GaugeCohomClass}
C^\infty\big(M,\U(1)\big)\to H^1(M,\Z),\quad u\mapsto
\Big[\frac{u^{-1}du}{2\pi i}\Big].
\end{equation}
Now, $u\in C^\infty\big(M,\U(1)\big)$ is mapped to 0 if and only
if $u=\exp(2\pi i g)$ for some $g\in C^\infty(M,\R)$, because then
\[
u^{-1}du = 2\pi i dg.
\]
Also note that there exists $g$ with $u=\exp(2\pi i g)$ precisely
if $u$ is homotopic to a constant map. Then we have the following
topological interpretation of \eqref{MapTorBundleIsomCond}.

\begin{prop}\label{MapTorIsomCond}
Let $L\to M$ be a Hermitian line bundle, and let $\widehat
f_1,\widehat f_2\in \cG_f(L)$. Define $u\in \cG(L)$ by requiring
that $\widehat f_2 = \widehat f_1\circ u$. Then the line bundles
$L_{\widehat f_1}$ and $L_{\widehat f_2}$ are isomorphic if and
only if
\begin{equation}\label{IsomCond:1}
\Big[\frac{ u^{-1}du }{2\pi i}\Big] \in
\im\big(\Id-f^*\big)\subset H^1(M,\Z).
\end{equation}
\end{prop}

\begin{proof}
Assume first that $L_{\widehat f_1}\cong L_{\widehat f_2}$. Then
there exists $\gf_t\in C^\infty\big(\R,\cG(L)\big)$ as in Lemma
\ref{MapTorBundleIsom} such that $\gf_{t+1}= \widehat
f_2^{-1}\circ \gf_t \circ \widehat f_1$. Now, since $\U(1)$ is
abelian, $\widehat f_1^{-1}\circ \gf_t \circ \widehat
f_1=\gf_t\circ f$. Hence,
\[
\gf_{t+1}= u^{-1}\circ \widehat f_1^{-1}\circ \gf_t\circ \widehat
f_1 = u^{-1} (\gf_t\circ f).
\]
This yields
\begin{equation}\label{IsomCond:2}
\gf_1^{-1}d\gf_1 = f^*(\gf_0^{-1}d\gf_0) - u^{-1}du
\end{equation}
Moreover, since $\gf_0$ is homotopic to $\gf_1$ we can find
$g:M\to \R$ such that
\[
 \gf_0^{-1}d\gf_0 =  \gf_1^{-1}d\gf_1  + 2\pi i dg .
\]
From this and \eqref{IsomCond:2} we see that \eqref{IsomCond:1}
is valid. Conversely, let us assume that \eqref{IsomCond:1} holds.
This means that there exist $\gf\in C^\infty\big(M,\U(1)\big)$ and
$g\in C^\infty(M,\R)$ such that
\[
u^{-1}du  =f^*(\gf^{-1} d \gf) - \gf^{-1} d \gf  + 2\pi i dg.
\]
For $t\in [0,1]$ define $\gf_t:= \gf \exp\big(-2\pi it g\big)$.
Then, upon adding a constant to $g$, we conclude that
\[
u = \gf_1^{-1}(\gf_0\circ f).
\]
This allows us to extend $\gf_t$ for all $t\in \R$ in such a way
that
\[
\gf_{t+1}= u^{-1}(\gf_t\circ f)= \widehat f_2^{-1}\circ \gf_t
\circ \widehat f_1,
\]
and we get the required isomorphism.
\end{proof}

\begin{remark}\label{MapTorTrivFlat}\quad\nopagebreak
\begin{enumerate}
\item As we have mentioned in Section \ref{FlatTrivSec}, line bundles
over $M_f$ are classified by the group $H^2(M_f,\Z)$. This group
fits into an exact sequence associated to the fiber bundle
$M\hookrightarrow M_f\to S^1$. One can use for example the
five-term exact sequence induced by the Leray-Serre spectral
sequence to obtain an exact sequence
\begin{equation}\label{WangSeq}
\begin{split}
\ldots\longrightarrow H^1(M) \xrightarrow{\Id-f^*} &H^1(M)
\longrightarrow H^2(M_f)\\
&\xrightarrow{i^*}H^2(M)\xrightarrow{\Id-f^*}H^2(M)\longrightarrow
\ldots
\end{split}
\end{equation}
For higher dimensional spheres as the base this is usually
referred as the \emph{Wang sequence}, see \cite[p. 254]{DavKir}.

Using the discussion preceding Proposition \ref{MapTorIsomCond}
we can give a geometric interpretation of the map $H^1(M) \to
H^2(M_f)$ in \eqref{WangSeq}. First we use \eqref{GaugeCohomClass}
to represent an element of $H^1(M,\Z)$ by a gauge transformation
$u\in C^\infty\big(M,\U(1)\big)$, and let $L_u\to M_f$ be the line
bundle defined by $u$ as in Remark \ref{MapTorBundleTriv}. Since
the isomorphism class of $L_u$ is independent of the homotopy
class of $u$ we get a well-defined map
\[
H^1(M,\Z) \to H^2(M_f,\Z),\quad \Big[\frac{u^{-1}du}{2\pi
i}\Big]\mapsto c_1(L_u).
\]
Now Proposition \ref{MapTorIsomCond} identifies the kernel of this
map and gives a geometric explanation for the exactness of the
\eqref{WangSeq} at $H^1(M)$. Also, from a geometric point of view,
exactness at $H^2(M_f,\Z)$ is immediate; this simply means that a
line bundle $L\to M_f$ is of the form $L=L_u$ if and only if it
restricts to the trivial line bundle over $M$. Concerning
exactness at $H^2(M)$ we observe that the restriction $L|_M$ of a
line bundle $L\to M_f$ has to satisfy $f^*c_1(L|_M)= c_1(L|_M)$.
Conversely, this is precisely the condition which enables us to
define a line bundle over $M_f$.
\item The sequence \eqref{WangSeq} also shows that there will in
general be flat line bundles over $M_f$ which are topologically
non-trivial. A gauge transformation $u\in
C^\infty\big(M,\U(1)\big)$ gives a flat bundle $L_u$ if and only
if $c_1(L_u)$ is a torsion element in $H^2(M_f,\Z)$. According to
Proposition \ref{MapTorIsomCond} and part (i) of this remark, this
is precisely the case if there exists $N\in \N$ such that
\[
N\cdot\Big[\frac{u^{-1}d u}{2\pi i}\Big]= \Big[\frac{u^{-N}d
u^N}{2\pi i}\Big] \in \im(\Id- f^*)\subset H^1(M,\Z).
\]
Already in the case that $M$ is the 2-torus $T^2$ one can easily
find a diffeomorphism $f:T^2\to T^2$ such that
\[
\coker\big (\Id - f^*:H^1(M,\Z)\to H^1(M,\Z)\big)
\]
contains torsion elements, see Remark \ref{TorusBundleFlatConnRem}
in Section \ref{TorusBundlesGen}.
\end{enumerate}
\end{remark}

\noindent\textbf{Flat Connections over Mapping Tori.} For
convenience we assume for the rest of this section that a flat
$\U(k)$-bundle over $M$ is necessarily trivializable. According
to Remark \ref{SurfaceBundleFlat} this is satisfied for example
if $M=\gS$ is a closed surface, which is the case we are
considering in Chapter \ref{3dimMapTor}.

Under this assumption a flat bundle $E\to M_f$ restricts to the
trivial bundle over the fiber, so that it is of the form
considered in Remark \ref{MapTorBundleTriv}, i.e.,
\[
E=E_u\quad u\in C^\infty\big(M,\U(k)\big),
\]
where $E_u$ is the mapping torus of the bundle isomorphism
$\widehat f_u:M\times \C^k\to M\times \C^k$. We can then identify
the space of sections of $E_u$ as
\[
C^\infty(M_f,E_u)= \bigsetdef{\gf_t:\R\to
C^\infty(M,\C^k)}{\gf_{t+1}= u^{-1}(\gf_t\circ f)}.
\]
More generally, $\widehat f_u$ induces a pullback
\begin{equation}\label{BundleIsomPullback}
\widehat f_u^* \ga = u^{-1}(f^*\ga),\quad \ga\in
\gO^\bullet(M,\C^k),
\end{equation}
and
\begin{equation}\label{FormsMapTor}
\gO^\bullet(M_f,E_u)= \bigsetdef{\ga_t\in \gO^1(M\times
\R,\C^k)}{\ga_{t+1}= \widehat f_u^*\ga_t }.
\end{equation}
In a similar way, a $\U(k)$-connection $A$ on $E_u$ can
equivalently be described as a Lie algebra valued 1-form
\[
A= a_t + b_tdt,\quad a_t\in
C^\infty\big(\R,\gO^1(M,\cu(k)\big),\quad  b_t \in
C^\infty\big(\R,C^\infty(M,\cu(k)\big),
\]
which is $\Z$-equivariant in the sense that
\begin{equation}\label{MapTorConnDef}
a_{t+1} = \widehat f_u^* a_t = u^{-1} (f^*a_t) u + u^{-1}du,\quad
b_{t+1} = \widehat f_u^* b_t = u^{-1} (b_t\circ f) u.
\end{equation}
The curvature of $A$ is given by
\[
F_{A} = d_{M\times \R} A + A\wedge A= d_M a_t + a_t\wedge a_t
+\big(d_M b_t - \pd_t a_t + [a_t,b_t]\big)\wedge dt.
\]
Hence, if $F_{a_t}\in C^\infty\big(\R,\gO^2(M,\cu(k)\big)$ denotes
the curvature associated to the path of connections $a_t$ over
$M$, we have
\[
F_{A} = F_{a_t} + \big(d_{a_t}b_t-\pd_t a_t\big)\wedge dt.
\]
Therefore, $A$ is flat if and only if
\begin{equation}\label{MapTorFlatDef}
F_{a_t} = 0,\quad\text{and}\quad \pd_t a_t=d_{a_t}b_t.
\end{equation}

Before we can describe the structure of the moduli space of flat
$\U(k)$-bundles on $M_f$ we need the following technical result.

\begin{lemma}\label{ModuliSpace:1}
Let $A=a_t+b_tdt$ be a flat connection on $E_u$. Then there
exists a gauge transformation $v\in C^\infty\big(M,\U(k)\big)$
such that the constant path $a_0$ defines a connection $A_0$ on
$E_v$, and there exists an isomorphism
\[
\gF:E_v\to E_u,\quad\text{with}\quad \gF^* A = A_0.
\]
\end{lemma}

\begin{proof}
Let $\gf_t:\R\to C^\infty\big(M,\U(k)\big)$ be the unique
solution to
\[
\pd_t \gf_t = - b_t \gf_t,\quad \gf_0\equiv e,
\]
where $e\in \U(k)$ is the identity matrix. Define $v:=u\gf_1\in
C^\infty\big(M,\U(k)\big)$. We claim that
\begin{equation}\label{ModuliSpace:2}
\gf_{t+1} = u^{-1} (\gf_t\circ f) v\quad\text{and}\quad
\pd_t(a_t\cdot \gf_t) = 0.
\end{equation}
\begin{proof}[Proof of \eqref{ModuliSpace:2}]
For the first assertion we recall from \eqref{MapTorConnDef} that
$b_t$ satisfies
\[
b_{t+1} = u^{-1}(b_t\circ f) u.
\]
From this and the definition of $\gf_t$, one easily checks that
both sides of the claimed equality are solutions of the initial
value problem
\[
\pd_t \widetilde \gf_t = - u^{-1} (b_t\circ f) u\,  \widetilde
\gf_t,\quad \widetilde \gf_0 = u^{-1}v,
\]
and hence agree. Now, using the identity $\pd_t \gf_t^{-1} =
-\gf_t^{-1}(\pd_t \gf_t)\gf_t^{-1}$, we compute that
\[
\begin{split}
\pd_t(a_t\cdot \gf_t) &= \pd_t(\gf_t^{-1}a_t\gf_t + \gf_t^{-1}d_M \gf_t)\\
&= -\gf_t^{-1}(\pd_t\gf_t) \gf_t^{-1}a_t \gf_t + \gf_t^{-1}(\pd_t
a_t)\gf_t +
\gf_t^{-1}a_t(\pd_t\gf_t)\\
&\qquad\qquad - \gf_t^{-1}(\pd_t\gf_t) \gf_t^{-1} d_M \gf_t +
\gf_t^{-1} d_M (\pd_t\gf_t)\\ &= \gf_t^{-1}(\pd_t a_t) \gf_t +
\gf_t^{-1}(b_ta_t)\gf_t - \gf_t^{-1}(a_tb_t)\gf_t\\&\qquad\qquad +
\gf_t^{-1} b_t d_M \gf_t - \gf_t^{-1}d_M(b_t\gf_t)\\
&= \gf_t^{-1}\big(\pd_ta_t - d_M b_t - [a_t,b_t]\big)\gf_t =
\gf_t^{-1}\big(\pd_ta_t - d_{a_t} b_t\big)\gf_t.
\end{split}
\]
Since we are assuming that $A$ is flat, condition
\eqref{MapTorFlatDef} yields the second part of
\eqref{ModuliSpace:2}.
\end{proof}

Having established \eqref{ModuliSpace:2} we continue with the
proof of Lemma \ref{ModuliSpace:1}. First of all we can use the
first formula in \eqref{ModuliSpace:2} and Remark
\ref{MapTorBundleIsomTriv}, to deduce that the family $\gf_t$
defines a bundle isomorphism
\[
\gF:E_v\xrightarrow{\cong} E_u.
\]
On the other hand, we can use \eqref{MapTorConnDef} and the
definition of $\gf_1$ to compute that
\[
f^* a_0 = a_1\cdot u^{-1} = (a_1\cdot \gf_1)\cdot v^{-1} =
a_0\cdot v^{-1},
\]
where in the last step we have used \eqref{ModuliSpace:2} to see
that $a_t\cdot \gf_t$ is constant. Also note that $a_0\cdot \gf_0
= a_0\cdot e = a_0$. This implies that the constant path $a_0$
defines a connection $A_0$ on $E_v$. Moreover, $\gF^*A$ is given
by
\begin{equation}\label{MapTorConnTrans}
\gF^*A = a_t\cdot \gf_t + \big(\gf_t^{-1}b_t \gf_t +
\gf_t^{-1}\pd_t\gf_t\big)dt,
\end{equation}
and thus, $\gF^*A=a_t\cdot\gf_t \equiv A_0$.
\end{proof}

\noindent\textbf{The Moduli Space of Flat Connections.} We can now
state geometric version of Proposition \ref{MapTorModuliAlg}.
Recall that we are assuming that a flat bundle over $M$ is
necessarily trivial. We let $\cF_M$ denote the space of flat
$\U(k)$-connections over $M$, and define
\[
\widehat \cM(M_f) := \bigsetdef{(a,v)\in \cF_M\times
C^\infty\big(M,\U(k)\big)}{\widehat f_v^*a = a }.
\]
Here, $\widehat f_v^*a$ is defined as in
\eqref{BundleIsomPullback}. Note the similarity of the definition
of $\widehat \cM(M_f)$ and
\[
\bigsetdef{(\ga,g)\in \Hom\big(\pi_1(M),\U(k)\big)\times
\U(k)}{g^{-1}\ga g = f^*\ga},
\]
which we considered in \eqref{MapTorModuliAlg:1} in the context of
the algebraic description of $\cM\big(M_f,\U(k)\big)$. There is a
natural action of $C^\infty\big(M,\U(k)\big)$ on $\widehat
\cM(M_f)$, given by
\[
(a,v)\cdot \gf := \big(a\cdot\gf,(\gf^{-1}\circ f)v\gf\big),\quad
(a,v)\in \widehat \cM(M_f),\quad \gf\in C^\infty\big(M,\U(k)\big).
\]
This action is well-defined since for $u:= (\gf^{-1}\circ f)v\gf$
we have
\[
\widehat f_u^* (a\cdot \gf) = (\widehat f_v^*a)\cdot
v^{-1}u(\gf\circ f) = a\cdot \gf.
\]
For notational brevity we use the following abbreviations for the
moduli spaces of flat bundles
\[
\cM(M):=\cM\big(M,\U(k)\big),\quad
\cM(M_f):=\cM\big(M_f,\U(k)\big).
\]
We then have the following analog of of Proposition
\ref{MapTorModuliAlg}.

\begin{prop}\label{MapTorModuliGeom}
With respect to the natural action of $C^\infty\big(M,\U(k)\big)$
on $\widehat \cM(M_f)$ we have
\[
\widehat \cM(M_f)\big/C^\infty\big(M,\U(k)\big) \cong \cM(M_f).
\]
Moreover, the projection $\widehat \cM(M_f) \to \cF_M$ onto the
first factor induces a surjection
\[
[i^*]: \cM(M_f) \to \Fix(f^*)\subset \cM(M).
\]
If we represent $[a]\in \Fix(f^*)\subset \cM(M)$ by $a\in \cF_M$,
then
\[
[i^*]^{-1}[a] \cong \bigsetdef{v\in
C^\infty\big(M,\U(k)\big)}{\widehat f_v^*a=a}\big/I(a),
\]
where $I(a)=\setdef{g\in \U(k)}{g^{-1}ag = a}$ denotes the
isotropy group of $a$ in $\U(k)$.
\end{prop}

\begin{proof}
Every element $(a,v)\in \cM(M_f)$ defines a flat $\U(k)$-bundle
$(E_v,A_0)$ over $M_f$. Here, $A_0$ is used as in Lemma
\ref{ModuliSpace:1} to denote the flat connection on $E_v$
induced by the constant path $a:\R\to \gO^1\big(M,\cu(k)\big)$.
We define
\[
\Psi:\widehat \cM(M_f) \to \cM(M_f),\quad (a,v)\mapsto [E_v,A_0].
\]
If $(a,v)\in \widehat \cM(M_f)$ and $\gf\in
C^\infty\big(M,\U(k)\big)$, then
\[
\gf  = v^{-1}(\gf\circ f) u,\quad \text{where}\quad
u:=(\gf^{-1}\circ f)v\gf .
\]
According to Lemma \ref{MapTorBundleIsom} this implies that
$E_u\cong E_v$. Moreover, the connection $A_0$ on $E_v$ pulls back
to the connection $A_0\cdot \gf$ on $E_u$, see
\eqref{MapTorConnTrans}. Thus,
\[
\Psi\big((a,v)\cdot \gf\big) = \Psi(a,v),
\]
and we get a well-defined map
\[
\overline \Psi: \widehat \cM(M_f)\big/C^\infty\big(M,\U(k)\big)
\to \cM(M_f).
\]
To see that $\overline \Psi$ is surjective, let $[E_u,A]\in
\cM(M_f)$. Then Lemma \ref{ModuliSpace:1} implies that there exist
$a\in \cF_M$ and $v\in C^\infty\big(M,\U(k)\big)$ in such a way
that $E_u$ is isomorphic to $E_v$, and the connection $A$ pulls
back to the connection $A_0$ defined by $a$. By definition, this
gives an element $(a,v)\in \widehat \cM(M_f)$ such that
$\Psi(a,v) = [E_u,A]$. To check injectivity, assume that
$\Psi(a,v) = \Psi(\widetilde a,\widetilde v)$. According to Lemma
\ref{MapTorBundleIsom} and \eqref{MapTorConnTrans} this means
that there exists $\gf_t:\R\to  C^\infty\big(M,\U(k)\big)$ such
that
\[
\gf_{t+1}= v^{-1} (\gf_t\circ f) \widetilde v,\quad a\cdot \gf_t
= \widetilde a, \quad\text{and}\quad \gf_t^{-1}\pd_t\gf_t = 0.
\]
In particular, $\gf_t\equiv \gf$ is independent of $t$, and
\[
\widetilde a = a\cdot \gf,\quad\text{and}\quad \widetilde v =
(\gf^{-1}\circ f) v \gf.
\]
Hence, $(\widetilde a,\widetilde v) = (a,v)\cdot \gf$, which
establishes injectivity. The rest of the proof is formally the
same as the proof of Proposition \ref{MapTorModuliAlg} and shall
be omitted.
\end{proof}

\section{Holomorphic Line Bundles over Riemann
Surfaces.}\label{HolomAspects}

In this section we discuss some aspects of closed surfaces
related to complex geometry. A general reference is \cite[Ch.'s
I--III]{W}. Moreover, a concise introduction can be found in
\cite[Sec. 1.4]{Nic:SW}.\\

\noindent\textbf{Complex Structures on Closed Surfaces.} Let
$\gS$ be a closed, oriented surface. Endow $\gS$ with a
Riemannian metric $g_\gS$, with volume form $\vol_\gS$ of unit
volume. The metric $g_\gS$ defines the structure of a complex
manifold on $\gS$ in the following way: Let $*$ be the Hodge star
operator on $\gS$. On $\gO^1(\gS)$ it satisfies $*^2 = -1$ and
thus gives an almost complex structure on $\gS$ such that
\begin{equation}\label{AlmostComplexDef}
\gO^{1,0}= \setdef{\ga\in \gO^1}{*\ga =
-i\ga}\quad\text{and}\quad \gO^{0,1}= \setdef{\ga\in \gO^1}{*\ga
= i\ga}.
\end{equation}
Denote by $P^{1,0}$ and $P^{0,1}$ the associated projections and
define the Dolbeault operators
\[
\pd:= P^{1,0}\circ d,\quad \bar\pd:=P^{0,1}\circ d.
\]
Using the Leibniz rule they extend to $\gO^\bullet\to\gO^\bullet$
and satisfy
\[
\pd^2=\bar\pd^2 = 0.
\]
This is because $\gO^{2,0}$ and $\gO^{0,2}$ are trivial on
2-dimensional almost complex manifolds. Therefore, the almost
complex structure is \emph{integrable}. This means that we can
find local coordinates $z=x+iy$ with holomorphic transition
functions such that $* dx = dy$. In these coordinates, the metric
$g_\gS$ is conformal to the standard metric on $\C$, i.e.,
\[
g_\gS = e^{2u(z,\bar z)} dz\otimes d\bar z,\quad u:U\subset \gS\to
\R.
\]
For a proof see \cite[Sec. 5.10]{Tay:I} or \cite[Thm
4.16]{McDSal}. The sheaf of holomorphic functions on $\gS$ is
given by
\[
\cO_\gS(U) = \ker \bar\pd|_U \subset C^\infty(U), \quad U\subset
\gS\text{ open}.
\]

\noindent\textbf{Holomorphic Line Bundles.} Now let $L\to\gS$ be a
Hermitian line bundle on $\gS$, and let $A$ be a unitary
connection on $L$ with associated covariant derivative
\[
d_A:\gO^0(\gS,L)\to \gO^1(\gS,L).
\]
Define the twisted Dolbeault operators by
\begin{equation}\label{DolbDef}
\begin{split}
&\pd_A:=P^{1,0}\circ d_A: \gO^0(\gS,L)\to \gO^{1,0}(\gS,L),\\
&\bar\pd_A:=P^{0,1}\circ d_A: \gO^0(\gS,L)\to \gO^{0,1}(\gS,L).
\end{split}
\end{equation}
As in the untwisted case the extension of $\bar\pd_A$ to
$\gO^{0,\bullet}(\gS,L)$ satisfies $\bar\pd_A^2 =0$. We can then
define a holomorphic structure on $L$ by requiring that its sheaf
of holomorphic sections is given by
\[
\cO(U,L_A):= \ker \bar\pd_A|_U.
\]
Since $\bar\pd_A(fs)=(\bar\pd f)s+f\bar\pd_As$, it is clear that
$\cO(L_A)$ is a sheaf of $\cO_\gS$-modules. Moreover, it follows
from elliptic theory that the space of global sections $\cO(\gS,
L_A):= \ker \bar\pd_A$ is finite dimensional.

\begin{remark*}
According to the above definition, every unitary connection on a
line bundle defines a holomorphic structure and one might ask
whether this is a suitable definition. Indeed, one can construct
a complex structure on the total space of $L$ such that the
projection $L\to\gS$ is holomorphic. A section $s\in
C^\infty(\gS,L)$ is then holomorphic as a map between $\gS$ and
$L$ if and only if $\bar\pd_A s=0$, see \cite[Thm. 2.1.53 \& Sec.
2.2.2]{DK}.
\end{remark*}

\noindent\textbf{The Riemann-Roch Theorem.} Although we will not
need it in the main body of the thesis, we digress briefly on the
famous Riemann-Roch Theorem in its version for line bundles, see
\cite[Thm 5.4.1]{JosRS} and \cite[Sec. 5.6]{JosRS}. Let
$K:=T\gS\to\gS$ be the tangent bundle of $\gS$ viewed as a complex
line bundle. The metric $g_\gS$ and the Levi-Civita connection
endow $K$ with a Hermitian metric, respectively, a holomorphic
structure. $K$ together with this structure is called the
\emph{canonical line bundle} of $\gS$.

\begin{theorem}[Riemann-Roch]\label{RiemannRoch}
Let $\gS$ be a closed Riemann surface of genus $g$. Let
$L_A\to\gS$ be a Hermitian line bundle endowed with the
holomorphic structure given by a unitary connection $A$. Then
\[
\dim \cO(\gS,L_A) - \dim \cO(\gS,K\otimes L_A^{-1}) = \deg L_A - g
+1.
\]
\end{theorem}

\begin{remark*}
The left hand side of the Riemann-Roch Theorem is the index of
the operator
\[
\bar\pd_A:\gO^0(\gS,L_A)\to \gO^{0,1}(\gS,L_A).
\]
To see this, note first that by definition $\ker
\bar\pd_A=\cO(\gS,L_A)$. To identify the cokernel, we note that
$\bar\pd^t_A = - *\pd_A *$. Recall that we are using the complex
\emph{linear} $*$ operator. Therefore, the $*$ operator maps the
kernel of $\bar\pd^t_A$ to
\[
\ker\big( \pd_A: \gO^{0,1}(\gS,L_A)\to \gO^2(\gS,L_A)\big).
\]
Observing that $K^{-1}=(T^*\gS)^{0,1}$, we can interpret $\pd_A$
as an operator
\[
\pd_A:\gO^0(\gS,K^{-1}\otimes L_A)\to \gO^{1,0}(\gS,K^{-1}\otimes
L_A).
\]
Since anti-holomorphic sections of a line bundle are in 1-1
correspondence to holomorphic sections of the dual bundle we get
that
\[
\dim(\ker \bar\pd^t_A) = \dim\cO(\gS,K\otimes L_A^{-1})
\]
which identifies the left hand side of the Riemann-Roch Formula
as an index. The right hand side is then the integral over the
corresponding index density as in the Atiyah-Singer Index Theorem
\ref{IndThm}, see \cite[Sec. 4.1]{BGV}.\qed
\end{remark*}

\noindent\textbf{Relation to the Signature Operator.} Since we are
usually dealing with the signature operator rather than the
Dolbeault operator we also want to mention how they can be
related. It follows from the definition \eqref{AlmostComplexDef}
of the almost complex structure on $\gS$ that
\[
\gO^{1,0} = \gO^+\cap \gO^1\quad\text{and}\quad \gO^{0,1} =
\gO^-\cap \gO^1,
\]
We thus have isometries
\[
\gF_+ : \gO^+ \to \gO^0\oplus\gO^{1,0},\quad \ga\mapsto \sqrt2
\ga_{[0]} + \ga_{[1]},
\]
and
\[
\gF_- : \gO^- \to \gO^{0,1}\oplus\gO^2,\quad \ga\mapsto \ga_{[1]}
+\sqrt2 \ga_{[2]}.
\]
Let $\bar\pd_A$ be the Dolbeault operator on
$\gO^\bullet(\gS,L_A)$. One checks that
\[
\gF_-\circ (d_A+d^t_A)\circ \gF_+^{-1} = \sqrt 2\bar\pd_A:
\gO^0\oplus\gO^{1,0} \to \gO^{0,1}\oplus\gO^2
\]
and
\[
\gF_+\circ (d+d^t)\circ \gF_-^{-1} = \sqrt 2\bar\pd^t_A:
\gO^{0,1}\oplus\gO^2 \to \gO^0\oplus\gO^{1,0}.
\]
This implies

\begin{lemma}\label{deRhamDolb}
The de Rham operator $d_A+d^t_A$ on $\gS$ with values in the line
bundle $L_A$ is unitary equivalent to
\[
\sqrt2(\bar\pd_A+\bar\pd^t_A):\gO^\bullet(\gS,L_A)\to
\gO^\bullet(\gS,L_A).
\]
Under this equivalence, the signature operator $D_A^+$
corresponds to
\[
\sqrt 2\bar\pd_A: \gO^0(\gS,L_A)\oplus\gO^{1,0}(\gS,L_A) \to
\gO^{0,1}(\gS,L_A)\oplus\gO^2(\gS,L_A)
\]
\end{lemma}

\begin{remark*}
With Lemma \ref{deRhamDolb} at hand one could derive the
Riemann-Roch Theorem from the Hirzebruch Signature Theorem and the
Gauss-Bonnet Theorem, or vice versa. Certainly, the relation
among these results becomes more complicated in higher dimensions.
\end{remark*}

\noindent\textbf{The Moduli Space of Holomorphic Line Bundles.}
We now want to give some remarks on the notion of equivalence of
holomorphic line bundles. First, recall from Section
\ref{FlatConnRep} that the group of gauge transformations
$\cG=C^\infty\big(\gS,\U(1)\big)$ acts on the space of Hermitian
connections $\cA(L)$ on $L$ via
\begin{equation*}
A\cdot u=A +u^{-1}du,\quad u\in \cG,\quad A\in \cA(L)
\end{equation*}
One easily checks that the associated twisted Dolbeault operators
satisfy
\[
\bar\pd_{A\cdot u} = \bar\pd_A + u^{-1}\bar\pd u = u^{-1}\bar\pd_A
u.
\]
Therefore, for every $U\subset \gS$ we get an isomorphism of
$\cO_\gS(U)$-modules
\begin{equation}\label{HolomEquiv}
\cO(U,L_A)\to \cO(U,L_{A\cdot u}),\quad s\mapsto u^{-1}s.
\end{equation}
Because of this, gauge equivalent connections on $L$ give rise to
equivalent holomorphic structures. The converse is certainly not
true. For this note that we can use any $f\in C^\infty(\gS,\C^*)$
to define a holomorphic structure on $L$ via
\begin{equation}\label{ActionCR}
(\bar\pd_A)_f:=f^{-1}\bar\pd_A f = \bar\pd_A + f^{-1}\bar\pd f.
\end{equation}
Via an isomorphism of the form \eqref{HolomEquiv}, this
holomorphic structure is equivalent to the one induced by $A$.
However, the underlying connection $A+f^{-1}df$ is in general not
unitary.

The relation between the moduli space of complex line bundles and
the moduli space of line bundles with connection can be made very
explicit. Since we want to avoid dealing with isomorphic
Hermitian line bundles which are not equal, we fix one Hermitian
line bundle $L\to\gS$ of degree 1 and let $L_k:=L^{\otimes k}$,
where a negative exponent means taking the tensor product of the
dual. We then let $\cA_k$ be the space of Hermitian connections on
$L_k$ and $\cA=\bigcup_k\cA_k$.

We now consider the complexified group $\cG^c:=
C^\infty(\gS,\C^*)$ of gauge transformations. It acts on the set
of holomorphic structures on $L_k$ via \eqref{ActionCR}. We want
to lift this to an action on $\cA_k$. The underlying idea why
this should be possible is a theorem of Chern that associates to
every holomorphic structure and Hermitian metric on a vector
bundle a unique compatible Hermitian connection, see \cite[Thm.
III.2.1]{W}. In the simple case at hand we have the following:

\begin{lemma}
Let $A\in \cA_k$ with associated Dolbeault operator $\bar\pd_A$,
and let $f\in \cG^c$. Define
\[
A\cdot f := A +f^{-1}\bar \pd f - \bar f^{-1}\pd \bar f.
\]
The $A\cdot f$ is a Hermitian connection satisfying
\[
\bar\pd_{A_f} = f^{-1}\bar\pd_A f.
\]
This defines an action of $\cG^c$ on $\cA_k$ which for $\cG\subset
\cG^c$ coincides with the standard action.
\end{lemma}

\begin{proof}
It is straightforward to check that $f^{-1}\bar \pd f - \bar
f^{-1}\pd \bar f$ is an imaginary valued 1-form on $\gS$. Thus,
the connection $A\cdot f$ is Hermitian. Since $P^{0,1}\big(\bar
f^{-1}\pd \bar f\big)=0$, the second assertion also follows. The
Leibniz rule applied to $\bar\pd$ and $\pd$ shows that
$A\cdot(fg)= (A\cdot f)\cdot g$ so that we indeed get an action of
$\cG^c$ on $\cA_k$. Moreover, a short computations shows that
\begin{equation}\label{ComplexGaugeAct}
f^{-1}\bar \pd f - \bar f^{-1}\pd \bar f = f^{-1}df -\pd\log |f|^2
\end{equation}
from which one readily finds that if $f$ takes values in $\U(1)$,
the action coincides with the usual one.
\end{proof}

\begin{dfn}\label{PicDef}
The moduli space of \emph{holomorphic structures} on $L_k$ is
defined as the quotient $\cA_k/\cG^c$. Moreover, we define the
\emph{Picard group} as
\[
\Pic(\gS) := \bigcup_{k\in\Z} \cA_k/\cG^c = \cA/\cG^c,
\]
where the group structure is induced by the tensor product of
line bundles with connection.
\end{dfn}

\begin{remark*}
We note without further details that we have defined $\cA_k/\cG^c$
and thus $\Pic(\gS)$ purely in differential geometric language.
The proof that the objects we obtain coincide with the ones
defined in holomorphic terms requires more work than sketched
here.
\end{remark*}

\noindent\textbf{Relation to the Moduli Space of Unitary
Connections.} We now have the ingredients to relate the moduli
space $\cA/\cG$ of line bundles with connections and the moduli
space $\Pic(\gS)$ of holomorphic line bundles. We start by
recalling a structure result for $\cA/\cG$, see \cite[Sec.
2.2.1]{DK}. Consider the map sending a connection to its
Chern-Weil representative,
\[
CW: \cA \to \gO^2(\gS,\R),\quad CW(A):=\frac i{2\pi}F_A.
\]
Note that no trace is involved since we are dealing with line
bundles. The image of $CW$ is easily seen to coincide with the
space $\gO^2_{cl}(\gS,\Z)$ of closed 2-forms with integral
periods, see Remark \ref{IntPeriod}. Moreover, $CW(A+ia)=CW(A)$
for every closed 1-form $a$. These considerations produce a short
exact sequence
\[
0\to\gO^1_{cl}(\gS,\R)\rightarrow \cA\xrightarrow{CW}
\gO^2_{cl}(\gS,\Z)\to 0.
\]
We also note that the action of $\cG$ on $\cA$ changes a
connection $A$ by a closed 1-form with integral periods. In
particular, the map $CW$ is invariant under the action of $\cG$.
Moreover, in the case at hand
\begin{equation}\label{JacobianTorus}
\gO^1_{cl}(\gS,\R)/\gO^1_{cl}(\gS,\Z) =
H^1(\gS,\R)/H^1(\gS,\Z)=H^1\big(\gS,\U(1)\big).
\end{equation}
Thus, taking quotients, we get the following exact sequence of
groups
\[
0\to H^1\big(\gS,\U(1)\big)\to \cA/\cG \xrightarrow{CW}
\gO_{cl}^2(\gS,\Z)\to 0.
\]

\begin{remark*}
This exact sequence generalizes to the case of an arbitrary closed
manifold $M$. However, the proof we sketched does not generalize
immediately. For this note that in the case of a surface $\gS$,
there are no line bundles which give torsion elements. Moreover,
the first homology group $H_1(\gS,\Z)$ is torsion-free, a fact we
used in the second equality of \eqref{JacobianTorus}. In the
general case, the Chern-Weil map does not capture possible
torsion. However, the universal coefficient theorem shows that
\[
\Tor\big(H^2(M,\Z)\big) = \Tor\big(H_1(M,\Z)\big).
\]
Therefore, the information about torsion is contained in
\[
H^1\big(M,\U(1)\big) = \Hom\big(H_1(M,\Z),\U(1)\big)
\]
which appears on the left hand side of the above sequence.
\end{remark*}

Returning now to holomorphic line bundles over surfaces we state
the following structure result.

\begin{prop}\label{PicModRel}
The natural projections
\[
\cA/\cG\to \Pic(\gS)\quad\text{ and } \quad \gO_{cl}^2(\gS,\Z)\to
H^2(\gS,\Z)
\]
fit into the following commutative diagram with exact rows
\[
\begin{CD}
0@>>>H^1\big(\gS,\U(1)\big) @>>> \cA/\cG @>{CW}>>
\gO_{cl}^2(\gS,\Z) @>>> 0 \\
@. @| @VVV @VVV @.\\
0@>>>H^1\big(\gS,\U(1)\big) @>>> \Pic(\gS) @>{c_1}>>
H^2(\gS,\Z)@>>> 0.
\end{CD}
\]
\end{prop}

\begin{proof}[Sketch of proof]
We can write every element $f\in\cG^c$ uniquely as
\[
f=\exp(\gf)u,\quad u\in\cG,\quad \gf\in C^\infty(\gS,\R).
\]
It follows from \eqref{ComplexGaugeAct} that for every $A\in\cA$
\[
A\cdot f = A + u^{-1}du + d\gf -2\pd \gf.
\]
Therefore, the moduli space $\Pic(\gS)=\cA/\cG^c$ can
alternatively be described as a quotient of $\cA/\cG$ by the
action
\[
[A]\cdot \gf := \big[A + d\gf - 2\pd \gf\big],\quad [A]\in
\cA/\cG,\quad \gf\in C^\infty(\gS,\R).
\]
With respect to this, the map $CW$ has the following equivariance
property:
\[
CW\big([A]\cdot \gf\big) = CW([A]) +\frac 1{\pi i} \bar\pd \pd
\gf = CW([A]) +\frac 1{2\pi} (\gD \gf)\vol_\gS.
\]
The latter equality is a consequence of the K\"{a}hler identities
but can also be checked directly in this simple case. Therefore,
we can equip the image of $CW$, i.e., the space
$\gO^2_{cl}(\gS,\Z)$, with a natural action of $C^\infty(\gS,\R)$
by defining
\[
\go\cdot\gf := \go +\frac 1{2\pi} (\gD \gf)\vol_\gS,\quad \gf\in
C^\infty(\gS,\R),\quad \go\in \gO^2_{cl}(\gS,\Z).
\]
It follows from the Hodge decomposition theorem that the quotient
of this action coincides with $H^2(\gS,\Z)$. The stabilizers
consist of the constant functions and thus agree with the
stabilizers of the action of $C^\infty(\gS,\R)$ on $\cA/\cG$.
Moreover, the action of $C^\infty(\gS,\R)$ on the fiber
$H^1\big(\gS,\U(1)\big)$ is trivial. Therefore, taking quotients
in the equivariant exact sequence
\[
0\to H^1\big(\gS,\U(1)\big)\to \cA/\cG \xrightarrow{CW}
\gO_{cl}^2(\gS,\Z)\to 0
\]
gives the exact sequence
\[
0\to H^1\big(\gS,\U(1)\big)\to \Pic(\gS) \xrightarrow{c_1}
H^2(\gS,\Z)\to 0
\]
and thus the requested diagram.
\end{proof}

\begin{remark*}
As in the case of Hermitian line bundles, there is a description
of $\Pic(\gS)$ in terms of \v{C}ech cohomology. Let $\cO^*_\gS$
be the sheaf of nowhere vanishing holomorphic functions on $\gS$.
Then
\[
\Pic(\gS) = H^1(\gS,\cO^*_\gS),
\]
see \cite[Lem. III.4.4]{W}. Moreover, the exponential sequence
$\Z\to\cO_\gS\to\cO^*_\gS$ gives rise to a long exact sequence in
cohomology
\[
...\to H^p(\gS,\Z) \to H^p(\gS,\cO_\gS) \to H^p(\gS,\cO_\gS^*)
\to H^{p+1}(\gS,\Z) \to...
\]
The sheaf $\cO_\gS$ is not fine so that the above sequence
contains much more information than its smooth version discussed
in Section \ref{FlatTrivSec}. In the case at hand, this sequence
is essentially the second line of the diagram in Proposition
\ref{PicModRel}. To see this, note that as $\gS$ is complex
1-dimensional, we have $H^2(\gS,\cO_\gS)=0$. Thus, the above
sequence reads
\[
H^1(\gS,\Z) \to H^1(\gS,\cO_\gS) \to \Pic(\gS) \to H^2(\gS,\Z)\to
0.
\]
Moreover, the space $H^1(\gS,\cO_\gS)$ is isomorphic to the
Dolbeault cohomology group $H^{0,1}(\gS)$. Via Hodge theory, the
latter can be identified with $H^1(\gS,\R)$. Using
\eqref{JacobianTorus} we arrive at the sequence of Proposition
\ref{PicModRel}.
\end{remark*}
\cleardoublepage
\chapter{Some Computations}\label{Comp}

Here, we include a computational discussion which will be used in
the main body of this thesis. We introduce basic Eta and Zeta
functions, and derive some of their values. In the second part of
this appendix, we briefly discuss the Dedekind sums and their
generalization which we use in Section \ref{TorusBundlesExp}. In
particular, we establish the relation among them which we need to
prove Theorem \ref{RhoHyp}.

\section{Values of Zeta and Eta Functions}\label{CompGamma}

\subsection{The Gamma and the Hurwitz Zeta Function}

We will need some standard facts concerning the Gamma function
and the generalization by Hurwitz of the Riemann Zeta function,
see for example \cite[Ch. 9]{Cohen} as a general reference. Recall
that the integral representation of the Gamma function is
\[
\gG(s) = \int_0^\infty e^{-t}t^{s-1}dt,\quad \Re(s)>0.
\]
It satisfies the functional equation
\begin{equation}\label{GammaFunctEqn}
\gG(s+1)=s\gG(s),
\end{equation}
which can be used to extend $\gG(s)$ meromorphically to the whole
plane. Then $\gG(s)$ has no zeros and only simple poles at
$s=0,-1,-2,\ldots$ The residues are given by
\begin{equation}\label{GammaPoles}
\Res\gG(s)\big|_{s=-n}  = \frac{(-1)^n}{n!},\quad n\in \N.
\end{equation}
The \emph{Hurwitz Zeta function} is defined for $q\in \R^+$ as
\begin{equation}\label{Hurwitz} \gz_q(s):=
\sum_{ n =0}^\infty \frac 1{(n+q)^s}, \quad \Re(s)>1.
\end{equation}
In particular, $\gz_1(s)$ is the Riemann Zeta function. Using the
Mellin transform one can derive basic properties of $\gz_q(s)$.
First, we note that
\[
\frac 1{(n+q)^s} = \frac 1{\gG(s)}\int_0^\infty t^{s-1}
e^{-t(n+q)}dt.
\]
The formula for the geometric series shows that for $t>0$,
\begin{equation}\label{HurwitzKernel}
\sum_{n=0}^\infty e^{-t(n+q)} = \frac {e^{-tq}}{1- e^{-t}},
\end{equation}
which clearly decays exponentially with $t$ as $t\to \infty$. \\

\noindent\textbf{Bernoulli Polynomials.} The behaviour as $t\to
0$ of \eqref{HurwitzKernel} is determined by the \emph{Bernoulli
polynomials} $B_n(x)$. Recall, e.g. from \cite[Sec. 9.1]{Cohen},
that they can be defined via the generating function
\begin{equation}\label{BernoulliPol}
\frac{te^{xt}}{e^t-1} = \sum_{n= 0}^\infty B_n(x) \frac{t^n}{n!},
\quad |t|<2\pi,\quad x\in \R.
\end{equation}
Comparing this with \eqref{Bernoulli}, we note that the Bernoulli
numbers with respect to the normalization we are using are given
by $B_n=B_n(0)$. Moreover, expanding $e^{xt}$ in
\eqref{BernoulliPol} and comparing coefficients of $t^n$ yields
that
\[
B_n(x) = \sum_{k=0}^n \begin{pmatrix} n \\ k \end{pmatrix} B_k
x^{n-k}.
\]
In particular,
\begin{equation}\label{BernoulliPolExpl}
B_0(x) = 1,\quad B_1(x) = x- \lfrac 12,\quad B_2(x) = x^2 -x
+\lfrac 16.
\end{equation}
From the definition \eqref{BernoulliPol}, one easily finds that
the Bernoulli polynomials have the symmetry property
\begin{equation}\label{BernoulliSymm}
B_n(1-x) = (-1)^n B_n(x).
\end{equation}

\noindent\textbf{Values of the Hurwitz Zeta Function.} It now
follows from \eqref{BernoulliPol} and \eqref{BernoulliSymm}
applied to \eqref{HurwitzKernel} that for $t\in (0,2\pi)$
\[
\sum_{n=0}^\infty e^{-t(n+q)} = \frac {e^{(1-q)t}}{e^t-1} =
\sum_{n= 0}^\infty (-1)^n B_n(q) \frac{t^{n-1}}{n!}.
\]
In particular, $\sum_{n=0}^\infty e^{-t(n+q)} = O(t^{-1})$ as
$t\to 0$. This implies that we can apply the Mellin transform to
\eqref{Hurwitz} and interchange summation and integration. Then
splitting the integral into $\int_0^1+ \int_1^\infty$ one easily
obtains that for $\Re(s)>1$
\begin{equation*}
\gG(s)\gz_q(s) = h(s) + \sum_{n=0}^\infty \frac{(-1)^n
B_n(q)}{(s+n-1)n!},
\end{equation*}
where $h(s)$ can be extended to a holomorphic function of $s\in
\C$. Therefore, $\gG(s)\gz_q(s)$ extends to a meromorphic
function on the whole plane. It has simple poles at the points $s=
1,0,-1,-2,\ldots$ with residues
\begin{equation}\label{ZetaBernoulliRes}
\Res\big(\gG(s)\gz_q(s)\big)\big|_{s=-n+1} = (-1)^n
\frac{B_n(q)}{n!},\quad n\in \N.
\end{equation}
Now, since $\gG(s)$ has no zeros, we can deduce that $\gz_q(s)$
extends to a meromorphic function on the whole plane which can
have only simple poles. Using \eqref{GammaPoles} and
\eqref{ZetaBernoulliRes} one finds that $\gz_q(s)$ has a simple
pole at $s=1$ with residue 1. The other poles are cancelled out
by the zeroes of $\gG(s)^{-1}$ and
\begin{equation}\label{ValueZeta}
\gz_q(-n) = -  \frac{B_{n+1}(q)}{n+1},\quad n\in \N.
\end{equation}

\subsection{An Eta Function and a ``Periodic'' Zeta Function}

\begin{dfn}\label{PeriodicBernoulli}
For $x\in \R$ let $[x]$ denote the largest integer less or equal
than $x$. We define the \emph{$n$-th periodic Bernoulli function}
as
\[
P_n(x):= \begin{cases}\hphantom{B_n\big(x} 0, &\text{if $x\in \Z$},\\
B_n\big(x- [x]\big), &\text{if $x\notin \Z$,}
\end{cases}
\quad \text{for $n$ odd,}
\]
and
\[
\quad P_n(x) := B_n\big(x- [x]\big),\quad \text{for $n$ even}.
\]
\end{dfn}

\begin{remark*}
The definition of $P_n(x)$ for odd $n$ is a bit artificial. Note
that \eqref{BernoulliSymm} implies that $B_n(1) = (-1)^n B_n(0)$,
so that for $n$ even, we have $B_n(1) = B_n(0)$. We note without
proof, that for odd $n$ with $n\neq 1$, one has $B_n(1) =
B_n(0)=0$ so that a distinction is unnecessary. However, when
working with $P_1(x)$, the above convention is convenient. Using
\eqref{BernoulliPolExpl} we note that explicitly,
\begin{equation}\label{ExplPeriodBern}
P_1(x)= \begin{cases}\hphantom{B_n\big(x} 0, &\text{if $x\in
\Z$},\\ x- [x]-\lfrac 12, &\text{if $x\notin \Z$,}
\end{cases}\quad P_2(x) = \big(x-[x]\big)^2 - \big(x-[x]\big)
+\lfrac 16.
\end{equation}
\end{remark*}

Most of our computations of Eta invariants in the main body of
this thesis will be based on the following result.

\begin{prop}\label{ZEtaCalc}
Let $q\in\R$, and define for $\Re(s)>1$
\[
\eta_q(s):= \sum_{\begin{smallmatrix} n \in\Z \\
n\neq q
\end{smallmatrix}}\frac{\sgn(n-q)}{|n -q|^s},
\quad\text{and}\quad
\widetilde \gz_q(s):= \sum_{\begin{smallmatrix} n \in\Z \\
n\neq q
\end{smallmatrix}}\frac 1{|n -q|^s}.
\]
\begin{enumerate}
\item The function $\eta_q(s)$ extends to a holomorphic function
for all $s\in \C$,  and
\[
\eta_q(0) = 2 P_1(q).
\]
\item The function $\widetilde \gz_q(s)$ extends to a meromorphic
function with only one simple pole at $s=1$. Moreover,
\[
\widetilde \gz_q(0)=\begin{cases}\,\, 0, &\text{\rm if $q\notin \Z$}, \\
-1, &\text{\rm if $q\in \Z$},
\end{cases}\qquad \gz_q(-1) = - P_2(q).
\]
\end{enumerate}
\end{prop}

\begin{proof}
Write $q_0 := q -[q]\in [0,1)$. Since the sums defining
$\eta_q(s)$ and $\widetilde \gz_q(s)$ converge absolutely for
$\Re(s)>1$, we can change the order of summation. Then, if $q \in
\Z$ so that $q_0 = 0$, we find that
\[
\eta_q(s) = \gz_1 (s) - \gz_1(s) =0,\quad \text{and}\quad
\widetilde \gz_q(s) = 2\gz_1(s).
\]
If $q\notin \Z$ we get
\[
\eta_q(s) = \sum_{n=1}^\infty \Big(\lfrac 1{n -q_0}\Big)^{s} -
\sum_{n=0}^\infty \Big(\lfrac 1{n +q_0}\Big)^{s} = \gz_{1-q_0}(s)
- \gz_{q_0}(s),
\]
and
\[
\widetilde \gz_q(s) = \sum_{n = 1}^\infty \Big(\lfrac 1{n
-q_0}\Big)^{s} + \sum_{n = 0}^\infty \Big(\lfrac 1{n
+q_0}\Big)^{s}= \gz_{1-q_0}(s) + \gz_{q_0}(s).
\]
Since $\gz_1 (s)$, $\gz_{1-q_0}(s)$ and $\gz_{q_0}(s)$ extend to
meromorphic functions on the whole plane, with only one simple
pole at $s=1$ of residue 1, we can extend $\eta_q(s)$ and
$\widetilde \gz_q(s)$ meromorphically. One sees that $\eta_q(s)$
has no pole, whereas $\widetilde \gz_q(s)$ has a simple pole at
$s=1$. Moreover, $\eta_q(s)$ vanishes if $q\in \Z$. Otherwise, we
deduce from \eqref{ValueZeta}, \eqref{BernoulliPolExpl} and
\eqref{BernoulliSymm} that
\[
\eta_q(0) = \gz_{1-q_0}(s) - \gz_{q_0}(s) = - B_1(1 -q_0) +
B_1(q_0) = 2 B_1(q_0).
\]
From the definition of $P_1(q)$ and $q_0$, part (i) follows.
Concerning part (ii), we first assume that $q_0=0$. Then
\[
\widetilde \gz_q(0) = - 2 B_1(1) = - 1,\quad\text{and}\quad
\widetilde \gz_q(-1) = - 2 \lfrac {B_2(1)}2 = - B_2(0).
\]
For $q_0\neq 0$, one finds that
\[
\widetilde \gz_q(0) = - B_1(1-q_0) - B_1(q_0) = B_1(q_0)-B_1(q_0)
= 0,
\]
and
\[
\widetilde \gz_q(-1) = -\lfrac 12\big ( B_2(1-q_0) +
B_2(q_0)\big) = - B_2(q_0).\qedhere
\]
\end{proof}

\begin{remark*}
Clearly, one could go on without difficulty, and determine more
values of $\eta_q(s)$ and $\widetilde \gz_q(s)$ in terms of the
periodic Bernoulli functions. However, Proposition \ref{ZEtaCalc}
covers all the cases we are interested in.
\end{remark*}

\section{Generalized Dedekind Sums}\label{CompDedekind}

When studying the Eta invariant for 2-dimensional torus bundles
over the circle, one naturally encounters versions of the
\emph{Dedekind sums}. In this section, we include some relevant
definitions and computations. We start to collect some facts
about finite Fourier series, see \cite[Ch. 7]{Beck}.

\subsection{Some Finite Fourier Analysis}\label{DedFour}

Let $c\in \Z$ with $c\neq 0$. In this section we will always use
the $c$-th root of unity $\xi := \exp\big(\lfrac{2\pi i}c\big)$.

\begin{dfn}
Assume that $f:\Z\to \C$ is \emph{$c$-periodic}, i.e.,
$f(k+c)=f(k)$ for all $k\in \Z$. The \emph{Fourier transform}
$\widehat f:\Z\to \C$ is defined as
\[
\widehat f(p) = \sum_{k=0}^{|c|-1} f(k)\xi^{-kp},\quad p\in \Z.
\]
\end{dfn}

\begin{remark}\label{FourierSignRel}
Since we allow $c$ to be negative, one has to be a bit careful
concerning signs. Let $\eps=\sgn(c)$ and denote by $\widehat
f^\eps$ the Fourier transform of $f$ with respect to the $c$-th
root of unity $\xi^{\eps}=\exp\big(\lfrac{2\pi i}{|c|}\big)$.
Then for all $p\in \Z$
\[
\widehat f^\eps(p) = \sum_{k=0}^{|c|-1} f(k)(\xi^\eps)^{-k p} =
\widehat f(\eps p).
\]
\end{remark}

Since $\xi^{-kp}$ is $c$-periodic in $p$, the Fourier transform is
again $c$-periodic. Moreover, we can shift the sum by any $m\in
\Z$, i.e.,
\[
\widehat f(p) = \sum_{k=m}^{|c|+m-1} f(k)\xi^{-kp}.
\]
This implies that if $g(k):= f(k+m)$, then
\begin{equation}\label{FourierTrans}
\widehat g(p) = \sum_{k=0}^{|c|-1} f(k+m)\xi^{-kp} =
\sum_{k=0}^{|c|-1} f(k)\xi^{-(k-m)p} = \xi^{mp}\widehat f(p).
\end{equation}
Furthermore, if $a\in \Z$ with $\gcd(a,c)=1$, then
$\setdef{ak}{k=0,\ldots, |c|-1}$ is a representation system of
$\Z$ modulo $c$, so that
\[
\widehat f(p) = \sum_{k=0}^{|c|-1} f(ak) \xi^{-akp},\quad
\gcd(a,c)=1.
\]
Hence, if $d\in Z$ is an inverse of $a$ modulo $c$, i.e., $ad
\equiv 1\,(c)$, then $g(k):= f(ak)$ satisfies
\begin{equation}\label{FourierMult}
\widehat g(p) = \sum_{k=0}^{|c|-1} f(ak) \xi^{-kp} =
\sum_{k=0}^{|c|-1} f(k) \xi^{-kdp} = \widehat f(dp).
\end{equation}
The finite geometric series yields that
\[
\sum_{k=0}^{|c|-1} \xi^{kp} =
\begin{cases} |c|, &\text{if $p\equiv 0\, (c)$},\\
\hphantom{|}0&\text{otherwise}. \end{cases}
\]
From this one easily deduces the \emph{Fourier inversion formula}
\begin{equation}\label{FourierInversion}
f(k) = \lfrac 1{|c|} \sum_{p=0}^{|c|-1} \widehat f(p)\xi^{pk},
\end{equation}
see \cite[Thm. 7.2]{Beck}. Moreover, if $g$ is another
$c$-periodic function, one has the convolution formul{\ae}
\begin{equation}\label{Convolution}
(f* g)(k) := \sum_{l=0}^{|c|-1} f(l)g(k-l) = \lfrac 1{|c|}
\sum_{p=0}^{|c|-1} \widehat f(p)\widehat g(p) \xi^{pk},
\end{equation}
and
\begin{equation}\label{Convolution:Alt}
(\widehat f* \widehat g)(p) =  |c| \sum_{k=0}^{|c|-1}   f(k) g(k)
\xi^{-pk},
\end{equation}
see \cite[Thm 7.10]{Beck}.

The facts we have collected so far are sufficient for the
application to generalized Dedekind sums in Section
\ref{CompDedekindSec} below. Yet, we need to compute the Fourier
transform for one particular class of functions, which form the
building blocks of generalized Dedekind sums. First, we introduce
some notation. Note that fixing a pair $(a,c)$ with
$\gcd(a,c)=1$, and an inverse $d$ of $a$ modulo $c$ is the same
as fixing a matrix
\[
M = \begin{pmatrix} a &b\\ c& d
\end{pmatrix}\in \SL_2(\Z),\quad c\neq 0.
\]
Here, $b$ is uniquely determined by requiring that $ad - bc =1$.
Moreover, for $x,y\in \R$ we define
\begin{equation}\label{x'Def}
\begin{pmatrix} x'\\ y' \end{pmatrix} := M^t \begin{pmatrix} x\\ y
\end{pmatrix}= \begin{pmatrix} ax +cy \\ bx+ dy \end{pmatrix}.
\end{equation}

\begin{prop}\label{DedekindPrep}
Let $P_1$ be the first periodic Bernoulli function, see
Definition \ref{PeriodicBernoulli}. Fix
$M=\left(\begin{smallmatrix} a &b\\ c& d
\end{smallmatrix}\right)\in \SL_2(\Z)$ with $c\neq 0$, let
$x,y \in \R$ and $x'$ as in \eqref{x'Def}. Define a $c$-periodic
function by
\[
f_{x,y,a,c}(k) := P_1\left(a\frac {k+x}{c} + y \right),\quad k\in
\Z.
\]
Then,
\[
\widehat f_{x,y,a,c}(p) = \sgn(c) \begin{cases} \hphantom{\frac
12\big( i\cot\big(\frac{\pi}c}
P_1(x'), &\text{if $p\equiv 0\,(c)$,}\\
\frac 12\Big( i\cot\big(\frac{\pi dp}{c}\big) -
\gd(x')\Big)\xi^{d[x']p}, &\text{otherwise.}
\end{cases}
\]
Here, $\gd$ is the characteristic function of $\R\setminus \Z$,
i.e., $\gd(x') = 0$ if $x'\in \Z$ and $\gd(x') = 1$ if $x'\notin
\Z$. Moreover, as in Definition \ref{PeriodicBernoulli}, the
expression $[x']$ refers to the largest integer less or equal
than $x'$.
\end{prop}

Before we give the proof of Proposition \ref{DedekindPrep} let us
collect some special cases.

\begin{cor}\quad
\begin{enumerate}
\item If $f(k):= f_{0,0,1,c}(k)= P_1\big(\frac k{c}\big)$, then
\begin{equation}\label{P1Cot}
\widehat f(p) = \begin{cases} \hphantom{\frac i2 \cot} 0,
&\text{if $p\equiv 0\,(c)$,}\\
\frac i2 \cot\big(\frac{\pi p}{|c|}\big), &\text{otherwise.}
\end{cases}
\end{equation}
In particular,
\begin{equation}\label{SumFullRepSys}
\sum_{k=1}^{|c|-1} P_1\big(\lfrac k{c}\big) = 0.
\end{equation}
\item Let $x\in \R$, and $f_x(k):= f_{x,0,1,c}(k)= P_1\big(\frac
{k+x}{c}\big)$. Then
\begin{equation}\label{P1TransCot}
\widehat f(p) = \sgn(c) \begin{cases} \hphantom{\frac 12\big(
i\cot\big(\frac{\pi}c}
P_1(x), &\text{if $p\equiv 0\,(c)$,}\\
\frac 12\Big( i\cot\big(\frac{\pi p}{c}\big) -
\gd(x)\Big)\xi^{[x]p}, &\text{otherwise.}
\end{cases}
\end{equation}
\end{enumerate}
\end{cor}

\begin{proof}[Proof of Proposition \ref{DedekindPrep}]
The proof consists of proving the special cases \eqref{P1Cot} and
\eqref{P1TransCot} first. The general case then follows using
\eqref{FourierTrans} and \eqref{FourierMult}. Formula
\eqref{P1Cot} is standard, see \cite[Lem 7.3]{Beck}.
Nevertheless, we sketch a proof, since most computations we
encounter are deduced from this formula. We assume first that
$c>0$. For $k\in \{0,\ldots,c-1\}$ we have
\[
f(k) = P_1\big(\lfrac kc \big)= \begin{cases}
\hphantom{ \frac kc\,} 0, &\text{if $k=0$,}\\
 \frac kc -\frac 12, &\text{otherwise}.
\end{cases}
\]
Therefore, for $p\equiv 0\,(c)$,
\[
\widehat f(p) = \sum_{k=1}^{c-1} \big( \lfrac kc -\lfrac 12\big)
= \lfrac 1c \lfrac {c(c-1)}2 - \lfrac {c-1}2 =0,
\]
which is also the claim in \eqref{SumFullRepSys}. Now, for $p$
not divisible by $c$, we have
\[
\widehat f(p) = \sum_{k=1}^{c-1} \big( \lfrac kc -\lfrac 12\big)
\xi^{-kp} =  - \lfrac 12 \sum_{k=1}^{c-1}\xi^{-kp} +
\frac{\pd}{\pd x}\Big|_{x=-2\pi i p}
\sum_{k=1}^{c-1}\exp\big(\lfrac {kx}c\big).
\]
For $x\notin 2\pi i c \Z$, the formula for the finite geometric
series states that
\[
\sum_{k=1}^{c-1}\exp\big(\lfrac {kx}c\big) =
\frac{1-e^x}{1-e^{x/c}} -1.
\]
Using this one verifies that
\[
\frac{\pd}{\pd x}\Big|_{x=-2\pi i p}
\sum_{k=1}^{c-1}\exp\big(\lfrac {kx}c\big) = -\frac 1{1-
\xi^{-p}}= - \frac {\xi^{p/2}}{\xi^{p/2}- \xi^{-p/2}} = \lfrac 12
\big(i\cot\big(\lfrac {\pi p }c\big) - 1\big).
\]
Moreover, for $p$ not divisible by $c$, one has
$\sum_{k=1}^{c-1}\xi^{-kp} = -1$. Therefore, in this case
\[
\widehat f(p) = \lfrac 12 + \lfrac 12 \big(i\cot\big(\lfrac {\pi p
}c\big) - 1\big) = \lfrac i2\cot\big(\lfrac {\pi p }c\big),
\]
which proves \eqref{P1Cot} for $c>0$. For general $c\neq 0$, let
$\eps=\sgn(c)$. Then $P_1\big(\lfrac kc\big) = \eps P_1\big(\lfrac
k{|c|}\big)$, so that we deduce from the case $c>0$ and Remark
\ref{FourierSignRel} that
\[
\widehat f(p) = \begin{cases} \hphantom{\frac i2 \cot} 0,
&\text{if $p\equiv 0\,(c)$,}\\
\frac i2 \eps   \cot\big(\frac{\pi \eps p}{|c|}\big),
&\text{otherwise.}
\end{cases}
\]
Since the cotangent is an odd function, we obtain \eqref{P1Cot}
for $c<0$ as well.

Concerning \eqref{P1TransCot}, we assume first that $x\in [0,1)$,
and again that $c>0$. Then
\begin{equation}\label{DedekindPrep:2}
f_x(k) = \begin{cases}\hphantom{P_1} P_1\big(\frac xc\big),
&\text{if $k\equiv
0\, (c)$},\\
P_1\big(\frac kc\big) + \frac xc, &\text{otherwise.}
\end{cases}
\end{equation}
Therefore, for $p\equiv 0\, (c)$,
\[
\widehat f_x(p) = P_1\big(\lfrac xc\big) + \sum_{k=1}^{c-1}\lfrac
xc + \sum_{k=1}^{c-1} P_1\big(\lfrac kc\big) = P_1\big(\lfrac
xc\big) + \lfrac {c-1}c x,
\]
where we have used \eqref{SumFullRepSys}. Now, one easily
verifies that for $x\in [0,1)$
\[
P_1\big(\lfrac xc\big) + \lfrac {c-1}c x= P_1(x),
\]
which implies \eqref{P1TransCot} for the case $p\equiv 0\,(c)$. If
$p$ is not divisible by $c$ we use \eqref{DedekindPrep:2} and
\eqref{P1Cot} to deduce that
\[
\widehat f_x(p) = P_1\big(\lfrac xc\big)+ \lfrac xc
\sum_{k=1}^{c-1} \xi^{-kp} + \sum_{k=1}^{c-1} P_1\big(\lfrac
kc\big) \xi^{-kp}  = P_1\big(\lfrac xc\big) -\lfrac xc + \lfrac
i2\cot\big(\lfrac{\pi p}c\big).
\]
Since $P_1\big(\lfrac xc\big) -\lfrac xc =-\lfrac 12 \gd(x)$,
formula \eqref{P1TransCot} follows for $x\in [0,1)$. For
arbitrary $x\in \R$ we write $x= [x]+ x_0$, so that
\[
f_x(k) = f_{x_0} \big(k + [x]\big),\quad x_0\in [0,1).
\]
Since we have already proved \eqref{P1TransCot} for $f_{x_0}$, we
can use \eqref{FourierTrans} to get the required formula for
$f_x$. For this note that in the case $p\equiv 0\,(c)$ we have
$\xi^{[x]p}=1$. The case $c<0$ follows as before, using Remark
\ref{FourierSignRel} and the fact that $f_x$ is odd.

Now for the general formula of Proposition \ref{DedekindPrep}
observe that
\[
f_{x,y,a,c}(k) = P_1\big(a\lfrac {k+x}c + y \big) = P_1\big(\lfrac
{ak+ax+ cy}c\big) =P_1\big(\lfrac {ak+x'}c\big) = f_{x'}(ak),
\]
where $x'=ax+ cy$. Then \eqref{FourierMult} implies that
\[
\widehat f_{x,y,a,c}(p) = \widehat f_{x'}(dp),
\]
so that Proposition \ref{DedekindPrep} follows immediately from
\eqref{P1TransCot}.
\end{proof}

\subsection{Relation Among Some Dedekind Sums}\label{CompDedekindSec}

Let $M= \left(\begin{smallmatrix} a &b\\ c& d
\end{smallmatrix}\right)\in \SL_2(\Z)$ with $c\neq 0$.
Recall that the classical \emph{Dedekind sums} are defined by
\begin{equation}\label{DedDef}
s(a,c) := \sum_{k=1}^{|c|-1} P_1\big(\lfrac
{ak}{c}\big)P_1\big(\lfrac k{c}\big),
\end{equation}
see e.g., \cite[p. 128]{Beck}. Since the first periodic Bernoulli
function is odd, we can replace $c$ with $|c|$ in both
denominators of \eqref{DedDef}. Moreover, using that
$\setdef{ap}{p=0,\ldots,|c|-1}$ is a representation system for
$\Z$ modulo $c$, and that $ad\equiv 1\,(c)$, one can replace $a$
with $d$, so that
\begin{equation}\label{DedSymm}
s(a,c) = s(a,|c|)= s(d,c) = s(d,|c|).
\end{equation}
There are many more relations among different Dedekind sums but
their discussion would lead to far afield. We only mention that
\eqref{P1Cot} and \eqref{Convolution:Alt} easily imply the
classical \emph{cotangent formula}
\begin{equation*}
s(a,c) = \lfrac 1{|c|}\sum_{p=1}^{|c|-1} \lfrac
i2\cot\big(\lfrac{dp}{|c|}\big)\lfrac i2\cot\big(-
\lfrac{p}{|c|}\big) = \lfrac 1{4|c|}\sum_{p=1}^{|c|-1}
\cot\big(\lfrac{dp}{c}\big)\cot\big(\lfrac{p}{c}\big),
\end{equation*}
The Dedekind sums $s(a,c)$ were generalized in several ways. The
generalization we are interested in was considered in \cite{Die59,
Mey57, Rad64}.

\begin{dfn}\label{GenDedDef}
For $x,y\in \R$ define
\[
s_{x,y}(a,c) := \sum_{k=0}^{|c|-1} P_1\big(a\lfrac
{k+x}{c}+y\big)P_1\big(\lfrac {k+x}{c}\big).
\]
\end{dfn}

Again, there are several relations among generalized Dedekind sums
for different values of $(x,y)$ and $(a,c)$. For example,
$s_{x,y}(a,c)$ depends on $(x,y)$ only modulo $\Z^2$. For $y$
this is immediate and for $x$ this is because for $m\in \Z$, one
has
\begin{equation}\label{GenDedDep}
s_{x+m,y}(a,c) = \sum_{k=-m}^{|c|-m-1} P_1\big(a\lfrac
{k+x}{c}+y\big)P_1\big(\lfrac {k+x}{c}\big) = \sum_{k=0}^{|c|-1}
P_1\big(a\lfrac {k+x}{c}+y\big)P_1\big(\lfrac {k+x}{c}\big).
\end{equation}
We will not include a separate treatment of other relations but
give an ad hoc explanation whenever we are using them. However,
what we want to single out is the following straightforward
consequence of Proposition \ref{DedekindPrep} and
\eqref{Convolution:Alt}:

\begin{lemma}\label{GenDedCot}
Let $x'$ be as in \eqref{x'Def}, and let $\gd$ be the
characteristic function of $\R\setminus\Z$. Then
\[
\begin{split}
s_{x,y}(a,c) = \lfrac 1{|c|} P_1(x)P_1(x') - \lfrac
1{4|c|}\sum_{p=1}^{|c|-1} \Big(i \cot&\big(\lfrac{\pi p}{c}\big)
- \gd(x)\Big) \\ &\times\Big(i \cot\big(\lfrac{\pi dp}{c}\big) +
\gd(x')\Big) \xi^{([x]-d[x'])p}.
\end{split}
\]
\end{lemma}

We will now make a simplifying assumption, which renders the
following formul{\ae} a bit more transparent.

\begin{assumption}
From now on, $(x,y)\in \R^2$ will always be chosen in such a way
that $x\in [0,1)$ and
\begin{equation}\label{IntegerAssump}
\begin{pmatrix} m\\ n
\end{pmatrix}:= (\Id - M^t) \begin{pmatrix} x\\ y
\end{pmatrix} = \begin{pmatrix} x-x'\\ y-y'
\end{pmatrix}\in \Z^2.
\end{equation}
\end{assumption}

Note that under this assumption we have $x- x'= m\in \Z$ and
$[x]=0$, which yields that $[x']=-m$ and $x'-[x'] = x$. Thus, the
cotangent formula of Lemma \ref{GenDedCot} simplifies to
\begin{equation}\label{GenDedCotAssump}
s_{x,y}(a,c) = \lfrac 1{|c|} P_1(x)^2 - \lfrac
1{4|c|}\sum_{p=1}^{|c|-1} \Big(i \cot\big(\lfrac{\pi p}{c}\big) -
\gd(x)\Big)\Big(i \cot\big(\lfrac{\pi dp}{c}\big) + \gd(x)\Big)
\xi^{dmp}.
\end{equation}
We then have the following relation between the generalized
Dedekind sums and the classical Dedekind sums, see also
\cite[Sec. 8]{Mey57} and \cite[Sec. 4]{Mey61}.

\begin{prop}\label{GenDedRel}
Under the above assumption, let $r\in \{0,\ldots |c|-1\}$ with
$m\equiv r\, (c)$. Then
\[
\begin{split}
s_{x,y}(a,c) - s(a,c) = \lfrac 1{|c|} \big(P_2(x)&-\lfrac 16\big)
+ \sum_{k=1}^{|c|-r} P_1\big(\lfrac{dk}{|c|}\big)  +   \lfrac 12
P_1\big(\lfrac {dm}{|c|}\big)
\\  & +\lfrac 12 \gd(x)\Big(
P_1\big(\lfrac{m}{|c|}\big) - P_1\big(\lfrac{dm}{|c|}\big)\Big) +
\lfrac 1{4}\gd(x)\big(1-\gd(\lfrac{m}c)\big).
\end{split}
\]
\end{prop}

\begin{proof}
The cotangent formula \eqref{GenDedCotAssump} shows that
\begin{equation}\label{GenDedRel:1}
\begin{split}
s_{x,y}(a,c) &= \lfrac 1{|c|} P_1(x)^2 + \lfrac
1{4|c|}\sum_{p=1}^{|c|-1} \cot\big(\lfrac{\pi p}{c}\big)
\cot\big(\lfrac{\pi dp}{c}\big) \xi^{dmp} + \lfrac
1{4|c|}\gd(x) \sum_{p=1}^{|c|-1}\xi^{dmp}\\
&\qquad+ \lfrac 1{2|c|}\gd(x) \sum_{p=1}^{|c|-1} \lfrac i2
\cot\big(\lfrac{\pi dp}{c}\big)\xi^{dmp} - \lfrac 1{2|c|}\gd(x)
\sum_{p=1}^{|c|-1} \lfrac i2 \cot\big(\lfrac{\pi
p}{c}\big)\xi^{dmp} .
\end{split}
\end{equation}
Let $\Tilde x = 0$ and $\Tilde y= -\frac mc$, so that in the
notation of \eqref{x'Def} we have $\Tilde x' = -m$, and $[\Tilde
x] - d[\Tilde x'] = dm$. Applying Lemma \ref{GenDedCot} to
$s_{\Tilde x,\Tilde y}(a,c)$ yields
\[
s_{\Tilde x,\Tilde y}(a,c) = \lfrac 1{4|c|}\sum_{p=1}^{|c|-1}
\cot\big(\lfrac{\pi p}{c}\big) \cot\big(\lfrac{\pi dp}{c}\big)
\xi^{dmp},
\]
which is exactly the first sum in \eqref{GenDedRel:1}. On the
other hand, by definition,
\[
s_{\Tilde x,\Tilde y}(a,c) = \sum_{k=1}^{|c|-1}
P_1\big(\lfrac{ak-m}{c}\big) P_1\big(\lfrac k{c}\big) =
\sum_{k=1}^{|c|-1} P_1\big(\lfrac{ak+m}{c}\big) P_1\big(\lfrac
k{c}\big) = \sum_{k=1}^{c-1} P_1\big(\lfrac{k+m}{c}\big)
P_1\big(\lfrac {dk}{c}\big),
\]
where we have first used that $P_1$ is odd and summed over $-k$,
and then summed over $ak$ instead of $k$. Comparing the two
expressions for $s_{\Tilde x,\Tilde y}(a,c)$ we find that
\begin{equation}\label{GenDedRel:2}
\lfrac 1{4|c|}\sum_{p=1}^{|c|-1} \cot\big(\lfrac{\pi p}{c}\big)
\cot\big(\lfrac{\pi dp}{c}\big) \xi^{dmp} = \sum_{k=1}^{|c|-1}
P_1\big(\lfrac{k+m}{c}\big) P_1\big(\lfrac {dk}{c}\big)
\end{equation}
Using \eqref{P1Cot} and the Fourier inversion formula
\eqref{FourierInversion} we can immediately identify two other
terms in \eqref{GenDedRel:1}, namely
\begin{equation}\label{GenDedRel:3}
\lfrac 1{|c|}\sum_{p=1}^{|c|-1} \lfrac i2 \cot\big(\lfrac{\pi
dp}{c}\big)\xi^{dmp} = P_1\big(\lfrac{m}{|c|}\big) ,\quad
\lfrac1{|c|}\sum_{p=1}^{|c|-1} \lfrac i2 \cot\big(\lfrac{\pi
p}{c}\big)\xi^{dmp} = P_1\big(\lfrac{dm}{|c|}\big).
\end{equation}
Moreover,
\[
\sum_{p=1}^{|c|-1}\xi^{dmp} = |c|\big(1-\gd(\lfrac{dm}c)\big) - 1
= |c|\big(1-\gd(\lfrac{m}c)\big) - 1,
\]
where in the last equality we have used that $\gcd(d,c)=1$.
Employing this together with \eqref{GenDedRel:2} and
\eqref{GenDedRel:3}, we can rewrite \eqref{GenDedRel:1} as
\begin{equation}\label{GenDedRel:4}
\begin{split}
s_{x,y}(a,c) = & \lfrac 1{|c|} P_1(x)^2 + \sum_{k=1}^{c-1}
P_1\big(\lfrac{k+m}{c}\big) P_1\big(\lfrac {dk}{c}\big)\\
&\qquad + \lfrac 1{4}\gd(x)\Big(1- \gd(\lfrac{m}c) - \lfrac
1{|c|}\Big) + \lfrac 12 \gd(x)\Big( P_1\big(\lfrac{m}{|c|}\big) -
P_1\big(\lfrac{dm}{|c|}\big)\Big).
\end{split}
\end{equation}
To find the claimed formula for $s_{x,y}(a,c) - s(a,c)$, let us
first study the difference
\begin{equation}\label{GenDedRel:5}
\sum_{k=1}^{c-1} P_1\big(\lfrac{k+m}{c}\big) P_1\big(\lfrac
{dk}{c}\big) - s(a,c) = \sum_{k=1}^{c-1}P_1\big(\lfrac
{dk}{|c|}\big)\Big(P_1\big(\lfrac{k+m}{|c|}\big) -
P_1\big(\lfrac{k}{|c|}\big)\Big),
\end{equation}
where we have used that we can replace $c$ with $|c|$ in the
first term and that according to \eqref{DedSymm} we have
$s(a,c)=s(d,|c|)$. Now, with $r\in\{0,\ldots,|c|-1\}$ satisfying
$m\equiv r\,(c)$, we get
\[
P_1\big(\lfrac{k+m}{|c|}\big) - P_1\big(\lfrac{k}{|c|}\big)
=\begin{cases}
\hphantom{P_1}\lfrac r{|c|}, &\text{if $k<|c|-r$},\\
P_1\big(\lfrac{r}{|c|}\big), &\text{if $k=|c|-r$},\\
\,\lfrac r{|c|}-1, &\text{if $k>|c|-r$.}
\end{cases}
\]
Hence, straightforward manipulations yield that
\eqref{GenDedRel:5} equals
\[
\begin{split}
\big(\lfrac r{|c|} - 1\big)\sum_{k=1}^{|c|-1}&P_1\big(\lfrac
{dk}{|c|}\big) + \sum_{k=1}^{|c|-r-1}P_1\big(\lfrac
{dk}{|c|}\big) - \big(\lfrac r{|c|} -
1\big)P_1\big(\lfrac{d(|c|-r)}{|c|}\big) +
P_1\big(\lfrac{d(|c|-r)}{|c|}\big) P_1\big(\lfrac {r}{|c|}\big) \\
&= \sum_{k=1}^{|c|-r}P_1\big(\lfrac {dk}{|c|}\big) +
P_1\big(\lfrac{d(|c|-r)}{|c|}\big)\Big(P_1\big(\lfrac
{r}{|c|}\big) -\lfrac {r}{|c|}\Big)\\
&= \sum_{k=1}^{|c|-r}P_1\big(\lfrac {dk}{|c|}\big) + \lfrac 12
P_1\big(\lfrac {dm}{|c|}\big).
\end{split}
\]
Note that we have used \eqref{SumFullRepSys} to see that
$\sum_{k=1}^{|c|-1}P_1\big(\lfrac {dk}{|c|}\big)=0$. Combining
this computation with \eqref{GenDedRel:4}, we arrive at
\[
\begin{split}
s_{x,y}(a,c) -s(a,c)= & \lfrac 1{|c|} P_1(x)^2 - \lfrac
1{4|c|}\gd(x) + \sum_{k=1}^{|c|-r}P_1\big(\lfrac {dk}{|c|}\big) +
\lfrac 12 P_1\big(\lfrac {dm}{|c|}\big)  \\
&\qquad + \lfrac 1{4}\gd(x)\big(1-\gd(\lfrac{m}c)\big)  + \lfrac
12 \gd(x)\Big( P_1\big(\lfrac{m}{|c|}\big) -
P_1\big(\lfrac{dm}{|c|}\big)\Big),
\end{split}
\]
which is the claimed formula, since $P_1(x)^2 - \lfrac 1{4}\gd(x)
= x^2-x = P_2(x) -\lfrac 16$.
\end{proof}

\cleardoublepage
\chapter{Local Variation of the Eta Invariant}\label{AppEtaVar}

We include some more analytical details concerning the heat
operator. Specifically, we will give some remarks concerning
families of heat operators, and use this to give a proof of
Proposition \ref{EtaDiffSF}.

\section{More Results on the Heat Operator}\label{MoreHeat}

In Chapter \ref{Signature} we have defined the heat operator
$e^{-tH}$ for a formally self-adjoint elliptic differential
operator $H$ of order 2 with positive definite leading symbol. We
have done this by explicitly writing it as the limit of integral
operators with smooth kernels, which were defined in terms of a
spectral decomposition of $H$. In Lemma \ref{BasicKernelEst}, we
have derived a basic estimate for the heat kernel for large times.
Concerning the small time behaviour of the heat operator, it is
useful to have a description of $e^{-tH}$ in terms of the
resolvent $(H-z)^{-1}$ for $z\notin\spec(H)$. In particular, the
proof of the asymptotic expansion in Theorem \ref{HeatTrace} as
in \cite{Gil} uses this description.

\subsection{Expression via the Resolvent}

\begin{figure}[htbp]
\centering
\includegraphics[width=0.6\linewidth]{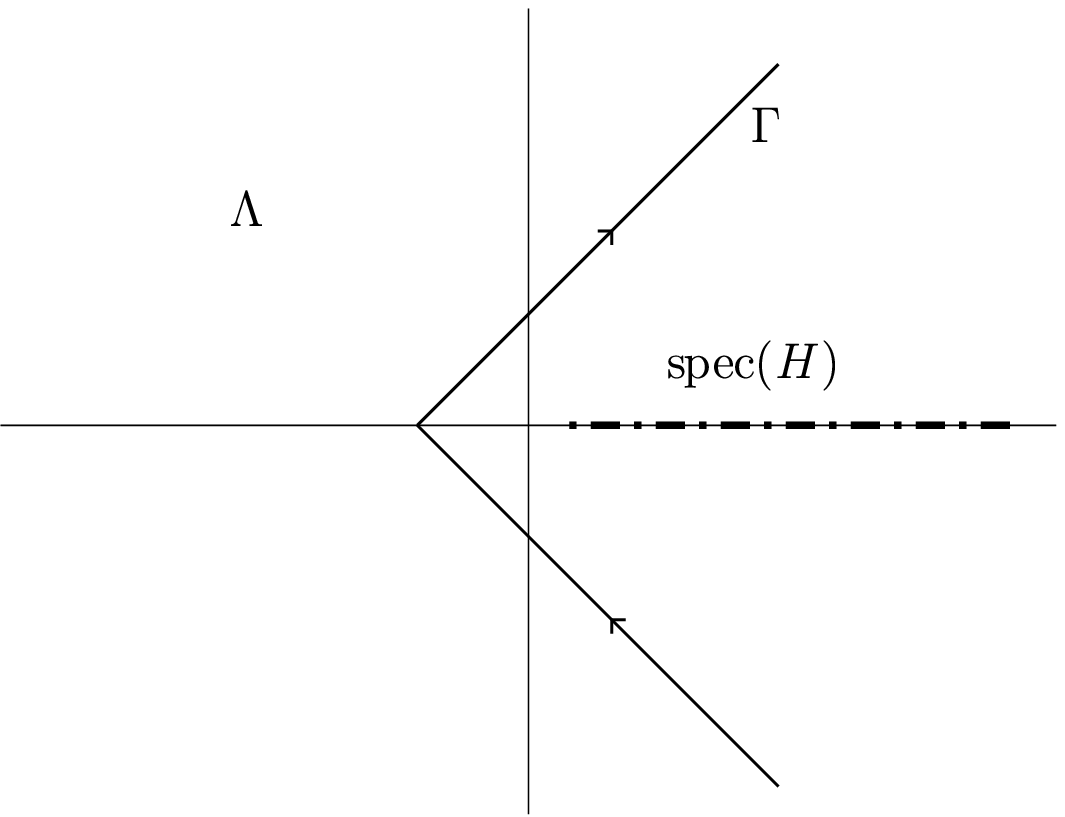}
\caption{The contour $\gG$}\label{Fig:Contour1}
\end{figure}

Let $E$ be a vector bundle over a closed manifold $M$ of
dimension $m$, and let $H$ be a formally self-adjoint elliptic
differential operator of order 2 with positive definite leading
symbol. As before, for $s,s'\in \R$, let $\sB(L^2_s,L^2_{s'})$
denote the space of bounded operators from $L_s^2$ to $L^2_{s'}$
endowed with the operator norm $\|.\|_{s,s'}$. Although this is
not really necessary, we assume for simplicity that $H$ is
non-negative, i.e., $\spec(H)\subset [0,\infty)$, which is
certainly true for an operator of the form $H=D^2$. Then the
elliptic estimate \eqref{Garding} implies that for every $s\ge 0$,
there exists a constant $C$ such that
\begin{equation}\label{SobolevScale}
C^{-1} \big\| (\Id+H)^{s/2} \gf\big\|_{L^2}\le \|\gf\|_{L^2_s} \le
C \big\| (\Id+H)^{s/2} \gf\big\|_{L^2},\quad \gf\in C^\infty(M,E).
\end{equation}
Consider the region
\begin{equation}\label{RegionDef}
\gL :=\bigsetdef{z\in \C}{\Re(z)+ 1\le |\Im(z)|}.
\end{equation}
Then there exists a constant $C>0$ such that
\begin{equation*}
\dist(z,\spec(H))\ge C|z|\quad\text{for all $z\in \gL$}.
\end{equation*}
This, together with \eqref{SobolevScale} and the spectral theorem,
implies that for $s\ge 0$, $z\in\gL$, and $\gf\in C^\infty(M,E)$
\[
\begin{split}
\|(H-z)^{-1}\gf\|_{L_s^2}&\le C_1
\|(\Id+H)^{s/2}(H-z)^{-1}\gf\|_{L^2}\\
&\le C_2|z|^{-1}\, \|(\Id+H)^{s/2}\gf\|_{L^2} \le  C_3 |z|^{-1}\,
\|\gf\|_{L^2_s}
\end{split}
\]
where we have used that $(\Id+H)$ and $(H-z)^{-1}$ commute. By
duality, this also holds for $s<0$. Thus, for every $s\in \R$,
there exists a constant $C$ such that
\begin{equation}\label{BasicResolventEst:1}
\|(H-z)^{-1}\|_{s,s} \le C |z|^{-1},\quad z\in \gL.
\end{equation}
We now consider the contour $\gG:=\pd \gL$, oriented as the
boundary, i.e., in such a way that $[0,\infty)$ lies in the
interior, see Figure \ref{Fig:Contour1}. Then the Cauchy formula
implies that

\begin{equation}\label{HeatResolventRel}
e^{-tH} = \frac 1{2\pi i} \int_\gG e^{-tz} (H-z)^{-1}dz.
\end{equation}

Because of \eqref{BasicResolventEst:1}, this expression converges
in the operator norm in $L^2_s(M,E)$ for every $s\in\R$. To get
norm estimates in $\sB(L_s^2,L_{s+l}^2)$ for $l>0$ we need the
following basic resolvent estimate.

\begin{prop}
Let $s\in \R$, $k\in\N$, and $0\le l\le 2k$. Then there exists a
constant $C$ such that for all $z\in \gL$
\begin{equation}\label{BasicResolventEst:2}
\|(H-z)^{-k}\|_{s,s+l} \le C |z|^{l/2-k}.
\end{equation}
\end{prop}
\begin{proof}
For $z\in \gL$, the spectrum of the operator $(H-z)^{-1}$ is not
contained in the negative real line, so that $(H-z)^{-r}$ can be
defined for every $r\ge 0$ by the spectral theorem. As in
\eqref{BasicResolventEst:1} one then finds
\[
\|(H-z)^{-r}\|_{s,s}\le C|z|^{-r}.
\]
Then,
\[
\begin{split}
\|(H-z)^{-k}\|_{s,s+l} &\le C_1
\big\|(\Id+H)^{l/2}(H-z)^{-k}\big\|_{s,s} \\
&\le C_1 \big\|\big((\Id+H)(H-z)^{-1}\big)^{l/2}\big\|_{s,s}
\big\|(H-z)^{l/2-k}\big\|_{s,s}\\
&\le C_2 |z|^{l/2-k} \Big(1 +
\big\|(\Id+z)(H-z)^{-1}\big\|_{s,s}^k\Big)\\
&\le C_2 |z|^{l/2-k} \Big(1+ C_3\lfrac{|z+1|}{|z|}\Big)^k.
\end{split}
\]
The last factor can be bounded, independently of $z\in\gL$. This
proves \eqref{BasicResolventEst:2}.
\end{proof}

An important consequence of \eqref{BasicResolventEst:2} is that
$e^{-tH}$ does indeed solve the heat equation. We give the
following summary of well-known facts, see in particular
\cite[Lem. 1.7.5]{Gil}.

\begin{prop}\label{HeatSolve}
Let $H$ be a non-negative operator in $\sP^2_{s,e}(M,E)$ and let
$s \in\R$.
\begin{enumerate}
\item The one-parameter family $(e^{-tH})_{t\in(0,\infty)}$ is a smooth
family\footnote{See Definition \ref{SmoothFamilySmoothing}
below.} of smoothing operators. If $\gf\in C^\infty(M,E)$, then
\[
(\lfrac d{dt} + H)e^{-tH}\gf = 0,\quad \text{and}\quad \lim_{t\to
0}\|e^{-tH}\gf - \gf\|_{L^2_s}=0.
\]
\item The collection $e^{-tH}$ forms a semi-group, i.e.,
\[
e^{-(t+t')H} = e^{-tH}e^{-t'H},\quad t,t'>0.
\]
\item For $l\ge 0$, there exists a constant $C>0$
such that for $t\in (0,1)$
\begin{equation}\label{BasicSmallTimeEst}
\|e^{-tH}\|_{s,s+l} \le C t^{-l/2}.
\end{equation}
\item Let $c>0$ be smaller than the smallest non-zero eigenvalue of
$H$ and let $t_0>0$. Then there exists a constant $C>0$ such that
for all $t\ge t_0$
\begin{equation}\label{BasicLargeTimeEst}
\|e^{-tH}-P_0\|_{s,s}\le C e^{-ct},
\end{equation}
where $P_0$ is the projection onto $\ker H$.
\end{enumerate}
\end{prop}

\begin{figure}[htbp]
\centering
\includegraphics[width=0.6\linewidth]{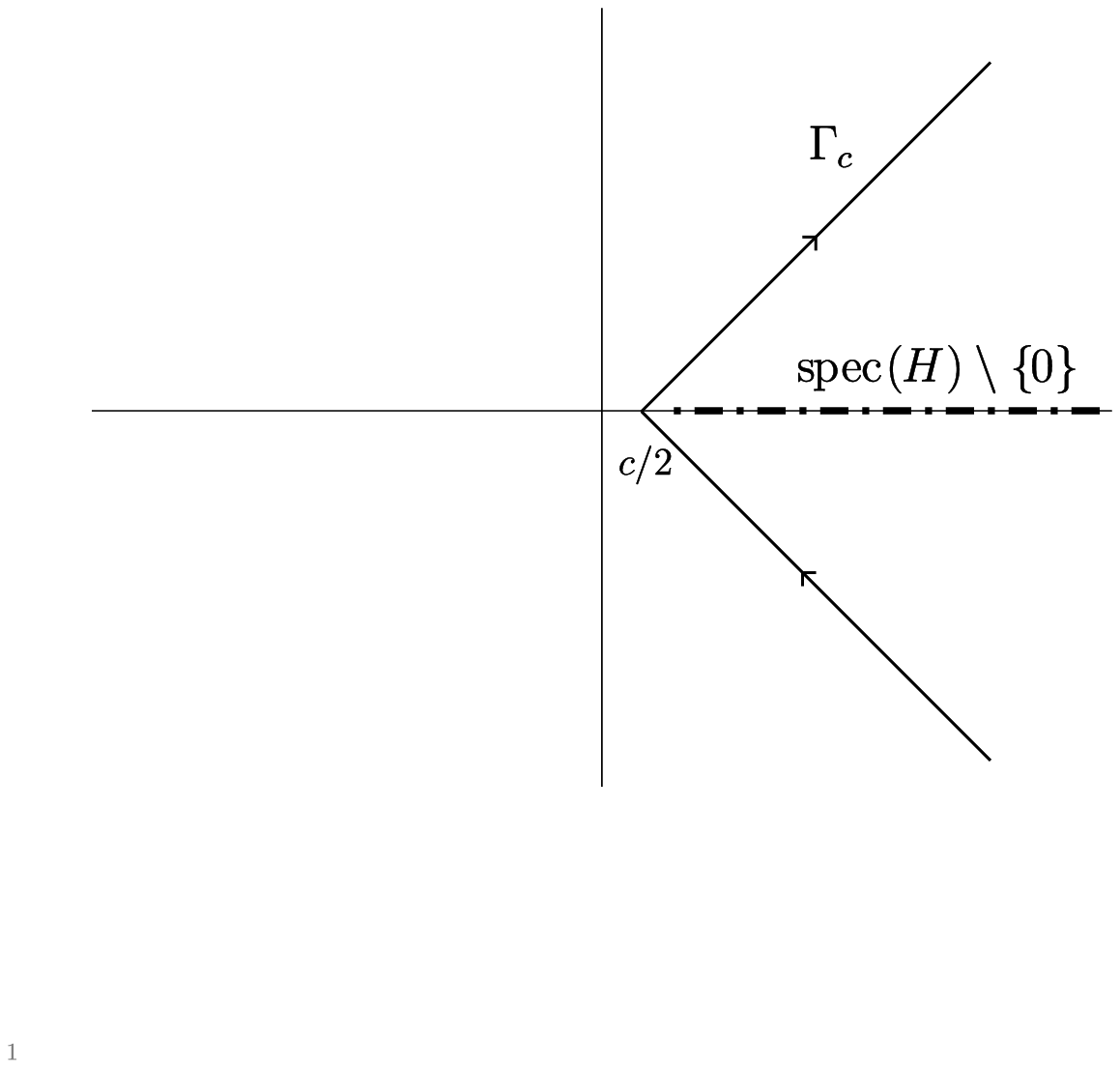}
\caption{The contour $\gG_c$}\label{Fig:Contour2}
\end{figure}

\begin{proof}[Sketch of proof]
The assertion (i) is \cite[Lem. 1.7.5]{Gil}. We skip the proof
since it uses the same ideas as the proof of part (iii). Part
(ii) is also standard: Let $\eps>0$, and consider the contour
$\gG' := \gG -\eps$. Then, using Cauchy's formula and the
resolvent equation, one finds that
\[
\begin{split}
e^{-t'H}e^{-t H} &= -\frac1{4\pi^2} \int_{\gG'}\int_\gG
e^{-t'z'}e^{-t z }(H-z)^{-1}(H-z')^{-1}dz dz'\\
&= -\frac1{4\pi^2} \int_{\gG'}\int_\gG \frac{e^{-tz
-t'z'}}{z-z'}\big[(H-z)^{-1}-(H-z')^{-1}\big]dz dz'.
\end{split}
\]
It now follows from standard ``$\C$-valued'' complex analysis,
that first term is equal to $e^{-(t'+t)H}$, while the second term
vanishes. To prove \eqref{BasicSmallTimeEst} choose $k\in \N$
with $l\le 2k$. We integrate \eqref{HeatResolventRel} by parts to
find that
\begin{equation}\label{HeatResolventRel:2}
e^{-tH} = \frac 1{2\pi i}\frac {(k-1)!}{t^{k-1}} \int_\gG e^{-tz}
(H -z)^{-k}dz.
\end{equation}
Now let $t\in (0,1)$ and substitute $\gz=tz$ in
\eqref{HeatResolventRel:2}. Then we can use Cauchy's theorem to
change integration over $t\gG$ back to integration over $\gG$.
This shows that for $t\in(0,1)$
\[
e^{-tH} = \frac 1{2\pi i}\frac {(k-1)!}{t^k} \int_\gG e^{-\gz}
(H-\gz/t)^{-k}d\gz.
\]
Now \eqref{BasicResolventEst:2} applied to $z=\gz/t\in \gL$ shows
that
\[
\|e^{-tH}\|_{s,s+l} \le C_1 t^{-k} \Big(\int_\gG
|e^{-\gz}|\,|\gz/t|^{l/2-k}
 d\gz\Big)\le C_2 |t|^{-l/2}.
\]
This proves \eqref{BasicSmallTimeEst}. Concerning (iii) we only
note that the large time estimate can be easily deduced from
Lemma \ref{BasicKernelEst}. Alternatively it can be proved using
\eqref{HeatResolventRel:2} by considering the contour
\[
\gG_c := \bigsetdef{z\in \C}{\Re(z)-c/2=|\Im(z)|},
\]
and employing estimates corresponding to the ones in
\eqref{BasicResolventEst:2}, see Figure \ref{Fig:Contour2}
\end{proof}

%

\subsection{Perturbed Operators}

Let $E$ be a vector bundle over a closed manifold $M$ of
dimension $m$, and let $K$ be an integral operator on $L^2(M,E)$
with smooth kernel $k(x,y)\in C^\infty(M\times M,E\boxtimes E^*)$.
Then $K$ is a smoothing operator and thus,
\[
K\in \sB\big(L_s^2,L_{s'}^2\big)\quad\text{for all $s, s'\in \R$}.
\]
Moreover, the $C^k$-norms of $k(x,y)$ can be controlled by the
operator norms of $K$. More precisely, for each $k\in \N$, we can
find a constant $C$ as in \cite[Lem. 1.2.7]{Gil} such that
\begin{equation}\label{KernOpEst}
\|k(x,y)\|_{C^k}\le C \,\|K\|_{-l,l},\quad l> k+m/2.
\end{equation}

\begin{remark*}
Conversely, the Schwartz Kernel Theorem (see e.g. \cite[Sec.
4.6]{Tay:I}) ensures that an operator $K$ on $L^2(M,E)$, which
satisfies $\|K\|_{-l,l}<\infty$ for $l> k+m/2$, has a kernel
$k(x,y)$ of class $C^k$ such that \eqref{KernOpEst} holds.
\end{remark*}

As before, assume that $H\in \sP_{s,e}^2(M,E)$ is non-negative. If
$K$ is a symmetric smoothing operator on $C^\infty(M,E)$, the
operator $H+K$ will in general not be an elliptic differential
operator. Nevertheless, it follows from standard perturbation
theory that $H+K$ has all the properties, we have obtained in
Theorem \ref{EllOpSpec} for the unperturbed case. In particular,
we get a well-defined heat operator $e^{-t(H+K)}$. We need to
have a control on the difference $e^{-t(H+K)} - e^{-tH}$. The
following result is basically \cite[Prop. 9.46]{BGV}, only that
we have changed it into a statement about operators rather than
kernels.

\begin{prop}\label{Volterra}
Let $K$ be a symmetric smoothing operator on $C^\infty(M,E)$, and
let $H\in \sP_{s,e}^2(M,E)$ be non-negative. For $k\ge 1$ and
$t>0$ define inductively
\[
K_0(t) := e^{-tH},\quad K_k(t):= \int_0^te^{-(t-s)H}K K_{k-1}(s)
ds.
\]
Then, if $l\in\N$, there exists a constant $C>0$ such that for all
$k\ge 1$
\[
\|K_k(t)\|_{-l,l} \le \frac {C^kt^k}{k!}\,\|K\|_{-l,l}^k,
\]
Moreover, for each $N\ge 1$ there exists $C>0$ such that for all
$t\in (0,1)$
\[
\Big\|e^{-t(H+K)} -e^{-tH} - \sum_{k=1}^N (-1)^k
K_k(t)\Big\|_{-l,l} \le Ct^{N+1}.
\]
\end{prop}

\begin{remark*}
Before we sketch the proof we want to point out that for $k\ge 1$
the operator $K_k(t)$ can also be described as follows. Let
\[
\gD_k:=\bigsetdef{(s_1,\ldots,s_k)\in \R^k}{0\le s_1\le\ldots\le
s_k =1}
\]
be the standard $k$-simplex. Then
\begin{equation}\label{VolterraAlternative}
K_k(t)= t^k\int_{\gD_k} e^{-t(1-s_k) H}Ke^{-t(s_k-s_{k-1})H}\dots
Ke^{-t(s_2-s_1)  H} K e^{-ts_1  H} ds.
\end{equation}
\end{remark*}

\begin{proof}[Sketch of proof]
We use \eqref{VolterraAlternative} to prove the estimate on
$K_k(t)$. First it follows from \eqref{BasicSmallTimeEst} and
\eqref{BasicLargeTimeEst} that $\|e^{-t H}\|_{s,s}$ can be
bounded, independently of $t$. Thus,
\[
\|K e^{-ts_1  H}\|_{-l,l} \le \|K\|_{-l,l} \|e^{-ts_1
H}\|_{-l,-l} \le C \|K\|_{-l,l},
\]
On the other hand,
\[
\big\|e^{-t(1-s_k) H}Ke^{-t(s_k-s_{k-1})}\dots Ke^{-t(s_2-s_1)
H}\big\|_{l,l}\le C^{k} \|K\|_{l,l}^{k-1},
\]
Since $\|K\|_{l,l}\le C \|K\|_{-l,l}$ we deduce,
\[
\|K_k(t)\|_{-l,l} \le t^k \vol(\gD_k) C^{k} \|K\|_{-l,l}^k\le
\frac {C^kt^k}{k!}\,\|K\|_{-l,l}^k,
\]
where we have used that $\vol(\gD_k)= \lfrac 1{k!}$. Now the
definition of $K_k(t)$ shows that
\[
\lfrac d{dt} K_k(t)= K K_{k-1}(t) -H K_k(t).
\]
Thus, for all $N\ge 1$,
\[
\big(\lfrac d{dt} +H+K\big)\sum_{k=0}^N (-1)^kK_k(t) = (-1)^N K
K_N(t).
\]
Using the estimate on $\|K_k(t)\|_{-l,l}$ we see that
\[
\Big\|\big(\lfrac d{dt} +H+K\big)\sum_{k=0}^N
(-1)^kK_k(t)\Big\|_{-l,l}\le C_N t^N
\]
Thus, $\sum_{k=0}^N (-1)^kK_k(t)$ is an approximate solution to
the heat equation in terms of $H+K$ in the sense of \cite[Sec.
2.4]{BGV}. This implies that
\[
\Big\|e^{-t(H+K)} -e^{-tH} - \sum_{k=1}^N (-1)^k
K_k(t)\Big\|_{-l,l} \le Ct^{N+1}.\qedhere
\]
\end{proof}

We also need a result on the heat trace asymptotics when we
perturb the operator $H$ by a smoothing operator $K$. We first
borrow the following from \cite[Prop. 2.47]{BGV}.

\begin{prop}\label{TraceAsympSmoothingOp}
Let $K$ be a smoothing operator on $L^2(M,E)$. Then there exists
an asymptotic expansion of the form
\[
\Tr(e^{-tH}K)\sim \Tr(K) + \sum_{n=0}^\infty a_n t^n,\quad\text{as
$t\to 0$.}
\]
\end{prop}

\begin{remark*}
The proof in \cite{BGV} relies on the explicit description of the
heat kernel by geometrically constructed approximations. This
method to obtain the heat trace asymptotics is different from the
one presented in \cite{Gil}, which is the one we are following in
our presentation. However, a rough idea to prove Proposition
\ref{TraceAsympSmoothingOp} with the methods already established
is the following: For each $N\in \N$ let
\[
K_N(t):= \sum_{n=0}^N \frac{(-t)^n}{n!}H^nK.
\]
Since $K$ is a smoothing operator, $K_N(t)$ is a smoothing
operator as well. Applying the heat equation to $K_N(t)$ yields
\[
\big(\lfrac d{dt}+H\big) K_N(t) = \frac{(-t)^N}{N!}H^{N+1}K.
\]
Since we can bound $\|H^{N+1}K\|_{-l,l}$ for fixed $N$, this
shows that $K_N(t)$ is an approximate solution to the heat
equation. For $t\to 0$ it converges strongly to $K$ which implies
that
\[
\big\|e^{-tH}K-K_N(t)\big\|_{-l,l}\le C t^{N+1}.
\]
Using \eqref{KernOpEst} and the expression of the trace in terms
of kernels, we can then estimate
\[
\big|\Tr\big(e^{-tH}K - K_N(t)\big)\big|\le C t^{N+1}.
\]
From this one finds that the assertion holds with
\[
a_n:= \frac{(-1)^n}{n!}\Tr(H^nK).
\]
\end{remark*}

The following result is what we were aiming at in this section.
We will give a version which suffices for our considerations,
although a more general statement should be possible.

\begin{prop}\label{TraceAsympPerturb}
Let $H$ and $K$ be as before, and assume in addition that $K$ and
$H$ commute. Let $D$ be an auxiliary formally self-adjoint
differential operator of order $d$, and let $K'$ be a symmetric
smoothing operator on $C^\infty(M,E)$. Then there exists an
asymptotic expansion
\[
\Tr\big((D+K')e^{-t(H+K)} \big) \sim \sum_{n=0}^\infty
t^{\frac{n-m-d}{2}}a_n ,\quad\text{as }t\to 0.
\]
Moreover, if we denote by $a_n(D,H)$ the coefficients of the
asymptotic expansion of $\Tr(De^{-tH })$, we have
\[
a_n = a_n(D,H),\quad\text{for $n\le m$}.
\]
\end{prop}

\begin{proof}[Sketch of proof]
According to Theorem \ref{HeatTrace} it suffices to check that
\[
\Tr\big((D+K')e^{-t(H+K)} \big) - \Tr\big(De^{-tH}\big)
\]
has an asymptotic expansion as $t\to 0$ as a power series in $t$.
The assumption that $K$ and $H$ commute simplifies the situation
considerably, since then
\[
e^{-t(H+K)} = e^{-tH}\sum_{n=0}^\infty \frac{(-t)^n}{n!}K^n,
\]
where the series converges in every $L_{-l,l}^2$, because $K$ is
smoothing and thus bounded. This reduces the claim to the study
of terms of the form
\begin{equation}\label{TraceAsympPerturb:1}
\Tr\big(D e^{-tH}K^n\big)\quad\text{and}\quad \Tr(K' e^{-tH}K^n),
\end{equation}
where $n\ge 1$. Now the trace property shows that both are of the
form $\Tr\big(e^{-tH}\tilde K\big)$ for some smoothing operator
$\tilde K$. We can thus apply Proposition
\ref{TraceAsympSmoothingOp} to deduce that the terms in
\eqref{TraceAsympPerturb:1} are indeed asymptotic to power series
in $t$.
\end{proof}

\subsection{Variation of the Heat Operator}

Now let $(K_u)_{u\in U}$ be a $p$-parameter family of smoothing
operators. The constant $C$ in the estimate \eqref{KernOpEst} is
independent of $u$. This implies that if $K_u$ depends smoothly
on $u$ with respect to $\|.\|_{-l,l}$ for each $l\in\N$, then the
family of kernels $k_u(x,y)$ will depend smoothly on $u$ with
respect to all $C^k$-norms, where $l>k+m/2$. This motivates the
following definition.

\begin{dfn}\label{SmoothFamilySmoothing}
Let $K_u$ be a one-parameter family of operators with smooth
kernels. Then $K_u$ is called a \emph{smooth family of smoothing
operators}, if for all $l\in \N$, the assignment $u\mapsto K_u$ is
smooth in $\sB\big(L_{-l}^2,L_l^2\big)$.
\end{dfn}

We can now formulate the following version of Duhamel's formula,
compare with \cite[Thm. 2.48]{BGV}.

\begin{theorem}\label{DuhamelThm}
Consider a one-parameter family $H_u$ in $\sP^2_{s,e}(M,E)$ of
non-negative operators, and assume that $H_u$ is smooth in the
sense of Definition \ref{SmoothFamily}. Then $e^{-tH_u}$ is a
smooth family of smoothing operators, and
\begin{equation}\label{Duhamel}
\lfrac{d}{du} e^{-tH_u} = -\int_0^t
e^{-(t-s)H_u}(\lfrac{d}{du}H_u)e^{-sH_u}ds.
\end{equation}
\end{theorem}

\begin{remark*}
Note that in contrast to the discussion in \cite[Sec. 2.7]{BGV}
the smoothness in $u$ does not follow from our description
\eqref{HeatKernel} of the heat kernel since in general the
eigenvalues and eigenvectors of $H_u$ will not depend smoothly on
$u$. However, we can use the description of $e^{-tH_u}$ in terms
of the resolvent to obtain the result.
\end{remark*}

\begin{proof}
First note that the basic resolvent estimates
\eqref{BasicResolventEst:1} and \eqref{BasicResolventEst:2} can
be made uniform in $u$ since the constants appearing there depend
on $u$ only through the elliptic estimate, which can be made
locally uniform in $u$.

We now want to prove that $(H_u-z)^{-1}$ varies smoothly with
$u$. Without loss of generality we consider an interval around 0.
First of all, let $\gG$ be the contour defined as the boundary of
\eqref{RegionDef}. Then for all $z\in\gG$,
\begin{equation}\label{ResolventFormula}
H_u-z = (H_0-z)\big(\Id + (H_0-z)^{-1}(H_u-H_0)\big)=:
(H_0-z)(\Id + T_u).
\end{equation}
It follows from our smoothness assumption that $T_u$ is a bounded
operator on each $L_l^2$ for every choice of  $l\in\Z$, and that
the assignment
\[
\R\to\sB(L_l^2,L_l^2),\quad u\mapsto T_u
\]
is smooth. Moreover, we can choose $\gd>0$ such that for all
$u\in(-\gd,\gd)$
\[
\|T_u\|_{n,n} \le \lfrac 12 \quad \text{for } n\in\{-l,\ldots,l\}.
\]
Using the Neumann series, this ensures that $\Id + T_u$ is
invertible in each $\sB(L_n^2,L_n^2)$. Furthermore, the
assignment $u\mapsto (\Id +T_u)^{-1}$ is differentiable with
\[
\lfrac d{du} (\Id +T_u)^{-1} =(\Id +T_u)^{-1} \big(\lfrac
d{du}T_u\big)(\Id +T_u)^{-1}.
\]
Inductively, one finds that $(\Id +T_u)^{-1}$ is smooth in $u$.
Now \eqref{ResolventFormula} shows that for all $z\in\gG$
\[
(H_u-z)^{-1} = (\Id +T_u)^{-1} (H_0-z)^{-1}.
\]
This yields that the family $(H_u-z)^{-1}$ depends smoothly on $u$
as a map to $\sB(L_{n-2}^2,L_n^2)$, where as before
$n\in\{-l,\ldots,l\}$. Hence, for $k$ large enough, the assignment
\[
(-\gd,\gd)\to \sB(L_{-l}^2,L_l^2),\quad u\mapsto (H_u-z)^{-k}
\]
is well-defined and smooth in $u$. Moreover, one easily finds that
\begin{equation}\label{ResolventPowerDiff}
\lfrac d{du} (H_u-z)^{-k} = - \sum_{n=1}^k
(H_u-z)^{-n}\big(\lfrac d{du}H_u\big)(H_u-z)^{n-k-1}.
\end{equation}
Since the basic resolvent estimate \eqref{BasicResolventEst:2}
can be made uniform in $u$, we deduce from
\eqref{HeatResolventRel:2} that $e^{-tH_u}$ is differentiable and
that we can differentiate under the integral using
\eqref{ResolventPowerDiff}. Inductively, one finds that
$e^{-tH_u}$ is a smooth map to $\sB(L_{-l}^2,L_l^2)$.

Having established the smoothness of $e^{-tH_u}$ in $u$, we can
prove Duhamel's formula \eqref{Duhamel} along the same lines as
in \cite[Sec. 2.7]{BGV}: First, the heat equation implies that
\[
(\lfrac d{dt} +  H_u)\lfrac d{du}e^{-tH_u} =-(\lfrac d{du}H_u)
e^{-tH_u}.
\]
On the other hand, one finds that
\[
(\lfrac d{dt} +  H_u) \int_0^t
e^{-(t-s)H_u}(\lfrac{d}{du}H_u)e^{-sH_u}ds = (\lfrac d{du}H_u)
e^{-tH_u}.
\]
Thus,
\begin{equation}\label{Duhamel:3}
\lfrac{d}{du} e^{-tH_u} + \int_0^t
e^{-(t-s)H_u}(\lfrac{d}{du}H_u)e^{-sH_u}ds
\end{equation}
solves the heat equation. We are thus left to show that
\eqref{Duhamel:3} converges to 0 in $L^2(M,E)$ as $t\to 0$.
Concerning the first term, we use
\[
\lfrac{d}{du} e^{-tH_u} = - \frac  1{2\pi i} \int_\gG e^{-tz}
(H_u-z)^{-1}\big(\lfrac d{du}H_u\big)(H_u-z)^{-1}dz
\]
and the resolvent estimate \eqref{BasicResolventEst:1}. Then, as
in the proof of \eqref{BasicSmallTimeEst}, one easily finds small
time estimates
\[
\|\lfrac{d}{du} e^{-tH_u}\gf\|_{L^2}\le C t
\|\gf\|_{L^2},\quad\gf\in L^2(M,E),
\]
which are uniform in $u$ and $t\in (0,1)$. This shows that
$\lfrac{d}{du} e^{-tH_u}$ converges to 0 as $t\to 0$ in
$\sB(L^2,L^2)$. For the second term in \eqref{Duhamel:3}, we can
argue as in the proof of Proposition \ref{Volterra} to get
uniform bounds on the integrand as $t\to 0$. This implies our
assertion.
\end{proof}

In a similar way one can show the following result.

\begin{prop}\label{ProjGreenSmooth}
Let $D_u$ be a smooth family in $\sP_{s,e}^d(M,E)$ and assume
further that $\dim(\ker D_u)$ is constant.
\begin{enumerate}
\item The projection $P_u$ onto $\ker(D_u)$ is a smooth family of
smoothing operators.
\item Denote by $G_u$ the family of Green's operators of $D_u^2$,
defined by
\[
G_u|_{\ker(P_u)}= \big(D_u|_{\ker(P_u)}\big)^{-1},\quad
G_u|_{\im(P_u)} =0.
\]
Then, for every $s\in \R$, the family $G_u$ is smooth in $u$ as a
map to $\sB(L_s^2,L_{s+d}^2)$.
\end{enumerate}
\end{prop}

\begin{proof}[Sketch of proof]
Since $D_u$ varies smoothly with $u$, one can verify---using the
smooth dependence of the resolvent on $u$---that for compact
$u$-intervals, the non-zero eigenvalues are uniformly bonded away
from 0. Then, according to the spectral theorem, the projection
$P_u$ onto $\ker D_u$ is given by
\begin{equation}\label{ContIntProj}
P_u= \frac 1{2\pi i}\int_{\gG} (D_u-z)^{-1} dz,
\end{equation}
where $\gG$ is the clockwise oriented boundary of a small disk
$B$ such that
\[
B\cap\spec(D_u) = \{0\},\quad\text{for all $u$.}
\]
Integrating by parts, we get
\[
P_u= \frac 1{2\pi i}\int_{\gG} (-z)^{k-1} (D_u-z)^{-k} dz.
\]
Resolvent estimates as in \eqref{BasicResolventEst:1} and
\eqref{BasicResolventEst:2} and the formula
\eqref{ResolventPowerDiff} for the derivative of the resolvent
with respect to $u$ then yield that $P_u$ is a smooth family of
smoothing operators, see \cite[Prop. 9.10]{BGV} for a related
discussion.

Concerning the family of Green's operators, we can take the same
small disk $B$ and orient the boundary $\gG$ counter-clockwise.
Then
\begin{equation}\label{GreensOp}
G_u= \frac 1{2\pi i}\int_{\gG} z^{-1}(D_u-z)^{-1} dz.
\end{equation}
Now the uniform bound $\|(D_u-z)^{-1}\|_{s,s+d}\le C$ and the
formula for the derivative of the resolvent show that $G_u$ is a
smooth family of operators in $\sB(L_s^2,L_{s+d}^2)$.
\end{proof}

\begin{remark}\label{SmoothDepRem}\quad\nopagebreak
\begin{enumerate}
\item Note that in \eqref{GreensOp} we cannot integrate by parts
to increase the regularity of $G_u$. This is, of course, already
clear in the case of a single operator.
%
%
\item If there exists a constant $c$ such that $c\notin
\spec(|D_u|)$ for all $u$ in some interval $[-\gd,\gd]$, we can
use \eqref{ContIntProj}---with $\gG$ being a circle of radius $c$
around 0---to deduce that the spectral projection onto all
eigenvalues of norm less than $c$ is smooth $u$, see also
\cite[Prop. 9.10]{BGV}.
\item Theorem \ref{DuhamelThm} and Proposition
\ref{ProjGreenSmooth} carry over verbatim, if we consider a smooth
$p$-parameter families of formally self-adjoint elliptic
operators. This is important for the discussion of fiber bundles.
\end{enumerate}
\end{remark}

We also need the following consequence of Theorem \ref{DuhamelThm}

\begin{lemma}
Let $H_u$ be a smooth one-parameter family of non-negative
operators in $\sP^2_{s,e}(M,E)$, and let $D_u$ be an auxiliary
smooth one-parameter family of formally self-adjoint differential
operators of order $d$. Then $D_ue^{-tH_u}$ is a differentiable
family of trace-class operators and, if $D_u$ commutes with $H_u$,
\begin{equation}\label{DuhamelCor}
\lfrac d{du} \Tr\big(D_ue^{-tH_u}\big) = \Tr\big(\lfrac
{dD_u}{du} e^{-tH_u} \big) - t \Tr\big(D_u
\lfrac{dH_u}{du}e^{-tH_u}\big).
\end{equation}
\end{lemma}

\begin{proof}
Let $K_u:= D_ue^{-tH_u/2}$ and $L_u:= e^{-tH_u/2}$. Then Theorem
\ref{DuhamelThm} implies that $L_u$ is a smooth family of
smoothing operators, and the same is true for $K_u$, since by
Lemma \ref{SmoothFamilyMap} the operators
$D_u\in\sB(L_{l+d}^2,L_l^2)$ depend smoothly on $u$ for every
$l\in\Z$. We now observe the following
\begin{enumerate}
\item $K_u$ and $L_u$ are smooth families of operators in
$\sB(L^2,L^2)$, and
\[
\lfrac d{du}(K_uL_u) = \lfrac{dD_u}{du} e^{-tH_u} + D_u
\lfrac{d}{du} e^{-tH_u},
\]
where we are using the pairing
\[
\sB(L^2_d,L^2)\times \sB(L^2,L^2_d) \to \sB(L^2,L^2),\quad
(S,T)\mapsto ST,
\]
to differentiate $D_ue^{-tH_u/2}$.
\item  $K_u$ and $L_u$ are trace-class operators. Expressing
the trace in terms of the kernels it follows from
\eqref{KernOpEst} that they depend continuously on $u$ with
respect to the trace norm.
\end{enumerate}
Using the H\"older inequality $|\Tr(ST)|\le \Tr|S|\, \|T\|_{0,0}$,
these observations imply that $\Tr\big(D_ue^{-tH_u}\big)$ is
differentiable. If $D_u$ and $H_u$ commute, we can use
\eqref{Duhamel} to compute that
\[
\begin{split}
\lfrac d{du} \Tr\big(D_ue^{-tH_u}\big) &=
\Tr\big(\lfrac{dD_u}{du} e^{-tH_u}\big) + \Tr\big(D_u
\lfrac{d}{du} e^{-tH_u}\big) \\
&= \Tr\big(\lfrac{dD_u}{du} e^{-tH_u}\big) -
\int_0^t\Tr\big(D_ue^{-(t-s)H_u}(\lfrac{d}{du}H_u)e^{-sH_u}\big)ds\\
&= \Tr\big(\lfrac{dD_u}{du} e^{-tH_u}\big) - \int_0^t
\Tr\big(D_u(\lfrac{d}{du}H_u)e^{-tH_u}\big)ds.\\ &=\Tr\big(\lfrac
{dD_u}{du} e^{-tH_u} \big) - t \Tr\big(D_u
\lfrac{dH_u}{du}e^{-tH_u}\big).\qedhere
\end{split}
\]
\end{proof}

\section{Parameter Dependent Eta Invariants}

As a consequence of the above discussion, we can now describe a
proof of the variation formula for the Eta invariant.

\subsection{Large Time Behaviour}

To study parameter dependent Eta invariants we let $H_u=D_u^2$,
where $D_u$ is a smooth family in $\sP^1_{s,e}(M,E)$. We can then
state the following parameter dependent version of Proposition
\ref{BasicTraceProp}, see also \cite[Lem. A.14]{DF}.

\begin{lemma}\label{Duhamel:2:Lem}
Under the assumptions of Theorem \ref{DuhamelThm}, the
one-parameter family $D_u e^{-tD_u^2}$ is a differentiable family
of trace-class operators, and
\begin{equation}\label{Duhamel:2}
\lfrac{d}{du} \Tr(D_u e^{-tD_u^2}) = \big(1+2t\lfrac d{dt}\big)
\Tr\big(\lfrac{dD_u}{du} e^{-tD_u^2}\big).
\end{equation}
Moreover, if $\dim(\ker D_u)$ is constant, then for $t_0>0$, there
exist constants $c$ and $C$ such that for $t\ge t_0$
\begin{equation}\label{UniformTraceEst}
\big| \Tr(D_u e^{-tD_u^2}) \big| \le Ce^{-ct},
\end{equation}
locally uniform in $u$.
\end{lemma}

\begin{proof}
By \eqref{DuhamelCor} we have
\[
\begin{split}
\lfrac{d}{du} \Tr(D_u e^{-tD_u^2}) &= \Tr\big(\lfrac {dD_u}{du}
e^{-tD_u^2} \big) - t \Tr\big(D_u
\lfrac{dD_u^2}{du}e^{-tD_u^2}\big)\\ &= \Tr\big(\lfrac {dD_u}{du}
e^{-tD_u^2} \big) - 2 t \Tr\big(\lfrac{dD_u}{du} D_u^2
e^{-tD_u^2}\big)\\
&= \Tr\big(\lfrac {dD_u}{du} e^{-tD_u^2} \big) + 2 t\lfrac d{dt}
\Tr\big(\lfrac{dD_u}{du}e^{-tD_u^2}\big).
\end{split}
\]
Proposition \ref{BasicTraceProp} shows that for $t_0>0$, we can
find constants $C(u)$ such that for $t\ge t_0$
\[
\big| \Tr(D_u e^{-tD_u^2}) \big| \le C(u)e^{-t\gl_0(u)/2},
\]
where $\gl_0(u)$ is the smallest non-zero eigenvalue of $D_u^2$.
The proof of Proposition \ref{BasicTraceProp} shows that $C(u)$
depends continuously on $u$ so that we can find $C$ with $C(u)\le
C$ for compact $u$ intervals. Moreover, the eigenvalues of
$D_u^2$ vary continuously. Hence if $\dim(\ker D_u)$ is constant,
the non-zero eigenvalues of $D_u^2$ have a uniform positive lower
bound on compact $u$ intervals. This proves the second assertion.
\end{proof}

The above lemma has the following consequence

\begin{cor}\label{EtaFctnDer}
If $\dim(\ker D_u)$ is constant, then
\begin{equation*}
\frac{d}{du} \int_{1}^\infty \Tr\big(D_ue^{-tD_u^2}\big)
t^{\frac{s-1}2}dt= -2 \Tr\big(\lfrac{dD_u}{du}e^{-D_u^2}\big) - s
\int_{1}^\infty \Tr\big(\lfrac{dD_u}{du}e^{-tD_u^2}\big)
t^{\frac{s-1}2}dt,
\end{equation*}
and both sides are holomorphic for all $s\in\C$.
\end{cor}

\begin{proof}
Let $P_u$ be the orthogonal projection onto $\ker D_u$. Since
$\dim(\ker D_u)$ is constant we know from Proposition
\ref{ProjGreenSmooth} that $P_u$ depends smoothly on $u$.
Moreover,
\[
\begin{split}
\lfrac{d}{du} D_u &=
\lfrac{d}{du}\big((\Id-P_u)D_u(\Id-P_u)\big)\\
&= - \lfrac{dPu}{du}D_u(\Id-P_u) +
(\Id-P_u)\lfrac{dD_u}{du}(\Id-P_u) - (\Id-P_u)D_u \lfrac{dPu}{du}
\end{split}
\]
From this one deduces that
$\Tr\big(\lfrac{dD_u}{du}e^{-tD_u^2}\big)$ satisfies an estimate
of the form \eqref{UniformTraceEst}. Thus, for fixed $T>1$ we can
use \eqref{Duhamel:2} to differentiate under the integral.
Integrating by parts, we find that
\[
\begin{split}
\frac d{du} \int_{1}^Tt^{\frac{s-1}2} \Tr\big(D_ue^{-tD_u^2}\big)
dt &= \int_{1}^T t^{\frac{s-1}2} \big(1+2t\lfrac
d{dt}\big)\Tr\big(\lfrac{dD_u}{du} e^{-tD_u^2}\big) dt\\
&= 2\Big[T^{\frac{s+1}2}\Tr\big(\lfrac{dD_u}{du} e^{-TD_u^2}\big)-
 \Tr\big(\lfrac{dD_u}{du} e^{- D_u^2} \big)\Big]\\
&\qquad\quad - s \int_{1}^T
\Tr\big(\lfrac{dD_u}{du}e^{-tD_u^2}\big) t^{\frac{s-1}2}dt.
\end{split}
\]
The uniform estimates on $\Tr\big(D_ue^{-tD_u^2}\big)$ and
$\Tr\big(\lfrac{dD_u}{du}e^{-tD_u^2}\big)$ show that both sides
are holomorphic for $s\in \C$ and allow us to take $T\to \infty$.
This proves the result.
\end{proof}

\subsection{Small Times and Meromorphic Extension}

To extend Corollary \ref{EtaFctnDer} to small times $t\in [0,1]$,
we need the following parameter dependent version of Theorem
\ref{HeatTrace}, see \cite[Lem. 1.9.3]{Gil}. A detailed proof can
also be found in \cite[Thm.'s 3.3 \& 3.4]{BO86}.

\begin{theorem}\label{ParamHeatTrace}
Let $A_u$ be an auxiliary smooth family of formally self-adjoint
operators of order $a$. There is an asymptotic expansion, locally
uniform in $u$, such that
\begin{equation}\label{ParamHeatTraceExp}
\Tr(A_u e^{-tD_u^2}) \sim \sum_{n=0}^\infty
t^{\frac{n-m-a}{2}}a_n(A_u,D_u) ,\quad\text{as }t\to 0.
\end{equation}
The $a_n(A_u,D_u)$ are integrals over quantities locally
computable from the total symbols of $A_u$ and $D_u$. Moreover,
the functions $a_n(A_u,D_u)$ are smooth in $u$ and
\eqref{ParamHeatTraceExp} can be differentiated term by term,
i.e.,
\[
\lfrac{d}{du} \Tr(A_u e^{-tD_u^2}) \sim \sum_{n=0}^\infty
t^{\frac{n-m-a}{2}}\lfrac {d}{du}a_n(A_u,D_u) ,\quad\text{as }t\to
0.
\]
\end{theorem}

\begin{prop}\label{EtaDerApp}
If $\dim(\ker D_u)$ is constant, the meromorphic extension of
$\eta(D_u,s)$ is continuously differentiable in $u$, and
\begin{equation*}
\lfrac{d}{du}\eta(D_u) = - \lfrac2{\sqrt\pi}\,
a_m(\lfrac{dD_u}{du},D_u^2),
\end{equation*}
where $a_m(\lfrac{dD_u}{du},D_u^2)$ is the constant term in the
asymptotic expansion of
\[
\sqrt t\,\Tr\big(\lfrac{dD_u}{du}e^{-tD_u^2}\big),\quad\text{as
}t\to 0.
\]
\end{prop}

\begin{proof}
According to Theorem \ref{ParamHeatTrace} we have asymptotic
expansions, which can be differentiated in $u$,
\[
\Tr\big(D_ue^{-tD_u^2}\big) \sim \sum_{n=0}^\infty
t^{\frac{n-m-1}{2}}a_n(u) ,\quad\text{as }t\to 0,
\]
and
\[
\Tr\big(\lfrac{dD_u}{du}e^{-tD_u^2}\big) \sim \sum_{n=0}^\infty
t^{\frac{n-m-1}{2}} a_n'(u) ,\quad\text{as }t\to 0.
\]
As remarked in Theorem \ref{HeatTrace} they can be differentiated
in $t$ as well so that \eqref{Duhamel:2} implies
\begin{equation}\label{AsympCoeffDiff}
\lfrac d{du} a_n(u) = (n-m)\,  a_n'(u).
\end{equation}
Let $N\in\N$ be fixed, and let
\[
r_N(t,u):= \Tr\big(D_ue^{-tD_u^2}\big) - \sum_{n=0}^N
t^{\frac{n-m-1}{2}}a_n(u),
\]
and
\[
r_N'(t,u):= \Tr\big(\lfrac{dD_u}{du}e^{-tD_u^2}\big) -
\sum_{n=0}^N t^{\frac{n-m-1}{2}}a_n'(u).
\]
Then $r_N$ and $r_N'$ satisfy estimates, locally uniform in $u$,
\begin{equation}\label{UniformRemainderBounds}
|r_N(t,u)|\le C t^N,\quad |r_N'(t,u)|\le Ct^N,\quad \text{as
}t\to 0.
\end{equation}
Moreover, \eqref{Duhamel:2} and \eqref{AsympCoeffDiff} imply that
\begin{equation}\label{RemainderDiff}
\lfrac d{du}r_N(t,u) = (1+2t\lfrac d{dt}) r_N'(t,u).
\end{equation}
Now for all $s\in\C$ with $\Re(s)> m- (N+1)$ and
$s\notin\setdef{m-n}{n\in\N}$
\[
\int_0^1t^{\frac{s-1}2} \Tr\big(D_ue^{-tD_u^2}\big) dt =
\sum_{n=0}^N \frac {2a_n(u)}{s+n-m} +
\int_0^1t^{\frac{s-1}2}r_N(t,u)dt.
\]
Using \eqref{UniformRemainderBounds} we can differentiate under
the integral to deduce from \eqref{AsympCoeffDiff} and
\eqref{RemainderDiff} that
\[
\begin{split}
\frac d{du}\int_0^1 &t^{\frac{s-1}2} \Tr\big(D_ue^{-tD_u^2}\big)
dt = \sum_{n=0}^N \frac {2(n-m)\, a_n'(u)}{s+n-m} +
\int_0^1t^{\frac{s-1}2}(1+2t\lfrac d{dt}) r_N'(t,u) dt \\
&= \sum_{n=0}^N \frac {2(n-m)\, a_n'(u)}{s+n-m} + 2r_N'(1,u) -
s\int_0^1 t^{\frac{s-1}2}r_N'(t,u) dt\\
&= 2\Tr\big(\lfrac{dD_u}{du}e^{-D_u^2}\big)+ \sum_{n=0}^N \frac
{-2s\, a_n'(u)}{s+n-m} - s\int_0^1 t^{\frac{s-1}2}r_N'(t,u) dt,
\end{split}
\]
since
\[
\frac{2(n-m)\, a_n'(u)}{s+n-m} -2a_n'(u) = \frac{-2s\,
a_n'(u)}{s+n-m}.
\]
On the other hand,
\[
\int_0^1 t^{\frac{s-1}2}
\Tr\big(\lfrac{dD_u}{du}e^{-D_u^2}\big)dt = \sum_{n=0}^N \frac
{2\, a_n'(u)}{s+n-m} +\int_0^1 t^{\frac{s-1}2}r_N'(t,u) dt.
\]
This shows that the meromorphic extension to $\C$ of $\int_0^1
t^{\frac{s-1}2}\Tr\big(D_ue^{-tD_u^2}\big) dt$ is continuously
differentiable in $u$ away from the poles, with derivative given
by the meromorphic extension of
\[
2\Tr\big(\lfrac{dD_u}{du}e^{-D_u^2}\big) -s \int_0^1
t^{\frac{s-1}2} \Tr\big(\lfrac{dD_u}{du}e^{-D_u^2}\big)dt.
\]
Since $\dim(\ker D_u)$ is constant we can apply Corollary
\ref{EtaFctnDer} to deduce that the meromorphic extension of the
Eta function $\eta(D_u,s)$ is continuously differentiable in $u$.
Moreover, for $N > m$ and a suitable function $h_N(u,s)$,
holomorphic for $\Re(s)> m-(N+1)$,
\[
\lfrac d{du}\eta(D_u,s) = \frac{-s}{\gG(\lfrac{s+1}2)}\Big(
\sum_{n=0}^N \frac {2\, a_n'(u)}{s+n-m} +h_N(u,s) \Big).
\]
In particular,
\[
\lfrac d{du}\eta(D_u) = - 2\, \gG(\lfrac 12)^{-1} a_m'(u)= -
\lfrac2{\sqrt\pi}\, a_m(\lfrac{dD_u}{du},D_u^2).\qedhere
\]
\end{proof}

\subsection{The Case of Varying Kernel Dimension}

If we want to drop the assumption on $\dim(\ker D_u)$, we have to
study the reduced $\xi$-invariant $[\xi(D_u)]\in \R/\Z$ as in
Definition \ref{EtaDef}.

%

\begin{prop}\label{RedEtaDerApp}
Let $D_u$ be a smooth family in $\sP^1_{s,e}(M,E)$. Then the
reduced $\xi$-function $[\xi(D_u)]\in \R/\Z$ is continuously
differentiable in $u$, and
\begin{equation*}
\lfrac d{du}[\xi(D_u)] = - \lfrac1{\sqrt\pi}\,
a_m(\lfrac{dD_u}{du},D_u^2),
\end{equation*}
where $a_m(\lfrac{dD_u}{du},D_u^2)$ is the constant term in the
asymptotic expansion of
\[
\sqrt t\,\Tr\big(\lfrac{dD_u}{du}e^{-tD_u^2}\big),\quad\text{as
}t\to 0.
\]
\end{prop}

The result is the same as \cite[Prop. A.17]{DF}. We include a
proof for completeness and sketch a few more details.

\begin{proof}
Choose $\gd$ small enough so that there exists $c\in (0,1)$ with
$c\notin\spec(|D_u|)$ for all $u\in(-\gd,\gd)$. Denote by
$\gl_u^i$ with $i=1,\ldots,i_0$ the finite number of eigenvalues
of $D_u$ with $|\gl_u^i|<c$, and let
\[
E_u(c):= \bigoplus_{i=1}^{i_0}\ker\big(D_u-\gl_u^i\big).
\]
Let $P_u(c)$ be the projection onto $E_u(c)$. According to Remark
\ref{SmoothDepRem} $P_u(c)$ is a smooth\footnote{Kato's selection
theorem (\cite[Sec. II.6]{K}), ensures that the eigenvalues
$\gl_u^i$ can be ordered in such a way that they are parametrized
by continuously differentiable functions. Nevertheless, the total
projection onto all eigenspaces spanned by the collection
$\gl_u^i$ depends smoothly on $u$.} family of finite rank
operators with smooth kernel. Thus
\[
D_u(c):= D_u\big(\Id-P_u(c)\big) + P_u(c)
\]
is a smooth perturbation of $D_u$ by finite rank operators with
smooth kernel. Note that for fixed $u$ we replace with 1 the
finite number of eigenvalues of $D_u$ which are of norm smaller
than $c$. Thus, the large eigenvalues of $D_u(c)$ are the same as
those of $D_u$. This implies that the Eta function is well-defined
for $\Re(s)>m$ and satisfies
\begin{equation}\label{RedEtaDerApp:0}
\eta\big(D_u(c),s\big) = \eta(D_u,s) + \dim\ker(D_u) +
\sum_{\gl_u^i\neq 0}\Big[ 1-
\lfrac{\sgn(\gl_u^i)}{|\gl_u^i|^s}\Big],\quad \Re(s)>m.
\end{equation}
Since the right hand side admits a meromorphic continuation to
$\C$, the same holds for $\eta\big(D_u(c),s\big)$. Moreover, $s=0$
is no pole, and the reduced $\xi$-invariant satisfies
\begin{equation}\label{RedEtaDerApp:0.5}
\big[\xi\big(D_u(c)\big)\big] = \big[\xi(D_u)\big].
\end{equation}
We now need to understand the variation of $\xi\big(D_u(c)\big)$.
Since $D_u(c)$ is invertible we will study $\eta\big(D_u(c)\big)$
instead. Note that we cannot directly apply Proposition
\ref{EtaDerApp} since $D_u(c)$ is in general not a family of
differential operators. However, we have already done the major
work and indicate the changes to be made:
\begin{enumerate}
\item Since $D_u(c)$ is a smooth perturbation of $D_u$ by symmetric
smoothing operators, the proof of Theorem \ref{DuhamelThm} goes
through verbatim, showing that $e^{-tD_u(c)^2}$ is a smooth family
of smoothing operators which satisfies \eqref{Duhamel}.
\item Then, as in Lemma \ref{Duhamel:2:Lem}, one finds that the
one-parameter family $D_u(c)e^{-tD_u(c)^2}$ is a differentiable
family of trace-class operators satisfying \eqref{Duhamel:2}.
Clearly, the uniform large-time estimate \eqref{UniformTraceEst}
and Corollary \ref{EtaFctnDer} also continue to hold.
\item Since $D_u(c)^2$ is of the form $D_u^2+K_u$, where $K_u$ is
a smooth family of smoothing operators which commutes with
$D_u^2$, we can apply Proposition \ref{TraceAsympPerturb} for
fixed $u$ to get asymptotic expansions
\begin{equation}\label{RedEtaDerApp:1}
\Tr\big(D_u(c)e^{-tD_u(c)^2}\big) \sim \sum_{n=0}^\infty
t^{\frac{n-m-1}{2}}a_n(u) ,\quad\text{as }t\to 0,
\end{equation}
and
\begin{equation}\label{RedEtaDerApp:2}
\Tr\Big(\big(\lfrac{d}{du}D_u(c)\big)e^{-tD_u(c)^2}\Big) \sim
\sum_{n=0}^\infty t^{\frac{n-m-1}{2}} a_n'(u) ,\quad\text{as
}t\to 0,
\end{equation}
where $a_m'(u)$ is equal to the constant term in the asymptotic
expansion as $t\to 0$ of $\sqrt
t\,\Tr\big(\lfrac{dD_u}{du}e^{-tD_u^2}\big)$. Moreover, one
deduces from Theorem \ref{ParamHeatTrace} and the proof of
Proposition \ref{TraceAsympPerturb} that \eqref{RedEtaDerApp:1}
and \eqref{RedEtaDerApp:2} are locally uniform in $u$ with
coefficients depending smoothly on $u$.
\end{enumerate}
Now the proof of Proposition \ref{EtaDerApp} carries over to the
situation at hand yielding that $\eta\big(D_u(c)\big)$ is
continuously differentiable in $u$ with
\[
\lfrac d{du} \eta\big(D_u(c)\big) = - \lfrac2{\sqrt\pi}\,
a_m(\lfrac{dD_u}{du},D_u^2).
\]
Since $\big[\xi(D_u)\big] = \lfrac 12 \eta\big(D_u(c)\big)$, the
proposition follows.
\end{proof}

From the proof of Proposition \ref{RedEtaDerApp} we can also
easily deduce the variation formula for the $\xi$-invariant, as
stated in Proposition \ref{EtaDiffSF}.

\begin{cor}\label{EtaDiffSFApp}
Let $D_u$ with $u\in [a,b]$ be a smooth one-parameter family of
operators in $\sP_{s,e}^1(M,E)$. Then
\begin{equation*}
\xi(D_b) - \xi (D_a) = \SF(D_u)_{u\in[a,b]} + \int_a^b \lfrac
d{du}[\xi(D_u)]du.
\end{equation*}
\end{cor}

\begin{proof}
Without loss of generality we may assume that $[a,b]=[-\gd,\gd]$
for $\gd$ as in the proof of Proposition \ref{RedEtaDerApp}.
Moreover, let $D_u(c)$ be the family of perturbed operators
defined there. Then the fundamental theorem of calculus shows that
\[
\eta(D_{\gd}(c)) - \eta (D_{-\gd}(c)) =  \int_{-\gd}^\gd \lfrac
d{du}\eta(D_u(c)) du.
\]
We deduce from analytic continuation of \eqref{RedEtaDerApp:0} to
$s=0$ that
\[
\begin{split}
\xi(D_{\gd})  - \xi (D_{-\gd}) &= \int_{-\gd}^\gd \lfrac
d{du}\xi(D_u(c)) du  -
\#\bigsetdef{i\in\{1,\ldots,i_0\}}{\gl_\gd^i
<0}\\&\qquad\qquad\qquad\qquad\qquad\quad +
\#\bigsetdef{i\in\{1,\ldots,i_0\}}{\gl_{-\gd}^i <0}\\
&= \int_{-\gd}^\gd \lfrac d{du}\xi(D_u(c)) du +
\SF(D_u)_{u\in[-\gd,\gd]}.
\end{split}
\]
Since $\xi(D_u(c))$ is differentiable,
\[
\lfrac d{du}\xi(D_u(c)) = \lfrac d{du} \big[\xi(D_u(c))\big].
\]
Now the equality $\big[\xi\big(D_u(c)\big)\big] =
\big[\xi(D_u)\big]$ from \eqref{RedEtaDerApp:0.5} proves the
result.
\end{proof}

\cleardoublepage

\markboth{\textsc{Bibliography}}{\textsc{Bibliography}}

\bibliography{main}

\providecommand{\bysame}{\leavevmode\hbox to3em{\hrulefill}\thinspace}
\providecommand{\MR}{\relax\ifhmode\unskip\space\fi MR }
\providecommand{\MRhref}[2]{%
  \href{http://www.ams.org/mathscinet-getitem?mr=#1}{#2}
}
\providecommand{\href}[2]{#2}
\begin{thebibliography}{10}

\bibitem{AK01}
J.A. \'Alvarez~L\'opez and Y.A. Kordyukov, \emph{{Long time behavior of
  leafwise heat flow for Riemannian foliations}}, Compos. Math. \textbf{125}
  (2001), no.~2, 129--153.

\bibitem{And}
J.E. Andersen, \emph{{The Witten invariant of finite order mapping tori I}},
  {to appear in J. Reine Angew. Math.}

\bibitem{Ati87}
M.F. Atiyah, \emph{{The logarithm of the Dedekind $\eta$-function}}, Math. Ann.
  \textbf{278} (1987), 335--380.

\bibitem{AB64}
M.F. Atiyah and R.~Bott, \emph{The index problem for manifolds with boundary},
  Differential analysis (Bombay, 1964), Oxford Univ. Press, London, 1964,
  pp.~175--186.

\bibitem{ABP}
M.F. Atiyah, R.~Bott, and V.K. Patodi, \emph{On the heat equation and the index
  theorem}, Invent. Math. \textbf{19} (1973), 279--330.

\bibitem{ADS}
M.F. Atiyah, H.~Donnelly, and I.M. Singer, \emph{Eta invariants, signature
  defects of cusps, and values of ${L}$-functions}, Ann. of Math. \textbf{118}
  (1983), 131--177.

\bibitem{APS1}
M.F. Atiyah, V.K. Patodi, and I.M. Singer, \emph{Spectral asymmetry and
  riemannian geometry {I}}, Math. Proc. Camb. Phil. Soc. \textbf{77} (1975),
  43--69.

\bibitem{APS2}
\bysame, \emph{Spectral asymmetry and riemannian geometry {II}}, Math. Proc.
  Camb. Phil. Soc. \textbf{78} (1975), 405--432.

\bibitem{APS3}
\bysame, \emph{Spectral asymmetry and riemannian geometry {III}}, Math. Proc.
  Camb. Phil. Soc. \textbf{79} (1976), 71--99.

\bibitem{AS4}
M.F. Atiyah and I.M. Singer, \emph{The index of elliptic operators {IV}}, Ann.
  of Math. \textbf{93} (1971), 119--138.

\bibitem{BeSa}
S.~Bechtluft-Sachs, \emph{{The computation of $\eta$-invariants on manifolds
  with free circle action}}, J. Funct. Anal. \textbf{174} (2000), no.~2,
  251--263.

\bibitem{Beck}
M.~Beck and S.~Robins, \emph{{Computing the continuous discretely.
  Integer-point enumeration in polyhedra}}, {Undergraduate Texts in
  Mathematics}, {Springer}, {Berlin etc.}, 2007.

\bibitem{BGV}
N.~Berline, E.~Getzler, and M.~Vergne, \emph{Heat kernels and {D}irac
  operators}, Grundlehren der mathematischen Wissenschaften, vol. 298,
  Springer, Berlin etc., 1992.

\bibitem{B86}
J.M. Bismut, \emph{The {A}tiyah-{S}inger index theorem for families fo {D}irac
  operators: {T}wo heat equation proofs}, Invent. Math. \textbf{83} (1986),
  91--151.

\bibitem{B98}
\bysame, \emph{Local index theory, eta invariants and holomorphic torsion: a
  survey}, Surveys in {D}ifferential {G}eometry, vol. III, Int. Press, 1998,
  pp.~1--76.

\bibitem{BC89}
J.M. Bismut and J.~Cheeger, \emph{$\eta$-invariants and their adiabatic
  limits}, J. Amer. Math. Soc. \textbf{2} (1989), 33--70.

\bibitem{BC92}
\bysame, \emph{Transgressed {E}uler classes of {SL}$(2n,\mathbb{Z})$, vector
  bundles, adiabatic limits of eta invariants and special values of
  {L}-functions}, Ann. scient. \'Ec. Norm. Sup. \textbf{25} (1992), 335--391.

\bibitem{BF2}
J.M. Bismut and D.S. Freed, \emph{The analysis of elliptic families
  {II}:{D}irac operators, eta invariants, and the holonomy theorem of
  {W}itten}, Commun. Math. Phys. \textbf{107} (1986), 103--163.

\bibitem{BL95}
J.M. Bismut and J.~Lott, \emph{Flat vector bundles, direct images and higher
  real analytic torsion}, J. Amer. Math. Soc. \textbf{8} (1995), no.~2,
  291--363.

\bibitem{BLP}
B.~Boo\ss-Bavnbek, M.~Lesch, and J.~Phillips, \emph{Unbounded {F}redholm
  operators and spectral flow}, Canad. J. Math. \textbf{57} (2005), no.~2,
  225--250.

\bibitem{BW}
B.~Boo\ss-Bavnbek and K.P. Wojciechowski, \emph{{Elliptic boundary problems for
  {D}irac operators}}, {Mathematics: Theory and Applications}, {Birkh\"auser},
  {Boston, MA}, 1993.

\bibitem{BT}
R.~Bott and L.W. Tu, \emph{{Differential forms in algebraic topology}},
  {Graduate Texts in Mathematics}, vol.~{82}, {Springer}, {Berlin etc.}, 1982.

\bibitem{BO86}
T.P. Branson and B.~{\O}rsted, \emph{Conformal indices of riemannian
  manifolds}, Comp. Math. \textbf{60} (1986), no.~3, 261--293.

\bibitem{BL98}
J.~Br{\"u}ning and M.~Lesch, \emph{Spectral theory of boundary value problems
  for {D}irac type operators}, Geometric aspects of partial differential
  equations (Roskilde, 1998), Contempory Mathematics, vol. 242, Amer. Math.
  Soc., Providence, RI, 1999, pp.~203--215.

\bibitem{BunMa}
U.~Bunke and X.~Ma, \emph{Index and secondary index theory for flat bundles
  with duality}, Gil Juan (ed.) et al., Aspects of boudary problems in analysis
  and geometry, Operator Theory: Advances and Applications, vol. 151,
  Birkh\"auser, Basel, 2004, pp.~265--341.

\bibitem{Che87}
J.~Cheeger, \emph{$\eta$-invariants, the adiabatic approximation and conical
  singularities}, J. Differ. Geom. \textbf{26} (1987), 175--221.

\bibitem{CHS}
S.S. Chern, F.~Hirzebruch, and J.-P. Serre, \emph{The index of a fibered
  manifold}, Proc. Amer. Math. Soc. \textbf{8} (1957), no.~3, 587--596.

\bibitem{CS}
S.S. Chern and J.~Simons, \emph{Characteristic forms and geometric invariants},
  Ann. of Math. \textbf{99} (1974), 48--69.

\bibitem{Cohen}
H.~Cohen, \emph{{Number theory. Volume II: Analytic and modern tools}},
  {Graduate Texts in Mathematics}, vol. {240}, {Springer}, {Berlin etc.}, 2007.

\bibitem{Dai91}
X.~Dai, \emph{Adiabatic limits, nonmultiplicativity of signature, and {L}eray
  spectral sequence}, J. Am. Math. Soc. \textbf{4} (1991), no.~2, 265--321.

\bibitem{DF}
X.~Dai and D.S. Freed, \emph{{$\eta$}-invariants and determinant lines}, J.
  Math. Phys. \textbf{35} (1994), no.~10, 5155--5194.

\bibitem{DaiZha95}
X.~Dai and W.~Zhang, \emph{{Circle bundles and the Kreck-Stolz invariant}},
  Trans. Am. Math. Soc. \textbf{347} (1995), no.~9, 3587--3593.

\bibitem{DavKir}
J.~Davis and P.~Kirk, \emph{Lecture notes in algebraic topology}, Graduate
  Studies in Mathematics, vol.~35, Amer. Math. Soc., Providence, RI, 2001.

\bibitem{Ded}
R.~Dedekind, \emph{{Erl\"{a}uterungen zu den {F}ragmenten {XXVIII}}}, {Riemann:
  Gesammelte mathematische Werke und wissenschaftlicher Nachlass}, {Teubner},
  {Leipzig}, 2nd ed., 1892, Reprint: Springer, Berlin etc., 1990, pp.~466--478.

\bibitem{Die59}
U.~Dieter, \emph{{Das Verhalten der Kleinschen Funktionen log $\sigma_{g,h}
  (\omega_1, \omega_2)$ gegen\"uber Modultransformationen und verallgemeinerte
  Dedekindsche Summen}}, J. Reine Angew. Math. \textbf{201} (1959), 37--70.

\bibitem{Dom98}
D.~Dom\'inguez, \emph{Finiteness and tenseness theorems for riemannian
  foliations}, Amer. J. Math. \textbf{120} (1998), 1237--1276.

\bibitem{DK}
S.K. Donaldson and P.B. Kronheimer, \emph{The geometry of four-manifolds},
  Oxford Mathematical Monographs, Clarendon Press, Oxford, 1990.

\bibitem{DosSal}
S.~Dostoglou and D.~Salamon, \emph{{Instanton homology and symplectic fixed
  points}}, {Symplectic geometry (Warwick, 1990)}, {Lond. Math. Soc. Lect. Note
  Ser.}, vol. 192, {Camb. Univ. Press}, {Cambridge}, 1993, pp.~57--93.

\bibitem{FL}
M.~Farber and J.P. Levine, \emph{Jumps of the eta-invariant}, Math. Z.
  \textbf{223} (1996), no.~2, 197--246.

\bibitem{FarKra}
H.M. Farkas and I.~Kra, \emph{{Riemann surfaces}}, 2nd ed., {Graduate Texts in
  Mathematics}, vol.~{71}, {Springer}, {Berlin etc.}, 1992.

\bibitem{F95}
R.~Forman, \emph{{Spectral sequences and adiabatic limits}}, Commun. Math.
  Phys. \textbf{168} (1995), no.~1, 57--116.

\bibitem{FreVaf}
D.S. Freed and C.~Vafa, \emph{{Global anomalies on orbifolds.}}, Commun. Math.
  Phys. \textbf{117} (1988), no.~2, 349--351.

\bibitem{Fuj}
H.~Fujita, \emph{{On the functoriality of the {C}hern-{S}imons line bundle and
  the determinant line bundle}}, Commun. Contemp. Math. \textbf{8} (2006),
  no.~6, 715--735.

\bibitem{Get83}
E.~Getzler, \emph{Pseudodifferential operators on supermanifolds and the index
  theorem}, Commun. Math. Phys. \textbf{92} (1983), 163--178.

\bibitem{Get86}
\bysame, \emph{A short proof of the {A}tiyah-{S}inger index theorem}, Topology
  \textbf{25} (1986), 111--117.

\bibitem{Get94}
\bysame, \emph{The {B}argmann representation, generalized {D}irac operators and
  the index of pseudodeffierential operators on {$\R^n$}}, {Maeda, Yoshiaki
  (ed.) et al., Symplectic geometry and quantization}, {Contempory
  Mathematics}, vol. 179, Amer. Math. Soc., Providence, RI, 1994, pp.~63--81.

\bibitem{Gil79}
P.~Gilkey, \emph{The residue of the local eta function at the origin}, Math.
  Ann. \textbf{240} (1979), 183--189.

\bibitem{Gil81}
\bysame, \emph{The residue of the global eta function at the origin}, Adv.
  Math. \textbf{40} (1981), 290--307.

\bibitem{Gil}
P.B. Gilkey, \emph{Invariance theory, the heat equation, and the
  {A}tiyah-{S}inger index theorem}, 2nd ed., CRC Press, Boca Raton, 1995.

\bibitem{GLP}
P.B. Gilkey, J.V. Leahy, and J.~Park, \emph{{Spinors, spectral geometry, and
  Riemannian submersions}}, Lecture Notes Series, vol.~40, {Seoul National
  Univ.}, Seoul, 1998, electronic reprint.

\bibitem{Goet}
S.~Goette, \emph{Equivariant {$\eta$}-invariants and {$\eta$}-forms}, J. Reine
  Angew. Math. \textbf{526} (2000), 181--236.

\bibitem{Hir73}
F.~Hirzebruch, \emph{Hilbert modular surfaces}, L'Enseign. Math. \textbf{19}
  (1973), 183--281.

\bibitem{Iva}
N.~Ivanov, \emph{{Mapping class groups}}, {Daverman, R. J. (ed.) et al.,
  Handbook of geometric topology}, Elsevier, Amsterdam, 2002, pp.~523--633.

\bibitem{Jef92}
L.~Jeffrey, \emph{{Chern-Simons-Witten invariants of lens spaces and torus
  bundles, and the semiclassical approximation}}, Commun. Math. Phys.
  \textbf{147} (1992), 563--604.

\bibitem{JosRS}
J.~Jost, \emph{{Compact Riemann surfaces. An introduction to contemporary
  mathematics}}, 3rd ed., Universitext, Springer, Berlin etc., 2006.

\bibitem{K}
T.~Kato, \emph{Perturbation theory for linear operators}, Classics in
  Mathematics, Springer, Berlin etc., 1995, Reprint of the 1980 edition.

\bibitem{KK90}
P.~Kirk and E.~Klassen, \emph{Chern-{S}imons invariants of 3-manifolds and
  respresentation spaces of knot groups}, Math. Ann. \textbf{287} (1990),
  343--367.

\bibitem{KK94}
\bysame, \emph{Computing spectral flow via cup products}, J. Differ. Geom.
  \textbf{40} (1994), no.~3, 505--562.

\bibitem{KK97}
\bysame, \emph{The spectral flow of the odd signature operator and higher
  {M}assey products}, Math. Proc. Camb. Phil. Soc. \textbf{121} (1997), no.~2,
  297--320.

\bibitem{KKR}
P.~Kirk, E.~Klassen, and D.~Ruberman, \emph{Splitting the spectral flow and the
  {A}lexander matrix}, Comment. Math. Helv. \textbf{69} (1994), no.~3,
  375--416.

\bibitem{KL04}
P.~Kirk and M.~Lesch, \emph{The {$\eta$}-invariant, {M}aslov index, and
  spectral flow for {D}irac-type operators on manifolds with boundary}, Forum
  Math. \textbf{16} (2004), no.~4, 553--629.

\bibitem{KN}
S.~Kobayashi and K.~Nomizu, \emph{Foundations of differential geometry, vol.
  {I+II}}, Wiley Classics Library, John Wiley \& Sons Inc., New York, 1996.

\bibitem{Kom}
M.~Komuro, \emph{{On Atiyah-Patodi-Singer $\eta$-invariant for $S^1$-bundles
  over Riemann surfaces}}, J. Fac. Sci., Univ. Tokyo, Sect. I A \textbf{30}
  (1984), 525--548.

\bibitem{Lang}
S.~Lang, \emph{{Elliptic functions}}, 2nd ed., {Graduate Texts in Mathematics},
  vol. {112}, {Springer}, {Berlin etc.}, 1987.

\bibitem{L05}
M.~Lesch, \emph{{The uniqueness of the spectral flow on spaces of unbounded
  self-adjoint {F}redholm operators}}, Contemp. Math. \textbf{366} (2005),
  193--224.

\bibitem{MazMel}
R.~Mazzeo and R.B. Melrose, \emph{{The adiabatic limit, {H}odge cohomology and
  {L}eray's spectral sequence for a fibration}}, J. Differ. Geom. \textbf{31}
  (1990), 185--213.

\bibitem{McC}
J.~McCleary, \emph{{A user's guide to spectral sequences}}, 2nd ed., Cambridge
  Studies in Advanced Mathematics, vol.~58, {Camb. Univ. Press}, Cambridge,
  2001.

\bibitem{McDSal}
D.~McDuff and D.~Salamon, \emph{{Introduction to symplectic topology}}, 2nd
  ed., Oxford Mathematical Monographs, {Oxford Univ. Press}, New York, 1998.

\bibitem{McKSin}
H.~McKean and I.M. Singer, \emph{{Curvature and the eigenforms of the
  {L}aplacian}}, J. Differ. Geom. \textbf{1} (1967), 43--69.

\bibitem{Mey57}
C.~Meyer, \emph{{\"{U}ber einige Anwendungen Dedekindscher Summen}}, J. Reine
  Angew. Math. \textbf{198} (1957), 143--203.

\bibitem{Mey61}
\bysame, \emph{{Bemerkungen zu den allgemeinen Dedekindschen Summen}}, J. Reine
  Angew. Math. \textbf{205} (1961), 186--196.

\bibitem{Mey73}
W.~Meyer, \emph{Die {S}ignatur von {F}l\"achenb\"undeln}, Math. Ann.
  \textbf{201} (1973), 239--264.

\bibitem{MS}
J.W. Milnor and J.D. Stasheff, \emph{Characteristic classes}, Annals of
  Mathematics Studies, no.~76, Princeton Univ. Press, Princeton, NJ, 1974.

\bibitem{Mo}
P.~Molino, \emph{Riemannian foliations}, Progress in Mathematics, no.~73,
  Birkh\"auser, Basel, Boston, 1988.

\bibitem{Mor}
T.~Morifuji, \emph{{The $\eta$-invariant of mapping tori with finite
  monodromies}}, Topology Appl. \textbf{75} (1997), 41--49.

\bibitem{Mos}
J.~Moser, \emph{{On the volume elements on a manifold}}, Trans. Am. Math. Soc.
  \textbf{120} (1965), 286--294.

\bibitem{Mue84}
W.~M{\"u}ller, \emph{Signature defects of cusps of {H}ilbert modular varieties
  and values of ${L}$ series at $s=1$}, J. Differ. Geom. \textbf{20} (1984),
  55--119.

\bibitem{Nic98}
L.I. Nicolaescu, \emph{{Adiabatic limits of the Seiberg-Witten equations on
  Seifert manifolds.}}, Commun. Anal. Geom. \textbf{6} (1998), no.~2, 331--392.

\bibitem{Nic99}
\bysame, \emph{Eta invariants of {D}irac operators on circle bundles over
  {R}iemann surfaces and virtual dimensions of finite energy {S}eiberg-{W}itten
  moduli spaces}, Israel J. Math. \textbf{114} (1999), 61--123.

\bibitem{Nic:SW}
\bysame, \emph{Notes on {S}eiberg-{W}itten theory}, Graduate Studies in
  Mathematics, vol.~28, Amer. Math. Soc., Providence, RI, 2000.

\bibitem{ON}
B.~O'Neill, \emph{{Semi-Riemannian geometry. {W}ith applications to
  relativity}}, {Pure and Applied Mathematics}, vol. {103}, {Academic Press},
  {New York etc.}, 1983.

\bibitem{Phi}
J.~Phillips, \emph{Self-adjoint {F}redholm operators and spectral flow}, Canad.
  Math. Bull. \textbf{39} (1996), no.~4, 460--467.

\bibitem{Qui85a}
D.~Quillen, \emph{Determinants of {C}auchy-{R}iemann operators over a {R}iemann
  surface}, Funktsional Anal. i Prilozhen \textbf{19} (1985), no.~1, 37--41.

\bibitem{Qui85b}
\bysame, \emph{Superconnections and the {C}hern character}, Topology
  \textbf{24} (1985), 89--95.

\bibitem{Rad64}
H.~Rademacher, \emph{{Some remarks on certain generalized Dedekind sums.}},
  Acta Arith. \textbf{9} (1964), 7--105.

\bibitem{Ram}
S.~Ramanan, \emph{{Global calculus}}, Graduate Studies in Mathematics, vol.~65,
  Amer. Math. Soc., Providence, RI, 2005.

\bibitem{RS71}
D.B. Ray and I.M. Singer, \emph{{R-torsion and the Laplacian on Riemannian
  manifolds.}}, Adv. Math. \textbf{7} (1971), 145--210.

\bibitem{RS73}
\bysame, \emph{Analytic torsion for complex manifolds}, Ann. of Math.
  \textbf{98} (1973), 154--177.

\bibitem{RobSal}
J.~Robbin and D.~Salamon, \emph{{The spectral flow and the {M}aslov index}},
  Bull. Lond. Math. Soc. \textbf{27} (1995), no.~1, 1--33.

\bibitem{Sal}
D.~Salamon, \emph{{Seiberg-Witten invariants of mapping tori, symplectic fixed
  points, and Lefschetz numbers}}, Turk. J. Math. \textbf{23} (1999), no.~1,
  117--143.

\bibitem{Sco83}
P.~Scott, \emph{{The geometries of 3-manifolds}}, Bull. Lond. Math. Soc.
  \textbf{15} (1983), 401--487.

\bibitem{Scz78}
R.~Sczech, \emph{{Ein einfacher Beweis der Transformationsformel f\"ur $\log
  \eta (z)$}}, Math. Ann. \textbf{237} (1978), 161--166.

\bibitem{Shi}
G.~Shimura, \emph{{Introduction to the arithmetic theory of automorphic
  functions}}, {Publications of the Mathematical Society of Japan}, vol.~11,
  Princeton Univ. Press, Princeton, NJ, 1971.

\bibitem{Tay:I}
M.E. Taylor, \emph{Partial differential equations {I}}, Appl. Math. Sci., vol.
  115, Springer, Berlin etc., 1996.

\bibitem{Ton}
P.~Tondeur, \emph{{Geometry of foliations}}, {Monographs in Mathematics},
  {Birk\-h\"auser}, {Basel}, 1997.

\bibitem{W}
R.O. Wells, \emph{{Differential analysis on complex manifolds}}, {Graduate
  Texts in Mathematics}, vol.~{65}, {Springer}, {Berlin etc.}, 1980.

\bibitem{Wit85}
E.~Witten, \emph{Global gravitational anomalies}, Commun. Math. Phys.
  \textbf{100} (1985), 197--229.

\bibitem{Zha94}
W.~Zhang, \emph{{Circle bundles, adiabatic limits of $\eta$-invariants and
  Rokhlin congruences}}, Ann. Inst. Fourier \textbf{44} (1994), no.~1,
  249--270.

\bibitem{Z}
\bysame, \emph{{Lectures on {C}hern-{W}eil theory and {W}itten deformations}},
  {Nankai Tracts in Mathematics}, vol.~{4}, {World Scientific}, {Singapore},
  2001.

\end{thebibliography}
\addcontentsline{toc}{chapter}{Bibliography}
\bibliographystyle{amsplain}



\end{document}